\documentclass[12pt,twoside,leqno]{amsart}
\usepackage{amssymb,amsbsy,amsmath,amsfonts,amssymb,amscd,epsfig,times,graphics,color,xypic}
\definecolor{blue}{cmyk}{1.,1.,0.,0.43}
\definecolor{red}{cmyk}{0.,1.,1.,0.43}
\definecolor{green}{cmyk}{1.,0.,1.,0.43}
\usepackage[latin1]{inputenc}
\newcommand{\A}{\mathbb{A}}\newcommand{\C}{\mathbb{C}}
\newcommand{\B}{\mathbb{B}}
\newcommand{\F}{\mathbb{F}}
\newcommand{\K}{\mathbb{K}}\newcommand{\LL}{\mathbb{L}}
\newcommand{\N}{\mathbb{N}}
\newcommand{\R}{\mathbb{R}}
\newcommand{\X}{\mathbb{X}}\newcommand{\Z}{\mathbb{Z}}
\def\dl{{[\![}}\def\dr{{]\!]}}
\newtheorem{problem}{Problem}
\newtheorem{definition}{Definition}
\newtheorem{theorem}{Theorem}
\newtheorem{proposition}{Proposition}
\newtheorem{lemma}{Lemma}
\newtheorem{corollary}{Corollary}
\newtheorem{example}{Example}
\newtheorem{openquestion}{Open question}
\newtheorem{openproblem}{Open problem}
\newcommand{\blue}{\textcolor{blue}}
\newcommand{\green}{\textcolor{green}}
\newcommand{\red}{\textcolor{red}}

\begin{document}


\title[
Holomorphic extensions
and removable singularities
]{
{\large
Holomorphic extension of CR functions,
\\
$\,$ 
\\
envelopes of holomorphy,
\\
$\,$ 
\\
and removable singularities}
}

\author{Jo\"el Merker and Egmont Porten}


\date{\number\year-\number\month-\number\day}

\maketitle

\begin{center}
\begin{minipage}[t]{12.5cm}
\baselineskip =0.35cm
{\tiny

\noindent 

D\'epartement de Math\'ematiques et Applications, UMR 8553
du CNRS, \'Ecole Normale
Sup\'erieure, 45 rue d'Ulm, F-75230 Paris Cedex 05, 
France.
\\
{\tt merker@dma.univ-mrs.fr
http://www.cmi.univ-mrs.fr/$\sim$merker/index.html}
\\
Department of Engineering, Physics and
Mathematics, Mid Sweden University, Campus Sundsvall,
S-85170 Sundsvall, Sweden
\ \
{\tt Egmont.Porten@miun.se}

\bigskip
\bigskip

\centerline{\bf Table of contents (7 main parts)}

\medskip

{\bf I.~Introduction \dotfill I.}

{\bf II.~Analytic vector field systems, formal CR mappings and local
CR automorphism groups \dotfill II.}

{\bf III.~Sussman's orbit theorem, locally integrable 
systems of vector fields and CR functions \dotfill III.}

{\bf IV.~Hilbert transform and Bishop's equation in H\"older spaces
\dotfill IV.}

{\bf V.~Holomorphic extension of CR functions \dotfill V.}

{\bf VI.~Removable singularities \dotfill VI.}

\smallskip

\hfill
{\footnotesize\tt [0+3+7+1+19+7 diagrams]}






}\end{minipage}
\end{center}

\bigskip\bigskip\bigskip\bigskip\bigskip
\bigskip\bigskip

\begin{center}
\fbox{\large \green{\bf IMRS} 
\green{\sf International Mathematics Research Surveys}} 

\smallskip
\fbox{\red{Volume {\bf 2006}, Article ID 28295, 287~pages}}

\smallskip
\fbox{\blue{www.hindawi.com/journals/imrs}}
\end{center}

\newpage

\centerline{{\Large \bf I:} {\large \bf Introduction}}

\subsection*{1.1.~CR extension theory}
In the past decades, remarkable progress has been accomplished towards
the understanding of compulsory extendability of holomorphic
functions, of CR functions and of differential forms. These phenomena,
whose exploration is still active in current research, originate from
the seminal Hartogs-Bochner extension theorem.

In local CR extension theory, the most satisfactory achievement was
the discovery that, on a smooth embedded generic submanifold $M
\subset \C^n$, there is a precise correspondence between {\sl CR
orbits}\, of $M$ and families of small {\sl Bishop discs}\, attached
to $M$. Such discs cover a substantial part of the polynomial hull of
$M$, and in most cases, this part may be shown to constitute a global
{\sl one-sided neighborhood} $\mathcal{ V }^\pm (M)$ of $M$, if $M$ is
a hypersurface, or else a {\sl wedgelike domain}\, $\mathcal{ W}$
attached to $M$, if $M$ has codimension $\geqslant 2$. A local
polynomial approximation theorem, or a CR version of the {\sl
Kontinuit\"atssatz}\, (continuity principle) assures that 
CR functions automatically extend holomorphically to such domains
$\mathcal{ W}$, which are in addition contained in the envelope of
holomorphy of arbitrarily thin neighborhoods of $M$ in $\C^n$.

Tr\'epreau in the hypersurface (1986) case and slightly after Tumanov
in arbitrary codimension (1988) established a nowadays celebrated
extension theorem: {\it if $M \subset \C^n$ is a sufficiently smooth
($\mathcal{ C}^2$ or $\mathcal{ C}^{ 2, \alpha }$ suffices) generic
submanifold, then at every point $p\in M$ whose local CR orbit
$\mathcal{ O }_{ CR}^{ loc} (M, p)$ has maximal dimension equal to
$\dim M$, there exists a local wedge $\mathcal{ W}_p$ of edge $M$ at
$p$ to which continuous CR functions extend holomorphically}. Several
reconstructions and applications of this groundbreaking result,
together with surveys about the local Bishop equation have already
appeared in the literature.

Propagational aspects of CR extension theory are less known by
contemporary experts of several complex variables, but they lie deeper
in the theory. Using {\sc fbi} transform and concepts of microlocal
analysis, Tr\'epreau showed in 1990 that holomorphic extension to a
wedge propagates along curves whose velocity vector is
complex-tangential to $M$. His conjecture that extension to a wedge
should hold at every point of a generic submanifold $M \subset \C^n$
consisting of a single global CR orbit has been answered independently
by J\"oricke and by the first author in 1994, using tools introduced
previously by Tumanov. To the knowledge of the two authors, there is
no survey of these global aspects in the literature.

The first main objective of the present survey is to expose the
techniques underlying these results in a comprehensive and unified
way, emphasizing propagational aspects of embedded CR geometry and
discussing optimal smoothness assumptions. Thus, topics that are
necessary to build the theory from scratch will be selected and
accompanied with thorough proofs, whereas other results that are
nevertheless central in CR geometry will be presented in concise
survey style, without any proof.

The theory of CR extension by means of analytic discs combines various
concepts emanating mainly from three (wide) mathematical areas: {\sl
Harmonic analysis}, {\sl Partial differential equations} and {\sl
Complex analysis in several variables}. As the project evolved, we
felt the necessity of being conceptional, extensive and systematic in
the restitution of (semi)known results, so that various contributions
to the subject would recover a certain coherence and a certain
unity. With the objective of adressing to a younger audience, we
decided to adopt a style accessible to doctoral candidates working on
a dissertation. Parts~III, IV and~V present elementarily general CR
extension theory. Also, most sections of the text may be read
independently by experts, as quanta of mathematical information.

\subsection*{1.2.~Concise presentation of the contents}
The survey text is organized in six main parts. Actually, the
present brief introduction constitutes the first and shortest one.
Although the reader will find a ``conceptional summary-introduction''
at the beginning of each part, a few descriptive words explaining some
of our options governing the reconstruction of CR extension theory
(Parts~III, IV and V) are welcome.

The next Part~II is independent of the others and can be skipped in a
first reading. It opens the text, because it is concerned with
propagational aspects of analytic CR structures, better understood
than the smooth ones.

\smallskip\noindent$\bullet$
In Part~III, exclusively concerned with the smooth
category, Sussmann's orbit theorem and its consequences are first
explained in length. Involutive structures and embedded CR manifolds,
together with their elementary properties, are introduced. Structural
properties of finite type structures, of CR orbits and of CR functions
are presented without proofs. As a collection of background material,
this part should be consulted first.

\smallskip\noindent$\bullet$
In Part~IV, fundamental results about singular integral operators in
the complex plane are first surveyed. Explicit estimates of the norms
of the Cauchy, of the Schwarz and of the Hilbert transforms in the
H\"older spaces $\mathcal{ C}^{ \kappa, \alpha }$ are provided. They
are useful to reconstruct the main Theorem~3.7(IV), due to Tumanov,
which asserts the existence of unique solutions to a parametrized
Bishop-type equation with an optimal loss of smoothness with respect
to parameters. Following Bishop's constructive philosophy, the
smallness of the constants insuring existence is precised explicitly,
thanks to sharp norm inequalities in H\"older spaces. This part is
meant to introduce interested readers to further reading of Tumanov's
recent works about extremal (pseudoholomorphic) 
discs in higher codimension.

\smallskip\noindent$\bullet$
In Part~V, CR extension theory is first discussed in the hypersurface
case. A simplified proof of wedge extendability that treats both
locally minimal and globally minimal generic submanifolds on the same
footing constitutes the main Theorem~4.12(V): {\it If $M$ is a
globally minimal $\mathcal{ C}^{ 2, \alpha}$ $(0 < \alpha < 1)$
generic submanifold of $\C^n$ of codimension $\geqslant 1$ and of CR
dimension $\geqslant 1$, there exists a wedgelike domain $\mathcal{
W}$ attached to $M$ such that every continuous CR function $f\in
\mathcal{ C}_{ CR}^0 (M)$ possesses a holomorphic extension $F \in
\mathcal{ O} (\mathcal{ W}) \cap \mathcal{ C}^0 (M\cup \mathcal{ W}) $
with $F\vert_M = f$}. The figures are intended to share the geometric
insight of experts in higher codimensional geometry. 



\smallskip
In fact, throughout the text, diagrams (33 in sum)
facilitating readability (especially of Part~V) are included. Selected
open questions and open problems (16 in sum) are formulated. They are
systematically inserted in the right place of the architecture. The
sign ``[$*$]'' added after one or several bibliographical references
in a statement (Problem, Definition, Theorem, Proposition, Lemma,
Corollary, Example, Open question and Open problem, {\it e.g.}
Theorem~1.11(I)) indicates that, compared to the existing literature,
a slight modification or a slight improvement has been brought by the
two authors. Statements containing no bibliographical reference are
original and appear here for the first time.

\smallskip
We apologize for having not treated some central topics of CR geometry
that also involve propagation of holomorphicity, {\it exempli
gratia}\, the geometric reflection principle, in the sense of Pinchuk,
Webster, Diederich, Forn{\ae}ss, Shafikov and Verma. By lack of space,
embeddability of abstract CR structures, polynomial hulls, Bishop
discs growing at elliptic complex tangencies, filling by Levi-flat
surfaces, Riemann-Hilbert boundary value problems, complex Plateau
problem in K\"ahler manifolds, partial indices of analytic discs,
pseudoholomorphic discs, {\it etc.} are not reviewed either. 
Certainly, better experts will fill this gap in the near future.

\smallskip
To conclude this introductory presentation, we believe that, although
uneasy to build, surveys and syntheses play a decisive r\^ole in the
evolution of mathematical subjects. For instance, in the last decades,
the remarkable development of $\overline{ \partial }$ techniques and
of $L^2$ estimates has been regularly accompanied by monographs and
panoramas, some of which became landmarks in the field. Certainly,
the (local) method of analytic discs deserves to be known by a
wider audience; in fact, its main contributors have brought it to the
degree of achievement that opened the way to the present survey.

\subsection*{ 1.3.~Further readings}
Using the tools exposed and reconstructed in this survey, the research
article~\cite{ mp2006a} studies removable singularities on CR
manifolds of CR dimension equal to 1 and solves a delicate remaining
open problem in the field ({\it see} the Introduction there for
motivations). Recently also, the authors built in~\cite{ mp2006c} a
new, rigorous proof of the classical Hartogs extension theorem which
relies only on the basic local Levi argument along analytic discs,
hence avoids both multidimensional integral representation formulas
and the Serre-Ehrenpreis argument about vanishing of $\overline{
\partial}$ cohomology with compact support.

\newpage

\begin{center}
{\Large\bf II:~Analytic vector field systems
\\ 
and formal CR mappings}
\end{center}

\bigskip\bigskip\bigskip

\begin{center}
\begin{minipage}[t]{11cm}
\baselineskip =0.35cm
{\scriptsize

\centerline{\bf Table of contents}

\smallskip

{\bf 1.~Analytic vector field systems and Nagano's theorem 
\dotfill 7.}

{\bf 2.~Analytic CR manifolds, Segre chains and minimality 
\dotfill 19.}

{\bf 3.~Formal CR mappings, jets of Segre varieties and CR reflection
mapping \dotfill 28.}

\smallskip

\hfill
{\footnotesize\tt [3 diagrams]}

}\end{minipage}
\end{center}

\bigskip
\bigskip

{\small


According to the theorem of Frobenius, a system $\LL$ of local vector
fields having real or complex analytic coefficients enjoys the
integral manifolds property, provided it is closed under Lie bracket.
If the Lie brackets exceed $\LL$, considering the smallest analytic
system $\LL^{ \rm lie}$ containing $\LL$ which is closed under Lie
bracket, Nagano showed that through every point, there passes a
submanifold whose tangent space is spanned by $\LL^{ \rm
lie}$. Without considering Lie brackets, these submanifolds may also
be constructed by means of compositions of local flows of elements of
$\LL$. Such a construction has applications in real analytic
Cauchy-Riemann geometry, in the reflection principle, in formal CR
mappings, in analytic hypoellipticity theorems and in the problem of
local solvability and of local uniqueness for systems of first order
linear partial differential operators (Part~III).

For a generic set of $r\geqslant 2$ vector fields having analytic
coefficients, $\LL^{ \rm lie}$ has maximal rank equal to the dimension
of the ambient space.

The extrinsic complexification $\mathcal{ M}$ of a real algebraic or
analytic Cauchy-Riemann submanifold $M$ of $\C^n$ carries two pairs of
intrinsic foliations, obtained by complexifying the classical Segre
varieties together with their conjugates. The Nagano leaves of this
pair of foliations coincide with the extrinsic complexifications of
local CR orbits. If $M$ is (Nash) algebraic, its CR orbits are
algebraic too, because they are projections of complexified algebraic
Nagano leaves.

A complexified formal CR mapping between two complexified generic
submanifolds must respect the two pairs of intrinsic foliations that
lie in the source and in the target. This constraint imposes strong
rigidity properties, as for instance: convergence, analyticity or
algebraicity of the formal CR mapping, according to the smoothness of
the target and of the source. There is a combinatorics of various
nondegeneracy conditions that entail versions of the so-called {\sl
analytic reflection principle}. The concept of {\sl CR reflection
mapping} provides a unified synthesis of recent results of the
literature.

}

\section*{ \S1.~Analytic vector field systems and Nagano's theorem}

\subsection*{ 1.1.~Formal, analytic and (Nash) algebraic power series} 
Let $n\in \N$ with $n\geqslant 1$ and let ${\sf x} = ({\sf x}_1, \dots,
{\sf x}_n) \in \K^n$, where $\K = \R$ or $\C$. Let $\K \dl {\sf x}
\dr$ be the ring of formal power series in $( {\sf x}_1, \dots,\, {\sf
x}_n)$. An element $\varphi( {\sf x}) \in \K \dl {\sf x } \dr$ writes
$\varphi( {\sf x }) = \sum_{ \alpha \in \N^n}\, \varphi_\alpha \, {\sf
x }^\alpha$, with ${\sf x}^\alpha := {\sf x }_1^{ \alpha_1 } \cdots
{\sf x }_n^{ \alpha_n}$ and with $\varphi_\alpha \in \K$, for every
multiindex $\alpha := (\alpha_1, \dots, \alpha_n )\in \N^n$. We put
$\vert \alpha \vert := \alpha_1 + \cdots + \alpha_n$.

On the vector space $\K^n$, we choose once for all the {\sl maximum
norm} $\vert {\sf x} \vert := \max_{ 1\leqslant i\leqslant n} \vert
{\sf x}_i \vert$ and, for any ``radius'' $\rho_1$ satisfying $0 <
\rho_1 \leqslant \infty$, we define the {\sl open cube}
\[
\square_{\rho_1}^n 
:= 
\{ 
{\sf x} \in \K^n : \vert x \vert < \rho_1 \}
\] 
as a fundamental, concrete
open set. For $\rho_1 = \infty$, we identify of course
$\square_\infty^n$ with $\K^n$.

If the coefficients $\varphi_\alpha$ satisfy a Cauchy estimate of the
form $\vert \varphi_\alpha \vert \leqslant C \rho_2^{ -\vert \alpha
\vert}$, $C >0$, for every $\rho_2$ satisfying $0 < \rho_2 < \rho_1$,
the formal power series is $\K$-analytic ($\mathcal{ C}^\omega$) in
$\square_{\rho_1}^n$. It then defines a true point map $\varphi:
\square_{ \rho_1}^n \to \K$. Such a $\K$-analytic function $\varphi$
is called (Nash) $\K$-{\sl algebraic} if there exists a nonzero
polynomial $P({\sf X}, \Phi ) \in\K[ {\sf X}, \Phi ]$ in $(n+1)$
variables such that the relation $P({\sf x}, \varphi ( {\sf x}))
\equiv 0$ holds in $\K \dl {\sf x} \dr$, hence for all ${\sf x}$ in
$\square_{ \rho_1}^n$. The category of $\K$-algebraic functions and
maps is stable under elementary algebraic operations, under
differentiation and under composition. Implicit solutions of
$\K$-algebraic equations are $\K$-algebraic
(\cite{ ber1999}).

\subsection*{ 1.2.~Analytic vector field systems
and their integral manifolds} Let 
\[
\LL^0 
:= 
\{ L_a \}_{ 1\leqslant a
\leqslant r}, 
\ \ \
r\in \N, 
\ \ \ r\geqslant 1, 
\]
be a finite set of vector
fields $L_a = \sum_{ i=1 }^n\, \varphi_{a,i} ({\sf x}) \, \frac{
\partial }{ \partial {\sf x}_i}$, whose coefficients $\varphi_{ a, i}$
are algebraic or analytic in $\square_{ \rho_1}^n$. Let $\A_{ \rho_1}$
denote the ring of algebraic or analytic functions in $\square_{
\rho_1}^n$. The set of linear combinations of elements of $\LL^0$ with
coefficients in $\A_{\rho_1}$ will be denoted by $\LL$ (or $\LL^1$)
and will be called the $\A_{\rho_1}$-{\sl linear hull of} $\LL^0$.

If $p$ is a point of $\square_{ \rho_1}^n$, denote by $L_a (p)$ the
vector $\sum_{ i=1 }^n\, \varphi_{a,i} (p) \, \frac{ \partial }{
\partial {\sf x}_i} \big\vert_p$. It is an element of $T_p \square_{
\rho_1}^n \simeq \K^n$. Define the linear subspace
\[
\LL(p):= 
{\rm Span}_\K\,\{L_a(p):1\leqslant a\leqslant r\}
=
\{L(p):L\in\LL\}.
\]
No constancy of dimension, no linear independency assumption are
made.

\def\theproblem{1.3}\begin{problem}
Find local submanifolds $\Lambda$ passing through the origin
satisfying $T_q \Lambda \supset \LL (q)$ for every $q \in \Lambda$.
\end{problem}

\noindent
By the theorem of Frobenius (\cite{ stk2000}; original
article: \cite{ fr1877}), if the $L_a$ are linearly independent at
every point of $\square_{ \rho_1}^n$ and if the Lie brackets $\left[
L_a, L_{ a'} \right]$ belong to $\LL$, for all $a, a' = 1, \dots,
r$, then $\square_{ \rho_1}^n$ is foliated by $r$-dimensional
submanifolds $N$ satisfying $T_q N = \LL (q)$ for every $q\in N$.

\def\thelemma{1.4}\begin{lemma}
If there exists a local submanifold $\Lambda$ passing through the
origin and satisfying $T_q \Lambda \supset \LL (q)$ for every $q \in
\Lambda$, then for every two vector fields $L, L' \in \LL$,
the restriction to $\Lambda$ of the Lie bracket $\left[ L, L'
\right]$ is tangent to $\Lambda$.
\end{lemma}

Accordingly, set $\LL^1 := \LL$ and for $k\geqslant 2$, define $\LL^k$
to be the $\A_{ \rho_1}$-linear hull of $\LL^{ k-1} + \left[ \LL^1,
\LL^{k-1} \right]$. Concretely, $\LL^k$ is generated by $\A_{
\rho_1}$-linear combinations of iterated Lie brackets $\left[ L_1,
\left[ L_2, \dots, \left[ L_{ k-1}, L_k\right] \dots \right] \right]$,
where $L_1, L_2, \dots, L_{ k-1}, L_k \in \LL^1$. The Jacobi identity
insures (by induction) that $\left[ \LL^{k_1}, \LL^{ k_2} \right]
\subset \LL^{ k_1 + k_2}$. Define then $\LL^{ \rm lie} :=
\cup_{k\geqslant 1} \, \LL^k$. Clearly, $\left[ L, L' \right] \in
\LL^{ \rm lie}$, for every two vector fields $L, L' \in \LL^{\rm
lie}$.

\def\thetheorem{1.5}\begin{theorem}
{\rm ({\sc Nagano} \cite{ na1966, trv1992, ber1999,
bch2005})} There exists a unique local $\K$-analytic submanifold
$\Lambda$ of $\K^n$ passing through the origin which satisfies $\LL(q)
\subset T_q \Lambda = \LL^{\rm lie} (q)$, for every $q \in \Lambda$.
\end{theorem}

A discussion about what happens in the algebraic category is postponed
to \S1.12. In Frobenius' theorem, $\LL^{ \rm lie} = \LL$ and the
dimension of $\LL^{\rm lie} (p)$ is constant. In the above theorem,
the dimension of $\LL^{\rm lie} (q)$ is constant for $q$ belonging to
$\Lambda$, but in general, not constant for $p\in \square_{ \rho_1
}^n$, the function $p \mapsto \dim_\K \LL (p)$ being lower
semi-continuous.

Nagano's theorem is stated at the origin; it also holds at every point
$p \in \square_{ \rho_1}^n$. The associated local submanifold
$\Lambda_p$ passing through $p$ with the property that $T_q\Lambda =
\LL^{\rm lie} (q)$ for every $q\in \Lambda_p$ is called a 
(local) {\sl Nagano leaf}.

In the $\mathcal{ C }^\infty$ category, the consideration of $\LL^{\rm
lie}$ is insufficient. Part~III handles smooth vector field systems,
providing a different answer to the search of similar submanifolds
$\Lambda_p$.

\def\theexample{1.6}\begin{example}
{\rm
In $\R^2$, take $\LL^0 = \{ L_1, L_2\}$, where $L_1= \partial_{ {\sf
x}_1}$ and $L_2 = e^{ -1/ {\sf x}_1^2 } \, \partial_{ {\sf x
}_2}$. Then $\LL^{ \rm lie} ( 0)$ is the line $\R \partial_{ \sf x_1}
\vert_0$, while $\LL^{ \rm lie} (p) = \R \partial_{ {\sf x }_1}
\vert_p + \R\partial_{ {\sf x }_2} \vert_p$ at every point $p \not \in
\R \times \{ 0\}$. Hence, there cannot exist a $\mathcal{ C}^\infty$
curve $\Lambda$ passing through $0$ with $T_0 \Lambda = \R
\partial_{\sf x_1} \vert_0$ and $T_q \Lambda = \LL^{ \rm lie} (q)$ for
every $q\in \Lambda$.
}
\end{example}

\smallskip
\noindent
{\it Proof of Theorem~1.5.} (May be skipped in a first reading.) If
$n=1$, the statement is clear, depending on whether or not all vector
fields in $\LL^{\rm lie}$ vanish at the origin. Let $n\geqslant 2$. Since
$\LL(q) \subset \LL^{\rm lie} (q)$, the condition $T_q \Lambda =
\LL^{\rm lie} (q)$ implies the inclusion $\LL (q) \subset T_q
\Lambda$. Replacing $\LL$ by $\LL^{\rm lie}$ if necessary, we may
therefore assume that $\LL^{\rm lie} = \LL$ and we then have to prove
the existence of $\Lambda$ with $T_q \Lambda = \LL^{\rm lie} (q) = \LL
(q)$, for every $q\in \Lambda$.

We reason by induction, supposing that, in dimension $(n-1)$, for
every $\A_{\rho_1}$-linear system $\LL' =(\LL' )^{ \rm lie}$ of
vector fields locally defined in a neighborhood of the origin in $\K^{
n-1}$, there exists a local $\K$-analytic submanifold $\Lambda'$
passing through the origin and satisfying $T_{q'} \Lambda' = \LL'
(q')$, for every $q'\in \Lambda'$.

If all vector fields in $\LL = \LL^{\rm lie}$ vanish at $0$, we are
done, trivially. Thus, assume there exists $L_1 \in \LL$ with $L_1(0)
\neq 0$. After local straightening, $L_1 = \partial_{ {\sf x
}_1}$. Every $L \in \LL$ writes uniquely $L = a ({\sf x}) \,
\partial_{ {\sf x}_1} + \widetilde{ L}$, for some $a ({\sf x}) \in \K
\{ {\sf x} \}$, with $\widetilde{ L} = \sum_{ 2 \leqslant i\leqslant n}\, a_i(
{\sf x} )\, \partial_{ {\sf x}_i}$. Introduce the space $\widetilde{
\LL} := \{ \widetilde{ L} : L \in \LL\}$ of such vector fields. As
$\partial_{{\sf x}_1}$ belongs to $\LL$ and as $\LL$ is
$\A_{\rho_1}$-linear, $\widetilde{ L} = L - a ({\sf x}) \,
\partial_{{\sf x}_1}$ belongs to $\LL$. Since $[\LL, \LL] \subset
\LL$, we have $\big[ \widetilde{ \LL}, \widetilde{ \LL} \big] \subset
\LL$. On the other hand, we observe that the Lie bracket between two
elements of $\widetilde{ \LL}$ does not involve $\partial_{{\sf
x}_1}$:
\def\theequation{1.7}\begin{equation}
\aligned
{}
\big[
\widetilde{ L}_1, \widetilde{ L}_2 
\big]
&
=
\Big[
\sum_{2\leqslant i_2\leqslant n}\,
a_{i_2}^1\, \partial_{{\sf x}_{i_2}}, \
\sum_{2\leqslant i_1\leqslant n}\,
a_{i_1}^2\, \partial_{{\sf x}_{i_1}}
\Big]
\\
&
=
\sum_{2\leqslant i_1\leqslant n}\,
\Big(
\sum_{2\leqslant i_2\leqslant n}\,
\Big[
a_{i_2}^1\,
\frac{\partial a_{i_1}^2}{\partial x_{i_2}}
-
a_{i_2}^2\,
\frac{\partial a_{i_1}^1}{\partial x_{i_2}}
\Big]
\Big)
\partial_{{\sf x}_{i_1}}.
\endaligned
\end{equation}
We deduce that $\big[ \widetilde{ \LL}, \widetilde{ \LL} \big] \subset
\widetilde{ \LL }$. In other words, $\widetilde{ \LL}^{\rm lie} =
\widetilde{ \LL}$. Next, we define the restriction
\[
\LL' := \big\{ L' = \widetilde{ L}
\big\vert_{\{{\sf x}_1 = 0\}} :
\widetilde{ L} \in \widetilde{ \LL}
\big\},
\]
and we claim that $(\LL')^{\rm lie} = \LL'$ also holds true. Indeed,
restricting~\thetag{ 1.7} above to $\{ {\sf x}_1 = 0\}$, we observe
that
\[
\big[ \widetilde{ L}_1\big\vert_{\{{\sf x}_1 = 0\}}, \widetilde{
L}_2\big\vert_{\{{\sf x}_1 = 0\}} \big] = \big[ \widetilde{ L}_1,
\widetilde{ L}_2 \big] \big\vert_{\{{\sf x}_1 = 0\}},
\]
since neither $\widetilde{ L}_1$ nor $\widetilde{ L}_2$ involves
$\partial_{ {\sf x}_1}$. This shows that $\big[ \LL', \LL' \big]
\subset \LL'$, as claimed.

Since $(\LL' )^{\rm lie} = \LL'$, the induction assumption applies:
there exists a local $\K$-analytic submanifold $\Lambda'$ of
$\K^{ n-1}$ passing through the origin such that $T_{q'} \Lambda' =
\LL'( q')$, for every point $q'\in \Lambda'$. Let $d$ denote its
codimension. If $d=0$, {\it i.e.} if $\Lambda'$ coincides with an open
neighborhood of the origin in $\K^{n-1}$, it suffices to chose for
$\Lambda$ an open neighborhood of the origin in $\K^n$. Assuming
$d \geqslant 1$, we split the coordinates ${\sf x} = ({\sf x }_1, {\sf x}')
\in \K \times \K^{n-1}$ and we let $\rho_j( {\sf x}')=0$, $j=1,
\dots, d$, denote local $\K$-analytic defining equations for
$\Lambda'$. We claim that it suffices to choose for $\Lambda$ the
local submanifold of $\K^n$ with the same equations, hence
having the same codimension.

Indeed, since these equations are independent of ${\sf x}_1$, it is
first of all clear that the vector field $\partial_{ {\sf x }_1} \in
\LL$ is tangent to $\Lambda$. To conclude that every $L = a\,
\partial_{ {\sf x}_1} + \widetilde{ L} \in \LL$ is tangent to $\Lambda$,
we thus have to prove that every $\widetilde{ L} \in \widetilde{ \LL}$
is tangent to $\Lambda$.

Let $\widetilde{ L} = \sum_{ 2\leqslant i\leqslant n}\, a_i ({\sf x}, {\sf x}')
\, \partial_{{\sf x }_i} \in \widetilde{ \LL}$. As a preliminary
observation:
\[
({\rm ad}\, \partial_{{\sf x}_1})
\widetilde{L}
:=
\big[
\partial_{{\sf x}_1},\widetilde{L}
\big]
=
\sum_{2\leqslant i\leqslant n}\,
\frac{\partial a_i}{\partial {\sf x}_1}
({\sf x}_1,{\sf x}')\,
\frac{\partial}{\partial {\sf x}_i},
\]
and more generally, for $\ell \in \N$ arbitrary:
\[
({\rm ad}\, \partial_{{\sf x}_1})^\ell
\widetilde{L}
=
\sum_{2\leqslant i\leqslant n}\,
\frac{\partial^\ell a_i}{\partial {\sf x}_1^\ell}
({\sf x}_1,{\sf x}')\,
\frac{\partial}{\partial {\sf x}_i}.
\]
Since $\LL$ is a Lie algebra, we have $({\rm ad}\, \partial_{ {\sf
x}_1 })^\ell \widetilde{ L} \in \LL$. Since $({\rm ad}\, \partial_{
{\sf x}_1 })^\ell \widetilde{ L}$ does not involve $\partial_{{\sf
x}_1}$, according to its expression above, it belongs in fact to
$\widetilde{ \LL}$. Also, after restriction $({\rm ad}\, \partial_{
{\sf x}_1 })^\ell \widetilde{ L} \big \vert_{ {\sf x}_1 = 0} \in
\LL'$. By assumption, $\LL '$ is tangent to $\Lambda'$. We deduce
that, for every $\ell \in \N$, the vector field
\[
L_\ell':=({\rm ad}\, \partial_{{\sf x}_1})^\ell
\widetilde{L}\big\vert_{{\sf x}_1=0}
=
\sum_{2\leqslant i\leqslant n}\,
\frac{\partial^\ell a_i}{\partial {\sf x}_1^\ell}
(0,{\sf x}')\,
\frac{\partial}{\partial {\sf x}_i}
\]
is tangent to $\Lambda'$. Equivalently, $[L_\ell' \, \rho_j] ({\sf
x}') = 0$ for every ${\sf x}' \in \Lambda'$. Letting $({\sf x}_1,
{\sf x}') \in \Lambda$, whence ${\sf x}' \in \Lambda '$, we compute:
\[
\aligned
\big[
\widetilde{ L}\,\rho_j
\big]({\sf x}_1, {\sf x}')
&
=
\sum_{2\leqslant i\leqslant n}\,
a_i({\sf x}_1,{\sf x}')\,
\frac{\partial \rho_j}{\partial {\sf x}_i}({\sf x}')
\\
&
=
\sum_{2\leqslant i\leqslant n}\,
\sum_{\ell=0}^\infty\,
\frac{{\sf x}_1^\ell}{\ell!}\,
\frac{\partial^\ell a_i}{\partial {\sf x}_1^\ell}
(0,{\sf x}')\,
\frac{\partial\rho_j}{\partial {\sf x}_i}({\sf x}')
\ \ \ \ \ 
\text{\rm [Taylor development]}
\\
&
=
\sum_{\ell=0}^\infty\,
\frac{{\sf x}_1^\ell}{\ell!}\,
\big[L_\ell'\,\rho_j\big]({\sf x}')
=0,
\endaligned
\]
so $\widetilde{ L}$ is tangent to $\Lambda$. Finally, the property
$T_{{\sf x}_1, {\sf x}'} \Lambda = \LL ( {\sf x}_1, {\sf x}')$ follows
immediately from $T_{{\sf x}'} \Lambda' = \LL' ({\sf x}')$ and the
proof is complete (the Taylor development argument above was crucially
used, and this enlightens why the theorem does not hold in the
$\mathcal{ C}^\infty$ category).

\subsection*{ 1.8.~Free Lie algebras and generic sets of
$\K$-analytic vector fields} For a generic set of $r\geqslant 2$ vector
fields $\LL^0 = \{ L_a \}_{
1 \leqslant a \leqslant r}$, 
or after slightly perturbing any given set, one
expects that $\LL^{ \rm lie} (0) = T_0 \K^n$. Then the Nagano
leaf $\Lambda$ passing through $0$ is just an open neighborhood of $0$
in $\K^n$. Also, one expects that the dimensions
of the intermediate spaces $\LL^k (0)$ be maximal.

To realize this intuition, one has to count the maximal number of
iterated Lie brackets that are linearly independent in $\LL^k$, for $k
= 1, 2, 3, \dots$, modulo antisymmetry and Jacobi identity.

Let $r\geqslant 2$ and let $h_1, h_2, \dots, h_r$ be $r$ linearly
independent elements of a vector space over $\K$. The {\sl free Lie
algebra} ${\sf F}(r)$ of rank $r$ is the smallest (non-commutative,
non-associative) $\K$-algebra (\cite{ re1993}) having $h_1, h_2,
\dots, h_r$ as elements, with multiplication $(h, h') \mapsto h\, h'$
satisfying antisymmetry $0 = h \, h' + h' \, h$ and Jacobi identity $0
= h(h'\,h'') + h'' (h\, h') + h' (h''\, h)$. It is unique up to
isomorphism. The case $r= 1$ yields only ${\sf F} (1) = \K$. The
multiplication in ${\sf F} (r)$ plays the role of the Lie bracket in
$\LL^{ \rm lie}$. Importantly, {\it no linear relation exists between
iterated multiplications}, {\it i.e.} between iterated Lie brackets,
{\it except those generated by antisymmetry and Jacobi
identity}. Thus, ${\sf F} (r)$ is infinite-dimensional. Every
finite-dimensional Lie $\K$-algebra having $r$ generators embeds as a
subalgebra of ${\sf F} (r)$, {\it see} \cite{ re1993}.

Since the bracket multiplication is not associative, one must
carefully write some parentheses, for instance in $(h_1\, h_2)h_3$, or
in $h_1( h_2 (h_1\, h_2))$, or in $(h_1\, h_2)( h_3( h_5\,
h_1))$. Writing all such words only with the alphabet $\{ h_1, h_2,
\dots, h_r \}$, we define the {\sl length} of a word $\mathbf{ h}$ to be
the number of elements $h_{ i_\alpha }$ in it. For $\ell \in \N$ with
$\ell \geqslant 1$, let ${\sf W}_r^\ell$ be the set of words of length
equal to $\ell$ and let ${\sf W}_r = \bigcup_{ \ell \geqslant 1} {\sf
W}_r^\ell$ be the set of all words.

Define ${\sf F}_1(r)$ to be the vector space generated by $h_1,
h_2, \dots, h_r$ and for $\ell \geqslant 2$, define ${\sf F}_\ell (r)$ to
be the vector space generated by words of length $\leqslant \ell$.
This corresponds to $\LL^\ell$, except that in $\LL^\ell$, there might
exist special linear relations that are absent in the abstract case.
Thus, ${\sf F} (r)$ is a graded Lie algebra. The Jacobi identity
insures (by induction) that ${\sf F}_{ \ell_1} (r) {\sf F}_{ \ell_2}
(r) \subset {\sf F}_{ \ell_1 + \ell_2} (r)$, a property similar to
$\left[ \LL^{ k_1}, \LL^{ k_2} \right] \subset \LL^{ k_1 + k_2}$. It
follows that ${\sf F}_\ell (r)$ is generated by words of the form
\[
h_{i_1}(h_{i_2}(\dots(h_{i_{\ell'-1}}h_{i_{\ell'}})\dots)),
\]
where $\ell' \leqslant \ell$ and where $1 \leqslant i_1, i_2, \dots,
i_{ \ell' -1}, i_{ \ell'} \leqslant r$. For instance, $(h_1\,
h_2)(h_3( h_5\, h_1))$ may be written as a linear combination of such
simple words whose length is $\leqslant 5$. Let us denote by
\[
{\sf SW}_r 
= 
\bigcup_{\ell\geqslant 1}\
{\sf SW}_r^\ell 
\]
the set of these simple words, where ${\sf SW}_r^\ell$ denotes the set
of simple words of length $\ell$. Although it generates ${\sf F} (r)$
as a vector space over $\K$, we point out that it is not a basis of
${\sf F} (r)$: for instance, we have $h_1 (h_2 (h_1 h_2)) = h_2 (h_1
(h_1 h_2))$, because of an obvious Jacobi identity in 
which $(h_1h_2) (h_1h_2) = 0$ disappears. In fact, one
verifies that this is the only
Jacobi relation between simple words of length $4$, that
simple words of length $5$ have no
Jacobi relation, hence a basis of ${\sf F}_5 (2)$ is
\[
\aligned
&
h_1,
\ \ 
h_2,
\ \ \ \ \ 
h_1h_2,
\\
&
h_1(h_1h_2),
\ \ 
h_2(h_1h_2),
\\
&
h_1(h_1(h_1h_2)),
\ \ 
h_1(h_2(h_1h_2)),
\ \ 
h_2(h_2(h_2h_1)),
\\ 
&
h_1(h_1(h_1(h_1h_2))),
\ \
h_1(h_1(h_2(h_1h_2))),
\ \
h_1(h_2(h_2(h_2h_1))),
\\
& 
\ \ \ \ \ \ \ \ \ \ \ \ \ \ \ \ \ \ \ \ \ \ \
h_2(h_1(h_1(h_1h_2))),
\ \ 
h_2(h_2(h_1(h_2h_1))),
\ \
h_2(h_2(h_2(h_2h_1))).
\endaligned
\]

In general, what are the dimensions of the ${\sf F}_\ell (r)$\,? How
to find bases for them, when considered as vector spaces\,?

\def\thedefinition{1.9}\begin{definition}{\rm
A {\sl Hall-Witt basis} of ${\sf F}(r)$ is a linearly ordered (infinite)
subset ${\sf HW}_r = \bigcup_{ \ell \geqslant 1} {\sf HW}_r^\ell$ of the
set of simple words ${\sf SW}_r$ such that:

\begin{itemize}

\smallskip\item[$\bullet$]
if two simple words ${\bf h}$ and ${\bf h}'$ satisfy ${\sf
length}({\bf h}) < {\sf length}({\bf h'})$, then ${\bf h} < {\bf h'}$;

\smallskip\item[$\bullet$]
${\sf HW}_r^1 = \{ h_1, h_2, \dots, h_r \}$;

\smallskip\item[$\bullet$]
${\sf HW}_r^2 = \{ h_{ i_1} h_{ i_2} : \ 1\leqslant i_1 <
i_2 \leqslant r\}$;

\smallskip\item[$\bullet$]
${\sf HW}_r \backslash ( {\sf HW}_r^1 \cup {\sf HW}_r^2) = \{ {\bf
h}({\bf h'} {\bf h''}): \ {\bf h}, {\bf h'}, {\bf h''} \in {\sf HW}_r,
\ {\bf h'} < {\bf h''} \ {\rm and} \ {\bf h'} \leqslant {\bf h} < {\bf
h'}{\bf h''} \}$.

\end{itemize}\smallskip

}\end{definition}

A Hall-Witt basis essentially consists of the choice, for every $\ell
\geqslant 1$, of some (among many possible) finite subset ${\sf HW}_r^\ell$
of ${\sf SW}_r^\ell$ that generates the finite-dimensional quotient
vector space ${\sf F}_\ell ( r) / {\sf F}_{ \ell - 1} (r)$. To fix
ideas, an arbitrary linear ordering is added among the elements of the
chosen basis ${\sf HW}_r^\ell$ of the vector space ${\sf F}_\ell ( r)
/ {\sf F}_{ \ell - 1} (r)$. The last condition of the definition takes
account of the Jacobi identity.

\def\thetheorem{1.10}\begin{theorem}
{\rm (\cite{ bo1972, re1993})}
Hall-Witt bases exist and are bases of the free Lie algebra ${\sf F}
(r)$ of rank $r$, when considered as a vector space. The
dimensions ${\sf n }_\ell (r) - {\sf n}_{
\ell - 1} (r)$ of ${\sf F}_\ell (r) / {\sf
F}_{ \ell - 1} (r)$, or equivalently the cardinals of ${\sf HW 
}_r^\ell$, satisfy the induction relation
\[
{\sf n}_\ell(r)
-
{\sf n}_{\ell-1}(r)
=
\frac{1}{\ell}\,
\sum_{d\ {\rm divides}\ \ell}
\mu(d)\,r^{\ell/d},
\] 
where $\mu$ is the M\"obius function.

\end{theorem}

Remind that
\[
\mu(d)
=
\left\{
\aligned
&
1,\ {\rm if}\ d=1;
\\
&
0,\ {\rm if}\ d\
\text{\rm contains square integer factors};
\\
&
(-1)^\nu,\ {\rm if}\ d=p_1\cdots p_\nu\
\text{\rm is the product of}\ \nu \
\text{\rm distinct prime numbers}.
\endaligned
\right.
\]

\smallskip

Now, we come back to the system $\LL^0 = \{ L_a \}_{ 1\leqslant a
\leqslant r}$ of local $\K$-analytic vector fields of \S1.1, where
$L_a = \sum_{ i= 1}^n \, \varphi_{ a, i} ({\sf x}) \, \frac{ \partial
}{ \partial {\sf x}_i}$. If the vector space $\LL (0)$ has dimension
$< r$, a slight perturbation of the coefficients $\varphi_{ a, i}
({\sf x})$ of the $L_a$ yields a system ${\LL' }^0$ with $\LL' (0)$ of
dimension $=r$. By an elementary computation with Lie brackets, one
sees that a further slight perturbation yields a system ${\LL''}^0$
with $\LL{ ''} (0)$ of dimension $r + \frac{ r ( r- 1)}{ 2} = {\sf
n}_2 (r)$.

To pursue, any simple iterated Lie bracket $[ L_{
a_1}, [ L_{ a_2}, \dots [ L_{ a_{ \ell-1}}, L_{ a_\ell}]\dots ]]$ of
length $\ell$ is a vector field $\sum_{ i=1}^n\, A_{ a_1, a_2, \dots,
a_{ \ell - 1}, a_\ell}^i \, \frac{ \partial}{ \partial {\sf x}_i}$
having coefficients $A_{ a_1, a_2, \dots, a_{ \ell - 1}, a_\ell }^i$
that are universal polynomials in the jets
\[
J_{\sf x}^{\ell-1}\varphi({\sf x})
:=
\big(
\partial_{\sf x}^\alpha
\varphi_{a,i}({\sf x})
\big)_{1\leqslant a\leqslant r,\ 
1\leqslant i\leqslant n}^{
\alpha\in\N^n,\
\vert\alpha\vert\leqslant\ell-1}
\in\K^{ N_{ rn, n, \ell-1}}
\]
of order $(\ell - 1)$ of the coefficients of $L_1, L_2, \dots,
L_r$. Here, $N_{ rn, n, \ell - 1} = rn \, \frac{ (n+ \ell - 1)!}{ n! \
(\ell - 1)!}$ denotes the number of such independent partial
derivatives. A careful inspection of the polynomials $A_{ a_1, a_2,
\dots, a_{ \ell - 1}, a_\ell }^i$ enables to get the following
genericity statement, whose proof will appear elsewhere. It says in a
precise way that $\LL^{ \rm lie} (0) = T_0 \K^n$ with the maximal
freedom, for generic sets of vector fields.

\def\thetheorem{1.11}\begin{theorem}
{\rm (\cite{ gv1987, ge1988}, [$*$])} 
If $\ell_0$ denotes the smallest length
$\ell$ such that ${\sf n}_\ell (r) \geqslant n$, there exists a proper
$\K$-algebraic subset $\Sigma$ of the jet space $J_0^{ \ell_0 - 1}
\varphi = \K^{ N_{ rn, n, \ell_0-1}}$ such that for every collection
$\LL^0 = \{ L_a \}_{ 1\leqslant a\leqslant r}$ of $r$ vector fields
$L_a = \sum_{ i=1}^n\, \varphi_{ a, i} ({\sf x}) \, \frac{ \partial
}{\partial {\sf x}_i}$ such that $J_0^{ \ell_0 - 1} \varphi (0)$ does
not belong to $\Sigma$, the following two properties hold{\rm :}

\begin{itemize}

\smallskip\item[$\bullet$]
$\dim \LL^\ell (0) 
= 
{\sf n}_\ell (r)$, 
for every $\ell \leqslant \ell_0 - 1$, 

\smallskip\item[$\bullet$]
$\dim \LL^{\ell_0} (0)
= n$, hence $\LL^{ \rm lie} (0) = 
T_0 \K^n$. 

\end{itemize}

\smallskip

\end{theorem}

The number of divisors of $\ell$ being an ${\rm O} (\frac{ {\rm log}
\, \ell}{ {\rm log}\, 2})$, one verifies that ${\sf n }_\ell(r) - {\sf
n}_{ \ell-1}(r) = \frac{ 1}{ \ell }\, r^\ell + {\rm O} (r^{ \ell /2}\,
\frac{ {\rm log}\, \ell }{ {\rm log}\, 2})$. It follows that, for $r$
fixed, the integer $\ell_0$ of the theorem is equivalent to $\frac{
{\rm log}\, n}{ {\rm log}\, r}$ as $n \to \infty$.

\subsection*{1.12.~Local orbits of $\K$-analytic and of
(Nash) $\K$-algebraic systems} We now describe a second, more
concrete, simple and useful approach to the local Nagano Theorem~1.5.
It is inspired by Sussmann's Theorem~1.21(III) and does {\it not}\,
involve the consideration of any Lie bracket. Theorem~1.13 below will
be applied in \S2.11.

As above, consider a {\it finite}\, set 
\[
\LL^0 
:= 
\{ L_a \}_{ 1
\leqslant a \leqslant r}, 
\ \ \
r \in \N,
\ \ \
r \geqslant 1, 
\]
of nonzero vector fields defined in the cube $\square_{ \rho_1 }^n$
and having $\K$-analytic coefficients. We shall neither consider its
$\A_{ \rho_1 }$-linear hull $\LL$, nor $\LL^{ \rm lie}$. We will
reconstruct the Nagano leaf passing through the origin only by means
of the flows of $L_1, L_2, \dots, L_r$.

Referring the reader to \S1.3(III) for background, we denote the flow
map of a vector field $L \in \LL^0$ shortly by $({\sf t}, {\sf x})
\mapsto L_{\sf t} ( {\sf x}) = \exp ( {\sf t} L) ({\sf x})$. It is
$\K$-analytic. What happens in the algebraic category\,?

So, assume that the coefficients of all vector fields $L\in \LL^0$ are
$\K$-algebraic. Unfortunately, algebraicity fails to be preserved
under integration, so the flows are only $\K$-analytic, in general. To
get algebraicity of Nagano leaves, there is nothing else than
supposing that the flows are algebraic, which we will do
(second phrase of {\bf (5)} below).

Choose now $\rho_2$ with $0 < \rho_2 < \rho_1$. Let $k \in \N$ with $k
\geqslant 1$, let $L = ( L^1, \ldots, L^k) \in (\LL^0)^k$, let ${\sf
t} = ( {\sf t}_1, \ldots, {\sf t }_k) \in \K^k$ with $\vert {\sf t}
\vert < \rho_2$, {\it i.e.} ${\sf t} \in \square_{ \rho_2}^k$, and let
${\sf x} \in \square_{ \rho_2 }^n$. We shall adopt the contracted
notation 
\[
L_{\sf t }( {\sf x}) : = L_{ {\sf t }_k }^k ( \cdots( L_{
{\sf t}_1 }^1 ({\sf x} )) \cdots ) 
\]
for the composition of flow maps, whenever it is defined. In fact,
since $L_0 (0) = \exp (0 L ) (0) = 0$, it is clear that if we bound
the length $k\leqslant 2n$, then there exists $\rho_2 > 0$ sufficiently
small such that all maps $({\sf t}, {\sf x}) \mapsto L_{\sf t} ({\sf
x})$ are well-defined, with $L_{ \sf t} ({\sf x}) \in \square_{
\rho_1}^n$, at least for all ${\sf t}\in \square_{ \rho_2}^k$ and all
${\sf x } \in \square_{ \rho_2 }^n$. The reason why we may restrict to
consider only compositions of length $k \leqslant 2 n$ will appear
{\it a posteriori}\, in the proof of the theorem below. We shall be
concerned with rank properties of $({\sf t}, {\sf x})
\mapsto L_{\sf t} ( {\sf x})$.

Let $n \geqslant 1$, $m \geqslant 1$, $\rho_1 > 0$, $\sigma_1 > 0$ and
let $f : \square_{ \rho_1 }^n \to \square_{ \sigma_1 }^m$, ${\sf x}
\mapsto f ({\sf x})$, be a $\K$-algebraic or $\K$-analytic map between
two open cubes. Denote its Jacobian matrix by ${\rm Jac} \, (f) =
\big( \frac{ \partial f_j }{ \partial {\sf x}_i} ({\sf x}) \big)_{
1\leqslant i \leqslant n }^{ 1 \leqslant j \leqslant m}$. At a point
${\sf x} \in \square_{ \rho_1 }^n$, the map $f$ has rank $r$ if and
only if ${\rm Jac} \, f$ has rank $r$ at ${\sf x}$. Equivalently, by
linear algebra, there is a $r\times r$ minor that does not vanish at
${\sf x}$ but all $s \times s$ minors with $r + 1 \leqslant s \leqslant
n$ do vanish at ${\sf x}$.

For every $s \in \N$ with $1 \leqslant s \leqslant \min (n, m)$,
compute all the possible $s \times s$ minors $\Delta_1^{s \times s},
\dots, \Delta_{ N(s)}^{ s \times s}$ of ${\rm Jac}\, ( f)$. They are
universal (homogeneous of degree $s$) polynomials in the partial
derivatives of $f$, hence are all $\K$-algebraic or $\K$-analytic
functions. Let $e$ with $0 \leqslant e \leqslant \min (n, m)$ be the
maximal integer $s$ with the property that there exists a minor
$\Delta_{ \mu}^{ s\times s} ({\sf x})$, 
$1 \leqslant \mu \leqslant N(s)$, not
vanishing identically. Then the set
\[
\mathcal{R}_f
:=
\big\{ {\sf x} \in \square_{ \rho_1}^n : \, 
\Delta_\mu^{ s\times s} ({\sf x}) = 0, \, 
\mu 
= 
1, \dots, N(s) 
\big\}
\]
is a {\it proper}\, $\K$-algebraic or analytic subset of $\square_{
\rho_1}^n$. The principle of analytic continuation insures that
$\square_{ \rho_1 }^n \big \backslash \mathcal{ R}_f$ is open and
dense.

The integer $e$ is called the {\sl generic rank} of $f$. For every
open, connected and nonempty subset $\Omega \subset \square_{
\rho_1}^n$ the restriction $f \vert_\Omega$ has the same generic rank
$e$.

\def\thetheorem{1.13}\begin{theorem} 
{\rm (\cite{ me1999, me2001a, me2004a})} There exists an integer $e$
with $1\leqslant e \leqslant n$ and an $e$-tuple of vector fields
$L^*=( L^{*1 }, \ldots, L^{ *e}) \in (\LL^0)^{e}$ such that the
following six properties hold true.

\smallskip

\begin{itemize}
\item[{\bf (1)}] 
For every $k=1, \dots, e$, the map $({\sf t }_1, \dots, {\sf t}_k)
\mapsto L_{ {\sf t}_k }^{ *k}( \cdots (L_{{\sf t }_1 }^{ *1} (0))
\cdots)$ is of generic rank equal to $k$.

\smallskip

\item[{\bf (2)}] 
For every arbitrary element $L'\in \LL^0$, the map $({\sf
t}_1,\dots,{\sf t}_e,{\sf t}')\mapsto L_{ {\sf t}'}'(L_{{\sf
t}_e}^{*e}(\cdots (L_{{\sf t}_1}^{*1}(0)) \cdots))$ is of generic rank
$e$, hence $e$ is the maximal possible generic rank.

\smallskip

\item[{\bf (3)}] 
There exists an element ${\sf t }^*\in \square_{ \rho_2}^e$
arbitrarily close to the origin which is of the special form $({\sf t
}_1^*, \ldots, {\sf t }_{e-1 }^*,0)$, namely with ${\sf t}_{e}^*=0$,
and there exists an open connected neighborhood $\omega^*$ of ${\sf t
}_*$ in $\square_{ \rho_2}^e$ such that the map ${\sf t} \mapsto L_{
{\sf t}_e}^{ *e}( \cdots (L_{ {\sf t }_1}^{ *1 }(0 )) \cdots)$ is of
constant rank $e$ in $\omega^*$.

\smallskip

\item[{\bf (4)}]
Setting $L^* :=( L^{ *1}, \ldots, L^{*e})$, $K^* := (L^{ *e-1 },
\ldots,L^{*1})$ and $s^* := (-{\sf t }_{ e-1}^*, \ldots, - {\sf t
}_1^*)$, we have $K^*_{ s^*} \circ L_{ {\sf t }^*}^*
(0)=0$ and the map $\psi: \omega^* \to \square_{ \rho_1 }^n$ defined by
$\psi: {\sf t} \mapsto K_{ s^* }^*\circ L_{ \sf t}^*(0)$ is also of
constant rank equal to $e$ in $\omega^*$.

\smallskip

\item[{\bf (5)}] 
The image $\Lambda := \psi (\omega^*)$ is a piece of $\K$-analytic
submanifold passing through the origin enjoying the most important
property that every vector field $L' \in \LL^0$ is tangent to
$\Lambda$. If the flows of all elements of $\LL^0$ are algebraic,
$\Lambda$ is $\K$-algebraic.

\smallskip

\item[{\bf (6)}] 
Every local $\K$-algebraic or $\K$-analytic submanifold $\Lambda'$
passing trough the origin to which all vector fields $L' \in \LL^0$
are tangent must contain $\Lambda$ in a neighborhood of $0$.

\end{itemize} 

\smallskip

In conclusion, the dimension $e$ of $\Lambda$ is characterized by the
generic rank properties {\bf (1)} and {\bf (2)}.
\end{theorem}

Previously, $\Lambda$ was called {\sl Nagano leaf}. Since the above
statement is superseded by Sussmann's Theorem~1.21(III), we prefer to
call it the {\sl local $\LL$-orbit of} $0$, introducing in advance the
terminology of Part~III and denoting it by $\mathcal{ O }_{\LL^0}^{
loc} (\square_{\rho_1}^n, 0)$. The integer $e$ of the theorem is
$\leqslant n$, just because the target of the maps $({\sf t }_1,
\dots, {\sf t}_k) \mapsto L_{ {\sf t}_k }^{ *k}( \cdots (L_{{\sf t }_1
}^{ *1} (0)) \cdots)$ is $\K^n$. It follows that in {\bf (4)} and {\bf
(5)} we need $2e - 1 \leqslant 2n -1$ compositions of flows to cover
$\Lambda$.

\smallskip

We quickly mention
an application about separate algebraicity.
In~\cite{ bm1949}, it is shown that a local $\K$-analytic function $g
: \square_{ \rho_1 }^n \to \K$ is $\K$-algebraic if and only if its
restriction to every 
affine coordinate segment is $\K$-algebraic. Call the
system $\LL^0$ {\sl minimal} at the origin if $\mathcal{ O}_{ \LL^0
}^{loc} ( \square_{ \rho_1 }^n, 0)$ contains a neighborhood of the
origin. Equivalently, the integer $e$ of Theorem~1.13 equals $n$.

\def\thetheorem{1.14}\begin{theorem}
{\rm (\cite{ me2001a})} If $\LL^0 = 
\{ L_a \}_{ 1\leqslant a \leqslant r}$ 
is minimal at $0$, a local
$\K$-analytic function $g : \square_{ \rho_1 }^n \to \K$ is
$\K$-algebraic if and only it its restriction 
to every integral curve of every
$L_a \in \LL^0$ is $\K$-algebraic.
\end{theorem}

\proof[Proof of Theorem~1.13.]
(May be skipped in a first reading.)
If all vector fields of $\LL^0$ vanish at the origin, $\Lambda = \{
0\}$. We now exclude this possibility. Choose a vector field $L^{ *1}
\in \LL^0$ which does not vanish at $0$. The map ${\sf t}_1 \mapsto L_{
{\sf t}_1 }^{ *1}(0)$ is of (generic) rank one at
every ${\sf t}_1 \in \square_{ \rho_2}^1$. If there exists $L'\in
\LL^0$ such that the map $({\sf t}_1, {\sf t}')\mapsto L_{ {\sf
t}'}'(L_{{\sf t}_1 }^{ *1 }(0))$ is of generic rank two, we choose one
such $L'$ and we denote it by $L^{ *2}$. Continuing in this way, we
get vector fields $L^{*1}, \dots, L^{ *e}$ satisfying properties {\bf
(1)} and {\bf (2)}, with $e\leqslant n$. 

Since the generic rank of the map $({\sf t }_1, \dots, {\sf t}_e)
\mapsto L_{ {\sf t}_e}^{ *e}( \cdots (L_{ {\sf t}_1 }^{*1} (0))
\cdots)$ equals $e$, and since this map is $\K$-analytic, there exists
a ${\sf t}^*\in \square_{ \rho_2 }^e$ arbitrarily close to the origin
at which its rank equals $e$. We claim that we can moreover choose
${\sf t}^*$ to be of the special form $( {\sf t}_1^*, \dots, {\sf
t}_{e-1 }^*, 0)$, {\it i.e.} with $t_e^*=0$. It suffices to apply the
following lemma to $\varphi ({\sf t}) := L_{ {\sf t}_{e-1}}^{ *e-1}(
\cdots (L_{ {\sf t}_1 }^{*1} (0)) \cdots)$ and to $L' := 
L^{ *e}$.

\def\thelemma{1.15}\begin{lemma}
Let $n \in\N$, $n\geqslant 1$, let $e \in \N$, $1\leqslant e\leqslant n$, let ${\sf
t} \in \square_{\rho_2}^{e-1}$ and let
\[
\square_{\rho_2}^{e-1}\ni
{\sf t} \mapsto \varphi ({\sf t}) = (\varphi_1( {\sf
t}), \dots, \varphi_n ({\sf t}))
\in\square_{\rho_1}^n
\]
be a $\K$-analytic map whose generic rank equals $(e-1)$. Let $L'$
be a $\K$-analytic vector field
and assume that the map $\psi: ( {\sf t}, {\sf t}')\mapsto L_{
{\sf t}'}'(\varphi( {\sf t}))$ has generic rank $e$. Then there exists
a point $( {\sf t}^*,0)$ arbitrarily close to the origin at which the
rank of $\psi$ is equal to $e$.
\end{lemma}

\proof
Choose ${\sf t}^\sharp \in \square_{ \rho_2}^{ e-1}$ arbitrarily close
to zero at which $\varphi$ has maximal rank, equal to $(e-1)$. Since
the rank is lower semi-continuous, there exists a connected
neighborhood $\omega^\sharp$ of ${\sf t}^\sharp$ in $\square_{
\rho_2}^{ e-1}$ such that $\varphi$ has rank $(e-1)$ at every point of
$\omega^\sharp$. By the constant rank theorem, $\Pi := \varphi(
\omega^\sharp)$ is then a local $\K$-analytic submanifold of
$\square_{\rho_1}^n$ passing through the point $\varphi ({\sf t
}^\sharp )$. To complete the proof, we claim that there exists ${\sf
t}^* \in \omega^\sharp$ arbitrarily close to ${\sf t}^\sharp$ such
that the map $( {\sf t}, {\sf t}') \mapsto L_{{\sf t}'}' (\varphi (
{\sf t}))$ has rank $e$ at $({\sf t}^*, 0)$.

Let us reason by contradiction, supposing that at all points of the
form $({\sf t}^*, 0)$, for ${\sf t}^* \in \omega^\sharp$, the map
$\psi: ( {\sf t}, {\sf t}') \mapsto L_{{\sf t}'}' (\varphi( {\sf t}
))$ has rank equal to $(e-1)$. Pick
an arbitrary ${\sf t}^* \in \omega^\sharp$. 
Reminding that when ${\sf t}' = 0$, we have $L_{\sf t'}' =
L_0 ' = {\rm Id}$, we observe that $\psi( {\sf t}, 0) \equiv \varphi
({\sf t})$. Consequently, the partial derivatives of $\psi$ with
respect to the variables ${\sf t}_i$, $i=1, \dots, e-1$ at an
arbitrary point $({\sf t}^*, 0)$, with ${\sf t}^* \in \omega^\sharp$,
coincide with the $(e-1)$ linearly independent vectors $\frac{
\partial \varphi }{ \partial {\sf t}_i} ({\sf t}^*) \in \K^n$,
$i=1, \dots, e-1$. In fact, the tangent space to $\Pi$ at the point
$\psi ({\sf t}^*, 0) = \varphi ( {\sf t}^*)$ is generated by these
$(e-1)$ vectors.

Reminding the fundamental property $\left. \frac{\partial }{ \partial
{\sf t}'} L_{{\sf t}'} ' ({\sf x}) \right \vert_{{\sf t'}=0} = L' ({\sf
x})$, we deduce [from our assumption that the map $( {\sf t}, {\sf
t}') \mapsto L_{{\sf t}'}' (\varphi( {\sf t} ))$ has rank equal to
$(e-1)$] that the vector
\[
\left.
\frac{\partial }{\partial {\sf t}'}
L_{{\sf t}'}'(\varphi ({\sf t}))
\right\vert_{{\sf t}'=0}
=
L'(\varphi({\sf t}))
\]
must be linearly dependent with the $(e-1)$ vectors $\frac{ \partial
\varphi }{ \partial {\sf t}_i} ({\sf t})$, $i=1, \dots, e-1$, for
every ${\sf t} \in \omega^\sharp$. Equivalently, the vector field
$L'$ is tangent to the submanifold $\Pi$. It follows that the local
flow of $L'$ necessarily stabilizes $\Pi$: if ${\sf x} = \varphi ({\sf
t}) \in \Pi$, ${\sf t} \in \omega^\sharp$, then $L_{ {\sf t}'}' ({\sf
x}) \in \Pi$, for all ${\sf t}' \in \square_{ \rho( {\sf t})}^1$, where
$\rho ({\sf t}) >0$ is sufficiently small. Set $\Omega^\sharp := \{
({\sf t}, {\sf t}') : {\sf t } \in \omega^\sharp, {\sf t}' \in
\square_{\rho ({ \sf t })^1 }\}$. It is a nonempty connected open subset
of $\square_{ \rho_2}^e$. We have thus deduced that $\psi
(\Omega^\sharp)$ is contained in the $(e - 1)$-dimensional submanifold
$\Pi$. This constraint entails that $\psi$ is of rank $\leqslant e-1$
at every point of $\Omega^\sharp$. However,
$\psi\vert_{\Omega^\sharp}$ being $\K$-analytic and of generic rank
equal to $e$, by 
assumption, 
it should be of rank $e$ at every point of an open dense
subset of $\Omega^\sharp$. This is the desired contradiction which
proves the lemma.
\endproof

\smallskip

Hence, there exists ${\sf t}^* =( {\sf t }_1^*, \dots, {\sf t }_{
e-1}^*,0) \in \square_{ \rho_2 }^e$ arbitrarily close to the origin at
which the rank of ${\sf t} \mapsto L_{ {\sf t}_e }^{ *e} (\cdots (L_{
{\sf t }_1}^{ *1} (0)) \cdots)$ is maximal (hence locally constant)
equal to $e$, so we get the constant rank property {\bf (3)}, for a
sufficiently small neighborhood $\omega^*$ of ${\sf t}^*$.

In {\bf (4)}, the property $K_{s^*}^*\circ L_{t^*}^*(0) = 0$ is obvious,
using ${\sf x} \equiv L_0 ({\sf x}) \equiv L_{-{\sf t}} \circ L_{
\sf t} ( {\sf x})$:
\[
L_{-{\sf t}_1^*}^{*1}\circ\cdots 
\circ L_{-{\sf t}_{e-1}^*}^*\circ
L_{0}^{*e}\circ L_{{\sf t}_{e-1}^*
}^*\circ \cdots \circ L_{{\sf t}_1^*}^*({\sf
x})\equiv{\sf x}.
\]
Since the map ${\sf x} \mapsto K_{ s^* }^*( {\sf x})$ is a local
diffeomorphism, it is clear that the map $\psi : {\sf t} \mapsto K_{
{\sf s}^*}^* \circ L_{\sf t}^* (0)$ is also of constant rank $e$ in
$\omega^*$. Thus, we obtain {\bf (4)}, and moreover, by the constant
rank theorem, the image $\Lambda := \psi (\omega^*)$ constitutes a
local $\K$-analytic submanifold of $\K^n$ passing through the
origin. If the flows of elements of $\LL^0$ are all $\K$-algebraic,
clearly $\psi$ and $\Lambda$ are also $\K$-algebraic.

It remains to check that every vector field $L'\in \LL^0$ is tangent to
$\Lambda$. As a preliminary, denote by $L_{{\sf t'}}' (\varphi ({\sf
t}))$, ${\sf t} \in \square_{\rho_2}^e$, ${\sf t}' \in
\square_{\rho_2}^1$, the map appearing in 
{\bf (2)}, where $L' \in \LL^0$ is
arbitrary. Reasoning as in the lemma above, we see that $L'$ is
necessarily tangent to some local submanifold $\Pi$ obtained as the
local image of an open connected set where $\varphi$ has maximal
locally constant rank. It follows that the flows and the multiple
flows of elements of $\LL^0$ stabilize this submanifold. We deduce a
generalization of {\bf (2)}: for $k\leqslant 2n$, for $L'\in 
(\LL^0)^k$, for
${\sf t}' \in \square_{ \rho_2}^k$, the map $({\sf t}, {\sf t}')
\longmapsto L_{ {\sf t'} }'(L_{{\sf t}_e}^{ *e}( \cdots (L_{{\sf t}_1
}^{ *1} (0)) \cdots))$ is of generic rank $e$.

In particular, for every $L'\in \LL^0$, the map $({\sf t}', {\sf s},
{\sf t}) \longmapsto L_{{\sf t}'}' \circ K_{\sf s}^* \circ L_{\sf t}^*
(0)$ is of generic rank $e$. In fact, the restriction $\psi: {\sf t}
\mapsto K_{{\sf s}^*}^* \circ L_{\sf t}^* (0)$ of this map to the open
set $\{ (0, {\sf s}^*, {\sf t}): {\sf t} \in \omega^*\}$ is already of
rank $e$ at every point and its image is the local submanifold
$\Lambda$, by the above construction. So the map $({\sf t}', {\sf t} )
\longmapsto L_{{\sf t}'} ' \circ K_{{\sf s}^*}^* \circ L_{\sf t}^* (0)$
must be of rank $e$ at every point. In particular, the vector
\[
\left. \frac{ \partial }{\partial {\sf t}'} \,
L_{{\sf t}'} ' \circ
K_{{\sf s}^*}^* \circ L_{\sf t}^* (0) 
\right \vert_{{\sf t}'=0}
=
L'
\big(
K_{{\sf s}^*}^* \circ L_{\sf t}^* (0) 
\big)
\] 
must
necessarily be tangent to $\Lambda$ at the point $K_{{\sf s}^*}^*
\circ L_{\sf t}^* (0) \in \Lambda$. Thus, {\bf (5)} is proved.

Take $\Lambda'$ as in {\bf (6)}. The local 
flows of all vector $L' \in \LL^0$ stabilize
$\Lambda'$. Shrinking $\rho_2$ if necessary, 
all the maps $({\sf t}, {\sf x}) 
\longmapsto L_{\sf t} ({\sf x})$ have range in $\Lambda'$. 
So $\Lambda \subset \Lambda'$, proving
{\bf (6)}.
\endproof 

\section*{\S2.~Analytic CR manifolds, Segre chains and minimality}

\subsection*{2.1.~Local Cauchy-Riemann submanifolds of $\C^n$}
Let $(z_1, \dots, z_n) = (x_1 + i y_1, \dots, x_n + i y_n)$ denote the
canonical coordinates on $\C^n$. As before, we use the maximum norms
$\vert x \vert = \max_{1 \leqslant k \leqslant n} \, \vert x_k \vert$,
$\vert y \vert = \max_{ 1\leqslant k\leqslant n} \, \vert y_k \vert$
and $\vert z\vert = \max_{1 \leqslant k \leqslant n} \, \vert
z_k\vert$, where $\vert z_k \vert= (x_k^2 + y_k^2 )^{1 /2}$. If
$\rho>0$, we denote by $\Delta_\rho^n = \{z \in\C^n: \, \vert z\vert <
\rho\}$ the open polydisc of radius $\rho$ centered at the origin, not
to be confused with $\square_\rho^{ 2n} = \{ x+iy \in \C^n : \vert x
\vert, \vert y \vert < \rho \}$.

Let $J$ denote the complex structure of $T \C^n$, acting on real
vectors as if it were multiplication by $\sqrt{ -1}$. Precisely, if
$p$ is any point, $T_p \C^n$ is spanned by the $2n$ vectors
$\frac{ \partial}{ \partial x_k} \big\vert_p$, $\frac{
\partial}{ \partial y_k} \big\vert_p$, $k=1, \dots, n$, and $J$ acts
as follows: $J \frac{ \partial}{ \partial x_k} \big\vert_p =
\frac{ \partial }{\partial y_k} \big\vert_p $; $J
\frac{ \partial}{ \partial y_k } \big\vert_p = - \frac{
\partial }{\partial x_k} \big\vert_p$.

Choose the origin as a center point and consider a real
$d$-codimensional local submanifold $M$ of $\C^n \simeq \R^{ 2n}$
passing through the origin, defined by $d$ Cartesian equations $r_1(
x, y) = \cdots= r_d (x,y) = 0$, where the differentials $dr_1, \dots,
dr_d$ are linearly independent at the origin. The functions $r_j$ are
assumed to be of class\footnote{ Background about H\"older classes
appears in Section~1(IV). } $\mathcal{ C }^\mathcal{ R}$, where
$\mathcal{ R} = (\kappa, \alpha)$, $\kappa \geqslant 1$, $0 \leqslant
\alpha \leqslant 1$, $\mathcal{ R} = \infty$, $\mathcal{ R} = \omega$
or $\mathcal{ R} = \mathcal{ A} lg$. Accordingly, $M$ is said to be
of class $\mathcal{ C}^{\mathcal{ A}lg}$ (real algebraic), $\mathcal{
C}^\omega$ (real analytic), $\mathcal{ C}^\infty$ or $\mathcal{ C}^{
\kappa, \alpha}$.

For $p\in M$, the smallest $J$-invariant subspace of the tangent space
$T_pM$ is given by $T_p^cM := T_pM \cap J T_pM$ and is called the {\sl
complex tangent space to $M$ at $p$}.

\def\thedefinition{2.2}\begin{definition}
{\rm
The submanifold $M$ is called{\rm :}

\smallskip

\begin{itemize}
\item[$\bullet$]
{\sl holomorphic} if $T_p^cM = T_pM$ at every point $p\in M$;

\smallskip

\item[$\bullet$]
{\sl totally real} if $T_p^cM = \{0\}$ at every point $p\in M$;

\smallskip

\item[$\bullet$]
{\sl generic} if $T_pM+JT_pM=T_p\C^n$ at every point $p\in M$;

\smallskip

\item[$\bullet$]
{\sl Cauchy-Riemann} (CR for short) if
the dimension of $T_p^cM$ is equal to a fixed constant at
every point $p\in M$.
\end{itemize} 
}
\end{definition}

For fundamentals about Cauchy-Riemann (CR for short) structures, we
refer the reader to \cite{ ch1989, ja1990, ch1991, bo1991,
ber1999, me2004a}. Here, we only summarize some
elementary useful properties. The two $J$-invariant spaces $T_pM \cap
JT_pM$ and $T_pM+ JT_pM$ are of even real dimension. We denote by
$m_p$ the integer $\frac{ 1}{ 2}\dim_\R (T_pM \cap JT_pM)$ and call it
the {\sl CR dimension of $M$ at $p$}. If $M$ is CR, $m_p \equiv m$ is
constant. Holomorphic, totally real and generic submanifolds are CR,
with $m= n - \frac{ 1}{ 2} \, d$, $m=0$ and $m=n-d$ respectively. If
$M$ is totally real and generic, $\dim_\R \, M=n$ and $M$ is called
{\sl maximally real}. We denote by $c_p$ the integer $n- \frac{ 1}{ 2}
\dim_\R( T_pM + JT_pM)$ and call it the {\sl holomorphic codimension
of $M$ at $p$}. It is constant if and only if $M$ is CR. Holomorphic,
totally real, generic and Cauchy-Riemann submanifolds are all CR and
have constant holomorphic codimensions $c = \frac{ 1}{ 2}\, d$, $c = d
- n$, $c = 0$ and $c = d-n+m$ respectively. Submanifolds of class
$\mathcal{ C}^{\kappa, \alpha}$ or $\mathcal{ C}^\infty$ will be
studied in Part~III.

Let $M$ or be a real algebraic ($\mathcal{ C}^{\mathcal{ A}lg}$) or
analytic ($\mathcal{ C }^\omega$) submanifold of $\C^n$ of (real)
codimension $d$ and let $p_0 \in M$. There exist complex algebraic or
analytic coordinates centered at $p_0$ and $\rho_1 >0$ such that $M$
is locally represented as follows.

\def\thetheorem{2.3}\begin{theorem}
{\rm (\cite{ ch1989, bo1991,
ber1999, me2004a})}

\smallskip

\begin{itemize}
\item[$\bullet$] 
If $M$ is {\rm holomorphic}, letting $m= n - \frac{ 1}{ 2} \, d
\geqslant 0$ and $c := \frac{ 1}{ 2}\, d$, then $m+ c = n$ and $M =
\left\{ (z, w_1) \in \Delta_{ \rho_1}^m \times \Delta_{ \rho_1}^c
: w_1 = 0 \right\}$.

\smallskip

\item[$\bullet$]
If $M$ is {\rm totally real}, letting $d_1 = 2n - d \geqslant 0$ and $c = d
- n \geqslant 0$, then $d_1 + c = n$ and $M = \left\{ (w_1, w_2) \in
\square_{ \rho_1 }^{ 2 d_1} \times \Delta_{ \rho_1}^c : {\rm Im }\,
w_1 = 0, w_2 = 0 \right\}$.

\smallskip

\item[$\bullet$]
If $M$ is {\rm generic}, letting $m = d - n$, then $m + d = n$ and 
\[
M = \left\{ (z, w) \in \Delta_{ \rho_1}^m \times (\square_{\rho_1 }^d
+ i\R^d) : {\rm Im}\, w = \varphi ( z, \bar z, {\rm Re} \, w)
\right\},
\]
for
some $\R^d$-valued algebraic or analytic map $\varphi$ satisfying
$\varphi (0) = 0$ whose power series converges normally in $\Delta_{
2\rho_1}^m \times \Delta_{ 2\rho_1}^m \times \square_{2\rho_1}^d$.

\smallskip

\item[$\bullet$]
If $M$ is {\rm Cauchy-Riemann}, letting $m = {\rm CRdim}\, M$, $c =
d-n+m \geqslant 0$, and $d_1 = 2n - 2m -d \geqslant 0$, then $m + d_1 + c = n$
and
\[
\aligned
M=
\big\{
(z, w_1, w_2) 
&
\in \Delta_{\rho_1}^m \times \left( \square_{ \rho_1 }^{d_1} + i \R^{
d_1} \right) \times \Delta_{ \rho_1}^c:
\\
&
{\rm Im}\, w_1 = \varphi_1 ( z,
\bar z, {\rm Re} \, w_1), \ w_2 = 0 \big\}, 
\endaligned
\]
for some $\R^{ d_1}$-valued algebraic or analytic map $\varphi_1$
satisfying $\varphi_1 (0) = 0$ whose power series converges normally
in $\Delta_{ 2\rho_1 }^m \times \Delta_{ 2\rho_1}^m \times \square_{2
\rho_1 }^{ d_1}$.
\end{itemize}

\smallskip

A further linear change of coordinates may yield $d\varphi (0) = 0$
and $d\varphi_1 (0) = 0$.
\end{theorem}

A CR algebraic or analytic manifold being generic in some local
complex manifold of (smaller) dimension $n-c$, called its {\sl
intrinsic complexification}, in most occasions, questions, results and
articles, one deals with generic manifolds. In this chapter, all
generic submanifolds will be of positive codimension $d\geqslant 1$ and of
positive CR dimension $m\geqslant 1$.

\subsection*{ 2.4.~Algebraic and analytic generic submanifolds and
their extrinsic complexification} Let $M$ be generic, represented by
${\rm Im}\, w= \varphi ( z, \bar z, {\rm Re}\, w)$. The implicit
function theorem applied to the vectorial equation $\frac{ w-\bar w}{
2i}= \varphi \left( z, \bar z, \frac{ w+\bar w}{ 2} \right)$, enables
to solve the variables $\bar w \in \C^d$, or the variables $w \in
\C^d$. This yields the so-called {\em complex defining equations} for
$M$, most useful in applications, as stated just below. Recall that,
given a power series $\Phi (t) = \sum_{ \gamma \in \N^n}\,
\Phi_\gamma\, t^\gamma$, $t \in \C^n$, $\Phi_\gamma \in \C$, $\gamma
\in \N^n$, one defines the series $\overline{ \Phi} (t) := \sum_{
\gamma \in \N^n }\, \overline{ \Phi }_\gamma\, t^\gamma$ by
conjugating only its complex coefficients. Then $\overline{ \Phi(t)}
\equiv \overline{ \Phi} (\bar t)$, a frequently used property.

\def\thetheorem{2.5}\begin{theorem} 
{\rm (\cite{ ber1999, me2004a})} A local generic real
algebraic or analytic $d$-codimensional generic submanifold $M \cap
\Delta_{\rho_1}^n$ may be represented by $\bar w = \Theta ( \bar z,z,
w)$, or equivalently by $w = \overline{ \Theta} (z,\bar z,\bar w)$,
for some complex algebraic or analytic $\C^d$-valued map $\Theta$
whose power series converges normally in $\Delta_{ 2\rho_1}^m \times
\Delta_{ 2\rho_1}^m \times \Delta_{ 2\rho_1}^d$, with $\rho_1
>0$. Here, $\Theta$ and $\overline{ \Theta}$ satisfy the two {\rm
(}equivalent by conjugation{\rm )} vectorial functional equations{\rm
:}
\[
\left\{
\aligned
\bar w \equiv 
& \
\Theta(\bar z, z, \overline{\Theta}(z,\bar z, \bar w)), 
\\
w\equiv 
& \
\overline{\Theta}(z,\bar z, 
\Theta(\bar z, z, w)).
\endaligned\right.
\]
Conversely, if such a $\C^d$-valued map $\Theta$ satisfies the
above, the set $M := \{(z,w)\in \Delta_{
\rho_1}^n : \, \bar w = \Theta (\bar z, z, w)\}$ is a {\rm real} local
generic submanifold of codimension $d$.
\end{theorem}

The coordinates $(z, w) \in \C^m \times \C^d$ will also be
denoted by $t
\in \C^n$. Let $\tau = (\zeta, \xi) \in \C^m \times \C^d$ be new
independent complex variables. Define the {\sl extrinsic
complexification $\mathcal{ M} = (M )^c$ of $M$} to be the complex
algebraic or analytic $d$-codimensional submanifold of $\C^n\times
\C^n$ defined by the vectorial equation $\xi - \Theta ( \zeta, t) = 0$
(the map $\Theta$ being analytic, we may indeed substitute $\zeta$
for $\bar z$ in its power series). We also write $\tau = (\overline{
t})^c$. Observe that $M$ identifies with the intersection $\mathcal{ M}
\cap \{ \tau = \bar t\}$.

\def\thelemma{2.6}\begin{lemma}
{\rm (\cite{ me2004a, me2005})} There exists an
invertible $d \times d$ matrix $a (t, \tau )$ of algebraic or analytic
power series converging normally in $\Delta_{ 2 \rho_1 }^n \times
\Delta_{ 2 \rho_1 }^n$ such that $w- \overline{ \Theta}( z, \tau)
\equiv
a(t, \tau)\, [ \xi - \Theta( \zeta,t) ]$.
\end{lemma}

Thus, $\mathcal{ M}$ is equivalently defined by $w - \overline{
\Theta} ( z, \tau) = 0$.

\subsection*{2.7.~Complexified Segre varieties and complexified CR 
vector fields} Let $\tau_ p, t_p \in \Delta_{ \rho_1 }^n$ be fixed and
define the {\sl complexified Segre varieties} $\mathcal{ S}_{\tau_p}$
and the {\sl complexified conjugate Segre varieties} $\underline{
\mathcal{ S}}_{ t_p}$ by:
\[
\left\{
\aligned
\mathcal{ S}_{\tau_p}
& \
:=\left\{
(t,\, \tau)\in 
\Delta_{\rho_1}^n \times \Delta_{ \rho_1}^n : \ 
\tau= \tau_p, \
w= \overline{\Theta}( z, \, \tau_p)
\right\}
\ \ \ \ \ {\rm and}
\\
\underline{\mathcal{ S}}_{t_p}
& \
:=\left\{ (t,\, \tau)\in
\Delta_{\rho_1}^n \times \Delta_{ \rho_1}^n
: \ 
t= t_p, \ 
\xi = 
\Theta (\zeta,\, t_p)
\right\}.
\endaligned\right.
\]
Geometrically, $\mathcal{ S}_{\tau_p} = \mathcal{ M} \cap \{ \tau =
\tau_p\}$ and $\underline{ \mathcal{ S}}_{t_p} = \mathcal{ M} \cap \{
t = t_p\}$. We draw a diagram.

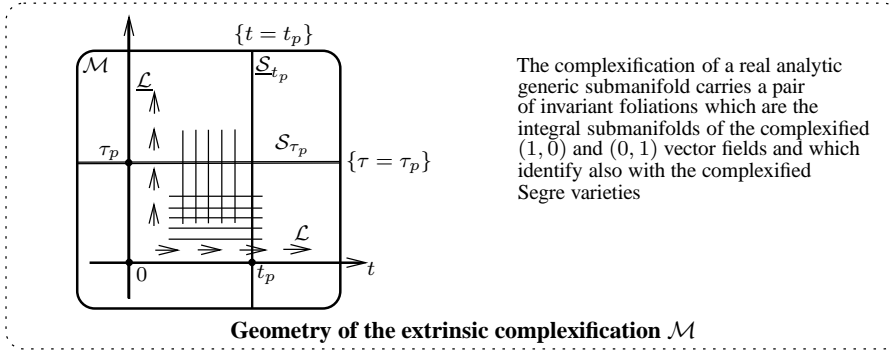
\begin{center}
\input complexification.pstex_t
\end{center}

We warn the reader that 
\[
\dim_\C \, \mathcal{ M} - \dim_\C \, 
\mathcal{ S}_{ \tau_p } - \dim_\C
\, \underline{ \mathcal{ S }}_{ t_p}
= d \geqslant 1,
\]
so that the ambient codimension $d$ of the unions of $\mathcal{
S}_{\tau_p}$ and of $\underline{ \mathcal{ S}}_{ t_p}$ is invisible in
this picture; one should imagine for instance that $\mathcal{M}$
is the three-dimensional physical space equipped with a pair of
foliations by horizontal orthogonal real lines.

Next, define two collections of complex vector fields:
\[
\left\{
\aligned
\mathcal{L}_k:= & \ 
{\partial \over\partial z_k}+
\sum_{j=1}^d\,
\frac{\partial\overline{\Theta}_j}{ \partial z_k}
(z, \zeta, \xi)\,
\, {\partial\over\partial w_j}, \ \ \ \ 
k=1,\dots,m, \ \ \ \ \ \ {\rm and} \\
\underline{\mathcal{L}}_k:= & \
{\partial\over\partial \zeta_k}+\sum_{j=1}^d\,
\frac{ \partial \Theta_j}{ \partial \zeta_k}
(\zeta, z, w)\, \,
{\partial\over\partial \xi_j}, \ \ \ \
k=1,\dots,m.
\endaligned\right.
\]
One verifies that $\mathcal{ L}_k \left( w_j- \overline{ \Theta}_j (z,
\zeta, \xi) \right) \equiv 0$, which shows that the $\mathcal{ L}_k$
are tangent to $\mathcal{ M}$. Similarly, $\underline{ \mathcal{ L}}_k
\left( \xi_j - \Theta_j (\zeta, z, w) \right) \equiv 0$, so the
$\underline{ \mathcal{ L}}_k$ are also tangent to $\mathcal{ M}$. In
addition, $[\mathcal{ L}_k,\, \mathcal{ L}_{k'}] =0$ and $[\underline{
\mathcal{ L}}_k,\, \underline{ \mathcal{ L}}_{ k'}] =0$ for $k,\, k'=
1,\dots, m$, so the theorem of Frobenius applies. In fact, the
$m$-dimensional integral submanifolds of the two collections $\{
\mathcal{ L}_k \}_{ 1\leqslant k\leqslant m}$ and $\{ \underline{ \mathcal{
L}}_k \}_{ 1\leqslant k\leqslant m}$ are the $\mathcal{ S}_{\tau_p}$ and the
$\underline{ \mathcal{ S}}_{ t_p}$. In summary, $\mathcal{ M}$ carries
a fundamental pair of foliations.

Observe that the vector fields $\mathcal{ L}_k$ are the
complexifications of the vector fields $L_k := \frac{ \partial
}{\partial z_k} + \sum_{ j=1}^d\, \frac{ \partial \overline{
\Theta}_j}{ \partial z_k} (z, \bar z, \bar w) \, \frac{ \partial}{
\partial w_j}$, $k=1, \dots, m$, that generate the holomorphic tangent
bundle $T^{ 1, 0} M$. A similar observation applies to the vector
fields $\underline{ \mathcal{ L}}_k$.

In general (unless $M$ is Levi-flat), the total collection $\{
\mathcal{ L}_k, \underline{ \mathcal{ L}}_k \}_{ 1\leqslant k\leqslant
m}$ does {\it not}\, enjoy the Frobenius property. In fact, the
noncommutativity of this system of $2m$ vector fields is at the very
core of Cauchy-Riemann geometry.

To apply Theorem~1.13, introduce the ``multiple'' flows of the two
collections $\{ \mathcal{ L}_k \}_{1 \leqslant k \leqslant m}$ and $\{
\underline{ \mathcal{L }}_{ k} \}_{ 1\leqslant k \leqslant m}$. If
$p\in \mathcal{ M }$ has coordinates $(z_p, w_p, \zeta_p, \xi_p) \in
\Delta_{ \rho_1 }^m \times \Delta_{ \rho_1}^d \times \Delta_{
\rho_1}^m \times \Delta_{ \rho_1}^d$ satisfying $w_p = \overline{
\Theta} ( z_p, \zeta_p, \xi_p)$ and $\xi_p = \Theta ( \zeta_p, z_p,
w_p)$ and if $z_1 := (z_{ 1, 1}, \dots, z_{1, m}) \in\C^m$ is a small
``multitime'' parameter, define the ``multiple'' flow of $\mathcal{
L}$ by:
\def\theequation{2.8}\begin{equation}
\aligned
\mathcal{L}_{z_1}(z_p, w_p, \zeta_p, \xi_p) 
& 
:=
\exp\left(z_1\mathcal{L}\right)(p)
\\
&
:=
\exp\left(z_{1, 1}\mathcal{L}_1(\cdots(\exp(z_{1, m}\mathcal{L}_m(
p)))\cdots)\right) 
\\
&
:=\left(
z_p+z_1, \overline{\Theta}(z_p+z_1, \zeta_p, \xi_p), 
\zeta_p, \xi_p\right).
\endaligned
\end{equation}
Of course, $\mathcal{ L}_{ z_1}(p) \in \mathcal{ M}$. Similarly, for
$p\in \mathcal{ M}$ and $\zeta_1 \in \C^m$, defining:
\def\theequation{2.9}\begin{equation}
\underline{\mathcal{L}}_{\zeta_1}(z_p,w_p,\zeta_p,\xi_p):=
(z_p,w_p,\zeta_p+\zeta_1,\Theta(\zeta_p+\zeta_1,z_p,w_p)), 
\end{equation}
we have $\underline{ \mathcal{ L}}_{\zeta_1}(p) \in \mathcal{
M}$. Clearly, $(p, z_1) \mapsto \mathcal{ L}_{z_1} (p)$ and $(p,
\zeta_1) \mapsto \underline{ \mathcal{ L}}_{\zeta_1} (p)$ are complex
algebraic or analytic local maps.

\subsection*{2.10.~Segre chains}
Let us start from $p = 0$ being the origin and move vertically along
the complexified conjugate Segre variety $\underline{ \mathcal{ S}
}_0$ of a height $z_1\in \C^m$, namely let us consider the point
$\underline{ \mathcal{ L }}_{ z_1}(0)$, which we shall also denote by
$\underline{ \Gamma}_1 (z_1)$. We have $\underline{ \Gamma}_1
(0)=0$. Let $z_2 \in \C^m$. Starting from the point $\underline{
\Gamma}_1 (z_1)$, let us move horizontally along the complexified
Segre variety of a length $z_2 \in \C^m$, namely let us consider the
point
\[
\underline{\Gamma}_2(z_1,z_2):=
\mathcal{L}_{z_2}(\underline{\mathcal{L}}_{z_1}(0)).
\]
Next, define $\underline{ \Gamma}_3 (z_1, z_2, z_3) := \underline{
\mathcal{ L}}_{ z_3}( \mathcal{ L}_{ z_2} (\underline{ \mathcal{ L}
}_{ z_1} (0 )))$, and then
\[
\underline{ \Gamma}_4 (z_1,
z_2, z_3, z_4) := \mathcal{ L}_{ z_4}( \underline{ \mathcal{ L}}_{
z_3}( \mathcal{ L}_{ z_2} (\underline{ \mathcal{ L}}_{ z_1} (0 )))),
\]
and so on. We draw a diagram:

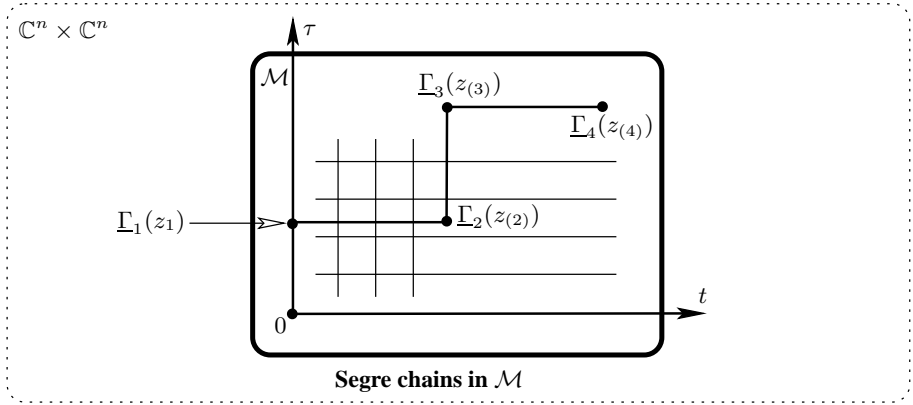
\begin{center}
\input new-orbit.pstex_t
\end{center}

By induction, for every $k\in \N$, $k\geqslant 1$, we obtain a local
complex algebraic or analytic map $\underline{ \Gamma }_k (z_1,
\dots,z_k)$, valued in $\mathcal{ M}$, defined for sufficiently
small $z_1, \dots, z_k \in \C^m$ which satisfies $\underline{
\Gamma}_k(0, \dots, 0) =0$. The abbreviated notation $z_{ (k)}
:=(z_1, \dots, z_k) \in \C^{ mk}$ will be used. The map
$\underline{ \Gamma}_k$ is called the {\sl $k$-th conjugate Segre
chain} (\cite{ me2004a, me2005}).

If we had conducted this procedure by starting with $\mathcal{ L}$
instead of starting with $\underline{ \mathcal{ L }}$, we would have
obtained maps $\Gamma_1( z_1):= \mathcal{ L}_{ z_1}(0)$, $\Gamma_2
(z_{(2)}) := \underline{ \mathcal{ L}}_{ z_2} (\mathcal{
L}_{z_1}(0))$, {\it etc.}, and generally $\Gamma_k(z_{(k)})$. The map
$\Gamma_k$ is called the {\sl $k$-th Segre chain}.

There is a symmetry relation between $\Gamma_k$ and $\underline{
\Gamma}_k$. Indeed, let $\overline{ \sigma}$ be the antiholomorphic
involution of $\C^n\times \C^n$ defined by $\overline{ \sigma} (t,
\tau):=(\bar \tau, \bar t)$. Since we have $w = \overline{ \Theta}( z,
\zeta, \xi )$ if and only if $\xi = \Theta( \zeta, z, w)$, this
involution is a bijection of $\mathcal{M}$. Applying $\overline{
\sigma}$ to the definitions~\thetag{ 2.8} and~\thetag{ 2.9} of the
flows of $\mathcal{ L}$ and of $\underline{ \mathcal{ L}}$, one may
verify that $\overline{ \sigma} ( \mathcal{ L}_{ z_1} (p)) =
\underline{ \mathcal{ L }}_{ \bar z_1} (\overline{ \sigma} (p))$. It
follows the general symmetry relation $\overline{ \sigma} \left(
\Gamma_k( z_{ (k) }) \right) = \underline{ \Gamma}_k \left( \overline{
z_{ (k) }} \right)$. Thus, $\Gamma_k$ and $\underline{ \Gamma}_k$ have
the same behavior.

\subsection*{2.11.~Minimality of 
$\mathcal{ M}$ at the origin and complexified local CR orbits} Since
$\Gamma_k (0) = \underline{ \Gamma }_k (0) = 0$, for every integer $k
\geqslant 1$, there exists $\delta_k >0$ sufficiently small such that
$\Gamma_k( z_{ (k) })$ and $\underline{ \Gamma }_k (z_{ (k)})$ are
well defined and belong to $\mathcal{ M}$, at least for all $z_{ (k)}
\in \Delta_{ \delta_k}^{ mk}$. To fiw ideas, it will be convenient to
consider that $\Delta_{ \delta_k}^{ mk}$ is the precise domain of
definition of $\Gamma_k$ and of $\underline{ \Gamma}_k$. We aim to
apply the procedure of Theorem~1.13 to the system $\LL^0 := \big\{
\mathcal{ L}_1, \dots, \mathcal{ L}_m, \, \underline{ \mathcal{ L}}_1,
\dots, \underline{ \mathcal{ L}}_m \big\}$.

However, there is a slight (innocuous) difference: each multitime
${\sf t} = ({\sf t}_1, \dots, {\sf t}) \in \K^k$ had scalar components
${\sf t}_i \in \K$, whereas now each $z_{ (k) } = (z_1, \dots, z_k)
\in \C^{m k}$ has vectorial components $z_i \in \C^m$. It is easy to
see that both $\Gamma_1$ and $\underline{ \Gamma}_1$ are of constant
rank $m$. Also, both $\Gamma_2$ and $\underline{ \Gamma}_2$ are of
constant rank $2m$, since $\mathcal{ L}_1, \dots, \mathcal{ L}_m, \,
\underline{ \mathcal{ L}}_1, \dots, \underline{ \mathcal{ L}}_m$ are
linearly independent at the origin. However, when passing to
(conjugate) Segre chains of length $\geqslant 3$, it is necessary to
speak of {\it generic}\, ranks and to introduce some combinatorial
integers $e_k \geqslant 1$. Justifying examples may be found in~\cite{
me1999, me2004a}.

\def\thetheorem{2.12}\begin{theorem} 
{\rm (\cite{ ber1996, ber1999, me1999, me2001a, me2004a})} There
exists an integer $\nu_0$ with 
$1 \leqslant \nu_0 \leqslant d$ and, for
$k= 3, \dots, \nu_0 + 1$, integers $e_k$ with 
$1 \leqslant e_k \leqslant
m$ such that the following nine properties hold true.

\begin{itemize}

\smallskip\item[{\bf (1)}]
For every $k = 3, \dots, \nu_0 + 1$, the two maps $\Gamma_k$ and
$\underline{ \Gamma}_k$ are of generic rank equal to $2m + e_3 +
\cdots + e_k$. In the special case $\nu_0 = 1$, the $e_k$ are
inexistent\footnote{ One may set
$e_1 := m$ and $e_2 := m$ in any case.} and nothing is stated.

\smallskip\item[{\bf (2)}]
For every $k \geqslant \nu_0 + 1$, both $\Gamma_k$ and $\underline{
\Gamma}_k$ are of fixed, stabilized generic rank equal to $2m + e$,
where
\[
e 
:= 
e_3 +
\cdots + e_{ \nu_0} \leqslant d.
\]

\smallskip\item[{\bf (3)}]
Setting $\mu_0 := 2\nu_0 + 1$, there exist two points $z_{ (\mu_0 )}^*
\in \Delta_{ \delta_{\mu_0} }^{ m\mu_0}$ and $\underline{ z}_{ (\mu_0
)}^* \in \Delta_{ \delta_{ \mu_0} }^{ m\mu_0 }$ satisfying $\Gamma_{
\mu_0} (z_{ (\mu_0 )}^* ) = 0$ and $\underline{ \Gamma}_{ \mu_0} (
\underline{ z}_{ (\mu_0 )}^* ) = 0$ which are arbitrarily close to the
origin in $\Delta_{ \delta_{ \mu_0} }^{ m \mu_0}$ such that $\Gamma_{
\mu_0}$ and $\underline{ \Gamma}_{ \mu_0}$ are of constant rank $2m +
e$ in neighborhoods $\omega^*$ and $\underline{ \omega}^*$ of $z_{
(\mu_0)}^*$ and of $\underline{ z}_{ (\mu_0)}^*$. The images $\Gamma_{
\mu_0} (\omega^*)$ and $\underline{ \Gamma}_{ \mu_0} (\underline{
\omega}^*)$ then constitute two pieces of local $\K$-algebraic or
analytic submanifold of dimension $2m + e$ contained in $\mathcal{
M}$.

\smallskip\item[{\bf (4)}]
Both $\Gamma_{ \mu_0} (\omega^*)$ and $\underline{ \Gamma}_{ \mu_0}
(\underline{ \omega}^*)$ enjoy the most important property that all
vector fields $\mathcal{ L}_1, \dots, \mathcal{ L}_m, \, \underline{
\mathcal{ L}}_1, \dots, \underline{ \mathcal{ L}}_m$ are tangent to
$\Gamma_{ \mu_0} (\omega^*)$ and to $\underline{ \Gamma}_{ \mu_0}
(\underline{ \omega}^*)$.

\smallskip\item[{\bf (5)}]
$\Gamma_{ \mu_0} (\omega^*)$ and $\underline{ \Gamma}_{ \mu_0}
(\underline{ \omega}^*)$ coincide together in a neighborhood of $0$ in
$\mathcal{ M}$. 

\smallskip\item[{\bf (6)}]
Denoting by 
\[
\mathcal{ O}_{ \mathcal{ L}, \underline{ \mathcal{ L}}}
(\mathcal{ M}, 0) 
\]
this common local piece of complex analytic
submanifold of $\mathcal{ M}$, it is algebraic provided that the flows
of $\big\{ \mathcal{ L}_1, \dots, \mathcal{ L}_m, \, \underline{
\mathcal{ L}}_1, \dots, \underline{ \mathcal{ L}}_m \big\}$ are
themselves algebraic.

\smallskip\item[{\bf (7)}]
Every local complex analytic or algebraic submanifold $\mathcal{ N}
\subset \mathcal{ M}$ passing through the origin to which $\mathcal{
L}_1, \dots, \mathcal{ L}_m, \, \underline{ \mathcal{ L}}_1, \dots,
\underline{ \mathcal{ L}}_m$ are all tangent must contain $\mathcal{
O}_{ \mathcal{ L}, \underline{ \mathcal{ L}}} (\mathcal{ M}, 0)$ in a
neighborhood of the origin.

\smallskip\item[{\bf (8)}]
The integers $\nu_0$, $e_3, \dots, e_{ \mu_0}$ and $e$ are
biholomorphic invariants of $\mathcal{ M}$.

\smallskip\item[{\bf (9)}]
$\Gamma_{ \mu_0} (\omega^*)$ and $\underline{ \Gamma}_{ \mu_0}
(\underline{ \omega}^*)$ also coincide {\rm (}in a neighborhood of the
origin{\rm )} with the Nagano leaf of the system $\big\{ \mathcal{
L}_1, \dots, \mathcal{ L}_m, \, \underline{ \mathcal{ L}}_1, \dots,
\underline{ \mathcal{ L}}_m \big\}$, as it was constructed in
Theorem~1.5.

\end{itemize}\smallskip

\end{theorem}

As in \cite{ me2004a, me2005} (with different notations), 
the integer $\nu_0$ 
will be called the {\sl Segre type of $M$}. 

The ``orbit notation'' $\mathcal{
O }_{ \mathcal{ L}, \underline{ \mathcal{ L} }} ( \mathcal{ M }, 0)$
anticipates the presentation and the notation of Section~1(III).
We will abandon Lie brackets and Nagano leaves.

The complex vector fields $L_k := \frac{ \partial }{ \partial z_k} +
\sum_{ j=1 }^d\, \frac{\partial \overline{ \Theta}_j }{ \partial z_k}
(z, \bar z, \bar w) \, \frac{\partial }{ \partial w_j}$, $k=1, \dots,
m$, are tangent to $M$ of equations $w_j = \overline{ \Theta}_j ( z,
\bar z, \bar w)$, $j=1, \dots, d$; their conjugates $\overline{ L}_k$
are also tangent to $M$; it follows that the real and imaginary parts
${\rm Re}\, L_k$ and ${\rm Im}\, L_k$ are also tangent to $M$. We may
then apply Theorem~1.13 to the system $\{ {\rm Re}\, L_k, {\rm Im}\,
L_k \}_{ 1\leqslant k \leqslant m}$, getting a certain real 
analytic local submanifold
$\mathcal{ O}_{L, \overline{ L}} ( M, 0)$ of $M$
passing through the origin. It will be
called the {\sl local CR orbit of the
origin in} $M$ (terminology of Part~III).

The relation between $\mathcal{ O }_{ \mathcal{ L}, \underline{
\mathcal{ L} }} ( \mathcal{ M }, 0)$ and $\mathcal{ O}_{L,
\overline{ L}} ( M, 0)$ is as follows (\cite{ ber1996, me1999,
me2001a, me2004a}). Let $\pi_t (t, \tau) := t$ and $\pi_\tau (t, \tau)
:= \tau$ denote the two canonical projections associated to the
product $\Delta_{ \rho_1 }^n \times \Delta_{ \rho_1}^n$. Let
$\underline{ A} := \big\{ 
(t, \tau) \in \Delta_{ \rho_1}^n \times \Delta_{
\rho_1 }^n : \tau = \bar t \big\}$ be the antiholomorphic diagonal.
Observe that $\pi_t ( \underline{ A } \cap \mathcal{ M}) = M$. 

\smallskip

\begin{itemize}

\item[$\bullet$]
The extrinsic complexification $\big[ \mathcal{ O}_{ L, \overline{ L
}} ( M, 0) \big]^c = \mathcal{ O}_{ \mathcal{ L}, \underline{
\mathcal{ L}}} ( \mathcal{ M}, 0)$.

\smallskip

\item[$\bullet$]
The projection $\pi_t \big( \underline{ A} \cap \mathcal{
O}_{\mathcal{ L}, \underline{ \mathcal{ L}}} (\mathcal{ M}, 0) \big) =
\mathcal{ O}_{L, \overline{ L}} ( M, 0)$.

\end{itemize}

\smallskip

Concerning smoothness, a striking subtelty happens: if $M$ is real
algebraic, although the local multiple flows of $\mathcal{ L}$ and of
$\underline{ \mathcal{ L}}$ are complex algebraic (thanks to their
definitions~\thetag{ 2.8} and~\thetag{ 2.9}), the flows of ${\rm Re}
\, L_k$ and of ${\rm Im}\, L_k$ are only real analytic in general.

\def\theexample{2.13}\begin{example}{\rm 
{\rm (\cite{ me2004a})} For the real algebraic hypersurface of
$\C^2$ defined by ${\rm Im}\, w = \sqrt{ 1 + z\bar z} - 1$, the vector
field $L := \frac{\partial }{\partial z} + i\bar z
\sqrt{ 1+z\bar z} \, \frac{\partial }{\partial w}$ generates $T^{
1, 0} M$ and the flow of $2\, {\rm Re}\, L$ involves the transcendent
function ${\rm Arcsh}$.
}\end{example}

\def\thetheorem{2.14}\begin{theorem}
{\rm (\cite{ ber1996, me2001a})}
The local CR orbit $\mathcal{ O}_{ L, \overline{ L}} (M, 0)$ is
real algebraic if $M$ is.
\end{theorem}

For the proof, assuming $M$ to be real algebraic, it is
impossible, because of the example, to apply the second phrase of
Theorem~1.13 {\bf (5)} to the system $\{ {\rm Re}\, L_k, {\rm Im}\,
L_k \}_{ 1\leqslant k \leqslant m}$. Fortunately, this phrase applies
to the complexified system $\{ \mathcal{ L}_k, \underline{ \mathcal{
L}}_k \}_{ 1 \leqslant k \leqslant m}$, whence $\mathcal{ O}_{
\mathcal{ L}, \underline{ \mathcal{ L}}} ( \mathcal{ M}, 0)$ is
algebraic, and then the local CR orbit $\mathcal{ O}_{L, \overline{
L}} ( M, 0) = \pi_t \left( \underline{ A} \cap \mathcal{ O}_{\mathcal{
L}, \underline{ \mathcal{ L}}} (\mathcal{ M}, 0) \right)$ is real
algebraic.

\def\thedefinition{2.15}\begin{definition}
{\rm 
The generic submanifold $M$ or its extrinsic complexification
$\mathcal{ M}$ is said to be {\sl minimal at the origin} if $\mathcal{
O}_{ L, \overline{ L}} (M, 0)$ contains a neighborhood of $0$ in $M$,
or equivalently if $\mathcal{ O}_{\mathcal{ L}, \underline{ \mathcal{
L}}} (\mathcal{ M}, 0)$ contains a neighborhood of $0$ in $\mathcal{
M}$.
}\end{definition}

The minimality at the origin of the algebraic or analytic complexified
local generic submanifold $\mathcal{ M}= (M )^c$ is a
biholomorphically invariant property{\rm ;} it neither depends on the
choice of defining equations nor on the choice of a conjugate pair of
systems of complex vector fields $\{ \mathcal{ L }_k \}_{ 1 \leqslant
k \leqslant m}$ and $\{ \underline{ \mathcal{ L} }_k \}_{ 1 \leqslant
k \leqslant m}$ spanning the tangent space to the two foliations.

Minimality at $0$ reads $e = d$ in Theorem~2.12. For a
hypersurface $M$, namely with $d = 1$, minimality at $0$ is equivalent
to $\nu_0 = 2$.

\subsection*{2.16.~Projections of the submersions $\Gamma_{ \mu_0}$
and $\underline{ \Gamma}_{ \mu_0}$} Let $\mu_0 = 2 \nu_0 + 1$ as in
Theorem~2.12. If $\mathcal{ M}$ is minimal at the origin, the two
local holomorphic maps
\[
\Gamma_{\mu_0} 
\ \ 
{\rm and} \ \ \, 
\underline{ \Gamma}_{\mu_0} \ : \ \ \ \ \
\Delta_{\delta_{\mu_0}}^{m\mu_0}
\longrightarrow \mathcal{ M}
\]
satisfy $\Gamma_{ \mu_0} ( z_{ (\mu_0)}^*) =0$ and $\underline{
\Gamma}_{ \mu_0} ( \underline{ z}_{ (\mu_0)}^*) =0$ and they are
submersive at $z_{ (\mu_0)}^*$ and at $\underline{ z}_{ (\mu_0)}^*$.

Consider the two projections $\pi_t (t, \tau) := t$ and $\pi_\tau ( t,
\tau) := \tau$ and four compositions $\pi_t \left( \Gamma_{\mu_0}
(z_{(\mu_0)}) \right)$, $\pi_t \left( \underline{ \Gamma}_{ \mu_0} (
z_{(\mu_0)}) \right)$ and $\pi_\tau \left( \Gamma_{ \mu_0}
(z_{(\mu_0)}) \right)$, $\pi_\tau \left( \underline{ \Gamma}_{ \mu_0}
( z_{(\mu_0)}) \right)$. Since $\mu_0 = 2 \nu_0 + 1$ is odd, observe
that the composition 
$\underline{ \Gamma}_{ 2\nu_0 + 1} = \underline{ \mathcal{ L}}
(\cdots)$ ends with a $\underline{ \mathcal{ L}}$ and that $\Gamma_{
2\nu_0 + 1} = \mathcal{ L} (\cdots)$ ends with a $\mathcal{ L}$.
According to the two definitions of the flow maps, the coordinates
$(\zeta_p, \xi_p)$ are untouched in~\thetag{ 2.8} and the coordinates
$(z_p, w_p)$ are untouched in~\thetag{ 2.9}. It follows that
\[
\left\{
\aligned
\pi_t \left(
\underline{ \Gamma}_{2\nu_0+1}(z_{(2\nu_0+1)})
\right)
& \
\equiv
\pi_t \left(
\underline{ \Gamma}_{2\nu_0}(z_{(2\nu _0)})
\right) 
\ \ \ \ \
{\rm and} \\
\pi_\tau \left(
\Gamma_{2\nu_0+1}(z_{(2\nu_0+1)})
\right)
& \
\equiv
\pi_\tau \left(
\Gamma_{2\nu_0}(z_{(2\nu_0)})
\right).
\endaligned\right.
\]

\def\thecorollary{2.17}\begin{corollary} 
{\rm (\cite{ me1999, ber1999, me2004a})} If $M$ is minimal at the
origin, there exists a integer $\nu_0 \leqslant d +1$
{\rm (}the {\rm Segre
type} of $M$ at the origin{\rm )} 
and there exist points $\underline{ z}_{
(2 \nu_0 )}^* \in \C^{ 2m \nu_0}$ and $z_{ (2 \nu_0 )}^* \in \C^{ 2m
\nu_0}$ arbitrarily close to the origin, such that the two maps
\[
\left\{
\aligned
\Delta_{\delta_{2\nu_0}}^{m2\nu_0}
\ni z_{(2\nu_0)} 
& \
\longmapsto \pi_t \left(\underline{ 
\Gamma}_{2\nu_0} ( 
z_{ (2\nu_0)})\right) \in \C^n 
\ \ \ \ \ {\rm and} \\
\Delta_{\delta_{2\nu_0}}^{m2\nu_0}
\ni z_{(2\nu_0)} 
& \
\longmapsto \pi_\tau \left( 
\Gamma_{2\nu_0} (
z_{ (2\nu_0)})\right) \in \C^n
\endaligned\right.
\]
are of rank $n$ and send $\underline{ z}_{ (2
\nu_0)}^*$ and $z_{(2\nu_0)}^*$ to the origin.
\end{corollary}

\section*{ \S3.~Formal CR mappings, jets of Segre varieties \\
and CR reflection mapping}

\subsection*{ 3.1.~Complexified CR mappings respect pairs of foliations}
Let $n'\in \N$ with $n'\geqslant 1$ and let $M' \subset \C^{n'}$ be a
second algebraic or analytic generic submanifold of codimension
$d'\geqslant 1$ and of CR dimension $m' = n' - d'\geqslant 1$. Let $p'
\in M'$. There exist local coordinates $t' = (z', w') \in \C^{ m'}
\times \C^{ d'}$ centered at $p'$ in which $M'$ is represented by
$\bar w ' = \Theta ' (\bar z', t')$, or equivalently by $w' =
\overline{ \Theta} ' ( z', \overline{ t}')$. If $(\overline{ t}')^c =
\tau' = (\zeta', \xi') \in \C^{ m'} \times \C^{ d'}$, the extrinsic
complexification is represented by $\xi' = \Theta ' (\zeta', t')$, or
equivalently by $w' = \overline{ \Theta}' ( z', \tau')$. We shall
denote by $0'$ the origin of $\C^{ n'}$.

Let $t\in \C^n$ and let $h(t) = (h_1(t), \dots, h_{n'}(t)) \in \C \dl
t \dr^{ n'}$ be a formal power series mapping with no constant term,
{\it i.e.} $h(0) = 0'$; it may also be holomorphic namely $h (t) \in
\C \{ t\}^{ n'}$, or even (Nash) algebraic. We have $(\overline{
h(t)})^c = \overline{ h} ( (\overline{ t})^c) = \overline{ h}
(\tau)$. Define $h^c (t, \tau) := (h(t), \overline{ h} (\tau))$.

Set $r (t, \tau) := \xi - \Theta ( \zeta, t)$, set $\overline{ r}
(\tau, t) := w - \overline{ \Theta} ( z, \tau)$, set $r '( t', \tau ')
:= \xi' - \Theta' (\zeta', t')$ and set $\overline{ r} '( \tau', t')
:= w' - \overline{ \Theta} '( z', \tau')$. We say that the power
series mapping $h$ is a {\sl formal CR mapping} from $(M,0)$ to
$(M',0')$ if there exists a $d' \times d$ matrix of formal power
series $b(t, \bar t)$ such that
\[
r' \left( h(t),
\overline{ h} (\bar t) \right)
\equiv b(t, \bar t) \, r(t, \bar t)
\]
in $\C \dl t, \bar t \dr^{d'}$. By complexification, it follows that
$r' \left( h (t), \overline{ h} (\tau) \right) \equiv
b(t, \tau) \, r (t, \tau)$ in $\C \dl t, \tau
\dr^{d'}$, namely $h^c ( t, \tau) = (h (t), \overline{ h} (\tau))$
maps $(\mathcal{ M}, 0)$ formally to $(\mathcal{ M}', 0')$. By
Lemma~2.6, there exist two complex analytic invertible matrices $a(t,
\tau)$ and $a' \left( t', \tau' \right)$ satisfying~:
\[
\left\{
\aligned
r(t,\tau) 
&
\equiv
a (t, \tau) \, 
\overline{r}(\tau,t), 
\ \ \ \ \ \ \ \ \ \ \ \ \ \ \ \ \
r'\left(t',\tau' \right) 
\equiv
a'\left(t',\tau'\right) \, 
\overline{r}' \left(\tau',t'\right), \\
\overline{r}(\tau,t) 
&
\equiv \
\overline{ a} (\tau, t) \, 
r(t,\tau), 
\ \ \ \ \ \ \ \ \ \ \ \ \ \ \ \ 
\overline{r}'(\tau',t') 
\equiv \
\overline{ a}'\left(\tau', t' \right) \, 
r'\left(t',\tau'\right),
\endaligned\right.
\]
in $\C \dl t, \tau \dr^{ d}$ and in $\C \dl t', \tau ' \dr^{ d'}$. So,
to define a complexified formal CR mapping $h^c : (\mathcal{ M }, 0)
\mapsto_{ \mathcal{ F }} (\mathcal{ M }', 0')$, we get four vectorial
formal identities, each one implying the remaining three:
\[
\left\{
\aligned
r' \left( h(t),
\overline{ h} (\tau) \right) 
&
\equiv b(t, \tau) \, r(t, \tau), 
\ \ \ \ \ \ \ \ \ \ \ 
r' \left( 
h (t),\overline{ h}(\tau)
\right) 
\equiv \overline{ c}(\tau,t) \, 
\overline{r} (\tau, t), 
\\
\overline{ r}' \left(
\overline{ h} (\tau), h(t) \right)
&
\equiv
\overline{ b} (\tau, t) \,
\overline{ r} (\tau, t),
\ \ \ \ \ \ \ \ \ \ \ 
\overline{r}' 
\left(
\overline{ h} (\tau), h(t) \right) 
\equiv 
c(t,\tau) \, 
r (t, \tau).
\endaligned\right.
\]
Here, we have set $c(t, \tau) := \overline{ b}( \tau, t) \, a
(t, \tau)$. 

These identities are independent of the choice of local coordinates
and of local complex defining equations for $(M,0)$ and for
$(M',0')$. Since $h$ is not a true point-map, we write $h: (M, 0)
\to_{ \mathcal{ F}} (M', 0')$, the index $\mathcal{ F}$ being the
initial of {\sl F}ormal. If $h$ is convergent, it is a true point-map
from a neighborhood of $0$ in $M$ to a neighborhood of $0'$ in $M'$.

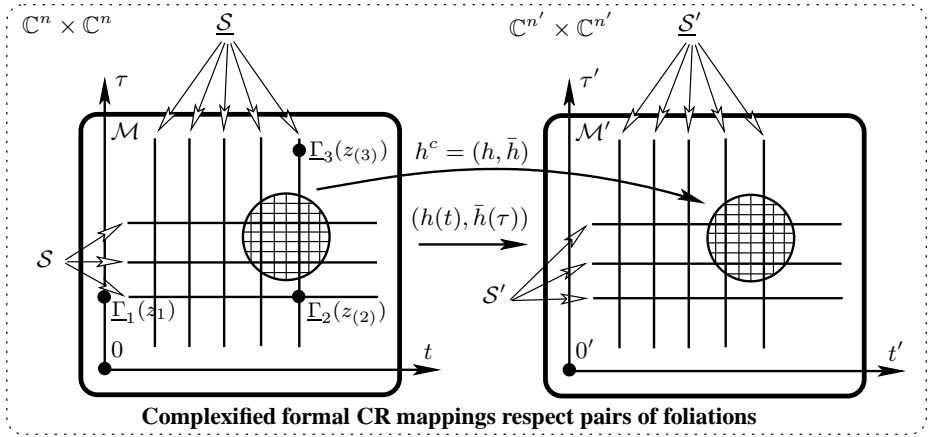
\begin{center}
\input stabilization-foliations.pstex_t
\end{center}

If $h$ is holomorphic in a polydisc $\Delta_{\rho_1}^n$, $\rho_1 >0$,
its extrinsic complexification $h^c$ sends both the $n$-dimensional
coordinate spaces $\{ t = {\rm cst.} \}$ and $\{ \tau = {\rm cst.} \}$
to the $n'$-dimensional coordinate spaces $\{ t' = {\rm cst.} \}$ and
$\{ \tau' = {\rm cst.} \}$.

Equivalently, $h^c$ maps complexified (conjugate) Segre varieties of
the source to complexified (conjugate) Segre varieties of the
target. {\it Some strong rigidity properties are due to the fact that
$h^c=(h, \bar h)$ must respect the two pairs of Segre foliations}.

The most important rigidity feature, 
called the {\em reflection principle}\footnote{ Other rigidity
phenomena are: parametrization of CR automorphism groups by a jet of
finite order, finiteness of their dimension, genericity of
nonalgebraizable CR submanifolds, genericity of CR submanifolds having
no infinitesimal CR automorphisms, {\it etc.}}, says that
the smoothness of $M$, $M'$ governs the smoothness of $h$:

\begin{itemize}

\smallskip\item[$\bullet$]
suppose that $M$ and $M'$ are real analytic and that $h (t) \in \C \dl
t \dr^{ n'}$ is only formal; statement: under suitable assumptions, $h
(t) \in \C \{ t\}^{ n'}$ is in fact convergent.

\smallskip\item[$\bullet$]
suppose that $M$ and $M'$ are real algebraic and that 
$h (t) \in \C \dl t \dr^{ n'}$ is only formal; statement: under
suitable assumptions, $h (t)$ is complex algebraic.

\end{itemize}\smallskip

After a mathematical phenomenon has been observed in a special, well
understood situation, the research has to focus attention on the
finest, the most adequate, the necessary and sufficient conditions
insuring it to hold true.

In this section, we aim to expose various possible assumptions for the
reflection principle to hold. Our goal is to provide a synthesis by
gathering various nondegeneracy assumptions which imply
reflection. For more about history, for other results, for complements
and for different points of view we refer to \cite{ pi1975, le1977,
we1977, we1978, pi1978, df1978, dw1980, df1988, br1988, br1990,
dp1993, dp1995, dp1998, ber1999, sh2000, ber2000, me2001a, me2002,
hu2001, sh2003, dp2003, ro2003, mmz2003b, er2004, me2005}.

The main theorems will be presented in \S3.19 and in \S3.22 below,
after a long preliminary. In these results, $M$ will always be assumed
to be minimal at the origin. Corollary~2.17 says already how to use
concretely this assumption: to show the convergence or the
algebraicity of a formal CR mapping $h : (M, 0) \mapsto_{ \mathcal{
F}} (M', 0')$, it suffices to establish that for every $k \in \N$, the
formal maps $z_{ (k)} \longmapsto_{ \mathcal{ F}} \, h \left(
\pi_t\left( \underline{ \Gamma }_k (z_{ (k)}) \right) \right)$ are
convergent or algebraic.

Before surveying recent results about the reflection principle
(without any indication of proof), we have to analyze thoroughly the
geometry of the target $\mathcal{ M}'$ and to present the
nondegeneracy conditions both on $\mathcal{ M}'$ and on $h$. Of
course, everything will also be meaningful for sufficiently smooth
($\mathcal{ C}^\infty$ or $\mathcal{ C}^\kappa$) local CR mappings, by
considering Taylor series.

These conditions are classical in local analytic geometry and they may
already be illustrated here with a plain formal map $h(t) \in \C \dl t
\dr^{ n'}$, not necessarily being CR.

\def\thedefinition{3.2}\begin{definition}{\rm 
A formal power series mapping $h : (\C^n, 0 ) \mapsto_{ \mathcal{
F}} (\C^{n'}, 0')$ with components $h_{i'} (t) \in \C\dl t\dr$,
$i'=1, \dots, n'$, is called

\begin{itemize}

\smallskip\item[{\bf (1)}]
{\sl invertible} if $n'=n$ and ${\rm det}\, ([\partial h_{i_1}/
\partial t_{i_2 }](0))_{1 \leqslant i_1, i_2\leqslant n}\neq 0$;

\smallskip\item[{\bf (2)}]
{\sl submersive} if $n\geqslant n'$ and there exist integers
$1\leqslant i (1)< \cdots < i(n') \leqslant n$ such that ${\rm det}\, 
([\partial h_{ i_1'}/ \partial t_{i (i_2' )}](0))_{1 \leqslant i_1',
i_2' \leqslant n'}\neq 0$;

\smallskip\item[{\bf (3)}]
{\sl finite} if the ideal generated by the components
$h_1(t), \dots, h_{n'}(t)$ is of finite codimension in $\C \dl t\dr$;
this implies $n'\geqslant n$;

\smallskip\item[{\bf (4)}]
{\sl dominating} if $n\geqslant n'$ and there exist integers 
$1\leqslant i(1)< \cdots < i({n'}) \leqslant n$ such that
${\rm det}([\partial h_{ i_1'}/\partial t_{ i(i_2')}] (t))_{1\leqslant
i_1',i_2' \leqslant n'}\not \equiv 0$ in $\C\dl t\dr$;

\smallskip\item[{\bf (5)}]
{\sl transversal} if there does not exist a nonzero power series
$G(t_1', \dots, t_{n'}')\in \C \dl t_1',\dots, t_{n'}'\dr$ such that
$G(h_1(t), \dots, h_{n'}(t))\equiv 0$ in $\C \dl t\dr$.
\end{itemize}\smallskip

}\end{definition}

It is elementary to see that invertibility implies submersiveness
which implies domination. Furthermore, if a formal power series is
either invertible, submersive or dominating, then it is
transversal. Philosophically, the ``distance'' between finite and
dominating or transversal is large, whereas the ``distance'' between
invertible and submersive or finite is ``small''.

\subsection*{ 3.3.~Jets of Segre varieties and Segre mapping} 
The target $M'$ concentrates all geometric conditions that are central
for the reflection principle. With respect to $\mathcal{ M}'$, the
complexified conjugate Segre variety associated to a fixed $t'$ is
$\underline{ \mathcal{ S }}_{t'}' := \{ (\zeta', \xi ') \in \C^{ n'} :
\, \xi' = \Theta ' (\zeta', t')\}$. Here, $\zeta'$ is a parametrizing
variable. For $k' \in \N$, define the {\sl morphism of $k'$-th jets of
complexified conjugate Segre varieties by}:
\[
\varphi_{k'}'(\zeta', t')
:=
J_{\tau'}^{k'}\underline{\mathcal{S}}_{t'}'
:=
\left( \zeta',\
\left({1\over\beta'!}\, \partial_{\zeta'}^{\beta'}
\Theta_{j'}'(\zeta', t')\right)_{
1\leqslant j'\leqslant d', 
\ 
\beta'\in \N^{m'},\ 
\vert\beta'\vert\leqslant
k'}\right).
\]
It takes values in $\C^{m'+ N_{d', m', k'}}$, with $N_{ d', m', k'} :=
d' \frac{ (m'+k' ) !}{ m'! \ k'!}$. If $k_1 ' \leqslant k_2 '$, we
have of course $\pi_{ k_2', k_1'} \circ \varphi_{ k_2'}' = \varphi_{
k_1' }'$.

As observed in~\cite{ dw1980}, the properties of this morphism govern
the various reflection principles. We shall say (\cite{ me2004a,
me2005}) that $M'$ (or equivalently $\mathcal{ M}'$) is:

\smallskip

\begin{itemize}
\item[{\bf (nd1)}]
{\sl Levi non-degenerate at the origin} if $\varphi_1 '$ is of rank
$m'+n'$ at $(\zeta', t')= (0', 0')$;

\smallskip\item[{\bf (nd2)}]
{\sl finitely nondegenerate at the origin} if there exists an integer
$\ell_0'$ such that $\varphi_{k'}'$ is of rank $n'+ m'$ at $(\zeta',
t')= (0', 0')$, for $k' = \ell_0'$, hence for all $k'\geqslant \ell_0'$;

\smallskip\item[{\bf (nd3)}]
{\sl essentially finite at the origin} if there exists an integer
$\ell_0'$ such that $\varphi_{k'} '$ is a finite holomorphic map at
$(\zeta', t')= (0', 0')$, for $k' = \ell_0'$, hence for all $k'\geqslant
\ell_0'$;

\smallskip\item[{\bf (nd4)}]
{\sl Segre nondegenerate at the origin} if there exists an integer
$\ell_0'$ such that the restriction of $\varphi_{ k'}'$ to the
complexified Segre variety $\mathcal{ S }_0'$ (of complex dimension
$m'$) is of generic rank $m'$, for $k' = \ell_0'$, hence
for all $k'\geqslant \ell_0'$;

\smallskip\item[{\bf (nd5)}]
{\sl holomorphically nondegenerate} if there exists an integer
$\ell_0'$ such that the map $\varphi_{ k'}'$ is of maximal possible
generic rank, equal to $m' +n'$, for $k' = \ell_0'$, hence
for all $k' \geqslant \ell_0'$.

\end{itemize}\smallskip

\def\thetheorem{3.4}\begin{theorem}
{\rm (\cite{ me2004a})}
These five conditions are biholomorphically invariant and{\rm :}
$
\text{\bf (nd1)} \Rightarrow 
\text{\bf (nd2)} \Rightarrow 
\text{\bf (nd3)} \Rightarrow 
\text{\bf (nd4)} \Rightarrow 
\text{\bf (nd5)}
$.
\end{theorem}

Being not punctual, the last condition {\bf (nd5)} is the finest:
as every condition of maximal generic rank, 
it propagates from any small open subet to
big connected open sets, thanks to 
the principle of analytic continuation. 
Notably, if a connected real analytic $M'$ is
holomorphically nondegenerate ``at'' a point, it is automatically
holomorphically nondegenerate ``at'' every point (\cite{ st1996,
ber1999, me2004a}).

\smallskip

To explain the (crucial) biholomorphic invariance of the jet map
$\varphi_{ k'}'$, consider a local biholomorphism $t' \mapsto h '(t')=
t''$, where $t', t'' \in \C^{ n'}$, that fixes the origin, $h_{i'}'
(t') \in \C \{ t'\}$, $h_{i'}' (0') = 0'$, for $i ' =1,\dots, n'$.
Splitting the coordinates $t'' = (z'', w'') \in \C^{ m'} \times \C^{
d'}$, the image $M''$ may be similarly represented by $\bar w'' =
\Theta '' ( \bar z'', t'')$ and there exists a $d' \times d'$ matrix
$b'(t', \tau')$ of local holomorphic functions such that
\[
r'' \big( h'(t'), \overline{
h}' ( \tau' ) \big) \equiv b '(t', \tau')\, r' (t', \tau')
\] 
in $\C \{ t', \tau'\}^{ d'}$, where $r_{j'}'(t', \tau') := \xi_{j'}' -
\Theta_{j'}' (\zeta', t')$ and $r_{j'}'' \left( t'', \tau'' \right) :=
\xi_{j'}'' - \Theta_{j'}'' \left( \zeta'', t'' \right)$, for $j'=1,
\dots , d'$. Setting $h'(t') := (f'(t'), g'(t'))\in \C \{ t'\}^{ m'}
\times \C\{ t'\}^{d'}$ and replacing $\xi'$ by $\Theta ' (\zeta', t')$
in the above equation, the right hand side vanishes identically 
(since $r' (t', \tau') = 
\xi' - \Theta' (\zeta', t')$ by definition) and
we obtain the following formal identity in
$\C\{ \zeta', t'\}^{ d'}$:
\[
\overline{ g}' \left( \zeta',\Theta'(\zeta',t') \right)
\equiv \Theta''\big(\overline{ f}'
(\zeta',\Theta'(\zeta',t')),
h'(t')\big).
\]
Some algebraic manipulations conduct to the
following.

\def\thelemma{3.5}\begin{lemma}
{\rm (\cite{ me2004a, me2005})}
For every $j'=1, \dots, d'$ and every $\beta' \in \N^{m'}$, there
exists a {\rm universal} rational map $Q_{j', \beta'}'$ whose
expression depends neither on $\mathcal{ M}'$, nor on $h'$, nor on
$\mathcal{ M}''$, such that the following identities in $\C \{
\zeta', t'\}$ hold true\,{\rm :}
\[
\small
\aligned
{}
&
\frac{1}{\beta'!}\, 
{\partial^{\vert \beta'\vert} \Theta_{j'}''\over
\partial (\zeta'')^{\beta'}}\left(
\overline{ f}' (\zeta',\Theta'(\zeta',t')), \
h'(t')\right)\equiv
\\
&
\equiv
Q_{j',\beta'}'\left(
\left(\partial_{\zeta'}^{\beta_1'}\Theta_{j_1'}'(\zeta',t')
\right)_{1\leqslant j_1'\leqslant d', \,
\vert \beta_1'\vert \leqslant \vert \beta'\vert}, \,
\left(\partial_{\tau'}^{\alpha_1'}
\overline{ h}_{i_1'}'(\zeta',\Theta'(\zeta',t'))
\right)_{1\leqslant i_1'\leqslant n', 
\vert \alpha_1' \vert\leqslant 
\vert \beta' \vert}\right) \\
&
=:
R_{j', \beta'}' \left(\zeta',\ 
\left(\partial_{\zeta'}^{\beta_1'}\Theta_{j_1'}(\zeta',t')
\right)_{1\leqslant j_1'\leqslant d', \,
\vert \beta_1'\vert \leqslant \vert \beta'\vert}
\right),
\endaligned
\]
where the last line defines $R_{j', \beta'}'$ by forgetting the jets
of $\overline{ h}'$. Here, the $Q_{j', \beta'}'$ are holomorphic in a
neighborhood of the constant jet 
\[
\left(( \partial_{\zeta'}^{
\beta_1'} \Theta_{j_1'}' (0,0) )_{ 1 \leqslant j_1' \leqslant d', \, \vert
\beta_1' \vert \leqslant \vert \beta'\vert}, \,(\partial_{\tau'}^{
\alpha_1' } \overline{ h}_{ i_1'} '(0,0) )_{1 \leqslant i_1'\leqslant n,
\vert \alpha_1' \vert \leqslant \vert \beta' \vert}\right). 
\]
Some symmetric relations hold after replacing $\Theta'$, $\Theta''$,
$\zeta'$, $t'$, $\overline{ f}'$, $h'$ by $\overline{ \Theta}'$,
$\overline{ \Theta}''$, $z'$, $\tau'$, $f'$, $\overline{ h}'$.
\end{lemma}

The existence of $R_{j', \beta'}'$ says that the following diagram
is commutative~:
\[
\diagram 
(\mathcal{M}',0') \rto^{(h')^c} 
\dto_{J_\bullet^{k'}\underline{\mathcal{S}}_\bullet'}
& (\mathcal{M}',0')
\dto^{J_\bullet^{k'}\underline{\mathcal{S}}_\bullet''} \\
\C^{m+N_{d',m',k'}} \rto^{ R_{k'}'((h')^c)} & 
\C^{m+N_{d',m',k'}}
\enddiagram,
\]
where the biholomorphic map $R_{k'}' ((h')^c)$, which depends on
$(h')^c$, is defined by its components $R_{j', \beta'}'$ for $j'=1,
\dots, d'$ and $\vert \beta' \vert \leqslant k'$. Thanks to the
invertibility of $h'$, the map $R_{k'}' ((h')^c)$ is also checked to
be invertible, and then the invariance of the five nondegeneracy
conditions {\bf (nd1)}, {\bf (nd2)}, {\bf (nd3)}, {\bf (nd4)} and {\bf
(nd5)} is easily established (\cite{ me2004a}).

\smallskip

We now present the {\sl Segre mapping} of $M'$. By developing the
series $\Theta_{ j'} ' (\zeta', t')$ in powers of $\zeta'$, we may
write the equations of $\mathcal{ M}'$ under the form $\xi_{ j'}' =
\sum_{\gamma ' \in \N^{ m'}}\, (\zeta')^{ \gamma'} \, \Theta_{j',
\gamma'}' (t')$ for $j' = 1, \dots, d'$. In terms of such a
development, the {\sl infinite Segre mapping of $M'$} is defined to be
the mapping
\[
\mathcal{Q}_\infty': 
\ \ 
\C^{n'}\ni t' 
\longmapsto
(\Theta_{j',\gamma'}'(t'))_{
1\leqslant j'\leqslant d',\ 
\gamma'\in\N^{m'}}
\in 
\C^\infty.
\]
Let $k'\in\N$. For finiteness reasons, 
it is convenient to truncate this infinite
collection and to define the 
{\sl $k'$-th Segre mapping of $M'$} by
\[
\mathcal{Q}_{k'}': 
\ \ 
\C^{n'}
\ni 
t' 
\longmapsto
(\Theta_{j',\gamma'}'
(t'))_{1\leqslant j'\leqslant d',\, 
\vert\gamma'\vert\leqslant k'}
\in 
\C^{N_{d',n',k'}},
\]
where $N_{d', n', k'} =d ' \, \frac{ (n' + k')! }{ n'! \ k'!}$. If
$k_2' \geqslant k_1'$, we have $\pi_{ k_2', k_1'} \big[ \mathcal{ Q}_{
k_2'}' (t') \big] = \mathcal{ Q}_{k_1'}'(t')$. One verifies (\cite{
me2004a}) the following characterizations.

\smallskip

\begin{itemize}
\item[{\bf (nd1)}]
$M'$ is {\sl Levi non-degenerate at the origin} if and only if
$\mathcal{ Q}_1 '$ is of rank $n'$ at $t'= 0'$.

\smallskip\item[{\bf (nd2)}]
$M'$ is {\sl finitely nondegenerate at the origin} if and only if
there exists an integer $\ell_0'$ such that $\mathcal{ Q}_{ k'}'$ is
of rank $n'$ at $t'= 0'$, for all $k' \geqslant \ell_0'$.

\smallskip\item[{\bf (nd3)}]
$M'$ is {\sl essentially finite at the origin} if there exists an
integer $\ell_0'$ such that $\mathcal{ Q}_{ k'}'$ is a finite
holomorphic map at $t'= 0'$, for all $k' \geqslant \ell_0'$.

\smallskip\item[{\bf (nd4)}]
$M'$ is {\sl Segre nondegenerate at the origin} if there exists an
integer $\ell_0'$ such that the restriction of $\mathcal{ Q}_{k'}'$ to
the complexified Segre variety $\mathcal{ S }_{0'}'$ (of complex
dimension $m'$) is of generic rank $m'$, for all $k' \geqslant \ell_0'$.

\smallskip\item[{\bf (nd5)}]
$M'$ is {\sl holomorphically nondegenerate} if there exists an integer
$\ell_0'$ such that the map $\mathcal{ Q}_{k'}'$ is of maximal
possible generic rank, equal to $n'$, for all $k' \geqslant \ell_0'$.

\end{itemize}\smallskip

\subsection*{3.6.~Essential holomorphic dimension and Levi multitype} 
Assume now that $M'$ is not nececessarily local, but connected.
Denote by $\ell_{ M'}'$ the smallest integer $k'$ such that the
generic rank of the jet mappings $(t', \tau') \mapsto J_{ \tau'}^{ k'}
\underline{ \mathcal{S }}_{ t'} '$ does not increase after $k'$ and
denote by $m'+ n_{M'}' \leqslant m '+ n'$ the (maximal) generic rank of
$(t',\tau' )\mapsto J_{\tau'}^{ \ell_{ M'}'} \underline{ \mathcal{ S
}}_ { t'}'$. Since $w' \mapsto \Theta ' (\zeta', z', w')$ is of rank
$d'$ according to Theorem~2.5, the (generic) rank of the zero-th
order jet map satisfies
\[
{\rm genrk}_\C\left((t ',\tau')\mapsto
J_{ \tau'}^0\underline{\mathcal{S}}_{t'}'
=
(\zeta',\Theta'(\zeta',z',w'))
\right)
=
m'+d'
=
n'.
\] 
Thus, $d' \leqslant n_{M'}' \leqslant n'$. It is natural to call
$n_{ M'}'$ the {\sl essential holomorphic dimension of $M'$}
because of the following.

\def\theproposition{3.7}\begin{proposition}
{\rm (\cite{ me2001a, me2004a})}
Locally in a neighborhood of a Zariski-generic point $p'\in M'$, the
generic submanifold $M'$ is biholomorphically equivalent to the
product $\underline{M }_{p' } '\times \Delta^{ n'- n_{M' }'}$, of a
generic submanifold $\underline{ M }_{p' }'$ of codimension $d'$ in
$\C^{ n_{ M'}'}$ by a complex polydisc $\Delta^{ n'- n_{M' }'}$.
\end{proposition}

Generally speaking, we may define $\lambda_{0, M'}' := {\rm genrk }_\C
\left( (t', \tau') \mapsto J_{ \tau' }^0 \underline{ \mathcal{S }}_{
t'}' \right) - m' = d'$ and for every $k'=1, \dots, \ell_{M' }'$,
\[
\lambda_{k',M'}':={\rm genrk}_\C\left(
(t',\tau')\mapsto 
J_{\tau'}^{k'}\underline{\mathcal{S}}_{t'}'
\right)
-
{\rm genrk}_\C\left(
(t',\tau')\mapsto 
J_{\tau'}^{k'-1}\underline{\mathcal{S}}_{t'}'
\right).
\]
One verifies (\cite{ me2004a}) that $\lambda_{ 1,M'}' \geqslant 1,\dots,
\lambda_{ \ell_{M'}' , M '}'\geqslant 1$.
With these definitions, we have the relations
\[
{\rm genrk}_\C\left(
(t',\tau')\mapsto 
J_{\tau'}^{k'}\underline{\mathcal{S}}_{t'}'
\right)=
m'
+
\lambda_{0,M'}'
+
\lambda_{1,M'}'
+\cdots+
\lambda_{k',M'}',
\]
for $k'= 0,1, \dots, \ell_{M'}'$ and
\[
{\rm genrk}_\C\left(
(t',\tau')\mapsto 
J_{\tau'}^{k'}\underline{\mathcal{S}}_{t'}'
\right)
=
m'+d'+\lambda_{1,M'}'+\cdots+\lambda_{\ell_{M'}',M'}'
=
m'
+
n_{M'}',
\]
for all $k'\geqslant \ell_{M'}'$. It follows that 
\[
\ell_{M'}'\leqslant \lambda_{1,{M'}}'
+\cdots+
\lambda_{\ell_{M'}',M'}'=
n_{M'}'
-
d'
\leqslant 
m'.
\]

\def\thetheorem{3.8}\begin{theorem}
{\rm (\cite{ me2004a})} 
Let $M'$ be a {\rm connected} real algebraic or analytic generic
submanifold in $\C^{n'}$ of codimension $d'\geqslant 1$ and of CR
dimension $m'= n'-d' \geqslant 1$. Then there exist well defined
integers $n_{M'}' \geqslant d'$, $\ell_{ M'}' \geqslant 0$, $\lambda_{
0, M'}' \geqslant 1$, $\lambda_{ 1, M'}' \geqslant 1$, $\dots$,
$\lambda_{ \ell_{M'}', M'}' \geqslant 1$ and a proper real algebraic
or analytic subvariety $E'$ of $M'$ such that for every point $p' \in
M' \backslash E'$ and for every system of coordinates $(z',w')$
vanishing at $p'$ in which $M'$ is represented by defining equations
$\bar w_{j'} = \Theta_{ j'} '(\bar z',t')$, $j'=1, \dots,d'$, then the
following four properties hold{\rm :}

\begin{itemize}

\smallskip\item[$\bullet$]
$\lambda_{0,M'}'= d'$, $d'\leqslant n_{M'}' \leqslant n'$ and $\ell_{ M'}' \leqslant
n_{ M'}' -d'$.

\smallskip\item[$\bullet$]
For every $k' =0, 1, \dots, \ell_{ M'}'$, the mapping of $k'$-th order
jets of the conjugate complexified Segre varieties $(t',\tau') \mapsto
J_{ \tau'}^{ k'} \underline{ \mathcal{S }}_{t'}' $ is of rank equal to
$m'+ \lambda_{0, M'}' +\cdots+ \lambda_{k',M'}'$ at $(t_{p'}', \bar
t_{p'}')= (0',0')$.

\smallskip\item[$\bullet$]
$n_{ M'}' =d' + \lambda_{1, M'}' +\cdots+ \lambda_{ \ell_{M'}' ,M'}'$
and for every $k'\geqslant \ell_{ M'}'$, the mapping of $k'$-th order jets
of the conjugate complexified Segre varieties $(t',\tau') \mapsto J_{
\tau'}^{ k'} \underline{ \mathcal{S}}_{ t'}'$ is of rank equal to $n_{
M'}'$ at $(0',0')$.

\smallskip\item[$\bullet$]
There exists a local complex algebraic or analytic change of
coordinates $t''=h'(t')$ fixing $p'$ such that the image $M_{p'}'':
=h'(M')$ is locally in a neighborhood of $p'$ the product $\underline{
M}_{p'}'' \times \Delta^{n'- n_{ M' }'}$ of a real algebraic or
analytic generic submanifold of codimension $d'$ in $\C^{ n_{ M'}'}$
by a complex polydisc $\Delta^{ n'-n_{ M'}'}$. Furthermore, at the
central point $\underline{p}' \in \underline{ M}_{ p'} '' \subset \C^{
n_{M'} '}$, the generic submanifold $\underline{ M}_{p'} ''$ is $\ell_{
M'} '$-finitely nondegenerate, hence in particular its essential
holomorphic dimension $n_{ \underline{ M }_{p'}''}'$ coincides with $n_{
M'}'$.
\end{itemize}

In particular, $M'$ is holomorphically nondegenerate if and only if
$n_{ M'}' = n'$ and in this case, $M'$ is finitely nondegenerate at
every point of the Zariski-open subset $M' \backslash E'$.
\end{theorem}

\subsection*{3.9.~CR-horizontal nondegeneracy conditions} 
As in \S3.1, let $h = h(t) \in \C \dl t\dr^{ n'}$ be a formal
CR mapping $(M, 0) \to_{ \mathcal{ F}} (M',0')$. Decompose $h
(t) = (f(t), g(t)) \in \C \dl t \dr^{ m'} \times \C \dl t \dr^{ d'}$,
as in the splitting $t ' = (z', w') \in \C^{ m'} \times \C^{
d'}$. Replacing $w$ by $\overline{ \Theta} (z, \tau)$ in the
fundamental identity $\overline{ r}' \left( \overline{ h} (\tau ), h
(t) \right) \equiv \overline{ b}( \tau, t) \, \overline{ r} (\tau,
t)$, the right hand side vanishes identically (since $\overline{ r}
(\tau, t) = w - \overline{ \Theta} (z, \tau)$ by definition), and we
get a formal identity in $\C \dl z, \tau \dr^{ d'}$:
\[
g
\big(z, \overline{ \Theta} (z, \tau) \big) \equiv
\overline{ \Theta} ' \left( f(z, \overline{
\Theta} (z, \tau)), \overline{ h}( \tau) \right).
\]
Setting $\tau := 0$, we get $g \left(z, \overline{ \Theta} (z, 0)
\right) \equiv \overline{ \Theta}' \big( f(z, \overline{ \Theta} (z,
0)), 0 \big)$. In other words, $h \vert_{ \mathcal{ S }_0 }$ maps
$\mathcal{ S }_0$ formally to $\mathcal{ S}_{ 0'} '$. The restriction $h
\vert_{ \mathcal{ S }_0}$ coincides with the formal map:
\[
\C^m \ni z \longmapsto_{\mathcal{ F}} \
\left(
f\left(
z, \overline{ \Theta} (z, 0)
\right), \ \
\overline{ \Theta}'\left(
f\left(
z, \overline{ \Theta} (z, 0)
\right), \ 0
\right)
\right)\in \C^{m'} \times \C^{ d'}.
\]
The rank properties of this formal map are the same as those of its
{\sl CR-horizontal part}:
\[
\C^m \ni z \longmapsto_{\mathcal{ F}} \
f\left(
z, \overline{ \Theta} (z, 0)
\right)\in \C^{m'}.
\]
The formal CR mapping $h$ is said (\cite{ me2004a}) to be:

\smallskip

\begin{itemize}
\item[{\bf (cr1)}]
{\sl CR-invertible} at the origin if $m'=m$ and if its CR-horizontal
part is a formal equivalence at $z=0$;

\smallskip

\item[{\bf (cr2)}]
{\sl CR-submersive} at the origin if $m' \leqslant m$ and if its
CR-horizontal part is a formal submersion at $z=0$;

\smallskip

\item[{\bf (cr3)}]
{\sl CR-finite} at the origin if $m' = m$ and if its CR-horizontal
part is a finite formal map at $z=0$, namely the quotient ring $\C \dl
z \dr / \left( f_{k'} ( z, \overline{ \Theta} ( z, 0))_{1\leqslant k'
\leqslant m'} \right)$ is finite-dimensional (the requirement $m' = m$
is necessary for the reflection principle below);

\smallskip

\item[{\bf (cr4)}]
{\sl CR-dominating} at the origin if $m' \leqslant m$ and if there exist
integers $1\leqslant k(1) < \cdots < k( {m'}) \leqslant m$ such that the
determinant ${\rm det}([\partial \phi_{ k_1'}/ \partial z_{
k(k_2' )}] ( z ))_{1 \leqslant k_1', k_2' \leqslant m'} \not \equiv 0$
does not vanish identically in $\C \dl z \dr$, where $\phi_{
k'} ( z) := f_{ k'} \big( z, \overline{ \Theta} (z, 0) \big)$;

\smallskip

\item[{\bf (cr5)}]
{\sl CR-transversal} at the origin if there does not exist a nonzero
formal power series $F' ( f_1, \dots, f_{m'}) \in \C \dl f_1, \dots,
f_{m'} \dr$ such that $F '(\phi_1( z), \dots, \phi_{m'} (
z))\equiv 0$ in $\C \dl z \dr$, where $\phi_{ k'} ( z) := f_{ k'}
\big( z, \overline{ \Theta} (z, 0) \big)$.
\end{itemize}

\smallskip
\noindent
One verifies (\cite{ me2004a}) biholomorphic invariance and the four
implications:
\[
\text{\bf (cr1)} \Rightarrow
\text{\bf (cr2)} \Rightarrow
\text{\bf (cr3)} \Rightarrow
\text{\bf (cr4)} \Rightarrow
\text{\bf (cr5)},
\]
provided that $m'=m$ in the second and in the third. By far, 
CR-transversality is the most general nondegeneracy condition.

\subsection*{3.10.~Nondegeneracy conditions for CR mappings}
This subsection explains how to synthetize the combinatorics of
various formal reflection principles published in the last decade.

As in \S3.1, let $h^c: (\mathcal{M}, 0) \to_{\mathcal{ F}}
(\mathcal{M}', 0)$ be a complexified formal CR mapping between two
formal, analytic or algebraic complexified generic submanifolds of
equations $0 = r (t, \tau):= \xi - \Theta (\zeta, t)$ and $0 = r' (t',
\tau'):= \xi ' - \Theta ' ( \zeta', t')$. By hypothesis, $r' (h(t),
\overline{ h} (\tau)) \equiv b(t, \tau)\, r(t, \tau)$. Denoting $h =
(f, g)\in \C^{ m'} \times \C^{ d'}$, replacing $\xi$ by $\Theta
(\zeta, t)$ in $r' (h(t), \overline{ h} (\tau)) \equiv b(t, \tau)\,
r(t, \tau)$ and developing $\Theta ' (\bar f, h) = \sum_{ \gamma ' \in
\N^{ m'}}\, \bar f^{\gamma '}\, \Theta_{ \gamma' }' (h)$, we start
with the following fundamental power series identity in $\C\dl \zeta,
t\dr^{ d'}$:
\[
\aligned
\bar g(\zeta,\Theta(\zeta,t))
&
\equiv
\Theta'
\big(
\bar f(\zeta,\Theta(\zeta,t)),h(t)
\big)
\\
&
\equiv
\sum_{\gamma'\in\N^{m'}}\,
\bar f(\zeta,\Theta(\zeta,t))^{\gamma'}\,
\Theta_{\gamma'}'(h(t)).
\endaligned
\]
Consider the $m$ complex vector fields $\underline{ \mathcal{ L}}_1,
\dots, \underline{ \mathcal{L }}_m$ tangent to $\mathcal{M}$ that were
defined in \S2.7. For every $\beta=(\beta_1, \dots, \beta_m) \in
\N^m$, define the multiple derivation $\underline{ \mathcal{
L}}^\beta= \underline{ \mathcal{ L}}_1^{ \beta_1} \cdots \underline{
\mathcal{ L}}_m^{ \beta_m}$. Applying them to the above $d'$ scalar
equations, observing that they do not differentiate the variables $t =
(z, w)$, we get, without writing the arguments:
\[
\underline{\mathcal{L}}^\beta\bar g_{j'}
- 
\sum_{\gamma'\in\N^{m'}}\,
\underline{\mathcal{L}}^\beta(\bar f^{\gamma'})\,
\Theta_{j',\gamma'}'(h)\equiv 0,
\]
for all $\beta \in \N^m$, all $j' =1, \dots,d'$ and all $(t, \tau)\in
\mathcal{ M}$. 

\def\thelemma{3.11}\begin{lemma}
{\rm (\cite{ me2004a, me2005})}
For every $i'=1,\dots,n'$ and every $\beta\in\N^m$, there exists a
polynomial $P_{i',\beta}$ in the jet $J_\tau^{\vert \beta \vert} \bar
h(\tau)$ with coefficients being power series in $(t,\tau)$ which
depend only on the defining functions $\xi_j-\Theta_j(\zeta,t)$ of
$\mathcal{M}$ and which can be computed by means of some combinatorial
formula, such that
\[
\underline{\mathcal{L}}^\beta \bar h_{i'}(\tau)
\equiv
P_{i',\beta}
\left(
t,\tau,J_\tau^{\vert \beta \vert} \bar h(\tau)
\right).
\]
\end{lemma}

\smallskip
\noindent
{\bf Convention 3.12.} Let $k,l \in\N$. On the complexification
$\mathcal{M}$, equipped with either the coordinates $(z, \tau)$ or
$(\zeta, t)$, which correspond to either replacing $w$ by $\overline{
\Theta} ( z, \tau)$ or $\xi$ by $\Theta (\zeta, t)$, we shall identify
(notationally) a power series written under the complete form
\[
R(t,\tau,J^k
h(t),J^l\bar h(\tau)),
\]
with a power series written under one of the following four forms:
\[
\aligned
\bullet & \ \ 
R\left(t,\zeta,\Theta(\zeta,t),J^k h(t), J^l\bar
h(\zeta,\Theta(\zeta,t))\right),\\
\bullet & \ \ 
R\left(t,\zeta,J^k h(t),
J^l\bar h(\zeta,\Theta(\zeta,t))\right),\\
\bullet & \ \
R\left(z,\overline{\Theta}(z,\tau),\tau,J^k
h(z,\overline{\Theta}(z,\tau)),
J^l\bar h(\tau)\right),\\
\bullet & \ \
R\left(z,\tau, J^kh(z,\overline{\Theta}(z,\tau)),
J^l\bar h(\tau)\right).
\endaligned
\]

\smallskip

Thanks to the lemma and to the convention, we may therefore write:
\def\theequation{3.13}\begin{equation}
\underline{\mathcal{L}}^\beta 
\big[
\bar g_{j'}(\tau)
-
\Theta_{j'}'(\bar f(\tau), h(t))
\big]
=:
R_{j',\beta}'
\big(
t,\tau,J_\tau^{\vert \beta \vert} 
\bar h(\tau): \, 
h(t)
\big)
\equiv
0,
\end{equation}
for $j'=1, \dots, d'$. Remind that $h(t)$ is not differentiated, since
the derivations $\underline{ \mathcal{ L }}^\beta$ involve only $\frac{
\partial }{ \partial \tau_i}$, $i= 1, \dots, n$. This is why we write
$h(t)$ after ``:''. Furthmerore, the identities ``$\equiv 0$'' are
understood ``on $\mathcal{ M}$'', namely as formal power series
identities in $\C\dl \zeta, t\dr$ after replacing $\xi$ by $\Theta(
\zeta,t)$ or equivalently, as a formal power series identities in
$\C \dl z , \tau\dr$ after replacing $w$ by $\overline{ \Theta }(z,
\tau)$.

To understand the reflection principle, it is important to observe
immediately that the smoothness of the power series $R_{j',\beta}'$ is
the minimum of the two smoothnesses of $M$ and of $M'$. For instance,
the power series $R_{j',\beta}'$ are all complex analytic if $M$ is
real analytic and if $M'$ is real algebraic, even if the power series
CR mapping $h(t)$ was assumed to be purely formal and 
nonconvergent. By a careful inspection of the application of the chain
rule in the development of the above equations~\thetag{ 3.13} ({\it
cf.} Lemma~3.11), we even see that each $R_{j',\beta}'$ is relatively
polynomial with respect to the derivatives of positive order $(
\partial_\tau^\alpha \bar h (\tau ))_{ 1\leqslant \vert \alpha \vert
\leqslant \vert \beta \vert}$.

\subsection*{3.14.~Nondegeneracy conditions for formal CR mappings} 
In the equations~\thetag{3.13}, we replace $h(t)$ by a new independent
variable $t' \in \C^{n'}$, we set $(t, \tau) = (0,0)$, and we define the
following collection of power series
\[
\Psi_{j',\beta}'(t')
:=
\Big[\underline{\mathcal{L}}^\beta\bar g_{j'}-
\sum_{\gamma'\in\N^{m'}}\, 
\underline{\mathcal{L}}^\beta(\bar f^{\gamma'})\, 
\Theta_{j',\gamma'}'(t')
\Big]_{
t=\tau=0}, 
\]
for $j'=1, \dots, d'$ and $\beta \in \N^m$. Here, if $\beta = 0$, we
mean that $\Psi_{j', 0 }' (t')= - \Theta_{ j' }'( 0, t')$. According
to~\thetag{ 3.13}, an equivalent definition is:
\[
\Psi_{j',\beta}'(t'):=
R_{j',\beta}'\big(0,0,J_\tau^{
\vert\beta\vert}\bar h(0):\,t'\big).
\]
Now, just before introducing five new nondegeneracy conditions, we
make a crucial heuristic remark. When $n = n'$, $m = m'$, 
$M = M'$ and $h = {\rm
Id}$, writing $T'$ instead of $t'$ the special variable above in order
to avoid confusion, we get for $j' = 1, \dots, d'$ and $\beta ' \in
\N^{ m'}$:
\[
\aligned
\Psi_{j',\beta'}'(T')
& \
=
\Big[\underline{
\mathcal{L}'}^{\beta'}\xi_{j'}'
-
\sum_{\gamma'\in\N^{m'}}\, 
\underline{\mathcal{
L}'}^{\beta'}(\zeta')^{\gamma'}\, 
\Theta_{j',\gamma'}'(T')
\Big]_{
t'=\tau'=0'}
\\
& \ 
=
\left[\underline{
\mathcal{L}'}^{\beta'} 
\Theta_{j'}'(\zeta',t')
-
\beta'!\,
\Theta_{j',
\beta'}'(T')
\right]_{t'=\tau'=0'}
\\
& \
=
\beta'!\, 
\Big(
\Theta_{j',\beta'}'(0')
-
\Theta_{j',\beta'}'(T')
\Big).
\endaligned
\]
Consequently, up to a translation by a constant, we recover with
$\Psi_{j', \beta'}' (T')$ the components of the infinite Segre
mapping $\mathcal{ Q }_\infty'$ of $M'$. Hence the
next definition generalizes the concepts introduced before.

\def\thedefinition{3.15}\begin{definition}
{\rm 
The formal CR mapping $h : 
(M, 0) \to_{ \mathcal{ F}} (M', 0')$ is called

\begin{itemize}

\smallskip\item[{\bf (h1)}]
{\sl Levi-nondegenerate} at the origin if the mapping
\[
t'\mapsto \left(
R_{j',\beta}'(0,0,J_\tau^{\vert \beta \vert} \bar h(0) :\,
t')
\right)_{1\leqslant j'\leqslant d', \, \vert
\beta \vert \leqslant 1}
\]
is of rank $n'$ at $t' =0'$;

\smallskip\item[{\bf (h2)}]
{\sl finitely nondegenerate} at the origin if there exists an integer
$\ell_1$ such that the mapping
\[
t'\mapsto \left(
R_{j',\beta}'(0,0,J_\tau^{\vert 
\beta \vert} \bar h(0) :\, t')
\right)_{1\leqslant j'\leqslant d', \, \vert
\beta \vert \leqslant k}
\]
is of rank $n'$ at $t' = 0'$, for $k= \ell_1$, hence for every $k
\geqslant \ell_1$;

\smallskip\item[{\bf (h3)}]
{\sl essentially finite} at the origin if there exists an integer
$\ell_1$ such that the mapping
\[
t'\mapsto \left(
R_{j',\beta}'(0,0,J_\tau^{\vert \beta \vert} \bar h(0) :\,
t')
\right)_{1\leqslant j'\leqslant d', \, \vert
\beta \vert \leqslant k}
\]
is locally finite at $t' = 0'$, for $k= \ell_1$, hence for every $k
\geqslant \ell_1$;

\smallskip\item[{\bf (h4)}]
{\sl Segre nondegenerate} at the origin if there exist an integer
$\ell_1$, integers ${j_*' }^1, \dots, {j_*' }^{ n'}$ with $1 \leqslant
{j_*' }^{i' }\leqslant d'$ for $i'=1, \dots, n'$ and multiindices
$\beta_*^1, \dots, \beta_*^{n'}$ with $\vert \beta_*^{i' } \vert
\leqslant \ell_1$ for $i'= 1, \dots, n'$, such that the determinant
\[
{\rm det}\, \left(
\frac{\partial R_{{j_*'}^{i_1'},\beta_*^{i_1'}}'}{
\partial t_{i_2'}'}\left(
z,\overline{\Theta}(z,0),0,0,
J^{\vert \beta_*^{i_1'} \vert} \bar h(0) :\,
h(z,\overline{\Theta}(z,0))
\right)
\right)_{
1\leqslant i_1',i_2'\leqslant n'}
\]
does not vanish identically in $\C\dl z \dr$;

\smallskip\item[{\bf (h5)}]
{\sl holomorphically nondegenerate} at the origin if there exists an
integer $\ell_1$, integers ${j_*' }^1, \dots, {j_*'}^{ n'}$ with
$1\leqslant {j_*' }^{ i' }\leqslant d'$ for $i'=1, \dots, n'$ and multiindices
$\beta_*^1, \dots, \beta_*^{n'}$ with $\vert \beta_*^{i'} \vert \leqslant
\ell_1$ for $i' =1, \dots, n'$, such that the determinant
\[
{\rm det}\, \left(
\frac{\partial R_{{j_*'}^{i_1'},\beta_*^{i_1'}}'}{
\partial t_{i_2'}'}\left(
0,0,0,0,
J^{\vert\beta_*^{i_1'}\vert} 
\bar h(0) :\, 
h(t)
\right)
\right)_{
1\leqslant i_1',i_2'\leqslant n'}
\]
does not vanish identically in $\C\dl t \dr$.

\end{itemize}\smallskip

}\end{definition}

The nondegeneracy of the formal mapping $h$ requires
the same nondegeneracy on the target $(M', 0')$.

\def\thelemma{3.16}\begin{lemma}
{\rm (\cite{ me2004a})} Let $h: (M,0) \to_{\mathcal{ F}} (M', 0')$ be
a formal CR mapping.

\begin{itemize}

\smallskip\item[{\bf (1)}]
If $h$ is Levi-nondegenerate at $0$, then
$M'$ is necessarily Levi-nondegenerate at $0'$. 

\smallskip\item[{\bf (2)}]
If $h$ is finitely nondegenerate at $0$, then 
$M'$ is necessarily finitely nondegenerate at $0'$.

\smallskip\item[{\bf (3)}]
If $h$ is essentially finite at $0$, then 
$M'$ is necessarily essentially finite at $0'$.

\smallskip\item[{\bf (4)}]
If $h$ is Segre nondegenerate at $0$, then 
$M'$ is necessarily Segre nondegenerate at $0'$.

\smallskip\item[{\bf (5)}]
If $h$ is holomorphically nondegenerate at $0$, then 
$M'$ is necessarily holomorphically nondegenerate at $0'$.

\end{itemize}\smallskip

\end{lemma}

We now show that CR-transversality of the mapping $h$ insures that it
enjoys exactly the same nondegeneracy condition as the target
$(M', 0')$. 

\def\thetheorem{3.17}\begin{theorem}
{\rm (\cite{ me2004a})} Assume that the
formal CR mapping $h : (M, 0) \to_{ \mathcal{ F}}
(M', 0')$ is CR-transversal at $0$. Then the following five
implications hold:

\begin{itemize}

\smallskip\item[{\bf (1)}]
If $M'$ is Levi nondegenerate at $0'$, then $h$ is finitely
nondegenerate at $0$.

\smallskip\item[{\bf (2)}]
If $M'$ is finitely nondegenerate at $0'$, then $h$ is finitely
nondegenerate at $0$.

\smallskip\item[{\bf (3)}]
If $M'$ is essentially finite at $0'$, then $h$ is 
essentially finite at $0$.

\smallskip\item[{\bf (4)}]
If $M'$ is Segre nondegenerate at $0'$, then $h$ is Segre
nondegenerate at $0$.

\smallskip\item[{\bf (5)}]
If $M'$ is holomorphically nondegenerate, and if moreover $h$
is transversal at $0$, then $h$ is holomorphically nondegenerate at
$0$.

\end{itemize}\smallskip

\end{theorem}

The above five implications also hold under the assumption that $h$ is
either CR-invertible, or CR-submersive, or CR-finite with $m=m'$ or
CR-dominating: this provides at least 20 more
(less refined) versions of the theorem, some of which 
appear in the literature.

\smallskip

Other relations hold true between the nondegeneracy conditions on $h$
and on the generic submanifolds $(M, 0)$ and $(M', 0')$. 
We mention some, concisely. As above,
assume that $h : (M, 0) \mapsto_{ \mathcal{ F} } (M', 0)$ is a formal
CR mapping. Since $dh_0 ( T_0^c M) \subset T_0^c M'$, a linear map
$dh_0^{\rm trv} : T_0 M / T_0^cM \to T_0 M' / T_0^c M'$ is
induced. Assume $d' = d$ and $m' = m$. The next statement may be
interpreted as a kind of Hopf Lemma for CR mappings.

\def\thetheorem{3.18}\begin{theorem}
{\rm (\cite{ br1990, er2004})} If $M$ is minimal at $0$ and if
$h$ is CR-dominating at $0$, then $d h_0^{ \rm trv} : T_0 M / T_0^cM \to
T_0 M' / T_0^c M'$ is an isomorphim.
\end{theorem}

An open question is to determine whether the condition that the
jacobian determinant ${\rm det} \, \big( \frac{\partial h_i}{
\partial t_j } (t) \big)_{ 1\leqslant i, j \leqslant n}$ does not vanish
identically in $\C \dl t\dr$ is sufficient to insure that $d h_0^{ \rm
trv} : T_0 M / T_0^c M \to T_0 M' / T_0^c M'$ is an isomorphism. A
deeper understanding of the constraints between various nondegeneracy
conditions on $h$, $M$ and $M'$ would be desirable.

\subsection*{ 3.19.~Classical versions of the reflection principle}
Let $h: (M, 0) \to_{ \mathcal{ F}} (M', 0)$ be a formal power series
CR mapping between two generic submanifolds. Assume that
$M$ is minimal at $0$.

\def\thetheorem{3.20}\begin{theorem}
{\rm (\cite{ ber1999, me2004a, me2005})} 
If $M$ and $M'$ are real analytic, if $h$ is either Levi
nondegenerate, or finitely nondegenerate, or essentially finite, or
Segre nondegenerate at the origin, then $h (t)$ is convergent, namely
$h(t ) \in \C \{ t\}^{ n'}$. If moreover, $M$ and $M'$ are algebraic,
then $h$ is algebraic.
\end{theorem}

If one puts separate nondegeneracy conditions on $h$ and on $M'$, as
in Theorem~3.17, one obtains a combinatorics of possible statements,
some of which appear in the literature.

\smallskip

If $h$ is finitely nondegenerate (level {\bf (2)}), the (paradigmatic)
proof yields more information.

\def\thetheorem{3.21}\begin{theorem} 
{\rm (\cite{ ber1999, me2005})}
As above, let $h: ( M, 0) \to ( M', 0')$ be a formal power series CR
mapping. Assume that $M$ is minimal at $0$ and let $\nu_0$ be the
integer of Corollary~2.17. Assume also that $h$ is $\ell_1$-finitely
nondegenerate at $0$. Then there exists a $\C^{ n'}$-valued power
series mapping $H(t, J^{ 2\nu_0 \ell_1})$ which is constructed
algorithmically by means of the defining equations of $(M, 0)$ and of
$( M', 0')$, such that the power series identity
\[
h(t)
\equiv 
H(t,J^{2\nu_0\ell_1}h(0))
\]
holds in $\C\dl t \dr^{n'}$. If $M$ and $M'$ are real analytic {\rm
(}resp. algebraic{\rm )}, $H$ is holomorphic {\rm (}resp. complex
algebraic{\rm )} in a neighborhood of $0\times J^{2\nu_0\ell_1} h
(0)$.
\end{theorem}

In~\cite{ ber1999, gm2004}, 
the above formula $h(t) \equiv H(t, J^{2 \nu_0 \ell_1}
h(0))$ is studied horoughly in the case where $M' = M$ and
$h$ is a local holomorphic automorphism of $(M, 0)$ close to the
identity.

\smallskip

At level {\bf (5)}, namely with a holomorphically nondegenerate target
$(M', 0')$, the reflection principle is much more delicate.
It requires the introduction of a new object, whose regularity
properties hold in fact {\it without any nondegeneracy 
assumption on the target $(M', 0')$}.

\subsection*{ 3.22.~Convergence of the reflection mapping}
The {\sl reflection mapping} associated to $h$ and to the system of
coordinates $(z', w')$ is~:
\[
\mathcal{ R}_h'(\tau',t
):= 
\xi' 
-
\Theta'(\zeta',h(t))\in 
\C\dl\tau',t\dr^{d'}.
\]
Since $h$ is formal, it is only a formal power series mapping. As
argued in the introduction of~\cite{ me2005}, it is the most
fundamental object in the analytic reflection principle. In the case
of CR mappings between essentially finite hypersurfaces, the analytic
regularity of the reflection mapping is equivalent to the extension of
CR mappings as correspondences, as studied in \cite{ dp1995,
sh2000, sh2003, dp2003}. Without nondegeneracy
assumption on $(M', 0')$, the reflection mapping enjoys regularity
properties from which all analytic reflection principles may be
deduced. Here is the very main theorem of this
Section~3.

\def\thetheorem{3.23}\begin{theorem} 
{\rm (\cite{ me2001b, bmr2002, me2005})} If $M$ is
minimal at the origin and if $h$ is either 
CR-invertible, or CR-submersive, or CR-finite, or CR-dominating, or
CR-transversal, then for every system of coordinates $(z', w') \in
\C^{ m'} \times \C^{ d'}$ in which the extrinsic complexification
$\mathcal{ M}'$ is represented by $\xi ' = \Theta '( \zeta', t')$, the
associated CR-reflection mapping is convergent, namely $\mathcal{
R}_h' (\tau', t) \in \C\{ \tau', t\}^{ d'}$.
\end{theorem}

If the convergence property holds in one such system of coordinates,
it holds in all systems of coordinates (\cite{ me2005};
Proposition~3.26 below). Further, if we develope $\Theta ' (\zeta',
t') = \sum_{ \gamma ' \in \N^{ m'}} \, (\zeta' )^{ \gamma'} \,
\Theta_{ \gamma ' }' (t')$, the convergence of $\mathcal{ R}_h ' (
\tau', t)$ has a concrete signification.

\def\thecorollary{3.24}\begin{corollary}
All the components $\Theta_{ \gamma '}' (h(t))$ of the reflection
mapping are convergent, namely $\Theta_{ \gamma '}' (h(t)) \in \C \{
t\}^{ d'}$ for every $\gamma' \in \N^{ m'}$.
\end{corollary}

Conversely (\cite{ me2001b, me2005}), if $\Theta_{ \gamma '}'
(h(t)) \in \C \{ t\}^{ d'}$ for every $\gamma' \in \N^{ m'}$, an
elementary application of the Artin approximation Theorem~3.28
(below) yields
Cauchy estimates: there exist $\rho >0$, $\sigma >0$ and $C>0$ so that
$\vert \Theta_{ \gamma'}' (h(t)) \vert < C \, (\rho)^{ - \vert \gamma '
\vert}$, for every $t \in \C^n$ with $\vert t\vert < \sigma$.
It follows that $\mathcal{ R}_h' (\tau', t)
\in \C \{ \tau', t\}^{ d'}$.

\smallskip

Taking account of the nondegeneracy conditions {\bf (ndi)} and {\bf
(crj)}, several corollaries may be deduced from the theorem. Most of
them are already expressed by Theorem~3.20, {\it except notably}\, the
delicate case where $(M',0')$ is holomorphically nondegenerate.

\def\thecorollary{3.25}\begin{corollary}
{\rm (\cite{ me2001b, me2005})} If $M$ is minimal at the
origin, if $(M', 0')$ is holomorphically nondegenerate and if $h$ is
either CR-invertible and invertible, or CR-submersive and submersive,
or CR-finite and finite with $m' = m$, or CR-dominating and
dominating, or CR-transversal and transversal, then $h (t) \in \C \{
t\}^{ n'}$ is convergent.
\end{corollary}

It is known (\cite{ st1996}) that $(M',0')$ is holomorphically
degenerate if and only if there exists a nonzero $(1,0)$
vector field $X' = \sum_{ i' = 1}^{ n'} \, a_{i'} ' (t') \frac{
\partial }{ \partial t_{ i' }'}$ having holomorphic coefficients which
is tangent to $(M',0')$. In the corollary above, holomorphic
nondegeneracy is optimal for the convergence of a formal equivalence:
if $M'$ is holomorphically degenerate, if $(s', t')\longmapsto \exp
(s' X')(t')$ denotes the local flow of $X'$, where $s'\in \C$, $t'\in
\C^{n'}$, there indeed exist (\cite{ ber1999, me2005})
nonconvergent power series $\varpi' ( t') \in \C\dl t' \dr$ such that
$t' \mapsto_{ \mathcal{ F }} \exp (\varpi' (t') X') (t')$ is a
nonconvergent formal equivalence of $M'$.

\smallskip

The invariance of the reflection mapping is crucial.

\def\theproposition{3.26}\begin{proposition}
{\rm (\cite{ me2002, me2004a, me2005})} The
convergence of the reflection mapping is a biholomorphically invariant
property. More precisely, if $t'' = \phi ' (t')$ is a local
biholomorphism fixing $0'$ and transforming $(M', 0')$ into a generic
submanifold $(M'', 0')$ of equations $\bar w_{ j'} '' = \Theta_{ j'} ''
(\bar z'', t'')$, $j' = 1, \dots, d'$, the composed reflection
mapping of $\phi' \circ h : (M, 0) \to_{ \mathcal{ F}} (M'', 0')$
defined by
\[
\aligned
\mathcal{ R}_{\phi'\circ h}''(\tau'',t)
:= 
&\
\xi'' 
-
\Theta''
\big(
\zeta'',\phi'(h(t))
\big)
\\
=
&\
\xi''
-
\sum_{\gamma'\in\N^{m'}}\,
(\zeta'')^{\gamma'}\,
\Theta_{\gamma'}''
\big(
\phi'(h(t))
\big)
\endaligned
\]
has components $\Theta_{ \gamma' }'' \big( \phi' (h(t )) \big)$ given
by formulas
\[
\Theta_{\gamma'} '' 
\big(\phi'(h(t))\big) 
\equiv
S_{\gamma'}'\left(
\big(
\Theta_{\gamma_1'}'(h(t))
\big)_{ 
\gamma_1'\in\N^{m'}}
\right),
\] 
where the local holomorphic functions $S_{ \gamma' }'$ depend only on
the biholomorphism $t'' = \phi' (t')$ {\rm (}they have an infinite
number of variables, but the necessary Cauchy estimates insuring
convergence are automatically satisfied{\rm )}.
\end{proposition}

A few words about the proof of the main Theorem~3.23. Although the
classical reflection principle deals only with the ``reflection
identities''~\thetag{ 3.13}, to get the most adequate version of the
reflection principle, it is unavoidable to understand the
symmetry between the variables $t$ and the 
variables $\tau = (\bar t)^c$.

The assumption that $h^c$ maps formally $( \mathcal{ M}, 0)$ to
$(\mathcal{ M}', 0')$ is equivalent to each one of the following two
formal identities:
\[
\left\{
\aligned
\overline{ g} (\tau) 
& \
=
\sum_{ \gamma' \in \N^{m'}}\, 
\overline{ f} (\tau)^{\gamma '} \,
\Theta_{\gamma' } ' (h(t)), \\
g (t) 
& \
=
\sum_{\gamma' \in\N^{ m'}} \, 
f(t)^{\gamma '} \, 
\overline{ \Theta}_{ 
\gamma'} ' \left(\overline{ h} (\tau) \right),
\endaligned\right.
\]
on $\mathcal{ M}$, namely after replacing either $w$ by $\overline{
\Theta} (z, \tau)$ or $\xi$ by $\Theta (\zeta, t)$. The symmetry may
be pursued by considering the two families of derivations:
\[
\left\{
\aligned
\mathcal{ L}^\beta :=
& \ 
( \mathcal{ L }_1 )^{ \beta_1}(
\mathcal{ L }_1 )^{ \beta_2} \cdots ( \mathcal{L }_m)^{ \beta_m}
\ \ \ \ \
{\rm and} \\
\underline{ \mathcal{ L}}^\beta :=
& \ 
(\underline{ \mathcal{ L }}_1 )^{ \beta_1}
(\underline{ \mathcal{ L
}}_1 )^{ \beta_2} \cdots ( 
\underline{ \mathcal{L }}_m)^{ \beta_m},
\endaligned\right.
\]
where $\beta = (\beta_1, \beta_2, \dots, \beta_m) \in \N^m$. Applying
them to the two formal identities above,
if we respect the completeness of the combinatorics, we will get
{\it four}\, families of {\sl reflection identities}. The first pair
is obtained by applying $\underline{ \mathcal{ L}}^\beta$ to the
two formal identities above:
\[
\left\{
\aligned
\underline{ \mathcal{ L}}^\beta \, 
\overline{ g} (\tau) 
& \
=
\sum_{ \gamma' \in \N^{m'}}\,
\underline{ \mathcal{ L}}^\beta \big[ 
\overline{ f} (\tau)^{\gamma '}\big] \,
\Theta_{\gamma' } ' (h(t)), \\
0
& \
=
\sum_{\gamma '\in \N^{ m'}} \, 
f(t)^{\gamma '} \, 
\underline{ \mathcal{ L}}^\beta \big[
\overline{ \Theta}_{ 
\gamma'} ' \left(\overline{ h} (\tau) \right)\big].
\endaligned\right.
\]
The second pair is obtained by applying $\mathcal{
L}^\beta$, permuting the two lines:
\[
\left\{
\aligned
\mathcal{ L}^\beta 
g (t)
& \
=
\sum_{\gamma' \in \N^{ m'}} \,
\mathcal{ L}^\beta \big[ 
f(t)^{\gamma '}\big] \, 
\overline{ \Theta}_{ 
\gamma'} ' \left(\overline{ h} (\tau) \right), 
\\
0
& \
=
\sum_{ \gamma' \in \N^{m'}}\, 
\overline{ f} (\tau)^{\gamma '} \,
\mathcal{ L}^\beta \left[
\Theta_{\gamma' } ' (h(t))\right]. \\
\endaligned\right.
\]
We immediately see that these two pairs are conjugate line by line.
In each pair, we notice a crucial difference between the first and the
second line: whereas it is $\overline{ g}$ and the power $\overline{
f}^{ \gamma'}$ (or $g$ and $f^{ \gamma'}$) that are differentiated in
each first line, in each second line, only the components
$\overline{ \Theta}_{ \gamma'} ' ( \overline{ h} )$ (or $\Theta_{
\gamma'}' (h)$) of the reflection mapping, which are the right
invariant functions, are differentiated. In a certain sense, it is
forbidden to differentiate $\overline{ g}$ and $\overline{ f}^{
\gamma'}$ (or $g$ and $f^{ \gamma'}$), because the components $(f, g)$
of $h$ need not enjoy a reflection principle. In fact, in the proof
of the main Theorem~3.23, one has to play constantly with the four
reflection identities above.

Since we cannot summarize here the long and refined proof, we only
formulate the main technical proposition. Denote by $J_t^\ell \psi$
the $\ell$-th jet of a power series $\psi (t) \in \C\dl t\dr^{ d'}$,
for instance $J_t^\ell 
\Theta_{ \gamma'} '( h)$ for some $\gamma ' \in\N^{
m'}$. Remind that $\underline{ \Gamma }_k$ and $\Gamma_k$ are
(conjugate) Segre chains.
Let $N_{ d', n, \ell} := 
d'\, \frac{ (n + \ell) ! }{ n! \ \ell !}$.

\def\theproposition{3.27}\begin{proposition}
{\rm (\cite{ me2005})}
For every $k\in \N$ and every $\ell \in \N$,
the following two properties hold{\rm :}
\begin{itemize}
\smallskip\item[$\bullet$]
if $k$ is odd, for every $\gamma' \in \N^{ m'}${\rm :}
\[
\left[J_t^\ell\Theta_{\gamma'}'(h)\right]
\left(
\Gamma_k\left(
z_{ (k)}\right) 
\right)\in\C\{z_{(k)}\}^{N_{d',n,\ell}}~;
\]
\smallskip\item[$\bullet$]
if $k$ is even, for every $\gamma' \in \N^{ m'}${\rm :} 
\[
\left[J_\tau^\ell \overline{\Theta}_{ 
\gamma'}'(\overline{ h})\right]
\left(
\Gamma_k\left(
z_{ (k)}\right) 
\right)\in\C\{z_{(k)}\}^{N_{d',n,\ell}}.
\]
\end{itemize}
\end{proposition}

With $\ell = 0$ and $k = 2\nu_0$, thanks to Corollary~2.17, we deduce
from this main proposition that $\Theta_{ \gamma'} ' (h(t)) \in \C\{
t\}^{ \gamma'}$ for every $\gamma' \in \N^{ m'}$.
This yields Theorem~3.23.

\smallskip

The main tool in the proof of this proposition is an approximation 
theorem saying that a formal power series mapping
that is a solution of some analytic equations may be corrected so as
to become convergent and still a solution.

\def\thetheorem{3.28}\begin{theorem}
{\rm ({\sc Artin} \cite{ ar1968, jopf2000})} Let $\K = \R$ or
$\C$, let $n \in \N$ with $n \geqslant 1$, let
${\sf x} = ({\sf x}_1, \dots,\, {\sf x}_n) \in \K^n$, let $m \in \N$,
with $m \geqslant 1$, let ${\sf y} =({\sf y}_1, \dots ,\, {\sf y }_m ) \in
\K^n$, let $d \in \N$ with $d \geqslant 1$ and let $R_1 ( {\sf x}, {\sf y
}), \dots, \, R_d({\sf x},\, {\sf y})$ be an arbitrary collection of
formal power series in $\K \{{\sf x},\, {\sf y} \}$ that vanish at the
origin, namely $R_j(0,\, 0) =0$, $j=1, \dots,\, d$. Assume that there
exists a formal mapping $h( {\sf x}) = (h_1( {\sf x }), \dots,\, h_m(
{\sf x})) \in \K\dl {\sf x} \dr^m$ with $h(0) =0$ such that
\[
R_j\left({\sf x} ,\, h({\sf x})\right)\equiv 0 \ \ 
{\rm in} \ \ \K\dl {\sf x} \dr, \ \ \ \
{\rm for} \ \
j=1,\dots,\, d.
\]
Let $\mathfrak{m}( {\sf x}) := {\sf x}_1 \K\dl {\sf x} \dr + \cdots+
{\sf x}_n \K \dl {\sf x} \dr$ be the maximal ideal of $\K \dl {\sf x}
\dr$. For every integer $N \geqslant 1$, there exists a convergent power
series mapping $h^N( {\sf x}) \in \K\{ {\sf x} \}^m$ such that
\[
R_j\left({\sf x},\, h^N({\sf x})\right)\equiv 0 \ \ 
{\rm in} \ \ \K \dl {\sf x} \dr, \ \ \ \
{\rm for} \ \ 
j=1,\dots,\, d,
\]
that approximates $h({\sf x})$ to order $N-1${\rm :}
\[
h^N({\sf x}) \equiv h({\sf x}) \ \
{\rm mod} \, \left(\mathfrak{m}(x)^N \right).
\]
\end{theorem}

As an application of the main Theorem~3.23, an approximation property
for formal CR mappings holds.

\def\thetheorem{3.29}\begin{theorem}
{\rm (\cite{ me2005})} 
Under the assumptions of Theorem~3.23, for every integer $N \geqslant 1$,
there exists a convergent power series mapping ${\sf H}^N (t) \in \C\{
t\}^{ n'}$ with ${\sf H}^N (t) \equiv h(t) \ {\rm mod} \, (\mathfrak{
m} (t))^N$ {\rm (}whence $H(0) = 0${\rm )}, that induces a local
holomorphic map from $(M, 0)$ to $(M', 0')$.
\end{theorem}

\def\thecorollary{3.30}\begin{corollary}
{\rm (\cite{ me2001b, me2005})} Assume that $n' = n$, that $d'
= d$, that $M$ is minimal at the origin, and that $h : (M, 0) \to_{
\mathcal{ F}} (M', 0')$ is a formal {\rm (}invertible{\rm )}
equivalence. Then $M$ and $M'$ are biholomorphically equivalent.
\end{corollary}

It is known (\cite{ st1996, ber1999, gm2004}) that a minimal
holomorphically nondegenerate real analytic generic submanifold of
$\C^n$ has finite-dimensional local holomorphic automorphism
group. Unique determination by a jet of finite order follows from a
representation formula, as in Theorem~3.21. More
generally:

\def\thecorollary{3.31}\begin{corollary}
{\rm (\cite{ me2001b, bmr2002, me2005})} Assume that
$m' = m$ and $d' = d$, that $(M, 0)$ is minimal at the origin and that
$(M', 0)$ is holomorphically nondegenerate. There exists an integer
$\kappa = \kappa (m, d)$ such that, if two local biholomorphisms $h^1,
h^2 : (M, 0) \to (M', 0)$ have the same $\kappa$-th jet at the origin,
then $h^1 = h^2$.
\end{corollary}

From an inspection of the proof, Theorem~3.29 holds without the
assumption that $(M, 0)$ is minimal, but with the assumption that its
CR orbits have constant dimension in a neighborhood of $0$.
However, the case where CR orbits have arbitrary dimension
is delicate.

\def\theopenquestion{3.32}\begin{openquestion}
Does formal equivalence coincide with biholomorphic equivalence in the
category of real analytic generic local submanifolds of $\C^n$ whose
CR orbits have non-constant dimension\,?
\end{openquestion}

\subsection*{ 3.33.~Algebraicity of the reflection mapping}
We will assume that both $M$ and $M'$ are algebraic. Remind that
Theorem~3.20 shows the algebraicity of $h$ under some hypotheses. A
much finer result is as follows. It synthetizes all existing results
(\cite{ we1977, ss1996, cms1999, ber1999, za1999}) about algebraicity
of local holomorphic mappings.

\def\thetheorem{3.34}\begin{theorem}
{\rm (\cite{ me2001a})} If $h$ is a local holomorphic map $(M, 0) \to
(M', 0')$, if $M$ and $M'$ are algebraic, if $M$ is minimal at the
origin and if $M'$ is the smallest {\rm (}for inclusion{\rm )} local
real algebraic manifold containing $h(M)$, then the reflection mapping
$\mathcal{ R}_h' (\tau', t)$ is algebraic.
\end{theorem}

Trivial examples (\cite{ me2001a}) show that the algebraicity of
$\mathcal{ R}_h'$ need not hold if $M'$ is not the smallest one.

In fact, Theorem~3.34 also holds (with the same proof) if one assumes
only that the source $M$ is minimal at a Zariski-generic point: it
suffices to shrink $M$ and the domain of definition of $h$ around such
points, getting local algebraicity of $\mathcal{ R }_h'$ there, and
since algebraicity is a global property, $\mathcal{ R }_h'$ is
algebraic everywhere.

An equivalent formulation of Theorem~3.34 uses the
concept of transcendence degree, studied in~\cite{ pu1990, cms1999,
me2001a}. With $n_{ M' }'$ being the essential holomorphic
dimension of $(M', 0')$ defined in \S3.6, set $\kappa_{ M'}' := n' -
n_{ M'}'$. Observe that $(M', 0')$ is holomorphically nondegenerate
precisely when $\kappa_{ M'}' = 0$. Denote by $\C [t ]$ the ring of
complex polynomials of the variable $t \in \C^n$ and by $\C(t)$ its
quotient field. Let $t' = h (t )$ be a local holomorphic mapping as in
Theorem~3.34. and let $\C(t) (h_1 (t), \ldots, h_{n'} (t))$ be the
field generated by the components of $h$.

\def\thetheorem{3.35}\begin{theorem}
{\rm (\cite{ me2001a})} With the same assumptions as in Theorem~3.34,
the transcendence degree of the field extension $\C(t) \to \C(t)
(h(t))$ is less than or equal to $\kappa_{ M'}'$.
\end{theorem}

\def\thecorollary{3.36}\begin{corollary}
{\rm (\cite{ cms1999, za1999, me2001a})}
If $M$ is minimal at a Zariski-generic point and if the real algebraic
target $M'$ does not contain any complex algebraic curve, then the
local holomorphic mapping $h$ is algebraic.
\end{corollary}

However, in case $h$ is only a formal CR mapping, it is impossible to
shift the central point to a nearby minimal point. Putting the
simplest rank assumption (invertibility) on $h$, we may thus
formulate delicate problems for the future.

\def\theopenquestion{3.37}\begin{openquestion}
Let $h$ be a formal equivalence between two real analytic generic
submanifolds of $\C^n$ which are minimal at a Zariski-generic point.

\begin{itemize}

\smallskip\item[$\bullet$] 
Is the reflection mapping convergent\,?

\smallskip\item[$\bullet$]
Is $h$ uniquely determined by a jet of finite order when the target
is holomorphically nondegenerate\,?

\smallskip\item[$\bullet$]
Is $h$ convergent under the assumption that the real analytic target
$M'$ does not contain any complex analytic curve\,?

\end{itemize}\smallskip

\end{openquestion}

For $M'$ algebraic containing no complex algebraic curve and $M$
minimal at $0$, the third question has been settled in~\cite{
mmz2003b}. However, the assumption of algebraicity of $M'$ is strongly
used there, because these authors deal with the transcendence degree
of the field extension $\C (t) \to \C (t) (h(t))$, a concept which is
meaningless if $M'$ is real analytic. For further (secondary)
results and open questions, we refer to~\cite{ bmr2002, ro2003}. This
closes up our survey of the formal/algebraic/analytic reflection
principle.

A generic submanifold $M \subset \C^n$ is called {\sl locally
algebraizable} at one of its points $p$ if there exist local
holomorphic coordinates centered at $p$ in which it is Nash algebraic.
Unlike partial results, the following question remains up to now
unsolved.

\def\theopenproblem{3.38}\begin{openproblem}
{\rm (\cite{ hu2001, hjy2001, ji2002, gm2004, fo2004})} Formulate a
{\rm necessary and sufficient} condition for the local
algebraizability of a real analytic hypersurface $M \subset \C^n$ in
terms of a basis of the {\rm (}differential{\rm )} algebra of its
Cartan-Hachtroudi-Chern invariants.
\end{openproblem}

To conclude, we would like to mention that the complete theory of CR
mappings may be transferred to systems of partial differential
equations having finite-dimensional Lie symmetry group.
This aspect will be treated in subsequent publications (\cite{
me2006a, me2006b}).

\newpage

\begin{center}
{\Large\bf III:~Systems of vector fields and CR functions}
\end{center}

\bigskip\bigskip\bigskip

\begin{center}
\begin{minipage}[t]{11cm}
\baselineskip =0.35cm
{\scriptsize

\centerline{\bf Table of contents}

\smallskip

{\bf 1.~Sussman's orbit theorem and structural properties of orbits 
\dotfill 50.}

{\bf 2.~Finite type systems and their genericity (openess and
density) \dotfill 63.}

{\bf 3.~Locally integrable CR structures \dotfill 74.}

{\bf 4.~Smooth generic submanifolds and their CR orbits \dotfill 87.}

{\bf 5.~Approximation and uniqueness principles \dotfill 106.}

\smallskip

\hfill 
{\footnotesize\tt [7 diagrams]}

}\end{minipage}
\end{center}

\bigskip
\bigskip

{\small


Beyond the theorems of Frobenius and of Nagano, Sussmann's theorem
provides a means, valid in the smooth category, to construct all the
integral manifolds of an arbitrary system of vector fields, as orbits
of the pseudo-group actions of global flows. The fundamental
properties of such orbits: lower semi-continuity of dimension, local
flow box structure, propagation of embeddedness, intersection with a
transversal curve in the one-codimensional 
case, are essentially analogous, but
different from the ones known in foliation theory. Orbits possess
wide applications in Control Theory, in sub-Riemannian Geometry, in the
Analysis of Linear Partial Differential Equations and in
Cauchy-Riemann geometry.

Let $L_j f = g_j$, $j=1, \dots, \lambda$, be a linear PDE system with
unknown $f$, where $g$ is smooth and where $\{ L_k \}_{ 1\leqslant k
\leqslant r}$ is an involutive (in the sense of Frobenius) system of
smooth vector fields on $\R^n$ having {\it complex-valued}\,
coefficients. Since Lewy's celebrated discovery of an example of a
single equation $L f = g$ in $\R^3$ without any solution, a major
problem in the Analysis of PDE's is to find adequate criterions for
the existence of local solutions. Condition (P) of Nirenberg-Treves
has appeared to be necessary and sufficient to insure local
integrability of a single equation of principal type having simple
characteristics. The problem of characterizing systems of several
linear first order PDE's having maximal space of solutions is not yet
solved in full generality; several fine questions remain open.

Following Treves, to abstract the notion of systems involving several
equations, an {\sl involutive structure} on a smooth $\mu$-dimensional
real manifold $M$ is a $\lambda$-dimensional complex subbundle
$\mathcal{ L}$ of $\C \otimes TM$ satisfying $\left[\mathcal{ L},
\mathcal{ L} \right] \subset \mathcal{ L}$. The automatic
integrability of smooth almost complex structures (those with
$\mathcal{ L} \oplus \overline{ \mathcal{ L}} = \C \otimes TM$) and
the classical (non)integrability theorems for smooth abstract CR
structures (those with $\mathcal{ L} \cap \overline{ \mathcal{ L}} =
\{ 0 \}$) are inserted in this general framework.

Beyond such problematics, it is of interest to study the analysis and
the geometry of subbundles $\mathcal{ L}$ whose space of solutions is
maximal, viz the preceding question is assumed to be solved,
optimally: in a neighborhood of every point of $M$, there exist $(\mu
- \lambda)$ local complex valued functions $z_1, \dots, z_{\mu -
\lambda}$ having linearly independent differentials which are
solutions of $\mathcal{ L} z_k = 0$. Such involutive structures are
called {\sl locally integrable}. Some representative examples are
provided by the bundle of anti-holomorphic vector fields tangent to
various embedded generic submanifolds of $\C^n$. According to a
theorem due to Baouendi-Treves, every local solution of $\mathcal{ L}
f = 0$ may be approximated sharply by polynomials in a set of
fundamental solutions $z_1, \dots, z_{ \mu - \lambda}$, in the
topology of functional spaces as $\mathcal{ C}^{\kappa, \alpha}$, $L_{
loc}^{\sf p}$, or $\mathcal{ D}'$.

In a locally integrable structure, the Sussmann orbits of the vector
fields ${\rm Re}\, L_k, {\rm Im}\, L_k$ are then of central importance
in analytic and in geometrical questions. They show up propagational
aspects, as for instance: the support of a function or distribution
solution $f$ of $\mathcal{ L} f = 0$ is a union of orbits. The
approximation theorem also yields an elegant proof of uniqueness in
the Cauchy problem. Further propagational aspects will be studied in
the next chapters, using the method of analytic discs.
Sections~3, 4 and~5 of this chapter and the remainder of the memoir
are focused on embedded generic submanifolds.

} 

\section*{ \S1.~Sussmann's theorem and structural properties
\\
of orbits}

\subsection*{ 1.1.~Integral manifolds of a system of vector fields}
Ordinary differential equations in the modern sense emerged in the
seventieth century, concomitantly with the infinitesimal calculus.
Nowadays, in contemporary mathematics, the abstract study of vector
fields is inserted in several broad areas of research, among
which we perceive the following. 

\smallskip

\begin{itemize}

\item[$\bullet$]
{\sf Control Theory:} 
controllability of vector fields on $\mathcal{ C}^\infty$ and real
analytic manifolds; nonholonomic systems; sub-Riemannian geometry
(\cite{ gv1987, bel1996}).

\smallskip

\item[$\bullet$]
{\sf Dynamical systems:} singularities of real or complex vector
fields and foliations; normal forms and classification; phase
diagrams; Lyapunov theory; Poincar\'e-Bendixson theory; theory of
limit cycles of polynomial and analytic vector fields; small divisors
(\cite{ ar1978, ar1988}).

\smallskip

\item[$\bullet$]
{\sf Lie-Cartan theory:} infinitesimal symmetries of differential
equations; classification of local Lie group actions; Lie algebras of
vector fields; representations of Lie algebras; exterior differential
systems; Cartan-Vessiot-K\"ahler theorem; Janet-Riquier theory;
Cartan's method of equivalence (\cite{ ol1995, stk2000}).

\smallskip

\item[$\bullet$]
{\sf Numerical analysis:} systems of (non)linear ordinary differential
equations; methods of: Euler, Newton-Cotes, Newton-Raphson,
Runge-Kutta, Adams-Bashforth, Adams-Moulton (\cite{ de1996}).

\smallskip

\item[$\bullet$]
{\sf PDE theory:}
Local solvability of linear partial differential equations; uniqueness
in the Cauchy problem; propagation of singularities; {\sc fbi}
transform and control of wave front set
(\cite{ es1993, trv1992}).
\end{itemize}

\smallskip

\noindent
To motivate the present Part~III, let us expose informally two dual
questions about systems of vector fields. Consider a set $\LL$ of
local vector fields defined on a domain of $\R^n$. Frobenius' theorem
provides local foliations by submanifolds to which every element of
$\LL$ is tangent, provided $\LL$ is closed under Lie
brackets. However, for a generic set $\LL$, the condition $\left[ \LL,
\LL \right] \subset \LL$ fails and in addition, the tangent spaces
spanned by elements of $\LL$ are of varying dimension. To surmount
these imperfections, two inverse options present themselves:

\smallskip

\begin{itemize}

\item[{\bf Sub:}]
find the subsystems $\LL ' \subset \LL$ which satisfy Frobenius'
condition $\left[ \LL ', \LL ' \right] \subset \LL '$ and which are
maximal in an appropriate sense;

\smallskip

\item[{\bf Sup:}]
find the supsystems $\LL ' \supset \LL$ which have integral manifolds
and which are minimal in an appropriate sense.

\end{itemize}

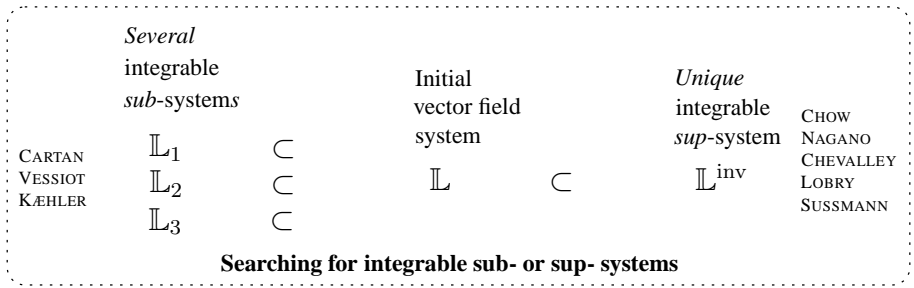
\begin{center}
\input sub-sup-distributions.pstex_t
\end{center}

The first problem {\bf Sub} is answered by the Cartan-Vessiot-K\"ahler
theorem, thanks to an algorithm which provides all the minimal
Frobenius-integrable subsystems $\LL'$ of $\LL$ (we recommend \cite{
stk2000} for a presentation). Generically, there are infinitely many
solutions and their cardinality is described by means of a sequence of
integers together with the so-called {\sl Cartan character} of
$\LL$. In the course of the proof, the Cauchy-Kowalevskaya
integrability theorem, valid only in the analytic category, is heavily
used. It was not a serious restriction at the time of \'E.~Cartan,
but, in the second half of the twentieth century, the progress of the
Analysis of PDE showed deep new phenomena in the differentiable
category. Hence, one may raise the:

\def\theopenproblem{1.2}\begin{openproblem}
Find versions of the Cartan-Vessiot-K\"ahler theorem for systems of
vector fields having smooth non-analytic coefficients.
\end{openproblem}

The Cauchy characteristic subsystem of $\LL$ (\cite{ stk2000}) is
always involutive, hence the smooth Frobenius theorem applies to
it\footnote{ We are grateful to Stormark for pointing out this
observation}. However, for intermediate systems, the question is wide
open. Possibly, this question is related to some theorems about local
solvability of smooth partial differential equations ({\it cf.}
Section~3) that were established to understand the Hans Lewy
counterexample (\S3.1).

The second problem {\bf Sup} is already answered by Nagano's theorem
(Part~II), though only in the analytic category, with a unique
integrable minimal supsystem $\LL^{\rm lie} \supset \LL$. In the
general smooth category, the stronger Chevalley-Lobry-Stefan-Sussmann
theorem, dealing with flows of vector fields instead of Lie brackets,
shows again that there is a unique integrable sup-system of $\LL$
which has integral manifolds. As this theorem will be central in this
memoir, it will be exposed thoroughly in the present Section~1.

\subsection*{1.3.~Flows of vector fields and their regularity}
Let $\K = \R$ or $\C$. Let $D$ be a open connected subset of
$\K^n$. Let ${\sf x} = ({\sf x}_1, \dots, {\sf x}_n )\in D$. Let
$L=\sum_{ i=1 }^n \, a_i({\sf x})\, \frac{ \partial}{ \partial {\sf
x}_i}$ be a vector field defined over $D$. Throughout this section, we
shall assume that its coefficients $a_i$ are either $\K$-analytic (of
class $\mathcal{ C}^\omega$), of class $\mathcal{ C }^\infty$, or of
class $\mathcal{ C }^{ \kappa, \alpha}$, where $\kappa \geqslant 1$
and $0 \leqslant \alpha \leqslant 1$ ({\it see}\, Section~1(IV) for
background about H\"older classes).

By the classical Cauchy-Lipschitz theorem,
through each point ${\sf x}_0 \in D$, there passes a unique local {\sl
integral curve} of the vector field $L$, namely a local solution ${\sf
x}( {\sf t})= ({\sf x}_1( {\sf t} ), \dots, {\sf x}_n ( {\sf t}))$ of
the system of ordinary differential equations:
\[
d{\sf x}_1({\sf t})/d{\sf t}
=
a_1({\sf x}({\sf t})),\dots\dots,\ d{\sf
x}_n({\sf t})/d{\sf t}
=
a_n({\sf x}({\sf t})),
\]
which satisfies the initial condition ${\sf x}(0)={\sf x}_0$. This
solution is defined at least for small ${\sf t} \in \K$ and is
classically denoted by ${\sf t} \mapsto \exp ({\sf t} L) ({\sf x}_0)$,
because it has the local pseudogroup property 
\[
\exp (
{\sf t}'L)
\big(
\exp
({\sf t}L)({\sf x}_0)
\big)
= 
\exp 
\big(
({\sf t}+{\sf t}')L
\big)
({\sf x}_0),
\]
whenever the composition is defined. Denote by $\Omega_{ {\sf x}_0}$
the largest connected open set containing the origin in $\K$ in which
$\exp ({\sf t} L) ({\sf x}_0)$ is defined. One shows that the union of
various $\Omega_{ {\sf x}_0}$, for ${\sf x}_0$ running in $D$, is an
open connected set $\Omega_L$ of $\K \times \K^n$ which contains
$\{0\} \times D$. Some regularity with respect to both variables ${\sf
t}$ and ${\sf x}_0$ is got automatically.

\def\thetheorem{1.4}\begin{theorem}
{\rm (\cite{ la1983}, [$*$])}
The global flow $\Omega_L \ni ( {\sf t}, {\sf x}_0) \mapsto \exp (
{\sf t} \, L) ({\sf x}_0) \in D$ of a vector field $L = \sum_{
i=1}^n\, a_i({\sf x}) \, \partial_{ {\sf x}_i}$ defined in the domain
$D$ has exactly the {\rm same smoothness} as $L$, namely it is
$\mathcal{ C}^\omega$, $\mathcal{ C}^\infty$ or $\mathcal{ C
}^{\kappa, \alpha}$.
\end{theorem}

As a classical corollary, a local straightening property holds~: in a
neighborhood of a point at which $L$ does not vanish, there exists a
$\mathcal{ C}^\omega$, $\mathcal{ C}^\infty$ or $\mathcal{ C}^{\kappa,
\alpha}$ change of coordinates ${\sf x}'= {\sf x}'({\sf x})$ in which
the transformed vector field is the unit positive vector field
directed by the ${\sf x}_1'$ lines, viz $L'= \partial / \partial {\sf
x}_1'$.

Up to the end of this Section~1, we will work with $\K = \R$.

\subsection*{1.5.~Searching integral manifolds of a system of
vector fields} Let $M$ be a smooth paracompact {\it real}\, manifold,
which is $\mathcal{ C }^\omega$, $\mathcal{ C}^\infty$ or $\mathcal{
C}^{ \kappa+1, \alpha }$, where $\kappa \geqslant 1$, $0 \leqslant
\alpha \leqslant 1$. Let $\LL:=\{L_a \}_{a \in A}$ be a collection of
vector fields defined on open subsets $D_a$ of $M$ and having
$\mathcal{ C }^\omega$, $\mathcal{ C }^\infty$ or $\mathcal{ C }^{
\kappa, \alpha}$ coefficients, where $A$ is an arbitrary set. It is no
restriction to assume that $\cup_{ a \in A} \, D_a = M$, since
otherwise, it suffices to shrink $M$. Call $\LL$ a {\sl system of
vector fields on $M$}.

\def\theproblem{1.6}\begin{problem}
Find submanifolds $N$ of $M$ such that each element of $\LL$ is
tangent to $N$.
\end{problem}

\noindent
To analyze this (still imprecise) problem, let $\F_M$ denote the
collection of all $\mathcal{ C}^\omega$, $\mathcal{ C}^\infty$ or
$\mathcal{ C}^{ \kappa, \alpha}$ functions defined on open subsets
of $M$, and call the system $\LL$ of vector fields {\sl $\F_M$-linear}
if every combined vector field $fK + gL$ belongs to $\LL$, whenever $f,g
\in \F_M$ and $K,L\in\LL$. Here, $fK + gL$ is defined in the
intersection of the domains of definition of $f$, $g$, $K$ and $L$. To
study the problem, it is obviously no restriction to assume that $\LL$
is $\F_M$-linear.

For $p\in M$ arbitrary, define
\[
\LL(p) 
:=
\{L(p):\ 
L\in\LL\}.
\]
Since $\LL$ is $\F_M$-linear, this is a linear subspace of
$T_pM$. So Problem~1.6 is to find submanifolds $N$ satisfying $T_p N
\supset \LL (p)$, for every $p\in N$. Notice that an appropriate
answer should enable one to find {\it all}\, such submanifolds. Also,
suppose that $N_1$ and $N_2$ are two solutions with $N_2 \subset
N_1$. Then the problem with the pair $(M, N)$ is exactly the same as
the problem with the pair $(N_1, N_2)$. Hence a better formulation.

\def\theproblem{1.6'}\begin{problem}
Find {\rm all} the submanifolds $N \subset M$ {\rm of
smallest dimension} that satisfy
$T_p N \supset \LL (p)$, for every $p \in N$.
\end{problem}

The classical {\sl Frobenius theorem} (\cite{ fr1877, sp1970,
ber1999, bo1991, ch1991, stk2000, 
trv1992}) provides an answer in the (for us simplest) case
where $\LL$ is closed under Lie brackets and is of constant dimension:
{\it every point $p \in M$ admits an open neighborhood foliated by
submanifolds $N$ satisfying $T_q N = \LL (q)$, for every $q \in
N$}. The global properties of these submanifolds were not much studied
until C.~Ehresmann and G.~Reeb endeavoured to understand them (birth
of foliation theory). A line with irrational slope in the $2$-torus
$(\R / \Z)^2$ shows that it is necessary to admit submanifolds $N$ of
$M$ which are not closed. Let $\mathcal{A}_M$ denote the manifold
structure of $M$.

\def\thedefinition{1.7}\begin{definition}{\rm
An {\sl immersed submanifold}\, of $(M, \mathcal{ A}_M)$ is a subset
of $N$ of $M$ equipped with its own smooth manifold structure
$\mathcal{ A}_N$, such that the inclusion map $i: (N,\mathcal{ A}_N)
\to (M, \mathcal{ A}_M)$ is smooth, immersive and injective.
}\end{definition}

Thus, to keep maximally open Problem~1.6', one should seek 
immersed submanifolds and make no assumption about closedness under
Lie brackets. For later use, recall that
an immersed submanifold $N$ of $M$ is {\sl embedded} if its own
manifold structure coincides with the manifold inherited from
the inclusion $N\subset M$. It is well known (\cite{ cln1985})
that an immersed submanifold $N$ is embedded {\rm if and only if} for
every point $p \in N$, there exists a neighborhood $U_p$ of $p$ in $M$
such that the pair $(U_p, N\cap U_p)$ is diffeomorphic to
$(\R^{\dim M}, \R^{ \dim N})$.

\subsection*{1.8.~Maximal strong integral manifolds property}
In order to understand Problem~1.6', for heuristic reasons, it will be
clever to discuss the differences between the two possibilities $\LL
(p) = T_p N$ and $\LL (p) \varsubsetneq T_p N$. Consider an arbitrary
$\F_M$-linear system of vector fields $\widehat{ \LL}$ containing
$\LL$, for instance $\LL$ itself. Let $p \in M$ and define the linear
subspace $\widehat{ \LL} (p) := \{\widehat{ L}(p): \, \widehat{ L} \in
\widehat{ \LL} \}$.

\def\thedefinition{1.9}\begin{definition}{\rm 
An immersed submanifold $N$ of $M$ is said to be{\rm :}

\begin{itemize}

\smallskip\item[$\bullet$]
a {\sl strong $\widehat{ \LL}$-integral manifold} if $T_q N =
\widehat{ \LL} (q)$, at every point $q \in N$;

\smallskip\item[$\bullet$] 
a {\sl weak $\LL$-integral manifold} if $T_q N \supset
\widehat{ \LL} (q)$, at every point $q \in N$.
\end{itemize}\smallskip

}\end{definition}

In advance, the answer (Theorem~1.21 below) to Problem~1.6' states
that it is possible to construct a {\it unique}\, system of vector
fields $\widehat{ \LL}$ containing $\LL$, whose {\it strong}\,
integral manifolds coincide with the smallest {\it weak}\,
$\LL$-integral manifolds $N$. Further definitions are needed.

A system of vector fields $\widehat{ \LL}$ is said to have the {\sl
strong integral manifolds property} if for every point $p \in M$,
there exists a strong $\widehat{ \LL}$-integral submanifold $N$
passing through $p$. A {\sl maximal strong $\widehat{ \LL}$-integral
manifold} $N$ is an immersed $\widehat{ \LL}$-integral manifold with
the property that every connected strong $\widehat{ \LL}$-integral
manifold which intersects $N$ is an open submanifold of $N$. Thus,
through a point $p\in M$, there passes at most one maximal strong
$\widehat{ \LL}$-integral submanifold. Finally, the system $\widehat{
\LL}$ has the {\sl maximal strong integral manifolds} property if,
through every point $p\in M$, there passes a maximal strong $\widehat{
\LL}$-integral manifold. The $\F_M$-linear systems $\widehat{ \LL}$
containing $\LL$ are ordered by inclusion. We then admit that
Problem~1.6' is essentially reduced to:

\def\theproblem{1.6''}\begin{problem} 
How to construct the {\rm (}a posteriori unique{\rm )} smallest {\rm
(}for inclusion{\rm )} $\F_M$-linear system of vector fields
$\widehat{ \LL}$ containing $\LL$ which has the maximal strong
integral manifolds property\,?
\end{problem}

\subsection*{1.10.~Taking account of the Lie brackets}
Here is a basic geometric observation inspired by Frobenius' and
Nagano's theorems.

\def\thelemma{1.11}\begin{lemma} 
Assume the $\F_M$-linear system $\widehat{ \LL}$ has the strong
integral manifolds property. Then for every two vector fields
$\widehat{L}, \, \widehat{L}'\in\widehat{ \LL}$ and for every $p$ in
the intersection of their domains, the Lie bracket $\big[ \widehat{L},
\widehat{L}' \big] (p)$ belongs to $\widehat{ \LL} (p)$.
\end{lemma}

\proof
Indeed, let $N$ be a strong $\widehat{ \LL}$-integral manifold, namely
satisfying $TN = \widehat{ \LL} \big\vert_N$. If $\widehat{L}, \,
\widehat{L}'\in\widehat{ \LL}$, the two restrictions $\widehat{ L}
\big\vert_N$ and $\widehat{ L}' \big\vert_N$ are tangent to $N$. Hence
the restriction to $N$ of the Lie bracket $\big[ \widehat{ L},
\widehat{ L}' \big]$ is also tangent to $N$.
In conclusion, at every $p\in N$,
we have $\big[ \widehat{L}, \widehat{L}' \big] (p) \in T_p N =
\widehat{ \LL} (p)$.
\endproof

So it is a temptation to believe that the smallest system $\LL^{\rm
lie}$ of vector fields containing $\LL$ which is closed under Lie
brackets does enjoy the maximal integral manifolds property. However,
just after the statement of Nagano's theorem (Part~II), we have
already learnt by means of Example~1.6(II) that in the $\mathcal{ C
}^\infty$ and $\mathcal{ C}^{ \kappa, \alpha}$ categories, the
consideration of $\LL^{\rm lie}$ is inappropriate.

\subsection*{1.12.~Transport of a vector field by the
flow of another vector field} To understand why $\LL^{\rm lie}$ is
insufficient, it will be clever to recall one of the classical
definitions of the Lie bracket between two vector fields. Let $p \in
M$ and let $K$ be a vector field defined in a neighborhood of
$p$. Denote by $K (q)$ the {\sl value} of $K$ at a point $q$ (this is
a vector in $T_q M$), by $g_* (K)$ the push-forward of $K$ by a local
diffeomorphism $g$, and by $q \mapsto K_{ \sf s} (q)$ [instead of $q
\mapsto \exp( {\sf s} K) (q)$] the local diffeomorphism at time ${\sf
s}$ induced by the flow of $K$. If $L$ is a second vector field
defined in a neighborhood of $p$, the {\sl Lie bracket between $K$ and
$L$ at $p$} is defined by:
\def\theequation{1.13}\begin{equation}
[K,L] (p) 
:=
\lim_{{\sf s}\to 0}\,\left(
\frac{L(p)-(K_{\sf s})_*(L(K_{-{\sf s}}(p)))
}{{\sf s}}\right).
\end{equation}
Observe that for every fixed ${\sf s}
\neq 0$, the two vectors $L (p)$ and $(K_{ \sf s})_* 
(L( K_{ -{\sf s}} (p )))$ belong $T_pM$. 

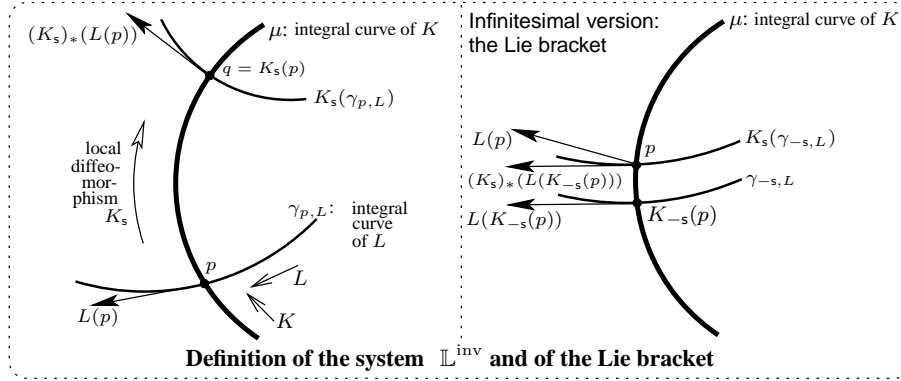
\begin{center}
\input S-invariant-new.pstex_t
\end{center}

We explain how to read the right hand side of the diagram. In it: the
integral curve of $K$ passing through $p$ is denoted by $\mu$; the
integral curve of $L$ passing through the point $K_{ -{\sf s}} (p)$
for ${\sf s}$ very small is denoted by $\gamma_{ -{\sf s}, L}$; its
image by the local diffeomorphism $K_{\sf s}$ is denoted by $K_{\sf s}
(\gamma_{ -{\sf s}, L})$; the vector $L(K_{ -{\sf s}} (p))$ is tangent
to $\gamma_{ -{\sf s}, L}$ at the point $K_{-{\sf s}}(p)$; the vector
$(K_{\sf s})_*(L(K_{ -{\sf s}}(p)))$, transported by the
differerential of $K_{\sf s}$, is in general distinct from the vector
$L(p)$; in fact, the difference $L(p) - (K_{\sf s} )_*(L(K_{- {\sf s}
}(p)))$ divided by ${\sf s}$, tends to $[ K, L] (p)$ as ${\sf s} \to
0$.

Essentially, $\LL^{ \rm lie}$ collects all vector fields obtained by
taking infinitesimal differences~\thetag{ 1.13} between vectors $L (p)$
and transported vectors $(K_{\sf s} )_*( L(K_{ -{\sf s}} (p)))$,
and then iterating this processus to absorb all multiple Lie brackets.

As suggested in the left hand side of the diagram, {\it instead of
taking the infinitesimal differences, it is more general to collect
all the vectors of the form $(K_{ \sf s} )_* (L(p ))$}. This is the
clue of Sussmann's theorem. In fact, the system $\widehat{ \LL}$ which
is sought for in Problem~1.6'' should not only contain $\LL^{ \rm
lie}$, but should also collect all the vector fields of the form $(K_{
\sf s })_* (L)$, where ${\sf s}$ is {\it not}\, an infinitesimal.

\def\thelemma{1.14}\begin{lemma}
Let $\widehat{ \LL}$ be a $\F_M$-linear system of vector fields
containing $\LL$ which has the strong integral manifolds property. Let
$p\in M$, let $K,L\in \LL$ be two arbitrary vector fields defined in a
neighborhood of $p$ and let $q = K_{ \sf s}(p)$ be a point in the
integral curve of $K$ issued from $p$, with ${\sf s} \in \R$
small. Then the linear subspace $\widehat{ \LL} (q)$ necessarily
contains the transported vector $(K_{ \sf s} )_* (L(p))$.
\end{lemma}

\proof
Let $N$ be a strong $\widehat{ \LL}$-integral manifold passing through
$p$. As $\widehat{ \LL}(r) = T_rN$ at every point $r\in N$, and as
$\LL$ is contained in $\widehat{ \LL}$, it follows that the restricted
vector field $K\vert_N$ is tangent to $N$. Consequently, the integral
curve of $K$ issued from $p$ is locally contained in $N$, hence the
point $q = K_{\sf s} (p)$ belongs to $N$.

Moreover, as $\LL$ is contained in $\widehat{ \LL}$, the vector $L(p)$
is tangent to $N$ at $p$. The differential $(K_{ \sf s})_*$ being a
linear isomorphism between $T_p N$ and $T_q N$, it follows that
the vector $(K_{ \sf s} )_* (L(p))$ belongs to the tangent 
space $T_qN$, which coincides with $\widehat{ \LL} (q)$ by assumption.
\endproof

\subsection*{ 1.15.~The smallest $\LL$-invariant system of vector fields
$\LL^{\rm inv}$} Based on this crucial observation, we may introduce
the smallest $\F_M$-linear system of vector fields $\LL^{\rm inv}$
(``inv'' abbreviates ``invariant'') containing $\LL$ which contains
all vectors of the form $(K_{ \sf s} )_* (L)$, whenever $K, L \in \LL$
and ${\sf s} \in \R$. It follows that $(K_{\sf s})_* \left( \LL^{\rm
inv} (p) \right) = \LL^{\rm inv} \left( K_{\sf s}(p) \right)$: the
distribution of linear subspaces $p \mapsto \LL^{\rm inv} (p) \subset
T_p M$ is {\it invariant}\, under the local flow maps.

In~\cite{ su1973}, it is shown that $\LL^{\rm inv}$ is concretely and
finitely generated as stated in Lemma~1.16 below. At first, some more
notation is needed to denote the composition of several local
diffeomorphisms of the form $K_{\sf s}$. Let $\X$ denote the system of
all tangent vector fields to $M$, defined on open subsets of $M$. Let
$k\in\N$ with $k\geqslant 1$ and let $K= (K^1, \dots, K^k)\in \X^k$ be a
$k$-tuple of vector fields defined in their domains of
definition. If ${\sf s} =({\sf s}_1, \dots, {\sf s}_k) \in \R^k$ is a
$k$-tuple of ``time'' parameters, we will denote by $K_{\sf s} (p)$
the point
\[
K_{{\sf s}_1}^1\left(\cdots 
(K_{{\sf s}_k}^k(p))\cdots \right)
:=
\exp\left({\sf s}_1K^1
\left(\cdots(\exp({\sf s}_k K^k(p)))\cdots\right)\right),
\]
whenever the composition is defined. The $k$-tuple $({\sf s}_1, \dots,
{\sf s}_k)$ will also be called a {\sl multitime parameter}. For
${\sf s}$ fixed, the map $p\mapsto K_{\sf s} (p)$ is a local
diffeomorphism between two open subsets of $M$. Its local inverse is
the map $p\mapsto \widetilde{K}_{ -\widetilde{\sf s}}(p)$, where
$\widetilde{ K} := (K^k, \dots, K^1) \in \LL^k$ and $\widetilde{ \sf
s} := ( {\sf s}_k, \dots, {\sf s}_1)$. Moreover, if we define $({\sf
s}, {\sf s}') := ( {\sf s}_1, \dots, {\sf s}_k, {\sf s}_1', \dots,
{\sf s}_{k'}')$ for general ${\sf s} = ({\sf s}_1, \dots, {\sf s}_k)
\in \R^k$ and ${\sf s}' = ( {\sf s}_1', \dots, {\sf s}_{ k'}') \in
\R^{ k'}$, we have $K_{ {\sf s}' }'\circ K_{\sf s} = (K', K)_{({\sf
s}', {\sf s})}$. 

After shrinking the domains of definition, the
composition of local diffeomorphisms $K_{\sf s}$ is clearly
associative, where it is defined. It follows that the set of local
diffeomorphisms $K_{\sf s}$ constitutes a {\sl pseudogroup of local
diffeomorphisms}. Here, the term ``pseudo'' stems from the fact that
the domains of definitions have to be adjusted; not all compositions
are allowed.

\def\thelemma{1.16}\begin{lemma}
{\rm (\cite{ su1973})}
The system $\LL^{\rm inv}$ is generated by the $\F_M$-linear
combinations of all vector fields of the form $(K_{\sf s})_*(L)$, for
all $L\in \LL$, all $k$-tuples $K= (K^1, \dots, K^k) \in \LL^k$
of elements of $\LL$ and all multitime parameters ${\sf s} = ({\sf
s}_1, \dots, {\sf s}_k) \in \R^k$.
\end{lemma}

The definitions and the above reasonings show that the $\F_M$-linear
system $\LL^{\rm lie}$ is a subsystem of the $\F_M$-linear system
$\LL^{\rm inv}$ (of course, every system is contained in $\X$):
\[
\boxed{
\LL\subset \LL^{\rm lie} \subset \LL^{\rm inv}\subset \X}.
\]
 
In general, at a fixed point $p\in M$, the inclusions $\LL(p)\subset
\LL^{\rm lie}(p)\subset \LL^{\rm inv}(p)\subset 
\X (p) = T_pM$ may be all strict.

\def\theexample{1.17}\begin{example}
{\rm
On $\R^4$, consider the system $\LL$ generated by the three
vector fields 
\[
\frac{\partial }{\partial {\sf x}_1}, 
\ \ \ \ \ \ \ 
{\sf x}_1 \,
\frac{\partial }{\partial {\sf x}_2}, 
\ \ \ \ \ \ \ 
e^{- 1 / {\sf x}_1^2} \, 
\frac{\partial }{\partial {\sf x}_3}.
\]
Then it may be checked that
\[
\left\{
\aligned
\LL (0) 
&
=
\R \, \partial_{{\sf x}_1}, \\
\LL^{\rm lie} (0) 
&
=
\R \, \partial_{{\sf x}_1} \oplus 
\R \, \partial_{{\sf x}_2}, \\
\LL^{\rm inv} (0) 
&
=
\R \, \partial_{{\sf x}_1} \oplus 
\R \, \partial_{{\sf x}_2} \oplus 
\R \, \partial_{{\sf x}_3}, \\
\X (0) 
&
=
\R \, \partial_{{\sf x}_1} \oplus 
\R \, \partial_{{\sf x}_2} \oplus 
\R \, \partial_{{\sf x}_3} \oplus 
\R \, \partial_{{\sf x}_4}.
\endaligned\right.
\]

}
\end{example}

\def\thetheorem{1.18}\begin{theorem}
{\rm (\cite{ na1966, su1973})}
In the $\mathcal{ C}^\omega$, $\mathcal{ C}^\infty$ and $\mathcal{
C}^{\kappa, \alpha}$ categories, the system $\LL^{\rm inv}$ is the
smallest one containing $\LL$ that has the maximal strong integral
manifolds property. In the $\mathcal{ C}^\omega$ category, $\LL^{\rm
inv} = \LL^{\rm lie}$.
\end{theorem}

Further structural properties remain to be explained.

\subsection*{1.19.~$\LL$-orbits} 
The maximal strong integral manifolds of $\LL^{\rm inv}$
may be defined directly by means of $\LL$, without refering to
$\LL^{\rm inv}$, as follows. Two points $p,q\in M$ are said to be {\sl
$\LL$-equivalent} if there exists a local diffeomorphism of the form
$K_{\sf s}$, $K= (K^1, \dots, K^k)$, ${\sf s} = ({\sf s}_1, \dots,
{\sf s}_k)$, $k\in \N$, with $K_{\sf s}(p) =q$. This clearly defines
an equivalence relation on $M$. The equivalence classes are called the
{\sl $\LL$-orbits} of $M$ and will be denoted either by $\mathcal{
O}_\LL(p)$ or shortly by $\mathcal{O }_\LL$, when the reference to one
point of the orbit is superfluous.

Concretely, two points $p,q \in M$ belong to the same $\LL$-orbit if
and only if there exist a continuous curve $\gamma : [0,1]\to M$ with
$\gamma(0) = p$ and $\gamma(1) =q$ together with a partition of the
interval $[0,1]$ by numbers $0= {\sf s}_0 < {\sf s}_1<{\sf s}_2
<\cdots<{\sf s}_k =1$ and vector fields $K^1,\dots, K^k\in \LL$ such
that for each $i = 1, \dots, k$, the restriction of $\gamma$ to the
subinterval $[{\sf s}_{i-1},{\sf s}_i]$ is an integral curve of
$K^i$. Such a curve will be called a {\sl piecewise integral curve of}
\, $\LL$.

Let $p\in M$. Then its $\LL$-orbit $\mathcal{ O }_\LL (p)$ may be
equipped with the finest topology which makes all the maps ${\sf s}
\mapsto K_{ \sf s} (p)$ continuous, for all $k \geqslant 1$, all $K = (K^1,
\dots, K^k) \in \LL^k$ and all multitime parameters ${\sf s} = ({\sf
s}_1, \dots, {\sf s}_k)$. This topology is independent of the choice
of a central point $p$ inside a given orbit (\cite{ su1973}). Since
the maps $\R^k \ni {\sf s} \mapsto K_{\sf s} (p) \in M$ are already
continuous, the topology of $\mathcal{ O}_\LL (p)$ is always finer
than the topology induced by the inclusion $\mathcal{ O }_\LL (p)
\subset M$. It follows that the inclusion map from $\mathcal{ O }_\LL
(p)$ into $M$ is continuous. In particular, $\mathcal{ O }_\LL (p)$ is
Hausdorff.

\subsection*{ 1.20.~Precise statement of the orbit theorem} 
We now state in length the fundamental theorem of Sussmann, based on
preliminary versions due to Hermann (\cite{ he1963}), to Nagano (\cite{
na1966}) and to Lobry (\cite{ lo1970}). It describes $\LL$-orbits as
immersed submanifolds {\bf (1)}, {\bf (2)} enjoying the everywhere
accessibility conditions {\bf (3)}, {\bf (4)}, together with a local
flow-box property {\bf (5)}, useful in applications.

\def\thetheorem{1.21}\begin{theorem}
{\rm ({\sc Sussmann}
\cite{ su1973, trv1992, bm1997, ber1999,
bch2005}, [$*$])} The following five properties hold true.

\smallskip

\begin{itemize}
\item[{\bf (1)}]
Every $\LL$-orbit $\mathcal{ O}_\LL$, equipped with the finest
topology which makes all the maps ${\sf s} \mapsto K_{\sf s} (p)$
continuous, admits a {\rm unique differentiable structure} with the
property that $\mathcal{ O}_\LL$ is an {\rm immersed} submanifold of
$M$, of class $\mathcal{ C}^\omega$, $\mathcal{ C}^\infty$ or
$\mathcal{ C}^{\kappa, \alpha}$.

\smallskip

\item[{\bf (2)}]
With this topology, each $\LL$-orbit $\mathcal{ O}_\LL$ is
simultaneously a {\rm (}connected{\rm )} maximal weak integral
manifold of $\LL$ and a {\rm (}connected{\rm )} maximal strong
integral manifold of the $\LL$-invariant $\F_M$-linear system
$\LL^{\rm inv}${\rm ;} thus, for every point $p \in M$, it holds $T_p
\mathcal{ O}_\LL (p) = \LL^{\rm inv} (p)$, whence in particular $ \dim
\LL^{ \rm inv} (q) = \dim \mathcal{ O }_\LL (p)$ is constant for all
$q$ belonging to a given $\LL$-orbit $\mathcal{ O}_\LL (p)$.

\smallskip

\item[{\bf (3)}] 
For every $p\in M$, $k\geqslant 1$, $K\in \LL^k$, ${\sf s} \in \R^k$ such
that $K_{\sf s} (p)$ is defined, the differential map $(K_{\sf s})_*$
makes a linear isomorphism from $T_p \mathcal{ O}_\LL (p) = \LL^{\rm
inv}(p)$ onto $T_{K_{\sf s} (p)} \mathcal{ O}_\LL(p) = \LL^{\rm inv}
(K_{\sf s} (p))$.

\smallskip

\item[{\bf (4)}]
For every $p,q\in M$ belonging to the {\rm same} $\LL$-orbit, there
exists an integer $k\geqslant 1$, there exists a $k$-tuple of vector fields
$K= (K^1, \dots, K^k)\in \LL^k$ and there exists a multitime ${\sf
s}^* = ({\sf s}_1^*, \dots, {\sf s}_k^*) \in \R^k$ such that $p =
K_{{\sf s}^*} (q)$ and such that the differential at ${\sf s}^*$ of
the map
\def\theequation{1.22}\begin{equation}
\R^k\ni{\sf s}
\mapsto 
K_{\sf s}(q)
\in
\mathcal{O}_\LL(p) 
\end{equation}
is of rank equal to $\dim \mathcal{ O}_\LL(p)$.

\smallskip

\item[{\bf (5)}] 
For every $p \in M$, there exists an open connected neighborhood $V_p$
of $p$ in $M$ and there exists a $\mathcal{ C}^\omega$, $\mathcal{
C}^\infty$ or $\mathcal{ C}^{ \kappa,\alpha}$ diffeomorphism
\def\theequation{1.23}\begin{equation}
\square^e\times 
\square^{n-e}\ni 
({\sf s}, {\sf r}) 
\longmapsto 
\varphi({\sf s}, {\sf r})\in V_p,
\end{equation}
where $e = \dim \mathcal{ O }_\LL (p)$, where $\square =\{ {\sf x} \in
\R: \, \vert {\sf x} \vert < 1\}$, such that{\rm :}

\smallskip

\begin{itemize}
\item[$\bullet$]
$\varphi(0,0)=p${\rm ;}

\smallskip

\item[$\bullet$] 
the plaque $\varphi \left( \square^e \times \{0 \} \right)$ is an open
piece of the $\LL$-orbit of $p${\rm ;}

\smallskip

\item[$\bullet$] 
each plaque $\varphi \left( \square^e \times \{ {\sf r} \} \right)$ is
contained in a single $\LL$-orbit{\rm ;} and{\rm :}
\smallskip

\item[$\bullet$] 
the set of ${\sf r} \in \square^{ n -e}$ such that
$\varphi \big( \square^e \times \{ {\sf r} \} \big)$ is
contained in the same $\LL$-orbit
$\mathcal{ O}_\LL (p)$ is either finite or
countable.

\end{itemize}
\end{itemize}
\end{theorem}

In general, for ${\sf r} \neq 0$, the $e$-dimensional plaques $\varphi
(\square^e \times \{ {\sf r} \})$ have positive codimension in the
nearby orbits. We draw a diagram, in which $e= \dim \mathcal{
O}_\LL (p) = 1$, with the nearby $\LL$-orbits $\mathcal{ O}_2$,
$\mathcal{ O}_2'$, $\mathcal{ O}_3$ and $\mathcal{ O}_3'$ having
dimensions $2$, $2$, $3$ and $3$.

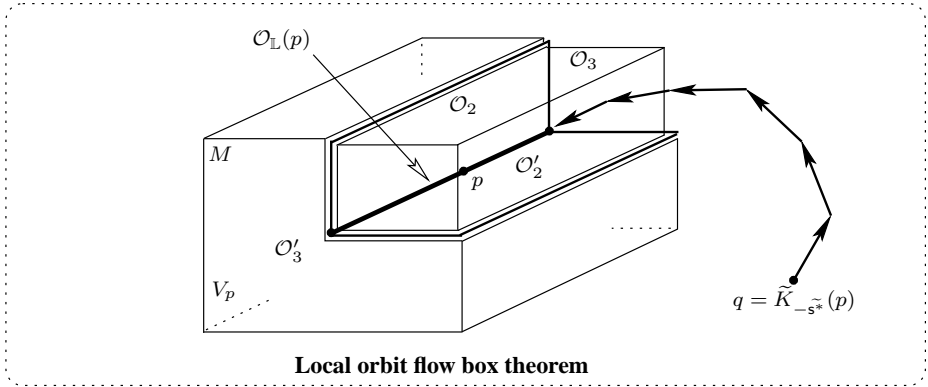
\begin{center}
\begin{picture}(0,0)%
\includegraphics{flowbox.pstex}%
\end{picture}%
\setlength{\unitlength}{3947sp}%
\begingroup\makeatletter\ifx\SetFigFont\undefined
\def\x#1#2#3#4#5#6#7\relax{\def\x{#1#2#3#4#5#6}}%
\expandafter\x\fmtname xxxxxx\relax \def\y{splain}%
\ifx\x\y   
\gdef\SetFigFont#1#2#3{%
  \ifnum #1<17\tiny\else \ifnum #1<20\small\else
  \ifnum #1<24\normalsize\else \ifnum #1<29\large\else
  \ifnum #1<34\Large\else \ifnum #1<41\LARGE\else
     \huge\fi\fi\fi\fi\fi\fi
  \csname #3\endcsname}%
\else
\gdef\SetFigFont#1#2#3{\begingroup
  \count@#1\relax \ifnum 25<\count@\count@25\fi
  \def\x{\endgroup\@setsize\SetFigFont{#2pt}}%
  \expandafter\x
    \csname \romannumeral\the\count@ pt\expandafter\endcsname
    \csname @\romannumeral\the\count@ pt\endcsname
  \csname #3\endcsname}%
\fi
\fi\endgroup
\begin{picture}(5799,2424)(1718,-2782)
\put(4656,-1518){\makebox(0,0)[lb]{\smash{\SetFigFont{8}{9.6}{rm}{\color[rgb]{0,0,0}$p$}%
}}}
\put(3023,-2205){\makebox(0,0)[lb]{\smash{\SetFigFont{8}{9.6}{rm}{\color[rgb]{0,0,0}$V_p$}%
}}}
\put(3001,-1359){\makebox(0,0)[lb]{\smash{\SetFigFont{8}{9.6}{rm}{\color[rgb]{0,0,0}$M$}%
}}}
\put(4927,-1427){\makebox(0,0)[lb]{\smash{\SetFigFont{8}{9.6}{rm}{\color[rgb]{0,0,0}$\mathcal{O}_2'$}%
}}}
\put(5270,-763){\makebox(0,0)[lb]{\smash{\SetFigFont{8}{9.6}{rm}{\color[rgb]{0,0,0}$\mathcal{O}_3$}%
}}}
\put(4514,-1045){\makebox(0,0)[lb]{\smash{\SetFigFont{8}{9.6}{rm}{\color[rgb]{0,0,0}$\mathcal{O}_2$}%
}}}
\put(3417,-1952){\makebox(0,0)[lb]{\smash{\SetFigFont{8}{9.6}{rm}{\color[rgb]{0,0,0}$\mathcal{O}_3'$}%
}}}
\put(6296,-2271){\makebox(0,0)[lb]{\smash{\SetFigFont{8}{9.6}{rm}{\color[rgb]{0,0,0}$q=\widetilde{K}_{-\widetilde{{\sf s}^*}}(p)$}%
}}}
\put(3278,-632){\makebox(0,0)[lb]{\smash{\SetFigFont{8}{9.6}{rm}{\color[rgb]{0,0,0}$\mathcal{ O}_\LL(p)$}%
}}}
\put(3544,-2689){\makebox(0,0)[lb]{\smash{\SetFigFont{9}{10.8}{rm}{\color[rgb]{0,0,0}{\bf Local orbit flow box theorem}}%
}}}
\end{picture}

\end{center}

Property {\bf (4)} is crucial: the maps~\thetag{ 1.22} of rank $\dim
\mathcal{ O}_\LL (p)$ are used to define the differentiable structure
on $\mathcal{ O}_\LL (p)$; they are also used to obtain the local
orbit flow box property {\bf (5)}, as follows.

Let $p\in M$ and choose $q\in \mathcal{ O}_\LL (p)$ with $q\neq p$, to
fit with the diagrams ($q = p$ would also do). Assuming that {\bf (4)}
holds, set $e := \dim \mathcal{ O}_\LL (p)$, introduce an open subset
$T_e$ in some $e$-dimensional affine subspace passing through ${\sf
s}^*$ in $\R^k$ so that the restriction of the map~\thetag{ 1.23} to
$T_e$ still has rank $e$ at ${\sf s} = {\sf s}^*$. Introduce also an
$(n-e)$-dimensional local submanifold $\Lambda_p$ passing through $p$
with $T_p \Lambda_p \oplus T_p \mathcal{ O}_\LL (p) = T_p M$ and set
$\Lambda_q := \widetilde{ K}_{ - \widetilde{ \sf s}} (
\Lambda_p)$. Notice that $T_q \Lambda_q \oplus T_q \mathcal{ O}_\LL
(p) = T_q M$, since the multiple flow map $K_{\sf s} (\cdot)$
stabilizes $\mathcal{ O}_\LL (p)$. Then, as one of the possible maps
$\varphi$ whose existence is claimed in {\bf (5)}, we may choose a
suitable restriction of:
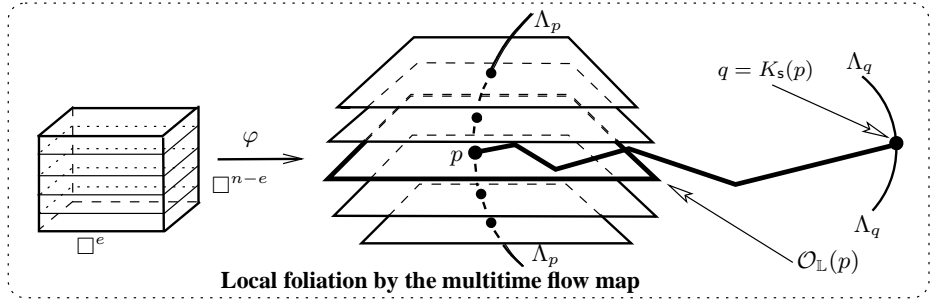
\[
T_e \times \Lambda_q \ni
({\sf s}, {\sf r}) \longmapsto 
K_{\sf s} ({\sf r}) \in M.
\]

\begin{center}
\input leaves.pstex_t
\end{center}

\subsection*{ 1.24.~Characterization of embedded $\LL$-orbits}
A smooth manifold $N$ together with an immersion $i : N \to M$ is
called {\sl weakly embedded} if for every manifold $P$, every smooth
map $\psi : P \to M$ with $\psi( P) \subset N$, then $\psi : P \to N$
is in fact smooth (\cite{ sp1970}; the diagram is also
borrowed).

\begin{center}
\begin{picture}(0,0)%
\includegraphics{eight.pstex}%
\end{picture}%
\setlength{\unitlength}{3947sp}%
\begingroup\makeatletter\ifx\SetFigFont\undefined
\def\x#1#2#3#4#5#6#7\relax{\def\x{#1#2#3#4#5#6}}%
\expandafter\x\fmtname xxxxxx\relax \def\y{splain}%
\ifx\x\y   
\gdef\SetFigFont#1#2#3{%
  \ifnum #1<17\tiny\else \ifnum #1<20\small\else
  \ifnum #1<24\normalsize\else \ifnum #1<29\large\else
  \ifnum #1<34\Large\else \ifnum #1<41\LARGE\else
     \huge\fi\fi\fi\fi\fi\fi
  \csname #3\endcsname}%
\else
\gdef\SetFigFont#1#2#3{\begingroup
  \count@#1\relax \ifnum 25<\count@\count@25\fi
  \def\x{\endgroup\@setsize\SetFigFont{#2pt}}%
  \expandafter\x
    \csname \romannumeral\the\count@ pt\expandafter\endcsname
    \csname @\romannumeral\the\count@ pt\endcsname
  \csname #3\endcsname}%
\fi
\fi\endgroup
\begin{picture}(5799,1674)(1189,-2473)
\put(1800,-2387){\makebox(0,0)[lb]{\smash{\SetFigFont{10}{12.0}{rm}{\color[rgb]{0,0,0}{\bf An immersion of the real line in $\R^ 2$ that is not weakly embedded}}%
}}}
\end{picture}

\end{center}

\def\theproposition{1.25}\begin{proposition}
{\rm (\cite{ bel1996, bm1997, bch2005})}
Each $\LL$-orbit is countable at infinity {\rm (}second countable{\rm
)} and weakly embedded in $M$.
\end{proposition}

As the multiple flows are diffeomorphisms, embeddability 
propagates. 

\def\theproposition{1.26}\begin{proposition}
{\rm (\cite{ bel1996, bm1997, bch2005})} 
Let $\mathcal{ O}_\LL$ be an $\LL$-orbit in $M$ and let $e:= \dim
\mathcal{ O}_\LL$. The following three conditions are equivalent{\rm :}

\smallskip

\begin{itemize}
\item[$\bullet$]
$\mathcal{ O}_\LL$ is an embedded submanifold of $M${\rm ;}

\smallskip

\item[$\bullet$]
for {\rm every} point $p\in \mathcal{ O}_\LL$, there exists a
straightening map $\varphi$ as in~\thetag{ 1.23} with $\mathcal{
O}_\LL \cap \varphi \left( \square^e \times \square^{ n- e} \right) =
\varphi \left( \square^e \times \{ 0 \} \right)${\rm ;}

\smallskip

\item[$\bullet$]
there {\rm exists at least one point} at which the
preceding property holds.
\end{itemize}

Conversely, $\mathcal{ O}_\LL$ is {\rm not embedded} in $M$ if and
only if for every $p\in \mathcal{ O}_\LL$ and for every local
straightening map $\varphi$ centered at $p$ as in~\thetag{ 1.23}, the
set of ${\sf r} \in \square^{ n - e}$ such that $\varphi \big(
\square^{ n-e} \times \{ {\sf r} \} \big)$ is contained in $\mathcal{
O}_\LL = \mathcal{ O}_\LL (p)$ is infinite {\rm (}nonetheless
countable{\rm )}.

\end{proposition}

\subsection*{1.27.~Local $\LL$-orbits and their smoothness} 
For $U$ running in the collection of all nonempty open connected
subsets of $M$ containing $p$, consider the localized $\LL
\vert_U$-orbit of $p$ in $U$, denoted by $\mathcal{O }_\LL(U, p)$. If
$p \in U_2 \subset U_1$, then $\mathcal{O }_\LL( U_2,p) \subset
\mathcal{ O }_\LL (U_1, p)\cap U_2$, so the dimension of $\mathcal{
O}_\LL(U, p)$ decreases as $U$ shrinks. Consequently, the localized
$\LL$-orbit $\mathcal{ O}_\LL (U,p)$ stabilizes and defines a unique
piece of local\footnote{ In certain references, local $\LL$-orbits are
considered as germs. Knowing by experience that the language of germs
becomes misleading when several quantifiers are involved in complex
statements, we will always prefer to speak of local submanifolds of a
certain small size.} $\LL$-integral submanifold through $p_0$, called
the {\sl local $\LL$ orbit of $p_0$} and denoted by $\mathcal{
O}_\LL^{ loc}(p)$. In the CR context, this concept will be of interest
in Parts~V and VI. Sometimes, $\LL$-orbits (in $M$) are called
{\sl global}, to distinguish them and to emphasize their nonlocal,
nonpointwise nature.

From the flow regularity Theorem~1.4 and from Theorem~1.21, it follows:

\def\thelemma{1.28}\begin{lemma}
Global and local $\LL$-orbits are as smooth as $\LL$, {\it i.e.}
$\mathcal{ C }^\omega$, $\mathcal{ C }^\infty$ or $\mathcal{ C }^{
\kappa, \alpha}$. Furthermore, 
\[
T_p\mathcal{O}_\LL^{loc}(p)
\subset
T_p\mathcal{O}_\LL(M,p)
=
\LL^{\rm inv}(p), 
\]
for every $p\in M$. This inclusion may be strict in the smooth
categories $\mathcal{ C}^\infty$ and $\mathcal{ C}^{\kappa, \alpha}$,
whereas, in the $\mathcal{ C}^\omega$ category, local and global CR
orbits have the same dimension.
\end{lemma}

In the $\mathcal{ C}^{ \kappa, \alpha}$ category, the maximal integral
curve of an arbitrary element of $\LL$ is $\mathcal{ C}^{ \kappa + 1,
\alpha}$, trivially because the right hand sides of the equations $d
{\sf x}_k ({\sf t}) / d {\sf t} = a_k ( {\sf x} ({\sf t}))$, $k=1,
\dots, n$, are $\mathcal{ C }^{ \kappa, \alpha}$. May it be induced
that general $\LL$-orbits are $\mathcal{ C }^{ \kappa + 1, \alpha}$?
Trivially yes if $\dim \LL = 1$ at every point.

Another instance is as follows. Let $r \in \N$ with $1\leqslant r
\leqslant n-1$ and let $\LL^0 = \{ L_a \}_{ 1\leqslant a \leqslant r}$
be a system of $\mathcal{ C}^{ \kappa, \alpha}$ vector fields defined
in a neighborhood of the origin in $\R^n$ that are linearly
independent there. Consider the system $\LL$ generated by linear
combinations of elements of $\LL^0$. Achieving Gaussian elimination
and a linear change of coordinates, we may assume that $r$ generators
of $\LL$, still denoted by $L_1, \dots, L_r$, take the form $L_i =
\frac{ \partial}{ \partial {\sf x}_i} + \sum_{ j= 1}^{ n -r} \, a_{
ij} ({\sf x}, {\sf y})\, \frac{ \partial }{\partial {\sf y}_j}$, $i=1,
\dots, r$, with $( {\sf x}, {\sf y} ) = ( {\sf x}_1, \dots, {\sf x}_r,
{\sf y}_1, \dots, {\sf y}_{ n-r})$ and with $a_{ij}({\sf x}, {\sf y})$
of class $\mathcal{ C}^{ \kappa, \alpha}$ in a neighborhood of the
origin.
 
We claim that {\it if $\mathcal{ O}_\LL (0)$ has {\rm (minimal
possible)} dimension $r$, then it is $\mathcal{ C}^{ \kappa+1,
\alpha}$}. This happens in particular if $\LL$ is
Frobenius-integrable.

Indeed, the local graphed equations of $\mathcal{ O}_\LL (0)$ must
then be of the form ${\sf y}_j = h_j ({\sf x})$, $j = 1, \dots, n-r$,
with the $h_j$ of class at least $\mathcal{ C}^{ \kappa, \alpha}$,
thanks to the lemma above. Observe that the $L_i$ are tangent to this
submanifold if and only if the $h_j$ satisfy the complete system of
partial differential equations $\frac{ \partial h_j}{ \partial {\sf
x}_i} ({\sf x}) = a_{ij} ({\sf x}, h( {\sf x}))$, for $i=1, \dots, r$,
$j=1, \dots, n-r$, implying directly that the $h_j$ are $\mathcal{
C}^{ \kappa+1, \alpha}$. In general, this argument shows that if $\dim
\mathcal{ O}_\LL (p)$ coincides with $\dim \LL (p)$, the orbit is
$\mathcal{ C}^{ \kappa + 1, \alpha}$ at $p$.

\def\theexample{1.29}\begin{example}{\rm
However, this improvement of smoothness is untrue when $\dim \LL (p) +
1 \leqslant \dim \mathcal{ O}_\LL (p) \leqslant n-1$.

Indeed, pick the function $\chi_{\kappa, \alpha} = \chi_{\kappa,
\alpha} ({\sf z})$ of ${\sf z} \in \R$ equal to zero for ${\sf z}
\leqslant 0$ and, for ${\sf z}\geqslant 0$, defined by:
\[
\chi_{\kappa, \alpha}({\sf z})
=
\left\{
\aligned
&
{\sf z}^{\kappa+\alpha}, 
\ \ \ \ \ \ \ \ \ \ \
{\rm if}
\ \
0<\alpha\leqslant 1,
\\
&
{\sf z}^\kappa/{\rm log}\, {\sf z},
\ \ \ \ \ \ \,
{\rm if}
\ \
\alpha=0.
\endaligned
\right.
\]
This function
is $\mathcal{ C}^{\kappa, \alpha}$ on $\R$, but for
$(\lambda, \beta) > (\kappa, \alpha)$,
it is not $\mathcal{ C}^{ \lambda, \beta}$ 
in any neighborhood of
the origin.
Then
in $\R^4$
equipped with coordinates $( {\sf x}, {\sf y}, {\sf z}, {\sf t})$,
consider the hypersurface $\Sigma$ of
equation:
\[
0
=
{\sf t}
-
\chi_{\kappa+1,\alpha}({\sf y})
\,
\chi_{\kappa,\alpha}({\sf z}),
\] 
Then $\Sigma$ is $\mathcal{ C}^{ \kappa, \alpha }$, not better.
The two vector fields $L_1:= \frac{ \partial
}{\partial {\sf x} }$ and $L_2:= \frac{ \partial }{\partial {\sf y}} +
\left[ {\sf x} \, \chi_{ \kappa, \alpha} (-{\sf y}) \right]
\frac{ \partial }{ \partial {\sf z}} + \left[ \chi_{ \kappa + 1,
\alpha}' ({\sf y})\, \chi_{ \kappa, \alpha} ({\sf z}) \right]
\frac{\partial }{\partial {\sf t}}$ have $\mathcal{ C}^{ \kappa,
\alpha}$ coefficients and are tangent to $\Sigma$. We claim
that $\Sigma$ is the local $\{ L_1, L_2\}$-orbit of the origin. 

Otherwise, there would exist a local two-dimensional submanifold $
\big\{ {\sf z} = g({\sf x}, {\sf y}), \ {\sf t} = h ({ \sf x}, {\sf
y}) \big\}$ with $L_1$ and $L_2$ tangent to it. Then $[L_1, L_2] =
\chi_{ \kappa, \alpha} (-{\sf y}) \, \frac{ \partial}{ \partial {\sf
z}}$ should also be tangent. However, at points $\left( 0, {\sf y},
g(0, {\sf y}), h (0, {\sf y}) \right)$, with ${\sf y}$ negative and
arbitrarily small, $L_1$, $L_2$ and $L_3$ are equal to the three
linearly independent vectors $\frac{ \partial }{ \partial {\sf x}}$,
$\frac{ \partial }{ \partial {\sf y}}$ and $\chi_{ \kappa, \alpha}
(-{\sf y}) \, \frac{ \partial}{ \partial {\sf z }}$.
\qed

}\end{example}

\section*{ \S2.~Finite type system and their genericity 
(openess and density)}

\subsection*{2.1.~Systems of vector
fields satisfying $\LL^{ \rm lie} = \LL$} Let $M$ be a $\mathcal{
C}^\kappa $ ($1\leqslant \kappa \leqslant \infty$) connected manifold
of dimension $n\geqslant 1$. By $\X$, denote the system of all vector
fields defined on open subsets of $M$ (it is a sheaf). Let
\[
\LL^0 = \{
L_a \}_{ 1\leqslant a\leqslant r}, 
\ \ \ \ \ \ \ \
r\geqslant 1, 
\]
be a {\it finite}\, collection of $\mathcal{ C}^{ \kappa-1}$ vector
fields defined on $M$, namely $L_a \in \X (M)$. Unlike in the
$\mathcal{ C}^\omega$ category, in the $\mathcal{ C}^\kappa$ category,
$\X (M)$ is always nonempty and quite large, thanks to partitions of
unity. For this reason, we shall not work in the real analytic
category, except in some specific local situations. The set of linear
combinations of elements of $\LL^0$ with coefficients in $\mathcal{
C}^{ \kappa-1} (M, \R)$ will be denoted by $\LL$ (or $\LL^1$) and
called the $\mathcal{ C}^{ \kappa-1} (M)$-{\sl linear hull of}
$\LL^0$.

\def\thedefinition{2.2}\begin{definition}{\rm
A $\mathcal{ C}^{ \kappa -1} 
(M)$-linear system $\LL \subset \X$ is said to
be {\sl of finite type at a point $p\in M$} if $\LL^{\rm lie} (p) =
T_p M$.
}\end{definition}

If $\LL^{ \rm lie}$ is of finite type at every point, then $\LL^{ \rm
lie} = \LL^{ \rm inv} = \X$ and there is just one maximal
$\LL$-integral manifold in the sense of Sussmann: $M$ itself.

In 1939, Chow had already shown that the equality $\LL^{ \rm lie} =
\X$ implies the {\sl everywhere accessibility condition}: every two
points of $M$ may be joined by integral curves of $\LL$. In 1967,
H\"ormander established that every second order partial differential
operator $P := L_1^2 + \cdots + L_r^2 + R_1 + R_0$ on a domain $\Omega
\subset \R^n$ whose top order part is a sum of squares of $\mathcal{
C}^\infty$ vector fields $L_a$, $1\leqslant a \leqslant r$, such that
$\LL^{ \rm lie} = \X$ is $\mathcal{ C}^\infty$-hypoelliptic, namely $P
f \in \mathcal{ C}^\infty$ implies $f \in \mathcal{ C}^\infty$. Vector
field systems satisfying $\LL^{ \rm lie} = \X$ have been further
studied by workers in hypoelliptic partial differential equations and
in nilpotent Lie algebras: M\'etivier, Stein, Mitchell, Stefan, Lobry
and others.

In the next Parts~V and VI, we will focus on propagational
aspects that are enjoyed by the (more general) smooth systems $\LL$
that satisfy $\LL^{\rm inv} = \X$, but possibly $\LL^{ \rm lie} (p)
\neq \X (p)$ at every $p \in M$. Nevertheless, for completeness, we
shall survey in the present section some classical geometric
properties of finite type systems.

\subsection*{ 2.3.~Lie bracket flags, weights, privilegied
coordinates and distance estimate} Define $\LL^1 := \LL$ and by
induction, for $s\in \N$ with $2\leqslant s \leqslant \kappa$, define
$\LL^s$ to be the $\mathcal{ C}^{ \kappa - s}$-linear hull of $\LL^{
s-1} + \left[ \LL^1, \LL^{ s-1} \right]$. Concretely, $\LL^s$ is
generated over $\mathcal{ C}^{\kappa - s}$ by iterated Lie brackets of
length $\leqslant s$ of the form:
\[
L_\alpha
=
[L_{\alpha_1},[L_{\alpha_2},\dots,
[L_{\alpha_{\ell-1}},L_{\alpha_\ell}]\dots]],
\ \ \ \ \ \ \ 
1\leqslant\ell\leqslant s.
\]
Jacobi's identity insures that $[ \LL^{ s_1}, \LL^{ s_2} ] \subset
\LL^{ s_1 + s_2}$.

Denote $\LL^s (p) := {\rm Vect}_\R \, \{ L(p) : \, L \in \LL^s\}$.
Clearly, $\LL$ is of finite type at $p\in M$ if and only it there
exists an integer $d(p) \leqslant \kappa$ with $\LL^{ d(p)} (p) = T_p M$.
The smallest $d(p)$ is sometimes called the {\sl degree of
non-holonomy} of $\LL$ at $p$. Other authors call it the {\sl type of}
$\LL$ at $p$, which we will do. The function $p \mapsto d (p) \in [1,
\kappa] \cup \{ \infty\}$ is upper-semi-continuous: $d(q) \leqslant d(p)$
for $q$ near $p$.

Combinatorially, at a finite type point, it is of interest to
introduce the {\sl Lie bracket flag}:
\[
\{0\}
\subset
\LL^1(p)
\subset
\LL^2(p)
\subset
\cdots
\subset
\LL^s(p)
\subset
\cdots
\subset
\LL^{d(p)}(p)
=
T_pM.
\]
Then a finite type point $p$ is called {\sl regular} if the integers
$n_s(q) := \dim \LL^s (q)$ remain constant in some neighborhood of
$p$. It is elementary to verify (\cite{ bel1996}) that, at such a
regular point, the dimensions are strictly increasing:
\[
0
<
n_1(p)
<
n_2(p)
<
\cdots
<
n_{d(p)}(p)
=
n.
\]

Fix $p$, not necessarily regular. A local coordinate system $(x_1,
x_2, \dots, x_n)$ centered at $p$ is {\sl linearly adapted at $p$} if:
\[
\left\{
\aligned
\LL^1(p)
&
=
{\rm Vect}_p
\left(
\frac{\partial}{\partial x_1},
\dots,
\frac{\partial}{\partial x_{n_1(p)}}
\right),
\\
\LL^2(p)
&
=
{\rm Vect}_p
\left(
\frac{\partial}{\partial x_1},
\dots,
\frac{\partial}{\partial x_{n_1(p)}},
\dots,
\frac{\partial}{\partial x_{n_2(p)}}
\right),
\\
\dots\dots
&
\dots\dots\dots\dots
\dots\dots\dots\dots
\dots\dots\dots\dots
\dots\dots\dots
\\
\LL^{d(p)}(p)
&
=
{\rm Vect}_p
\left(
\frac{\partial}{\partial x_1},
\dots,
\frac{\partial}{\partial x_{n_1(p)}},
\dots,
\frac{\partial}{\partial x_{n_2(p)}},
\dots
\frac{\partial}{\partial x_{n_{d(p)}(p)}}
\right).
\endaligned\right.
\]
Let us assign {\sl weights} $w_i$ to such linearly adapted coordinates
$x_i$ as follows: the first group $(x_1, \dots, x_{ n_1 (p)})$ being
linked to $\LL^1 (p)$, their weights will all be equal to one: $w_1 =
\dots = w_{ n_1 (p)} =1$. The second group $(x_{ n_1 (p)+ 1}, \dots,
x_{ n_2 (p)})$, linked to the quotient $\LL^2 (p) / \LL^1 (p)$, will
be assigned uniform weight two: $w_{ n_1 (p) + 1} = \dots = w_{ n_2
(p)} = 2$, and so on, until $w_{ n_{ d(p)-1}(p) +1} = \dots = w_{
n_{d(p)}(p)} = d(p)$.

Provided $\LL$ is of finite type at every point, we claim that the
original finite collection $\LL^0$ produces what is called a {\sl
sub-Riemannian metric}; then by means of weights, the topology
associated to this metric may be compared to the manifold topology of
$M$ in a highly precise way.

Indeed, let us define the (infinitesimal) {\sl sub-Riemannian length}
of a vector $v_p \in \LL^1 (p)$ by:
\[
\vert\!\vert
v_p
\vert\!\vert_{\LL^0}
:=
\inf
\big\{
(u_1^2+\cdots+u_m^2)^{1/2}:\,
v_p
=
u_1\,L_1(p)
+
\cdots
+
u_r\,L_r(p)
\big\}.
\]
For $v_p \not \in \LL^1 (p)$, we set $\vert\!\vert v_p
\vert\!\vert_{\LL^0} = \infty$. The length of a piecewise $\mathcal{
C}^1$ curve $\gamma (t)$, $t\in [ 0, 1]$, will be the integral:
\[
{\rm length}_{\LL^0}(\gamma)
:=
\int_0^1\,
\vert\!\vert
d\gamma(t)/dt
\vert\!\vert_{\LL^0}\,
dt.
\]
Finally, the distance associated to the finite collection $\LL^0$ is:
\[
d_{\LL^0}(p,q)
:=
\inf_{\gamma:\ p\to q}\,
{\rm length}_{\LL^0}(\gamma).
\]
Assume for instance $d(p) = 2$, so that $n_2(p) = n$. If the
coordinates are linearly adapted, the tangent space $T_p M$ then
splits in the ``horizontal'' space, the $(x_1, \dots, x_{
n_1(p)})$-plane, together with a (not unique) 
``vertical'' space generated {\it e.g.} by the
remaining coordinates. It is then classical that the distance from $p$
to a point of coordinates $(x_1, \dots, x_n)$ close to $p$ enjoys the
estimate:
\[
d_{\LL^0}
\big(
p,(x_1,\dots,x_n)
\big)
\asymp
\vert x_1\vert
+
\cdots
+
\vert x_{n_1(p)}\vert
+
\vert x_{n_1(p)+1}\vert^{1/2}
+
\cdots
+
\vert x_n\vert^{1/2}.
\]
Here, the abbreviation $\Phi \asymp \Psi$ means that there exists $C>1$
with $C^{ -1} \, \Psi < \Phi < C \, \Psi$. Notice that the successive
exponents coincide with the weights $w_1, \dots, w_{ n_1(p)}, 
w_{ n_1 (p) +1},
\dots,
w_n$. In particular, to reach a point of coordinates $(0, \dots, 0,
\varepsilon, \dots, \varepsilon)$, it is necessary to flow along
$\LL^0$ during a time $\sim {\rm cst.}\, 
\varepsilon^{ 1/2}$. Observe that
$\vert x_1 \vert + \cdots + \vert x_n\vert$ is equivalent to the
distance from $p$ to $x$ induced by any Riemannian metric. Thus,
the modified distance $d_{ \LL^0}$ is just obtained
by replacing each $\vert x_i\vert$ by $\vert x_i\vert^{ 1/w_i}$, 
up to a multiplicative constant. 

To generalize such a quantitative comparison between the
$d_{\LL^0}$-distance and the underlying topology of $M$, linearly
adapted coordinates appear to be insufficient. For $\beta = (\beta_1,
\dots, \beta_r) \in \N^r$, denote by $L^\beta$ the $\vert \beta
\vert$-th order derivation $L_1^{\beta_1} L_2^{ \beta_2} \cdots
L_r^{\beta_r}$. Beyond linearly adapted coordinates, one must introduce
{\sl privileged coordinates}, whose existence is assured by
the following.

\def\thetheorem{2.4}\begin{theorem}
{\rm (\cite{ bel1996})} There exist local
coordinates $(x_1, \dots,
x_n)$ centered at $p$ that
are {\rm privileged} in the
sense that each $x_i$ is of order exactly equal
to $w_i$ with respect to $\LL^0$-derivations, namely, for
$i=1, \dots, n${\rm :}
\[
\aligned
L^\gamma\,x_i\vert_p
&
=
0, 
\ \ \ \ \ 
\text{\rm for all}\
\gamma\
\text{\rm with}\
\vert\gamma\vert
\leqslant
w_i-1,
\\
L^{\beta_i^*}x_i\vert_p
&
\neq 0, 
\ \ \ \ \ 
\text{\rm for some}\
\beta_i^*\
\text{\rm with}\
\vert\beta_i^*\vert=w_i.
\endaligned
\]
\end{theorem}

Only if $d(p) = 2$, linearly adapted 
coordinates are automatically privileged (\cite{
bel1996}). As soon as
$d(p) \geqslant 3$, privileged systems are unavoidable.

\def\thetheorem{2.5}\begin{theorem}
{\rm (\cite{ bel1996})}
For $x$ in a neighborhood of $p$, the estimate{\rm :}
\[
d_{\LL^0}(p,(x_1,\dots,x_n))
\asymp
\vert x_1\vert^{1/w_1}
+
\cdots
+
\vert x_{n}\vert^{1/w_n}
\]
holds if and only if the coordinates are privileged.
\end{theorem}

For $\varepsilon >0$ small, define the anisotropic ball $B_{\LL^0} (p,
\varepsilon):= \{ x : \, d_{ \LL^0} (p, x) < \varepsilon \}$.

\def\thecorollary{2.6}\begin{corollary}
{\rm (\cite{ bel1996})}
There exist $C>1$ such that
\[
\frac{1}{C}\,
\prod_{i=1}^n\,
[-\varepsilon^{w_i},\varepsilon^{w_i}]
\subset
B_{\LL^0}(p,\varepsilon)
\subset
C\,
\prod_{i=1}^n\,
[-\varepsilon^{w_i},\varepsilon^{w_i}].
\]
\end{corollary}

\subsection*{ 2.7.~Local basis}
At a non-regular point, the integers $n_k(p)$, $k=1, \dots, d(p)$ are
not necessarily strictly increasing. Thus, it is necessary to express
the combinatorics of the Lie bracket flag with more precision, in
terms of what is sometimes called {\sl H\"ormander numbers} $m_i$,
$\ell_i$. From now on, we shall assume that the $r$ vector fields
$L_a$, $1\leqslant a \leqslant r$, are linearly independent at $p$ and
have $\mathcal{ C}^\infty$ or $\mathcal{ C}^\omega$ coefficients. In
both cases, the formal Taylor series of every coefficient exists.

In the flag
\[
\{0\}
\subset
\LL^1(p)
\subset
\LL^2(p)
\subset
\cdots
\subset
\LL^s(p)
\subset
\cdots
\subset
\LL^{d(p)}(p)
=
T_pM.
\]
let $m_1$ denote the smallest $k \geqslant 2$
such that the dimension of $\LL^k (p)$ is larger than the dimension of
$\LL^1 (p)$ (at a regular point, $m_1 = 2$) and set $\ell_1 := \dim
\LL^{ m_1} (p) - r \geqslant 1$. Similarly, let $m_2$ denote the smallest
$k \geqslant 1+ m_1$ such that the dimension of $\LL^k (p)$ is larger than
the dimension of $\LL^{ m_1} (p)$ (at a regular point, $m_2 = 3$) and
set $\ell_2 := \dim \LL^{ m_2} (p) - \dim \LL^{ m_1} (p)$. By
induction, let $m_{ j+1}$ denote the smallest $k \geqslant 1+ m_j$ such
that the dimension of $\LL^k (p)$ is larger than the dimension of
$\LL^{ m_j} (p)$ and set $\ell_{ j+1} := \dim \LL^{ m_{ j+1}} (p) -
\dim \LL^{ m_j} (p)$.

Since $p$ is a point of finite type, the process terminates until $m_h
= d(p)$ reaches the degree of non-holonomy at $p$, for a certain
integer $h\geqslant 1$. We thus have extracted the
interesting information, namely the strict flag
of linear spaces:
\[
\LL^1(p)
\subset
\LL^{m_1}(p)
\subset
\LL^{m_2}(p)
\subset
\cdots
\subset
\LL^{m_h}(p)
=
T_pM,
\]
with Lie bracket orders $1 < m_1 < m_2 < \cdots < m_h$, whose
successive dimensions may be listed parallelly:
\[
r 
< 
r+\ell_1 
<
r+\ell_1+\ell_2
< 
\cdots
<
r+\sum_{1\leqslant j\leqslant h}\,\ell_j.
\]

Next, let $x = (x_1, \dots, x_n)$ be linearly adapted coordinates,
vanishing at $p$. We shall denote them by $(y, s_1, s_2, \dots,
s_h)$, where $y\in \R^r$, $s_1\in \R^{ \ell_1}$, $s_2 \in \R^{
\ell_2}$, \ldots, $s_h \in \R^{ \ell_h}$. As in the preceding
paragraph, we assign weight $1$ to the $y$-coordinates, weight $m_1$
to the $s_1$-coordinates, weight $m_2$ to the $s_2$-coordinates,
\ldots, weight $m_h$ to the $s_h$-coordinates. The weight of a
monomial $x^\alpha = y^\beta \, s_1^{\gamma_1} \, s_2^{ \gamma_2}
\cdots s_h^{ \gamma_h}$ is obviously defined as $\vert \beta \vert +
m_1 \, \vert \gamma_1 \vert + m_2 \, \vert \gamma_2 \vert + \cdots +
m_h \, \vert \gamma_h \vert$. We say that a formal power series $a
(x) = a (y, s_1, \dots, s_h)$ is an ${\rm O} (\kappa)$ if all its
monomials have weight $\geqslant \kappa$. Also, $a(x)$ is called {\sl
weighted homogeneous of degree $\kappa$} if
\[
a
\left(
t y, 
t^{m_1}s_1, 
t^{m_2}s_2,
\dots,
t^{m_h}s_h
\right)
=
t^\kappa
a(y,s_1,s_2,\dots,s_h), 
\]
for all $t\in \R$. As in the case of $\R \dl z_1, \dots, z_n \dr$ with
all weights equal to $1$, every formal series $a (y, s_1, s_2, \dots,
s_h)$ may be decomposed as a countable sum of weighted homogeneous
polynomials of increasing degree.

Dually, we also assign weights to all the basic vector fields: $\frac{
\partial }{\partial y_a}$ will have weight $-1$, whereas for $j = 1,
\dots, m_h$, the $\frac{ \partial }{\partial s_{ jl}}$, $l = 1, \dots,
\ell_j$, will have weight $-m_j$. The weight of a monomial vector
field $x^\alpha \, \frac{ \partial }{\partial x_i}$ is defined to be
the sum the weights of $x^\alpha$
with the weight of $\frac{ \partial }{\partial
x_i}$. Every vector field having formal power series coefficients may
be decomposed as a countable sum of weighted homogeneous vector fields
having polynomial coefficients.

\def\thetheorem{2.8}\begin{theorem}
{\rm (\cite{ bel1996, ber1999})}
Assume the local vector fields $L_a$, $a = 1, \dots, r$, have
$\mathcal{ C}^\infty$ or $\mathcal{ C}^\omega$ coefficients and are
linearly independend at $p$. If the $\mathcal{ C}^\infty$ or
$\mathcal{ C}^\omega$ coordinates $x = (y, s_1, s_2, \dots, s_h)$
centered at $p$ are {\rm priveleged}, then each $L_a$ may be developed
as{\rm :}
\[
L_a
= 
\widehat{L}_a
+
{\rm O}(0),
\]
where each vector field{\rm :} 
\[
\widehat{L}_a
:=
\frac{\partial }{\partial y_a}
+
\sum_{1\leqslant j\leqslant 1}\,\sum_{1\leqslant
l\leqslant \ell_j}\, 
p_{a,j,l}
(y,s_1,\dots,s_{j-1})\,
\frac{\partial}{\partial s_{j,l}},
\]
is homogeneous of degree $-1$ and has as its coefficients some
polynomials $p_{ a, j, l} = p_{ a, j, l} (y, s_1, \dots, s_{ j-1})$
that are independent of $s_j$ and are homogeneous of degree $m_j -1$.
\end{theorem}

A crucial algebraic information is missing in this statement: what are
the nondegeneracy conditions on the $p_{ a, j, l}$ that insure that
the system is indeed of finite type at $p$ with the combinatorial
invariants $m_j$ and $\ell_j$\,? The real problem is to classify
vector field systems that are of finite type, up to local changes of
coordinates. At least, the following may be verified.

\def\thetheorem{2.9}\begin{theorem} 
{\rm (\cite{ bel1996, ber1999})}
The vector fields $\widehat{ L}_a$, $a= 1, \dots, r$, form a finite
type system $\widehat{ \LL}^0$ at $p$ having the {\rm same}
combinatorial invariants $m_j$ and $\ell_j$ and satisfying the same
distance estimate as $d_{\LL^0}$ in Theorem~2.5. Moreover, the linear
hull of $\widehat{ \LL}^0$ generates a Lie algebra $\widehat{
\LL}^{\rm lie}$ with the nilpotency property that all Lie brackets of
length $\geqslant m_h + 1$ all vanish.
\end{theorem}

\subsection*{ 2.10.~Finite-typisation of smooth systems 
of vector fields}
As previously, let $\LL^0 = \{ L_a \}_{ 1\leqslant a \leqslant r}$ be
a finite collection of $\mathcal{ C}^{ \kappa -1}$ vector fields
globally defined on a connected manifold $M$ of class $\mathcal{
C}^\kappa$ ($1\leqslant \kappa \leqslant \infty$) and of dimension
$n\geqslant 1$. Let $\LL$ be its $\mathcal{ C}^{ \kappa-1} (M)$-linear
hull. If $r = 1$, then $\LL^{ \rm lie} = \LL$, hence $\LL$ cannot be
of finite type, unless $n=1$. So we assume $n \geqslant 2$ and $r
\geqslant 2$. We want to perturb $\LL$ slightly to $\widetilde{ \LL}$
so as to get finite-typeness at every point: $\widetilde{ \LL}^{\rm
lie} (p) = T_p M$ at every $p\in M$. 
Since the composition of Lie brackets of 
length $\ell$ requires coefficients of
vector fields to be at least 
$\mathcal{ C}^\ell$, if $\kappa < \infty$, then necessarily
$\widetilde{ \LL}^{ \rm lie} = \widetilde{ \LL}^\kappa$ stops
at length $\kappa$.

At a central point, say the origin in $\K^n$, and for $\K$-analytic
vector field systems, the already presented Theorem~1.11(II) yields
small perturbations that are of finite type at $0$. Of course, the same
local result holds true for collections of vector fields that are
$\mathcal{ C}^\infty$, or even $\mathcal{ C}^{ \kappa -1}$ with
$\kappa$ large enough. Now, we want a global theorem.

What does it mean for $\widetilde{ \LL}$ to be close to $\LL$\,? A
vector field $L \in \X (M)$ may be interpreted as a section of the
tangent bundle, in particular a $\mathcal{ C}^{ \kappa-1}$ map $M \to
TM$. The most useful topology on the set $\mathcal{ C}^\lambda (M, N)$
of all $\mathcal{ C}^\lambda$ maps from a manifold $M$ to another
manifold $N$ ({\it e.g.} $N= TM$ with $\lambda = \kappa -1$) is the {\sl
strong Whitney topology}; it controls better than the
so-called {\sl weak topology}\,
the behaviour of maps at infinity in the noncompact case. Essentially,
$f, g \in \mathcal{ C}^\lambda (M, N)$ are
(strongly) close to each other if all
their partial derivatives of order $\leqslant \lambda$, computed in a
countable collection of charts $\varphi_\nu : U_\nu \to \R^n$ and
$\psi_\nu : V_\nu \to \R^m$ covering $M$ and $N$, $\nu \in \N$, are
$\varepsilon_\nu$-close, the smallness of $\varepsilon_\nu >0$
depending on the pair of charts $(\varphi_\nu, \psi_\nu)$. Precise
definitions may be found in the monograph~\cite{ hi1976}. We then
topologize this way the finite product 
$\X (M)^r$.

Already two vector fields may well be of finite type on a manifold of
arbitrary dimension, {\it e.g.} $\frac{ \partial }{\partial x_1}$ and
$\sum_{ i=2}^n \, x_1^{i-1}\, \frac{ \partial }{\partial x_i}$ on $\R^n$.

\def\thetheorem{2.11}\begin{theorem}
{\rm (\cite{ lo1970})} If the connected manifold $M$ of dimension
$n\geqslant 2$ is $\mathcal{ C}^{n + n^2}$, 
then the set of pairs of vector
fields $\LL^0 := ( K, L) \in \X(M)^2$ on $M$ whose $\mathcal{ C}^{ n^2
+ n -1}$-linear hull $\LL$ satisfies $\LL^{ n^2 + n} = \LL$, is open
and dense in the strong Whitney topology.
\end{theorem}

According to~\cite{ su1976}, the smoothness $M \in \mathcal{ C}^{ n +
n^2}$ in~\cite{ lo1970} was improved to $M \in \mathcal{ C}^{ 2n}$ in
Lobry's thesis (unpublished). We will summarize the demonstration in
the case $M \in \mathcal{ C}^{ 2n}$. However, since neither
$\mathcal{ C}^{2n}$ nor $\mathcal{ C}^{ n+ n^2}$ are optimal, we will
improve this result afterwards (Theorem~2.16 below).

\proof
Openness is no mystery. For denseness, we need some preliminary. If
$M$ and $N$ are two $\mathcal{ C}^\lambda$ manifolds, we denote by
$J^\lambda (M, N)$ the {\sl bundle of $\lambda$-th jets of $\mathcal {
C}^\lambda$ maps from $M$ to $N$}. We recall that, to a $\mathcal{
C}^\lambda$ map $f : M \to N$ is associated the $\lambda$-th jet map
$j^\lambda f : M \to J^\lambda (M, N)$, a continuous map that may be
considered as a kind of intrinsic collection of all partial
derivatives of $f$ up to order $\lambda$. Let $\pi : J^\lambda ( M, N)
\to M$ be the canonical projection, sending a jet to its base point.
For $p\in M$, the fiber $\pi^{ -1} (p)$ may be identified with $\R^{
N_{ m,n, \lambda}}$, where $N_{ m, n, \lambda} := m \, \frac{ (n+
\lambda) !}{ n! \ \lambda !}$ counts the number of partial derivatives
of order $\leqslant \lambda$ of maps $\R^n \to \R^m$. 

We will state a lemma which constitutes a special case of the jet
transversality theorem. This particular statement (Lemma~2.12 below)
generalizes the intuitively obvious statement that any $\mathcal{ C
}^0$ curve graphed over $\R \times \{ 0\}^2$ in $\R^3$ may always be
slightly perturbed to avoid a given fixed $\mathcal{ C }^1$ curve
$\Sigma$.

Call a subset $\Sigma \subset J^\lambda (M, N)$ {\sl algebraic in the
jet variables} if in {\it every}\, pair of local charts, it possesses
defining equations that are polynomials in the jet variables $f_{j,
\alpha}$, $1\leqslant j \leqslant m$, $\alpha \in \N^n$, $\vert \alpha
\vert \leqslant \lambda$, whose coefficients are {\it independent of
the coordinates $x \in M$}. Of course, after a local diffeomorphism
$x \mapsto \bar x ( x)$ of $M$, a general polynomial in the jet
variables which is independent of $x$ almost never remains independent
of $\bar x$ in the new coordinates. Nevertheless, in the sequel, we
shall only encounter special sets $\Sigma \subset J^\lambda (M, N) $
which, in any coordinate system, may be defined as zero sets of such
special polynomials.

For instance, taking $(x_k, y_j, y_{ j,k})$ as coordinates on $J^1 (M,
N)$, where $1\leqslant k \leqslant n = \dim M$ and $1 \leqslant j
\leqslant m = \dim N$, a change of coordinates $x \mapsto \bar x (x)$
induces $(x_k, y_j, y_{ j,k}) \longmapsto (\bar x_k, \bar y_j, \bar
y_{ j,k})$, where $\bar y_j = y_j$ is unchanged but the new jet
variables $y_{ j,k} = \sum_{ l= 1}^n \, \bar y_{ j,l} \,
\frac{\partial \bar x_l}{\partial x_k}$ involve the variables $x$ (or
$\bar x$). Nevertheless, the equations $\{ y_{ j,k} = 0, \,
1\leqslant j \leqslant m, \ 1 \leqslant k \leqslant n \}$ saying that
the first (pure) jet vanishes are equivalent to $\{ \bar y_{ j,k} = 0,
\, 1\leqslant j \leqslant m, \ 1 \leqslant k \leqslant n \}$, since
the invertible Jacobian matrix $\big( \frac{\partial \bar
x_l}{\partial x_k} \big)$ may be erased: vanishing properties in a jet
bundle are intrinsic~!

A theorem due to Whitney states that real algebraic sets are
stratified, {\it i.e.} are finite unions of geometrically smooth real
algebraic manifolds. The {\sl codimension} of $\Sigma$ is thus
well-defined.

\def\thelemma{2.12}\begin{lemma}
{\rm (\cite{ hi1976})} Assume $\Sigma \subset J^\lambda ( M, N)$ is
algebraic in the jet variables and of codimension $\geqslant 1 + \dim
M$. Then the set of maps $f\in \mathcal{ C}^\lambda (M, N)$ whose
$\lambda$-th prolongation $j^\lambda f : M \to J^\lambda (M, N)$ does not
meet $\Sigma$ at any point is open and dense in the strong Whitney
topology.
\end{lemma}

Although $j^\lambda f$ is only continuous, the fact that the bad set
$\Sigma$ is algebraic enables to apply the appropriate
version of Sard's theorem that is used in the jet
transversality theorem.

We shall apply the lemma by defining
a certain bad set $\Sigma$ which, if avoided, 
means that a pair of vector fields on 
$M$ is of finite type at every point. 

Assume $M \in \mathcal{ C}^{ 2n}$ and let
$(K, L) \in \X(M )^2$. Both vector fields have $\mathcal{ C}^{ 2 n
-1}$ coefficients. With $\lambda := 2n -1$, denote by $J^{2n-1} (\X
(M)^2)$ the fiber bundle of the $(2n-1)$-th jets of these pairs. In
some coordinates provided by a local chart $U \ni q \mapsto (x^1 (q),
\dots, x^n (q)) \in \R^n$, with $U \subset M$ open, we may write $K =
\sum_{ 1\leqslant i\leqslant n} \, K_i (x) \, \frac{ \partial }{\partial x^i}$
and $L = \sum_{ 1\leqslant i\leqslant n} \, L_i (x) \, \frac{ \partial }{
\partial x^i}$. In such a chart, the $(2n-1)$-th jet map $j^{ 2n-1}
(K, L) : U \longrightarrow J^{ 2n-1} ( \X^2(M) \vert_U)$ is concretely
given by:
\[
U\ni x
\longmapsto
\big(
\partial_x^\alpha\,K_i(x),\
\partial_x^\alpha\,L_i(x)
\big)_{\alpha\in\N^n,\
\vert\alpha\vert\leqslant 2n-1,\ 
1\leqslant i\leqslant n}.
\]
We denote by $K_{ i, \alpha}$ and $L_{ i, \alpha}$ the corresponding
jet variables. A $\mathcal{ C}^{ 2n}$
local diffeomorphism $x \mapsto \bar x = \bar x(x)$
induces a triangular transformation involving the chain rule between
these jets variables, with coefficients depending on the $2n$-th
jet of $\bar x (x)$, some of which are only $\mathcal{ C}^0$, which
might be unpleasant. Fortunately, our bad set $\Sigma$ will be shown
to be algebraic with respect to 
the jet variables $K_{ i, \alpha}$ and
$L_{ i, \alpha}$ {\it in any system of coordinates}.

Let $(K, L) \in \X (M)^2$. To write shortly iterated Lie brackets, we
denote ${\rm ad} (K) L := [ K, L]$, so that ${\rm ad} (K)^2 L = [ K, [
K, L]]$, ${\rm ad} (K)^3 L = [K, [ K, [ K, L]]]$ and so on. Also, we
set ${\rm ad}^0 (K) L := L$. Define a subset $\Sigma \subset J^{ 2n-1}
( \X^2 (M))$ as a union $\Sigma = \Sigma ' \cup \Sigma '' \cup
\Sigma'''$, where:

\begin{itemize}

\smallskip\item[$\bullet$]
firstly $\Sigma'$ is defined by the $2n$ equations $K_{i,0} = L_{i, 0} 
= 0$;

\smallskip\item[$\bullet$]
secondly, $\Sigma''$ is defined by requiring that all the $n\times n$
minors of the following $n \times (2n)$ matrix
\[
\left(
{\rm ad}^0(K)L\ \
{\rm ad}^1(K)L\ \ 
\cdots\cdots\
{\rm ad}^{2n-1}(K)L
\right),
\]
vanish;

\smallskip\item[$\bullet$]
thirdly, $\Sigma'''$ is defined similarly, after exchanging $K$ with
$L$.

\end{itemize}\smallskip

\def\thelemma{2.13}\begin{lemma}
In the vector space of real $n\times (2n)$ matrices, isomorphic to
$\R^{ 2n^2}$, the subset of matrices of rank $\leqslant (n-1)$ is a real
algebraic subset of codimension equal to $(n+1)$.
\end{lemma}

Without obtaining a complete explicit expression, it is easily
verified that ${\rm ad}^j (K) (L)$, $0 \leqslant j \leqslant 2n -1$, is a
universal polynomial in the jet variables $K_{ i, \alpha }$ and $L_{
i, \alpha}$. Under a local change of coordinates $x \mapsto \bar x
(x)$, if the two vector fields $K$ and $L$ transform to $\overline{
K}$ and to $\overline{ L}$ (push-forward), all the multiple Lie
brackets ${\rm ad}^j (K) L$ then transform to ${\rm ad}^j ( \overline{
K}) \overline{ L}$, thanks to the invariance of Lie
brackets. Geometrically, the vanishing of each of
the $n\times n$ minors defining $\Sigma''$ and $\Sigma '''$
means the linear dependence of a system of $n$ vectors, thus it is an
intrinsic condition. Consequently, although the jet variables $K_{ i,
\alpha}$ and $L_{ i, \alpha}$ are transformed in an unpleasant way
through diffeomorphisms, the sets $\Sigma'$, $\Sigma''$ and
$\Sigma'''$ may be defined by universal polynomials in the jet
variables $K_{ i, \alpha}$ and $L_{ i, \alpha}$, that are the same in
any system of local coordinates.

The lemma above and an inspection of a part of the complete expression
of the ${\rm ad}^j (K) (L)$, $0\leqslant j \leqslant 2n -1$ provides the
following information. Details will be skipped.

\def\thelemma{2.14}\begin{lemma}
The two subsets $\Sigma''$ and $\Sigma '''$ of $J^{ 2n-1} (\X (M)^2)$
are both algebraic in the jet variables and of codimension $(n+1)$
outside $\Sigma'$.
\end{lemma}

To conclude the proof of the theorem, we have to show that arbitrarily
close to $(K, L)$, there are pairs of finite type. Since $\Sigma'$ has
codimension $2n > \dim M$, a first application of the avoidance
Lemma~2.12 yields a perturbed pair, still denoted by $(K, L)$, with the
property that at every point $p\in M$, either $K(p) \neq 0$ or $L(p)
\neq 0$. Since $\Sigma''$ and $\Sigma'''$ both have codimension $n+1 >
\dim M$, a second application of the avoidance Lemma~2.12 yields a
perturbed pair such that the two collections of $2n$ vector fields
${\rm ad}^j (K) L$, $0 \leqslant j\leqslant 2n-1$, and ${\rm ad}^j (L) K$, $0
\leqslant j\leqslant 2n-1$, generate $TM$ at every point $p\in M$. The proof is
complete.
\endproof

To improve this theorem, let $r\geqslant 2$ and consider the set $\X (M)^r$
of collections of $r$ vector fields globally defined on $M$ that are
$\mathcal{ C}^{ \kappa -1}$ for some $\kappa \geqslant 2$ to be chosen
later. If $\LL^0 = \{ L_1, L_2, \dots, L_r\}$, is such a collection,
its elements may be expressed in a local chart $(x_1, \dots, x_n)$ as
$L_a = \sum_{ i= 1}^n\, \varphi_{ a, i} (x) \, \frac{ \partial }{
\partial x_i}$, for $a= 1, \dots, r$. Since the coefficients are
$\mathcal{ C}^{ \kappa - 1}$, it is possible to speak of $\LL^\lambda$
only for $\lambda \leqslant \kappa$. We want to determine the smallest
regularity $\kappa$ such that the set of 
$r$-tuples $\LL^0 \in \X (M)^r$ that are
of finite type at every point of $M$ is open and dense in $\X (M)^r$
for the strong Whitney topology.

As in \S1.8(II), let ${\sf n}_\kappa (r)$ denote the dimension of the
subspace ${\sf F}_\kappa (r)$ of the free Lie algebra ${\sf F} (r)$
that is generated as a real vector space by simple words (abstract Lie
brackets) of length $\leqslant \kappa$. Then ${\sf n}_\kappa (r)$ is the
maximal possible dimension of $\LL^\kappa (p)$ at a point $p\in M$. We
know that $\LL^\kappa$ is generated by simple iterated Lie brackets of
the form
\[
\big[
L_{a_1},
\big[
L_{a_2},\dots,
\big[L_{a_{\kappa-1}},L_{a_\lambda}\big]
\dots
\big]\big],
\]
for all $\lambda \leqslant \kappa$ and for certain (not all) $a_i$
with $1 \leqslant a_1, a_2, \dots, a_{ \lambda -1}, a_\lambda
\leqslant r$ that depend on the choice of a Hall-Witt basis
(Definition~1.9(II)) of $F_\kappa ( r)$.

We choose $\kappa$ minimal so that ${\sf n}_\kappa (r) \geqslant 2 \,
\dim M = 2n$. This fixes the smoothness of $M$. For $b= r+1, \dots,
{\sf n}_\kappa (r)$, we order linearly as $L_b = \sum_{ i= 1}^n\,
\psi_{ b,i} (x) \, \frac{ \partial }{ \partial x_i}$ the chosen
collection of iterated Lie brackets that generate $\LL^\kappa$. If
$\lambda = \lambda (b)$ denotes the length of $L_b$, namely $L_b \in
\LL^{ \lambda (b)}$ of the form $L_b = \big[ L_{ a_1}, \dots, \big[
L_{ a_{ \lambda (b)-1}}, L_{ a_{ \lambda (b)}} \big] \dots \big]$,
there are universal differential polynomials $A_{ a_1, \dots, a_{
\lambda(b) }}^i$ in the $(\lambda (b)- 1)$-th jet of the coefficients
$\varphi_{ a, i}$ such that $\psi_{ b, i} ( x) = A_{ a_1, \dots, a_{
\lambda (b)}}^i \big( J_x^{ \lambda (b) -1} \varphi (x) \big)$. Also,
in a fixed local system of coordinates, we form the $n
\times (2n)$ matrix
\[
\big(
\varphi_{1,i}\
\dots\
\varphi_{r,i}\ \
\psi_{r+1,i}\
\dots\
\psi_{2n,i}
\big)_{1\leqslant i\leqslant n}.
\]
Similarly as in the proof of the previous theorem, we define a ``bad''
subset $\Sigma$ of $J^{ \kappa - 1} ( \X (M)^r)$ by requiring that the
dimension of $\LL^\kappa (p)$ is $\leqslant (n-1)$ at every point $p\in
M$. This geometric condition is intrinsic and neither depends on the
choice of local coordinates nor on the choice of a Hall-Witt basis.
Concretely, in a local system of coordinates, $\Sigma$ is described as
the zero-set of all $n\times n$ minors of the above matrix. Thanks to
Lemma~2.13 and to an inspection of a portion of the explicit
expressions of the jet polynomials $A_{ a_1, \dots , a_{ \lambda (b) }
}^i \big( J_x^{ \lambda (b) -1} \varphi (x) \big)$, we may establish
the following assertion.

\def\thelemma{2.15}\begin{lemma}
The so defined subset $\Sigma = \big\{
\dim \LL^\kappa (p) \leqslant
n-1, \, \forall \, p \in M \big\}$ 
of $J^{ \kappa - 1} (\X (M)^r)$ is algebraic in
the jet variables and of codimension $(n+1)$.
\end{lemma}

Then an application of the avoidance Lemma~2.12 yields that, after an
arbitrarily small perturbation of $\LL^0$, still denoted by $\LL^0$,
we have $\LL^\kappa (p) = T_p M$ for every $p\in M$. Equivalently, the
type $d(p)$ of $p$ is finite at every point and satisfies $d(p) \leqslant
\kappa$.

\def\thetheorem{2.16}\begin{theorem}
Let $r\geqslant 2$ be an integer and assume that the connected
$n$-dimensional abstract manifold $M$ is $\mathcal{ C }^\kappa$, where
$\kappa$ is minimal with the
property that the dimension ${\sf n}_\kappa (r)$ of the vector
subspace ${\sf F}_\kappa (r)$, of the free Lie algebra ${\sf
F} (r)$ having $r$ generators, that is generated by all
brackets of length $\leqslant \kappa$, satisfies 
\[
{\sf n}_\kappa(r)\geqslant 2\dim M 
= 
2\,n. 
\]
Then the set of collections of $r$ vector fields $\LL^0 \in \X (M)^r$
that are of type $\leqslant \kappa$ at every point is open and dense
in $\X (M)^r$ for the strong Whitney topology.
\end{theorem}

A more general problem about finite-typisation of vector field
structures is concerned with general substructures of a given finite
type structure.

\def\theopenquestion{2.17}\begin{openquestion}
Given a finite type collection $\K^0 = \{ K_b \}_{ 1\leqslant b
\leqslant s}$, $s \geqslant 3$, of $\mathcal{ C}^{ \kappa -1}$ vector
fields on $M$ of class $\mathcal{ C }^\kappa$ with the property that
$\K^\kappa (p) = T_p M$ at every point and given a $\mathcal{ C }^{
\kappa - 1}$ subsystem $\LL^0 = \{ L_a \}_{ 1\leqslant a \leqslant
r}$, $2\leqslant r\leqslant s-1$, of the form $L_a = \sum_{ 1\leqslant
b \leqslant s}\, \psi_{a, b} \, K_b$, is it always possible to perturb
slightly the functions $\psi_{ a,b} : M \to \R$ so as to render
$\LL^0$ of finite type at every point\,? If so, what is the smallest
regularity $\kappa$, in terms of $r$, $s$ and the highest type of
$\K^0$ at points of $M$\,?
\end{openquestion}

Finally, we mention a result similar to Theorem~2.16 that is valid in the
$\mathcal{ C}^2$ category and does not use any Lie bracket. It is
based on Sussmann's orbit Theorem~1.21. The reference~\cite{
su1976} deals with several other genericity properties, motivated by
Control Theory.

\def\thetheorem{2.18}\begin{theorem}
{\rm (\cite{ su1976})} Assume $r\geqslant 2$ and $\kappa \geqslant 2$. The set
of collections $\LL^0 = \{ L_a \}_{ 1\leqslant a \leqslant r}$ of $r$ vector
fields on a connected $\mathcal{ C }^\kappa$ manifold $M$ so that $M$
consists of a single $\LL$-orbit, is open and dense in $\X (M)^r$
equipped with the strong Whitney $\mathcal{ C}^{\kappa - 1}$ topology.
\end{theorem}

\subsection*{ 2.19.~Transition}
The next Section~3 exposes the point of view of Analysis, where vector
field systems are considered as partial differential operators, until
we come back to the applications of the notion of orbits to CR
geometry in Section~4.

\section*{ \S3.~Locally integrable CR structures}

\subsection*{3.1.~Local insolvability of partial differential equations}
Until the 1950's, among analysts, it was believed and expected that
all linear partial differential equations having smooth coefficients
had local solutions (\cite{ trv2000}). In fact, elliptic, parabolic,
hyperbolic and constant coefficient equations were known to be locally
solvable. Although his thesis subject was to confirm this expectation
in full generality, in 1957, Hans Lewy (\cite{ lew1957}) exhibited a
striking and now classical counterexample of a $\mathcal{ C}^\infty$
function $g$ in a neighborhood of the origin of $\R^3$, such that
$\overline{ L} f = g$ has no local solution at all. Here, $\overline{
L} = \frac{\partial }{\partial \bar z} + z \frac{ \partial }{\partial
v}$ is the generator of the Cauchy-Riemann anti-holomorphic bundle
tangential to the Heisenberg sphere of equation $v = z\bar z$ in
$\C^2$, equipped with coordinates $(z, w) = (x + iy, u + iv)$.

From the side of Analysis, almost absent in the two grounding
works~\cite{ po1907} and~\cite{ ca1932} of Henri Poincar\'e and of
\'Elie Cartan, Lewy's discovery constituted the birth of smooth linear
PDE theory and of smoooth
Cauchy-Riemann geometry. Later, in 1971, the simpler
two-variables Mizohata equation $\frac{ \partial f}{\partial x} - i
x^k \, \frac{ \partial f}{\partial y} = g$ was shown by Grushin to be
non-solvable, if $k$ is odd, for certain $g$. One may verify that the
set of smooth functions $g$ for which Lewy's or Grushin's equation is
insolvable, even in the distributional sense, is generic in the sense
of Baire. For $k=1$, the {\sl Mizohata vector field} $\frac{ \partial
}{\partial x} - i \, x\, \frac{ \partial}{\partial y}$ intermixes the
holomorphic and antiholomorphic structures, depending on the sign of
$x$.

In 1973 answering a question of Lewy, Nirenberg (\cite{ ni1973})
exhibited a perturbation $\frac{ \partial }{\partial x} - i\, x ( 1 +
\varphi ( x, y)) \, \frac{\partial }{\partial y}$ of the Mizohata
vector field, where $\varphi$ is $\mathcal{ C}^\infty$ and null for $x
\leqslant 0$, such that the only local solutions of $Lf=0$ are the
constants. A year later, in~\cite{ ni1974}, he exhibited a
perturbation of the Lewy vector field having the same property. A
refined version is as follows. 

Let $\Omega$ be a domain in $\R^3$,
exhausted by a countable family of compact sets $K_j$, $j = 1, 2,
\dots$ with $K_j \subset {\rm Int} \, K_{ j+1 }$. If $f \in \mathcal{
C}^\infty ( \Omega, \C)$, define the Fr\'echet semi-norms $\rho_j (f)
:= \max_{ x \in K_j, \, \vert \alpha \vert \leqslant j} \, \left \vert
\partial_x^\alpha f (x) \right \vert$ and topologize $\mathcal{ C
}^\infty ( \Omega, \C)$ by means of the metric $d ( f, g) := \sum_{
j=1 }^{ \infty } \, \frac{ \rho_j (f-g)}{ 1 + \rho_j ( f-g
)}$. Consider the set
\[
\widehat{\bf L}
:=
\Big\{
L=
\sum_{j=1}^3\,a_j(x)\,
\frac{\partial}{\partial x_j}:a_j\in
\mathcal{C}^\infty(\Omega,\C)
\Big\},
\]
equipped with this topology for each coefficient $a_j$.

\def\thetheorem{3.2}\begin{theorem}
{\rm (\cite{ jt1982, ja1990})}
The set of $L \in \widehat{ \bf L}$ for which the solutions $u \in
\mathcal{ C}^1 ( \Omega, \C)$ of $L u = 0$ are the constants only, is
dense in $\widehat{ \bf L}$.
\end{theorem}

These phenomena and others were not suspected at the time of Lie, of
Poincar\'e, of \'E.~Cartan, of Vessiot and of Janet, when PDE theory
was focused on the algebraic complexity of systems of differential
equations having analytic coefficients. In 1959, H\"ormander
explained the behavior of the Lewy counter-example, as follows. The
references~\cite{ trv1970, trv1986, es1993,
trv2000} provide further survey informations about operators of
principal type, operators with multiple characteristics,
pseudodifferential operators, hypoelliptic operators, microlocal
analysis, {\it etc.}

Let $P = P (x, \partial_x) = \sum_{ \alpha \in \N^n, \, \vert \alpha
\vert \leqslant m} \, a_\alpha (x) \, \partial_x^\alpha$ be a linear
partial differential operator of degree $m$ having $\mathcal{
C}^\infty$ complex-valued coefficients $a_\alpha : \Omega \to \C$
defined in a domain $\Omega \subset \R^n$. Its {\sl symbol} $P (x,
\xi) := \sum_{ \vert \alpha \vert \leqslant m} \, a_\alpha (x) \, \left(
i\, \xi \right)^\alpha$ is a function from the cotangent $T^* \Omega
\equiv \Omega \times \R^n$ to $\C$. Its {\sl principal symbol} is the
homogeneous degree $m$ part $P_m (x, \xi) := \sum_{ \vert \alpha \vert
= m } \, a_\alpha (x) \left( i\, \xi \right)^\alpha$. The cone of
points $(x, \xi) \in \Omega \times (\R^n \backslash \{ 0 \} )$ such
that $P_m ( x, \xi) = 0$ is the {\sl characteristic set} of $P$, the
locus of the obstructions to existence as well as to regularity of
solutions $f$ of $P (x, \partial_x) f = g$.

The real characteristics of $P$ are called {\sl simple} if, at every
characteristic point $(x_0, \xi_0)$ with $\xi_0 \neq 0$, the
differential $d_\xi P_m = \sum_{ k=1 }^n\, \frac{ \partial P_m}{
\partial \xi_k} \, d\xi_k$ with respect to $\xi$ is nonzero. It
follows from homogeneity and from Euler's identity that the zeros of
$P$ are simple, so the characteristic set is a regular hypersurface of
$\Omega \times (\R^n \backslash \{ 0 \})$. One can show that this
assumption entails that the behaviour of $P$ is the same as that of
$P_m$: in a certain rigorous sense, lower order terms may be
neglected. In his thesis (1955), H\"ormander called such operators
{\sl of principal type}, a label that has stuck (\cite{ trv1970}).

Call $P$ {\sl solvable at a point} $x_0 \in \Omega$ if there exists a
neighborhood $U$ of $x_0$ such that for every $g\in \mathcal{
C}^\infty (U)$, there exists a distribution $f$ supported in $U$ that
satisfies $Pf = g$ in $U$. In 1955, H\"ormander had shown that a
principal type partial differential operator $P$ is locally solvable
if all the coefficients $a_\alpha (x)$, $\vert \alpha \vert = m$, of
its principal part $P_m$ are real-valued. On the contrary, 
if they are complex-valued, in 1959, he showed:

\def\thetheorem{3.3}\begin{theorem}
{\rm (\cite{ ho1963})}
If the quantity
\[
\sum_{k=1}^n\,
\frac{\overline{\partial P_m(x,\xi)}}{\partial \xi_k}\,
\frac{\partial P_m(x,\xi)}{\partial x_k}
\]
is nonzero at a characteristic point $(x_0, \xi_0) \in T^* \Omega$,
for some $\xi_0 \neq 0$, then $P$ is insolvable at $x_0$.
\end{theorem}

With $\overline{ P}_m (x, \partial_x) := \sum_{\vert \alpha \vert = m}
\, \overline{ a_\alpha (x)} \, \partial_x^\alpha$, denote by $C_{
2m-1} (x, \xi)$ the principal symbol of the commutator $\left[ P_m (
x, \partial_x), \overline{ P}_m (x, \partial_x) \right]$, obviously
zero if $P_m$ has real coefficients. The above necessary condition for
local solvability may be rephrased as: if $P$ is locally solvable at
$x_0$, then for all $\xi \in \R^n \backslash \{ 0\}$:
\[
P_m(x,\xi)
=
0
\ 
\Longrightarrow
\
C_{2m-1}(x,\xi)=0.
\]
This condition explained the non-solvability of the Lewy operator
appropriately.

\subsection*{ 3.4.~Condition {\bf (P)} of Nirenberg-Treves and local 
solvability} The geometric content of the above necessary
condition was explored and generalized by Nirenberg-Treves (\cite{
nt1963, nt1970, trv1970}). Recall that the {\sl
Hamiltonian vector field} associated to a function $f = f(x, \xi) \in
\mathcal{ C}^1 (\Omega \times \R^n)$ is $H_f := \sum_{ k=1 }^n \,
\left( \frac{ \partial f}{\partial \xi_k}\, \frac{ \partial }{\partial
x_k} - \frac{ \partial f}{\partial x_k}\, \frac{ \partial }{\partial
\xi_k} \right)$. A {\sl bicharacteristic} of the
real part $A (x,\xi)$ of $P_m (x, \xi)$ is an integral curve of $H_{
A}$, namely:
\[
\frac{dx}{dt}
=
{\rm grad}_\xi\, A(x,\xi), 
\ \ \ \ \ \ \
\frac{d\xi}{dt}
=
-
{\rm grad}_x\, A(x,\xi).
\]
It follows at once that the function $A(x,\xi)$ must be constant along
its bicharacteristics. When the constant is zero, a bicharacteristic
is called a {\sl null bicharacteristic}. In particular, null
bicharacteristics are contained in the characteristic set, which
explains the terminology.

Then H\"ormander's necessary condition may be interpreted as
follows. Let $B (x, \xi)$ be the imaginary part of $P_m (x, \xi)$. An
immediate computation shows that the principal symbol of $\left[ A(x,
\partial_x ), B(x, \partial_x) \right]$ is given by:
\[
C_1(x,\xi)
=
\sum_{k=1}^n\,
\left\{
\frac{\partial A}{\partial \xi_k}(x,\xi)\,
\frac{\partial B}{\partial x_k}(x,\xi)
-
\frac{\partial B}{\partial \xi_k}(x,\xi)\,
\frac{\partial A}{\partial x_k}(x,\xi)
\right\}.
\]
Equivalently, 
\[
C_1(x,\xi)
=
\left(
dB/dt
\right)(x,\xi).
\]
Theorem~3.3 says that the nonvanishing of $C_1$ at a characteristic
point entails insolvability. In fact, Nirenberg-Treves observed that
if the order of vanishing of $B$ along the null characteristic of $A$
is odd then insolvability holds. Beyond finite order of vanishing,
what appeared to matter is only the change of sign. Since the
equation $P f = g$ has the same solvability properties as $z \, Pf =
g$, for all $z\in \C \backslash \{ 0\}$, this led to the following:

\def\thedefinition{3.5}\begin{definition}{\rm
(\cite{ nt1963, nt1970}) A differential operator $P$
of principal type is said to {\sl satisfy condition} {\bf (P)}
if, for every $z\in \C \backslash \{ 0\}$, the function ${\rm Im} \,
(z\, P_m )$ does not change sign along 
the null bicharacteristic of ${\rm Re} \, (z\, P_m)$.
}\end{definition}

The next theorem has been shown for $P$ having $\mathcal{ C}^\omega$
coefficients and in certain cases for $P$ having $\mathcal{ C}^\infty$
coefficients by Nirenberg-Treves, and finally, in the general
$\mathcal{ C}^\infty$ category by Beals-Fefferman (sufficiency) and by
Moyer (necessity).

\def\thetheorem{3.6}\begin{theorem}
{\rm (\cite{ nt1963, trv1970, nt1970, befe1973, trv1986})} Condition
{\bf (P)} is necessary and sufficient for the local solvability in
$L^2$ of a principal type linear partial differential equation $P f =
g$.
\end{theorem}

Except for complex and strongly pseudoconvex structures, little is
known about solvability of $\mathcal{ C}^\infty$ systems of PDE's,
especially overdetermined ones (\cite{ trv2000}). In the sequel, only
vector field systems (order $m=1$), studied for themselves, will be
considered.

\subsection*{3.7.~Involutive and CR structures}
Following~\cite{ trv1981, trv1992, bch2005}, let $M$ be a $\mathcal{ C
}^\omega$, $\mathcal{ C}^\infty$ or $\mathcal{ C}^{\kappa, \alpha}$
$(\kappa \geqslant 2$, $0 < \alpha < 1$) paracompact Hausdorff second
countable abstract real manifold of dimension $\mu \geqslant 1$ and
let $\mathcal{ L}$ be a $\mathcal{ C}^\omega$, $\mathcal{ C}^\infty$
or $\mathcal{ C}^{\kappa -1, \alpha}$ complex vector subbundle of $\C
TM := \C \otimes TM$ of rank $\lambda$, with $1\leqslant \lambda
\leqslant \mu$. Denote by $\mathcal{ L}_p$ its fiber at a point $p\in
M$. Denote by $\mathcal{ T}$ the orthogonal of $\mathcal{ L}$ for the
duality between differential forms and vector fields. It is a vector
subbundle of $\C T^* M$, whose fiber at a point $p \in M$ is
$\mathcal{ L}_p^\bot = \left\{ \varpi 
\in \C T_p^* M : \varpi = 0 \ \, {\rm on}
\ \, \mathcal{ L}_p \right\}$. The {\sl characteristic set} $\mathcal{
C} := \mathcal{ T} \cap T^* M$ (real $T^* M$) is in general not a
vector bundle: the dimension of $\mathcal{ C }_p^0$ may vary with $p$,
as shown for instance by the bundle generated over $\R^2$ by the
Mizohata operator $\partial_x - ix \, \partial_y$.

From now on, we shall assume that the bundle $\mathcal{ L}$ is {\sl
formally integrable}, {\it i.e.} that $\left[ \mathcal{ L}, \mathcal{
L} \right] \subset \mathcal{ L}$. Then $\mathcal{ L}$ defines:

\smallskip

\begin{itemize}

\item[$\bullet$]
an {\sl elliptic structure} if $\mathcal{ C}_p = 0$ for all $p\in M$;

\smallskip

\item[$\bullet$]
a {\sl complex structure} of $\mathcal{ L}_p \oplus \overline{
\mathcal{ L}_p} = \C T_p M$ for all $p\in M$; 

\smallskip

\item[$\bullet$]
a {\sl Cauchy-Riemann} (CR for short) {\sl structure} if $\mathcal{
L}_p \cap \overline{ \mathcal{ L}_p} = \{ 0\}$ for all $p\in M$;

\smallskip

\item[$\bullet$]
an {\sl essentially real structure} if
$\mathcal{ L}_p = \overline{ \mathcal{ L}_p}$, 
for all $p\in M$.

\end{itemize}

\smallskip

In general, $\mathcal{ L}$ will be called an {\sl involutive
structure} if
$[ \mathcal{ L}, \mathcal{ L} ] = 
\mathcal{ L}$. Let us summarize basic linear algebra properties (\cite{
trv1981, trv1992, bch2005}). Every essentially real
structure is locally generated by real vector fields. Every complex
structure is elliptic. If $\mathcal{ L }$ is a CR structure (often
called {\sl abstract}), the characteristic set $\mathcal{ C}$ is in
fact a vector subbundle of $T^* M$ of rank $\mu - 2\lambda$; this
integer is the {\sl codimension} of the CR structure. A CR structure
is {\sl of hypersurface type} if its codimension equals $1$.




\subsection*{ 3.8.~Local integrability and generic submanifolds of
$\C^n$} The bundle $\mathcal{ L}$ is {\sl locally integrable} if every
$p\in M$ has a neighborhood $U_p$ in which there exist $\tau := \mu -
\lambda$ functions $z_1, \dots, z_\tau : U_p \to \C$ of class
$\mathcal{ C }^\omega$, $\mathcal{ C }^\infty$ or $\mathcal{ C}^{
\kappa, \alpha}$ whose differentials $dz_1, \dots, dz_\tau$ are
linearly independent and span $\mathcal{ T} \vert_{ U_p}$ (or
equivalently, are annihilated by sections of $\mathcal{ L}$). In other
words, the homogeneous PDE system $\mathcal{ L} f = 0$ has the best
possible space of solutions.

Here is a canonical example of locally integrable structure. Consider
a generic submanifold $M$ of $\C^n$ of class $\mathcal{ C}^\omega$,
$\mathcal{ C}^\infty$ or $\mathcal{ C}^{\kappa, \alpha}$, $\kappa
\geqslant 1$, $0 \leqslant \alpha \leqslant 1$, as defined in
\S2.1(II) and in \S4.1 below. Let $d\geqslant 0$ be its codimension
and let $m = n-d \geqslant 0$ be its CR dimension. Let $T^c M = TM
\cap JTM$ (a real vector bundle) and let $\C TM = \C \otimes
TM$. Define the two complex subbundles $T^{ 1, 0} M$ and $T^{ 0, 1} M
= \overline{ T^{ 1, 0} M}$ of $\C TM$ whose fibers at a point $p\in M$
are:
\[
\left\{
\aligned
{}
&
T_p^{1,0}M=\{X_p+iJX_p: \, 
X_p\in T_p^cM\}=
\{Z_p\in \C T_pM: \, 
JZ_p=-iZ_p\},\\
&
T_p^{0,1}M=\{X_p-iJX_p: \, 
X_p\in T_p^cM\}=
\{Z_p\in\C T_pM: \, 
JZ_p=iZ_p\}.
\endaligned\right.
\]
Geometrically, $T^{ 1, 0} M$ and $T^{ 0, 1} M$ are just the traces on
$M$ of the holomorphic and anti-holomorphic bundles $T^{ 1, 0} \C^n$
and $T^{ 0, 1} \C^n$, whose fibers at a point $p$ are $\sum_{ k=1 }^n
\, a_k \, \frac{ \partial }{\partial z_k} \big\vert_p$ and $\sum_{ k=1
}^n \, b_k \, \frac{ \partial }{\partial \bar z_k} \big\vert_p$. They
satisfy the Frobenius involutivity conditions $\left[ T^{ 1, 0} M, T^{ 1,
0} M \right] \subset T^{ 1, 0} M$ and $\left[ T^{ 0, 1} M, T^{ 0, 1} M
\right] \subset T^{ 0, 1} M$. More detailed background information
may be found in~\cite{ ch1991, bo1991, trv1992, ber1999}.

On such an embedded generic submanifold $M$, choose as structure
bundle $\mathcal{ L}$ just $T^{ 0, 1} M \subset \C TM$. Then clearly,
the $n$ holomorphic coordinate functions $z_1, \dots, z_n$ are
annihilated by the anti-holomorphic local sections $\sum_{ k=1 }^n \,
b_k \, \frac{ \partial }{\partial \bar z_k}$ of $T^{ 0, 1} M$ and they
have linearly independent differential, at every point of $M$. {\it A
generic submanifold embedded in $\C^n$ carries a locally integrable 
involutive structure}. Conversely:

\def\thelemma{3.9}\begin{lemma}
Every locally integrable CR structure is locally realizable as the
anti-holomorphic structure induced on a generic submanifold 
embedded 
$\C^n$.
\end{lemma}

\proof
Indeed, if a real $\mu$-dimensional $\mathcal{ C }^\omega$, $\mathcal{
C }^\infty$ or $\mathcal{ C }^{ \kappa, \alpha }$ manifold $M$ bears
a locally integrable CR structure, the map $Z = (z_1, \dots, z_\tau)$
produces an embedding of the open set $U_p$ as a local generic
submanifold $M := Z (U_p)$ of $\C^\tau$, with $Z_* ( \mathcal{ L} ) =
T^{ 0, 1} M$.
\endproof

A locally integrable CR structure is sometimes called locally {\sl
realizable} or locally
{\sl embeddable}.

\subsection*{ 3.10.~Levi form} 
Let $\mathcal{ L}$ be an involutive structure, not necessarily locally
integrable and let ${\sf c}_p \in \mathcal{ C}_p \subset T_p^* M$ be a
nonzero characteristic covector at $p$.

\def\thedefinition{3.11}\begin{definition}{\rm
The {\sl Levi form at $p$ in the characteristic codirection} ${\sf
c}_p \in \mathcal{ C}_p \backslash \{ 0\} \subset T_p^* M \backslash
\{ 0\}$ is the Hermitian form acting on two vectors $X_p, Y_p \in
\mathcal{ L} (p)$ as:
\[
\mathfrak{ L}_{ p, {\sf c}_p}
\left(
X_p, \overline{ Y}_p
\right)
:=
\frac{ 1}{2i}\,
{\sf c}_p 
\left(
\left[
X, \overline{ Y}
\right]
\right),
\]
where $X$, $Y$ are any two sections of $\mathcal{ L}$ defined in a
neighborhood of $p$ and satisfying $X(p) = X_p$, $Y (p) = Y_p$.
The resulting number is independent of the choice of
such extensions $X$, $Y$.
}\end{definition}

For the study of realizability of CR structures of codimension one,
nondegeneracy of the Levi form, especially positivity or negativity,
is of crucial importance. An abstract CR structure of hypersurface
type whose Levi form has a definite signe is said to be {\sl strongly
pseudoconvex}, since, after a possible rescaling of sign of a nonzero
characteristic covector, all the eigenvalues of its Levi form are
positive.

\subsection*{ 3.12.~Nonembeddable CR structures}
After Lewy's discovery, the first example of a smooth strictly
pseudoconvex CR structure in real dimension $3$ which is not locally
embeddable was produced by Nirenberg in 1973 (\cite{ ni1973}), {\it
cf.}~Theorem~3.3 above. For CR structures of hypersurface type,
Nirenberg's work has been generalized in higher dimension under the
assumption that the Levi form is neither positive nor negative, in any
characteristic codirection. Let $n \geqslant 2$ and let $\varepsilon_1 =
1$, $\varepsilon_k = -2$, $k = 2, \dots, n$.

\def\thetheorem{3.13}\begin{theorem}
{\rm (\cite{ jt1982, bch2005})} 
There exists a $\mathcal{ C}^\infty$ complex-valued function $g = g(x,
y, s)$ defined in a neighborhood of the origin in $\C^n \times \R$ and
vanishing to infinite order along $\{ x_1 = y_1 = 0\}$ such that the
vector fields{\rm :}
\[
\widehat{ L}_j
:=
\frac{\partial}{\partial \bar z_j}
-
i\,\varepsilon_j\,z_j\left(
1+g(x,y,s)
\right)
\frac{\partial}{\partial s},
\]
are pairwise commuting and such that every $\mathcal{ C}^1$ solution
$h$ of $\widehat{ L}_j h = 0$, $j=1,2, \dots, n$ defined in a
neighborhood of the origin must satisfy $\frac{ \partial h}{\partial
s} (0) = 0$.
\end{theorem}

This entails that the involutive 
structure spanned by the $\widehat{ L}_j$ is
not locally integrable at $0$. One may establish that the set of such
$g$ is generic. Crucially, the Levi-form is of signature $(n-1)$.

\def\theopenproblem{3.14}\begin{openproblem}
Find versions of generic non-embeddability for CR structures of
codimension $1$ having degenerate Levi-form. Find higher codimensional
versions of generic non-embeddability.
\end{openproblem}

\subsection*{ 3.15.~Integrability of complex structures and 
embeddability of strongly pseudoconvex CR structures} Let us now
expose positive results. Every formally integrable essentially real
structure $\mathcal{ L} = {\rm Re}\, \mathcal{ L}$ is locally
integrable, thanks to Frobenius' theorem; however the condition
$\left[ \mathcal{ L}, \mathcal{ L} \right] \subset \mathcal{ L}$
entails the involutivity of ${\rm Re}\, \mathcal{ L}$ only in this
special case. Also, every analytic formally integrable CR structure is
locally integrable: it suffices to complexify the coefficients of a
generating set of vector fields and to apply the holomorphic version
of Frobenius' theorem. For complex structure, the proof is due to
Libermann (1950) and to Eckman-Fr\"olicher (1951), {\it see}~\cite{
ah1972a, trv1981, trv1992}.

\def\thetheorem{3.16}\begin{theorem}
Smooth complex structures are locally integrable.
\end{theorem}

This deeper fact has a long history, which we shall review
concisely. On real analytic surfaces, isothermal coordinates where
discovered by Gauss in 1825--26, before he published his {\sl
Disquisitiones generales circa superficies curvas}. In the 1910's, by
a nontrivial advance, Korn and Lichtenstein transferred this theorem
to H\"older continuous metrics.

\def\thetheorem{3.17}\begin{theorem}
Let $ds^2 = E \, du^2 + 2\, F \, du dv + G \, dv^2$ be a $\mathcal{
C}^{ 0, \alpha}$ $(0 < \alpha < 1)$ Gaussian metric defined in some
neighborhood of $0$ in $\R^2$. Then there exists a $\mathcal{ C}^{ 1,
\alpha}$ change of coordinates $(u, v) \mapsto (\tilde{ u}, \tilde{
v})$ fixing $0$ and a $\mathcal{ C}^{ 0, \alpha}$ function $\tilde{
\lambda} = \tilde{ \lambda} ( \tilde{ u}, \tilde{ v })$ such that{\rm
:}
\[
\tilde{\lambda} 
\left(
d\tilde{u}^2
+
d\tilde{v}^2
\right)
=
E\,du^2+2\,F\,dudv+G\,dv^2.
\]
\end{theorem}

A modern proof of this theorem based on the complex notation and on
the $\overline{ \partial}$ formalism was provided by Bers (\cite{
be1957}) and by Chern (\cite{ ch1955}). In the monograph~\cite{
ve1962}, deeper weakenings of smoothness assumptions are provided.

As a consequence of this theorem, complex structures of class
$\mathcal{ C}^{ 0, \alpha}$ on surfaces may be shown to be locally
integrable. Let us explain in length this corollary.

At first, remind that an {\sl almost complex structure} on
$2n$-dimensional manifold $M$ is a smoothly varying field $J = (J_p)_{
p\in M}$ of endomorphisms of $T_p M$ satisfying $J_p \circ J_p = -
{\rm Id}$. Thanks to $J$, as in the standard complex case, one may
define $T_p^{ 0, 1} M := \left\{ X_p + i J_p X_p : \, X_p \in T_p M
\right\}$ and then the bundle $\mathcal{ L} := T^{ 0, 1} M$ is a
complex structure in the PDE sense of \S3.7 above. Conversely, given
a complex structure $\mathcal{ L} \subset \C TM$, then locally in some
neighborhood $U_p$ of an arbitrary point $p \in M$, there exist local
coordinates $(x_1, \dots, x_n, y_1, \dots, y_n)$ vanishing at $p$ so
that $n$ complex vector fields of the form:
\[
Z_j:=
\sum_{k=1}^n\,a_{k,j}\,\partial_{x_k}
+
i\,\sum_{k=1}^n\,b_{k,j}\,\partial_{y_k},
\]
with $a_{k,j} (0) = \delta_{k,j} = b_{k,j} (0)$, span $\mathcal{ L}
\vert_{ U_p}$. The associated almost complex structure is obtained by
declaring that, at a point of coordinates $(x, y)$, one has:
\[
\small
J\left(
\sum_{k=1}^n\,a_{k,j}\,\partial_{x_k}
\right)
=
\sum_{k=1}^n\,b_{k,j}\,\partial_{y_k}
\ \ \
{\rm and}
\ \ \
J\left(
\sum_{k=1}^n\,b_{k,j}\,\partial_{y_k}
\right)
=
-
\sum_{k=1}^n\,a_{k,j}\,\partial_{x_k}.
\]

\def\thelemma{3.18}\begin{lemma}
The bundle $\mathcal{ L} \subset \C TM$ satisfies $\left[ \mathcal{
L}, \mathcal{ L} \right] \subset \mathcal{ L}$ if and only if, for
every two vector fields $X$ and $Y$ on $M$, the Nijenhuis
expression{\rm :}
\[
N(X,Y)
:=
\left[
J\,X,J\,Y
\right]
-
J
\left[
X,J\,Y
\right]
-
J
\left[
J\,X,Y
\right]
-
\left[
X,Y
\right]
\]
vanishes identically.
\end{lemma}

The proof is abstract nonsense. Also, one verifies that $N ( f \, X,
g\, Y) = f\, g \, N (X, Y)$ for every two smooth local function $f$
and $g$: the expression $N$ is of tensorial character. In
symplectic and in almost complex geometry (\cite{ ms1995}), 
the following is settled.

\def\thedefinition{3.19}\begin{definition}{\rm
The almost complex structure is called {\sl integrable} if, in some
neighborhood $U_p$ of every point $p\in M$ there exist $n$
complex-valued functions $z_1, \dots, z_n : U_p \to \C$ of class at
least $\mathcal{ C}^1$ and having linearly independent
differentials such that $dz_k \circ J = i \circ dz_k$, for $k=1,
\dots, n$.
}\end{definition}

One verifies that it is equivalent to require $\mathcal{ L} z_k = 0$,
$k=1, \dots, n$: integrability of an almost complex structure
coincides with local integrability of $\mathcal{ L} = T^{ 0,
1} M$.

\smallskip

Now, we may come back to the integrability Theorem~3.16. To an
arbitrary Gaussian metric $g = ds^2$ as in Theorem~3.17, with $E>0$,
$G>0$ and $EG-F^2 >0$, are associated both a volume form and an almost
complex structure:
\[
d{\rm vol}_g
:=
\sqrt{EG-F^2}\,du\wedge dv
\ \ \ 
{\rm and}
\ \ \ 
J_g
:=
\frac{1}{\sqrt{EG-F^2}}
\left(
\begin{array}{cc}
-F&-G\\
E&F
\end{array}
\right).
\]
Conversely, given a volume form and an almost complex structure $J$ on
a surface, an associated Riemannian metric is provided by:
\[
g(\cdot,\cdot)
:=
d{\rm vol}(\cdot,J\cdot).
\]
According to Korn's and Lichtenstein's theorem, there exist
coordinates in which the metric is conformally flat, equal to $\lambda
\left( du^2 + dv^2 \right)$. In these coordinates, the associated
complex structure is obviously the standard one: $J \partial_u =
\partial_v$ and $J \partial_v = - \partial_u$. In fact, any local
change of coordinates $(u, v) \mapsto (\tilde{ u}, \tilde{ v})$ which
respects orthogonality of the curvilinear coordinates, {\it i.e.}
transforms the Gaussian isothermal metric to a similar one
$\tilde{\lambda} \left( d \tilde{ u}^2 + d\tilde{ v}^2 \right)$,
commutes with $J$, so that the map $u + i\, v \mapsto \tilde{ u} + i\,
\tilde{ v}$ is holomorphic. In conclusion:

\def\thetheorem{3.20}\begin{theorem}
$\mathcal{ C}^{ 0, \alpha}$ $(0 < \alpha < 1)$
complex structures are locally integrable.
\end{theorem}

The generalization to several variables of the theorem of Korn and
Lichtenstein is due to Newlander-Nirenberg, who solved a question
raised by Chern. The proof was modified and the smoothness assumption
was perfected by several mathematicians: Nijenhuis-Woolf (\cite{
nw1963}), Malgrange, Kohn, H\"ormander, Nirenberg, Treves (\cite{
trv1992}) and finally Webster (\cite{ we1989c}) who used the
Leray-Koppelman $\overline{ \partial}$ homotopy formula together with
the Nash-Moser rapidly convergent iteration scheme for solving
nonlinear functional equations.

\def\thetheorem{3.21}\begin{theorem} 
{\rm ($\mathcal{ C}^{ 2n, \alpha}$, $0 < \alpha < 1$: \cite{ nn1957};
$\mathcal{ C}^{1, \alpha}$, $0 < \alpha < 1$: \cite{ nw1963,
we1989c}; $\mathcal{ C}^\infty$: \cite{ trv1992})} Suppose that on the
real manifold $M$ of dimension $2n \geqslant 4$, the formally integrable
complex structure $\mathcal{ L}$ is $\mathcal{ C}^\infty$ or
$\mathcal{ C}^{ \kappa - 1, \alpha}$, $\kappa \geqslant 2$, $0 < \alpha <
1$. Then there exist local complex-valued coordinates $(z_1, \dots,
z_n)$ annihilated by $\mathcal{ L}$ which are $\mathcal{ C}^\infty$ or
$\mathcal{ C}^{\kappa, \alpha}$.
\end{theorem}

Finally, an elementary linear algebra argument (\cite{ trv1981,
trv1992, bch2005}) enables to deduce local
integrability of $\mathcal{ C}^\infty$ or $\mathcal{ C}^{ \kappa - 1,
\alpha}$ elliptic structures from the above theorem. In fact, elliptic
structures are shown to be locally isomorphic to $\C^{\tau} \times
\R^{ \lambda -\tau}$, equipped with $\frac{ \partial }{\partial \bar
z_i}$, $\frac{ \partial }{\partial t_j}$.

\def\theproblem{3.22}\begin{problem}
Is a formally integrable involutive structure having
positive-dimensional characteristic set locally integrable\,?
\end{problem}

Again the history is rich. Integrability results are known only for
strongly pseudoconvex CR structures of hypersurface type. Solving a
question raised by Kohn in 1965, Kuranishi (\cite{ ku1982}) showed in
1982 that $\mathcal{ C }^\infty$ strongly pseudoconvex abstract CR
structures of dimension $\geqslant 9$ are locally realizable. His delicate
proof involved a study of the Neumann operator in $L^2$ spaces, for
solving the tangential Cauchy-Riemann equations, together with the
Nash-Moser argument. In 1987, Akahori (\cite{ ak1987}) modified the
technique of Kuranishi and included the case of dimension $7$.

In 1989, to solve an associated linearized problem, instead of the
Neumann operator, Webster used the totally explicit integral operators
of Henkin.

\def\thetheorem{3.23}\begin{theorem}
{\rm (\cite{ we1989a, we1989b})}
Let $M$ be a strongly pseudoconvex $(2 n-1)$-dimensional CR manifold
of class $\mathcal{ C}^\mu$. Then $M$ admit, locally near each
point, a holomorphic embedding of class $\mathcal{ C}^\kappa$,
provided
\[
n\geqslant 4,
\ \ \ \ \ \ 
\kappa\geqslant 21,
\ \ \ \ \ \
\mu\geqslant 6\,\kappa+5\,n-3.
\]
\end{theorem}

The main new ingredient in his proof was Henkin' local homotopy
operator $\overline{ \partial}_M$ on a hypersurface $M \subset \C^n$:
\[
f
=
\overline{\partial}_M\,P(f)
+
Q(\overline{\partial}_M\,f), 
\ \ \ \ \ \ \
f\ {\rm a} \ (0,1)-{\rm form},
\]
known to hold for $n \geqslant 4$. For this reason, Webster suspected the
existence of refinements based on an insider knowledge of $\overline{
\partial}$ techniques. In 1994, using a modified homotopy formula
yielding better $\mathcal{ C}^\kappa$-estimates, Ma-Michel~\cite{
mm1994} improved smoothness:
\[
\kappa\geqslant 18, 
\ \ \ \ \ \ \
\mu\geqslant\kappa+13.
\]
Up to now, the five dimensional remains open. In fact, the solvability
of $\overline{ \partial }_M\, f = g$ for a $(0, 1)$-form on a
hypersurface of $\C^3$ requires a special trick which does not lead to
a homotopy formula. Nagel-Rosay~\cite{ nr1989} showed the
nonexistence of a homotopy formula in the $5$-dimensional case,
emphasizing an obstacle.

\def\theopenproblem{3.24}\begin{openproblem}
Find generalizations of the Kuranishi-Akahori-Webster-Ma-Michel
theorem to {\rm higher codimension}, using the integral formulas for
solving the $\overline{ \partial }_M$ due to Ayrapetian-Henkin.
Replace the assumption of strong pseudoconvexity by finer
nondegeneracy conditions, {\it e.g.} weak pseudoconvexity and finite
type in the sense of Kohn.
\end{openproblem}

\subsection*{ 3.25.~Local generators of locally 
integrable structures} Abandoning these deep problems of local
solvability and of local realizability, let us survey basic properties
of locally integrable structures. Thus, let $\mathcal{ L}$ be a
$\mathcal{ C }^\infty$ or $\mathcal{ C}^{\kappa - 1, \alpha}$ locally
integrable structure of rank $\lambda$ on a $\mathcal{ C}^{\kappa,
\alpha}$ or $\mathcal{ C}^\infty$ manifold $M$ of dimension
$\mu$. Denote by $\tau = \mu - \lambda$ the dimension of $\mathcal{ T}
= \mathcal{ L}^\bot$. Let $p\in M$ and let $\delta_p$ denote the
dimension of $\mathcal{ C}_p = \mathcal{ T} \cap T_p^* M$. Notice that
$(\tau - \delta_p) + (\tau - \delta_p) + \delta_p + (\lambda - \tau +
\delta_p) = \tau + \lambda = \mu$ just below.

\def\thetheorem{3.26}\begin{theorem}
{\rm (\cite{ trv1981, trv1992, bch2005})}
There exist real coordinates{\rm :}
\[
\left(
x_1,\dots,x_{\tau-\delta_p}, 
y_1,\dots,y_{\tau-\delta_p},
u_1,\dots,u_{\delta_p},
s_1,\dots,s_{\lambda-\tau+\delta_p}
\right),
\]
defined in a neighborhood $U_p$ of $p$ and vanishing at $p$, and there
exist $\mathcal{ C}^\infty$ or $\mathcal{ C}^{\kappa, \alpha}$
functions $\varphi_j = \varphi_j (x,y,u,s)$ with $\varphi_j (0) = 0$,
$d\varphi_j (0) = 0$, $j=1, \dots, \delta_p$, such that the
differentials of the $\tau$ functions{\rm :}
\[
\left\{
\aligned
z_k
&
:=
x_k+i\,y_k,
\ \ \ \ \ \ \ \
k
=
1,\dots,\tau-\delta_p,
\\
w_j
&
:=
u_j
+
i\,\varphi_j(x,y,u,s),
\ \ \ \ \ \ \ \
j
=
1,\dots,\delta_p
\endaligned\right.
\]
span $\mathcal{ T} \vert_{ U_p}$. 
\end{theorem}
 
Since $d\varphi_j (0) = 0$, there exist unique coefficients $b_{l,j} =
b_{ l, j} (x,y,u,s)$ such that the vector fields:
\[
K_j:=\sum_{l=1}^{
\delta_p}\,b_{l,j}\,\frac{\partial}{\partial u_l},
\ \ \ \ \ \ \ \
k=1,\dots,\delta_p,
\]
satisfy $K_{ j_1} (w_{ j_2} ) = \delta_{ j_1, j_2 }$, for $j_1, j_2 =
1, \dots, \delta_p$. Define then the $\lambda$ vector fields:
\[
\left\{
\aligned
\overline{L}_k
&
:=
\frac{\partial}{\partial \bar z_k} 
-
i\,\sum_{l=1}^{\delta_p}\,
\frac{\partial \varphi_l}{\partial \bar z_k}\,K_l, 
\ \ \ \ \ \ \ \ \
k=1,\dots,\tau-\delta_p, 
\\
L_j'
&
:=
\frac{\partial}{\partial s_j}
-
i\,\sum_{l=1}^{\delta_p}\,
\frac{\partial \varphi_l}{\partial s_j}\,K_l, 
\ \ \ \ \ \ \ \ \
j=1,\dots,\lambda-\tau+\delta_p. 
\endaligned\right.
\]
Clearly, $0 = \overline{ L}_{ k_1} (z_{ k_2}) = \overline{ L}_k (w_j)
= L_j' ( z_k) = L_{ j_1} '( w_{ j_2} )$, hence the structure bundle
$\mathcal{ L}\vert_{ U_p}$ is spanned by the $L_k, L_j'$. One may
verify the commutation relations (\cite{ trv1981, trv1992,
bch2005}):
\[
\aligned
0
&
=
\left[\overline{L}_{k_1},\overline{L}_{k_2}\right]
=
\left[\overline{L}_k,L_j'\right]
=
\left[L_{j_1}',L_{j_2}'\right], 
\\
0
&
=
\left[\overline{L}_k,K_j\right]
=
\left[L_{j_1}',K_{j_2}\right]
=
\left[K_{j_1},K_{j_2}\right].
\endaligned
\]

Remind that if an involutive structure $\mathcal{ L}$ is CR, then
$\delta_p$ is independent of $p$ and equal to the codimension $\mu - 2
\lambda =: \delta$. It follows that $\tau - \delta = \tau - \mu + 2
\lambda = \lambda$, or $\lambda - \tau + \delta = 0$: this means that
the variables $( s_1, \dots, s_{ \lambda - \tau + \delta_p })$
disappear.

\def\thecorollary{3.27}\begin{corollary}
In the case of a CR structure of codimension $\delta$, 
the local integrals are{\rm :} 
\[
\left\{
\aligned
z_k
&
:=
x_k+i\,y_k, 
\ \ \ \ \ \ \ \
k
=
1,\dots,\lambda,
\\
w_j
&
:=
u_j
+
i\,\varphi_j(x,y,u),
\ \ \ \ \ \ \ \
j
=
1,\dots,\delta,
\endaligned\right.
\]
and a local basis for the structure bundle $\mathcal{ L}
\vert_{ U_p}$ is{\rm :}
\[
\overline{L}_k
:=
\frac{\partial}{\partial \bar z_k} 
-
i\,\sum_{l=1}^\delta\,
\frac{\partial \varphi_l}{\partial \bar z_k}\,K_l, 
\ \ \ \ \ \ \ \ \
k=1,\dots,\lambda.
\]
\end{corollary}

We recover a generic submanifold embedded in $\C^\tau$ which is
graphed by the equations $v_j = \varphi_j ( x, y, u)$, as in
Theorem~2.3(II), or as in Theorem~4.2 below.

\subsection*{ 3.28.~Local embedding into a CR structure}
But in general, the coordinates $( s_1, \dots, s_{ \lambda - \tau +
\delta_p} )$ are present. A trick (\cite{ ma1992}) is to introduce
extra coordinates $( t_1, \dots, t_{ \lambda - \tau + \delta_p} )$ and
to define a new structure on the product $U_p \times \R^{ \lambda -
\tau + \delta_p}$ generated by the following local integrals:
\[
\left\{
\aligned
\tilde{z}_k
&
:=
z_k,
\ \ \ \ \ \ \ \ \
\ \ \ \ \ \ \ \ \
\ \ \ \ \ \ \ \ \
\ \ \ \ \ \ \ \ \
k=1,\dots,\tau-\delta_p,
\\
\tilde{z}_k
&
:=
s_{k-\tau+\delta_p}
+
i\,t_{k-\tau+\delta_p},
\ \ \ \ \ \ \ \ \ \
k=\tau-\delta_p+1,\dots,\lambda,
\\
\tilde{w}_j
&
:=
w_j,
\ \ \ \ \ \ \ \ \
\ \ \ \ \ \ \ \ \
\ \ \ \ \ \ \ \ \
\ \ \ \ \ \ \ \ \
j=1,\dots,\delta_p.
\endaligned\right.
\]
The associated structure bundle $\widetilde{ \mathcal{ L}}$ is spanned
by:
\[
\left\{
\aligned
\overline{\widetilde{L}}_k
&
:=
\overline{L}_k,
\ \ \ \ \ \ \ \ \
k=1,\dots,\tau-\delta_p,
\\
\overline{\widetilde{L}}_j'
&
:=
\frac{1}{2}\,
L_j'
+
\frac{i}{2}\,
\frac{\partial}{\partial t_j}, 
\ \ \ \ \ \ \ \ \
j=1,\dots,\lambda-\tau+\delta_p.
\endaligned\right.
\]
It is a CR structure of codimension $\delta_p$ on $U_p \times \R^{
\lambda - \tau + \delta_p}$. All analytico-geometric objects defined
in $U_p$ can be lifted to $U_p \times \R^{ \lambda - \tau +
\delta_p}$, just by declaring them to be independent of the extra
variables $(t_1, \dots, t_{ \lambda - \tau + \delta_p})$.

Such an embedding enables one to transfer elementarily several
theorems valid in embedded Cauchy-Riemann Geometry, to the more
general setting of locally integrable structures. For instance, this
is true of most of the theorems about holomorphic or CR extension of
CR functions presented Part~V. In addition, most of the results stated
in \S3, \S4 and \S5 below hold in locally integrable structures.

\subsection*{ 3.29.~Transition}
However, for reasons of space and because the possible generalizations
which we could state by applying this embedding trick would require
dry technical details, we will content ourselves to just mention these
virtual generalizations, as was done in~\cite{ mp1999}.
For further study of locally integrable
structures, we refer mainly to~\cite{ trv1992, bch2005} and to the
references therein.

In summary, in this Section~3, we wanted to show how our approach is
inserted into a broad architecture of questions about solvability of
partial differential equations, about the problem of realizability of
abstract CR structures, as well as into hypo-analytic structures.
Thus, even if some of the subsequent surveyed results ({\em exempli
gratia}: the celebrated {\sl Baouendi-Treves approximation}
Theorem~5.2) were originally stated in the context of locally
integrable structures, even though we could as well state them in this
context or at least in the context of locally embeddable CR structures
(as was done in~\cite{ mp1999}), we shall content ourselves to state
them in the context of {\sl embedded Cauchy-Riemann geometry}, just
because the very core of the present memoir is concerned by Several 
Complex Variables topics: analytic discs, envelopes of
holomorphy, removable singularities, {\it etc.}

\section*{ \S4.~Smooth generic submanifolds and their CR orbits}

\subsection*{ 4.1.~Definitions of CR submanifolds and local 
graphing equations} We begin by some coordinate-invariant geometric
definitions. Some implicit lemmas are involved (the reader is referred
to~\cite{ ch1989, ch1991, bo1991, ber1999}). Let $J$ denote the
complex structure of $T\C^n$ ({\it see} \S2.1(II)). A real connected
submanifold $M \subset \C^n$ of class at least $\mathcal{ C}^1$ is
called:

\begin{itemize}

\smallskip\item[$\bullet$]
{\sl Totally real} if $T_p M \cap J T_p M = \{ 0\}$ for every $p\in
M$. Then $M$ has codimension $d\geqslant n$ and is called {\sl maximally
real} if $d = n$. The complex vector subspace $H_p := T_p M + J T_p
M$ of $T_p \C^n$ has complex dimension $2n - d$. If ${\sf proj}_{ H_p}
(\cdot)$ denotes any $\C$-linear projection of $T_p \C^n$ onto $H_p$
and if $\mathcal{ U}_p$ is a small neighborhood of $p$ in $\C^n$, then
${\sf proj}_{ H_p} (M \cap \mathcal{ U}_p)$ is maximally real in
$H_p$.

\smallskip\item[$\bullet$]
{\sl Generic} if $T_p M + J T_p M = T_p \C^n$ for every $p\in M$. Then
$M$ has codimension $d\leqslant n$ and is maximally real only if $d =
n$. Then $T_p M \cap JT_p M$ is the maximal $\C$-linear subspace of
$T_p M$ and has complex dimension equal to the integer $m := n-d$,
called the {\sl CR dimension} of $M$. It is obviously constant, as $p$
runs in $M$.

\smallskip\item[$\bullet$]
{\sl Cauchy-Riemann} ({\sl CR} for short) if the maximal $\C$-linear
subspace $T_p M \cap JT_pM$ of $T_p M$ has constant dimension $m$
(necessarily $\leqslant n$) for $p$ running in $M$. If $M$ has
codimension $d$, the integer $c := d - n + m$ is called the {\sl
holomorphic codimension} of $M$. Then for $p$ running in $M$, the
complex vector subspaces $H_p := T_p M + J T_p M$ of $T_p \C^n$ have
constant complex codimension $c$, which justifies the terminology. If
${\sf proj}_{ H_p} (\cdot)$ denotes any $\C$-linear projection of $T_p
\C^n$ onto $H_p$, and if $\mathcal{ U}_p$ is a small open neighborhood
of $p$ in $\C^n$, then
\[
\widetilde{M}_p
:=
{\sf proj}_{H_p}(M\cap\mathcal{U}_p)
\]
is a generic submanifold of $\C^{ n-c}$.

\end{itemize}\smallskip

In \S2.1(II), we have graphed totally real, generic and, generally,
Cauchy-Riemann local submanifolds $M \subset \C^n$, but only in the
algebraic and in the analytic category. In the smooth category, the
intrinsic complexification $\{ w_2 = 0\}$ disappears, but $H_p = T_p M
+ J T_p M$ still exists, so that further graphing functions are
needed.

\def\thetheorem{4.2}\begin{theorem}
{\rm (\cite{ ch1989, ch1991, bo1991, ber1999, me2004a})} Let $M
\subset \C^n$ be a real submanifold of codimension $d$ and let $p \in
M$. There exist complex algebraic or analytic coordinates centered at
$p$ and $\rho_1 > 0$ such that $M$, supposed to be $\mathcal{
C}^{\mathcal{ R}}$, where $\mathcal{ R} = \infty$ or where $\mathcal{
R} = (\kappa, \alpha)$, $\kappa \geqslant 1$, $0 \leqslant \alpha
\leqslant 1$, is locally represented as follows{\rm :}

\smallskip

\begin{itemize}
\item[$\bullet$]
If $M$ is {\rm totally real}, letting $d_1 = 2n - d \geqslant 0$ and $c = d
- n \geqslant 0$, then $d_1 + c = n$ and
\[
\aligned
M = \big\{ (w_1, w_2) 
\in 
&
\left(
\square_{ \rho_1 }^{ d_1} \times i\R^{ d_1} \right) \times \C^c : 
\\
& \ \
{\rm Im }\, w_1 = \varphi_1 ({\rm Re} \, w_1), w_2 = \psi_2 ({\rm Re}
\, w_1) \big\},
\endaligned
\]
for some $\R^{ d_1}$-valued $\mathcal{ C }^{ \mathcal{ R }}$ map
$\varphi_1$ and some $\C^c$-valued $\mathcal{ C }^{ \mathcal{ R }}$
map $\psi_2$ satisfying $\varphi_1 (0) = 0$ and $\psi_2 (0) = 0$.

\smallskip

\item[$\bullet$]
If $M$ is {\rm generic}, letting $m = d - n$, then $m + d = n$ and 
\[
M = \left\{ (z, w) \in \Delta_{ \rho_1}^m \times \left( \square_{
\rho_1 }^d \times i \R^d \right) : {\rm Im}\, w = \varphi ( {\rm Re}\, z,
{\rm Im}\, z, {\rm Re} \, w) \right\},
\]
for some $\R^d$-valued $\mathcal{ C}^{ \mathcal{ R}}$ map $\varphi$
satisfying $\varphi (0) = 0$.

\smallskip

\item[$\bullet$]
If $M$ is {\rm Cauchy-Riemann}, letting $m = {\rm CRdim }\, M$, $c =
d-n+m \geqslant 0$, and $d_1 = 2n - 2m -d \geqslant 0$, then $m + d_1 + c = n$
and 
\[
\aligned
M = 
&
\big\{ (z, w_1, w_2)
\in 
\Delta_{ \rho_1 }^m
\times 
\left(
\square_{ \rho_1}^{ d_1} \times i\R^{ d_1}
\right)
\times 
\C^c : 
\\
&
\ \
{\rm Im}\,
w_1 = \varphi_1 ( {\rm Re}\, z, {\rm Im}\, z, {\rm Re} \, w_1),
\ 
w_2 = \psi_2 ({\rm Re}\, z, {\rm Im}\, z, { \rm Re }\, w_1) \big\},
\endaligned
\]
for some $\R^{ d_1 }$-valued $\mathcal{ C }^{ \mathcal{ R }}$ map
$\varphi_1$ with $\varphi_1 (0) = 0$ and some $\C^c$-valued $\mathcal{
C}^{ \mathcal{ R }}$ map $\psi_2$ with $\psi_2 (0) = 0$ which is CR
{\rm (definition in \S4.25 below)} on the local generic submanifold
$M^1 := \left\{ (z, w_1) \in \Delta_{ \rho_1 }^m \times \left(
\square_{ \rho_1 }^{ d_1} \times i \R^{ d_1 } \right) : {\rm Im }\,
w_1 = \varphi_1 ( {\rm Re }\, z, {\rm Im }\, z, {\rm Re }\, w_1)
\right\}$.
\end{itemize}

\smallskip

An adapted linear change of coordinates insures that the differentials
at the origin of the graphing maps all vanish.
\end{theorem}

A CR manifold $M$ being always locally graphed above a generic
submanifold of $\C^{ n-c}$, the remainder of this memoir will mostly
be devoted to the study of $\mathcal{ C}^\infty$ or $\mathcal{
C}^{\kappa, \alpha}$ {\it generic}\, submanifolds of $\C^n$. The above
local representation of a generic $M$ will be constantly used.

\subsection*{ 4.3.~CR vector fields}
Let $M$ be generic, of class at least $\mathcal{ C}^1$, represented by
$v = \varphi (x,y,u)$ as above, in coordinates $(z, w) = (x+iy,
u+iv)$. Sometimes, we shall also write $v = \varphi (z, u)$, being it
clear that $\varphi$ is {\it not}\, holomorphic with respect to $z$.
Here, we provide a description in coordinates of three useful
families of vector fields.

There exist $m$ anti-holomorphic vector fields defined in $\Delta_{
\rho_1}^m \times \left( \square_{ \rho_1}^d \times i\R^d \right)$ of
the form:
\[
\overline{L}_k' 
=
\frac{\partial}{\partial\bar z_k}
+
\sum_{j=1}^d\,a_{j,k}'(x,y,u)\,
\frac{\partial}{\partial\bar w_j},
\]
whose restrictions to $M$ span $T^{ 0, 1} M$. To compute the
coefficients $a_{ j,k}'$, the conditions $0 \equiv \overline{ L}_k '
\left( \varphi_j (x,y,u) - v_j \right)$ yield:
\[
2\,
\varphi_{j,\bar z_k}
=
\sum_{l=1}^d\,
\left(
i\,\delta_{j,l}
-
\varphi_{j,u_l}
\right)
a_{l,k}'.
\]
In matrix notation, the solution is: $a' = 2\, \left(i\, I_{ d\times
d} - \varphi_u \right)^{ -1} \cdot \varphi_{\bar z}$, with $a' = (a_{
j,k}')_{1\leqslant j \leqslant d}^{ 1\leqslant k \leqslant m}$ and $\varphi_{ \bar z} =
\left( \varphi_{ j, \bar z_k} \right)_{1\leqslant j \leqslant d}^{ 1\leqslant k \leqslant
m}$ both of size $d\times m$. Since $\left[ T^{ 0, 1} M, T^{ 0, 1} M
\right] \subset T^{ 0, 1} M$, the $\overline{ L}_k'$ commute. They are
extrinsic.

Also, there exist $m$ intrinsic sections of $\C TM$ of the form:
\[
\overline{L}_k
:=
\frac{\partial}{\partial \bar z_k}
+
\sum_{j=1}^d\,a_{j,k}(x,y,u)\,
\frac{\partial}{\partial u_j},
\]
written in the coordinates $(x, y, u)$ of $M$, which span the
structure bundle $T^{ 0, 1} M \subset \C TM$. Since $(x,y,u)$ are
coordinates on $M$, restricting the $\overline{ L}_k' \vert_M$ to $M$
amounts to just drop the terms $\frac{ i}{ 2}\, \frac{ \partial
}{\partial v_j}$ in each $\frac{\partial }{ \partial \bar w_j }$
appearing in $\overline{ L}_k'$. Hence:

\def\thelemma{4.4}\begin{lemma}
One has $a_{j,k} = \frac{1 }{ 2} \, a_{ j,k }'$.
\end{lemma} 

Another argument is to first introduce the $d$ vector fields:
\[
K_j 
=
\sum_{l=1}^d\,b_{l,j}(x,y,u)\,
\frac{\partial}{\partial u_l},
\]
that are uniquely determined by the conditions 
\[
\delta_{j_1,j_2}
=
K_{j_1}
\left(
w_{j_2}\vert_M
\right)
=
K_{j_1}\left(u_{j_2}+
i\varphi_{j_2}(x,y,u)\right).
\]
Equivalently, the coefficients $b_{ l,j}$
satisfy:
\[
\delta_{j_1,j_2}
=
\sum_{l=1}^d\,
\left(
\delta_{j_2,l}
+
i\,\varphi_{j_2,u_l}
\right)
b_{l,j_1},
\]
whence, in matrix notation: $b = \left( I_{ d\times d} + i\, \varphi_u
\right)^{-1}$. Here, $b = \left( b_{l, j} \right)_{1\leqslant l\leqslant d}^{
1\leqslant j \leqslant d}$ and $\varphi_u = \left( \varphi_{ j, u_l}
\right)_{1\leqslant j \leqslant d}^{ 1\leqslant l \leqslant d}$ both are $d\times d$
matrices.

Similarly as in Corollary~3.27, the $\overline{ L}_k$ defined above
span the structure bundle $\mathcal{ L} = T^{ 0, 1} M$ having the
local integrals $z_1, \dots, z_m, w_1 \vert_M, \dots, w_d \vert_M$, if
and only if they satisfy $0 = \overline{ L}_k \left( w_j \vert_M
\right) = \overline{ L}_k \left[ u_j + i \, \varphi_j (x,y,u) \right]
$. Seeking the $\overline{ L}_k$ under the form $\overline{ L}_k =
\frac{\partial }{\partial \bar z_k} + \sum_{ l=1 }^d \, c_{ l,k} (x,
y, u) \, K_l$, it follows from $\delta_{ j_1, j_2} = K_{ j_1 }\left(
u_{j_2}+ i \, \varphi_{j_2} (x,y,u) \right)$ that $c_{ j,k} = - i \,
\varphi_{ j, \bar z_k}$. Reexpressing explicitly the $K_l$ in terms of
the $\frac{ \partial}{ \partial u_j}$ as achieved
above, we finally get in matrix
notation $a = \left(i\, I_{ d\times d} - \varphi_u \right)^{ -1}
\cdot \varphi_{\bar z}$. This yields a second, more intrinsic
computation of the coefficients $a_{ j,k}$ and a second proof of $a_{
j,k} = \frac{ 1}{ 2} \, a_{ j,k} '$.

\smallskip

If $\chi$ is a $\mathcal{ C}^\infty$ or $\mathcal{ C}^{ \kappa,
\alpha}$ function on $M$, its differential may be computed as 
\[
d\chi 
=
\sum_{ k=1 }^m \, L_k (\chi) \, dz_k + \sum_{ k=1 }^m \, \overline{
L}_k (\chi) \, d\bar z_k + \sum_{ j=1 }^d \, K_j ( \chi) \,
dw_j\big\vert_M.
\]

\def\thelemma{4.5}\begin{lemma}
{\rm (\cite{ trv1981, br1987, trv1992,
bch2005})} The following relations hold{\rm :}
\[
\left\{
\aligned
&
L_{ k_1} ( z_{ k_2} ) = \delta_{ k_1, k_2},
\ \ 
L_k(w_j)=0,
\ \ 
K_j(z_k)=0,
\ \ 
K_{j_1}(w_{j_2}\vert_M)=\delta_{j_1,j_2}, 
\\ 
&
\left[L_{k_1},L_{k_2}\right]
=
\left[L_k,K_j\right]
=
\left[K_{j_1},K_{j_2}\right]
=
0.
\endaligned\right.
\]
\end{lemma}

\subsection*{ 4.6.~Vector-valued Levi form} 
Let $p\in M$ and denote by $\pi_p$ the projection $\C T_p M
\longrightarrow \C T_p M / \left( T_p^{ 1, 0} M \oplus T_p^{ 0, 1} M
\right)$.

\def\thedefinition{4.7}\begin{definition}{\rm
The {\sl Levi map} at $p$ is the Hermitian $\C^d$-valued form acting
on two vectors $X_p, Y_p \in T_p^{ 1, 0} M$ as:
\[
\left[
\aligned
\mathfrak{L}_p:
\ \ \ \ \
&
T_p^{1,0}M\times T_p^{1,0}M
\longrightarrow
\C T_pM/\left(
T_p^{1,0}M\oplus T_p^{0,1}M
\right)
\\
&
\mathfrak{L}_p(X_p,Y_p)
:=
\frac{1}{2i}\,\pi_p
\left(
\left[
X,\overline{Y}
\right](p)
\right),
\endaligned\right.
\]
where $X, Y$ are any two sections of $T^{ 1, 0} M$ defined in a
neighborhood of $p$ satisfying $X(p) = X_p$, $Y(p) = Y_p$. The
resulting number is independent of the choice of 
such extensions $X, Y$.
}\end{definition}

As $p$ varies, this yields a smooth bundle map. The Levi map
$\mathfrak{ L}_p$ is {\sl nondegenerate at} $p$ if its kernel is null:
$\mathfrak{ L}_p ( X_p, Y_p) = 0$ for every $Y_p$ implies $X_p = 0$.
On the opposite, $M$ is {\sl Levi-flat} if the kernel of $\mathfrak{
L}_p$ equals $T_p^{ 1, 0} M$ at every $p$. If $M$ is a hypersurface
$(d=1)$, one calls $\mathfrak{ L}_p$ the {\sl Levi form at} $p$. Then
$M$ is {\sl strongly pseudoconvex at} $p$ if the Levi form $\mathfrak{
L}_p$ is definite, positive or negative. These definitions agree with
the ones formulated in \S3.10 for abstract CR structures.

\def\thetheorem{4.8}\begin{theorem}
{\rm (\cite{ fr1977, ch1991})}
In a neighborhood $U_p$ of a point $p \in M$ in which the kernel of
the Levi map is of constant rank and defines a $\mathcal{ C}^\omega$,
$\mathcal{ C}^\infty$ or $\mathcal{ C}^{ \kappa - 1, \alpha}$ $(
\kappa \geqslant 2, 0 \leqslant \alpha \leqslant 1)$ distribution of
$m_1$-dimensional complex planes $P_q \subset T_q M$, $\forall \, q\in
U_p$, the distribution is Frobenius-integrable, hence $M$ is
$\mathcal{ C}^\omega$, $\mathcal{ C}^\infty$ or $\mathcal{ C}^{ \kappa
- 1, \alpha}$ foliated by complex manifolds of dimension $m_1$.
\end{theorem}

In particular, a Levi-flat generic submanifold of CR dimension $m$ is
foliated by $m$-dimensional complex manifolds. These observations go
back to Sommer (1959). In~\cite{ ch1991}, one founds a systematic
study of foliations by complex and by CR manifolds.

\subsection*{ 4.9.~CR orbits} 
Let $M \subset \C^n$ be generic and consider the system $\LL$ of
sections of $T^cM$. To apply the orbit Theorem~1.21, we need $M$ to
be at least $\mathcal{ C}^2$, in order that the flows are at least
$\mathcal{ C}^1$. By definition, a weak $T^cM$-integral submanifold
$S \subset M$ satisfies $T_p S \supset T_p^c M$, at every point $p \in
S$. Equivalently, $S$ has the same CR dimension as $M$ at every point.

In the theory of holomorphic extension exposed in Part~V, local and
global CR orbits will appear to be adequate objects of study. They
constitute one of the main topics of this memoir.

\def\theproposition{4.10}\begin{proposition}
{\rm (\cite{ trp1990, tu1990, tu1994a, me1994, jo1996, mp1999,
mp2002})} Let $M \subset \C^n$ be generic of class $\mathcal{
C}^\omega$, $\mathcal{ C}^\infty$ or $\mathcal{ C}^{\kappa, \alpha}$
with $\kappa \geqslant 2$ and $0 \leqslant \alpha \leqslant 1$.

\begin{itemize}

\smallskip\item[{\bf (a)}]
The {\rm (}global{\rm)} 
$T^cM$-orbits are called {\rm CR orbits}. They are denoted by
$\mathcal{ O}_{ CR}$ or by $\mathcal{ O}_{ CR} (M, p)$, if the
reference to one of point $p \in \mathcal{ O}_{ CR}$ is needed.

\smallskip\item[{\bf (b)}]
The local CR orbit of a point $p \in M$ is denoted by $\mathcal{ O}_{
CR}^{ loc} (M, p)$. It is a local submanifold embedded in $M$, closed
in a sufficiently small neighborhood of $p$ in $M$.

\smallskip\item[{\bf (c)}]
Local and global CR orbits are $\mathcal{ C}^\omega$, $\mathcal{
C}^\infty$ or $\mathcal{ C}^{ \kappa, \beta}$, for
every $\beta$ with $0 < \beta < \alpha$.

\smallskip\item[{\bf (d)}]
$M$ is partitioned in global CR orbits. Each global CR orbit is
injectively immersed and weakly embedded in $M$, is a CR submanifold
of $\C^n$ contained in $M$ and has the same CR dimension as $M$.

\smallskip\item[{\bf (e)}]
Every {\rm (}immersed{\rm )} CR submanifold $S \subset M$ having the
same CR dimension as $M$ contains the local CR orbit of each of its
points

\smallskip\item[{\bf (f)}]
CR orbits of the smallest possible real dimension $2m = 2\, {\rm
CRdim}\, M$ satisfy $T_p \mathcal{ O}_{ CR} = T_p^c \mathcal{ O}_{
CR}$ at every point, hence are complex $m$-dimensional submanifolds.

\end{itemize}\smallskip

\end{proposition}

According to Example~1.29, CR orbits should be $\mathcal{ C}^{ \kappa
-1 , \alpha}$, not smoother. But in generic submanifolds, they also
can be described as boundaries of small attached analytic discs
(\cite{ tu1990, tu1994a, me1994}) and the
$\mathcal{ C}^{ \kappa, \alpha - 0} = \bigcap_{ \beta <
\alpha} \, \mathcal{ C}^{ \kappa, \beta}$
smoothness of the solution in Theorem~3.7(IV) 
explains {\bf (c)}.

\smallskip
 
Let us summarize some structural properties of CR orbits, useful in
applications. A specialization of Theorem~1.21{\bf (4)} yields the
following.

\def\theproposition{4.11}\begin{proposition}
For every $p\in M$, there exist $k\in \N$, sections $L^1, \dots, L^k$
of $T^cM$ and parameters ${\sf s}_1^*, \dots, {\sf s}_k^* \in \R$ such
that $L_{{\sf s}_k^*}^k (\cdots ( L_{{\sf s}_1^*}^1 (p)) \cdots) = p$
and the map $({\sf s}_1, \dots, {\sf s}_k) \longmapsto L_{{\sf s}_k}^k
(\cdots ( L_{{\sf s}_1}^1 (p)) \cdots)$ is of rank $\dim \mathcal{
O}_{ CR} (M, p)$ at $({\sf s}_1^*, \dots, {\sf s}_k^*)$.
\end{proposition}

The dimension of any $\mathcal{ O}_{ CR}$ is equal to
$2m+e$, for some $e\in \N$ with $0 \leqslant e \leqslant d$. Denote:

\smallskip

\begin{itemize}

\item[$\bullet$]
$\mathcal{ O}_{ 2m+e} \subset
M$ the union of CR orbits of dimension $= 2m+e$;

\smallskip

\item[$\bullet$]
$\mathcal{ O}_{ 2m+e}^\geqslant \subset
M$ the union of CR orbits of dimension $\geqslant
2m+e$;

\smallskip

\item[$\bullet$]
$\mathcal{ O}_{ 2m+e}^\leqslant \subset
M$ the union of CR orbits of dimension $\leqslant
2m+e$.

\end{itemize}

\smallskip
The function $p \mapsto \dim \mathcal{ O}_{ CR} (M, p)$ is
lower semicontinuous. It follows that 
$\mathcal{ O}_{ 2m+e}^\geqslant$ is open in $M$ and that
$\mathcal{ O}_{ 2m+e}^\leqslant$ is closed in $M$. 

Let $p\in M$, let $\mathcal{ O}_p$ be a small piece of $\mathcal{ O}_{
CR} (M, p)$ passing through $p$, of dimension $2m+e_p$, for some
integer $e_p$ with $0 \leqslant e_p \leqslant d$, and let $H_p$ be a
local $\mathcal{ C}^\infty$ or $\mathcal{ C}^{\kappa, \alpha}$
submanifold of $M$ passing through $p$ and satisfying $T_p H_p \oplus
T_p \mathcal{ O}_p = T_p M$. Call $H_p$ a local {\sl
orbit-transversal}. Implicitly, $H_p = \emptyset$ if $e_p = d$. Then,
in a sufficiently small neighborhood of $p$:

\smallskip

\begin{itemize}

\item[$\bullet$]
lower semi-continuity:
$H_p \cap \mathcal{ O}_{ 2m + e_p + 1}^\geqslant = \emptyset$;

\smallskip

\item[$\bullet$]
equivalently: $H_p \cap \mathcal{ O}_{ 2m+ e_p}^\geqslant = H_p
\cap \mathcal{ O}_{ 2m+e_p}$;

\smallskip

\item[$\bullet$]
$H_p \cap \mathcal{ O}_{ 2m+ e_p}^\leqslant$ is closed.

\end{itemize}

\smallskip

\def\theproposition{4.12}\begin{proposition}
If $M$ is $\mathcal{ C}^\omega$, then at every point $p\in M$, for
every orbit-transversal $H_p$ passing through $p$, the closed set $H_p
\cap \mathcal{ O}_{ 2m + e_p}^\leqslant$ is a local real analytic
subset of $H_p$ containing $p$.
\end{proposition}

A {\sl CR-invariant} subset of $M$ is a union of CR orbits. A closed
(for the topology of $M$) CR-invariant subset is {\sl minimal} if it
does not contain any proper subset which is also a closed CR-invariant
subset.

\def\theproblem{4.13}\begin{problem}
Describe the possible structure of the decomposition of $M$ in CR
orbits.
\end{problem}

There are differences between embedded and locally embeddable generic
submanifolds, which we shall not discuss, assuming that $M$ is
embedded in $\C^n$ or in $P_n (\C)$. Also, the $\mathcal{ C }^\omega$
category enjoys special features.

Indeed, if $M$ is a connected real analytic hypersurface,
Proposition~4.12 entails that each minimal closed invariant subset of
$M$ is either an embedded complex hypersurface or an open orbit; also
if $M$ contains at least one CR orbit of maximal dimension $(2n-1)$
(hence an open subset of $M$), all its CR orbits of codimension one
are complex $(n - 1)$-dimensional {\it embedded}\, submanifold of $M$
(a real analytic subset of codimension one in $\R$ consists of
isolated points). In the smooth category things are different.

So, let $M$ be a {\it connected}\, $\mathcal{ C }^\infty$ or
$\mathcal{ C}^{ \kappa, \alpha}$ ($\kappa \geqslant 2, 0 \leqslant
\alpha \leqslant 1$) hypersurface of $\C^n$. Its CR-orbits are either
$(2n - 1)$-dimensional, {\it i.e.} open in $M$, or $(2n -
2)$-dimensional and $T^c M$-integral, hence complex $(n-
1)$-dimensional hypersurfaces immersed in $M$.

\def\theproposition{4.14}\begin{proposition}
{\rm (\cite{ jo1999a})} In the smooth category, the structure of every
minimal closed CR-invariant subset $C$ of $M$ has one of the following
types{\rm :}

\smallskip

\begin{itemize}

\item[{\bf (i)}]
$C = M$ consists of a single embedded open CR orbit{\rm ;}

\smallskip

\item[{\bf (ii)}]
$C = \bigcup_{ a \in A} \mathcal{ O}_{CR, a} = M$ is a union 
of complex hypersurfaces, each being dense in $C$, with ${\rm
Card}\, A = {\rm Card}\, \R${\rm ;}

\smallskip

\item[{\bf (iii)}]
$C = \bigcup_{ a \in A} \mathcal{ O}_{CR, a}$ has empty interior in
$M$ and is a union union of complex hypersurfaces, each
being dense in $C$, with ${\rm Card}\, A = {\rm Card}\, \R${\rm ;}

\smallskip

\item[{\bf (iv)}]
$C$ consists of a single complex hypersurface embedded in $M$.
\smallskip

\end{itemize}

Furthermore, the closure, with respect
to the topology of $M$, of every CR orbit of dimension $(2n-2)$
is a minimal closed CR-invariant subset of $M$.

\end{proposition}

These four options are known in foliation theory (\cite{ hh1983,
cln1985}). One has to remind that each CR orbit contained in $C$ is
dense in $C$. In the first two cases, the trace of $C$ on any
orbit-transversal is a full open segment; in the third, it is a Cantor
set; in the last, it is an isolated point. In the third case,
impossible if $M$ is real analytic, $C$ will be called an {\sl
exceptional minimal CR-invariant subset}, similarly as in foliation
theory. We shall see below that compactness of $M \subset \C^n$
imposes a strong restriction on the possible $C$'s.

We mention that an analog of Proposition~4.14 holds for connected
generic submanifolds of codimension $d\geqslant 2$, provided one puts
the restrictive assumption that all its CR orbits are of codimension
$\leqslant 1$, the only difference being that CR orbits of codimension
$1$ are not complex manifolds in this case.

The presence of CR orbits of codimension $\geqslant 2$ in $M$ may
produce minimal closed CR-invariant subsets with complicated
transversal structure, even in the real analytic category (\cite{
bm1997}). Also, in a bounded strongly pseudoconvex boundary ({\it
see}~\S1.15(V) for background), there may exist a CR orbit of
codimension one whose closure constitutes an exceptional minimal
CR-invariant subset.

\def\thetheorem{4.15}\begin{theorem}
{\rm (\cite{ jo1999a})}
There exists a bounded strongly pseudoconvex domain $\Omega \subset
\C^3$ with $\mathcal{ C }^\infty$ boundary and a compact $\mathcal{
C}^\infty$ submanifold $M \subset \partial \Omega$ of dimension $4$
which is generic in $\C^3$ such that{\rm :}

\smallskip

\begin{itemize}

\item[$\bullet$]
$M$ is $\mathcal{ C}^\infty$ foliated
by CR orbits of dimension $3${\rm ;}

\smallskip

\item[$\bullet$]
$M$ contains a compact exceptional minimal 
CR-invariant set, but no compact CR orbit.

\end{itemize}

\smallskip

\end{theorem}

\proof[Summarized proof] The main idea is to start with an example due
to Sacksteder, known in foliation theory, of a compact real analytic
$3$-dimensional manifold $\mathcal{ N}$ equipped with a $\mathcal{
C}^\omega$ foliation $\mathcal{ F}$ of codimension one which carries a
compact exceptional minimal set, but no compact leaf. According to
\cite{ hh1983}, there exists such a pair $(\mathcal{ N}, \mathcal{
F})$, together with a $\mathcal{ C}^\infty$ diffeomorphism $\varphi_1
: \mathcal{ N} \to B \times S^1$, where $B$ is some compact orientable
$\mathcal{ C}^\infty$ surface of genus $2$ embedded in $\R^3$, and
where $S^1$ is the unit circle. Let $B \ni b \mapsto {\bf n} (b) \in
T_b \R^3$ denote the $\mathcal{ C}^\infty$ unit outward {\it normal}\,
vector field to such a $B \subset \R^3$, and consider $\R^3$ 
to be embedded in $\C^3$. For $\varepsilon >0$ small, the map
\[
\varphi_2:\ 
B\times S^1\ni
(b,\zeta)
\longmapsto
b+{\bf n}(b)\cdot\varepsilon\zeta
\in\C^3
\]
may be seen to be a totally real $\mathcal{ C}^\infty$ embedding. By
results of Bruhat-Whitney and Grauert, one may approximate the
$\mathcal{ C}^\infty$ totally real embedding $\varphi_2 \circ
\varphi_1$ by a $\mathcal{ C}^\omega$ embedding $\varphi : \mathcal{
N} \to \C^3$ which is arbitrarily close to $\varphi_2 \circ \varphi_1$
in $\mathcal{ C}^1$ norm, hence totally real. Denote $N := \varphi
(\mathcal{ N})$. The transported foliation $F := \varphi_* (\mathcal{
F})$ being real analytic, one may then proceed to an intrinsic
complexification of all its totally real $2$-dimensional leaves,
getting some $5$-dimensional real analytic hypersurface $N^{ i_c}$
containing $N$, equipped with a foliation $F^{ i_c}$ of $N^{ i_c}$ by
$2$-dimensional complex manifolds, with $F^{ i_c} \cap N = F$. This
foliation $F^{ i_c}$ is closed in some tubular neighborhood $\Omega$
of $N$ in $\C^3$, say $\Omega := \{ z \in \C^3 : \, {\rm dist} (z, N)
< \delta \}$, with $\delta >0$ small. Since $N$ is totally real, the
boundary $\partial \Omega$ is strongly pseudoconvex (Grauert) and is
$\mathcal{ C}^\infty$. Clearly, the intersection $M := N^{ i_c} \cap
\partial \Omega$ is a $4$-dimensional $\mathcal{ C}^\infty$
submanifold. The intersections of the $2$-dimensional complex leaves
of $F^{ i_c}$ with $\partial \Omega$ show that $M$ is foliated by
strongly pseudoconvex $3$-dimensional boundaries, which obviously
consist of a single CR orbit. Thus a CR orbit of $M$ is just the
intersection of a global leaf of $F^{ i_c}$ with $\partial \Omega$. In
conclusion, letting ${\rm Exc}_{ \mathcal{ F}}$ be the minimal
exceptional set of Sacksteder's example, $M$ contains the exceptional
minimal CR-invariant set $\left[ \varphi_* ({\rm Exc}_{ \mathcal{ F}}
) \right]^{ i_c} \cap \partial \Omega$ and no compact CR orbit.
\endproof

\subsection*{ 4.16.~Global minimality and laminations by complex 
manifolds} The CR orbits being essentially independent bricks, it is
natural to define the class of CR manifolds which consist of only one
brick.

\def\thedefinition{4.17}\begin{definition}{\rm
A $\mathcal{ C }^\omega$, $\mathcal{ C }^\infty$ or $\mathcal{ C}^{
\kappa, \alpha }$ CR manifold $M$ is called {\sl globally minimal} if
if consists of a single CR orbit.
}\end{definition}

Each CR orbit of a CR manifold is a globally minimal immersed CR
submanifold of $\C^n$. To simplify the overall presentation and not to
expose superficial corollaries, almost all the theorems of Parts~V 
and~VI in this memoir will be formulated with a global minimality
$M$.

\def\thelemma{4.18}\begin{lemma}
{\rm (\cite{ gr1968, jo1995, bch2005})} A compact
connected $\mathcal{ C}^2$ hypersurface in $\C^n$ is necessarily
globally minimal.
\end{lemma}

\proof 
Otherwise, the closure of a CR-orbit of codimension one in $M$
would produce a compact CR-invariant subset $C$ which is a union of
immersed complex hypersurfaces, each dense in $C$. Looking at a point
of $C$ where the pluriharmonic function $r_i := {\rm Re}\, z_i$ (or
$s_i := {\rm Im}\, z_i$) is maximal, the maximum principle entails
that $r_i$ (or $s_i$) is constant on $C$, for $i = 1, \dots, n$, hence
$C = \{ pt.\}$, contradiction. 
\endproof

More generally, the same simple argument yields:

\def\thecorollary{4.19}\begin{corollary}
Any Stein manifold cannot contain a compact set which is laminated by
complex manifolds
\end{corollary}

In the projective space $P_n (\C)$, one expects compact orientable
connected $\mathcal{ C}^2$ hypersurfaces $M$ to be still globally
minimal, but arguments are far to be simple. In fact, there are
infinitely many compact projective algebraic complex hypersurfaces
$\Sigma$ in $P_n (\C)$. However, they cannot be contained in such an
$M \subset P_n (\C)$ since, otherwise, their complex normal bundle $T
P_n (\C)\vert_{\Sigma} / T \Sigma$, known to be never trivial, but
equal to the complexification of the trvial bundle $TM\vert_\Sigma /
T\Sigma$, would be trivialized.

Thus, according to Proposition~4.14 above, the very question is about
nonexistence of closed CR-invariant sets $C \subset M$ laminated by
complex hypersurfaces which either coincide with $M$ or are
transversally Cantor sets. If $M \subset P_n (\C)$ is real analytic,
it might only be Levi-flat.

Nonexistence of orientable Levi-flat hypersurfaces in $P_n (\C)$ was
expected, because they would divide the projective space in two
smoothly bounded pseudoconvex domains. In the real analytic case,
non-existence was verified by Lins-Neto for $n\geqslant 3$ and by
Ohsawa for $n=2$; in the smooth (harder) case, {\it see}~\cite{
si2000}.

\def\theopenquestion{4.20}\begin{openquestion}
Is any compact orientable connected $\mathcal{ C}^2$ 
hypersurface of $P_n (\C)$ globally minimal\,?
\end{openquestion}

So, the expected answer is yes. In fact, the question is a particular
case of a deep conjecture stemming from Hilbert's sixteen problem
about phase diagrams of vector fields having polynomial coefficients
on the two-dimensional projective space. This conjecture is inspired
by the Poincar\'e-Bendixson theory valid over the real numbers,
according to which the closure of each leaf of such a foliation over
$P_2 (\R)$ contains either a compact leaf or a singularity. In its
most general form, it says that $P_n (\C)$ cannot contain a compact
set laminated by $(n-1)$-dimensional complex manifold, unless it is a
trivial lamination, viz just a compact complex projective algebraic
hypersurface; however, nontrivial laminations by $(n-2)$-dimensional
complex manifolds may be shown to exist.

A related general open question is to find topologico-geometrical
criteria on open subsets of $P_n (\C)$ insuring the existence of
nonconstant holomorphic functions there.

\subsection*{ 4.21.~Finite-typisation of generic submanifolds}
Let $M$ be a connected $\mathcal{ C}^\kappa$ ($2 \leqslant \kappa
\leqslant \infty$) generic submanifold of $\C^n$ of codimension $d
\geqslant 1$ and of CR dimension $m = n - d \geqslant 1$. The
distribution of subspaces $p \mapsto T_p^cM$ of $TM$ is of constant
rank $2m$. We apply to the complex tangential bundle $T^c M$ the
concept of finite-type.

\def\thedefinition{4.22}\begin{definition}{\rm
A point $p\in M$ is said to be of {\sl finite type} if the system
$\LL$ of local sections of $T^c M$ defined in a neighborhood of $p$
satisfies $\LL^\kappa ( p) = T_p M$. The smallest integer $d( p)
\leqslant \kappa$ with $\LL^{ d( p)} ( p) = T_p M$ is called the {\sl
type} of $M$ at $p$.
}\end{definition}

We want now to figure out how, in general, a generic submanifold of
$\C^n$ must be globally minimal and in fact, of finite type at every
point. We equip with the strong Whitney topology the set ${}^\kappa
\mathcal{ G}_{d,m}^n$ of $\mathcal{ C }^\kappa$ ($2\leqslant \kappa
\leqslant \infty$) connected generic submanifolds $M \subset \C^n$ of
codimension $d\geqslant 1$ and of CR dimension $m = n - d \geqslant
1$. No restriction is made on the global topology.

As a model case, let $\kappa \geqslant 2$ and
consider $M$ to be rigid algebraic represented by
\[
w_j
=
\bar w_j
+
i\,P_j(z,\bar z)
=
\bar w_j
+
i\,\sum_{\vert\alpha\vert+\vert\beta\vert\leqslant\kappa}\,
p_{j,\alpha,\beta}\,
z^\alpha\bar z^\beta,
\]
where $\alpha, \, \beta\in \N^m$, 
where the polynomials $P_j$ are real,
$p_{ j, \alpha, \beta} = 
\overline{ p}_{ j, \beta, \alpha}$ and have no pluriharmonic term,
namely $0 \equiv P_j (z, 0 ) \equiv P_j (0, \bar z)$, for 
$j=1, \dots, d$. A basis of
$(1, 0)$ and of $(0, 1)$ vector fields is given by
\[
L_k
:=
\frac{\partial}{\partial z_k}
+
i\,\sum_{j=1}^d\,
P_{j,z_k}\,\frac{\partial}{\partial w_j}
\ \ \ \ \
\text{\rm and}
\ \ \ \ \
\overline{L}_k
:=
\frac{\partial}{\partial \bar z_k}
-
i\,\sum_{j=1}^d\,
\overline{P}_{j,\bar z_k}\,\frac{\partial}{\partial \bar w_j},
\]
for $k=1, \dots, n$. The Lie algebra $\LL^\kappa$ generated by all Lie
brackets of length $\leqslant \kappa$ of the system $\LL := \{ L_k,
\overline{ L}_k \}_{ 1\leqslant
\kappa \leqslant m}$ contains the subalgebra $\LL_{ CR, rigid}^\kappa$
generated by the only brackets of the form
\[
\left[
L_{\lambda_1},\dots,
\left[
L_{\lambda_a},
\left[
\overline{L}_{\mu_1},
\dots, 
\left[
\overline{L}_{\mu_b},
\left[
L_{k_1},\,\overline{L}_{k_2}
\right]
\right]
\dots
\right]
\right]
\dots
\right]
\]
where $2 + a + b \leqslant \kappa$ and
where $1 \leqslant \lambda_1, 
\dots, \lambda_a, \mu_1, \dots, \mu_b, 
k_1, k_2 \leqslant m$.
One computes 
\[
\left[
L_{k_1},\,\overline{L}_{k_2}
\right]
=
-i
\left(
\sum_{j=1}^d\,
P_{j,z_{k_1}\bar z_{k_2}}\,
\frac{\partial}{\partial w_j}
+
\sum_{j=1}^d\,
\overline{P}_{j,z_{k_1}\bar z_{k_2}}\,
\frac{\partial}{\partial \bar w_j}
\right),
\]
hence the above iterated Lie brackets are equal to
\[
-i
\left(
\sum_{j=1}^d\,
P_{j,z_{\lambda_1}\dots z_{\lambda_a}\bar z_{\mu_1}\dots\bar z_{\mu_b}
z_{k_1}\bar z_{k_2}}\,
\frac{\partial}{\partial w_j}
+
\sum_{j=1}^d\,
\overline{P}_{j,z_{\lambda_1}\dots z_{\lambda_a}\bar z_{\mu_1}\dots
\bar z_{\mu_b}z_{k_1}\bar z_{k_2}}\,
\frac{\partial}{\partial \bar w_j}
\right).
\]
This shows that all these brackets are linearly independent 
as functions of the
jets of the $P_j$. In fact, the number of such brackets is exactly
equal to the dimension of the space of polynomials $P (z, \bar z)$ of
degree $\leqslant \kappa$ having no pluriharmonic term, namely equal
to
\[
\frac{(2m+\kappa)!}{2m!\ \kappa!}
-
2\,\frac{(m+\kappa)!}{m!\ \kappa!}
+
1.
\]
For a general $\mathcal{ C}^\kappa$ submanifold $M$ (not necessarily
rigid), one verifies that the same collection of brackets is
independent in terms of the jets of the defining equation of $M$.
Generalizing slightly Lemma~2.13, we see that in the vector space of
$d \times (d + e)$ (real or complex) matrices, the subset of matrices
of rank $\leqslant d - 1$ is a real algebraic set of codimension equal
to $(e + 1)$. If we choose $\kappa$ large enough so that the dimension
of $\LL_{ CR, rigid}^\kappa$ is $\geqslant d + \dim M = 2 ( m + d ) =
2n$ (applying the previous assertion with $e := \dim M$), if we form
the $d \times (d + e')$, $e' \geqslant e$, matrix consisting of the
coordinates of the $\frac{ \partial }{\partial w}$-part of brackets of
length $\leqslant \kappa$ as above, then the set where this matrix is
of rank $\leqslant d-1$ is of codimension $\geqslant \dim M + 1$ in
the space of $\kappa$-th jets of the defining equations of
$M$. Consequently, the jet transversality theorem applies.

\def\thetheorem{4.23}\begin{theorem}
Let $n\geqslant 1$, $m\geqslant 1$ and $d\geqslant 1$ be integers
satisfying $m+d = n$ and let $\kappa$ be the minimal integer having
the property that
\[
\frac{(2m+\kappa)!}{2m!\ \kappa!}
-
2\,\frac{(m+\kappa)!}{m!\ \kappa!}
+
1
\geqslant
2(m+d)
=
2n.
\]
Then the set of $\mathcal{ C}^\kappa$ connected generic submanifolds
$M \subset \C^n$ of codimension $d$ and of CR dimension $m$ that are
of finite type $\leqslant \kappa$ at every point is open and dense in 
the set ${}^\kappa \mathcal{ G}_{ d, m}^n$ of all 
generic submanifolds.
\end{theorem}

In particular, a connected $\mathcal{ C}^4$ (resp. $\mathcal{ C}^3$,
resp. $\mathcal{ C}^2$) hypersurface in $\C^2$ (resp. in $\C^3$,
resp. in $\C^n$ for $n\geqslant 4$) is of finite type $4$ (resp. $3$,
resp. $2$, or equivalently, is not Levi-flat) at every point after an
arbitrarily small perturbation.

Similarly, if instead of the subalgebra $\LL_{ \text{\rm CR},
\text{\rm rigid }}$, one would have considered the (smaller)
subalgebra consisting of only the brackets $\big[ L_{\lambda_1},
\dots, \big[ L_{\lambda_a}, \big[ L_{k_1}, \, \overline{ L}_{ k_2}
\big] \big] \dots \big]$, where $2 + a \leqslant \kappa$ and where $1
\leqslant \lambda_1 , \dots, \lambda_a, k_1, k_2 \leqslant \kappa$,
one would have obtained finite type $\leqslant \kappa$, for $\kappa$
minimal satisfying $2m \, \frac{ (m + \kappa - 1)!}{ m! \ (\kappa -
1)!} \geqslant m^2 + 2 ( 2m + d)$. We also mention that the same
technique enables one to prove that, after an arbitrarily small
perturbation, $M$ is finitely nondegenerate at every point and of
finite nondegeneracy type $\leqslant \ell$, with $\ell$ minimal
satisfying $2 d \, \frac{ (\ell + m)!}{ \ell ! \ m!} \geqslant 4n -
1$. In particular, $\ell = 3$ when $m = d = 1$ while $\ell = 2$
suffices when $m = 1$ for all $d \geqslant 2$. Details are left to
the reader.

To conclude, we state the analog of Open question~2.17 for induced CR
structures.

\def\theopenquestion{4.24}\begin{openquestion}
{\rm (\cite{ js2004}, [$*$])}
Given a fixed generic submanifold $M$ of class $\mathcal{ C}^\kappa$
that is of finite type at every point, is it always possible to
perturb slightly a $\mathcal{ C}^\kappa$ submanifold $M_1$ of $M$ that
is generic in $\C^n$, of codimension $d_1 \geqslant 1$ and of CR
dimension $m_1 = n - d_1 \geqslant 1$, as a $\mathcal{ C}^\kappa$
submanifold $\widetilde{ M_1}$ of $M$ that is of finite type at every
point\,? If so, what is the smallest regularity $\kappa$ in terms of
$d$, $m$, $d_1$, $m_1$ and of the highest type at points of $M$\,?
\end{openquestion}

\subsection*{ 4.25.~Spaces of CR functions and of CR distributions}
A $\mathcal{ C}^1$ function $f: M \to \C$ is called {\sl
Cauchy-Riemann} ({\sl CR} briefly) if it is annihilated by every
section of $T^{ 0, 1} M$. Equivalently:

\smallskip

\begin{itemize}
\item[$\bullet$]
$df$ is $\C$-linear on $T^c M$;

\smallskip

\item[$\bullet$]
$df \wedge dt_1 \wedge \cdots \wedge dt_n \big\vert_M = 0$;

\smallskip

\item[$\bullet$]
$\int_M \, f \, \overline{ \partial} \omega = 0$, for every
$\mathcal{ C}^1$ form $\omega$ of type $( n, m-1)$ in $\C^n$
having compact support.

\end{itemize}

\smallskip

(Remind the local expression of $(r, s)$ forms: $\sum_{ I, J} \, a_{
I, J} \, dt^I \wedge d\bar t^J$, where $I = (i_1, \dots, i_r)$ and $J
= (j_1, \dots, j_s)$.) A (only) {\it continuous}\, function $f : M \to
\C$ is {\sl CR} if the last condition $\int_M \, f \overline{
\partial} \omega = 0$ holds. Further, Lebesgue-integrable CR
functions, CR measures, CR distributions and CR currents may be
defined as follows (\cite{ trv1981, trv1992, hm1998,
trp1996, jo1999b}).

Thanks to graphing functions, one may equip locally $M$ with some (in
fact many) volume form, or equivalently, some deformation of the
canonical $\dim M$-dimensional Legesgue measure defined on tangent
spaces. Let ${\sf p}$ be a real number with $1\leqslant {\sf p}
\leqslant \infty$. Since two such measures are multiple of each other,
it makes sense to speak of $L_{ loc}^{ \sf p}$ functions $M \to
\C$. In this setting, a $L_{ loc}^{ \sf p} (M)$ function $f$ is CR if
$\int_M \, f \, \overline{ \partial} \omega = 0$, for every $\mathcal{
C}^1$ form $\omega$ of type $(n, m-1)$ in $\C^n$ having compact
support.

A distribution ${\sf T}$ on $M$ is {\sl CR} if for every section
$\overline{ L}$ of $T^{ 0, 1} M$ defined in an open subset $U \subset
M$ and every $\chi \in \mathcal{ C}_c^\infty (U, \C)$, one has $\left<
{\sf T}, \overline{ L} (\chi) \right> = 0$.

A CR distribution of order zero on $M$ is called
a {\sl CR measure}. Equivalently, a CR measure is 
a continuous linear map $\omega' \mapsto 
\mu (\omega ')$ from compactly supported forms on $M$
of maximal degree $2m + d$ to $\C$, that is
CR in the weak sense, namely
$\mu (\overline{ \partial} \omega) = 0$, 
for every $\mathcal{ C}^1$ form $\omega$ of type $(n, m-1)$
having compact support.
Once a volume form ${\rm dvol}_M$ is fixed on $M$, 
the quantity $\mu\, {\rm dvol}_M$ is a CR (Borel) measure on $M$.

\subsection*{ 4.26.~Traces of CR functions on CR orbits}
A $\mathcal{ C}^1$ function $f : M \to \C$ is CR on $M$ if and
only if its restriction to every CR orbit of $M$ is CR (obvious). If
$f$ is $\mathcal{ C}^0$ or $L_{loc}^{\rm p}$, a similar but nontrivial
statement holds. By ``{\sl 
almost every CR orbit}'', we shall mean ``except a
union of CR orbits whose $\dim M$-dimensional measure vanishes''.

\def\thetheorem{4.27}\begin{theorem} {\rm ($d=1$: \cite{ jo1999b};
$d\geqslant 1$: \cite{ po1997, mp1999})}
Assume that $M$ is at least $\mathcal{ C}^3$ and let $f$ be a function
in $L_{loc}^{\sf p}(M)$ with $1\leqslant {\sf p} \leqslant \infty$. Then the
restriction $f\vert_{\mathcal{ O}_{ CR}}$ is in $L_{loc}^{\sf p}$ on
$\mathcal{ O}_{ CR}$, for almost every $\mathcal{ O}_{ CR}$.
Furthermore, $f$ is CR if and only if, for almost every CR orbit
$\mathcal{ O }_{ CR}$ of $M$, its restriction $f \vert_{\mathcal{
O}_{CR}}$ is CR.
\end{theorem}

The theorem also holds for $f$ continuous, with $f\vert_{\mathcal{
O}_{CR}}$ being CR for {\it every} CR orbit. Here, $\mathcal{ C
}^3$-smoothness is needed. Property {\bf (5)} of Sussman's orbit
Theorem~1.21 together with a topological reasoning yields
a covering by orbit-chart which is used in the proof.

\def\theproposition{4.28}\begin{proposition}
{\rm (\cite{ jo1999a, po1997, mp1999})} Assume $M$ is $\mathcal{
C}^\infty$ or $\mathcal{ C }^{\kappa, \alpha }$, with $\kappa
\geqslant 2$, $0 \leqslant \alpha \leqslant 1$ and let $\square := \{
x \in \R : \vert x \vert < 1\}$. There exists a countable covering
$\bigcup_{ k\in \N} U_k = M$ such that for each $k$, there exist $e_k
\in \N$ with $0 \leqslant e_k \leqslant d$ and a $\mathcal{ C}^{
\kappa - 1, \alpha}$ diffeomorphism{\rm :}
\[
\varphi_k : 
(s_k, t_k) \ni
\square^{ 2m+e_k}\times \square^{ d - e_k}
\longmapsto 
\varphi_k (s_k, t_k)\in U_k,
\]
such that{\rm :}

\smallskip

\begin{itemize}

\item[$\bullet$]
$\varphi_k \left( \square^{ 2m + e_k} \times \{ t_k^* \} \right)$ is
contained in a single CR orbit, for every fixed $t_k^* \in \square^{
d-e_k}${\rm ;}

\smallskip

\item[$\bullet$]
for each $p\in M$, there exists $k = k_p \in \N$ with $p \in U_{ k_p}$,
viz there exist $s_{ k_p, p}$ and $t_{ k_p, p }$ with $\varphi_{ k_p} (
s_{ k_p, p}, t_{ k_p, p}) = p$, such that $\varphi_{ k_p} \left(
\square^{ 2m + e_{ k_p}} \times \{ t_{ k_p, p}\} \right)$ is an open
piece of the CR orbit of $p$, {\it i.e.} $\dim \mathcal{ O}_{ CR} (M,
p) = 2m + e_{k_p}$.

\end{itemize}

\smallskip

\end{proposition}

In the proof of the theorem, $\mathcal{ C }^2$-smoothness of the maps
$\varphi_k$ (hence $\mathcal{ C}^3$-smoothness of $M$) is required to
insure that the pull-back $\varphi_k^* (T^c M \vert_{ U_k })$ is
$\mathcal{ C }^1$. However, we would like to mention that if $M$ is
$\mathcal{ C}^{2, \alpha}$ with $0 < \alpha < 1$ results of \cite{
tu1990, tu1994a, tu1996} and Theorem~3.7(IV) insuring the $\mathcal{
C}^{ 2, \beta}$-smoothness of local and global CR orbits, for every
$\beta < \alpha$, this would yield orbit-charts $\varphi_k$ of class
$\mathcal{ C}^{ 2, \beta}$, and then the above theorem holds true with
$M$ of class $\mathcal{ C}^{ 2, \alpha}$.

\subsection*{ 4.29.~Boundary values of holomorphic functions for 
functional spaces $\mathcal{ C}^{\kappa, \alpha}$, $\mathcal{ D}'$,
$L_{loc }^{\sf p}$} Let $M$ be a generic submanifold of $\C^n$ of
codimention $d \geqslant 1$ and of nonnegative CR dimension $m
\geqslant 0$ (we admit $m=0$). Assume $M$ is at least $\mathcal{
C}^1$. In appropriate coordinates $t = (z, w) = (x+ iy, u+ iv) \in
\C^n \times \C^m$ centered at one of its points $p$:
\[
M 
=\left\{
(z,w)\in\Delta_{\rho_1}^m
\times
\left(
\square_{\rho_1}^d\times i\R^d
\right):\
v=\varphi(x,y,u)
\right\},
\]
for some $\rho_1 > 0$, with $\varphi (0) = 0$ and $d\varphi (0) = 0$.
Let $\rho$ be a real number with $0 \leqslant \rho \leqslant \rho_1$.
The height function:
\[
\sigma(\rho)
:=
\max_{\vert x\vert,\vert y\vert,\vert u\vert\leqslant\rho}
\left\vert\varphi(x,y,u)\right\vert
\]
is continuous and tends to $0$, as $\rho$ tends to $0$. For every
$\rho\leqslant \rho_1$ and every $\sigma > \sigma (\rho)$, the boundary of
$M \cap \left[ \Delta_\rho^m \times \left( \square_\rho^d \times i
\square_\sigma^d \right) \right]$ is contained in the boundary
$\partial \left( \Delta_\rho^d \times \square_\rho^d \right)$ of the
horizontal space, times the vertical space $i \square_\sigma^d$.

Let $C$ be an open convex cone in $\R^d$ having vertex $0$. We shall
assume it to be {\sl salient}, namely contained in one side of some
hyperplane passing through the origin. Equivalently, its intersection
$C \cap S^{ d-1}$ with the unit sphere of $\R^d$ is open, contained in
some open hemisphere and convex in the sense of spherical geometry.

A {\sl local wedge of edge $M$ at $p$ directed by $C$} is an open set
of the form:
\def\theequation{4.30}\begin{equation}
\aligned
\mathcal{W}
=
\mathcal{W}(\rho,\sigma,C)
:=
\big\{
&
(x+iy,u+iv)
\in\Delta_\rho^m
\times
\square_\rho^d
\times
i\square_\sigma^d:
\\
& \
v-\varphi (x,y,u)\in C\big\},
\endaligned
\end{equation}
for some $\rho, \sigma >0$ satisfying $\rho \leqslant \rho_1$ and $\sigma >
\sigma (\rho)$. This type of open set is independent of the choice of
local coordinates and of local defining functions; in codimension
$d\geqslant 2$, it generalizes the notion of local side of a
hypersurface. Notice that $\mathcal{ W}$ is connected.

If there exists a
function $F$ that
is holomorphic in $\mathcal{ W}$ and
that extends continuously up to the
{\sl edge}
\[
M_\rho
:=
M
\cap \left[ \Delta_\rho^m 
\times \left( \square_\rho^d \times i
\R^d\right)\right]
\]
{\sl of the wedge} $\mathcal{ W}$, then the limiting values of $F$
define a continuous CR function on $M_\rho$.

A more general phenomenon holds. A function $F$, holomorphic in the
wedge $\mathcal{ W}$, has {\sl slow growth up to $M$}, if there exist
$k\in \N$ and $C>0$ such that
\[
\vert F(t)\vert
\leqslant
C\left\vert v-\varphi(x,y,u)\right\vert^{-k},
\ \ \ \ \ \ \ \ 
t=(x+iy,u+iv)\in\mathcal{ W}.
\]
Equivalently, $\vert F(t) \vert \leqslant C\left[ {\rm dist}(t,M)
\right]^{-k}$, with the same $k$ but a possibly different $C$. 
As in the cited references, we shall assume
$M$ to be $\mathcal{ C }^\infty$.

\def\thetheorem{4.31}\begin{theorem}
{\rm (\cite{ bct1983, ho1985, br1987,
ber1999})} If $F \in \mathcal{ O} \left( \mathcal{ W} (\rho,
\sigma, C) \right)$ has slow growth up to $M$, it possesses a boundary
value $b_M F$ which is a CR distribution on the edge $M \cap \left[
\Delta_\rho^m \times \left( \square_\rho^d \times i \square_\sigma^d
\right) \right]$ precisely defined by{\rm :}
\[
\aligned
\left<
{\sf b}_MF,\chi
\right>
:=
\lim_{C\ni\theta\to 0}\,
\int_{\Delta_\rho^m\times\square_\rho^d}\
F
&
\left(
x+iy,
u+i\varphi(x,y,u)+i\theta
\right)
\cdot
\\
&
\ \ \ \ \ \ \
\cdot
\chi(x,y,u)\,
dx\,dy\,du,
\endaligned
\]
where $\chi = \chi (x, y, u)$ is a $\mathcal{ C}^\infty$ function
having compact support in $\Delta_\rho^m \times \square_\rho^d$. 

\smallskip

\begin{itemize}

\item[{\bf (i)}]
The limit is independent of the way how $\theta \in C$ approaches $0
\in \R^d$.

\smallskip

\item[{\bf (ii)}]
If ${\sf b}_M F$ is $\mathcal{ C}^{ \lambda, \beta}$, $\lambda
\geqslant 0$, $0 \leqslant \beta \leqslant 1$, then $F$ extends as a
$\mathcal{ C}^{\lambda, \beta}$ function on $\mathcal{ W}' \cup \left(
M \cap \left[ \Delta_\rho^m \times \left( \square_\rho^d \times i
\square_{ \sigma'}^d \right) \right] \right)$, for every wedge
$\mathcal{ W}' = \mathcal{ W}' (\rho, \sigma', C')$ with $\sigma (
\rho ) < \sigma' \leqslant \sigma$ and with $C' \cap S^{ d-1} \subset
\subset C \cap S^{ d-1}$.

\smallskip

\item[{\bf (iii)}]
Finally, $F$ vanishes identically in the wege $\mathcal{ W}$ if
and only if ${\sf b}_M F$ vanishes on some nonempty open subset of
the edge $M_\rho$.
\end{itemize}
\end{theorem}

The integration is performed on the translation $M_\rho^\theta :=
M_\rho + (0, i\, \theta)$, drawn as follows.

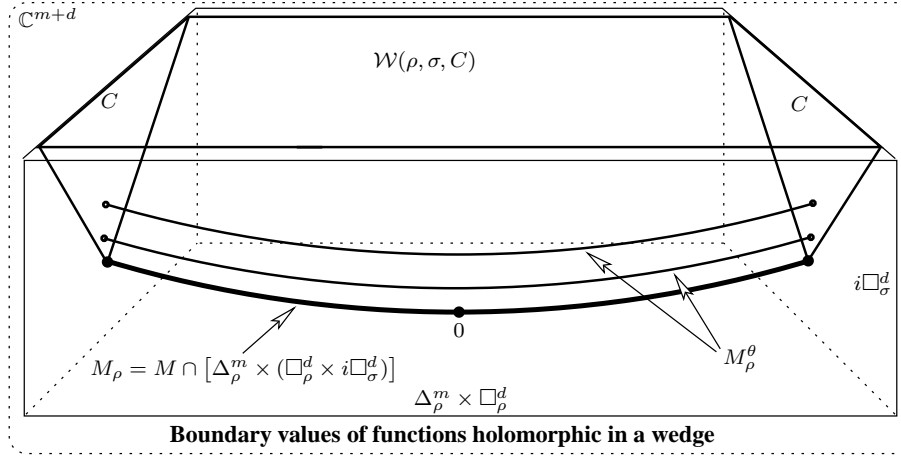
\begin{center}
\input bv.pstex_t
\end{center}

The proof is standard for $M \equiv \R^n$ (\cite{ ho1985}), the main
argument going back to Hadamard's finite parts. With technical
adaptations in the case of a general generic $M$, several integrations
by part are performed on a thin $(\dim M + 1)$-dimensional cycle
delimited by $M_\rho^0$ and $M_\rho^\theta$, taking advantage of
Cauchy's classical formula, until the rate of explosion of $F$ up to
the edge is dominated. The uniqueness property 
{\bf (iii)} requires analytic disc
methods (Part~V).

\smallskip

Boundary values in the $L^{\sf p}$ sense requires special
attention. At first, remind that a function $F$ holomorphic in the
unit disc $\Delta$ {\sl belongs to the Hardy class} $H^{ \sf p}
(\Delta)$ if the supremum:
\[
\left\vert\!\left\vert
F
\right\vert\!\right\vert_{H^{\sf p}(\Delta)}
:=
\sup_{0<r<1}
\left(
\int_{-\pi}^\pi\,
\left\vert
F(re^{it})
\right\vert^{\sf p}
\right)^{1/{\sf p}}
<
\infty
\]
is finite. According to Fatou and Privalov, such a function $F$ has
radial boundary values $f(e^{it}) := \lim_{ r\underset{ < }{
\rightarrow }1} F(re^{ it})$, for almost every $t \in [-\pi, \pi]$, so
that the boundary value $f$ belongs to $L^{\sf p} ( [ -\pi,
\pi])$. Furthermore, if $1 \leqslant {\sf p} < \infty$:
\[
\lim_{ r\underset{ < }{
\rightarrow }1}\,
\int_{-\pi}^\pi
\left\vert
F(re^{it})-f(e^{it})
\right\vert^{\sf p}
=
0.
\]
Consider a bounded domain $D \subset \C^n$ having boundary
of class at least $\mathcal{ C}^2$, defined by $D = \{ z \in
\C^n : \rho (z) < 0 \}$, with $\rho \in \mathcal{ C }^2$
satisfying $d\rho \neq 0$ on $\partial D$. For $\varepsilon >0$
small, let $D_\varepsilon := \{ z \in D : \rho (z) < - \varepsilon
\}$. The induced Euclidean measure on $\partial D_\varepsilon$
(resp. $\partial D$) is denoted by $d \sigma_\varepsilon$
(resp. $d\sigma$). Then the {\sl Hardy space} $H^{ \sf p} (D)$
consists of holomorphic functions $F \in \mathcal{ O} (D)$ having the
property that the supremum:
\[
\left\vert\!\left\vert
F
\right\vert\!\right\vert_{H^{\sf p}(D)}
:=
\sup_{\varepsilon>0}
\left(
\int_{\partial D_\varepsilon}\,
\left\vert
F(z)
\right\vert^{\sf p}\,
d\sigma_\varepsilon(z)
\right)^{1/{\sf p}}
<\infty
\]
is finite. The resulting space does not depend on the choice of
a defining function $\rho$ (\cite{ st1972}). Let $\mathbf{ n }_z$ be
the outward-pointing normal to the boundary at $z \in \partial D$.

\def\thetheorem{4.32}\begin{theorem}
{\rm (\cite{ st1972})}
If $F \in H^{\sf p} (D)$, for almost all $z \in \partial D$, the
normal boundary value $f (z) := \lim_{ \varepsilon \underset{ > }{
\rightarrow }0 } \, F \left( z - \varepsilon \mathbf{ n}_z \right)$
exists and defines a function $f$ which belongs to $L^{\sf p} (\partial
D)$. Furthermore, if $1 \leqslant {\sf p} < \infty${\rm :}
\[
\lim_{\varepsilon\underset{>}{\rightarrow}0}\,
\int_{\partial D}\,
\left\vert
F(z-\varepsilon\,{\bf n}_z)-f(z)
\right\vert^{\sf p}\,d\sigma(z)
=0.
\]
\end{theorem}

In arbitrary codimension, the notion of $L^{\sf p}$ boundary values
may be defined in the local sense as follows. Let
$M$ be generic, let $p\in M$ and let $\mathcal{ W} = \mathcal{ W}
(\rho, \sigma, C)$ be a local wedge of edge $M$ at $p$, as defined
by~\thetag{ 4.30}. A holomorphic function $F \in \mathcal{ O}
(\mathcal{ W})$ {\sl belongs to the Hardy space} $H_{ loc}^{\sf p} (
\mathcal{ W} )$ if for every cone $C' \subset \R^d$ with $C ' \cap S^{
d-1} \subset \subset C \cap S^{ d-1}$ and every $\rho ' < \rho$, the
supremum:
\[
\sup_{\theta'\in C'}\,
\int_{\Delta_{\rho'}^m
\times\square_{\rho'}^d}\,
\left\vert
F\left(x+iy,u+i\varphi(x,y,u)+i\theta'\right)
\right\vert^{\sf p}\,
dx\wedge dy \wedge du \ 
< 
\ \infty
\]
is finite. Up to shrinking cubes, polydiscs and cones, the resulting
space neither depends on local coordinates nor on the choice of
local defining equations.

\def\thetheorem{4.33}\begin{theorem}
{\rm ($d=1$: \cite{ st1972, jo1999b}; $d\geqslant 2$: \cite{
po1997})} If $F \in H^{\sf p}_{ loc} (\mathcal{ W})$, for almost
$(x,y, u+i\, \varphi (x,y,u )) \in M_\rho$ and for every cone $C'$
with $C' \cap S^{ d-1} \subset \subset C \cap S^{ d-1}$, the boundary
value{\rm :}
\[
f(x,y,u):=\lim_{C'\ni\theta' \to 0}\, 
F\left(x+iy,u+i\,\varphi(x,y,y)
+i\,\theta'\right)
\]
exists and defines a function $f$ which belongs to $L_{loc,CR}^{\sf p}
(M_\rho)$. Furthermore, if $1\leqslant {\sf p} < \infty$, for every $\rho'
< \rho${\rm :}
\begin{small}
\[
\aligned
\lim_{C'\ni\theta'\to 0}\,
\int_{\Delta_{\rho'}^m\times\square_{\rho'}^d}\,
\big\vert
F\big(x+iy,
&
u+i\,\varphi(x,y,y)+i\,\theta'\big)
-
\\
&
-
f(x,y,u)
\big\vert^{\sf p}\,
dx\wedge dy \wedge du 
=0.
\endaligned
\]
\end{small}
\end{theorem}

\subsection*{4.34.~Holomorphic extendability of CR functions in 
$\mathcal{ C}^{ \kappa, \alpha}$, $\mathcal{ D }'$, $L_{loc }^{\sf
p}$} In Part~V, we will
study sufficient conditions for the existence of wedges to
which CR functions and distributions extend holomorphically.

\def\thedefinition{4.35}\begin{definition}{\rm
A CR function of class $\mathcal{ C}^{ \kappa, \alpha}$ or $L_{ loc
}^{\sf p}$ $(1\leqslant {\sf p} < \infty)$ or a CR distribution $f$
defined on $M$ is {\sl holomorphically extendable} if there exists a
local wedge $\mathcal{ W} = \mathcal{ W} (\rho, \sigma, C)$ at $p$ and
a holomorphic function $F \in \mathcal{ O} ( \mathcal{ W})$ whose
boundary value ${\sf b}_M F$ equals $f$ on $M_\rho$ in the 
$\mathcal{ C}^{ \kappa, \alpha}$, $L^{\sf p}$ or
$\mathcal{ D}'$ sense.
}\end{definition}

\subsection*{ 4.36.~Local CR distributions supported by a local CR orbit}
Assume now that $M$, of class $\mathcal{ C}^\infty$ and represented as
in \S4.3, is {\it not}\, locally minimal at $p$. Equivalently,
$\mathcal{ O}_{ CR }^{loc} (M, p)$ is of dimension $2 m + e \leqslant
2m + d -1$. In a small neighborhood, $S := \mathcal{ O }_{ CR}^{ loc}
(M, p)$ is a closed connected CR submanifold of $M$ passing through
$p$ and having the same CR dimension as $M$. There exist local
holomorphic coordinates $( z, w) = (z, w_1, w_2) \in \C^m \times \C^{
e} \times \C^{ d - e}$ vanishing at $p$ in which $M$ is represented by
$v = \varphi (x, y,u)$ and $S$ is represented by the supplementary
(scalar) equation(s) $u_2 = \lambda_2 ( x, y, u_1)$, with $\varphi$ and
$\lambda_2$ of class $\mathcal{ C }^\infty$ satisfying $\varphi ( 0) =
0$, $d\varphi ( 0) = 0$, $\lambda_2 ( 0) = 0$ and $d \lambda_2 (0) =
0$. According to Theorem~4.2, the assumption that $S$ is CR and has
the same CR dimension as $M$ may be expressed as follows.

\def\theproposition{4.37}\begin{proposition}
Decomposing $\varphi = ( \varphi_1 , \varphi_2 )$
and defining{\rm :}
\[
\aligned
v_1
&
=
\varphi_1\left(x,y,u_1,\lambda_2(x,y,u_1)\right)
=:
\mu_1(x,y,u_1),
\\
v_2
&
=
\varphi_2\left(x,y,u_1,\lambda_2(x,y,u_1)\right)
=:
\mu_2(x,y,u_1).
\endaligned
\]
the map{\rm :}
\[
\psi_2(x,y,u_1)
:=
\lambda_2(x,y,u_1)+i\mu_2(x,y,u_1)
\]
is CR on the generic submanifold $v_1 = \mu_1 ( x, y, u_1)$ of
$\C^m \times \C^e$.
\end{proposition}

In a small neighborhood $U$ of $p$, the restrictions 
\[
dz_1
\big\vert_S,\dots, 
dz_m 
\big\vert_S, 
dw_1 
\big\vert_S, 
\dots, dw_e
\big\vert_S
\] 
span an $(m+e)$-dimensional subbundle of $\C T^*
S$. Denoting $dz := dz_1 \wedge \cdots \wedge dz_m$, $d\bar z := d\bar
z_1 \wedge \cdots \wedge d\bar z_m$ and $dw' := dw_1 \wedge \cdots
\wedge dw_e$, for $\chi \in \mathcal{ C}_c^\infty ( U, \C)$, consider
the (localized) distribution defined by:
\[
\left<
{\sf [S]}, \chi
\right>
:=
\int_{U\cap S}\,
\chi\cdot dz \wedge dw'\wedge d\bar z.
\]

\def\theproposition{4.38}\begin{proposition}
{\rm (\cite{ trv1992, ht1993})}
Then ${\sf [S]}$ is a nonzero local CR measure supported by $S \cap U$.
\end{proposition}

\proof
It is clear that ${\sf [ S ]}$ is supported by $S \cap U$ and is of
order zero. Let the $(0,1)$ vector fields $\overline{ L}_k$ and the
complex-transversal ones $K_j$ be as in \S4.3. Reminding
$d\chi = \sum_{ k=1 }^m \, L_k ( \chi) \, dz_k + \sum_{ k=1 }^m \,
\overline{ L}_k (\chi) \, d\bar z_k + \sum_{ j=1 }^d \, K_j ( \chi) \,
dw_j \big\vert_M$, we observe:
\[
\overline{L}_k(\chi)\,dz\wedge dw'\wedge d\bar z
=
\pm d
\left(\chi\cdot
dz\wedge dw'\wedge d\bar z_1 \wedge
\cdots\wedge\widehat{d\bar z_k}
\wedge\cdots\wedge d\bar z_m
\right).
\]
Replacing this volume form in the integrand:
\[
\aligned
\left<
\overline{L}_k{\sf [S]},\chi
\right>
:=
&
-
\left<
{\sf [S]},\overline{L}_k(\chi)
\right>
=
-\int_{S\cap U}\,
\overline{L}_k(\chi)\,dz
\wedge dw'\wedge d\bar z
\\
=
&
\pm\int_{S\cap U}\,d
\left(\chi\cdot
dz\wedge dw'\wedge d\bar z_1 \wedge
\cdots\wedge\widehat{d\bar z_k}
\wedge\cdots\wedge d\bar z_m
\right)
\endaligned
\]
and applying Stokes' theorem, we deduce $\left<
\overline{L }_k {\sf [S]}, \chi \right> = 0$, {\it i.e.}
${\sf [S]}$ is CR.
\endproof 

The last assertion of Theorem~4.31 and the vanishing of $[ S]$ on the
dense open set $U\backslash (S\cap U)$ entails the following.

\def\thecorollary{4.39}\begin{corollary}
{\rm (\cite{ trv1992, ht1993})} The nonzero local CR measure
${\sf [ S ]}$ does not extend holomorphically to any local wedge of
edge $M$ at $p$.
\end{corollary}

By means of this wedge nonextendable CR measure, one may construct
non-extendable CR functions of arbitrary smoothness. Indeed, let $M$
be a local generic submanifold with central point $p$, as represented
in \S4.3 and let $K_j$ be the complex-transversal vector fields
satisfying $K_{ j_1} ( w_{ j_2} ) = \delta_{ j_1, j_2}$ and
$[ K_{ j_1}, K_{ j_2} ] = 0$.

\def\theproposition{4.40}\begin{proposition}
{\rm (\cite{ bt1981, trv1981, br1990,
trv1992, ht1996, ber1999})} For every CR distribution ${\sf T}$ on $M$
and every $\kappa \in \N$, there exist an integer $\mu \in \N$ and a
local CR function $f$ of class $\mathcal{ C}^\kappa$ defined in some
neighborhood of $p$ such that{\rm :}
\[
{\sf T}
=
\left(
K_1^2+\cdots+K_d^2
\right)^{\mu}\,
f.
\]
\end{proposition}

Then with ${\sf T} := {\sf [ S ]}$ and for $\kappa \in \N$, an
associated CR function $f$ of class $\mathcal{ C }^\kappa$ is
also shown to be not holomorphically extendable to any local wedge of
edge $M$ at $p$. A Baire category argument (\cite{ br1990}) enables
to treat the $\mathcal{ C }^\infty$ case.

\def\thetheorem{4.41}\begin{theorem}
{\rm (\cite{ br1990, ber1999})} If $M$ is not locally minimal at $p$,
then for every $\kappa = 0, 1, 2, \dots, \infty$, there exists a CR
function $h$ of class $\mathcal{ C}^\kappa$ defined in a neighborhood
of $p$ which does not extend holomorphically to any local wedge of
edge $M$ at $p$.
\end{theorem}

\def\theopenproblem{4.42}\begin{openproblem}
Find criteria for the existence of CR distributions or functions
supported by a global CR orbit.
\end{openproblem}

In~\cite{ bm1997}, this question is dealt with in the case of CR
orbits of hypersurfaces which are immersed or embedded complex
manifolds.

\smallskip

To conclude this section, we give the general form of a CR
distribution supported by a local CR orbit $S = \mathcal{ O}_{ CR}^{
loc} (M, p)$. After restriction to $S$, the collection $K_S := (K_{
e+1}, \dots, K_d)$ of vector fields spans the normal bundle to $S$ in
$M$, in a neighborhood of $p$. Let ${\sf T}$ be a local CR
distribution defined on $M$ that is supported by $S$.

\def\thetheorem{4.43}\begin{theorem}
{\rm (\cite{ trv1992, bch2005})} There exist an integer
$\kappa \in \N$, and for all $\beta \in \N^{ d-e }$ with $\vert \beta
\vert \leqslant \kappa$, local CR distributions ${\sf T }_\beta^S$
defined on $S$ such
that{\rm :}
\[
\left<
{\sf T}, \chi
\right>
=
\sum_{\beta\in\N^{d-e},\,\vert\beta\vert\leqslant\kappa}\,
\left<
{\sf T}_\beta^S,\
(K_S)^\beta(\chi)\big\vert_S
\right>.
\]
\end{theorem}

\section*{ \S5.~Approximation and uniqueness principles}

\subsection*{ 5.1.~Approximation of CR functions and of CR distributions}
Let $M$ be a generic submanifold of $\C^n$. The following
approximation theorem has appeared to be a fundamental tool in
extending CR functions holomorphically (Part~V) and in removing their
singularities (Part~VI). It is also used naturally in the
proof of Theorem~4.43 just above as well as in the Cauchy uniqueness
principle Corollary~5.4 below. The statement is valid in the general
context of locally integrable structures $\mathcal{ L}$, but, as
explained in the end of Section~3, we decided to focus our attention
on embedded Cauchy-Riemann geometry.

\def\thetheorem{5.2}\begin{theorem}
{\rm (\cite{ bt1981, hm1998, jo1999b,
bch2005})} For every $p\in M$, there exists a neighborhood $U_p$ of
$p$ in $M$ such that for every function $f$ or distribution ${\sf T}$
as defined below, there exists a sequence of holomorphic polynomials
$\left( P_k (z) \right)_{ k \in \N}$ with{\rm :}

\smallskip

\begin{itemize}

\item[$\bullet$]
if $M$ is $\mathcal{ C}^{ \kappa+2, \alpha}$, with $\kappa \geqslant 0$,
$0\leqslant \alpha \leqslant 1$, if $f$ is a CR function of class $\mathcal{
C}^{\kappa, \alpha}$ on $M$, then $\lim_{ k\to \infty} \, \left\vert
\! \left\vert P_k - f \right\vert \! \right\vert_{\mathcal{
C}^{\kappa, \alpha} (U_p) } \to 0${\rm ;} in particular, continuous CR
functions on a $\mathcal{ C}^2$ generic submanifold are
approximable sharply by holomorphic polynomials{\rm ;}

\smallskip

\item[$\bullet$]
if $M$ is at least $\mathcal{ C}^2$, if $f$ is a $L_{ loc}^{\sf
p}$ CR function $(1\leqslant {\sf p} < \infty)$, then $\lim_{ k\to
\infty} \, \left\vert \! \left\vert P_k - f \right\vert \!
\right\vert_{ L_{ loc}^{\sf p} (U_p) } \to 0${\rm ;}

\smallskip

\item[$\bullet$]
if $M$ is $\mathcal{ C}^{ \kappa+2}$, if ${\sf T}$ is a CR
distribution of order $\leqslant \kappa$ on $M$, then $\lim_{ k\to \infty}
\, \left< P_k, \chi \right> = \left< {\sf T}, \chi \right>$ for every
$\chi \in \mathcal{ C}_c^\infty (U_p)$.
\end{itemize}

\smallskip

\end{theorem}

In~\cite{ hm1998, bch2005}, convergence in Besov-Sobolev
spaces $L_{s, loc}^{ \sf p}$ and in Hardy spaces $h^{\sf p }$,
frequently used as substitutes for the $L^{\sf p}$ spaces when $ 0 <
{\sf p} < 1$, is also considered, in the context of locally integrable
structure.

\smallskip

\proof
Let us describe some ideas of the proof, assuming for simplicity that
$M$ is $\mathcal{ C}^2$ and $f$ is $\mathcal{ C}^1$. In
coordinates $(t_1, \dots, t_n)$ vanishing at $p$, choose a local
maximally real $\mathcal{ C}^2$ submanifold $\Lambda_0$
contained in $M$, passing through $p$ and satisfying $T_p \Lambda_0 
= \{ {\rm
Re}\, t = 0\}$. Let $V_p$ be a small neighborhood of $p$, whose
projection to $T_p M$ is a $(2m+d)$-dimensional open ball. We may
assume that $\Lambda_0$ is contained in $V_p$ with boundary $B_0 :=
\overline{\Lambda_0} \cap \partial V_p$ being diffeomorphic to the
$(n-1)$-dimensional sphere. Consider a parameter $u\in \R^d$
satisfying $\vert u \vert < \delta$, with $\delta >0$ small. We may
include $\Lambda_0$ in a family $(\Lambda_u )_{ \vert u \vert <
\delta}$ of maximally real $\mathcal{ C}^2$ submanifolds of
$U_p$ with $\Lambda_u\big\vert_{ u=0 } = \Lambda_0$, whose boundary is
fixed: $\partial \Lambda_u \equiv \partial \Lambda_0 = B_0$, such that
the $\Lambda_u$ foliates a small neighborhood $U_p$ of $p$ in $M$. For
$t\in U_p$, there exists a $u = u (t)$ such that $t$ belongs to
$\Lambda_{ u(t)}$. We then introduce the entire functions:
\[
F_k(t)
:=
\left(
\frac{k}{\pi}
\right)^{n/2}\,
\int_{\Lambda_{u(t)}}\,
e^{-k(t-\tau)^2}\,
f(\tau)\,
d\tau_1\wedge\cdots\wedge d\tau_n,
\]
where $(t-\tau)^2 := \sum_{ j=1 }^{ n}\, (t_j - \tau_j )^2$ and where
$k\in \N$. Shrinking $V_p$ and $U_p$ if necessary, we may assume that
$\vert {\rm Im}\, (t - \tau) \vert \leqslant \frac{ 1}{ 2}\, \vert
{\rm Re}\, (t - \tau) \vert$ for all $t, \tau \in \Lambda_u \cap U_p$
and all $\vert u \vert < \delta$. Here, the $\mathcal{
C}^2$-smoothness assumption is used. With this inequality, the above
multivariate Gaussian kernel is easily seen to be an approximation of
the Dirac distribution at $\tau = t$ on $\Lambda_{
u(t)}$. Consequently $F_k (t)$ tends to $f( t)$ as $k\to
\infty$. Moreover, the convergence is uniform and holds in $\mathcal{
C}^0 (U_p)$.

We claim that the assumption that $f$ is CR insures that $F_k
(t)$ has the same value if the integration is performed on
$\Lambda_0$:
\def\theequation{5.3}\begin{equation}
F_k(t)
=
\left(
\frac{k}{\pi}
\right)^{n/2}\,
\int_{\Lambda_0}\,
e^{-k(t-\tau)^2}\,
f(\tau)\,
d\tau_1\wedge\cdots\wedge d\tau_n.
\end{equation}
Indeed, $\Lambda_{ u(t)}$ and $\Lambda_0$ bound a $(n+1)$-dimensional
submanifold $\Pi_t$ contained in $V_p$ with $\partial \Pi_t =
\Lambda_{ u(t)} - \Lambda_0$. Since $e^{ - k (t-\tau)^2 }$ is
holomorphic with respect to $\tau$ and since $d f(\tau) \wedge d\tau_1
\wedge \cdots \wedge d\tau_n \big \vert_M = 0$, because $f$ is
$\mathcal{ C}^1$ and CR, the $(n, 0)$ form $\omega := e^{ - k (t -
\tau)^2} \, f(\tau) \, d\tau = 0$ is closed: $d\omega = 0$. By an
application of Stokes' theorem, it follows that $0 = \int_{ \Pi_t} \,
d\omega = \int_{ \Lambda_{ u(t)}} \, \omega - \int_{ \Lambda_0}\,
\omega$, which proves the claim.

Finally, to approximate $f$ by polynomials on $U_p$ in the $\mathcal{
C}^0$ topology, in the above integral~\thetag{ 5.3} that is performed
on the fixed maximally real submanifold $\Lambda_0$, it suffices to
develop the exponential in Taylor series and to integrate term by
term. In other functional spaces, the arguments have to be adapted.
\endproof

As a consequence, uniqueness in the Cauchy problem holds. It may be
shown (\cite{ trv1981, trv1992}) that the trace of a CR
distribution on a maximally real submanifold always exists, in the
distributional sense.

\def\thecorollary{5.4}\begin{corollary}
{\rm (\cite{ trv1981, trv1992})}
If a CR function or distribution vanishes on a maximally real
submanifold $\Lambda$ of $M$, there exists an open neighborhood
$U_\Lambda$ of $\Lambda$ in $M$ in which it vanishes identically.
\end{corollary}

Since every submanifold $H$ of $M$ which is generic in $\C^n$ contains
small maximally real sumanifolds passing through every of its points,
the corollary also holds with $\Lambda$ replaced by such a $H$.

\proof
It suffices to localize the above construction in a neighborhood of an
arbitrary point $p \in \Lambda$ and to take for $\Lambda_0$ a
neighborhood of $p$ in $\Lambda$. The integral~\thetag{ 5.3} then
vanishes identically.
\endproof

\def\thecorollary{5.5}\begin{corollary}
{\rm (\cite{ trv1981, trv1992})} The support of a CR function
or distribution on $M$ is a closed CR-invariant subset of $M$.
\end{corollary}

\proof
By contraposition, if a CR function or distribution vanishes in a
neighborhood $U_p$ of a point $p$ in $M$, it vanishes identically in
the {\sl CR-invariant hull} of $U_p$, viz the union of CR orbits of
all points $q \in U_p$. The CR orbits being covered by concatenations
of CR vector fields, neglecting some technicalities, the main step
is to establish:

\def\thelemma{5.6}\begin{lemma}
Let $p\in M$, let $L$ be a section of $T^c M$ and let $q^* = \exp ( {\sf
s}^* L ) (p)$ for some ${\sf s}^* \in \R$. If a CR function or
distribution vanishes in a neighborhood of $p$, it vanishes also in a
neighborhood of $q$.
\end{lemma}

Indeed, we may construct a one-parameter family $\left( H_{\sf s}
\right)_{ 0 \leqslant {\sf s} \leqslant {\sf s}^*}$ of $\mathcal{ C}^2$
hypersurfaces of $M$ with $q^* \in H_{{\sf s}^*}$ and with $H_0$
contained in a small neighborhood of $p$ at which the CR function of
distribution vanishes already. As illustrated by the following
diagram, we can insure that at every point $q_{\sf s} = \exp ( {\sf s}
L) (p)$, the vector $L ( q_{\sf s})$ is nontangent to $H_{\sf s}$.

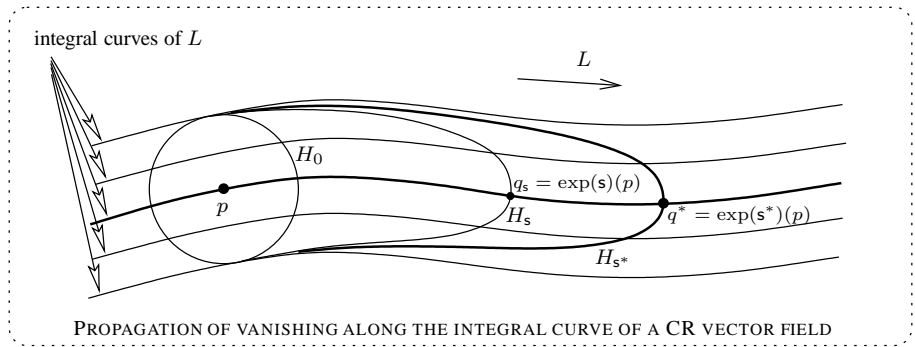
\begin{center}
\begin{picture}(0,0)%
\includegraphics{uniqueness.pstex}%
\end{picture}%
\setlength{\unitlength}{4144sp}%
\begingroup\makeatletter\ifx\SetFigFont\undefined
\def\x#1#2#3#4#5#6#7\relax{\def\x{#1#2#3#4#5#6}}%
\expandafter\x\fmtname xxxxxx\relax \def\y{splain}%
\ifx\x\y   
\gdef\SetFigFont#1#2#3{%
  \ifnum #1<17\tiny\else \ifnum #1<20\small\else
  \ifnum #1<24\normalsize\else \ifnum #1<29\large\else
  \ifnum #1<34\Large\else \ifnum #1<41\LARGE\else
     \huge\fi\fi\fi\fi\fi\fi
  \csname #3\endcsname}%
\else
\gdef\SetFigFont#1#2#3{\begingroup
  \count@#1\relax \ifnum 25<\count@\count@25\fi
  \def\x{\endgroup\@setsize\SetFigFont{#2pt}}%
  \expandafter\x
    \csname \romannumeral\the\count@ pt\expandafter\endcsname
    \csname @\romannumeral\the\count@ pt\endcsname
  \csname #3\endcsname}%
\fi
\fi\endgroup
\begin{picture}(5424,2049)(439,-1648)
\put(1693,-837){\makebox(0,0)[lb]{\smash{\SetFigFont{8}{9.6}{rm}{\color[rgb]{0,0,0}$p$}%
}}}
\put(3842, 37){\makebox(0,0)[lb]{\smash{\SetFigFont{8}{9.6}{rm}{\color[rgb]{0,0,0}$L$}%
}}}
\put(601,166){\makebox(0,0)[lb]{\smash{\SetFigFont{8}{9.6}{rm}{\color[rgb]{0,0,0}integral curves of $L$}%
}}}
\put(4385,-886){\makebox(0,0)[lb]{\smash{\SetFigFont{7}{8.4}{rm}{\color[rgb]{0,0,0}$q^*=\exp({\sf s}^*)(p)$}%
}}}
\put(3465,-693){\makebox(0,0)[lb]{\smash{\SetFigFont{7}{8.4}{rm}{\color[rgb]{0,0,0}$q_{\sf s}=\exp({\sf s})(p)$}%
}}}
\put(833,-1569){\makebox(0,0)[lb]{\smash{\SetFigFont{8}{9.6}{rm}{\color[rgb]{0,0,0}{\sc Propagation of vanishing along the integral curve of a CR vector field}}%
}}}
\put(2159,-534){\makebox(0,0)[lb]{\smash{\SetFigFont{8}{9.6}{rm}{\color[rgb]{0,0,0}$H_0$}%
}}}
\put(3419,-898){\makebox(0,0)[lb]{\smash{\SetFigFont{8}{9.6}{rm}{\color[rgb]{0,0,0}$H_{\sf s}$}%
}}}
\put(3954,-1146){\makebox(0,0)[lb]{\smash{\SetFigFont{8}{9.6}{rm}{\color[rgb]{0,0,0}$H_{{\sf s}^*}$}%
}}}
\end{picture}

\end{center}

It follows that the hypersurfaces $H_{\sf s}$ are generic in $\C^n$,
for every ${\sf s}$. Then the phrase after Corollary~5.4 applies to each
$H_{\sf s}$ from $H_0$ up to $H_{ {\sf s}^*}$, showing the propagation
of vanishing.
\endproof

\subsection*{ 5.7.~Unique continuation principles}
At least three unique continuation properties are known to be enjoyed by
holomorphic functions $h$ of several complex variables defined in a
domain $D\subset \C^n$. Indeed, we have $h \equiv 0$ in either of the
following three cases:

\smallskip

\begin{itemize}

\item[{\bf (ucp1)}]
the restriction of $h$ to some nonempty open subset of $D$ vanishes
identically;

\smallskip

\item[{\bf (ucp2)}]
the restriction of $h$ to some generic local submanifold $\Lambda$ of
$D$ vanishes identically;

\smallskip

\item[{\bf (ucp3)}]
there exists a point $p\in D$ at which the infinite jet of $h$
vanishes.

\end{itemize}

\smallskip

In Complex Analysis and Geometry, the {\bf (ucpi)} have deep influence
on the whole structure of the theory. Finer principles involving
tools from Harmonic Analysis appear in~\cite{ mp2006b}.

\def\theproblem{5.8}\begin{problem}
Find generalizations of the {\bf (ucpi)} to the category of embedded
generic submanifolds $M$.
\end{problem}

Since a domain $D$ of $\C^n$ trivially consists of a single CR orbit,
it is natural to assume that the given generic manifold $M$ is
globally minimal (although some meaningful questions arise without
this assumption, we prefer not to enter such technicalities). In this
setting, Corollaries~5.4 and~5.5 provide a complete generalization of
{\bf (ucp1)} and of {\bf (ucp2)}.

A version of {\bf (ucp3)} with the point $p$ in the boundary $\partial
D$ does not hold, even in complex dimension one. Indeed, the function
$\exp \left( e^{ i 5 \pi /4} / \sqrt{ w} \right)$ is holomorphic in
$\mathbb{ H}^+ := \{ w \in \C : \, {\rm Re}\, w >0 \}$, of class
$\mathcal{ C}^\infty$ on $\overline{ \mathbb{ H}}^+$ and flat at $w =
0$. The restriction of this function to the Heisenberg sphere ${\rm
Re}\, w = z\bar z$ of $\C^2$ provides a CR example.

To generalize rightly {\bf (ucp3)}, let $M$ be a $\mathcal{ C}^1$
generic submanifold of codimension $d\geqslant 1$ and of CR dimension
$m\geqslant 1$ in $\C^n$, with $n = m+d$. Let $\Sigma$ be a $\mathcal{
C}^1$ submanifold of $M$ satisfying:
\[
T_q^c M\oplus T_q\Sigma
=
T_qM, 
\ \ \ \ \ \ \ \ \ 
q\in\Sigma.
\]
Here, $\Sigma$ plays the r\^ole of the point $p$ in {\bf (ucp3)}.
Denote by $\mathcal{ O}_{ CR }^\Sigma$ the union of CR orbits of
points of $\Sigma$, {\it i.e.} the {\sl CR-invariant hull} of
$\Sigma$. It is an open subset of $M$. We say that a CR function $f :
M \to \C$ of class $\mathcal{ C}^1$ {\sl vanishes to infinite order
along} $\Sigma$ if for every $p\in \Sigma$, there exists an open
neighborhood $U_p$ of $p$ in $M$ such that for every $\nu \in \N$,
there exists a constant $C >0$ with
\[
\vert h(t)
\vert \leqslant C \left[ {\rm dist} (t, \Sigma) \right]^\nu, 
\ \ \ \ \ \ \ \ \ 
t\in U_p.
\]

\def\thetheorem{5.9}\begin{theorem}
{\rm (\cite{ ro1986b, bt1988}, [$*$])} Assume that $\Sigma$ is
the intersection with $M$ of some $d$-dimensional holomorphic
submanifold of $\C^n$. If a CR function of class $\mathcal{ C }^1$ 
vanishes to infinite order along $\Sigma$, then it vanishes
identically on the globally minimal generic submanifold $M$.
\end{theorem}

Assuming that $\Sigma$ is only a conic $d$-codimensional holomorphic
submanifold entering a wedge to which all CR functions of $M$ extend
holomorphically (Theorem~3.8(V)), the proof of this theorem may
be easily generalized.

\def\theopenquestion{5.10}\begin{openquestion}
{\rm (\cite{ ro1986b, bt1988})} Is the above unique
continuation true for $\Sigma$ merely $\mathcal{ C}^1$\,?
\end{openquestion}

To attack this question, one should start with $M$ being unit sphere
$S^3 \subset \C^2$ and $\Sigma \subset S^3$ being any $T^c
S^3$-transversal real segment which is nowhere locally the boundary of
a complex curve lying inside the ball.

\newpage

\begin{center}
{\Large\bf IV:~Hilbert transform and Bishop's equation in 
H\"older spaces}
\end{center}

\bigskip\bigskip\bigskip

\begin{center}
\begin{minipage}[t]{11cm}
\baselineskip =0.35cm
{\scriptsize

\centerline{\bf Table of contents}

\smallskip

{\bf 1.~H\"older spaces: basic properties \dotfill 112.}

{\bf 2.~Cauchy integral, Sokhotski\u\i-Plemelj
formulas and Hilbert transform \dotfill 114.}

{\bf 3.~Solving a local parametrized
Bishop equation with optimal loss of smoothness \dotfill 129.} 

{\bf 4.~Appendix: proofs of some lemmas \dotfill 146.} 

\smallskip

{\footnotesize\tt \hfill [1 diagram]}

}\end{minipage}
\end{center}

\bigskip
\bigskip

{\small

{\small


In complex and harmonic analysis, the spaces $\mathcal{ C}^{ \kappa,
\alpha}$ of fractionally differentiable maps, called {\sl H\"older
spaces}, are very flexible to generate inequalities and they yield
rather satisfactory norm estimates for almost all the classical
singular integral operators, especially when $0 < \alpha < 1$. For
instance, the Cauchy integral of a $\mathcal{ C}^{ \kappa, \alpha}$
function $f : \Gamma \to \C$ defined on a $\mathcal{ C}^{ \kappa + 1,
\alpha }$ Jordan curve $\Gamma$ of the complex plane produces a
sectionally holomorphic function, whose boundary values from one or
the either side are $\mathcal{ C}^{ \kappa, \alpha}$ on the curve. The
Sokhotski\u\i-Plemelj formulas show that the arithmetic mean of the
two (in general different) boundary values at a point of the curve is
given by the principal value of the Cauchy integral at that point.

Harmonic and Fourier analysis on the unit disc $\Delta$ is of
particular interest for geometric applications in Cauchy-Riemann
geometry. According to a theorem due to Privalov, the Hilbert
transform ${\sf T}$ is a bounded linear endomorphism of $\mathcal{
C}^{ \kappa, \alpha} (\partial \Delta, \R)$ with norm $\left\vert \!
\left\vert \! \left\vert {\sf T} \right\vert \! \right\vert \!
\right\vert_{ \kappa, \alpha}$ equivalent to $\frac{ C}{ \alpha (1 -
\alpha)}$, for some absolute constant $C > 0$. This operator produces
the harmonic conjugate ${\sf T} u$ of any real-valued function $u:
\partial \Delta \to \R$ on the unit circle, so that $u + i\, {\sf T}
u$ always extends holomorphically to $\Delta$. Bishop (1965),
Hill-Taiani (1978), Boggess-Pitts (1985) and Tumanov (1990) formulated
and solved a functional equation involving ${\sf T}$ in order to find
small analytic discs with boundaries contained in a generic
submanifold $M$ of codimension $d$ in $\C^n$.

In a general setting, this Bishop-type equation is of the form:
$$
U (e^{ i\, \theta}) = U_0 - {\sf T}
\left[ \Phi ( U(\cdot), \cdot, s) \right] (e^{ i\, \theta}), 
$$
where $U_0 \in \R^d$ is a constant vector, where $\Phi = \Phi ( u,
e^{i\, \theta}, s)$ is an $\R^d$-valued $\mathcal{ C}^{ \kappa,
\alpha}$ map, with $\kappa \geqslant 1$ and $0 < \alpha < 1$, where
$u\in \R^d$, where $e^{ i\, \theta} \in \partial \Delta$ and where $s
\in \R^b$ is an additional parameter which is useful in geometric
applications. Under some explicit assumptions of smallness of $U_0$
and of the first order jet of $\Phi$, the general solution $U = U (
e^{ i\, \theta}, s, U_0)$ is of class $\mathcal{ C}^{ \kappa, \alpha}$
with respect to $e^{i\, \theta}$ and in addition, for every $\beta$
with $0 < \beta < \alpha$, it is of class $\mathcal{ C}^{ \kappa,
\beta}$ with respect to all the variables $(e^{ i\, \theta}, s, U_0)$.
These smoothness properties are optimal.

}

\section*{1.~H\"older spaces: basic properties}

\subsection*{1.1.~Background on H\"older spaces}
Let $n\in \N$ with $n\geqslant 1$ and let ${\sf x} = ({\sf x}_1,
\dots, {\sf x}_n) \in \R^n$. On the vector space $\R^n$, we choose
once for all the {\sl maximum norm} $\vert {\sf x} \vert := \max_{
1\leqslant i\leqslant n} \vert {\sf x}_i \vert$ and, for any
``radius'' $\rho$ satisfying $0 < \rho \leqslant \infty$, we define
the {\sl open cube} $\square_\rho^n := \{ {\sf x} \in \R^n : \vert x
\vert < \rho\}$ as a fundamental, concrete open set. For $\rho =
\infty$, we identify $\square_\infty^n$ with $\R^n$.

Let $\kappa \in \N$ and let $\alpha \in \R$ with $0 \leqslant \alpha
\leqslant 1$. If $\K = \R$ or $\C$, a scalar function $f :
\square_\rho^n \to \K$ belongs to the {\sl H\"older class $\mathcal{
C}^{\kappa, \alpha} (\square_\rho^n, \K)$} if, for every multiindex
$\delta = (\delta_1, \dots, \delta_n) \in \N^n$ of {\sl length} $\vert
\delta \vert \leqslant \kappa$, the partial derivative $f_{{\sf
x}^\delta} ( {\sf x}) := \frac{\partial^{\vert \delta \vert} f}{
\partial {\sf x}^{ \delta_1} \cdots \partial {\sf x}^{\delta_n}}$ is
continuous in $\square_\rho^n$ and if, moreover, the quantity:
\[
\left\vert\!\left\vert
f
\right\vert\!\right\vert_{
\kappa,\alpha}:=
\sum_{0\leqslant\vert\delta\vert\leqslant\kappa}
\sup_{{\sf x}\in\square_\rho^n}
\left\vert 
f_{{\sf x}^\delta}({\sf x})
\right\vert
+
\sum_{\vert\delta\vert=\kappa}\,
\sup_{{\sf x}''\neq{\sf x}'\in\square_\rho^n}
\frac{
\left\vert
f_{{\sf x}^\delta}({\sf x}'')
-
f_{{\sf x}^\delta}({\sf x}')
\right\vert
}{
\vert 
{\sf x}''
-
{\sf x}'
\vert^\alpha}
\]
is finite (if $\alpha =0$, it is understood that the second sum is
absent). In case $f = (f_1, \dots, f_m)$ is a $\K^m$-valued mapping,
with $m\geqslant 1$, we simply define $\left\vert\!\left\vert f
\right\vert\!\right\vert_{ \kappa,\alpha}:= \max_{ 1\leqslant j \leqslant m}
\left\vert\!\left\vert f_j \right\vert\!\right \vert_{
\kappa,\alpha}$. This is coherent with the choice of the maximum norm
$\left\vert {\sf y} \right\vert := \max_{ 1\leqslant i \leqslant n} \left\vert
{\sf y}_i \right\vert$ on $\K^m$. For short, such a map will be said
to be {\sl $\mathcal{ C}^{\kappa, \alpha}$-smooth} or {\sl of class
$\mathcal{ C}^{\kappa, \alpha}$} and we write $f \in \mathcal{
C}^{\kappa, \alpha}$. One may verify $\left \vert\! \left\vert f_1
f_2 \right \vert\!\right\vert_{ \kappa, \alpha} \leqslant \left\vert\!
\left\vert f_1 \right \vert\! \right \vert_{ \kappa, \alpha} \cdot
\left\vert\! \left\vert f_2 \right \vert\!\right\vert_{ \kappa,
\alpha}$ and of course $\left \vert\! \left\vert \lambda_1 f_1 +
\lambda_2 f_2 \right \vert\!\right\vert_{ \kappa, \alpha} \leqslant \vert
\lambda_1 \vert \left\vert\! \left\vert f_1 \right \vert\! \right
\vert_{ \kappa, \alpha} + \vert \lambda_2 \vert \left\vert\!
\left\vert f_2 \right \vert\!\right\vert_{ \kappa, \alpha}$. If
$\kappa =0$ and $\alpha = 1$, the map $f$ is called {\sl
Lipschitzian}. The condition $\big\vert f(
e^{i\, \theta''}) - f(e^{ i\, \theta'}) 
\big\vert \leqslant C \cdot \left\vert \theta '' - \theta' 
\right\vert$ on the
unit circle was first introduced by Lipschitz in 1864 as sufficient
for the pointwise convergence of Fourier series.

Thanks to a uniform convergence argument, the space $\mathcal{
C}^{\kappa, \alpha} (\square_\rho^n, \K)$ is shown to be complete,
hence it constitutes a {\sl Banach algebra}. The space of functions
defined on the closure $\overline{ \square_\rho^n}$ also constitutes a
Banach algebra. If $\alpha$ is positive, thanks to a prolongation
argument, one may verify that $\mathcal{ C}^{\kappa, \alpha}
(\square_\rho^n, \K)$ identifies with the restriction
$\left. \mathcal{ C}^{ \kappa, \alpha} \left( \overline{
\square_\rho^n}, \K \right) \right \vert_{ \square_\rho^n}$.

H\"older spaces may also be defined on arbitrary convex open subsets.
More generally, on an arbitrary subset $\Omega \subset \R^n$, it is
reasonable to define the H\"older norms $\left\vert \! \left \vert
\cdot \right\vert \! \right\vert_{\kappa, \alpha}$, $0<\alpha\leqslant 1$,
only if ${\rm dist}_\Omega ({\sf x}'', {\sf x}') \leqslant C \cdot
\left\vert {\sf x}'' - {\sf x} '\right \vert$ for every two points
${\sf x}'', {\sf x}' \in \Omega$. This is the case for instance if
$\Omega$ is a domain in $\R^n$ having
piecewise $\mathcal{ C}^{1, 0}$
boundary.

Introducing the total order $(\kappa_1, \alpha_1) \leqslant (\kappa,
\alpha)$ defined by: $\kappa_1 < \kappa$, or: $\kappa_1 = \kappa$ and
$\alpha_1 \leqslant \alpha$, we verify that $\mathcal{ C}^{ \kappa,
\alpha}$ is contained in $\mathcal{ C}^{ \kappa_1, \alpha_1}$ and
that:
\begin{itemize}
\item[$\bullet$]
$\left \vert \! \left \vert f \right\vert\! \right \vert_{ \kappa, 0}
\leqslant \left \vert \! \left \vert f \right \vert \! \right \vert_{
\kappa, \alpha}$ for all $\alpha$ with $0 < \alpha \leqslant 1$ and for all
$\kappa \in \N$;
\item[$\bullet$]
$\left \vert \! \left \vert f \right\vert\! \right \vert_{ \kappa,
\alpha_1} \leqslant 3 \left \vert \! \left \vert f \right \vert \! \right
\vert_{ \kappa, \alpha_2 }$ for all $\alpha_1, \alpha_2$ with $0 <
\alpha_1 < \alpha_2 \leqslant 1$ and for all $\kappa \in \N$;
\item[$\bullet$]
$\left \vert \! \left \vert f \right\vert\! \right \vert_{ \kappa, 1}
\leqslant \left \vert \! \left \vert f \right \vert \! \right \vert_{
\kappa+ 1, 0}$, for all $\kappa \in \N$.
\end{itemize}
The first inequality above is trivial while the
third follows from~\thetag{ 1.3} below.
We explain the factor $3$ in the second inequality. Since
$\left \vert {\sf x}'' - {\sf x}' \right\vert^{-\alpha_1} \leqslant
\left \vert {\sf x}'' - {\sf x}' \right\vert^{-\alpha_2}$ only if
$\left \vert {\sf x}'' - {\sf x}' \right\vert \leqslant 1$, we 
may estimate:
\[
\sup_{0<\left\vert {\sf x}'' -{\sf x}' \right\vert \leqslant 1}
\frac{\left\vert f({\sf x}'')-f({\sf x}')\right\vert}{
\left\vert {\sf x}''-{\sf x}' \right\vert^{\alpha_1}}
\leqslant
\sup_{0<\left\vert {\sf x}'' -{\sf x}' \right\vert \leqslant 1}
\frac{\left\vert f({\sf x}'')-f({\sf x}')\right\vert}{
\left\vert {\sf x}''-{\sf x}' \right\vert^{\alpha_2}}
\leqslant
\left \vert \! \left \vert f
\right\vert\! \right \vert_{ 0, \alpha_2}.
\]
On the other hand, if $\left \vert {\sf x}'' - {\sf x}'\right\vert
>1$, we simply apply the (not fine) inequalities:
\[
\aligned
\frac{\left\vert f({\sf x}'')-f({\sf x}')\right\vert}{
\left\vert {\sf x}''-{\sf x}' \right\vert^{\alpha_1}}
\leqslant 
\left\vert f({\sf x}'')-f({\sf x}')\right\vert
\leqslant
2
\left \vert \! \left \vert f
\right\vert\! \right \vert_{ 0, 0}
\leqslant
2
\left \vert \! \left \vert f
\right\vert\! \right \vert_{ 0, \alpha_2}.
\endaligned
\]
Consequently: 
\[
\left \vert \! \left \vert f \right\vert\! \right \vert_{ 
0, \alpha_1}
=
\left \vert \! \left \vert f \right\vert\! \right \vert_{
0,0}
+
\sup_{{\sf x}''\neq {\sf x}'}
\frac{\left\vert f({\sf x}'')-f({\sf x}')\right\vert}{
\left\vert {\sf x}''-{\sf x}' \right\vert^{\alpha_1}}
\leqslant
3 \left \vert \! \left \vert f \right\vert\! \right
\vert_{ 0, \alpha_2},
\]
with a factor $3$. For general $\kappa \geqslant 1$, the desired
inequality follows:
\[
\aligned
\left \vert \! \left \vert f
\right\vert\! \right \vert_{ \kappa, \alpha_1}
=
\left \vert \! \left \vert f
\right\vert\! \right \vert_{ \kappa-1, 0}
+
\sum_{\vert\delta\vert=\kappa}
\left \vert \! \left \vert 
f_{{\sf x}^\delta}
\right\vert\! \right \vert_{0,\alpha_1}
&
\leqslant 
\left \vert \! \left \vert 
f
\right\vert\! \right \vert_{ \kappa-1, 0}
+
3
\sum_{\vert\delta\vert=\kappa}
\left \vert \! \left \vert 
f_{{\sf x}^\delta}
\right\vert\! \right \vert_{0,\alpha_2}
\\
&
\leqslant
3
\left \vert\!\left \vert f
\right\vert\!\right\vert_{ \kappa, \alpha_2}.
\endaligned
\]

In the sequel, sometimes, we might abbreviate $\mathcal{ C}^{\kappa,
0}$ by $\mathcal{ C}^\kappa$, a standard notation. However, we shall
never abbreviate $\mathcal{ C}^{0, \alpha}$ by $\mathcal{ C}^\alpha$,
in order to avoid the unpleasant ambiguity $\mathcal{ C}^{1, 0} \equiv
\mathcal{ C}^1 \equiv \mathcal{ C}^{ 0, 1}$. Without providing
proofs, let us state some fundamental structural properties of
H\"older spaces. Some of them are in~\cite{ kr1983}.

\smallskip

$\bullet$ 
The inclusions $\mathcal{ C}^{\lambda, \beta} \subset
\mathcal{ C}^{\kappa, \alpha}$ for $(\lambda, \beta) > (\kappa,
\alpha)$ are all strict. For instance, on $\R$, the function
$\chi_{\kappa, \alpha} = \chi_{\kappa, \alpha} ({\sf x})$ equal to
zero for ${\sf x} \leqslant 0$ and, for ${\sf x}\geqslant 0$:
\[
\chi_{\kappa, \alpha}({\sf x})
=
\left\{
\aligned
&
{\sf x}^{\kappa+\alpha}, 
\ \ \ \ \ \ \ \ \ \ \
{\rm if}
\ \
0<\alpha\leqslant 1,
\\
&
{\sf x}^\kappa/\text{\rm log} \, {\sf x},
\ \ \ \ \ \ \,
{\rm if}
\ \
\alpha=0,
\endaligned
\right.
\]
is $\mathcal{ C}^{\kappa, \alpha}$ in any neighborhood of
the origin, not better.

\smallskip

$\bullet$
If $0 < \alpha_1 < \alpha$, any uniformly bounded set of
functions in $\mathcal{ C}^{\kappa, \alpha}$ contains a sequence of
functions that converges in $\mathcal{ C}^{\kappa, \alpha_1}$-norm to
a function in $\mathcal{ C}^{\kappa, \alpha_1}$. This is 
a H\"older-space version of the Arzel\`a-Ascoli lemma.

\smallskip

$\bullet$ For $0 < \alpha \leqslant 1$, define the {\sl H\"older semi-norm}
(notice the wide hat):
\[
\left\vert\!\left\vert
f
\right\vert\!\right\vert_{\widehat{0,\alpha}}
:=
\sup_{{\sf x}''\neq{\sf x}'\in\square_\rho^n}
\frac{\left\vert
f({\sf x}'')-f({\sf x}')
\right\vert}
{\left\vert
{\sf x}''-{\sf x}'
\right\vert^\alpha}. 
\]
The constants satisfy $\left \vert\! \left \vert c \right \vert\!
\right\vert_{\widehat{ 0, \alpha}} =0 $ and, of course, we have $\left
\vert\! \left \vert f \right \vert\! \right\vert_{ 0, \alpha} \equiv
\left \vert\! \left \vert f \right \vert\! \right\vert_{ 0, 0} + \left
\vert\! \left \vert f \right \vert\! \right\vert_{\widehat{ 0,
\alpha}}$. As a function of $\alpha$, {\it
the semi-norm is logarithmically
convex}:
\[
\left \vert\! \left \vert
f
\right \vert\! \right\vert_{
\widehat{ 0, t\alpha_1 + (1-t) \alpha_2 }}
=
\left(
\left \vert\! \left \vert
f
\right \vert\! \right\vert_{
\widehat{ 0, \alpha_1}}
\right)^t
\cdot
\left(
\left \vert\! \left \vert
f
\right \vert\! \right\vert_{
\widehat{ 0, \alpha_2}}
\right)^{1-t}.
\]
Here, $0 < \alpha_1 < \alpha_2 \leqslant 1$ and $0 \leqslant t \leqslant 1$.

\smallskip

$\bullet$
Importantly, if $f$ is $\K^m$-valued, if $1\leqslant l\leqslant m$, from the
{\sl Taylor integral formula}:
\def\theequation{1.2}\begin{equation}
f_l\left({\sf x}''\right)
-
f_l\left({\sf x}'\right)
=
\int_0^1\, 
\sum_{i=1}^n\,
\frac{\partial f_l}{\partial {\sf x}_i}
\left({\sf x}' + {\sf s} ({\sf x}''-{\sf x}')
\right)
\left[
{\sf x}_i''-{\sf x}_i'
\right]
d{\sf s},
\end{equation}
follows the {\sl mean value inequality}:
\def\theequation{1.3}\begin{equation}
\aligned
\left\vert
f\left({\sf x}''\right)
-
f\left({\sf x}'\right)
\right\vert
&
=
\max_{1\leqslant l\leqslant m}\,
\left\vert
f_l\left({\sf x}''\right)
-
f_l\left({\sf x}'\right)
\right\vert
\\
&
\leqslant
\left\vert\!\left\vert
f
\right\vert\!\right\vert_{\widehat{ 1,0}}
\cdot
\left\vert
{\sf x}''
-
{\sf x}'
\right\vert,
\endaligned
\end{equation}
where ${\sf x}'', {\sf x}' \in \square_\rho^n$ are arbitrary, and
where
\[
\left\vert\!\left\vert
f
\right\vert\!\right\vert_{\widehat{1,0}}
:=
\max_{1\leqslant l\leqslant m}\,
\sum_{k=1}^n\,
\sup_{\left\vert{\sf x}\right\vert<\rho}
\left\vert
f_{l,{\sf x}_k}({\sf x})
\right\vert.
\]
This useful inequality also holds (by definition) if $f$ is merely
Lipschitzian, with $\left \vert \! \left \vert f \right \vert \!
\right \vert_{ \widehat{ 1,0} }$ replaced by $\left \vert \! \left
\vert f \right \vert \! \right \vert_{ \widehat{ 0,1} }$.

\smallskip

$\bullet$
If a function $f$ is $\mathcal{ C }^{ \kappa, 0}$, then for every
multiindex $\delta \in \N^n$ of length $\left \vert \delta \right
\vert \leqslant \kappa$, the partial derivative $f_{{\sf x}^\delta}$ is
$\mathcal{ C}^{\kappa- \left \vert \delta \right\vert, 0}$ and $\left
\vert \! \left \vert f_{{\sf x}^\delta} \right \vert \! \right
\vert_{\kappa - \left \vert \delta \right \vert, 0} \leqslant \left \vert
\! \left \vert f \right \vert \! \right \vert_{\kappa, 0}$.

\section*{ \S2.~Cauchy integral, Sokhotski\u\i-Plemelj 
formulas and Hilbert transform}

\subsection*{2.1.~Boundary behaviour of the
Cauchy integral} Let $\Omega$ be a domain in $\C$, let $z \in \Omega$
and let $\Gamma$ be a $\mathcal{ C}^1$-smooth simple closed curve
surrounding $z$ and oriented counterclockwise. Assume that its
interior domain (to which $z$ belongs) is entirely contained in
$\Omega$. In case $\Gamma$ is a circle, Cauchy (\cite{ ca1831})
established in 1831 the celebrated representation formula:
\[
f(z)
= 
\frac{1}{2\pi i}\,
\int_\Gamma\,
\frac{f(\zeta)\,d\zeta}{\zeta-z},
\]
valid for all functions $f \in \mathcal{ O} (\Omega)$ holomorphic in
$\Omega$. Remarkably, $\Gamma$ may be modified and deformed without
altering the value $f(z)$ of the integral.

The best proof of this formula is to derive it from the more general
{\sl Cauchy-Green-Pompeiu formula}, itself being an elementary
consequence of the Green-Stokes formula, which is valid for functions
$f$ of class only $\mathcal{ C}^1$ defined on the closure
of a domain $\Omega \subset
\C$ having $\mathcal{ C}^1$-smooth oriented boundary $\partial \Omega$
(\cite{ ho1973}): 
\[
f (z) = \frac{1}{2\pi i}\, \int_{\partial
\Omega} \, \frac{f(\zeta)\,d\zeta
}{\zeta-z} +\frac{ 1}{ 2 \pi i} \int
\!\! \int_\Omega \, 
\frac{\partial f /\partial \bar \zeta}{\zeta - z}
\, d \zeta \wedge d\bar \zeta.
\]
Indeed, for holomorphic $f$, one
clearly sees that the ``remainder'' double integral disappears.

The holomorphicity of the kernel $\frac{ 1}{ \zeta - z}$ enables then
to build concisely the fundamental properties of holomorphic functions
from Cauchy's formula: local convergence of Taylor series,
residue theorem, Cauchy uniform convergence theorem, maximum principle,
{\it etc.} (\cite{ ho1973}). Studying the Cauchy integral for
itself appeared therefore to be of interest and became a thoroughly
investigated subject in the years 1910--1960, under the influence of
Privalov.

If $z \in \Omega$ belongs to the exterior of $\Gamma$, {\it i.e.} to
the unbounded component of $\C \backslash \Gamma$, by a fundamental
theorem also due to Cauchy, the integral vanishes: $0 = \frac{1}{2\pi
i}\, \int_\Gamma\, \frac{f(\zeta)\,d\zeta}{\zeta-z}$. Thus, fixing the
countour $\Gamma$, as $z$ moves toward $\Gamma$, the Cauchy integral
is constant, either equal to $f(z)$ or to $0$. What happens when $z$
hits the curve $\Gamma$\,?

Denote by $\zeta_0$ a point of $\Gamma$ and by $\Delta(\zeta_0,
\varepsilon)$ the open disc of radius $\varepsilon > 0$ centered at
$\zeta_0$. If $\Gamma_\varepsilon$ denotes the complement $\Gamma
\backslash \Delta(\zeta_0, \varepsilon)$, introducing 
an arc of small circle contained in $\partial \Delta
(\zeta_0, \varepsilon)$ to join the two extreme points of
$\Gamma_\varepsilon$, it may be verified that
\def\theequation{2.2}\begin{equation}
\frac{1}{2}\,
f(\zeta_0)
=
\lim_{\varepsilon \to 0}\,
\frac{1}{2\pi i}\,
\int_{\Gamma_\varepsilon}\,
\frac{f(\zeta)\,d\zeta}{\zeta-\zeta_0}.
\end{equation}
Geometrically speaking, essentially one half of the circle $\partial
\Delta(\zeta_0, \varepsilon)$ of radius $\varepsilon$ centered at
$\zeta_0$ is contained in the domain $\Omega$. Consequently, the
``correct value'' of the Cauchy integral at a point $\zeta_0$ of the
curve $\Gamma$ is equal to the arithmetic mean:
\[
\frac{ 1}{ 2} 
\left(
\lim_{ z \to \zeta_0, \ z \ {\rm inside}}
+
\lim_{ z \to \zeta_0, \ z \ {\rm outside}}
\right)
=
\frac{ 1}{ 2} 
\left(
f(\zeta_0)+0
\right).
\]
Let us recall briefly why the excision of an
$\varepsilon$-neighborhood of $\zeta_0$ in the domain of integration
is necessary to provide this ``correct'' average value. Parametrizing
$\Gamma$ by a real number, the problem of giving a sense to the
singular integral $\int_\Gamma \, \frac{ f(\zeta) \, d\zeta}{\zeta -
\zeta_0}$ amounts to the following classical definition of the notion
of {\sl principal value} (\cite{ mu1953, ga1966, ek2000}).

\subsection*{ 2.3.~Principal value integrals}
Let $a, b \in \R$ with $a < b$ and let $f$ be a $\mathcal{
C}^1$-smooth real-valued function defined on the open segment
$(a,b)$. Pick ${\sf x} 
\in \R$ with $a < {\sf x} < b$ and consider the integral
$\int_a^b \frac{ d{\sf y}}{ 
{\sf y} - {\sf x}}$ whose integrand is singular. The two
integrals avoiding the singularity from the left and from the right,
namely:
\[
\aligned
\int_a^{{\sf x}-\varepsilon_1}\, 
\frac{d{\sf y}}{{\sf y}-{\sf x}}
&
=
{\rm log}(\varepsilon_1) 
-
{\rm log}({\sf x}-a) 
\ \ \ \ \ \ \ \ \ \ \ \ \ \ \
{\rm and}
\\
\int_{{\sf x}+\varepsilon_2}^b\,
\frac{d{\sf y}}{{\sf y}-{\sf x}}
&
=
{\rm log}(b-{\sf x})
-{\rm log}(\varepsilon_2)
\endaligned
\]
tend to $-\infty$, as $\varepsilon_1 \to 0^+$, and to $+\infty$ as
$\varepsilon_2\to 0^+$. Clearly, if $\varepsilon_2 = \varepsilon_1$
(or more generally, if $\varepsilon_1$ and $\varepsilon_2$ both depend
continuously on an auxiliary parameter $\varepsilon >0$ with $1 =
\lim_{ \varepsilon \to 0^+} \frac{ \varepsilon_2( \varepsilon)}{
\varepsilon_1 (\varepsilon)}$), the positive and the negative parts
compensate, so that the {\sl principal value}:
\[
{\rm p.v.}
\int_a^b\,
\frac{d{\sf y}}{{\sf y}-{\sf x}}
:= 
\lim_{\varepsilon\to 0^+}
\left(
\int_a^{{\sf x}-\varepsilon}+
\int_{{\sf x}+\varepsilon}^b
\right)
=
{\rm log}\,
\frac{b-{\sf x}}{{\sf x}-a}
\]
exists. Briefly, there is a key cancellation of infinite parts, thanks
to the fact that the singular kernel $\frac{ 1}{ {\sf y}}$ 
is odd. This is
why in~\thetag{ 2.2} above, the integration was performed over the
excised curve $\Gamma_\varepsilon$.

Generally, if $g: [a, b] \to \R$ is a real-valued function the {\sl
principal value} integral, defined by:
\[
\aligned
{\rm p.v.}
\int_a^b\,
\frac{g({\sf y})\,d{\sf y}}{{\sf y}-{\sf x}}
:= 
&
\lim_{\varepsilon\to 0^+}
\left(
\int_a^{{\sf x}-\varepsilon}+
\int_{{\sf x}+\varepsilon}^b
\right)
\\
=
&
\int_a^b\,
\frac{g({\sf y}) - g({\sf x})}{{\sf y}-{\sf x}}\,d{\sf y}
+
g({\sf x})\,
{\rm p.v.}
\int_a^b\,\frac{d{\sf y}}{{\sf y}-{\sf x}}\,d{\sf y}
\\
=
&
\int_a^b\,
\frac{g({\sf y}) - g({\sf x})}{{\sf y}-{\sf x}}\,d{\sf y}
+
g({\sf x})\,{\rm log}\,\frac{b-{\sf x}}{{\sf x}-a}
\endaligned
\]
exists whenever the quotient $\frac{g({\sf y}) - g({\sf x})}{{\sf
y}-{\sf x}}$ is integrable. This is the case for instance if $g$ is
of class $\mathcal{ C}^{1,0}$ or of class $\mathcal{ C}^{0,\alpha}$,
with $\alpha>0$, since $\int_0^1 \, {\sf y}^{ \alpha - 1} \, d{\sf y}
< \infty$. More is true.

\def\thetheorem{2.4}\begin{theorem}
{\rm (\cite{ mu1953, ve1962, dy1991, sme1988, 
ek2000}, [$*$])} Let $g: [a, b] \to \R$ be $\mathcal{ C}^{\kappa,
\alpha}$-smooth, with $\kappa \geqslant 0$ and $0 < \alpha < 1$. Then for
every ${\sf x}\in (a,b)$, the principal value integral
\[
G({\sf x})
:=
{\rm p.v.}
\int_a^b\, \frac{ g({\sf y})\,d{\sf y} }{ {\sf y}-{\sf x}}
\]
exists. In every closed segment $[a', b']$ contained in $(a,b)$, the
function $G ({\sf x})$ becomes $\mathcal{ C}^{ \kappa, \alpha}$-smooth
and enjoys the norm inequality $\left\vert \! \left\vert G
\right\vert \! \right\vert_{ \mathcal{ C}^{\kappa, \alpha}[a', b']}
\leqslant \frac{C}{\alpha (1- \alpha)}\, \left\vert \! \left\vert g
\right\vert \! \right\vert_{\mathcal{ C}^{\kappa, \alpha}[a,b]}$, for
some constant $C = C(\kappa, a, b, a', b')$. If $g$ together with its
derivatives up to order $\kappa$ vanish at the two extreme points $a$
and $b$, the function $G ({\sf x})$ is $\mathcal{ C}^{ \kappa,
\alpha}$-smooth over $[a,b]$ and enjoys the norm inequality
$\left\vert \! \left\vert G \right\vert \! \right\vert_{ \mathcal{
C}^{\kappa, \alpha}[a, b]} \leqslant \frac{C }{ \alpha (1- \alpha)}\,
\left\vert \! \left\vert g \right\vert \! \right\vert_{\mathcal{
C}^{\kappa, \alpha}[a,b]}$, for some constant $C= C (\kappa, a,b)$.
\end{theorem}

Notice the presence of the (nonremovable) factor $\frac{ 1}{
\alpha(1-\alpha)}$.

\subsection*{2.5.~General Cauchy integral}
Beginning with works of Sokhotski\u\i~\cite{ so1873}, of
Harnack~\cite{ ha1885} and of Morera~\cite{ mo1889}, the {\sl Cauchy
integral transform}:
\[
F(z)
:=
\frac{ 1}{ 2\pi i} \,
\int_\Gamma\, \frac{ f(\zeta) \, d\zeta}{ \zeta - z} 
\]
has been studied for itself, in the more general case where $\Gamma$
is an arbitrary closed or non-closed curve in $\C$ and $f$ is an
arbitrary smooth complex-valued function defined on $\Gamma$, not
necessarily holomorphic in a neighborhood of $\Gamma$ (precise rigorous
assumptions will follow; historical account may be found in \cite{
ga1966}). In Sokhotski\u\i's and in Harnack's
works, the study of the boundary behaviour of the Cauchy integral was
motivated by physical problems; its boundary properties find
applications to mechanics, to hydrodynamics and to elasticity
theory. Let us restitute briefly the connection to the notion of
logarithmic potential (\cite{ mu1953}).

Assuming $\Gamma$ and $f: \Gamma \to \R$ to be real-valued and of
class at least $\mathcal{ C}^{ 1, 0}$, parametrize $\Gamma$ by
arc-length $\zeta = \zeta (s)$, denote ${\bf r} (s) := \zeta (s) - z$
the radial vector from $z$ to $\zeta (s)$, denote $r = r(s) =
\left\vert {\bf r} (s) \right\vert$ its euclidean norm, denote ${\bf
t} (s) := \frac{ d{\bf r}}{ ds}$ the unit tangent vector field to
$\Gamma$ and denote ${\bf n} (s) := \frac{ d{\bf r}}{ ds} / \left\vert
\frac{ d{\bf r}}{ ds} \right\vert$ the unit normal vector field to
$\Gamma$. Puting $z= x+ iy$ and decomposing the Cauchy transform $F(
z) = U( x, y) + i V( x, y)$ in real and imaginary parts, the two
functions $U$ and $V$ are harmonic in $\C \backslash \Gamma$, since
$F$ is clearly holomorphic there. After elementary computations, one
shows that $U$ may be expressed under the form:
\[
U(x,y)
=
\frac{1}{2\pi}\,
\int_\Gamma\,
f\,\frac{\cos({\bf r},{\bf n})}{r}\,ds,
\]
which, physically, represents the {\sl potential of a double layer
with moment-density} $\frac{f}{2\pi}$. Also, $V$ may be expressed
under the form:
\[
V(x,y)
=
\frac{1}{2\pi}\,
\int_\Gamma\,\frac{df}{ds}\,{\rm log}\,r\,ds,
\]
which, in the case where $\Gamma$ consists of a finite number of
closed Jordan curves, represents the potential of a single layer with
moment-density $-\frac{ 1}{ 2\pi}\, \frac{df}{ ds}$.

\subsection*{2.6.~The Sokhotski\u\i-Plemelj 
formulas} Coming back to the mathematical study of the Cauchy
integral, we shall assume that the curve $\Gamma$ over which the
integration is performed is a {\it connected}\, curve of finite length
parametrized by arc length
\[
[a,b] \ni s \longmapsto \zeta(s)\in\Gamma,
\]
where $a<b$, where $\zeta(s)$ is of class $\mathcal{ C}^{\kappa+1,
\alpha}$ over the closed segment $[a,b]$, and where $\kappa \geqslant 0$,
$0 < \alpha < 1$. Topologically, we shall assume that $\Gamma = \zeta
[a, b]$ is either:
\begin{itemize}
\item[$\bullet$]
a {\sl Closed Jordan arc}, namely $\zeta : [a, b] \to \C$ is an
embedding;
\item[$\bullet$]
or a {\sl Jordan contour}, namely $\zeta : (a, b) \to \C$ is an
embedding, $\zeta ( a) = \zeta ( b)$, $\zeta$ extends as a
$\mathcal{ C}^{\kappa+1,
\alpha}$-smooth map on the quotient $[a, b]/ (a
\sim b)$ and $\Gamma = \zeta [ a, b ]$ is diffeomorphic to a circle.
\end{itemize}

Various more general assumptions can be made: $\Gamma$ consists of a
finite number of connected pieces, $\Gamma$ is piecewise smooth
(corners appear), $\Gamma$ possesses certain cusps, $\Gamma$ is only
Lipschitz, the length of $\Gamma$ is not finite, $f$ is $L^{\sf
p}$-integrable, $f$ is $L^{\sf p}_\alpha$, {\it i.e.} $f\in L^{\sf
p}(\Gamma)$ and $\int_\Gamma \, \left\vert f( s + h ) - f ( s)
\right\vert^{\sf p} \leqslant {\rm Cte} \, \left\vert h \right
\vert^\alpha$, $f$ belongs to certain Sobolev spaces, $f(\zeta) \, 
d\zeta$ is replaced by a measure $d\mu (\zeta)$, {\it etc.}, but
we shall not review the theory ({\it see}~\cite{ mu1953, 
ve1962, ga1966} and especially \cite{ dy1991}).

The natural orientation of the segment $[a, b]$ induced by the order
relation on $\R$ enables to orient the two semi-local sides of
$\Gamma$ in $\C$: the region on the left to $\Gamma$ will be called
the positive side (``$+$''), while the region to the right will be
called negative (``$-$''). In the case where $\Gamma$ is Jordan
contour, we assume that $\Gamma$ is oriented counterclockwise, so
that the positive region coincides with the bounded component of
$\C \backslash \Gamma$.

\def\thetheorem{2.7}\begin{theorem} {\rm (\cite{ mu1953,
ve1962, ga1966, dy1991, sme1988, ek2000}, [$*$])} 
Let $\Gamma$ be a $\mathcal{ C}^{\kappa+1,
\alpha}$-smooth closed Jordan arc or Jordan contour in $\C$ and let
$f: \Gamma \to \C$ be a $\mathcal{ C}^{\kappa, \alpha}$-smooth
complex-valued function.

\smallskip

\begin{itemize}
\item[{\bf (a)}]
If $\Gamma'$ is any closed portion of $\Gamma$ having no ends in
common with those of $\Gamma$, then for every $\zeta_1 \in \Gamma'$,
the Cauchy transform $F (z) := \frac{ 1}{2\pi i}\, \int_\Gamma \,
\frac{ f(\zeta) \, d\zeta}{\zeta - z}$ possesses {\rm (a priori
distinct)} limits $F^+ (\zeta_1)$ and $F^- (\zeta_1)$, when $z$ tends
to $\zeta_1$ from the positive or from the negative side.

\smallskip

\item[{\bf (b)}]
These two limits $F^+$ and $F^-$ are of class $\mathcal{ C}^{\kappa,
\alpha}$ on $\Gamma'$ with a norm estimate $\left\vert \! \left\vert
F^\pm \right\vert \! \right \vert_{ \mathcal{ C}^{\kappa, \alpha}
(\Gamma')} \leqslant \frac{C (\kappa, \Gamma', \Gamma)}{\alpha (1-
\alpha)} \,\left \vert \! \left\vert f \right\vert \! \right\vert_{
\mathcal{ C}^{ \kappa, \alpha} (\Gamma)}$, for some positive constant
where $C(\kappa, \Gamma', \Gamma)$.

\smallskip

\item[{\bf (c)}] 
Furthermore, if $\omega_+'$ and $\omega_-'$ denote an upper and a
lower open one-sided neighborhood $\Gamma'$ in $\C$, the two functions
$F^\pm: \omega_\pm ' \to \C$ defined by{\rm :}
\[
\left\{
\aligned
F^\pm(z)
&
:=
F(z)
\ \ \ \ \ \ \ \
\ \ \ \ \ 
{\rm if} 
\ \ 
z\in\omega_\pm', 
\\
F^\pm(z)
&
:=
F^\pm(\zeta_1)
\ \ \ \ \ \ \ \ \ \, 
{\rm if} 
\ \ 
z=\zeta_1\in\Gamma',
\endaligned\right.
\]
are of class $\mathcal{ C}^{\kappa, \alpha}$ in $\omega_\pm' \cup
\Gamma'$, with a similar norm estimate $\left\vert \! \left\vert
F^\pm \right\vert \! \right\vert_{ \mathcal{ C }^{\kappa, \alpha}
(\omega_\pm' \cup \Gamma')} \leqslant \frac{ C_1 ( \kappa, \Gamma',
\Gamma)}{ \alpha(1 -\alpha)} \,\left\vert \! \left\vert f \right\vert
\! \right\vert_{ \mathcal{ C }^{ \kappa, \alpha } (\Gamma)}$.

\smallskip

\item[{\bf (d)}]
Finally, at every point $\zeta_0$ of the curve $\Gamma$ not coinciding
with its ends, $F^+$ and $F^-$ satisfy the two {\rm
Sokhotski\u\i-Plemelj formulas:}
\[
\left\{
\aligned
F^+ (\zeta_0)
-
F^-(\zeta_0)
&
=
f(\zeta_0), 
\\
\frac{1}{2}\left[
F^+(\zeta_0)+F^-(\zeta_0)\right]
&
=
{\rm p.v.}\,
\frac{1}{2\pi i}
\int_\Gamma\,\frac{f(\zeta)}{\zeta-\zeta_0}\,d\zeta.
\endaligned
\right.
\]
\end{itemize}
\end{theorem}

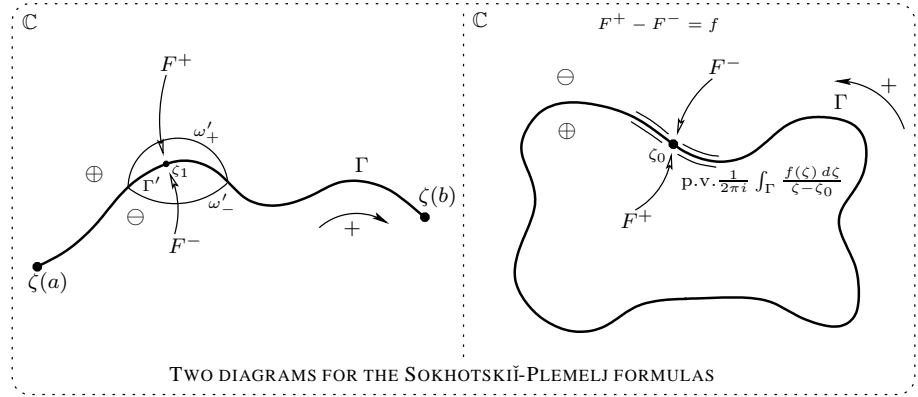
\begin{center}
\input sokhotskii-plemelj.pstex_t
\end{center}

Sometimes, $F$ is called {\sl sectionnally holomorphic}, as it is
discontinuous across $\Gamma$. Its jump across $\Gamma$ is provided by
the first formula above, while the arithmetic mean
$\frac{F^+ + F^- }{2}$ is given by the value of the Cauchy (singular)
integral at $\zeta_0 \in \Gamma$. Morera's classical theorem (\cite{
mo1889}) states that if $F^+$ and $F^-$ match up on the interior of
$\Gamma$, then the Cauchy integral is holomorphic in $\C$ minus the
endpoints of $\Gamma$. As is known (\cite{ sh1990}), this
theorem is also true for an arbitrary holomorphic function $F \in
\mathcal{ O} (\C \backslash \Gamma)$ which is not necessarily defined
by a Cauchy integral. 

\subsection*{2.8.~Less regular boundaries}
The boundary behaviour of the Cauchy transform at the two extreme
points $\gamma ( a)$ and $\gamma ( b)$ of a Jordan arc is studied
in~\cite{ mu1953}. We refer to~\cite{ dy1991} for a survey
presentation of the finest condition on $\Gamma$ (namely, it to be a
{\sl Carleson curve}) which insures that the Cauchy integral exists
and that the Sokhotski\u\i-Plemelj formulas hold true, almost
everywhere. Let us just mention what happens with the Cauchy integral
$F(z)$ in the limit case $\alpha = 0$.

If $\Gamma$ is (only) $\mathcal{ C}^{ 1, 0}$, if $f$ is (only)
$\mathcal{ C}^{ 0, 0}$, then for $\zeta_1$ in the interior of
$\Gamma$, {\it the limit $F^- (\zeta_1)$ exists if and only if the
limit $F^+ (\zeta_1)$ exists} (\cite{ mu1953}). However, generically,
none limit exists.

A more useful statement, valid in the case $\alpha = 0$, is as
follows. Assume $\Gamma$ to be $\mathcal{ C}^{ \kappa + 1, 0}$ with
$\kappa \geqslant 0$ and let $\Gamma'$ be a closed portion of the interior
of $\Gamma$. Parametrize $\Gamma'$ by a $\mathcal{ C}^{ \kappa+ 1,
0}$ map $\zeta' : [ a', b'] \to \Gamma'$. Extend $\zeta' = \zeta' (
s)$ as a a $\mathcal{ C}^{ \kappa + 1, 0}$ embedding $\zeta' (s,
\varepsilon)$ defined on $[ a', b'] \times (-\varepsilon_0,
\varepsilon_0)$, where $\varepsilon_0 >0$, with $\zeta' ( s, 0) \equiv
\zeta' (s)$ and with $\zeta' (s, \varepsilon)$ in the positive side of
$\Gamma'$ for $\varepsilon >0$. The family of curves
$\Gamma_\varepsilon ' := \zeta' ([ a', b'] \times \{ \varepsilon \})$
foliates a strip thickening of $\Gamma'$.

\def\thetheorem{2.9}\begin{theorem}
{\rm (\cite{ mu1953})} 
For every choice of a $\mathcal{ C}^{ \kappa + 1, 0}$ extension
$\zeta' (s, \varepsilon)$, and every $f\in \mathcal{ C}^{ \kappa, 0}
(\Gamma, \C)$, the difference from either side of the Cauchy transform
$F \vert_{ \Gamma_\varepsilon '} - F \vert_{ \Gamma_{
-\varepsilon'}'}$ tends to $f\vert_{ \Gamma'}$ in $\mathcal{ C}^{
\kappa, 0}$ norm as $\varepsilon \to 0${\rm :}
\[
\lim_{\varepsilon\to 0}\,
\sup_{s\in[a',b']}
\left\vert\!\left\vert
F(\zeta'(s,\varepsilon))
-
F(\zeta'(s,-\varepsilon))
-
f(\zeta'(s))
\right\vert\!\right\vert_{\kappa,0}
=0.
\]
\end{theorem}

To conclude, we state a criterion, due to Hardy-Littlewood, which
insures $\mathcal{ C}^{ \kappa, \alpha}$-smoothness
of holomorphic functions up to the boundary.

\def\thetheorem{2.10}\begin{theorem}
{\rm (\cite{ mu1953, ga1966}, [$*$])}
Let $\Gamma$ be a $\mathcal{ C}^{\kappa + 1, \alpha}$-smooth Jordan
contour, divinding the complex plane in two components $\Omega^+$ {\rm
(}bounded{\rm )} and $\Omega^-$ {\rm (}unbounded{\rm )}. If $f \in
\mathcal{ O} ( \Omega^\pm )$ satisfies the estimate $\vert
\partial_z^\kappa f (z) \vert \leqslant C \left( 1 - \vert z \vert
\right)^{1 - \alpha}$, for some $\kappa \in \N$, some $\alpha$ with $0
< \alpha < 1$, and some positive constant $C >0$, then $f$ is of class
$\mathcal{ C}^{\kappa, \alpha}$ in the closure $\overline{ \Omega^\pm
} = \Omega^\pm \cup \Gamma$.
\end{theorem}

\subsection*{2.11.~Functions and maps defined on the unit circle}
In the sequel, $\Omega$ will be the unit disc $\Delta := \{ \zeta \in
\C : \, \vert \zeta \vert < 1\}$ having as boundary the {\sl unit
circle} $\partial \Delta := \{\zeta \in \C : \vert \zeta \vert = 1\}$.
Consider a function $f: \partial \Delta \to \K$, where $\K= \R$ or
$\C$. Parametrizing $\partial \Delta$ by $\zeta = e^{ i\theta}$ with
$\theta \in \R$, such an $f$ will be considered as the function
\[
\R\ni\theta\longmapsto 
f(e^{ i\, \theta})\in\K.
\]
For $j\in \N$, we
shall write $f_{\theta^j} := \frac{ d^j f}{ d\theta^j }$.

Let $\alpha$ satisfy $0 < \alpha \leqslant 1$ and assume that $f\in
\mathcal{ C}^{0, \alpha}$. We define its $\mathcal{ C}^{ 0, \alpha}$
semi-norm (notice the
wide hat) precisely by:
\[
\left\vert\!\left\vert
f
\right\vert\!\right\vert_{\widehat{0,\alpha}}
:=
\sup_{\theta''\neq\theta'}\frac{\vert
f(e^{i\,\theta''})-f(e^{i\,\theta'})\vert}{ 
\left\vert \theta''-\theta'
\right\vert^\alpha}.
\]
Thanks to $2\pi$-periodicity, $\sup_{ \theta '' \neq \theta '}$ may be
replaced by $\sup_{ 0 < \left\vert \theta '' - \theta' \right \vert
\leqslant \pi}$. According to the definition of \S1.1,
the function $f$ is $\mathcal{ C}^{
\kappa, \alpha}$ if the quantity
\[
\left\vert\!\left\vert
f
\right\vert\!\right\vert_{\kappa,\alpha}
:=
\sum_{0\leqslant j\leqslant\kappa}\,
\left\vert\!\left\vert
f_{\theta^j}
\right\vert\!\right\vert_{0,0}
+
\left\vert\!\left\vert
f_{\theta^\kappa}
\right\vert\!\right\vert_{\widehat{0,\alpha}}
<\infty
\]
is finite. Besides H\"older spaces, we shall also consider the
Lebesgue spaces $L^{\sf p} (\partial \Delta)$, with ${\sf p} \in \R$
satisfying $1\leqslant {\sf p} \leqslant \infty$. As $\partial \Delta$ is
compact, the H\"older inequality entails the (strict) inclusions
$L^\infty (\partial \Delta) \subset L^{{\sf p}'} (\partial \Delta)
\subset L^{{\sf p}} (\partial \Delta) \subset L^1 (\partial \Delta)$,
for $1 < {\sf p} < {\sf p}' < \infty$.

\subsection*{ 2.12.~Fourier series of H\"older continuous functions}
If $f$ is at least of class $L^1$ on $\partial \Delta$, let
\[
\widehat{f}_k
:=
\frac{1}{2\pi i}
\int_{\partial \Delta}\,\zeta^{-k}\,
f(\zeta)\,
\frac{d\zeta}{\zeta}
\]
denote the $k$-th Fourier coefficient of $f$, where
$k\in \Z$. Given $n\in \N$, consider
the $n$-th partial sum of the Fourier series of $f$:
\[
{\sf F}_nf(e^{i\,\theta})
:=
\sum_{-n\leqslant k\leqslant n}\,
\widehat{f}_k\,e^{i\,k\,\theta}.
\]
We remind that Dini's (elementary) criterion:
\[
\int_0^\pi\,
\frac{ 
\left\vert
f(e^{i(\theta+t)})
+
f(e^{i(\theta-t)})
-
2\,f(e^{i\,\theta})
\right\vert}
{t}\,
dt
<
\infty
\]
insures the pointwise convergence $\lim_{ n\to \infty} \, {\sf F}_n
f (e^{ i\, \theta}) = f (e^{ i\, \theta})$. If $f$ is $\mathcal{ C}^{
0, \alpha}$ on $\partial \Delta$, with $0 < \alpha \leqslant 1$, the above
integral obviously converges at every $e^{ i\theta} \in \partial
\Delta$, so that we may identify $f$ with its (complete) Fourier
series:
\[
f(e^{i\,\theta})
=
{\sf F}f(e^{i\,\theta})
:=
\sum_{k\in\Z}\,
\widehat{f}_k\,e^{i\,k\,\theta}.
\]
In fact (\cite{ zy1959}), if $f\in \mathcal{ C}^{ \kappa, \alpha}$
with $\kappa \in \N$ and $0\leqslant \alpha\leqslant 1$, then $\big\vert
\widehat{ f}_k \big\vert \leqslant \frac{ \pi^{ 1+\alpha}}{ \vert k\vert^{
\kappa + \alpha}}\, \left\vert \! \left\vert f \right\vert \!
\right\vert_{ \widehat{ \kappa, \alpha}}$ for all $k\in \Z \backslash
\{ 0\}$. Also, if $f\in \mathcal{ C}^{ 0, \alpha}$, then $\sum_{ k\in
\Z} \, \big\vert \widehat{ f}_k \big \vert^c$ converges for $c >
\frac{ 2}{ 2\alpha +1}$. In 1913, Berste\u{\i}n proved absolute
convergence of $\sum_{ k\in \Z} \, \big\vert \widehat{ f}_k \big
\vert$ for $\alpha > 1/2$.

\subsection*{2.13.~Three Cauchy transforms in the unit disc} 
In the case $\Omega = \Delta$, our goal is to formulate Theorem~2.7
with more precision about the constant $C ( \kappa, \partial
\Omega)$. For $\eta \in \partial \Delta$ in the unit circle and $f \in
\mathcal{ C}^{\kappa, \alpha} (\partial \Delta, \C)$ with $\kappa \geqslant
0$, $0 < \alpha < 1$, as in \S2.6, we define:
\[
\aligned
{\sf C}^+f(\eta) 
&
:= 
\lim_{r\to 1^-}\,
\frac{1}{2\pi i}\,
\int_{\partial\Delta}
\frac{f(\zeta)}{
\zeta-r\eta}\,d\zeta,
\\
{\sf C}^0f(\eta)
&
:=
{\rm p.v.}\,\frac{1}{2\pi i}
\int_{\partial\Delta}\,
\frac{f(\zeta)}{\zeta-\eta}\,d\zeta,
\\
{\sf C}^-f(\eta) 
&
:= 
\lim_{r\to 1^+}\,
\frac{1}{2\pi i}\,
\int_{\partial\Delta}
\frac{f(\zeta)}{
\zeta-r\eta}\,d\zeta.
\endaligned
\]
The Sokhotski\u\i-Plemelj formulas hold: $f(\eta) = {\sf C}^+ f (\eta)
- {\sf C}^- f (\eta)$ and ${\sf C}^0 f (\eta) = \frac{ 1}{ 2} \left[
{\sf C}^+f(\eta)+{\sf C}^-f(\eta) \right]$. A theorem due to
Aleksandrov\footnote{We are grateful to Burglind J\"oricke who
provided the reference~\cite{ al1975}.} enables to obtain a precise
estimate of the $\mathcal{ C}^{ \kappa, \alpha}$ norms of these Cauchy
operators. To describe it, define:
\[
\mathcal{M}_0^\alpha
:=
\left\{
f\in\mathcal{C}^{0,\alpha}(\partial\Delta,\C):
\widehat{f}_0=0
\right\}.
\]
Then $\left\vert \! \left\vert \cdot \right\vert \! \right\vert_{
\widehat{ 0, \alpha }} $ is a norm on $\mathcal{ M}_0^\alpha$, since
only the constants $c$ satisfy $\left\vert \! \left\vert c \right\vert
\! \right\vert_{ \widehat{ 0, \alpha }} = 0$. For $p, q \in \R$ with
$0 < p, \, q < 1$, recall the definition $B (p, q) := \int_0^1 \, x^{
p-1} \, (1-x)^{ q-1} \, dx$ of the {\sl Euler beta function}.

\def\thetheorem{2.14}\begin{theorem}
{\rm (\cite{ al1975})}
The operator ${\sf C}^0 f (\eta) :=
{\rm p.v.} \, \frac{ 1}{ 2\pi i}\, \int_{
\partial \Delta} \, \frac{ f(\zeta)}{ \zeta - \eta} \, d\zeta$ is a
bounded linear endomorphism of $\mathcal{ M}_0^\alpha$ having norm{\rm
:}
\[
\big\vert\!\big\vert\!\big\vert
{\sf C}^0
\big\vert\!\big\vert\!\big\vert_{\widehat{0,\alpha}}
=
\frac{1}{2\pi}\,
B
\left(
\frac{\alpha}{2},\,
\frac{1-\alpha}{2}
\right).
\]
\end{theorem}

One may easily verify the two equivalences $B \left( \frac{ \alpha}{2},\,
\frac{1- \alpha }{2} \right) \sim \frac{ 2}{ \alpha}$ as $\alpha \to
0$ and $B \left( \frac{\alpha}{2},\, \frac{1- \alpha }{2} \right) \sim
\frac{ 2}{1- \alpha}$ as $\alpha \to 1$ as well as the two inequalities:
\[
\frac{1}{\alpha(1-\alpha)}
\leqslant
B\left( 
\frac{\alpha}{2},\,\frac{1-\alpha}{2}
\right)
\leqslant
\frac{4}{\alpha(1-\alpha)}.
\]
Thus, the nonremovable factor $\frac{ 1}{ \alpha ( 1- \alpha)}$ shows
what is the precise rate of explosion of the norm $\big\vert \!
\big\vert \! \big\vert {\sf C }^0 \big\vert \! \big\vert \!
\big\vert_{ \widehat{ 0, \alpha }}$ as $\alpha \to 0$ or
as $\alpha \to 1$.

Further, if $f \in \mathcal{ C}^{ 0, \alpha}$ does not necessarily
belong to $\mathcal{ M}_0^\alpha$, it is elementary to check that
$\big\vert \! \big\vert {\sf C}^0 f \big\vert \! \big\vert_{ 0,
0} \leqslant \frac{ C}{\alpha} \, \left\vert \! \left\vert f \right\vert \!
\right\vert_{ 0, \alpha}$, for some absolute constant $C>0$. It
follows that the (complete) operator norm $\big\vert \! \big\vert \!
\big\vert {\sf C }^0 \big\vert \! \big\vert \! \big\vert_{ 0,
\alpha }$ behaves like $\frac{ C}{ \alpha ( 1- \alpha)}$.

In conclusion, thanks to the Sokhotski\u\i-Plemelj
formulas ${\sf C}^+ f = \frac{ 1}{ 2} ( {\sf C}^0 f + f )$ and ${\sf
C}^- f = \frac{1}{ 2} ( {\sf C}^0 f - f)$, we deduce that there exists
an absolute constant $C_1> 1$ such that:
\[
\frac{1/C_1}{\alpha(1-\alpha)}\,
\leqslant
\big\vert\!\big\vert\!\big\vert
{\sf C}^{\sf b}
\big\vert\!\big\vert\!\big\vert_{0,\alpha}
\leqslant
\frac{C_1}{\alpha(1-\alpha)},
\]
where $0 < \alpha < 1$ and where
${\sf b} = -,0,+$.

Next, what happens with $f\in \mathcal{ C}^{ \kappa, \alpha}$, for
$\kappa \in \N$ arbitrary\,? For ${\sf b} = -, 0, +$, the ${\sf
C}^{\sf b}$ are bounded linear endomorphisms of $\mathcal{ C}^{
\kappa, \alpha}$ and similarly:

\def\thetheorem{2.15}\begin{theorem}
There exists an absolute constant $C_1 > 1$ such that 
if $\kappa \in \N$ and $0 < \alpha <1$,
for ${\sf b} = -, 0, +${\rm :}
\[
\frac{1/C_1}{\alpha(1-\alpha)}
\leqslant
\big\vert\!\big\vert\!\big\vert
{\sf C}^{\sf b}
\big\vert\!\big\vert\!\big\vert_{\kappa,\alpha}
\leqslant
\frac{C_1}{\alpha(1-\alpha)}.
\]
\end{theorem}

In other words, the constant $C_1$ is independent of $\kappa$. To
deduce this theorem from the estimates with $\kappa = 0$ (with
different absolute constant $C_1$) we proceed as follows, without
exposing all the rigorous details.

Inserting the Fourier series ${\sf F} ( f, e^{ i\, \theta})$ in the
integrals defining ${\sc C}^-$, 
${\sc C}^0$, ${\sc C}^+$ and integrating termwise (an operation which
may be justified), we get:
\[
\left\{
\aligned
{\sf C}^-f(e^{i\,\theta})
&
= 
-
\sum_{k<0}\,\widehat{f}_k\,
e^{i\,k\,\theta},
\\
{\sf C}^0f(e^{i\,\theta})
&
= 
-
\frac{1}{2}\,
\sum_{k<0}\,\widehat{f}_k\,
e^{i\,k\,\theta}
+
\frac{1}{2}\,
\widehat{f}_0
+
\frac{1}{2}\,
\sum_{k>0}\,\widehat{f}_k\,
e^{i\,k\,\theta},
\\
{\sf C}^+f(e^{i\,\theta})
&
=
\widehat{f}_0
+
\sum_{k>0}\,\widehat{f}_k\,
e^{i\,k\,\theta}.
\endaligned\right.
\]
If $\kappa \geqslant 1$, by differentiating termwise with respect to
$\theta$ these three Fourier representations of the 
${\sf C}^{\sf b}$, we see that
these operators commute with differentiation.

\def\thelemma{2.16}\begin{lemma}
For every $j\in \N$ with $0 \leqslant j\leqslant \kappa$ and for
${\sf b} = -, 0, +$, we have{\rm :}
\[
{\sf C}^{\sf b}
\left(
f_{\theta^j}
\right)
=
\left(
{\sf C}^{\sf b}f
\right)_{\theta^j}.
\]
\end{lemma}

Dealing directly with the principal value definition of ${\sf C}^0 f$,
another proof of this lemma for ${\sf C}^0$ would consist in
integrating by parts, deducing afterwards that ${\sf C}^-$ and ${\sf
C}^+$ enjoy the same property, thanks to the Sokhotski\u\i-Plemelj
formulas.

\smallskip

To establish Theorem~2.15, we introduce another auxiliary $\mathcal{
C}^{\kappa, \alpha}$ norm:
\[
\left\vert\!\left\vert
f
\right\vert\!\right\vert_{\kappa,\alpha}^\sim
:=
\sum_{0\leqslant j\leqslant \kappa}\,
\left\vert\!\left\vert
f_{\theta^j}
\right\vert\!\right\vert_{0,\alpha}
=
\left\vert\!\left\vert
f
\right\vert\!\right\vert_{\kappa,\alpha}
+
\sum_{0\leqslant j\leqslant\kappa-1}\,
\left\vert\!\left\vert
f_{\theta^j}
\right\vert\!\right\vert_{\widehat{0,\alpha}},
\]
which is equivalent to $\left\vert \! \left\vert \cdot \right\vert \!
\right\vert_{ \kappa, \alpha}$, thanks to the elementary inequalities
([$*$]):
\[
\left\vert\!\left\vert
f
\right\vert\!\right\vert_{\kappa,\alpha}
\leqslant
\left\vert\!\left\vert
f
\right\vert\!\right\vert_{\kappa,\alpha}^\sim
\leqslant
(1+\pi)
\left\vert\!\left\vert
f
\right\vert\!\right\vert_{\kappa,\alpha}.
\]
Notice that $\left\vert \! \left\vert \cdot \right\vert \!
\right\vert_{ 0, \alpha}^\sim = \left\vert \! \left\vert \cdot
\right\vert \! \right\vert_{ 0, \alpha}$. The next lemma applies to
${\sf L} = {\sf C}^-, {\sf C}^0, {\sf C}^+$ and to ${\sf L} = {\sf T}$, 
the Hilbert conjugation operator defined below.

\def\thelemma{2.17}\begin{lemma}
{\rm ([$*$])}
Let ${\sf L}$ be a bounded linear endomorphism of all the spaces
$\mathcal{ C}^{ \kappa, \alpha} ( \partial \Delta, \C)$
with $\kappa \in
\N$, $0 < \alpha < 1$, which commutes with differentiations, namely
${\sf L} \left( f_{\theta^j} \right) = \left( {\sf L} f
\right)_{\theta^j}$, for $j \in \N$. Assume that there
exist a contant $C_1 (\alpha) > 1$ depending on
$\alpha$ such that $C_1 (\alpha)^{-1} \leqslant \left\vert \! \left\vert \!
\left\vert {\sf L} \right\vert \! \right\vert \! \right\vert_{ 0,
\alpha} \leqslant C_1 (\alpha)$. Then for every $\kappa \in \N${\rm :}
\[
\left\{
\aligned
C_1(\alpha)^{-1}
\leqslant
&
\left\vert\!\left\vert\!\left\vert
{\sf L} \right\vert\!\right\vert\!\right\vert_{\kappa,
\alpha}^\sim \leqslant C_1(\alpha),
\\
(1+\pi)^{-1}\,C_1(\alpha)^{-1}
\leqslant
&
\left\vert\!\left\vert\!\left\vert
{\sf L}
\right\vert\!\right\vert\!\right\vert_{\kappa,\alpha}
\leqslant
(1+\pi)\,
C_1(\alpha).
\endaligned\right.
\]
\end{lemma}

\proof
Indeed, if $f \in \mathcal{ C}^{\kappa, \alpha}$, we 
develope a chain of (in)equalities:
\[
\aligned
\left\vert\!\left\vert
{\sf L}f
\right\vert\!\right\vert_{\kappa,\alpha}
&
\leqslant
\left\vert\!\left\vert
{\sf L}f
\right\vert\!\right\vert_{\kappa,\alpha}^\sim
=
\sum_{0\leqslant j\leqslant \kappa}\,
\left\vert\!\left\vert
({\sf L}f)_{\theta^j}
\right\vert\!\right\vert_{0,\alpha}
=
\sum_{0\leqslant j\leqslant \kappa}\,
\left\vert\!\left\vert
{\sf L}(f_{\theta^j})
\right\vert\!\right\vert_{0,\alpha}
\\
&
\leqslant
C_1(\alpha)\,
\sum_{0\leqslant j\leqslant\kappa}\,
\left\vert\!\left\vert
f_{\theta^j}
\right\vert\!\right\vert_{0,\alpha}
=
C_1(\alpha)\,
\left\vert\!\left\vert
f
\right\vert\!\right\vert_{\kappa,\alpha}^\sim
\\
&
\leqslant
(1+\pi)\,C_1(\alpha)\,
\left\vert\!\left\vert
f
\right\vert\!\right\vert_{\kappa,\alpha}.
\endaligned
\]
This yields the two majorations.
Minorations are obtained similarly.
\endproof

To conclude this paragraph, we state a T{\oe}plitz type theorem about
${\sf C}^+$, which will be crucial in solving Bishop's equation with
optimal loss of smoothness, as we will see in Section~3. A similar one
holds about ${\sf C}^-$, assuming $\phi \in H^\infty ( \overline{ \C}
\backslash \overline{ \Delta})$ instead, where $\overline{ \C}$ is the
Riemann sphere.

\def\thetheorem{2.18}\begin{theorem}
{\rm (\cite{ tu1994b}, [$*$])} There exists an
absolute constant $C_1 > 1$ such that for all $f \in \mathcal{ C}^{
\kappa, \alpha}$, $\kappa \in \N$, $0 < \alpha < 1$, and all $\phi \in
{\rm H }^\infty ( \Delta) := \mathcal{ O} ( \Delta) \cap L^\infty (
\Delta)${\rm :}
\[
\left\vert\!\left\vert
{\sf C}^+(f\overline{\phi})
\right\vert\!\right\vert_{\kappa, \alpha}
\leqslant
\frac{C_1}{\alpha(1-\alpha)}\,
\left\vert\!\left\vert
f
\right\vert\!\right\vert_{\kappa,\alpha}\,
\left\vert\!\left\vert
\phi
\right\vert\!\right\vert_{L^\infty}.
\]
\end{theorem}

Closely related to the Cauchy transform are the Schwarz and the
Hilbert transforms.

\subsection*{ 2.19.~Schwarz transform on the unit disc} 
Let $u \in L^1 ( \partial \Delta, \R)$ be real-valued. The {\sl
Schwarz transform} of $u$ is the function of $z \in \Delta$ defined
by:
\[
{\sf S} u (z) := 
\frac{ 1}{ 2\pi i}
\int_{\partial\Delta}
u(\zeta)
\left(
\frac{\zeta+z}{\zeta-z}
\right)
\frac{d\zeta}{\zeta}.
\]
Thanks to the holomorphicity of the kernel, ${\sf S} u(z)$ is a
holomorphic function of $z \in \Delta$. Decomposing it in real and
imaginary parts:
\[
{\sf S}u(z)
= 
{\sf P}u(z,\bar z)
+
i{\sf T}u(z,\bar z),
\]
we get the {\sl Poisson transform} of $u$:
\[
{\sf P}u(z,\bar z)
:=
\frac{1}{2\pi i}
\int_{\partial\Delta}\,
u(\zeta)\,{\rm Re}
\left(
\frac{\zeta+z}{\zeta-z}
\right)
\frac{d\zeta}{\zeta},
\]
together with the {\sl Hilbert transform} of $u$:
\[
{\sf T}u(z,\bar z)
:=
\frac{1}{2\pi i}
\int_{\partial\Delta}\,
u(\zeta)\,{\rm Im}
\left(
\frac{\zeta+z}{\zeta-z}
\right)
\frac{d\zeta}{\zeta}.
\]
Thanks to the harmonicity of the two kernels, ${\sf P} u$ and ${\sf T}
u$ are harmonic in $\Delta$. The power series of ${\sf C} u$, of ${\sf
P} u$ and of ${\sf T} u$ are given by:
\[
\left\{
\aligned
{\sf S}u(z) 
&
=
\widehat{u}_0
+
2\,\sum_{k>0}\,\widehat{u}_k\,z^k, 
\\
{\sf P}u(z,\bar z) 
&
=
\sum_{k<0}\,\widehat{u}_k\,\bar z^k
+
\widehat{u}_0
+
\sum_{k>0}\,\widehat{u}_k\,z^k,
\\
{\sf T}u(z,\bar z)
&
=
\frac{1}{i}\left(
-
\sum_{k<0}\,\widehat{u}_k\,\bar z^k
+
\sum_{k>0}\,\widehat{u}_k\,z^k
\right),
\endaligned\right.
\]
where $\widehat{ u}_k$ is the $k$-th Fourier coefficient of $u$. These
three series converge normally on compact subsets of $\Delta$.

\subsection*{2.20.~Poisson transform on the unit disc}
Let us first summarize the basic properties of the Poisson transform
(\cite{ ka1968, dr2002}). Setting $z = r\, e^{i\, \theta}$ with $0
\leqslant r <1$ and $\zeta = e^{ it}$, computing ${\rm Re}\, \big(
\frac{ \zeta+ z}{ \zeta - z} \big)$ and switching the convolution
integral, we obtain:
\[
{\sf P}u(r\,e^{i\,\theta})
=
\frac{1}{2\pi}\int_{-\pi}^\pi\,
P_r(t)\,u(e^{i(\theta-t)})\,dt
=
P_r*u\,(e^{i\,\theta}),
\]
where 
\[
P_r(t):=
\frac{1-r^2}{1-2r\,\cos t+r^2} 
\]
is the {\sl Poisson summability kernel}. It has
three nice properties:

\begin{itemize}
\smallskip\item[$\bullet$]
$P_r>0$ on $\partial \Delta$ for $0\leqslant r < 1$,
\smallskip\item[$\bullet$]
$\frac{1}{ 2\pi }\, \int_{ -\pi }^\pi\, P_r(t) \,dt = 1$
for $0 \leqslant r< 1$, and:
\smallskip\item[$\bullet$]
$\lim_{r\to1^-}\,P_r(t)=0$ for every $t\in [-\pi, \pi] 
\backslash \{ 0\}$.
\end{itemize}

\smallskip
\noindent
Consequently, $P_r$ is an approximation of the Dirac measure
$\delta_1$ at $1 \in \partial \Delta$. For this reason, the Poisson
convolution integral possesses excellent boundary value
properties.

\def\thelemma{2.21}\begin{lemma}
{\rm (\cite{ ka1968, dr2002})} Convergence in norm holds{\rm :}
\begin{itemize}
\smallskip\item[{\bf (i)}]
If $u\in L^{ \sf p}$ with $1\leqslant {\sf p} < \infty$ or $p =
\infty$ and $u$ is continuous, then $\lim_{ r\to 1^-}\, \left\vert \!
\left\vert P_r * u - u \right\vert \! \right\vert_{ L^{\sf p}} = 0$.
\smallskip\item[{\bf (ii)}]
If $u\in \mathcal{ C}^{\kappa, \alpha}$ with $\kappa \in \N$ and $0
\leqslant \alpha \leqslant 1$, including $\alpha = 0$ and $\alpha =
1$, then $\lim_{ r\to 1^-}\, \left\vert \! \left\vert P_r * u - u
\right\vert \! \right\vert_{ \kappa, \alpha} = 0$.
\end{itemize}
\end{lemma}

\smallskip
In $\mathcal{ C }^{ \kappa, \alpha}$, the pointwise convergence
$\lim_{ r\to 1^-}\, P_r * u (e^{ i\, \theta}) \to u ( e^{ i\,
\theta})$ follows obviously. However, in $L^{\sf p}$, from
convergence in norm one may only deduce pointwise convergence almost
everywhere for some sequence $r_k \to 1$ which depends on the
function. In $L^{\sf p}$, almost everywhere pointwise convergence was
proved by Fatou in 1906.

\def\thetheorem{2.22}\begin{theorem}
{\rm (\cite{ fa1906, ka1968, dr2002})}
If $u \in L^{\sf p}$ with $1\leqslant {\sf p} \leqslant \infty$, then
for almost every $e^{i\, \theta} \in \partial \Delta$, 
we have{\rm :}
\[
\lim_{r\to 1^-}\,
P_r*u\,(e^{i\,\theta})
=
u(e^{i\,\theta}).
\]
\end{theorem}

In summary, the Poisson transform ${\sf P} u$ yields a harmonic
extension to $\Delta$ of any function $u\in L^{\sf p} ( \partial
\Delta, \R)$ or $u\in \mathcal{ C}^{\kappa, \alpha} (\partial \Delta,
\R)$, with expected
boundary value ${\sf b}_{\partial \Delta} ( {\sf P} u) = u$
on $\partial \Delta$.

\subsection*{2.23.~Hilbert transform on the unit disc}
Next, we survey the fundamental properties of the Hilbert transform.
Again, $u$ is real-valued on $\partial \Delta$. Setting $z = r\,
e^{i\, \theta}$ with $0 \leqslant r < 1$ and $\zeta = e^{ it}$,
computing ${\rm Im}\, \big( \frac{ \zeta+ z}{ \zeta - z} \big)$ and
switching the convolution integral, we obtain:
\[
{\sf T}u(r\,e^{i\,\theta})
=
\frac{1}{2\pi}\int_{-\pi}^{\pi}\,
T_r(t)\,u(e^{i(\theta-t)})\,dt,
\]
where
\[
T_r(t)
:=
\frac{2\,r\,\sin t}{
1-2r\,\cos t+r^2}
\]
is the {\sl Hilbert kernel}. It is {\it not}\, a summability kernel,
being positive and negative with $L^1$ norm tending to $\infty$ as $r
\to 1^-$; for this reason, the Hilbert transform does not enjoy the
same nice boundary value properties as the Poisson transform:
H\"older classes are needed. 

Setting $r = 1$, the Poisson kernel ${\sf P}_1 (t)$ vanishes
identically and the Hilbert kernel tends to $\frac{ 2\, \sin t}{ 2 -
2\, \cos t} = \frac{ \cos t/2}{ \sin t/2}$. Near $t=0$, the function
$\cot (t/2)$ behaves like the function $2/t$, having
infinite $L^1$ norm. For
$u\in \mathcal{ C}^{ 0, \alpha} ( \partial \Delta, \R)$, it may be
verified that, as $z \to e^{ i\, \theta} \in \partial \Delta$, the
Hilbert transform ${\sf T} u (z)$ tends to
\[
\aligned
{\sf T}u(e^{i\,\theta})
:=
&\
{\rm p.v.}\,
\frac{1}{2\pi}\,
\int_{-\pi}^\pi\,
\frac{u(e^{i(\theta-t)})}{\tan(t/2)}\,dt
\\
=
&\
{\rm p.v.}\,
\frac{1}{2\pi i}\,
\int_{-\pi}^\pi\,
u(\zeta)\,
{\rm Im}\,
\Big(
\frac{\zeta+e^{i\,\theta}
}{
\zeta-e^{i\,\theta}} 
\Big)\,
\frac{d\zeta}{\zeta}.
\endaligned
\]
Since ${\rm Re} \big( \frac{ \zeta + e^{ i\, \theta}}{ \zeta - e^{ i\,
\theta}} \big) \equiv 0$ for $\zeta = e^{ i\, t} \in \partial \Delta$,
we get ${\rm Im} \, \big( \frac{ \zeta + e^{ i\, \theta}}{ \zeta - e^{
i\, \theta}} \big) = \frac{ 1}{ i}\, \frac{ \zeta + e^{ i\, \theta}}{
\zeta - e^{ i\, \theta}}$ so that we may rewrite
\[
i\,{\sf T}u(e^{i\,\theta})
=
{\rm p.v.}\,
\frac{1}{2\pi i}\int_{-\pi}^\pi\,
u(\zeta)\,
\frac{\zeta+e^{i\,\theta}
}{
\zeta-e^{i\,\theta}}\,
\frac{d\zeta}{\zeta}.
\]
Setting ${\sf P}_0 u := \frac{ 1}{ 2\pi i}\, \int_{ \partial \Delta}\,
u(\zeta) \frac{ d\zeta}{ \zeta} = \widehat{ u}_0$, the algebraic
relation $\frac{ 2}{ \zeta - e^{ i\, \theta}} - \frac{ 1}{ \zeta} =
\frac{\zeta+e^{i\,\theta} }{ \zeta-e^{i\,\theta}}\, \frac{ 1}{ \zeta}$
gives a fundamental relation between ${\sf C}^0$ and ${\sf T}$:
\[
2\,{\sf C}^0-{\sf P}_0
=
i\,{\sf T}.
\]
From Theorem~2.15, we deduce (${\sf P}_0$ is innocuous):

\def\thetheorem{2.24}\begin{theorem}
{\rm (\cite{ pri1916, hita1978, bo1991, ber1999}, [$*$])} There exist
an absolute constant $C_1 > 1$ such that if $\kappa \in\N$ and $0 <
\alpha < 1${\rm :}
\[
\frac{1/C_1}{\alpha(1-\alpha)}
\leqslant
\left\vert\!\left\vert\!\left\vert
{\sf T}
\right\vert\!\right\vert\!\right\vert_{\kappa,\alpha}
\leqslant
\frac{C_1}{\alpha(1-\alpha)}.
\]

\end{theorem}

It follows that at the level of Fourier series, ${\sf T}$
transforms $u ( e^{i\, \theta}) = {\sf F} u (e^{i\, \theta}) = \sum_{
k\in \Z} \, \widehat{ u}_k \, e^{ i\, k\, \theta}$ to
\[
{\sf T}u(e^{i\,\theta})
:=
\frac{1}{i}
\left(
-
\sum_{k<0}\,\widehat{u}_k\,e^{i\,k\,\theta}
+
\sum_{k>0}\,\widehat{u}_k\,e^{i\,k\,\theta}
\right).
\]
Notice that $(\widehat{ {\sf T} u})_0 = 0$. In fact, this formula
coincides with the series $\frac{ 1}{ i}\left( -\sum_{ k< 0}\,
\widehat{ u}_k \bar z^k + \sum_{ k>0}\, \widehat{ u}_k z^k \right)$,
written for $z \to e^{ i\, \theta}$, the limit existing provided $0 <
\alpha < 1$.

By termwise differentiation of the above formula, ${\sf
T} (u_{ \theta^j }) = ({\sf T} u)_{\theta^j}$ for $0 \leqslant
j\leqslant \kappa$, if $u\in \mathcal{ C}^{ \kappa, \alpha}$ (some
integrations by parts in the singular integral defining ${\sf T} u$
would yield a second proof of this property).

The Poisson transform ${\sf P} u$ of $u \in \mathcal{ C}^{ 0,\alpha}$
having boundary value ${\sf b}_{\partial \Delta} ({\sf P} u) = u$ and
the Schwarz transform being holomorphic in $\Delta$, we see that the
function $u + i \, {\sf T} u$ on $\partial \Delta$ extends
holomorphically to $\Delta$ as ${\sf S} u (z)$. So ${\sf T} u$ on
$\Delta$ is one of the Harmonic conjugates of $u$. In general, these
conjugates are defined up to a constant. The property $( \widehat{
{\sf T} u})_0 = 0$ means that ${\sf T} u (0) = 0$.

\def\thelemma{2.25}\begin{lemma}
The Hilbert transform ${\sf T} u$ on $\partial \Delta$
is the boundary value on $\partial \Delta$ of the
{\rm unique harmonic conjugate} in $\Delta$ of the harmonic
Poisson extension ${\sf P} u$, that vanishes at $0 \in \Delta$.
$$
\boxed{
\text{\rm For} \ 
u\in\mathcal{C}^{\kappa,\alpha}(\partial\Delta,\R), \
u+i\,{\sf T}u
\
\text{\rm extends holomorphically to} 
\
\Delta.
}
$$
Furthermore, ${\sf T} ({\sf T} u) = - u + \widehat{ u}_0$.
\end{lemma}

\subsection*{ 2.26.~Hilbert transform in $L^{\sf p}$ spaces} 
It is elementary to show that the study of the
principal value integral 
${\rm p.v.}\, \frac{ 1}{ 2\pi}\, 
\int_{ -\pi}^{ \pi}\,
\frac{ u( e^{ i\, (\theta -t)})}{
\tan (t/2)}\, dt$ 
is equivalent to the study of
the same singular convolution operator, in which $\cot (t/2)$ is replaced
by $2/t$. Similarly, one may define the {\sl Hilbert transform on the
real line}:
\[
{\sf H}f({\sf x})
:=
{\rm p.v.}\int_\R\,
\frac{f({\sf y})}{{\sf y}-{\sf x}}\,d{\sf y}.
\]
If $f$ is $\mathcal{ C}^{ 1, 0}$ on $\R$ and has compact support or
satisfies $\int_\R \, \vert f \vert < \infty$, replacing $f({\sf y})$
in the numerator by $[ f ({\sf y}) - f ({\sf x}) ] + f({\sf x})$ and
reasoning as in \S2.3, one straightforwardly shows the existence of
the above principal value.

Privalov showed that ${\sf H} f ({\sf x})$ exists for almost every
${\sf x} \in \R$ if $f \in L^1 (\R)$. A theorem due to M. Riesz
states that the two Hilbert transforms ${\sf H}$ on the real line and
${\sf T}$ on the unit circle are bounded endomorphisms of $L^{\sf p}$,
for $1 < {\sf p} < \infty$, namely if $f \in L^{ \sf p} (\R)$ and $u
\in L^{\sf p} (\partial \Delta)$, then:
\[
\left\vert\!\left\vert
{\sf H}f
\right\vert\!\right\vert_{L^{\sf p}(\R)}
\leqslant
C_{\sf p}
\left\vert\!\left\vert
f
\right\vert\!\right\vert_{L^{\sf p}(\R)}
\ \ \ \ \ \ \ \ \ \ \ \
{\rm and}
\ \ \ \ \ \ \ \ \ \ \ \
\left\vert\!\left\vert
{\sf T}u
\right\vert\!\right\vert_{L^{\sf p}(\partial\Delta)}
\leqslant
C_{\sf p}
\left\vert\!\left\vert
u
\right\vert\!\right\vert_{L^{\sf p}(\partial\Delta)},
\]
whith the same constant $C_{\sf p}$ (\cite{ zy1959}, Chapters~VII
and~XVI). In~\cite{ pi1972}, Zygmund's doctoral student Pichorides
obtained the best value of the constant $C_{\sf p}$: for $1< {\sf p}
\leqslant 2$, $C_{\sf p} = \tan \frac{ \pi }{ 2 {\sf p}}$, while, by a
duality argument, $C_{\sf p} = \cot \frac{ \pi}{ 2 {\sf p}}$ for
$2\leqslant {\sf p} < \infty$. The two elementary bounds $\tan \frac{
\pi }{ 2 {\sf p}} \leqslant \frac{ {\sf p}}{ {\sf p}-1}$ for $1 < {\sf
p} \leqslant 2$ and $\cot \frac{ \pi}{ 2{\sf p}} \leqslant {\sf p}$
for $2\leqslant {\sf p} < \infty$, yield:
\def\theequation{2.27}\begin{equation}
\left\vert\!\left\vert
{\sf H}f
\right\vert\!\right\vert_{L^{\sf p}(\R)}
\leqslant
\frac{{\sf p}^2}{{\sf p}-1}\,
\left\vert\!\left\vert
f
\right\vert\!\right\vert_{L^{\sf p}(\R)}
\ \ \ \ \ \ \ \ \
{\rm and}
\ \ \ \ \ \ \ \ \
\left\vert\!\left\vert
{\sf T}f
\right\vert\!\right\vert_{L^{\sf p}(\partial\Delta)}
\leqslant
\frac{{\sf p}^2}{{\sf p}-1}\,
\left\vert\!\left\vert
f
\right\vert\!\right\vert_{L^{\sf p}(\partial\Delta)},
\end{equation}
for $1 < {\sf p} < \infty$.
In $L^1$, the Hilbert transform is unbounded but, according to a
theorem due to Kolmogorov (\cite{ dy1991, dr2002}), it
sastisfies a {\sl weak inequality}:
\[
\mathfrak{m}
\left\{
{\sf H}f({\sf x})
>
a
\right\}
\leqslant 
\frac{C}{a}\,
\left\vert\!\left\vert
f
\right\vert\!\right\vert_{L^1},
\]
for every $a\in \R$ with $a>0$, where $\mathfrak{ m}$ is the Lebesgue
measure and where $C>0$ is some absolute constant.

\subsection*{ 2.28.~Pointwise convergence of
Fourier series} The boundedness of the Hilbert transform in $L^{ \sf
p}$ has a long history, closely related to the problem of pointwise
convergence of Fourier series. In 1913, before M.~Riesz
proved the estimates~\thetag{ 2.27}, using complex function theory and the
Riesz-Fischer theorem, Luzin showed that ${\sf H}$ is bounded in $L^2$
and formulated the celebrated conjecture that Fourier series of $L^2$
functions converge pointwise almost everywhere. This ``hypothetical
theorem'' was established by Carleson (\cite{ ca1966}) in 1966 and
slightly later by Hunt (\cite{ hu1966}) in $L^{\sf p}$ for $1 < {\sf
p} < \infty$. A complete self-contained restitution of these results
is available in~\cite{ dr2002}. Let us survey the main theorem.

The $n$-th partial sum of the Fourier series of a function $f$ on
$\partial \Delta$ is given by:
\[
{\sf F}_nf(e^{i\,\theta})
=
\frac{1}{2\pi}\,\int_{-\pi}^\pi\,
D_n(t)\,f(e^{i(\theta-t)})\,dt,
\]
where
\[
D_n(t)
:=
\frac{\sin (n+1/2)t}{\sin t/2}
\]
is the {\em Dirichlet kernel}, having unbounded $L^1$ norm $\left\vert
\! \left\vert D_n \right\vert \! \right\vert_{L^1} \sim \frac{
4}{\pi^2} \, {\rm log} \, n$. It is elementary to show that the
behaviour of this convolution integral, as $n\to \infty$, is
equivalent to the behaviour of the integral:
\[
\int_{-\pi}^\pi\,
\frac{\sin nt}{t}\,f(e^{i(\theta-t)})\,dt.
\]
Without loss of generality, $f$ is assumed to be real-valued, so that
the above integral is the imaginary part of the {\sl Carleson
integral}:
\[
{\sf C}_n(f,e^{i\,\theta})
:=
{\rm p.v.}\int_{-\pi}^\pi\,
\frac{e^{i\,n\,t}}{t}\,
f(e^{i(\theta-t)})\,dt.
\]
In, chapters~4, 5, 6, 7, 8, 9 and~10 of~\cite{ dr2002}, the main
proposition is to prove that the {\sl Carleson maximal sublinear
operator}:
\[
{\sf C}^*f(e^{i\,\theta})
:=
\sup_{n\in\N}\,
\left\vert
{\sf C}_n(f,e^{i\,\theta})
\right\vert
\]
is bounded from $L^{ \sf p}$ to $L^{ \sf p}$. The proof involves
dyadic partitions, changes of frequency, microscopic Fourier analysis
of $f$, choices of allowed pairs and seven exceptional sets. By an
elementary argument, one deduces that the {\sl maximal Fourier series
sublinear operator}:
\[
{\sf F}^*f(e^{i\,\theta})
:=
\sup_{n\in\N}\,
\left\vert
{\sf F}_nf(e^{i\,\theta})
\right\vert
\]
is bounded from $L^{ \sf p}$ to $L^{ \sf p}$. 

\def\thetheorem{2.29}\begin{theorem}
{\rm (\cite{ ca1966, hu1966, dr2002})} If $f \in L^{
\sf p}$ with $1 < {\sf p} < \infty$, there exists an absolute constant
$C>1$ such that{\rm :}
\[
\left\vert\!\left\vert
{\sf F}^*f
\right\vert\!\right\vert_{L^{\sf p}}
\leqslant
C\,\frac{{\sf p}^4}{({\sf p}-1)^3}\,
\left\vert\!\left\vert
f
\right\vert\!\right\vert_{L^{\sf p}}.
\]
\end{theorem}

Then by a standard argument, $\lim_{ n\to \infty}\, {\sf F}_n f (e^{
i\, \theta}) = f ( e^{ i\, \theta})$ almost everywhere.

\subsection*{ 2.30.~Transition}
Since the grounding article~\cite{ hita1978}, the nice behaviour of the
Hilbert transform in the H\"older classes (Theorem~2.24) is the main
reason why Bishop analytic discs have been constructed in the category
of $\mathcal{ C}^{\kappa, \alpha}$ generic submanifolds of $\C^n$
(\cite{ bpo1982, bpi1985, tu1990, trp1990, bo1991, brt1994, tu1994a,
me1994, trp1996, jo1996, ber1999}). Perhaps it is also interesting to
construct Bishop analytic discs in the Sobolev classes.

\section*{ \S3.~Solving a local parametrized
Bishop equation with optimal loss of smoothness} 

\subsection*{ 3.1.~Analytic discs attached
to a generic submanifold of $\C^n$} As in
Theorem~4.2(III), let
$M$ be a $\mathcal{ C}^{ \kappa, \alpha}$ local graphed generic
submanifold of equation $v = \varphi (x, y, u)$, where $\varphi$ is
defined for $\vert x + i \, y \vert < \rho_1$, $\vert u \vert <
\rho_1$, for some $\rho_1 >0$ and where $\varphi (0) = 0$, $d
\varphi (0) = 0$ and $\vert \varphi \vert < \rho_1$.

\def\thedefinition{3.2}\begin{definition}{\rm
An {\sl analytic disc} is a map
\[
\overline{ \Delta} \ni \zeta 
\longmapsto 
A ( \zeta) = ( Z (\zeta), W (\zeta))
\in \C^m \times \C^d 
\]
which is holomorphic in
the unit disc $\Delta$ and at least $\mathcal{ C}^{ 0, 0}$
in $\overline{ \Delta }$. It is {\sl attached to} $M$ if it sends
$\partial \Delta$ into $M$.

}\end{definition}
Thus, suppose that $( Z ( \zeta), W ( \zeta) )$ is attached to $M$ and
sufficiently small, namely $\vert [X + i \, Y] ( e^{ i\, \theta} )
\vert < \rho_1$, $\vert U ( e^{ i\, \theta} ) \vert < \rho_1$ and
$\vert V ( e^{ i\, \theta} ) \vert < \rho_1$ on $\partial \Delta$,
where $Z (\zeta) = X (\zeta) + i Y (\zeta)$ and $W ( \zeta) = U
(\zeta) + i V ( \zeta)$. Then clearly, the disc sends $\partial
\Delta$ to $M$ if and only if
\[
V(e^{i\, \theta})
=
\varphi
\left(X(e^{i\, \theta}),Y(e^{i\, \theta}),U(e^{i\, \theta})\right),
\]
for every $e^{i\, \theta} \in \partial \Delta$.
Thanks to the Hilbert transform, we claim that we may express
analytically the fact that the disc is attached to $M$.

At first, in order to guarantee the applicability of the harmonic
conjugation operator ${\sf T}$, all our analytic discs will $\mathcal{
C}^{ \kappa, \alpha}$ on $\overline{ \Delta}$, with $\kappa \in \N$
and $0 < \alpha < 1$. We let ${\sf T}$ act componentwise on maps $U =
(U^1, \dots, U^d) \in \mathcal{ C}^{ \kappa, \alpha} ( \partial
\Delta, \R^d)$, namely ${\sf T} U := \left( {\sf T} U^1, \dots, {\sf
T} U^d \right)$. We set $\left\vert \! \left\vert {\sf T} U
\right\vert \! \right\vert_{ \kappa, \alpha} := \max_{ 1\leqslant j
\leqslant d}\, \left\vert \! \left\vert {\sf T} U^j \right\vert \!
\right\vert_{ \kappa, \alpha}$. With a slight change of notation,
instead of ${\sf P} U (0)$, we denote by ${\sf P}_0 \, U := \frac{ 1}{
2\pi}\, \int_{ - \pi }^{ \pi} \, U ( e^{ i\, \theta} ) \, d\theta$ the
value at the origin of the Poisson extension ${\sf P}
U$. Equivalently, ${\sf P}_0 \, U = \widehat{ U}_0$ is the mean value
of $U$ on $\partial \Delta$. Here is a summary of the most useful
properties of ${\sf T}$.

\def\thelemma{3.3}\begin{lemma}
The $\R^d$-valued Hilbert transform ${\sf T}$ is a bounded linear
endomorphism of $\mathcal{ C}^{ \kappa, \alpha} ( \partial \Delta,
\R^d)$ with $\frac{ 1/C_1}{ \alpha ( 1 - \alpha)} \leqslant \left\vert
\! \left\vert \! \left\vert {\sf T} \right\vert \! \right\vert \!
\right\vert_{ \kappa, \alpha} \leqslant \frac{ C_1}{ \alpha ( 1 -
\alpha)}$ satisfying ${\sf T} ({\rm cst}) = 0$ and
\[
{\sf T}({\sf T}U)
=
-
U
+
{\sf P}_0\,U.
\]
\end{lemma}

In the sequel, we shall rather use the mild modification ${\sf T}_1$
of ${\sf T}$ defined by:
\[
{\sf T}_1U(e^{i\,\theta})
:=
{\sf T}U(e^{i\,\theta})
-
{\sf T}U(1).
\]
In fact, ${\sf T}_1$ is uniquely determined by the normalizing
condition ${\sf T}_1 U (1) = 0$. Then ${\sf T}_1$ is also bounded:
$\frac{ 1/C_1}{ \alpha ( 1 - \alpha)} \leqslant \left\vert \! \left\vert \!
\left\vert {\sf T}_1 \right\vert \! \right\vert \! \right\vert_{
\kappa, \alpha} \leqslant \frac{ C_1}{ \alpha ( 1 - \alpha)}$, also
annihilates constants: ${\sf T}_1 ({\rm cst}) = 0$ and
\[
{\sf T}_1({\sf T}_1U)
=
-
U
+
U(1).
\]
Furthermore, most importantly: 

\def\thelemma{3.4}\begin{lemma}
If $U \in \mathcal{ C}^{ \kappa, \alpha} ( \partial \Delta, \R^d)$,
then $U ( e^{ i\, \theta }) + i\, {\sf T}_1 U (e^{i \, \theta })$
extends as a holomorphic map $\Delta \to \C^d$ which is $\mathcal{
C}^{ \kappa, \alpha}$ in the closed disc $\overline{ \Delta }$.
\end{lemma}

To check that the extension is $\mathcal{ C}^{ \kappa, \alpha}$ in
$\overline{ \Delta }$, one may introduce the Poisson integral
formula and apply Lemma~2.21 {\bf (ii)}.

\smallskip

If $A = (Z, W)$ is an analytic disc attached to $M$, we set $U_0 :=
U(1)$ and $V_0 := V(1)$. Since $W$ is holomorphic, necessarily $V (
e^{ i\, \theta} ) = {\sf T}_1 U ( e^{ i\, \theta}) + V_0$. Applying
${\sf T}_1$ to both sides, we get ${\sf T}_1 V (e^{i\,\theta}) = -
U(e^{i\,\theta}) + U_0$ (the left and the right hand sides vanish at
$e^{ i\, \theta} = 1$). Applying ${\sf T}_1$ to $ V(e^{i\, \theta}) =
\varphi \left(X(e^{i\, \theta}),Y(e^{i\, \theta}),U(e^{i\,
\theta})\right)$ above and reorganizing, we obtain that $U$ satisfies
a functional equation\footnote{The origin of this equation may be
found in the seminal article~\cite{ bi1965} of Bishop. Since then, it
has been further exploited in \cite{ pi1974a, pi1974b, hita1978,
befo1978, we1982, bg1983, bpo1982, kw1982, r1983, hita1984, kw1984,
bpi1985, fr1985, trp1986, fo1986, tu1988, ai1989, tu1990, trp1990,
bo1991, dh1992, tu1994a, tu1994b, brt1994, cr1994, gl1994, me1994,
hukr1995, jo1995, trp1996, tu1996, jo1996, jo1997, me1997, po1997,
mp1998, tu1998, hu1998, cr1998, jo1999a, jo1999b, mp1999, ber1999,
po2000, mp2000, tu2001, da2001, ds2001, me2002, mp2002, hm2002,
po2003, js2004, me2004c}.} involving the Hilbert transform:
\def\theequation{3.5}\begin{equation}
\boxed{
U
(e^{i\,\theta})
= 
-
{\sf T}_1
\left[
\varphi(X(\cdot),Y(\cdot),U(\cdot))
\right](e^{i\,\theta})
+ 
U_0.}
\end{equation}
Here, the map $U : \partial \Delta \to \R^d$ is the unknown, whereas
the holomorphic map $Z = X+i\,Y : \partial \Delta \to \C^m$ and the
constant vector $U_0$ are given data.

Conversely, given $X + i\, Y$ and $U_0$, assume that $U \in \mathcal{
C }^{ \kappa, \alpha}$ satisfies the above functional equation. Set $V
( e^{ i\, \theta }) := {\sf T}_1 U ( e^{ i\, \theta}) + V_0$, where
$V_0 := \varphi (X (1), Y(1), U(1))$. Then $U (e^{ i\theta }) + i\, V(
e^{ i\, \theta})$ extends as a $\mathcal{ C}^{ \kappa, \alpha}$ map
$\overline{ \Delta} \ni \zeta \mapsto W ( \zeta) \in \C^d$ which is
holomorphic in $\Delta$. If $\vert [ X + i\, Y ] (e^{ i\, \theta})
\vert < \rho_1$, $\vert U ( e^{ i\, \theta}) \vert < \rho_1$ and
$\vert V ( e^{ i\, \theta}) \vert < \rho_1$, the disc $A := (Z, W)$ is
attached to $M$.

\smallskip

Bishop (1965) in the $\mathcal{ C}^{ \kappa, 0}$ classes and then
Hill-Taiani (1978), Boggess-Pitts
(1985) in the H\"older classes $\mathcal{ C}^{\kappa,
\alpha}$ established existence and uniqueness of the solution $U$ to
the fundamental functional equation~\thetag{ 3.5}.

\def\thetheorem{3.6}\begin{theorem}
{\rm (\cite{ bi1965, hita1978, bpi1985})} If $M$ is at
least $\mathcal{ C}^{ 1, \alpha}$, shrinking $\rho_1$ if necessary,
there exists $\rho_2$ with $0 <\rho_2 < \rho_1$ such that whenever the
data $Z \in \mathcal{ C}^{ 0, \alpha} ( \overline{ \Delta}, \C^m) \cap
\mathcal{ O} ( \Delta, \C^m)$ and $U_0 \in \R^d$ satisfy $\vert Z (e^{
i\, \theta}) \vert < \rho_2$ on $\partial \Delta$ and $\vert U_0 \vert
< \rho_2$, there exists a unique solution $U \in \mathcal{ C}^{ 0,
\beta} (\partial \Delta, \R^d )$, $0 < \beta < \alpha^2$, to the
Bishop-type functional equation~\thetag{ 3.5} above such that $\vert U
( e^{ i\theta} ) \vert < \rho_1$ on $\partial \Delta$ and such that in
addition $\vert V ( e^{ i\theta} ) \vert < \rho_1$ on $\partial
\Delta$, where
\[
V ( e^{ i\, \theta }) := {\sf
T}_1 U ( e^{ i\, \theta}) + \varphi (X (1), Y(1), U(1)).
\]
Consequently, the disc $(Z, U + i\, V)$ is attached to $M$.
\end{theorem}

Notice the (substantial) loss of smoothness, occuring also in~\cite{
brt1994, ber1999}, which is due to an application of a general
implicit function theorem in Banach spaces. The main theorem of this
chapter (\cite{ tu1990, tu1996}) refines the preceding result
with a negligible loss of smoothness, provided the graphing map
$\varphi$ belongs to the H\"older space $\mathcal{ C}^{\kappa,
\alpha}$. In the geometric applications (Parts~V and~VI), it is
advantageous to be able to solve a Bishop equation like~\thetag{ 3.5}
which involves supplementary parameters. Thus, instead of $\varphi$,
we shall consider an $\R^d$-valued $\mathcal{ C}^{ \kappa, \alpha}$
map $\Phi = \Phi ( u, e^{ i\, \theta}, s)$, where $s$ is a parameter.
For fixed $s$, we shall denote by $\Phi \vert_s$ the map 
$\square_{ \rho_1 }^d \times \partial \Delta \ni (u, e^{ i \, \theta })
\longmapsto \Phi \left( u, e^{ i\, \theta}, s \right) \in \R^d$. In
accordance with Section~1, we set:
\[
\left\vert\!\left\vert
\Phi_u
\right\vert\!\right\vert_{0,0}
:= 
\max_{1\leqslant j\leqslant d}\,
\left(
\sum_{1\leqslant l\leqslant d}\,
\left\vert\!\left\vert
\Phi_{u_l}^j
\right\vert\!\right\vert_{0,0}
\right),
\]
and similarly $\left\vert \! \left\vert \Phi_\theta \right\vert \!
\right\vert_{ 0,0 }
= \max_{1\leqslant j \leqslant d} \,
\big\vert \! \big\vert \Phi_\theta^j \big\vert \!
\big\vert_{ 0,0 }$.

\def\thetheorem{3.7}\begin{theorem}
{\rm (\cite{ tu1990, tu1996}, [$*$])} Let $\Phi = \Phi \left(
u, e^{ i\, \theta}, s \right)$ be an $\R^d$-valued map of class
$\mathcal{ C}^{ \kappa, \alpha}$, $\kappa \geqslant 1$, $\boxed{0 < \alpha
< 1}$, defined for $u\in \R^d$, $\vert u \vert < \rho_1$, $\theta \in
\R$ and $s\in \R^b$, $\vert s \vert < \sigma_1$, where $0 < \rho_1 <
1$ and $0 < \sigma_1 < 1$. Assume that on its domain of definition
$\square_{ \rho_1}^d \times \partial \Delta \times \square_{
\sigma_1}^b$, the map $\Phi$ and its derivatives with respect to $u$
and to $\theta$ satisfy the inequalities {\rm (}nothing is required
about $\Phi_s${\rm )}{\rm :}
\[
\left\vert\!\left\vert
\Phi
\right\vert\!\right\vert_{0,0}
\leqslant {\sf c}_1, 
\ \ \ \ \
\left\vert\!\left\vert
\Phi_u
\right\vert\!\right\vert_{0,0}
\leqslant {\sf c}_2, 
\ \ \ \ \
\left\vert\!\left\vert
\Phi_\theta
\right\vert\!\right\vert_{0,0}
\leqslant {\sf c}_3,
\]
for some small positive constants ${\sf c}_1$, 
${\sf c}_2$ and ${\sf c}_3$ such that
\def\theequation{3.8}\begin{equation}
{\sf c}_1 
\leqslant
C\,\alpha\,\rho_1,
\ \ \ \ \ \ \ \ \ 
{\sf c}_2
\leqslant
C^2
\alpha^2
\left[
1
+
\sup_{\vert s\vert<\sigma_1}
\left\vert\!\left\vert
\Phi\vert_s
\right\vert\!\right\vert_{1,\alpha}
\right]^{-2},
\ \ \ \ \ \
{\sf c}_3
\leqslant
\rho_1^2\,{\sf c}_2,
\end{equation}
where $0 < C < 1$ is an absolute constant. Then for every fixed
$U_0$ satisfying $\vert U_0 \vert < \rho_1 /16$ and every fixed $s \in
\square_{ \sigma_1 }^b$, the parameterized local Bishop-type
functional equation{\rm :}
\[
\boxed{
U(e^{i\, \theta})
= 
-
{\sf T}_1\left[\Phi\left(
U(\cdot),\cdot,s
\right)\right](e^{i\,\theta})
+
U_0
}
\]
has a unique solution{\rm :}
\[
\partial\Delta\ni 
e^{i\,\theta}
\longmapsto 
U(e^{i\,\theta},s,U_0)
\in\R^d,
\]
with $\left\vert \! \left\vert U \right\vert \! \right\vert_{ 0,
0}\leqslant \rho_1 /4$ which is of class $\mathcal{ C}^{
\kappa, \alpha}$ on $\partial \Delta$. Furthermore, this solution is
of class $\mathcal{ C}^{ \kappa, \alpha - 0} = \bigcap_{ \beta <
\alpha} \, \mathcal{ C}^{ \kappa, \beta}$ with respect to all the
variables, including parameters, namely the complete map
\[
\partial\Delta\times\square_{
\sigma_1}^b\times\square_{\rho_1/16}^d
\ni 
\left(e^{i\,\theta},s,U_0\right)
\longmapsto 
U(e^{i\,\theta},s,U_0)
\in\R^d
\]
is $\mathcal{ C}^{ \kappa, \alpha-0 }$.
\end{theorem}

Since the assumptions involve only the $\mathcal{ C}^{ 1, \alpha}$
norm of $\Phi$, we notice that the theorem is also true with $\Phi \in
\mathcal{ C}^{ \kappa, \alpha - 0}$, provided $\kappa \geqslant 2$, except
that the solution will only be $\mathcal{ C}^{ \kappa, \alpha - 0}$
with respect to $e^{ i\, \theta}$: it suffices to apply the theorem by
considering that $\Phi \in \mathcal{ C}^{ \kappa, \beta}$, with $\beta
< \alpha$ arbitrary, getting a solution that is $\mathcal{ C}^{
\kappa, \beta -0}$ with respect to all variables and concluding from
$\bigcap_{ \beta < \alpha} \, \mathcal{ C}^{ \kappa, \beta - 0} =
\mathcal{ C}^{ \kappa, \alpha - 0}$.
 
The main purpose of this section is to provide a thorough proof of the
theorem. In the sequel, $C$, $C_1$, $C_2$, $C_3$ and $C_4$ will denote
positive absolute constants. We may assure that they all will be
$\geqslant 10^{ -5}$ and $\leqslant 10^5$.

\smallskip

The smallness of $\left\vert\! \left\vert \Phi \right\vert\!
\right\vert_{0,0}$, of $\left\vert\! \left\vert \Phi_u \right\vert\!
\right\vert_{0,0}$ and of $\left\vert\! \left\vert \Phi_\theta
\right\vert\! \right\vert_{0,0}$ guarantee the smallness of
$\left\vert\! \left\vert \Phi \vert s \right\vert\!
\right\vert_{1, \alpha / 2}$ by virtue of an elementary observation.

\def\thelemma{3.9}\begin{lemma}
{\rm ([$*$])}
Let ${\sf x} \in \square_\rho^n$, $n\in \N$, $n\geqslant 1$, 
where $0 <\rho_i \leqslant
\infty$, and let $f =
f ({\sf x})$ be $\mathcal{ C}^{ 0, \alpha}$ function with values in
$\R^d$, $d\geqslant 1$. If $\left\vert \! \left\vert f
\right\vert \! \right\vert_{0,0} \leqslant {\sf c}$, for some quantity
${\sf c} >0$, then{\rm :}
\[
\left\vert\!\left\vert
f
\right\vert\!\right\vert_{\widehat{0,\alpha/2}}
\leqslant
{\sf c}^{1/2}
\left[
2
+
\left\vert\!\left\vert
f
\right\vert\!\right\vert_{\widehat{0,\alpha}}
\right].
\]
\end{lemma}

\noindent
We apply this inequality to $\Phi_u \vert_s$ and to $\Phi_\theta
\vert_s$, pointing out that for any $\beta$ with $0 < \beta \leqslant
\alpha$, by definition:
\begin{small}
\[
\aligned
\left\vert\!\left\vert
\Phi_u\vert_s
\right\vert\!\right\vert_{\widehat{0,\beta}}
&
=
\max_{1\leqslant j\leqslant d}
\left(
\sum_{l=1}^d\,
\frac{
\big\vert
\Phi_{u_l}^j(u'',e^{i\,\theta''},s)
-
\Phi_{u_l}^j(u',e^{i\,\theta'},s)
\big\vert
}{
\left\vert
(u'',\theta'')
-
(u',\theta')
\right\vert^\beta
}
\right),
\\
\left\vert\!\left\vert
\Phi_\theta\vert_s 
\right\vert\!\right\vert_{\widehat{0,\beta}}
&
=
\max_{1\leqslant j\leqslant d}
\frac{
\big\vert
\Phi_\theta^j(u'',e^{i\,\theta''},s)
-
\Phi_\theta^j(u',e^{i\,\theta'},s)
\big\vert
}{
\left\vert
(u'',\theta'')
-
(u',\theta')
\right\vert^\beta
}.
\endaligned
\]
\end{small}

\def\thelemma{3.10}\begin{lemma}
{\rm ([$*$])}
Independently of $s$, we have{\rm :}
\[
\aligned
\left\vert\!\left\vert
\Phi_u\vert_s
\right\vert\!\right\vert_{\widehat{0,\alpha/2}}
&
\leqslant
{\sf c}_2^{1/2}
\left[
2
+
\left\vert\!\left\vert
\Phi\vert_s
\right\vert\!\right\vert_{
1,\alpha}
\right],
\\
\left\vert\!\left\vert
\Phi_\theta\vert_s
\right\vert\!\right\vert_{\widehat{0,\alpha/2}}
&
\leqslant
{\sf c}_3^{1/2}
\left[
2
+
\left\vert\!\left\vert
\Phi\vert_s
\right\vert\!\right\vert_{
1,\alpha}
\right],
\\
\left\vert\!\left\vert
\Phi\vert_s
\right\vert\!\right\vert_{1,\alpha/2}
&
\leqslant
{\sf c}_1
+
{\sf c}_2
+
{\sf c}_3
+
\left(
{\sf c}_2^{1/2}
+
{\sf c}_3^{1/2}
\right)
\left[
2
+
\left\vert\!\left\vert
\Phi\vert_s
\right\vert\!\right\vert_{
1,\alpha}
\right].
\endaligned
\]
\end{lemma}

\noindent 
The presence of the squares in
the inequalities of Theorem~3.7 anticipates the roots
${\sf c}_2^{1/2}$ and ${\sf c}_3^{ 1/2}$ above. These two lemmas and
the next involve dry computations with H\"older norms. The
detailed proofs are postponed to Section~4.

\def\thelemma{3.11}\begin{lemma}
{\rm ([$*$])} If $U \in \mathcal{ C}^{ 1, \beta} ( \partial \Delta,
\R^d)$ with $0 < \beta \leqslant \alpha$ satisfies $\vert U (e^{ i\,
\theta}) \vert < \rho_1$ on $\partial \Delta$, then for every fixed $s
\in \square_{ \sigma_1}^b$, we have{\rm :}
\[
\aligned
\left\vert\!\left\vert
\Phi(U(\cdot),\cdot,s)
\right\vert\!\right\vert_{
\mathcal{ C}^{1,\beta}(\partial \Delta)}
&
\leqslant
\left\vert\!\left\vert
\Phi
\right\vert\!\right\vert_{0,0}
+
\left\vert\!\left\vert
\Phi_\theta
\right\vert\!\right\vert_{0,0}
+
\left\vert\!\left\vert
\Phi_\theta\vert_s
\right\vert\!\right\vert_{\widehat{0,\beta}}
\left[
1
+
\left(
\left\vert\!\left\vert
U
\right\vert\!\right\vert_{
\widehat{1,0}}
\right)^\beta
\right]
+
\\
&
+
\left\vert\!\left\vert
\Phi_u
\right\vert\!\right\vert_{0,0}
\left\vert\!\left\vert
U
\right\vert\!\right\vert_{1,\beta}
+
\left\vert\!\left\vert
\Phi_u\vert_s
\right\vert\!\right\vert_{\widehat{0,\beta}}
\left[
\left\vert\!\left\vert
U
\right\vert\!\right\vert_{
\widehat{1,0}}
+
\left(
\left\vert\!\left\vert
U
\right\vert\!\right\vert_{
\widehat{1,0}}
\right)^{1+\beta}
\right].
\endaligned
\]
\end{lemma}

Remind $\left\vert\!\left\vert U \right\vert\!\right\vert_{ \widehat{
1,0}} = \sup_{\theta}\, \left\vert U_\theta (e^{ i\, \theta})
\right\vert$. We then introduce the map:
\[
U\longmapsto 
\mathfrak{F}(U)
:=
U_0
-
{\sf T}_1\left[\Phi(U(\cdot),\cdot,s)\right](e^{i\,\theta}).
\]
To construct the solution $U$, we endeavour a Picard
iteration process, setting $U_k\vert_{ k=0} := U_0$ with $\vert U_0
\vert < \rho_1 / 16$ and $U_{ k+1} := \mathfrak{ F} ( U_k)$, for $k\in
\N$, whenever $\mathfrak{ F} ( U_k)$ may be defined, {\it i.e.}
whenever $\left\vert \! \left\vert U_k \right\vert \! \right\vert_{ 0,
0} < \rho_1$. We shall first work in $\mathcal{ C}^{ 1, \alpha /2}
\subset \mathcal{ C}^{ \kappa, \alpha}$.

\def\thelemma{3.12}\begin{lemma}
If we choose the absolute constant $C < 1$ appearing in the theorem
sufficiently small, then independently of $s$, the sequence $U_k$
satisfies the uniform boundedness estimate{\rm :}
\[
\left\vert\!\left\vert 
U_k 
\right\vert\!\right\vert_{ 
1, \alpha /2} 
\leqslant 
\rho_1/4
< 
\rho_1,
\]
hence it is defined for every $k\in \N$ and 
each $U_k$ belongs to $\mathcal{ C}^{
1, \alpha/2} ( \partial \Delta)$.
\end{lemma}

\proof
By Theorem~2.24, there exists an absolute constant $C_1 >1$ (not
exactly the same), such that
\[
\left\vert\!\left\vert\!\left\vert 
{\sf T}_1 
\right\vert\!\right\vert\!\right\vert_{1,\alpha/2}
\leqslant 
C_1/\alpha.
\]
Majorating by means of the $\mathcal{ C }^{0, \alpha / 2}$-norm:
\[
\left\vert\!\left\vert
\mathfrak{F}(U_k)
\right\vert\!\right\vert_{1,\alpha/2}
\leqslant
\vert U_0\vert
+
\left\vert\!\left\vert\!\left\vert
{\sf T}_1
\right\vert\!\right\vert\!\right\vert_{1,\alpha/2}
\,
\left\vert\!\left\vert
\Phi(U_k(\cdot),\cdot,s)
\right\vert\!\right\vert_{\mathcal{C}^{1,\alpha/2}(\partial\Delta)}.
\]
Assume by induction that $U_k$ is $\mathcal{ C}^{ 1, \alpha/2}$ and
satisfies $\left\vert \! \left\vert U_k \right\vert \! \right\vert_{
1, \alpha /2} \leqslant \rho_1 /4$ (this holds for $k=0$). Clearly $U_{k+1}
= \mathfrak{ F} (U_k)$ is $\mathcal{ C}^{ 1, \alpha/2}$. Thanks to
Lemma~3.11, and to the trivial majoration $( \left\vert \! \left\vert
U_k \right\vert \! \right\vert_{ \widehat{ 1, 0}} )^{ \alpha/2} \leqslant
(\rho_1/4 )^{\alpha /2} < 1$:
\[
\aligned
\left\vert\!\left\vert
\Phi(U_k(\cdot),\cdot,s)
\right\vert\!\right\vert_{\mathcal{C}^{1,\alpha/2}(\partial\Delta)}
&
\leqslant
\left\vert\!\left\vert
\Phi
\right\vert\!\right\vert_{0,0}
+
\left\vert\!\left\vert
\Phi_\theta
\right\vert\!\right\vert_{0,0}
+
2
\left\vert\!\left\vert
\Phi_\theta\vert_s
\right\vert\!\right\vert_{\widehat{0,\alpha/2}}
+
\\
& \ \ \ 
+
\left\vert\!\left\vert
\Phi_u
\right\vert\!\right\vert_{0,0}
\left\vert\!\left\vert
U_k
\right\vert\!\right\vert_{1,\alpha/2}
+
2
\left\vert\!\left\vert
\Phi_u\vert_s
\right\vert\!\right\vert_{\widehat{0,\alpha/2}}
\left\vert\!\left\vert
U_k
\right\vert\!\right\vert_{
1,\alpha/2}.
\endaligned
\]
Using then the assumptions~\thetag{ 3.8} of the theorem together with 
Lemma~3.10:
\[
\aligned
\left\vert\!\left\vert
U_{k+1}
\right\vert\!\right\vert_{1,\alpha/2}
&
\leqslant
\rho_1/16
+
C_1\,\alpha^{-1}
\left[
{\sf c}_1
+
{\sf c}_3
+
4\,{\sf c}_3^{1/2}
\left(
1
+
\left\vert\!\left\vert
\Phi\vert_s
\right\vert\!\right\vert_{1,\alpha}
\right)
+
\right.
\\
&
\left.
+
\left\vert\!\left\vert
U_k
\right\vert\!\right\vert_{1,\alpha/2}
\left(
{\sf c}_2
+
4\,{\sf c}_2^{1/2}
\left(
1
+
\left\vert\!\left\vert
\Phi\vert_s
\right\vert\!\right\vert_{1,\alpha}
\right)
\right)
\right].
\endaligned
\]
Using the two trivial majorations 
${\sf c}_3 \leqslant C \alpha \rho_1$ and
${\sf c}_2 \leqslant C \alpha$ together with the main 
assumptions~\thetag{ 3.8}
to majorate ${\sf c}_2^{ 1/2}$ and ${\sf c }_3^{ 1/2}$, we get:
\[
\left\vert\!\left\vert
U_{k+1}
\right\vert\!\right\vert_{1,\alpha/2}
\leqslant
\rho_1/16
+
C_1\,\rho_1\,6\,C
+
\left\vert\!\left\vert
U_k
\right\vert\!\right\vert_{1,\alpha/2}
\,C_1\,5\,C.
\]
Choosing $C \leqslant \frac{ 1}{ 16\, C_1\, 6}$
(whence $C \leqslant \frac{ 1}{ 2\, C_1\, 5}$), we finally
get:
\[
\left\vert\!\left\vert
U_{k+1}
\right\vert\!\right\vert_{1,\alpha/2}
\leqslant
\rho_1/8
+
(1/2)\,
\left\vert\!\left\vert
U_k
\right\vert\!\right\vert_{1,\alpha/2}.
\]
By immediate induction, the assumption $\vert U_0 \vert < \rho_1 /16$
and these (strict) inequalities entail that $\left\vert \! \left\vert
U_k \right\vert \! \right\vert_{1, \alpha /2} \leqslant 
\rho_1 /4$ for every $k\in \N$, as claimed.
\endproof

\def\thelemma{3.13}\begin{lemma}
{\rm (\cite{ tu1990}, [$*$])} For every $\beta$ with $0 < \beta \leqslant
\alpha$ and every fixed $s \in \square_{ \rho_1 }^b$, if two maps $U^j
\in \mathcal{ C}^{ 1, 0} (\partial \Delta, \R^d)$ with $\left\vert \!
\left\vert U^j \right\vert \! \right\vert_{ 0, 0} 
< \rho_1/3$ for $j=1,
2$ are given, the following inequality holds{\rm :}
\[
\left\vert\!\left\vert
\Phi(U^2(\cdot),\cdot,s)
-
\Phi(U^1(\cdot),\cdot,s)
\right\vert\!\right\vert_{\mathcal{C}^{0,\beta}(\partial\Delta)}
\leqslant
{\sf C}\,
\left\vert\!\left\vert
U^2
- 
U^1
\right\vert\!\right\vert_{\mathcal{C}^{0,\beta}(\partial\Delta)},
\]
where
\[
{\sf C}
=
\left\vert\!\left\vert
\Phi_{\vert s}
\right\vert\!\right\vert_{1,\beta}\,2
\left[
1
+
\left(
\left\vert\!\left\vert
U^1
\right\vert\!\right\vert_{\widehat{1,0}}
\right)^\beta
+
\left(
\left\vert\!\left\vert
U^2
\right\vert\!\right\vert_{\widehat{1,0}}
\right)^\beta
\right].
\]
\end{lemma}

Again, the (latexnically lengthy) 
proof is postponed to Section~4. 

\def\thelemma{3.14}\begin{lemma} 
If we choose the absolute constant $C$ of
the theorem sufficiently small, then
independently of $s$, the map{\rm :}
\[
U\longmapsto 
\mathfrak{F}(U)
:=
U_0
-
{\sf T}_1\left[\Phi(U(\cdot),\cdot,s)\right]
(e^{i\,\theta}),
\]
restricted to the set of those $U \in \mathcal{ C}^{ 1, \alpha/2}
(\partial \Delta, \R^d)$ that satisfy $\left \vert \! \left\vert U
\right\vert \! \right\vert_{ 1, \alpha /2} \leqslant
\rho_1 /4$, is a contraction{\rm :}
\[
\left\vert\!\left\vert
\mathfrak{F}(U^2)
-
\mathfrak{F}(U^1)
\right\vert\!\right\vert_{0,\alpha/2}
\leqslant
\frac{1}{2}\,
\left\vert\!\left\vert
U^2
-
U^1
\right\vert\!\right\vert_{0,\alpha/2}.
\]
\end{lemma}

\proof
Let $U^j \in \mathcal{ C}^{ 1, \alpha/2}$ with $\left\vert \!
\left\vert U^j \right\vert \! \right\vert_{ 1, \alpha/2} \leqslant
\rho_1/4$ for $j=1, 2$. In particular, $\left\vert \! \left\vert U^j
\right\vert \! \right\vert_{ 0, 0} < \rho_1 /3$, so Lemma~3.13
applies. In the majorations below, to pass to the fourth line, we use
the assumption $\rho_1 < 1$, which enables us to majorate simply by
$3$ the three terms in the brackets of the third line:
\[
\aligned
&
\left\vert\!\left\vert
\mathfrak{F}(U^2)
-
\mathfrak{F}(U^1)
\right\vert\!\right\vert_{0,\alpha/2}
=
\left\vert\!\left\vert
{\sf T}_1
\left[
\Phi\left(U^2(\cdot),\cdot,s\right)
-
\Phi\left(U^1(\cdot),\cdot,s\right)
\right]
\right\vert\!\right\vert_{0,\alpha/2}
\\
&
\leqslant
\left\vert\!\left\vert\!\left\vert
{\sf T}_1
\right\vert\!\right\vert\!\right\vert_{0,\alpha/2}\,
\left\vert\!\left\vert
\Phi\left(U^2(\cdot),\cdot,s\right)
-
\Phi\left(U^1(\cdot),\cdot,s\right)
\right\vert\!\right\vert_{0,\alpha/2}
\\
&
\leqslant
\frac{C_1}{\alpha}\,
\left\vert\!\left\vert
\Phi_{\vert s}
\right\vert\!\right\vert_{1,\alpha/2}\,2
\left[
1
+
\left(
\left\vert\!\left\vert
U^1
\right\vert\!\right\vert_{1,\alpha/2}
\right)^{\alpha/2}
+
\left(
\left\vert\!\left\vert
U^2
\right\vert\!\right\vert_{1,\alpha/2}
\right)^{\alpha/2}
\right]
\left\vert\!\left\vert
U^2
-
U^1
\right\vert\!\right\vert_{0,\alpha/2}
\\
&
\leqslant
\frac{C_2}{\alpha}\,
\left\vert\!\left\vert
\Phi_{\vert s}
\right\vert\!\right\vert_{1,\alpha/2}
\left\vert\!\left\vert
U^2
-
U^1
\right\vert\!\right\vert_{0,\alpha/2},
\endaligned
\]
where $C_2 > 1$ is absolute. Then we apply Lemma~3.10 to majorate
$\left\vert \! \left\vert \Phi\vert_s \right\vert \! \right\vert_{ 1,
\alpha /2}$, we use the three trivial majorations ${\sf c}_1, {\sf c}_2,
{\sf c}_3 \leqslant C \alpha$ and we majorate ${\sf c}_2^{ 1/2}, {\sf
c}_3^{ 1/2}$ by means of~\thetag{ 3.8}, dropping $\rho_1 < 1$ in ${\sf
c}_3^{ 1/2}$, which yields:
\[
\aligned
\left\vert\!\left\vert
\Phi\vert_s
\right\vert\!\right\vert_{1,\alpha/2}
&
\leqslant
{\sf c}_1
+
{\sf c}_2
+
{\sf c}_3
+
\left(
{\sf c}_2^{1/2}
+
{\sf c}_3^{1/2}
\right)
\left[
2
+
\left\vert\!\left\vert
\Phi_{\vert s}
\right\vert\!\right\vert_{
1,\alpha}
\right]
\\
&
\leqslant
3\,C\alpha
+
2\,C\alpha
+
4\,C\alpha
=
9\,C\alpha.
\endaligned
\]
Then we conclude that
\[
\left\vert\!\left\vert
\mathfrak{F}(U^2)
-
\mathfrak{F}(U^1)
\right\vert\!\right\vert_{0,\alpha/2}
\leqslant
C\,C_3
\left\vert\!\left\vert
U^2
-
U^1
\right\vert\!\right\vert_{0,\alpha/2}.
\]
Choosing the absolute constant $C$ of
the theorem $\leqslant \frac{ 1}{ 2C_3}$ yields the desired contracting factor
$\frac{ 1}{ 2}$.
\endproof

The fixed point theorem then entails that our sequence $U_k$ converges
in $\mathcal{ C}^{ 0, \alpha /2}$-norm towards some map $U \in
\mathcal{ C}^{ 0, \alpha /2} ( \partial \Delta, \R^d)$. More is true:

\def\thelemma{3.15}\begin{lemma}
For every fixed parameters $(s, U^0)$, this solution $U = U \left(
e^{i\, \theta}, s, U_0 \right) = \lim_{ k\to \infty} \, U_k$ belongs
in fact to $\mathcal{ C}^{ 1, \alpha / 2} (\partial \Delta, \R^d)$ and
satisfies $\left\vert \! \left\vert U \right\vert \! \right\vert_{ 1,
\alpha/2} \leqslant \rho_1 /4$.
\end{lemma}

\proof
Indeed, since $\left\vert \! \left\vert U_k \right\vert \!
\right\vert_{ 1, \alpha/2} \leqslant \rho_1 /4$ is bounded, it is possible
thanks to the Arzel\`a-Ascoli lemma to extract some subsequence
converging in $\mathcal{ C}^{ 1, 0 } (\partial \Delta)$ to a map,
still denoted by $U = U \left( e^{i\, \theta}, s, U_0 \right)$, which
is $\mathcal{ C}^{ 1, 0 }$ on $\partial \Delta$. Still denoting by
$U_k$ such a subsequence, we observe that the uniform convergence
$\left\vert \! \left\vert U_k - U \right\vert \! \right\vert_{ 1, 0}
\to 0$ plus the boundedness $\left\vert \! \left\vert U_k \right\vert
\! \right\vert_{ \widehat{ 1, \alpha/2}} \leqslant \rho_1 /4$
entail immediately that the following majoration holds:
\begin{small}
\[
\frac{
\big\vert
U_{\theta}(e^{i\,\theta''})
-
U_{\theta}(e^{i\,\theta'})
\big\vert}{
\left\vert
\theta''
-
\theta'
\right\vert^{\alpha/2}}
=
\lim_{k\to\infty}\
\frac{
\big\vert
U_{k,\theta}(e^{i\,\theta''})
-
U_{k,\theta}(e^{i\,\theta'})
\big\vert}{
\left\vert
\theta''
-
\theta'
\right\vert^{\alpha/2}}
\leqslant
\frac{\rho_1}{4},
\]
\end{small}
for arbitrary $0 < \left\vert \theta'' - \theta ' \right\vert \leqslant
\pi$. Consequently, $U$ belongs to $\mathcal{ C}^{ 1, \alpha /2}$.
Passing to the limit in $\left\vert \! \left\vert U_k \right\vert \!
\right\vert_{ 1, \alpha/2} \leqslant \rho_1 /4$, we also 
deduce $\left\vert \! \left\vert U \right\vert \!
\right\vert_{ 1, \alpha/2} \leqslant \rho_1 /4$.
\endproof

The next crucial step is to study the regularity of the solution $U =
U \left( e^{i\, \theta}, s, U_0 \right)$ with respect to $(s, U_0)$.

\def\thelemma{3.16}\begin{lemma}
The solution $U = U (e^{ i\, \theta}, s, U_0)$ satisfies a Lipschitz
condition with respect to the parameters $s$ and $U_0$.
\end{lemma}

\proof
Consider two parameters $s^1, s^2 \in \square_{ \sigma_1}^b$ and
define $U^j := U( e^{i\, \theta}, s^j, U_0)$ for $j=1, 2$. 
Then substract the two corresponding Bishop equations, 
insert two innocuous opposite terms and majorate:
\[
\aligned
\left\vert\!\left\vert
U^2
-
U^1
\right\vert\!\right\vert_{0,\alpha/2}
\leqslant
\left\vert\!\left\vert\!\left\vert
T_1
\right\vert\!\right\vert\!\right\vert_{0,\alpha/2}
&
\left[
\left\vert\!\left\vert
\Phi\left(U^2(\cdot),\cdot,s^2\right)
-
\Phi\left(U^2(\cdot),\cdot,s^1\right)
\right\vert\!\right\vert_{0,\alpha/2}
+
\right.
\\
&
\
\left.
+
\left\vert\!\left\vert
\Phi\left(U^2(\cdot),\cdot,s^1\right)
-
\Phi\left(U^1(\cdot),\cdot,s^1\right)
\right\vert\!\right\vert_{0,\alpha/2}
\right].
\endaligned
\]
To majorate the difference in the second line, we again apply
Lemma~3.13. To majorate the difference in the first line, we
apply:

\def\thelemma{3.17}\begin{lemma}
{\rm ([$*$])} Let $\beta$ with $0 < \beta \leqslant \alpha$, let $U \in
\mathcal{ C}^{1, 0} ( \partial \Delta, \R^d)$ with $\left\vert \!
\left\vert U \right\vert \! \right\vert_{ 0, 0} < \rho_1$ and
let two parameters $s^1, s^2 \in \square_{ \sigma_1}^b$. Then
\begin{small}
\[
\left\vert\!\left\vert
\Phi\left(U(\cdot),\cdot,s^2\right)
-
\Phi\left(U(\cdot),\cdot,s^1\right)
\right\vert\!\right\vert_{0,\beta}
\leqslant
\left\vert
s^2
-
s^1
\right\vert
\left(
\left\vert\!\left\vert
\Phi
\right\vert\!\right\vert_{1,0}
+
\left\vert\!\left\vert
\Phi
\right\vert\!\right\vert_{1,\beta}
\left[
1
+
\left(
\left\vert\!\left\vert
U
\right\vert\!\right\vert_{\widehat{1,0}}
\right)^\beta
\right]
\right).
\]
\end{small}
\end{lemma}

\noindent
With $\beta := \frac{\alpha}{2}$, we thus obtain:
\begin{small}
\[
\aligned
\left\vert\!\left\vert
U^2
-
U^1
\right\vert\!\right\vert_{0,\alpha/2}
&
\leqslant
\frac{C_1}{\alpha}
\left[
\left\vert
s^2
-
s^1
\right\vert
\left(
\left\vert\!\left\vert
\Phi
\right\vert\!\right\vert_{1,0}
+
\left\vert\!\left\vert
\Phi
\right\vert\!\right\vert_{1,\alpha/2}
\left[
1
+
\left(
\left\vert\!\left\vert
U^2
\right\vert\!\right\vert_{
1,0}
\right)^{\frac{\alpha}{2}}
\right]
\right)
+
\right.
\\
&
\left.
+
\sup_{\vert s\vert<\sigma_1}\,
\left\vert\!\left\vert
\Phi_{\vert s}
\right\vert\!\right\vert_{1,\alpha/2}\,2
\left[
1
+
\left(
\left\vert\!\left\vert
U^1
\right\vert\!\right\vert_{\widehat{1,0}}
\right)^{\frac{\alpha}{2}}
+
\left(
\left\vert\!\left\vert
U^2
\right\vert\!\right\vert_{\widehat{1,0}}
\right)^{\frac{\alpha}{2}}
\right]
\right]
\left\vert\!\left\vert
U^2
-
U^1
\right\vert\!\right\vert_{0,\alpha/2}.
\endaligned
\]
\end{small}

\noindent
Then we apply the majoration
of Lemma~3.10 to $\left\vert\!\left\vert
\Phi_{ \vert s} \right\vert \! \right\vert_{ 1, \alpha /2}$ and we use
$\rho_1 < 1$ to majorate by $1$ the terms $\left\vert \! \left\vert
U^j \right\vert \! \right\vert_{ \widehat{1, 0}} \leqslant \rho_1 /4$:
\[
\left\vert\!\left\vert
U^2
-
U^1
\right\vert\!\right\vert_{0,\alpha/2}
\leqslant
C_1\alpha^{-1}
\left\vert
s^2
-
s^1
\right\vert
3\,
\left\vert\!\left\vert
\Phi
\right\vert\!\right\vert_{1,\alpha/2}
+
C\,C_2
\left\vert\!\left\vert
U^2
-
U^1
\right\vert\!\right\vert_{0,\alpha/2}.
\]
Setting ${\sf K} := C_1\alpha^{-1} \, 3 \, \left\vert \! \left\vert
\Phi \right\vert \! \right\vert_{1, \alpha /2}$, requiring $C
\leqslant \frac{ 1}{ 2C_2}$ and reorganizing we obtain that $U (
e^{i\, \theta}, s, U_0)$ is Lipschitzian with respect to $s$:
\[
\left\vert\!\left\vert
U^2
- 
U^1 
\right\vert\!\right\vert_{0,0}
\leqslant
\left\vert\!\left\vert
U^2
- 
U^1 
\right\vert\!\right\vert_{0,\alpha/2} 
\leqslant
2\,{\sf K}
\left\vert 
s^2 
- 
s^1 
\right\vert.
\]
The proof that $U ( e^{i\, \theta}, s, U_0)$ is Lipschitzian
with respect to $U_0$ is similar and even simpler.
\endproof

In summary, the solution $U = U ( e^{i\, \theta}, s, U_0)$ is
$\mathcal{ C}^{ 1, \alpha /2}$ with respect to $e^{ i\, \theta}$ and
Lipschitzian with respect to {\it all}\, the variables $\left( e^{ i\,
\theta}, s, U_0 \right)$.

Consequently, according to a theorem due to Rademacher (\cite{ ra1919,
fe1969}), the partial derivatives $U_{ s_k}$, $k=1, \dots, b$ and $U_{
U_0^m}$, $m=1, \dots, d$ exist in $L^\infty$. We then have to
differentiate the Bishop-type equation of Theorem~3.7 with respect to
$\theta$, to $s_k$ and to $U_0^m$. However, the linear operator ${\sf
T}_1$ is {\it not}\, bounded in $L^\infty$; in fact, according
to~\thetag{ 2.27}, $\left\vert \! \left\vert \! \left\vert {\sf T}
\right\vert \! \right\vert \! \right\vert_{ L^{\sf p}} \sim {\sf p}$
as ${\sf p} \to \infty$. So we need more information.

\def\thelemma{3.18}\begin{lemma}
There exists a null-measure subset $\mathfrak{ N} \subset \square_{
\sigma_1}^b \times \square_{ \rho_1/16}^d$ and there exists a quantity
${\sf K} >0$ such that at every $(s, U_0) \not \in \mathfrak{ N}$, for
every $k= 1, \dots, b$ and for every $m = 1, \dots, d${\rm :}

\begin{itemize}

\smallskip\item[{\bf (i)}]
the partial derivatives $U_{ s_k} ( e^{i\, \theta}, s, U_0)$ 
and $U_{ U_0^m} ( e^{i\, \theta}, s, U_0 )$ exist {\rm for
every} $e^{ i\, \theta} \in \partial \Delta${\rm ;}

\smallskip\item[{\bf (ii)}]
the maps $e^{i\, \theta} \mapsto U_{ s_k }(e^{i\, \theta}, s,U_0)$ and
$e^{i\, \theta} \mapsto U_{ U_0^m }( e^{i\, \theta}, s,U_0)$ are
$\mathcal{ C}^{ 0, \alpha /2}$ on $\partial \Delta$ and satisfy the
uniform inequality
\[
\left\vert\!\left\vert 
U_{s_k}
(\cdot,s,U_0)
\right\vert\!\right\vert_{
\mathcal{C}^{0,\alpha/2}(\partial\Delta)}
\leqslant
{\sf K}
\ \ \ 
{\it and}
\ \ \
\left\vert\!\left\vert 
U_{U_0^m}
(\cdot,s,U_0)
\right\vert\!\right\vert_{
\mathcal{C}^{0,\alpha/2}(\partial\Delta)}
\leqslant
{\sf K}.
\]

\end{itemize}\smallskip 
\end{lemma}

\proof
Since $U$ is almost everywhere differentiable with respect to all its
arguments, there exist a subset $\mathfrak{ F} \subset \square_{
\sigma_1}^b \times \square_{ \rho_1 / 16}^d \times \partial \Delta$
having {\sl full measure}, namely its complement has null measure,
such that for every $(e^{ i\, \theta}, s, U_0) \in \mathfrak{ F}$, all
partial derivatives $U_\theta, \, U_{ s_k }, U_{ U_0^m}$ exist at
$(e^{ i\, \theta}, s, U_0)$. Since $\mathfrak{ F}$ has full measure,
there exists a null measure subset $\mathfrak{ N} \subset \square_{
\sigma_1}^b \times \square_{ \rho_1 /16 }^d$ such that for every $(s,
U_0) \not\in \mathfrak{ N}$, the slice
\[
\mathfrak{F}_{s,U_0}
:=
\big(
\partial\Delta\times\{s\}\times\{U_0\}
\big)
\cap
\mathfrak{F}
\]
is a subset of $\partial \Delta$ having full measure, so that
$U_\theta, \, U_{ s_k}, U_{ U_0^m}$ exist at $(e^{ i\, \theta}, s,
U_0)$ with $e^{ i\, \theta} \in \mathfrak{ F}_{ s, U_0}$.

Fix $(s, U_0) \not\in \mathfrak{ N}$. We will treat only the partial
derivatives with respect to the $s_k$, arguments being similar for the
$U_{ U_0^m}$. In the end of the proof of Lemma~3.17, we have shown:
\[
\left\vert\!\left\vert 
U^2-U^1 
\right\vert\!\right\vert_{0,\alpha/2} 
\leqslant
{\sf K} 
\left\vert 
s^2-s^1 
\right\vert, 
\]
for some (not the same) quantity ${\sf K} >0$. Fix $k\in \{ 1, 2,
\dots, b\}$, take $s^2$ and $s^1$ with $s_k^2 \neq s_k^1$ but $s_l^2 =
s_l^1$ for $l\neq k$. The inequality above says that for every $e^{
i\, \theta}, \, e^{ i\, \theta'}, \, e^{ i\, \theta''} \in \partial
\Delta$ with $0 < \vert \theta '' - \theta ' \vert \leqslant \pi$,
we have
\[
\aligned
&
\left\vert
\frac{
U(e^{i\,\theta},s^2,U_0)
-
U(e^{i\,\theta},s^1,U_0)
}{
s_k^2-s_k^1
}
\right\vert
+
\\
&
\ \ \ \
+
\left\vert
\frac{
U(e^{i\,\theta''},s^2,U_0)
-
U(e^{i\,\theta''},s^1,U_0)
}{
s_k^2-s_k^1
}
-
\right.
\\
&
\ \ \ \ \ \ \ \ \ \ \
\left.
-
\frac{
U(e^{i\,\theta'},s^2,U_0)
-
U(e^{i\,\theta'},s^1,U_0)
}{
s_k^2-s_k^1
}
\right\vert
\Big/
\vert\theta''-\theta'\vert^{\alpha/2}
\leqslant
{\sf K}.
\endaligned
\]
Assume $e^{ i\, \theta}, \, e^{ i\, \theta'}, \, e^{ i\, \theta''} \in
\mathfrak{ F}_{ s^1, U_0}$, let $s_k^2 \to s_k^1$ (the limits of the
quotients above exist) and rename $s^1$ by $s$:
\[
\big\vert
U_{s_k}(e^{i\,\theta},s,U_0)
\big\vert
+
\frac{
\big\vert
U_{s_k}(e^{i\,\theta''},s,U_0)
-
U_{s_k}(e^{i\,\theta'},s,U_0)
\big\vert}{
\vert
\theta''-\theta'
\vert^{\alpha/2}}
\leqslant
{\sf K}.
\]
This inequality says that $U_{ s_k} (\cdot, s, U_0)$ is $\mathcal{
C}^{ 0, \alpha /2}$ almost everywhere on $\partial \Delta$. The next
extension lemma concludes the proof.
\endproof

\def\thelemma{3.19}\begin{lemma}
Let $n\geqslant 1$, let ${\sf x} \in \R^n$, let $m\geqslant 1$, let
${\sf y} \in \R^m$, let $\rho >0$, let $\sigma >0$, and let $f = f(
{\sf x}, {\sf y})$ be a function defined {\rm (}only{\rm )} in a
full-measure subset $\mathfrak{ F} \subset \square_\rho^n \times
\square_\sigma^m$, so that there exists a null-measure subset
$\mathfrak{ N} \subset \square_\sigma^m$ with the property that for
every ${\sf y} \not \in \mathfrak{ N}$, the slice $\mathfrak{ F}_{\sf
y} := \big( \square_\rho^n \times \{ {\sf y} \} \big) \times
\mathfrak{ F}$ has full measure in $\square_\rho^n$.
Assume that for every ${\sf y} \not\in
\mathfrak{ N}$, every 
${\sf x}, {\sf x}', {\sf x}'' \in 
\mathfrak{ F}_{ \sf y}$, we have
\[
\vert
f({\sf x},{\sf y})
\vert
+
\frac{
\big\vert
f({\sf x}'',{\sf y})
-
f({\sf x}',{\sf y})
\big\vert}{
\vert
{\sf x}''-{\sf x}'
\vert^{\beta}}
\leqslant
{\sf K},
\]
for some $\beta$ with $0 < \beta \leq 1$ and some quantity ${\sf K} >
0$. Then for every ${\sf y} \not\in \mathfrak{ N}$, the function
${\sf x} \mapsto f ( {\sf x}, {\sf y})$ admits a unique continuous
prolongation to $\square_\rho^n$, still denoted by $f (\cdot, {\sf
y})$, that is $\mathcal{ C}^{ 0, \beta}$ in $\square_\rho^n$ with
\[
\left\vert\!\left\vert
f(\cdot,{\sf y})
\right\vert\!\right\vert_{\mathcal{C}^{0,\beta}(\square_\rho^n)}
\leqslant
{\sf K}.
\]
\end{lemma}

Thus, for every $(s, U_0) \not \in \mathfrak{ N}$,
the partial derivatives $U_{ s_k}$, $U_{ U_0^m}$ belong to $\mathcal{
C}^{ 0, \alpha /2} (\partial \Delta, \R^d)$. Since the operator ${\sf
T}_1$ is linear and bounded in $\mathcal{ C}^{ 0, \alpha /2}$, we may
differentiate the $d$ scalar Bishop-type equations
of Theorem~3.7 with
respect to $\theta$, to $s_k$, $k=1, \dots, b$ and to $U_0^m$, $m=1,
\dots, d$, which yields, for $j=1, \dots, d$:
\def\theequation{3.20}\begin{equation}
\left\{
\aligned
U_\theta^j
(e^{i\,\theta})
&
=
-
{\sf T}_1
\left[
\sum_{1\leqslant l\leqslant d}\,\Phi_{u_l}^j
\left(
U(\cdot),\cdot,s
\right)
U_\theta^l(\cdot)
+
\Phi_\theta^j
\left(
U(\cdot),\cdot,s
\right)
\right]
(e^{i\,\theta}),
\\
U_{s_k}^j
(e^{i\,\theta})
&
=
-
{\sf T}_1
\left[
\sum_{1\leqslant l\leqslant d}\,\Phi_{u_l}^j
\left(
U(\cdot),\cdot,s
\right)
U_{s_k}^l(\cdot)
+
\Phi_{s_k}^j
\left(
U(\cdot),\cdot,s
\right)
\right]
(e^{i\,\theta}),
\\
U_{U_0^m}^j
(e^{i\,\theta})
&
=
\delta_m^j
-
{\sf T}_1
\left[
\sum_{1\leqslant l\leqslant d}\,\Phi_{u_l}^j
\left(
U(\cdot),\cdot,s
\right)
U_{U_0^m}^l(\cdot)
\right]
(e^{i\,\theta}),
\endaligned\right.
\end{equation}
for every $e^{ i\, \theta} \in \partial \Delta$, provided $(s, U_0)
\not \in \mathfrak{ N}$. In the first line, we noticed that $({\sf
T}_1 V)_\theta = {\sf T} (V_\theta)$. We observe that as $U$ is
Lipschitzian, as $\Phi \in \mathcal{ C}^{ \kappa, \alpha}$ and as
$\kappa \geqslant 1$, the composite functions $\Phi_{u_l}^j, \Phi_\theta^j,
\Phi_{ s_i}^j \left( U ( e^{ i\, \theta}, s, U_0), e^{ i\, \theta},
s\right)$ are of class $\mathcal{ C}^{ 0, \alpha}$ with respect to all
variables.

We notice that in each of the three linear systems of Bishop-type
equations~\thetag{ 3.20} above, $t:= (s, U_0)$ is a parameter and there
appears the same matrix coefficients:
\[
p_l^j(e^{i\,\theta},t)
:= 
\Phi_{u_l}^j
\left(
U(e^{i\,\theta},s,U_0),e^{i\,\theta},s
\right),
\]
for $j,l= 1, \dots, d$. For any $\beta$ with $0 < \beta \leqslant \alpha$,
in order to be coherent with the definition of
$\left\vert \! \left\vert \Phi_u\vert_s \right\vert \! \right\vert_{
\widehat{0, \beta }}$ given after
Lemma~3.9, we set:
\[
\left\vert\!\left\vert 
p\vert_t
\right\vert\!\right\vert_{\widehat{0,\beta}}
:= 
\max_{1\leqslant j\leqslant d} 
\left( 
\sum_{1\leqslant l\leqslant d}\, 
\big\vert\!\big\vert 
p_l^j\vert_t
\big\vert\!\big\vert_{\widehat{0,\beta}} 
\right).
\]
We also set{\rm :}
\[
\left\vert\!\left\vert 
p 
\right\vert\!\right\vert_{0,0}
:=
\max_{1\leqslant j\leqslant d}
\left(
\sum_{1\leqslant l\leqslant d}\,
\big\vert\!\big\vert
p_l^j
\big\vert\!\big\vert_{0,0}
\right)
=
\left\vert\!\left\vert 
\Phi_u
\right\vert\!\right\vert_{0,0}.
\]
With these definitions, it is easy to check the inequality:
\[
\left\vert\!\left\vert 
p\vert_t
\right\vert\!\right\vert_{0,\beta}
\leqslant
\left\vert\!\left\vert 
\Phi_u\vert_s
\right\vert\!\right\vert_{0,\beta}
\left[
1
+
\left(
\left\vert\!\left\vert 
U
\right\vert\!\right\vert_{1,0}
\right)^\beta
\right].
\]
As $\left\vert \! \left\vert U \right\vert \! \right\vert_{ 1,
0} \leqslant \rho_1 / 4 < 1$, 
with $\beta := \alpha$, we deduce:
\[
\left\vert\!\left\vert 
p\vert_t
\right\vert\!\right\vert_{ 
0,\alpha} 
\leqslant 
2\left\vert\!\left\vert 
\Phi_u\vert_s 
\right\vert\!\right\vert_{0,\alpha} 
\leqslant 
2 
\left\vert\!\left\vert 
\Phi\vert_s 
\right\vert\!\right\vert_{ 
1,\alpha}.
\]
Taking $\sup_s$ and then $\sup_{U_0}$, adding $1$,
squaring and inverting:
\[
\left[
1
+
\sup_s
\left\vert\!\left\vert
\Phi\vert_s
\right\vert\!\right\vert_{1,\alpha}
\right]^{-2}
\leqslant
4\left[
1
+
\sup_{s,U_0}
\left\vert\!\left\vert
p\vert_t
\right\vert\!\right\vert_{0,\alpha}
\right]^{-2}.
\]
Consequently, the main assumption of the next Proposition~3.21,
according to which:
\[
\left\vert\!\left\vert
p
\right\vert\!\right\vert_{0,0}
\leqslant
C^2\alpha^2
\left[
1
+
\sup_t\,
\left\vert\!\left\vert
p\vert_t
\right\vert\!\right\vert_{0,\alpha}
\right]^{-2},
\]
for some positive absolute constant $C < 1$, is superseded by one of
the main assumptions, made in the theorem, 
according to which:
\[
\left\vert\!\left\vert
\Phi_u
\right\vert\!\right\vert_{0,0}
\leqslant
C^2\alpha^2
\left[
1
+
\sup_s\,
\left\vert\!\left\vert
\Phi\vert_s
\right\vert\!\right\vert_{0,\alpha}
\right]^{-2},
\]
for some ({\it a priori} distinct) positive absolute constant $C < 1$.

The following proposition applies to the three systems~\thetag{ 3.20}
and suffices to conclude the proof of Theorem~3.7 in the case $\kappa
= 1$. The case $\kappa \geqslant 2$ shall be discussed afterwards.

\def\theproposition{3.21}\begin{proposition}
{\rm (\cite{ tu1996}, [$*$])} Let $c\in \N$ with $c\geqslant 1$, let
$\tau_1 = (\tau_{ 1, 1}, \dots, \tau_{ 1, c}) \in \R^c$ with $0 <
\tau_{ 1, i} \leqslant \infty$, $i=1, \dots, c$, and denote by
$\square_{\tau_1}^c$ the polycube $\{ t\in \R^c: \vert t_i \vert <
\tau_{ 1,i}\}$. Consider vector-valued and matrix-valued H\"older
data{\rm :}
\[
\aligned
q
&
=
\left(q^j(e^{i\,\theta},t) 
\right)^{
1\leqslant j\leqslant d}
\ \in \
\mathcal{C}^{0,\alpha}
\left(
\partial\Delta\times\square_{\tau_1}^c,\R^d
\right),
\\
p
&
=
\left(p_l^j(e^{i\,\theta},t) 
\right)_{1\leqslant l\leqslant d}^{
1\leqslant j\leqslant d} \
\in\
\mathcal{C}^{0,\alpha}
\left(
\partial\Delta\times\square_{\tau_1}^c,\R^{d\times d}
\right),
\endaligned
\]
with $0 < \alpha < 1$. Suppose that a bounded measurable map{\rm :}
\[
u 
= 
\left( u^j (e^{i\, \theta}, t) \right)^{1\leqslant
j\leqslant d}\
\in\
L^\infty ( \partial \Delta \times
\square_{ \tau_1}^c, \R^d), 
\]
is $\mathcal{ C}^{ 0, \alpha /2}$ on $\partial \Delta$ for every fixed
$t$ not belonging to some null-measure subset $\mathfrak{ N}$ of
$\square_{ \tau_1}^c$ and suppose that it satisfies the system of {\rm
linear} Bishop-type equations in $\mathcal{ C}^{ 0, \alpha /2} (\partial
\Delta, \R^d)${\rm :}
\def\theequation{3.22}\begin{equation}
u^j
=
{\sf T}_*
\left(
\sum_{1\leqslant l\leqslant d}\,
p_l^j\,u^l
\right)
+
q^j,
\end{equation}
for $j=1, \dots, d$, where ${\sf T}_* = {\sf T}$ or ${\sf T}_* = {\sf
T}_1$. Assume that the norm of the matrix $p$ satisfies{\rm :}
\[
\left\vert\!\left\vert
p
\right\vert\!\right\vert_{0,0}
=
\max_{1\leqslant j\leqslant d}
\left(
\sum_{1\leqslant l\leqslant d}\,
\big\vert\!\big\vert
p_l^j
\big\vert\!\big\vert_{0,0}
\right)
\leqslant
{\sf c}_4,
\]
for some small positive constant ${\sf c}_4 < 1 /16$ such that
\def\theequation{3.23}\begin{equation}
{\sf c}_4
\leqslant
C^2\,
\alpha^2
\left[
1
+
\sup_{\vert t\vert<\tau_1}
\left\vert\!\left\vert
p\vert_t
\right\vert\!\right\vert_{0,\alpha}
\right]^{-2},
\end{equation}
where $C < 1$ is a positive absolute constant. Then,
after a correction of $u$ on $\mathfrak{ N}${\rm :}

\begin{itemize}

\smallskip\item[{\bf (i)}] 
on its full domain of definition $\partial \Delta \times \square_{
\tau_1}^c$, the corrected map $u ( e^{ i\, \theta}, t)$ is $\mathcal{
C}^{ 0, \alpha - 0} = \bigcap_{ 0 < \beta < \alpha} \, \mathcal{ C}^{
0, \beta}$ and furthermore{\rm :}

\smallskip\item[{\bf (ii)}] 
for every fixed $t$, the map $e^{i\, \theta} \mapsto u ( e^{ i\,
\theta}, t)$ is $\mathcal{ C}^{ 0, \alpha}$ on $\partial \Delta$.

\end{itemize}
\end{proposition}

In general, the Hilbert transform ${\sf T}$ does not preserve
$\mathcal{ C}^{ 0, \alpha}$ smoothness with respect to parameters, so
that the above solution $u = u (e^{ i\, \theta}, t)$ is not better
than $\mathcal{ C}^{ 0, \alpha - 0}$.

\def\theexample{3.24}\begin{example}
{\rm (\cite{ tu1996})
If $s\in \R$ with $\vert s \vert < 1$ is a parameter, the
function:
\[
\aligned
u(e^{i\,\theta},s)
&
:=
\vert s\vert^\alpha\ \ {\rm if}\
-\pi\leqslant \theta\leqslant -\vert s\vert^{1/2},
\\
&
:=
\theta^{2\alpha}\ \ \ {\rm if}
-\vert s\vert^{1/2}\leqslant \theta\leqslant 0,
\\
&
:=
\theta^\alpha\ \ \ {\rm if}\
0\leqslant \theta\leqslant\vert s\vert,
\\
&
:=
\vert s\vert^\alpha\ \ {\rm if}\
\vert s\vert\leqslant \theta\leqslant \pi,
\endaligned
\]
is $2\pi$-periodic with respect to $\theta$ and $\mathcal{ C}^{ 0,
\alpha}$ with respect to $(e^{i\, \theta}, s)$. As the function
$\cot( t/2) - 2/t$ is $\mathcal{ C}^{ 0, 0}$ on $[ -\pi, \pi]$, the
regularity properties of the singular integral ${\sf T}u
(e^{i\,\theta}) = {\rm p.v.}\, \frac{1}{ \pi} \int_{- \pi }^\pi\,
\frac{u(e^{i (\theta-t)})}{ \tan(t/2) }\,dt$ are the same as those
of{\rm :}
\[
\widetilde{\sf T}u(e^{i\,\theta})
:=
{\rm p.v.}\,\frac{1}{\pi}\int_{-\pi}^\pi\,
\frac{u(e^{i(\theta-t)})}{t}\,dt.
\]
However $\widetilde{ \sf T}u (1)$ involves the term $\vert
s\vert^\alpha \, {\rm log} \, \vert s\vert$ which is $\mathcal{ C}^{
0, \alpha -0}$ but not $\mathcal{ C}^{ 0, \alpha}$:
\[
\aligned
\widetilde{\sf T}u(1)
&
=
\frac{1}{\pi}
\left(
\int_{-\vert s\vert^{1/2}}^{-\vert s\vert}\,
\frac{\vert s\vert^\alpha}{t}\,dt
+
\int_{-\vert s\vert}^0\,
\frac{(-t)^\alpha}{t}\,dt
+
\int_{0}^{\vert s\vert^{1/2}}\,
\frac{t^{2\alpha}}{t}\,dt
\right)
\\
&
=
\frac{1}{2\pi}
\left(
\vert s\vert^\alpha\,{\rm log}\,\vert s\vert
-
\frac{\vert s\vert^\alpha}{\alpha}
\right).
\endaligned
\]
}\end{example}

\proof[Proof of the proposition]
We shall drop the indices, writing $u$, $p$ and $q$, without
arguments. Assume that $t \not \in \mathfrak{ N}$. For future
majorations, it is necessary to have ${\sf P}_0 u =0$. If this is not
the case, we set $u' := u - {\sf P}_0 u$ in order that ${\sf P}_0 \,
u' = 0$. Since $u$ satisfies either $u = {\sf T} (pu) + q$ or $u =
{\sf T} (pu) - {\sf T} [pu] (1) + q$, it follows that $u'$ satisfies
similar equations: either $u' = {\sf T} (pu') + q'$, with $q' := q -
{\sf P}_0 u$ or $u' = {\sf T} (pu') + q'$, with $q' := q - {\sf P}_0 u
- {\sf T}(pu)(1)$. Notice that $p$ is unchanged. It then suffices to
establish the improvements of smoothness {\bf (i)} and {\bf (ii)} for
$u'$. Equivalently, we may assume that $\widehat{ u}_0 = {\sf P}_0 \,
u = 0$ in the proposition.

For $t\not\in \mathfrak{ N}$, the function $\phi$ is $\mathcal{ C}^{
0, \alpha /2}$ on $\partial \Delta$. Applying ${\sf T}$ either to the
equation $u = {\sf T} (pu) + q$ or to the equation $u = {\sf T} (pu) -
{\sf T} [pu] (1) + q$, we get the {\it same}\, equation for both:
\[
{\sf T}u
=
-
pu
+
{\sf P}_0(pu)
+
{\sf T}q.
\]
As $\widehat{ u}_0 = 0$, we may write $u ( e^{i\, \theta} ) = \sum_{
k< 0} \, \widehat{ u}_k \, e^{ i\, k\, \theta} + \sum_{ k> 0} \,
\widehat{ u}_k \, e^{i\, k\, \theta} = : \overline{ \phi} + \phi$,
where $\phi$ extends holomorphically to $\Delta$. In fact, $\phi$ is
determined up to a imaginary constant $i A$ and we choose $A := -
{\sf P}_0 (pu)/2$, so that:
\def\theequation{3.25}\begin{equation}
\phi
=
\frac{u+i\,{\sf T}u
-
i\,{\sf P}_0(pu)}{2}.
\end{equation}
Equivalently:
\[
\left\{
\aligned
u
&
=
\phi
+
\overline{\phi},
\\
{\sf T}u
&
=
{\sf P}_0(pu)
-
i(\phi-\overline{\phi}).
\endaligned\right.
\]
Substituting, we rewrite~\thetag{ 3.25} as:
\[
-
i(\phi-\overline{\phi})
=
-
p(\phi+\overline{\phi})
+
{\sf T}q,
\]
or under the more convenient form:
\[
\phi
=
\overline{\phi}
+
P\,\overline{\phi}
+
Q,
\]
where the $d\times d$-matrix 
$P := -2\, i\, p \, ( I + i\, p)^{ -1}$ and the
$d$-vector $Q := i\, (I+i\, p)^{ -1}\, {\sf T} q$ both belong to
$\mathcal{ C}^{ 0, \alpha}$.

First of all, we establish {\bf (ii)} before any correction of $u$.

\def\thelemma{3.26}\begin{lemma}
For $t\not \in \mathfrak{ N}$, the map $e^{ i\, \theta} \mapsto u (e^{
i\, \theta}, t)$ is $\mathcal{ C}^{ 0, \alpha}$ on $\partial \Delta$.
\end{lemma}

\proof 
By assumption, the map $e^{ i\,
\theta} \mapsto \phi ( e^{ i\, \theta}, t)$ is $\mathcal{ C}^{ 0,
\alpha/2}$ on $\partial \Delta$. Since ${\sf C}^+$
is bounded in $\mathcal{ C}^{ 0, \alpha/2}$, we may apply ${\sf
C}^+$ to the vectorial equation $\phi = \overline{ \phi} + P\,
\overline{ \phi} + Q \, \overline{ \phi}$, noticing that ${\sf C}^+
(\phi) = \phi$ and that ${\sf C}^+ (\overline{ \phi}) = {\sf P}_0 (
\overline{ \phi})$, since, by construction, $\phi$ extends
holomorphically to $\Delta$. We thus get:
\def\theequation{3.27}\begin{equation}
\phi
=
{\sf P}_0\,\overline{\phi}
+
{\sf C}^+
\left[
P\,\overline{\phi}
+
Q
\right].
\end{equation}
Remind that ${\sf P}_0 \, \psi = \frac{ 1}{ 2\pi} \, \int_{-\pi}^\pi
\, \psi(e^{ i\, \theta}) \, d\theta$, so that $\left\vert \!
\left\vert {\sf P}_0 \, \psi \right\vert \! \right\vert_{ 0, 0} \leqslant
\left\vert \! \left\vert \psi \right\vert \! \right\vert_{ 0, 0}$.
Notice\footnote{ In Lemma~3.15 above, for $t\not \in \mathfrak{ N}$,
the map $e^{ i\,\theta} \mapsto u (e^{ i\, \theta}, t)$ was shown to
be $\mathcal{ C}^{ 0, \alpha /2}$ on $\partial \Delta$ in order to
insure that ${\sf T}u (\cdot, t)$ and $\phi (\cdot, t)$ are both
bounded on $\partial \Delta$, so that Theorem~2.24 may be applied in
the next phrase. In~\cite{ tu1996}, $u (\cdot, t)$ is only shown to
be in $L^\infty (\partial \Delta)$, but then ${\sf T} u (\cdot, t)$
and $\phi$ are not necessarily bounded.} that $\left\vert \!
\left\vert \phi (\cdot, t) \right\vert \! \right\vert_{ 0, 0} <
\infty$ for $t \not \in \mathfrak{ N}$. Taking the $\mathcal{ C}^{
0,\alpha}$ norm to~\thetag{ 3.27} and applying crucially Theorem~2.18
in its $\R^d$-valued version, as in~\cite{ tu1996}, we get:
\[
\aligned
\left\vert\!\left\vert
\phi
\right\vert\!\right\vert_{0,\alpha}
&
\leqslant
\left\vert\!\left\vert
{\sf P}_0\,\overline{\phi}
\right\vert\!\right\vert_{0,0}
+
\left\vert\!\left\vert
{\sf C}^+
\left[
P\,\overline{\phi}
+
Q
\right]
\right\vert\!\right\vert_{0,\alpha}
\\
&
\leqslant
\left\vert\!\left\vert
\phi
\right\vert\!\right\vert_{0,0}
+
\frac{C_1}{\alpha(1-\alpha)}
\left(
\left\vert\!\left\vert
P
\right\vert\!\right\vert_{0,\alpha}\,
\left\vert\!\left\vert
\phi
\right\vert\!\right\vert_{0,0}
+
\left\vert\!\left\vert
Q
\right\vert\!\right\vert_{0,\alpha}
\right)
\\
&
<
\infty,
\endaligned
\]
whence $\phi = \phi ( \cdot, t)$ is $\mathcal{ C}^{ 0, \alpha}$ on
$\partial \Delta$, as claimed.
\endproof

Next, we shall establish {\bf (i)} before any correction of $u$. To
this aim, let $t^1, t^2 \not\in \mathfrak{ N}$ and simply denote by
$u^1, u^2$, by $\phi^1, \phi^2$, by $P^1, P^2$ and by $Q^1, Q^2$ the
functions on $\partial \Delta$ with these two values of $t$. In the
theorem, to establish that $u$ is $\mathcal{ C}^{ 0, \alpha - 0}$, it
is natural to show that $u$ is $\mathcal{ C}^{ 0, \beta}$ for every
$\beta < \alpha$ arbitrarily close to $\alpha$.

\def\thelemma{3.28}\begin{lemma} 
For every $\beta$ with $\frac{ \alpha}{ 2} < \beta < \alpha$, every
two $t^1, t^2 \not \in \mathfrak{ N}$, we have $\left \vert \! \left
\vert u^2 - u^1 \right \vert \! \right \vert_{ 0, 0} \leqslant 4\, {\sf
K}_{\alpha, \beta} \, \left\vert t^2 - t^1 \right\vert^\beta$, for
some positive quantity ${\sf K}_{ \alpha, \beta} < \infty$.
\end{lemma}

We shall obtain ${\sf K}_{ \alpha, \beta}$ involving a nonremovable
factor $1/ (\alpha - \beta)$. This will confirm the optimality of
$\mathcal{ C}^{ 0, \alpha - 0}$-smoothness of $u$ with respect to 
the parameter $t$.

\proof
Since ${\sf P}_0 \, u = {\sf P_0} ({\sf T} u ) =0$, applying ${\sf
P}_0$ to the conjugate of~\thetag{ 3.25}, we get ${\sf P}_0 \,
\overline{ \phi} = \frac{ i}{ 2}\, {\sf P}_0 (pu)$, so that we may
rewrite~\thetag{ 3.27} under the form:
\[
\phi
=
\frac{i}{2}\,{\sf P}_0(pu)
+
{\sf C}^+
\left[
P\,\overline{\phi}
+
Q
\right].
\]
We may then organize the difference $\phi^2 - \phi^1$
as follows:
\[
\aligned
\phi^2
- 
\phi^1
&
=
\frac{i}{2}\,{\sf P}_0\left(p^2(u^2-u^1)\right)
+
\frac{i}{2}\,{\sf P}_0\left((p^2-p^1)u^1\right)
+
\\
&
\ \ \ \ \ 
+
{\sf C}^+\left((P^2-P^1)\overline{\phi}^2
+
Q^2-Q^1\right)
+
{\sf C}^+
\left(
P^1(\overline{\phi}^2-\overline{\phi}^1)
\right)
\\
&
=:
E_1+E_2+E_3+E_4.
\endaligned
\]
To majorate these four $E_i$'s, we shall need various
preliminaries.

\smallskip 

Firstly, to majorate $E_1$, we first observe the elementary
inequality:
\def\theequation{3.29}\begin{equation}
\left\vert\!\left\vert
u^2
-
u^1
\right\vert\!\right\vert_{0,0}
\leqslant 2\, 
\left\vert\!\left\vert
\phi^2
-
\phi^1
\right\vert\!\right\vert_{0,0}.
\end{equation}
Also, we notice that:
\[
\left\vert\!\left\vert
p^2
\right\vert\!\right\vert_{0,0}
=
\left\vert\!\left\vert
p(\cdot,t^2)
\right\vert\!\right\vert_{\mathcal{C}^{0,0}(\partial\Delta)}
\leqslant
\left\vert\!\left\vert
p
\right\vert\!\right\vert_{0,0}
\leqslant
{\sf c}_4.
\]
The majoration of $E_1$ is then easy:
\[
\left\vert\!\left\vert
E_1
\right\vert\!\right\vert_{0,0}
\leqslant
\frac{1}{2}\,
\left\vert\!\left\vert
p^2
\right\vert\!\right\vert_{0,0}\,
\left\vert\!\left\vert
u^2-u^1
\right\vert\!\right\vert_{0,0}
\leqslant
4\,{\sf c}_4
\left\vert\!\left\vert
\phi^2-\phi^1
\right\vert\!\right\vert_{0,0}.
\]

\smallskip

Secondly, to majorate $E_2$, let again $\beta$ with $\frac{ \alpha}{
2} < \beta < \alpha$. Using the inequality $\left\vert \! \left\vert
p \right\vert \! \right\vert_{ 0, \beta} \leqslant 3 \left\vert \!
\left\vert p \right\vert \! \right\vert_{ 0, \alpha}$ proved in the
beginning of Section~1, we may majorate $E_2$:
\[
\left\vert\!\left\vert
E_2
\right\vert\!\right\vert_{0,0}
\leqslant
\frac{1}{2}\,
\left\vert\!\left\vert
u^1
\right\vert\!\right\vert_{0,0}\,
\left\vert\!\left\vert
p^2-p^1
\right\vert\!\right\vert_{0,0}
\leqslant
\frac{3}{2}\,
\left\vert\!\left\vert
u^1
\right\vert\!\right\vert_{0,0}\,
\left\vert\!\left\vert
p
\right\vert\!\right\vert_{0,\alpha}\,
\left\vert
t^2-t^1
\right\vert^\beta.
\]

\smallskip
\noindent
Thirdly, to majorate $E_3$, we need:

\def\thelemma{3.30}\begin{lemma}
Let $f = f ({\sf x}, {\sf y})$ be a $\mathcal{ C}^{0, \alpha}$
function, defined in the product open cube $\square_\rho^n \times
\square_\rho^m$, where $0 < \alpha < 1$ and $\rho >0$. Let $\beta$
with $0 < \beta < \alpha$. For ${\sf x}', {\sf x}'' \in
\square_\rho^n$ arbitrary, the $\mathcal{ C}^{ 0, \alpha -
\beta}$-norm of the function $\square_\rho^m \ni {\sf y} \longmapsto
f( {\sf x}'', {\sf y}) - f( {\sf x}', {\sf y}) \in \R$ satisfies{\rm
:}
\[
\left\vert\!\left\vert
f({\sf x}'',\cdot)
-
f({\sf x}',\cdot)
\right\vert\!\right\vert_{0,\alpha-\beta}
\leqslant
4
\left\vert\!\left\vert
f
\right\vert\!\right\vert_{0,\alpha}
\left\vert
{\sf x}''
-
{\sf x}'
\right\vert^{\beta}.
\]
\end{lemma}

\noindent
Again, assume $\frac{ \alpha}{2} < \beta < \alpha$. Thanks to
Theorem~2.18 and to the lemma above, we may majorate $E_3$:
\[
\aligned
\left\vert\!\left\vert
E_3
\right\vert\!\right\vert_{0,0}
\leqslant
\left\vert\!\left\vert
E_3
\right\vert\!\right\vert_{0,\alpha-\beta}
&
\leqslant
\frac{C_1}{\alpha - \beta}\,
\left(
\left\vert\!\left\vert
P^2-P^1
\right\vert\!\right\vert_{0,\alpha-\beta}\,
\left\vert\!\left\vert
\phi^2
\right\vert\!\right\vert_{0,0}
+
\left\vert\!\left\vert
Q^2-Q^1
\right\vert\!\right\vert_{0,\alpha-\beta}
\right)
\\
&
\leqslant
\frac{C_2}{\alpha - \beta}\,
\left(
\left\vert\!\left\vert
P
\right\vert\!\right\vert_{0,\alpha}\,
\left\vert\!\left\vert
\phi^2
\right\vert\!\right\vert_{0,0}
+
\left\vert\!\left\vert
Q
\right\vert\!\right\vert_{0,\alpha}
\right)
\left\vert
t^2-t^1
\right\vert^\beta.
\endaligned
\]

Fourthly, to majorate $E_4$, we need a statement whose proof is
similar to that of Lemma~3.10.

\def\thelemma{3.31}\begin{lemma}
As $\left\vert\!\left\vert p \right\vert\!\right\vert_{0,0} \leqslant {\sf
c}_4$, then independently of $t \not \in \mathfrak{ N}$, we have{\rm
:}
\[
\vert\!\vert
p_{\vert_t}
\vert\!\vert_{0,\alpha/2}
\leqslant
{\sf c}_4
+
{\sf c}_4^{1/2}
\left[
2
+
\sup_t\,
\left\vert\!\left\vert
p\vert_t
\right\vert\!\right\vert_{0,\alpha}
\right].
\]
\end{lemma}

\noindent
Reminding the main assumption ${\sf c}_4 \leqslant C^2 \alpha^2 \left[ 1 +
\sup_t \left\vert \! \left\vert p\vert_t \right\vert \!
\right\vert_{ 0, \alpha} \right]^{ -2}$, whence ${\sf c}_4 \leqslant
C\alpha$, we deduce:
\[
\vert\!\vert
p_{\vert_t}
\vert\!\vert_{0,\alpha/2}
\leqslant
3\,C\alpha.
\]
Choosing then $C < 1/6$, we may insure that $\vert \! \vert p_{
\vert_t} \vert \! \vert_{ 0, \alpha /2} \leqslant 1/2$ for every fixed
$t$. It follows in particular that we may develope $P = -2\, i\, p\,
\sum_{ k\in \N}\, (-ip)^k$ and deduce the norm inequality:
\[
\left\vert\!\left\vert
P_{\vert_t}
\right\vert\!\right\vert_{0,\alpha/2}
\leqslant
\frac{
2\,
\vert\!\vert
p_{\vert_t}
\vert\!\vert_{0,\alpha/2}}{
1
-
\vert\!\vert
p_{\vert_t}
\vert\!\vert_{0,\alpha/2}}
\leqslant
4\,
\vert\!\vert
p_{\vert_t}
\vert\!\vert_{0,\alpha/2},
\]
valid for every fixed $t$. Then again thanks to Theorem~2.18 and
thanks to the previous inequalities:
\[
\aligned
\left\vert\!\left\vert
E_4
\right\vert\!\right\vert_{0,0}
\leqslant
\left\vert\!\left\vert
E_4
\right\vert\!\right\vert_{0,\alpha/2}
&
\leqslant
C_1\alpha^{-1}\,
\left\vert\!\left\vert
P^1
\right\vert\!\right\vert_{0,\alpha/2}\,
\left\vert\!\left\vert
\phi^2-\phi^1
\right\vert\!\right\vert_{0,0}
\\
&
\leqslant
C_2\alpha^{-1}\,
\vert\!\vert
p\vert_{t^1}
\vert\!\vert_{0,\alpha/2}\,
\left\vert\!\left\vert
\phi^2-\phi^1
\right\vert\!\right\vert_{0,0}
\\
&
\leqslant
C_3\,C
\left\vert\!\left\vert
\phi^2-\phi^1
\right\vert\!\right\vert_{0,0}
\\
&
\leqslant
4^{-1}\,
\left\vert\!\left\vert
\phi^2-\phi^1
\right\vert\!\right\vert_{0,0},
\endaligned
\]
provided $C < \frac{ 1}{ 4 C_3}$. We then insert these four
majorations in the inequality
\[
\left\vert\!\left\vert
\phi^2-\phi^1
\right\vert\!\right\vert_{0,0}
\leqslant
\left\vert\!\left\vert
E_1
\right\vert\!\right\vert_{0,0}
+
\left\vert\!\left\vert
E_2
\right\vert\!\right\vert_{0,0}
+
\left\vert\!\left\vert
E_3
\right\vert\!\right\vert_{0,0}
+
\left\vert\!\left\vert
E_4
\right\vert\!\right\vert_{0,0},
\]
and we get after reorganization:
\[
\left\vert\!\left\vert
\phi^2-\phi^1
\right\vert\!\right\vert_{0,0}
\left(
1
-
4\,{\sf c}_4
-
4^{-1}
\right)
\leqslant
{\sf K}_{\alpha,\beta}\,
\left\vert
t^2-t^1
\right\vert^\beta,
\]
for some quantity ${\sf K}_{\alpha, \beta} < \infty$ whose precise
expression is:
\[
K_{\alpha,\beta}
:=
\frac{C_2}{\alpha-\beta}
\left(
\left\vert\!\left\vert
P
\right\vert\!\right\vert_{0,\alpha}\,
\left\vert\!\left\vert
\phi^2
\right\vert\!\right\vert_{0,0}
+
\left\vert\!\left\vert
Q
\right\vert\!\right\vert_{0,\alpha}\,
\right)
+
\frac{3}{2}\,
\left\vert\!\left\vert
u^1
\right\vert\!\right\vert_{0,0}\,
\left\vert\!\left\vert
p
\right\vert\!\right\vert_{0,\alpha}.
\]
It suffices then to remind that ${\sf c}_4 < 1/16$ in the assumptions
of the theorem to insure the uniform H\"older condition:
\[
\left\vert\!\left\vert
\phi^2-\phi^1
\right\vert\!\right\vert_{0,0}
\leqslant
2\,
{\sf K}_{\alpha,\beta}\,
\left\vert
t^2-t^1
\right\vert^\beta,
\]
valid for $t^1, t^2 \not \in \mathfrak{ N}$. From~\thetag{ 3.29}, we
conclude that $\left\vert \! \left\vert u^2-u^1 \right\vert \!
\right\vert_{0,0} \leqslant 4\, {\sf K}_{ \alpha, \beta}\, \left\vert
t^2-t^1 \right\vert^\beta$. The proof of Lemma~3.28 is complete.
\endproof

\smallskip

Then the correction of $u$ is provided by the following statement.

\def\thelemma{3.32}\begin{lemma}
{\rm ([$*$])} Let $f = f({\sf x}, {\sf y}) : \square_\rho^n \times (
\square_\rho^m \backslash \mathfrak{ N} ) \to \R$ be a measurable
$L^\infty$ map defined only for ${\sf y}$ not belonging to some
null-measure subset $\mathfrak{ N} \subset \square_\rho^m$ and let
$\beta$ with $0 < \beta < \alpha$. If the map ${\sf x} \mapsto f
({\sf x}, {\sf y})$ is $\mathcal{ C}^{ 0, \beta}$ for every ${\sf y}
\not \in \mathfrak{ N}$ and if there exists ${\sf K} < \infty$ such
that{\rm :}
\[
\sup_{{\sf x}}\,
\left\vert
f({\sf x},{\sf y}^2)
-
f({\sf x},{\sf y}^1)
\right\vert
\leqslant
{\sf K}
\left\vert
{\sf y}^2
-
{\sf y}^1
\right\vert^\beta,
\]
for every two ${\sf y}^1, {\sf y}^2 \not \in \mathfrak{ N}$, 
then $f$ may be extended as a $\mathcal{ C}^{ 0,
\beta}$ map in the full domain $\square_\rho^n \times \square_\rho^m$.
\end{lemma}
 
In sum, still denoting by $u$ the $\mathcal{ C}^{ 0, \alpha -0}$
extension of $u$ through $\mathfrak{ N}$, property {\bf (i)} of the
proposition is proved. To obtain that $u$ is $\mathcal{ C}^{ 0,
\alpha}$ with respect to $e^{ i\, \theta}$, we re-apply the reasoning
of Lemma~3.26 to this extension.

The proof of Proposition~3.21 is complete.
\endproof

In conclusion, the functions $U_\theta^j$, $U_{s_k}^j$ and
$U_{U_0^m}^j$ are $\mathcal{ C}^{ 0, \alpha - 0}$ with respect to
$(e^{i\, \theta}, s, U_0)$ and $\mathcal{ C}^{ \alpha, 0}$ with
respect to $e^{ i\, \theta}$. Thus, the theorem is achieved if
$\kappa = 1$.

\smallskip

If $\kappa = 2$, the composite functions $\Phi_{ u_l}^j,
\Phi_\theta^j, \Phi_{ s_i}^j \left( U (e^{ i\, \theta}, s, U_0), e^{
i\, \theta}, s \right)$ are then of class $\mathcal{ C}^{ 1, \alpha -
0}$ with respect to $(e^{i\, \theta}, s, U_0)$ and of class $\mathcal{
C}^{ 1, \alpha}$ with respect to $e^{ i\, \theta}$. We then apply the
next lemma to the three families of Bishop-type 
vector equations~\thetag{ 3.20}.

\def\thelemma{3.33}\begin{lemma} 
Let $t \in \square_{ \tau_1}^c$ be a parameter with $c \in \N$,
$0< \tau_{1,i} \leqslant \infty$, 
$i=1, \dots, c$, and consider vector-valued and matrix-valued H\"older
data $q = \left( q^j (e^{ i\, \theta}, t) \right)^{ 1\leqslant j\leqslant d}$
and $p = \left( p_l^j (e^{ i\, \theta}, t) \right)_{ 1\leqslant l\leqslant d}^{
1\leqslant j\leqslant d}$ that are $\mathcal{ C}^{ 1, \alpha - 0}$ with respect
to $(e^{ i\, \theta}, t)$ and $\mathcal{ C}^{ 1, \alpha}$ with respect
to $e^{ i\, \theta}$. Suppose that a given map $u = \left( u^j (e^{
i\, \theta}, t) \right)^{ 1\leqslant j\leqslant d}$ which is $\mathcal{ C}^{ 0,
\alpha - 0}$ with respect to $(e^{ i\, \theta}, t)$ and $\mathcal{
C}^{ 0, \alpha}$ with respect to $e^{ i\, \theta}$ satisfies the
linear Bishop-type equation $u = {\sf T}_* (pu) + q$, where ${\sf T}_*
= {\sf T}$ or ${\sf T}_* = {\sf T}_1$. Assume that the norm of the
matrix $p$ satisfies the same inequality as in the proposition{\rm :}
$\left\vert \! \left\vert p \right\vert \! \right\vert_{0,0} \leqslant
{\sf c}_4$, for some small positive constant ${\sf c}_4 
\leqslant C^2\,
\alpha^2 \left[ 1 + \sup_t \left\vert\!\left\vert p\vert_t
\right\vert\!\right\vert_{0,\alpha} \right]^{-2}$, where $0 < C < 1$
is some absolute constant. Then $u$ is $\mathcal{ C}^{ 1, \alpha}$
with respect to $e^{ i\, \theta}$ and satisfies the Lipschitz
condition
\[
\left\vert\!\left\vert
u(\cdot,t^2)
-
u(\cdot,t^1)
\right\vert\!\right\vert_{0,\alpha/2}
\leqslant
{\sf K}\,
\left\vert
t^2
-
t^1
\right\vert,
\]
for some quantity ${\sf K} < \infty$. Furthermore, there exists a
null-measure subset $\mathfrak{ N} \subset \square_{ \tau_1}^c \times
\square_{ \rho_1 /16}^d$ such that at every $t \not \in \mathfrak{ N}$,
for every $l= 1, \dots, c${\rm :}
\begin{itemize}

\smallskip\item[{\bf (i)}]
the partial derivative $u_{ t_l} ( e^{ i\, \theta}, t)$ exists {\rm
for every} $e^{ i\, \theta} \in \partial \Delta${\rm ;}
\smallskip\item[{\bf (ii)}]
the map $e^{i\, \theta} \mapsto u_{ t_l} ( e^{ i\, \theta}, t)$
is $\mathcal{ C}^{ 0, \alpha/2}$ on
$\partial \Delta$.
\end{itemize}
\end{lemma}

\proof
The fact that $u$ is $\mathcal{ C}^{ 1, \alpha}$ with respect to
$e^{i\, \theta}$ is proved as in Lemma~3.26. For the Lipschitz
condition, the reasoning is simpler than Lemma~3.17, due to the
linearity of $u = {\sf T}_* (pu) + q$. Indeed, for two parameters $t^1,
t^2 \in \square_{\tau_1}^c$, if we take the $\mathcal{ C}^{ 0, \alpha
/2}$-norm of the difference:
\[
u^2-u^1
=
{\sf T}_*
\left(
p^2(u^2-u^1)
\right)
+
{\sf T}_*
\left(
(p^2-p^1)u^1
\right)
+
q^2-q^1,
\]
we get:
\[
\left\vert\!\left\vert
u^2-u^1
\right\vert\!\right\vert_{0,\alpha/2}
\leqslant
C_1\alpha^{-1}\,
\left\vert\!\left\vert
p^2
\right\vert\!\right\vert_{0,\alpha/2}\,\left\vert\!\left\vert
u^2-u^1
\right\vert\!\right\vert_{0,\alpha/2}
+
{\sf K}
\left\vert
t^2-t^1
\right\vert,
\]
and after substraction, taking account of Lemma~3.31:
\[
\left\vert\!\left\vert
u^2-u^1
\right\vert\!\right\vert_{0,\alpha/2}
\leqslant
2\,{\sf K}
\left\vert
t^2-t^1
\right\vert.
\]
Then the sequel of the reasoning is already known.
\endproof

So for $l=1, \dots, c$, the partial derivatives $u_{ t_l}$
exist almost everywhere and
they satisfy:
\[
u_{t_l}
=
{\sf T}_*(p\,u_{t_l})
+
q_l,
\]
with the same matrix $p$, where $q_l := {\sf T}_* ( p_{ t_l} u) + q_{
t_l}$ is $\mathcal{ C}^{ 0, \alpha - 0}$ with respect to
$(e^{ i\, \theta}, t)$ and $\mathcal{ C}^{ 0, \alpha}$ with respect to
$e^{ i\, \theta}$.

\def\thelemma{3.34}\begin{lemma}
Proposition~3.21 holds true if more generally, $q$ is
only assumed to be $\mathcal{ C}^{ 0, \alpha - 0}$ with respect to 
$(e^{ i\, \theta}, t)$ and $\mathcal{ C}^{ 0, \alpha}$ with
respect to $e^{ i\, \theta}$, with the same conclusion.
\end{lemma}

(It suffices only to inspect the majoration of $E_3$.) Consequently,
with $u = U_\theta, U_{ s_k}, U_{ U_0^m}$ in the three
equations~\thetag{ 3.20}, we have verified that the partial derivatives
$u_\theta$, $u_{ s_k}$ and $u_{ U_0^m}$ exist everywhere, are
$\mathcal{ C}^{ 0, \alpha - 0}$ with respect to $(e^{i\, \theta}, s,
U_0)$ and are $\mathcal{ C}^{ 0, \alpha}$ with respect to $e^{ i\,
\theta}$. In summary, the theorem is achieved if $\kappa = 2$.

Needless to say, we have clarified how to cover the case of general
$\kappa \geqslant 2$ by pure logical induction. In conclusion, the proof
of Theorem~3.7 is complete.
\endproof

\def\theopenproblem{3.35}\begin{openproblem}
Solve parametrized Bishop-type equations
in Sobolev spaces.
\end{openproblem}

\section*{\S4.~Appendix: proofs of some lemmas}

\subsection*{4.1.~Proofs of Lemmas~3.9 and~3.10}
Let ${\sf x} \in \R^n$, $n\geqslant 1$, with $\vert {\sf x} \vert <
\rho$, where $0 <\rho \leqslant \infty$. Assuming $\left\vert \!
\left\vert f \right\vert \! \right\vert_{ 0, 0} \leqslant {\sf c}$, we
estimate:
\[
\aligned
\left\vert\!\left\vert
f
\right\vert\!\right\vert_{\widehat{0,\alpha/2}}
\leqslant
\sup_{{\sf x}''\neq{\sf x}'}\
\frac{\left\vert
f({\sf x}'')-f({\sf x}')
\right\vert}{
\left\vert
{\sf x}''-{\sf x}'\right\vert^{\alpha/2}}
=
\max\ ({\sf A},\,{\sf B})
\leqslant
{\sf A}
+
{\sf B},
\endaligned
\]
where ${\sf A} := \sup_{0< \left\vert {\sf x}'' - {\sf x}' \right\vert
< {\sf c}^{ 1/\alpha}}$ and ${\sf B} := \sup_{ \left\vert {\sf x}'' -
{\sf x}' \right\vert \geqslant {\sf c}^{ 1/ \alpha}}$ satisfy:
\[
\aligned
{\sf A}
&
=
\sup_{0<\left\vert 
{\sf x}''-{\sf x}'\right\vert
<{\sf c}^{1/\alpha}}
\left(
\frac{\left\vert 
f({\sf x}'')-f({\sf x}')
\right\vert}{
\left\vert
{\sf x}''-{\sf x}'\right\vert^\alpha}\,
\left\vert
{\sf x}''-{\sf x}'
\right\vert^{\alpha/2}
\right)
\leqslant
\left\vert\!\left\vert
f
\right\vert\!\right\vert_{\widehat{0,\alpha}}\,
{\sf c}^{1/2},
\\
{\sf B}
&
=
\sup_{\left\vert 
{\sf x}''-{\sf x}'
\right\vert\geqslant
{\sf c}^{1/\alpha}}\,
\frac{\left\vert 
f({\sf x}'')-f({\sf x}')
\right\vert}{
\left\vert
{\sf x}''-{\sf x}'\right\vert^{\alpha/2}}
\leqslant
\frac{2
\left\vert\!\left\vert
f
\right\vert\!\right\vert_{0,0}}{{\sf c}^{1/2}}
\leqslant
2\,{\sf c}^{1/2}.
\endaligned
\]
Lemma~3.9 is proved.
\qed

\smallskip

Applying this to ${\sf x} = ( u, \theta)$, we
deduce:
\[
\aligned
\left\vert\!\left\vert
\Phi_u\vert_s
\right\vert\!\right\vert_{\widehat{0,\alpha/2}}
&
\leqslant
{\sf c}_2^{1/2}
\left[
2
+
\left\vert\!\left\vert
\Phi_u\vert_s
\right\vert\!\right\vert_{\widehat{0,\alpha}}
\right]
\leqslant
{\sf c}_2^{1/2}
\left[
2
+
\left\vert\!\left\vert
\Phi\vert_s
\right\vert\!\right\vert_{1,\alpha}
\right],
\\
\left\vert\!\left\vert
\Phi_\theta\vert_s
\right\vert\!\right\vert_{\widehat{0,\alpha/2}}
&
\leqslant
{\sf c}_3^{1/2}
\left[
2
+
\left\vert\!\left\vert
\Phi_\theta\vert_s
\right\vert\!\right\vert_{\widehat{0,\alpha}}
\right]
\leqslant
{\sf c}_3^{1/2}
\left[
2
+
\left\vert\!\left\vert
\Phi\vert_s
\right\vert\!\right\vert_{1,\alpha}
\right].
\endaligned
\]
Consequently:
\[
\aligned
\left\vert\!\left\vert
\Phi\vert_s
\right\vert\!\right\vert_{1,\alpha/2}
&
=
\left\vert\!\left\vert
\Phi\vert_s
\right\vert\!\right\vert_{0,0}
+
\left\vert\!\left\vert
\Phi_u\vert_s
\right\vert\!\right\vert_{0,0}
+
\left\vert\!\left\vert
\Phi_\theta\vert_s
\right\vert\!\right\vert_{0,0}
+
\left\vert\!\left\vert
\Phi_u\vert_s
\right\vert\!\right\vert_{\widehat{0,\alpha/2}}
+
\left\vert\!\left\vert
\Phi_\theta\vert_s
\right\vert\!\right\vert_{\widehat{0,\alpha/2}}
\\
&
\leqslant
{\sf c}_1
+
{\sf c}_2
+
{\sf c}_3
+
\left(
{\sf c}_2^{1/2}
+
{\sf c}_3^{1/2}
\right)
\left[
2
+
\left\vert\!\left\vert
\Phi
\right\vert\!\right\vert_{1,\alpha}
\right].
\endaligned
\]
Lemma~3.10 is proved.
\qed

\subsection*{4.2.~Proof of Lemma~3.11}
We shall abbreviate $\sup_{ 0 < \left\vert \theta '' - \theta '
\right\vert \leqslant \pi}$ by $\sup_{ \theta '' \neq \theta '}$. By
definition:
\begin{small}
\[
\aligned
&
\left\vert \! \left\vert \Phi \big( U(\cdot), \cdot, s
\big) \right\vert \! \right\vert_{ \mathcal{ C}^{ 1, \beta} (\partial
\Delta)}
=
\\
&
=
\sup_\theta\,
\left\vert
\Phi\big(
U(e^{i\,\theta}),e^{i\,\theta},s
\big)
\right\vert
+
\\
& \ \ \
+
\sup_\theta
\Big\vert
\sum_{1\leqslant l\leqslant d}\,
\Phi_{u_l}
\big(
U(e^{i\,\theta}),e^{i\,\theta},s
\big)
U_\theta^l(e^{i\,\theta})
+
\Phi_\theta
\big(
U(e^{i\,\theta}),e^{i\,\theta},s
\big)
\Big\vert
+
\\
& \ \ \
+ 
\sup_{\theta''\neq\theta'}
\Big\vert
\sum_{1\leqslant l\leqslant d}\,
\Phi_{u_l}
\big(
U(e^{i\,\theta''}),e^{i\,\theta''},s
\big)
U_\theta^l(e^{i\,\theta''})
+
\Phi_\theta
\big(
U(e^{i\,\theta''}),e^{i\,\theta''},s
\big)
-
\\
&
\ \ \ 
-
\sum_{1\leqslant l\leqslant d}\,
\Phi_{u_l}
\big(
U(e^{i\,\theta'}),e^{i\,\theta'},s
\big)
U_\theta^l(e^{i\,\theta'})
-
\Phi_\theta
\big(
U(e^{i\,\theta'}),e^{i\,\theta'},s
\big)
\Big\vert
\left\vert
\theta''-\theta'
\right\vert^{-\beta}
\endaligned
\]
\end{small}

\noindent
Majorating and inserting some appropriate new terms whose sum is zero:
\begin{small}
\[
\aligned
\leqslant
&
\left\vert\!\left\vert
\Phi
\right\vert\!\right\vert_{0,0}
+
\left\vert\!\left\vert
\Phi_u
\right\vert\!\right\vert_{0,0}
\left\vert\!\left\vert
U
\right\vert\!\right\vert_{\widehat{1,0}}
+
\left\vert\!\left\vert
\Phi_\theta
\right\vert\!\right\vert_{0,0}
+
\\
&
+
\sup_{\theta''\neq\theta'}
\Big\vert
\sum_{1\leqslant l\leqslant d}\,
\Phi_{u_l}
\big(
U(e^{i\,\theta''}),e^{i\,\theta''},s
\big)
\big[
U_\theta^l(e^{i\,\theta''})
-
U_\theta^l(e^{i\,\theta'})
\big]
\Big\vert
\left\vert
\theta''-\theta'
\right\vert^{-\beta}
+
\\
&
+
\sup_{\theta''\neq\theta'}
\Big\vert
\sum_{1\leqslant l\leqslant d}\,
\Big[
\Phi_{u_l}
\big(
U(e^{i\,\theta''}),e^{i\,\theta''},s
\big)
-
\Phi_{u_l}
\big(
U(e^{i\,\theta''}),e^{i\,\theta'},s
\big)
\Big]
U_\theta^l(e^{i\,\theta'})
\Big\vert
\left\vert
\theta''-\theta'
\right\vert^{-\beta}
+
\\
&
+
\sup_{\theta''\neq\theta'}
\Big\vert
\sum_{1\leqslant l\leqslant d}\,
\Big[
\Phi_{u_l}
\big(
U(e^{i\,\theta''}),e^{i\,\theta'},s
\big)
-
\Phi_{u_l}
\big(
U(e^{i\,\theta'}),e^{i\,\theta'},s
\big)
\Big]
U_\theta^l(e^{i\,\theta'})
\Big\vert
\left\vert
\theta''-\theta'
\right\vert^{-\beta}
+
\\
&
+
\sup_{\theta''\neq\theta'}
\Big\vert
\Phi_\theta
\big(
U(e^{i\,\theta''}),
e^{i\,\theta''},s
\big)
-
\Phi_\theta
\big(
U(e^{i\,\theta''}),
e^{i\,\theta'},s
\big)
\Big\vert
\left\vert
\theta''-\theta'
\right\vert^{-\beta}
+
\\
&
+\sup_{\theta''\neq\theta'}
\Big\vert
\Phi_\theta
\big(
U(e^{i\,\theta''}),
e^{i\,\theta'},s
\big)
-
\Phi_\theta
\big(
U(e^{i\,\theta'}),
e^{i\,\theta'},s
\big)
\Big\vert
\left\vert
\theta''-\theta'
\right\vert^{-\beta}.
\endaligned
\]
\end{small}

\noindent
Majorating:
\begin{small}
\[
\aligned
\leqslant
&
\left\vert\!\left\vert
\Phi
\right\vert\!\right\vert_{0,0}
+
\left\vert\!\left\vert
\Phi_u
\right\vert\!\right\vert_{0,0}
\left\vert\!\left\vert
U
\right\vert\!\right\vert_{\widehat{1,0}}
+
\left\vert\!\left\vert
\Phi_\theta
\right\vert\!\right\vert_{0,0}
+
\\
&
+
\left\vert\!\left\vert
\Phi_u
\right\vert\!\right\vert_{0,0}
\left\vert\!\left\vert
U
\right\vert\!\right\vert_{\widehat{1,\beta}}
+
\left\vert\!\left\vert
\Phi_u\vert_s
\right\vert\!\right\vert_{\widehat{0,\beta}}
\left\vert\!\left\vert
U
\right\vert\!\right\vert_{\widehat{1,0}}
+
\left\vert\!\left\vert
\Phi_u\vert_s
\right\vert\!\right\vert_{\widehat{0,\beta}}
\big(
\left\vert\!\left\vert
U
\right\vert\!\right\vert_{\widehat{1,0}}
\big)^\beta
\left\vert\!\left\vert
U
\right\vert\!\right\vert_{\widehat{1,0}}
+
\\
&
+
\left\vert\!\left\vert
\Phi_\theta\vert_s
\right\vert\!\right\vert_{\widehat{0,\beta}}
+
\left\vert\!\left\vert
\Phi_\theta\vert_s
\right\vert\!\right\vert_{\widehat{0,\beta}}
\big(
\left\vert\!\left\vert
U
\right\vert\!\right\vert_{\widehat{1,0}}
\big)^\beta,
\endaligned
\]
\end{small}

\noindent
which yields the lemma, noticing that $\left\vert \! \left\vert \Phi_u
\right\vert \! \right\vert_{ 0,0} \big( \left\vert \! \left\vert U
\right\vert \! \right\vert_{ \widehat{1,0}} + \left\vert \! \left\vert
U \right\vert \! \right\vert_{ \widehat{1,\beta}} \big) \leqslant
\left\vert \! \left\vert \Phi_u \right\vert \! \right\vert_{ 0,0}
\left\vert \! \left\vert U \right\vert \! \right\vert_{1, \beta}$.
\qed

\subsection*{4.3.~Proof of Lemma~3.13}
We need two preparatory lemmas.

\def\thelemma{4.4}\begin{lemma}
Let $n\geqslant 1$, let ${\sf x} \in \R^n$, let $m\geqslant 1$, let ${\sf y} \in
\R^m$, let $\rho >0$ and let $f = f ({\sf x}, {\sf y})$ be a $\in
\mathcal{ C}^{1, \alpha}$ map, with $0 < \alpha \leqslant 1$, defined in
the strip $\{ ({\sf x}, {\sf y}) \in \R^m \times \R^n : \left \vert
{\sf x} + {\sf y} \right \vert < \rho\}$ and valued in $\R^d$, 
$d\geqslant 1$. If four vertices
$({\sf x}', {\sf y}')$, $({\sf x}'', {\sf y}')$, $({\sf x}', {\sf
y}'')$ and $({\sf x}'', {\sf y}'')$ of a rectangle belong to the
strip, then{\rm :} 
\[
\left\vert
f({\sf x}'', {\sf y}'')
-
f({\sf x}', {\sf y}'')
-
f({\sf x}'', {\sf y}')
+
f({\sf x}', {\sf y}')
\right\vert
\leqslant
\left\vert\!\left\vert
f
\right\vert\!\right\vert_{\widehat{1,\alpha}}
\left\vert
{\sf x}''
-
{\sf x}'
\right\vert
\left\vert
{\sf y}''
-
{\sf y}'
\right\vert^\alpha.
\]
A similar inequality holds by reversing the r\^oles of the
variables ${\sf x}$ and ${\sf y}$.
\end{lemma}

\proof
We apply twice the Taylor integral formula~\thetag{ 1.2}
with respect to the variable ${\sf x}$ and we majorate:
\[
\aligned
&
\left\vert
\left(
f({\sf x}'',{\sf y}'')
-
f({\sf x}',{\sf y}'')
\right)
-
\left(
f({\sf x}'',{\sf y}')
-
f({\sf x}',{\sf y}')
\right)
\right\vert
\leqslant
\\
& 
\ \ \
\leqslant
\int_0^1
\sum_{1\leqslant i\leqslant n}
\left\vert
\frac{\partial f}{\partial {\sf x}_i}
\left(
{\sf x}'+{\sf s}({\sf x}''-{\sf x}'),{\sf y}''
\right)
-
\frac{\partial f}{\partial {\sf x}_i}
\left(
{\sf x}'+{\sf s}({\sf x}''-{\sf x}'),{\sf y}'
\right)
\right\vert
\left\vert
{\sf x}_i''-{\sf x}_i'
\right\vert
d{\sf s}
\\
& 
\ \ \
\leqslant
\left\vert\!\left\vert
f
\right\vert\!\right\vert_{\widehat{1,\alpha}}
\left\vert
{\sf y}''-{\sf y}'
\right\vert^\alpha
\left\vert
{\sf x}''-{\sf x}'
\right\vert. \qed
\endaligned
\]

\def\thelemma{4.5}\begin{lemma}
Let $n\geqslant 1$, let ${\sf x} \in \R^n$, let $\rho >0$ and let $H = H
({\sf t})$ be a $\mathcal{ C}^{1, \alpha}$ map, with $0 < \alpha \leqslant
1$, defined in the open cube $\square_\rho^n = \{ \vert {\sf t} \vert
< \rho \}$ and valued in $\R^d$. Let ${\sf x}$, ${\sf y}$, ${\sf z}$
with $\vert {\sf x} \vert, \vert {\sf y} \vert, \vert {\sf z} \vert
\leqslant \rho/3$, so that the four points ${\sf x}$, ${\sf y}$, ${\sf z}$
and ${\sf x} + {\sf y} - {\sf z}$ constitute the vertices of a
parallelogram which is contained in $\square_\rho^n$. Then{\rm :}
\[
\left\vert
H({\sf x} + {\sf y} - {\sf z})
-
H({\sf x})
-
H({\sf y})
+
H({\sf z})
\right\vert
\leqslant
4
\left\vert\!\left\vert
H
\right\vert\!\right\vert_{\widehat{1,\alpha}}
\left\vert
{\sf y}-{\sf z}
\right\vert
\left\vert
{\sf x}-{\sf z}
\right\vert^\alpha.
\]
A similar inequality holds after exchanging ${\sf x}$ and
${\sf y}$.
\end{lemma}

\proof 
To estimate the second difference of $H$, we introduce a new map
\[
f({\sf x},{\sf y}) 
:= 
H({\sf x}+{\sf y}),
\]
of $({\sf x}, {\sf y}) \in \R^n \times \R^n$, whose domain is the
strip $\{ \vert {\sf x} + {\sf y} \vert < \rho\}$. Let ${\sf x}$,
${\sf y}$, ${\sf z}$ with $\vert {\sf x} \vert, \vert {\sf y} \vert,
\vert {\sf z} \vert \leqslant \rho/3$. Fixing ${\sf x}' \in \R^n$
arbitrary, there exist unique ${\sf y}'$, ${\sf x}''$ and ${\sf y}''$
solving the linear system:
\[
\left\{
\aligned
{\sf x}' + {\sf y}' 
&
= 
{\sf z}, 
\ \ \ \ \ \ \ \ \ \ \ 
{\sf x}'' +
{\sf y}''
= {\sf x} + {\sf y} - {\sf z},
\\
{\sf x}' + {\sf y}'' 
&
= 
{\sf x}, 
\ \ \ \ \ \ \ \ \ \ \ 
{\sf x}'' + {\sf y}'
= 
{\sf y}.
\endaligned\right.
\]
In fact, ${\sf y}' = {\sf z} - {\sf x}'$, ${\sf x}'' = {\sf y} - {\sf
z} + {\sf x}'$ and ${\sf y}'' = {\sf x} - {\sf x}''$. Taking the norms
$\vert \cdot \vert$ of the four equations above, we see that the
rectangle $({\sf x}', {\sf y}')$, $({\sf x}'', {\sf y}')$, $({\sf x}',
{\sf y}'')$, $({\sf x}'', {\sf y}'')$ is contained in the strip $\{
\vert {\sf x} + {\sf y} \vert < \rho\}$. Applying then Lemma~4.4
(with $m = n$), we get:
\[
\aligned
&
\left\vert
H({\sf x} + {\sf y} - {\sf z})
-
H({\sf x})
-
H({\sf y})
+
H({\sf z})
\right\vert
=
\\
&
\ \ \ \ \
=
\left\vert
f({\sf x}'',{\sf y}'')
-
f({\sf x}',{\sf y}'')
-
f({\sf x}'',{\sf y}')
+
f({\sf x}',{\sf y}')
\right\vert
\\
&
\ \ \ \ \
\leqslant
\left\vert\!\left\vert
f
\right\vert\!\right\vert_{\widehat{1,\alpha}}
\left\vert
{\sf x}''
-
{\sf x}'
\right\vert
\left\vert
{\sf y}''
-
{\sf y}'
\right\vert^\alpha
\\
&
\ \ \ \ \
=
\left\vert\!\left\vert
f
\right\vert\!\right\vert_{\widehat{1,\alpha}}
\left\vert
{\sf y}
-
{\sf z}
\right\vert
\left\vert
{\sf x}
-
{\sf z}
\right\vert^\alpha.
\endaligned
\]
We claim that $\left\vert \! \left\vert f \right\vert \! \right\vert_{
\widehat{1, \alpha }} \leqslant 4 \left\vert \! \left\vert H \right\vert \!
\right\vert_{ \widehat{1, \alpha }}$, which will 
conclude. Carefulness and rigor are
required. In fact, to estimate:
\[
\small
\aligned
\left\vert\!\left\vert
f
\right\vert\!\right\vert_{ \widehat{1, \alpha}}
=
\sum_{i=1}^n\,
\sup_{
({\sf x}'',{\sf y}'')
\neq
({\sf x}',{\sf y}')}\,
\frac{
\left\vert
f_{{\sf x}_i}({\sf x}'',{\sf y}'')
-
f_{{\sf x}_i}({\sf x}',{\sf y}')
\right\vert
+
\left\vert
f_{{\sf y}_i}({\sf x}'',{\sf y}'')
-
f_{{\sf y}_i}({\sf x}',{\sf y}')
\right\vert
}{
\left\vert
({\sf x}'',{\sf y}'')
-
({\sf x}',{\sf y}')
\right\vert^\alpha
},
\endaligned
\]
we shall first transform the denominator. 
By definition:
\[
\left\vert
({\sf x}'',{\sf y}'')
- 
({\sf x}',{\sf y}') 
\right\vert 
= 
\max\left(
\left\vert 
{\sf x} '' 
- 
{\sf x}' 
\right\vert, 
\left\vert 
{\sf y} ''
- 
{\sf y}' 
\right\vert\right). 
\]
If we set $a := \left \vert {\sf x}'' - {\sf x}' \right \vert$ and $b
:= \left \vert {\sf y}'' - {\sf y} ' \right \vert$ and if we invert
the inequality $(a+b)^\alpha \leqslant 2\, \max\left( a^\alpha, b^\alpha
\right)$, noticing $2^\alpha \leqslant 2$, we obtain:
\[
\small
\aligned
\frac{1}{\left\vert({\sf x}'',{\sf y}'')
-({\sf x}',{\sf y}')\right\vert^\alpha}
=
\frac{1}{
\max\left(
a^\alpha, b^\alpha
\right)
}
\leqslant
\frac{2}{
(a+b)^\alpha
}
&
=
\frac{2}{
\left(
\left\vert
{\sf x}''-{\sf x}'
\right\vert
+
\left\vert
{\sf y}''-{\sf y}'
\right\vert
\right)^\alpha
}
\\
&
\leqslant
\frac{2}{
\left\vert
{\sf x}''+{\sf y}''
-
{\sf x}'-{\sf y}'
\right\vert^\alpha
}.
\endaligned
\]
Replacing the denominator above, we then transform 
and majorate the numerator: 
\[
\small
\aligned
\left\vert\!\left\vert
f
\right\vert\!\right\vert_{\widehat{1,\alpha}}
&
\leqslant
2\,
\sum_{i=1}^n\,
\sup_{
({\sf x}'',{\sf y}'')
\neq
({\sf x}',{\sf y}')}
\left(
\frac{
\left\vert
H_{{\sf t}_i}({\sf x}''+{\sf y}'')
-
H_{{\sf t}_i}({\sf x}'+{\sf y}')
\right\vert
+
\left\vert
H_{{\sf t}_i}({\sf x}''+{\sf y}'')
-
H_{{\sf t}_i}({\sf x}'+{\sf y}')
\right\vert
}{
\left\vert
{\sf x}''+{\sf y}''
-
{\sf x}'-{\sf y}'
\right\vert^\alpha
}
\right)
\\
&
\leqslant
4\,\sum_{i=1}^n\,
\sup_{{\sf t}''\neq{\sf t}'}
\left(
\frac{
\left\vert 
H_{{\sf t}_i}({\sf t''})
-
H_{{\sf t}_i}({\sf t'})
\right\vert}{
\left\vert
{\sf t}''-{\sf t}'
\right\vert^\alpha
}
\right)
\\
&
=
4
\left\vert\!\left\vert
H
\right\vert\!\right\vert_{ \widehat{1, \alpha}}.
\endaligned
\]
This completes the proof of Lemma~4.5.
\endproof

We can now state a slightly simplified version of Lemma~3.13.

\def\thelemma{4.6}\begin{lemma}
{\rm (\cite{ tu1990}, [$*$])} Let $u \in \R^d$, $d\geqslant 1$, let $\rho_1
>0$ and let $\Psi = \Psi (u)$ be a $\mathcal{ C}^{ 1, \alpha}$ map,
with $0 < \alpha \leqslant 1$, $u\in \R^d$, defined in the cube $\{ \vert u
\vert <\rho_1 \}$ and valued in $\R^d$. Let $U^1, U^2 \in \mathcal{
C}^{ 1, 0} (\partial \Delta, \R^d)$ with $\left\vert U^j (e^{ i\,
\theta}) \right\vert < \rho_1 /3$ on $\partial \Delta$, for $j=1,
2$. For every $\beta$ with $0 < \beta \leqslant \alpha$ the following
inequality holds{\rm :}
\[
\left\vert\!\left\vert
\Psi(U^2(\cdot))
-
\Psi(U^1(\cdot))
\right\vert\!\right\vert_{\mathcal{C}^{0,\beta}(\partial\Delta)}
\leqslant
{\sf D}\,
\left\vert\!\left\vert
U^2
- 
U^1
\right\vert\!\right\vert_{0,\beta},
\]
with
\[
{\sf D}
=
\left\vert\!\left\vert
\Psi
\right\vert\!\right\vert_{1,\beta}
\left[
1
+
2
\left(
\left\vert\!\left\vert
U^1
\right\vert\!\right\vert_{\widehat{1,0}}
\right)^\beta
+
2
\left(
\left\vert\!\left\vert
U^2
\right\vert\!\right\vert_{\widehat{1,0}}
\right)^\beta
\right].
\]
\end{lemma}

\proof
Firstly and obviously:
\[
\left\vert\!\left\vert
\Psi(U^2)
-
\Psi(U^1)
\right\vert\!\right\vert_{0,0}
\leqslant
\left\vert\!\left\vert
\Psi
\right\vert\!\right\vert_{\widehat{1,0}}
\left\vert\!\left\vert
U^2
-
U^1
\right\vert\!\right\vert_{0,0}.
\]
Secondly, we have $\left\vert\! \left\vert \Psi(U^2) - \Psi(U^1)
\right\vert \! \right\vert_{ \widehat{0, \beta}} = \sup_{ 0<
\left\vert \theta'' - \theta' \right\vert \leqslant\pi} \big( {\sf Q}/
\left\vert \theta''- \theta' \right \vert^\beta \big)$, where:
\[
{\sf Q}
:=
\big\vert
\Psi
\big(
U^2(e^{i\,\theta''})
\big)
-
\Psi
\big(
U^1(e^{i\,\theta''})
\big)
-
\Psi
\big(
U^2(e^{i\,\theta'})
\big)
+
\Psi
\big(
U^1(e^{i\,\theta'})
\big)
\big\vert.
\]
To majorate ${\sf Q}$, we start by inserting the term $\Psi \big[
U^1(e^{i\,\theta''}) + U^2(e^{i\,\theta'}) - U^1(e^{i\, \theta'})
\big]$, well-defined, thanks to the assumption $\left\vert \!
\left\vert U^j \right\vert \! \right\vert_{ 0, 0} < \rho_1 /3$:
\[
\aligned
{\sf Q}
\leqslant
&
\big\vert
\Psi
\big(
U^2(e^{i\,\theta''})
\big)
-
\Psi\big[
U^1(e^{i\,\theta''}) 
+ 
U^2(e^{i\,\theta'}) 
- 
U^1(e^{i\, \theta'})
\big]
\big\vert
+
\\
&
+
\big\vert
\Psi\big[
U^1(e^{i\,\theta''}) 
+ 
U^2(e^{i\,\theta'}) 
- 
U^1(e^{i\, \theta'})
\big]
-
\Psi\big(
U^1(e^{i\,\theta''}) 
\big)
-
\\
&
\ \ \ \ \ \ \ \ \ \ \ \ 
\ \ \ \ \ \ \ \ \ \ \ \ 
\ \ \ \ \ \ \ \ \ \ \ \ 
\ \ \ \ \
-
\Psi\big(
U^2(e^{i\,\theta'}) 
\big)
+
\Psi\big(
U^1(e^{i\,\theta'}) 
\big)
\big\vert.
\endaligned
\]
To estimate the second absolute value, we apply Lemma~4.5 with ${\sf
x} = U^1 (e^{ i\, \theta ''})$, with ${\sf y} = U^2 ( e^{ i\,
\theta'})$ and with ${\sf z} = U^1 ( e^{ i\, \theta'})$:
\[
\aligned
{\sf Q}
\leqslant
&
\left\vert\!\left\vert
\Psi
\right\vert\!\right\vert_{\widehat{1,0}}
\big\vert
[U^2-U^1](e^{i\,\theta''})
-
[U^2-U^1](e^{i\,\theta'})
\big\vert
+
\\
&
+
4\,
\left\vert\!\left\vert
\Psi
\right\vert\!\right\vert_{\widehat{1,\beta}}
\big\vert
U^2(e^{i\,\theta'})
-
U^1(e^{i\,\theta'})
\big\vert
\big\vert
U^1(e^{i\,\theta''})
-
U^1(e^{i\,\theta'})
\big\vert^\beta.
\endaligned
\]
We then achieve the remaining majorations:
\[
\aligned
{\sf Q}
\leqslant
&
\left\vert\!\left\vert
\Psi
\right\vert\!\right\vert_{\widehat{1,0}}
\left\vert\!\left\vert
U^2-U^1
\right\vert\!\right\vert_{\widehat{0,\beta}}
\left\vert
\theta''-\theta'
\right\vert^\beta
+
\\
&
+
4\,
\left\vert\!\left\vert
\Psi
\right\vert\!\right\vert_{\widehat{1,\beta}}
\left\vert\!\left\vert
U^2-U^1
\right\vert\!\right\vert_{0,0}
\left(
\left\vert\!\left\vert
U^1
\right\vert\!\right\vert_{\widehat{1,0}}
\right)^\beta
\left\vert
\theta''-\theta'
\right\vert^\beta.
\endaligned
\]
Reminding that $\left\vert\! \left\vert
U^2-U^1 \right\vert \! \right\vert_{ 0, 0} + \left\vert\! \left\vert
U^2-U^1 \right\vert \! \right\vert_{ \widehat{ 0,\beta}} =
\left\vert\! \left\vert U^2-U^1 \right\vert \! \right\vert_{
0,\beta}$ and summing, we obtain:
\[
\left\vert\!\left\vert
\Psi(U^2)
-
\Psi(U^1)
\right\vert\!\right\vert_{0,\beta}
\leqslant
\left\vert\!\left\vert
\Psi
\right\vert\!\right\vert_{1,\beta}
\left[
1
+
4
\left(
\left\vert\!\left\vert
U^1
\right\vert\!\right\vert_{\widehat{1,0}}
\right)^\beta
\right]
\left\vert\!\left\vert
U^2-U^1
\right\vert\!\right\vert_{0,\beta}.
\]
A similar inequality holds with $\big( \left\vert\! \left\vert U^2
\right\vert\! \right\vert_{\widehat{1,0}} \big)^\beta$ instead of
$\big( \left\vert\! \left\vert U^1 \right \vert\!
\right\vert_{\widehat{1,0}} \big)^\beta$. Taking the arithmetic mean,
we find the symmetric quantity ${\sf D}$ of the lemma. The proof is
complete.
\endproof

\proof[Proof of Lemma~3.13.]
By definition:
\begin{small}
\[
\aligned
{\sf R}
:=
&
\left\vert\!\left\vert
\Phi\left(U^2(\cdot),\cdot,s\right)
-
\Phi\left(U^1(\cdot),\cdot,s\right)
\right\vert\!\right\vert_{\mathcal{C}^{0,\beta}(\partial\Delta)}
\\ 
=
&
\sup_\theta\,
\left\vert
\Phi\big(U^2(e^{i\,\theta}),e^{i\,\theta},s\big)
-
\Phi\big(U^1(e^{i\,\theta}),e^{i\,\theta},s\big)
\right\vert
+
\\
&
+
\sup_{\theta''\neq\theta'}\,
\Big\vert
\Phi\big(U^2(e^{i\,\theta''}),e^{i\,\theta''},s\big)
-
\Phi\big(U^1(e^{i\,\theta''}),e^{i\,\theta''},s\big)
-
\\
&
\ \ \ \ \ \ 
\ \ \ \ \ \ 
\ \ \ \ \ \ 
-
\Phi\big(U^2(e^{i\,\theta'}),e^{i\,\theta'},s\big)
+
\Phi\big(U^1(e^{i\,\theta'}),e^{i\,\theta'},s\big)
\Big \vert
\left\vert
\theta''-\theta'
\right\vert^{-\beta}.
\endaligned
\]
\end{small}

\noindent
In the numerator, we insert $-\Phi\big( U^2 (e^{i\, \theta''}),
e^{i\,\theta'}, s\big) + \Phi\big (U^1( e^{i\, \theta''}),e^{i\,
\theta'},s \big)$ plus its opposite:
\begin{small}
\[
\aligned
{\sf R}
\leqslant
&
\left\vert\!\left\vert
\Phi
\right\vert\!\right\vert_{\widehat{1,0}}
\left\vert\!\left\vert
U^2-U^1
\right\vert\!\right\vert_{0,0}
+
\\
&
+
\sup_{\theta''\neq\theta'}\,
\Big\vert
\Phi\big(U^2(e^{i\,\theta''}),e^{i\,\theta''},s\big)
-
\Phi\big(U^1(e^{i\,\theta''}),e^{i\,\theta''},s\big)
-
\\
&
\ \ \ \ \ \ 
\ \ \ \ \ \ 
\ \ \ \ \ \ 
-
\Phi\big(U^2(e^{i\,\theta''}),e^{i\,\theta'},s\big)
+
\Phi\big(U^1(e^{i\,\theta''}),e^{i\,\theta'},s\big)
\Big\vert
\left\vert
\theta''-\theta'
\right\vert^{-\beta}
+
\\
&
+
\sup_{\theta''\neq\theta'}\,
\Big\vert
\Phi\big(U^2(e^{i\,\theta''}),e^{i\,\theta'},s\big)
-
\Phi\big(U^1(e^{i\,\theta''}),e^{i\,\theta'},s\big)
-
\\
&
\ \ \ \ \ \ 
\ \ \ \ \ \ 
\ \ \ \ \ \ 
-
\Phi\big(U^2(e^{i\,\theta'}),e^{i\,\theta'},s\big)
+
\Phi\big(U^1(e^{i\,\theta'}),e^{i\,\theta'},s\big)
\Big\vert
\left\vert
\theta''-\theta'
\right\vert^{-\beta}.
\endaligned
\]
\end{small}

\noindent
To majorate the first $\sup_{\theta '' \neq \theta '}$, we apply
Lemma~4.4 with ${\sf x}' = U^1 ( e^{ i\, \theta''})$, ${\sf y}' = e^{
i\, \theta'}$, ${\sf x}'' = U^2 ( e^{ i\, \theta''})$ and ${\sf y}'' =
e^{ i\, \theta'}$, where $s$ is considered as a dumb parameter. To
majorate the second $\sup_{\theta '' \neq \theta '}$, we apply
Lemma~4.6 to $\Psi ( u):= \Phi \big( u, e^{ i\, \theta}, s \big)$,
where $(e^{ i\, \theta}, s)$ are considered as dumb parameters. We
get:
\begin{small}
\[
\aligned
{\sf R}
\leqslant
&
\left\vert\!\left\vert
\Phi
\right\vert\!\right\vert_{\widehat{1,0}}
\left\vert\!\left\vert
U^2-U^1
\right\vert\!\right\vert_{0,0}
+
\\
&
+
\left\vert\!\left\vert
\Phi
\right\vert\!\right\vert_{\widehat{1,\beta}}
\sup_\theta\,
\left\vert
U^2(e^{i\,\theta})
-
U^1(e^{i\,\theta})
\right\vert
+
\\
&
+
\left\vert\!\left\vert
\Phi
\right\vert\!\right\vert_{1,\beta}
\left[
1
+
2\left(
\left\vert\!\left\vert
U^1
\right\vert\!\right\vert_{\widehat{1,0}}
\right)^\beta
+
2\left(
\left\vert\!\left\vert
U^2
\right\vert\!\right\vert_{\widehat{1,0}}
\right)^\beta
\right]
\left\vert\!\left\vert
U^2-U^1
\right\vert\!\right\vert_{0,\beta}.
\endaligned
\]
\end{small}

\noindent
To conclude, we use $\left\vert \! \left\vert \Phi \right\vert \!
\right\vert_{ \widehat{ 1, 0}} +\left\vert \! \left\vert \Phi
\right\vert \! \right\vert_{ \widehat{ 1, \beta}} \leqslant\left\vert \!
\left\vert \Phi \right\vert \! \right\vert_{1, \beta}$ and we get the
term ${\sf C}$ of Lemma~3.13.
\endproof

With these techniques, the proofs of Lemmas~3.17, 3.19, 3.30, 3.31
and~3.32 are easily guessed and even simpler.

\newpage

\begin{center}
{\Large\bf V:~Holomorphic extension of CR functions}
\end{center}

\bigskip\bigskip\bigskip

\begin{center}
\begin{minipage}[t]{11cm}
\baselineskip =0.35cm
{\scriptsize

\centerline{\bf Table of contents}

\smallskip

{\bf 1.~Hartogs theorem, jump formula and domains having the
holomorphic extension property \dotfill 153.}

{\bf 2.~Tr\'epreau's theorem, deformations of Bishop discs and
propagation on hypersurfaces \dotfill 166.}

{\bf 3.~Tumanov's theorem, deformations of
Bishop discs and propagation on generic submanifolds \dotfill 178.}

{\bf 4.~Holomorphic extension on globally minimal generic 
submanifolds \dotfill 192.}

\smallskip

\hfill 
{\footnotesize\tt [19 diagrams]}

}\end{minipage}
\end{center}

\bigskip
\bigskip

{\small


The method of analytic discs is rooted in the very birth of the theory
of functions of several complex variables. The discovery by Hurwitz
and Hartogs of the compulsory extension of holomorphic functions
relied upon an application of Cauchy's integral formula along a family
of analytic discs surrounding an illusory singularity. Since
H.~Cartan, Thullen, Behnke and Sommer, various versions of this
argument were coined ``{\sl Kontinuit\"attsatz}'' or 
``{\sl Continuity principle}''.
 
The removal of compact singularities culminated in the so-called {\sl
Hartogs-Bochner theorem}, usually proved by means of integral formulas
or thanks to the resolution of a $\overline{ \partial}$
problem with compact support. Contradicting 
all expectations, a subtle example
due to Forn{ae}ss (1998), shows that on a non-pseudoconvex domain, the
disc method may fail to fill in the domain, if the discs are
required to stay inside the domain.

Nevertheless, it is of the highest prize to build constructive methods
in order to describe significant parts of the envelope of holomorphy
of a domain, of a CR manifold, as well as the polynomial hull of
certain compact sets. In such problems, analytic discs with boundary
in the domain, the CR manifold or the
compact set remain the most adequate tools.

The precise existence Theorem~3.7(IV) 
for the solutions of Bishop's equation
that was established in the previous Part~IV may now be applied
systematically to a variety of geometric situations. In this respect,
we just followed Bishop's genuine philosophy that required to ensure
an explicit control of the size of solutions in terms of the size of
data, instead of appealing to some general, imprecise version of the
implicit function theorem.

Thanks to the jump theorem, holomorphic extension of CR functions
defined on a hypersurface $M$ is equivalent to the extension of the
functions that are holomorphic in one of the two sides to the other
side. Tr\'epreau's original theorem (1986) states that such an
extension holds at a point $p$ if and only if there does not exist a
local complex hypersurface $\Sigma$ of $\C^n$ with $p\in \Sigma
\subset M$. A deeper phenomenon of propagation (Tr\'epreau, 1990)
holds: if CR functions extend holomorphically to one side at a point
$p$, a similar extension holds at every point of the CR orbit of $p$
in $M$. By means of deformations of attached Bishop discs, there is an
elementary (and folklore) proof that contains both the local and the
global extension theorems on hypersurfaces, yielding a satisfactory
understanding of the phenomenon.

On a generic submanifold $M$ of $\C^n$ of higher codimension, the
celebrated {\sl Tumanov extension theorem} (1988) states that CR
functions defined on $M$ extend holomorphically to a local wedge of
edge $M$ at a point $p$ if the local CR orbit of $p$ contains a
neighborhood of $p$ in $M$. A globalization of this statement,
obtained independently by J\"oricke and the first author in 1994, states
that the same extension phenomenon holds if $M$ consists of a single
CR orbit, {\it i.e.} is {\sl globally minimal}. Both proofs heavily
relied on the local Tumanov theorem and on a precise control of the
propagation of directions of extension.
 
A clever proof that treats both locally minimal and globally minimal
generic submanifolds on the same footing constitutes the main
Theorem~4.12 of the present Part~V: {\it if $M$ is a globally minimal
$\mathcal{ C}^{ 2, \alpha}$ $(0 < \alpha < 1)$ generic submanifold of
$\C^n$ of codimension $\geqslant 1$ and of CR dimension $\geqslant 1$,
there exists a wedgelike domain $\mathcal{ W}$ attached to $M$ such
that every continuous CR function $f\in \mathcal{ C}_{ CR}^0 (M)$
possesses a holomorphic extension $F \in \mathcal{ O} (\mathcal{ W})
\cap \mathcal{ C}^0 (M\cup \mathcal{ W}) $ with $F\vert_M = f$}. This
basic statement as well as the techniques underlying its proof will be
the very starting point of the study of removable singularities in
Parts~VI and in~\cite{ mp2006a}.

}

\section*{\S1.~Hartogs theorem, jump formula \\
and domains having the extension property}

\subsection*{1.1.~Hartogs extension theorem: brief
history} \footnote{Further historical information may be found
in~\cite{ ra1986, str1988, fi1991, ra2002, hm2002}.} 
In 1897, Hurwitz showed
that a function holomorphic in $\C^2 \backslash \{ 0\}$ extends
holomorphically through the origin. In his thesis (1906), Hartogs
generalized this discovery, emphasizing the compulsory holomorphic
extendability of functions that are defined on the nowadays celebrated
{\sl Hartogs skeleton} (diagram below). The main argument is to apply
Cauchy's integral formula along families of analytic discs having
their boundary inside the domain and whose interior goes outside the
domain. In fact, the thinness of an embedded circle in $\C^n$ ($n
\geqslant 2$) offers much freedom to include illusory singularities
inside a disc.

In 1924, Osgood stated the ultimate generalization of the discovery of
Hurwitz and Hartogs: {\it if $\Omega \subset \C^n$ $(n\geqslant 2)$ is a
domain and if $K \subset \Omega$ is any compact such that $\Omega
\backslash K$ {\sl connected}, then $\mathcal{ O} (\Omega \backslash
K) = \mathcal{ O} (\Omega) \vert_{ \Omega \backslash K}$}. This
statement is nowadays called the {\sl Hartogs-Bochner} theorem.
Although the proof of Osgood was correct for geometrically simple
complements $\Omega \backslash K$, as for instance spherical shells,
it was incomplete for general $\Omega \backslash K$. In fact,
unpleasant topological and monodromy obstructions occur for general
$\Omega \backslash K$ when pushing analytic discs. In 1998, Forn{ae}ss
exhibited certain domains in which discs are forced to first
leave some intermediate domain $\Omega \backslash K_1$, $K_1 \subset
K$, before $\Omega$ may be filled in.

In the late 1930's, a rigorous proof of Osgood's general statement was
obtained by Fueter, by means of a generalization of the classical
Cauchy-Green-Pompeiu integral formula to several variables, in the
context of complex and quaternionic functions (Moisil 1931, Fueter
1935). In 1943, Martinelli simplified the formal treatment of
Fueter. Then Bochner observed that the same result holds more
generally if one assumes given on $\partial \Omega$ just a CR
function.

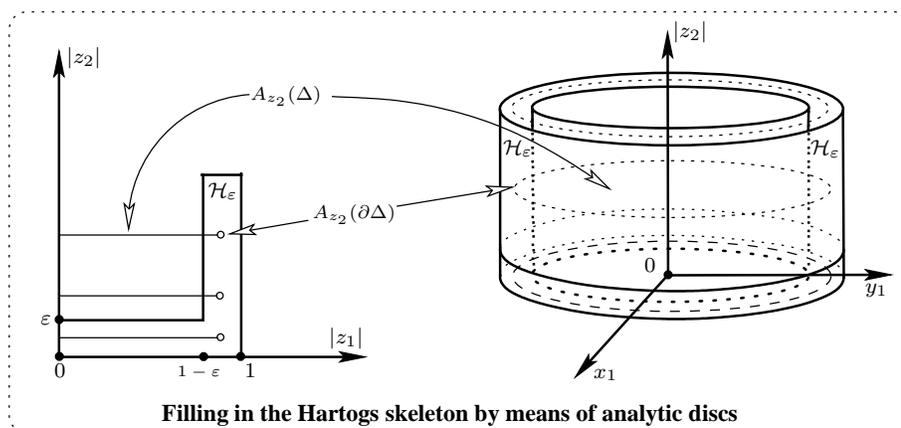
\begin{center}
\input pot.pstex_t
\end{center}

\subsection*{1.2.~Hartogs domain}
Consider the {\sl $\varepsilon$-Hartogs skeleton} (pot-looking)
domain:
\[
\mathcal{H}_\varepsilon
:=
\{
(z_1,z_2)\in\C^2:
\vert z_1\vert<1,
\vert z_2\vert<\varepsilon
\}
\, \bigcup\
\{
1-\varepsilon<\vert z_1\vert<1,
\vert z_2\vert<1
\}.
\]
We draw two diagrams: in $(\vert z_1 \vert, \vert z_2 \vert)$ and in
$(x_1,y_1, \vert z_2 \vert )$ coordinates.

\def\thelemma{1.3}\begin{lemma}
Every holomorphic function $f\in \mathcal{ O} ( \mathcal{
H}_\varepsilon)$ extends holomorphically to the bidisc $\Delta \times
\Delta$, the convex hull of $\mathcal{ H}_\varepsilon$.
\end{lemma}

\proof
Letting $\delta$ with $0 < \delta < \varepsilon$, for every $z_2 \in
\C$ with $\vert z_2 \vert < 1$, the analytic disc 
\[
\Delta\ni\zeta
\longmapsto
A_{z_2} (\zeta)
:=
\left(
[1-\delta]\zeta,z_2
\right)
\in\C^2
\]
has its boundary $A_{ z_2} (\partial \Delta)$ contained in $\mathcal{
H}_\varepsilon$, the domain where the function $f$ is defined. Thus,
we may compute the Cauchy integral
\[
F(z_1,z_2)
:=
\frac{1}{2\pi i}\,
\int_{\partial\Delta}\,
\frac{f(A_{z_2}(\zeta))}{\zeta-z_1}\,d\zeta.
\]
Differentiating under the sum, this extension $F$ is seen to be
holomorphic. In addition, for $\vert z_2 \vert < \varepsilon$, it
coincides with $f$. Obviously, the discs $A_{ z_2} (\Delta)$ fill 
in the hole of the domain $\mathcal{ H }_\varepsilon$.
\endproof

\subsection*{1.4.~Bounded domains in $\C^n$ and Hartogs-Bochner 
extension phenomenon} Let $\Omega$ be a connected open subset of
$\C^n$, a {\sl domain}. We assume it to be {\sl bounded}, {\it i.e.}
$\overline{ \Omega}$ is compact and that its {\sl boundary} $\partial
\Omega := \overline{ \Omega} \backslash \Omega$ is a hypersurface of
$\C^n$ of class at least $\mathcal{ C}^1$. By means of a partition of
unity, one can construct a real-valued function $r$ defined on $\C^n$
such that $\Omega = \{ z : r(z) < 0\}$ and $\partial \Omega = \{ z : r
(z) = 0\}$, with $dr(z) \neq 0$ for every $z \in \partial
\Omega$. Then $\partial \Omega$ is orientable.

Extensions of the above disc argument led to the most
general\footnote{ Often, some authors consider instead a compact $K
\subset \Omega$ with 
$\Omega \backslash K$
connected and state that $\mathcal{ O} (\Omega \backslash K) =
\mathcal{ O} (\Omega) \big\vert_{ \Omega \backslash K}$; a technical
check shows that the two statements are equivalent.} form of the
Hartogs theorem: {\it if $\Omega$ is a bounded domain in $\C^n$
($n\geqslant 2$) having connected boundary $\partial \Omega$, then
every function holomorphic in a neighborhood of $\partial \Omega$
uniquely extends as a function holomorphic in $\Omega$}. There are
three classical methods of proof:

\begin{itemize}

\smallskip\item[$\bullet$]
using the Bochner-Martinelli kernel;

\smallskip\item[$\bullet$]
using solutions of $\overline{ \partial} u = v$
having compact support;

\smallskip\item[$\bullet$]
pushing analytic discs, in successive
Hartogs skeletons.

\end{itemize}

The first two are rigorously established and we shall review the first
in a while. For almost one hundred years, it has been a folklore belief
that the third method could be accomplished somehow.
Let us be precise.

\subsection*{1.5.~Forn{ae}ss' counterexample and a disc theorem} 
Thus, let $\Omega$ be a bounded domain of $\C^2$ having connected
$\mathcal{ C }^1$ boundary. For $\delta >0$ small, consider the
one-sided neighborhood of $\partial \Omega$ defined by:
\[
\widetilde{\Omega}_\delta
:=
\{ z \in \Omega 
:{\rm dist}\, 
(z, \partial \Omega)
< \delta \}.
\]
The complement $\Omega \backslash \widetilde{ \Omega}_\delta$ is a
compact hole. Remind that the bidisc $\Delta^2$ is the convex hull of
the Hartogs skeleton $\mathcal{ H }_\varepsilon$.
Following 
\cite{ f1998}, we say that {\sl $\Omega$ can be filled in by
analytic discs} if for every $\delta >0$, there exist a finite
sequence of subdomains of $\Omega$ having $\mathcal{ C}^1$ boundary,
$\widetilde{ \Omega}_\delta = \Omega_1 \subset \Omega_2 \subset \cdots
\subset \Omega_k = \Omega$ and for each $j = 1, \dots, k-1$, an
$\varepsilon_j >0$ and a univalent holomorphic map $\Phi_j$ defined in
a neighborhood of $\overline{ \Delta }^2$ such that:

\begin{itemize}

\smallskip\item[{\bf (1)}]
$\Omega_{j+1} \subset \Omega_j \cup
\Phi_j ( \Delta^2) \subset \Omega$;

\smallskip\item[{\bf (2)}]
$\Phi_j(\mathcal{ H}_\varepsilon )\subset\Omega_j$;

\smallskip\item[{\bf (3)}]
$\Omega_j \cap \Phi_j ( \Delta^2)$ is connected;

\smallskip\item[{\bf (4)}]
$\Omega_{j+1}\cap\Phi_j(\Delta^2)$ is connected.
\end{itemize}\smallskip 

For such domains, by pushing analytic discs in the embedded Hartogs
figure, taking account of connectedness, we have $\mathcal{ O}
(\Omega_{ j+1}) \big\vert_{ \Omega_j} = \mathcal{ O} (\Omega_j)$. Then
by induction, uniquely determined holomorphic extension holds from
$\Omega_1$ up to $\Omega$. Importantly, the intermediate domains are
required to be all contained in $\Omega$.

In 1998, Forn{ae}ss \cite{ f1998} constructed a topologically
strange domain $\Omega \subset \C^2$ that cannot be filled in this
way. This example shows that the requirement that $\Omega_j \subset
\Omega_{ j + 1} \subset \Omega$ is too stringent.

Nevertheless, taking care of monodromy and working in the envelope of
holomorphy of $\Omega$, one may push analytic discs by allowing them
to wander in the outside, in order to get the general Hartogs theorem
stated above. As a preliminary, one perturbs and smoothes out the
boundary. Denote by $\vert \! \vert z \vert \! \vert := \big( \vert
z_1 \vert^2 + \cdots + \vert z_n \vert^2 \big)^{ 1/2}$ the Euclidean
norm of $z = (z_1, \dots, z_n) \in \C^n$ and by $\B^n ( p, \delta ) :=
\big\{ \vert \! \vert z - p \vert \! \vert < \delta \big\}$ the open
ball of radius $\delta >0$ centered at a point $p$.

\def\thetheorem{1.6}\begin{theorem}
{\rm (\cite{ mp2006c})}
Let $M \Subset \C^n$ {\rm (}$n \geqslant 2${\rm )} be a {\rm
connected} $\mathcal{ C }^\infty$ hypersurface bounding a domain
$\Omega_M \Subset \C^n$. Suppose to fix ideas that $2 \leqslant {\rm
dist }\, \big( 0 , \overline{ \Omega }_M \big) \leqslant 5$ and assume
that the restriction $r_M := r\vert_M$ of the distance function $r (z)
= \vert \! \vert z \vert \! \vert$ to $M$ is a Morse function having
only a finite number $\kappa$ of critical points $\widehat{ p}_\lambda
\in M$, $1 \leqslant \lambda \leqslant \kappa$, located on different
sphere levels{\rm :}
\[
2\leqslant 
\widehat{r}_1
:=
r(\widehat{p}_1)
<
\cdots
<
\widehat{r}_\kappa
:=
r(\widehat{p}_\kappa)
\leqslant
5
+
\text{\rm diam}
\big(\overline{\Omega}_M\big).
\]
Then there exists $\delta_1 >0$ such that for every $\delta$ with $0 <
\delta < \delta_1$, the {\rm (}tubular{\rm )} neighborhood 
\[
\mathcal{
V }_\delta (M)
:=
\cup_{p\in M}\, 
\B^n(p,\delta)
\] 
enjoys the global Hartogs extension property into 
$\Omega_M${\rm
:}
\[
\mathcal{O}
\big( 
\mathcal{V}_\delta( 
M)
\big)
=
\mathcal{O}
\big(
\Omega_M\cup
\mathcal{V}_\delta( 
M)
\big)
\big\vert_{\mathcal{V}_\delta(M)},
\]
by $''$pushing$''$ analytic discs inside
a finite number of Hartogs figures, without using neither the
Bochner-Martinelli kernel, nor solutions of some
auxiliary $\overline{ \partial }$ problem.
\end{theorem}

\subsection*{1.7.~Hartogs-Bochner theorem via the
Bochner-Martinelli kernel} By $\mathcal{ O} (C)$, where $C \subset
\C^n$ is closed, we mean $\mathcal{ O} (\mathcal{ V} (C))$ for some
open neighborhood $\mathcal{ V} (C)$ of $C$. Here is the general
statement.

\def\thetheorem{1.8}\begin{theorem}
{\rm (\cite{ hele1984, he1985, ra1986})}
Let $\Omega$ be a bounded domain in $\C^n$ having {\rm connected}
boundary. Then for every neighborhood $U$ of $\partial \Omega$ in
$\C^n$ and every holomorphic function $f\in \mathcal{ O} (U)$, there
exists a function $F \in \mathcal{ O} ( \overline{ \Omega})$ with $F
\vert_{\partial \Omega} = f \vert_{ \partial \Omega}$.
\end{theorem}

In the thin neighborhood $U$ of the not necessarily smooth boundary
$\partial \Omega$, by means of a partion of unity, one may construct a
connected boundary $\partial \Omega_1 \subset U$ close to $\partial
\Omega$ which is $\mathcal{ C}^1$, or $\mathcal{ C}^\infty$, or even
$\mathcal{ C}^\omega$, using Whitney approximation (\cite{ hi1976}; in
addition, one may assure that $r(z) \big\vert_{ \partial \Omega_1}$ is
as in Theorem~1.6, whence both statements are equivalent). Then the
restriction $F \vert_{ \partial \Omega_1}$ is CR on $\partial
\Omega_1$ and the previous theorem is a consequence of the next.

\def\thetheorem{1.9}\begin{theorem}
{\rm (\cite{ ra1986, he1985})} Let $\Omega$ be a bounded domain in
$\C^n$ $(n\geqslant 2)$ having connected $\mathcal{ C}^{ \kappa,
\alpha}$ boundary, with $1\leqslant \kappa \leqslant \infty$, $0
\leqslant \alpha \leqslant 1$. Then for every CR function $f :
\partial \Omega \to \C$ of class $\mathcal{ C}^{\kappa, \alpha}$,
there exists a function $F \in \mathcal{ O} ( \Omega ) \cap \,
\mathcal{ C}^{\kappa, \alpha} (\overline{ \Omega})$ with $F
\vert_{\partial \Omega} = f$.
\end{theorem}

\noindent
Some words about the proof. With $\zeta, z \in \C^n$, consider the
{\sl Bochner-Martinelli} kernel:
\[
{\sf BM}(\zeta,z)
:=
\frac{(n-1)!}{(2\pi i)^n}\,
\vert 
\zeta-z
\vert^{-2n}\,
\sum_{j=1}^n\,
\left(
\bar\zeta_j
-
\bar z_j
\right)
d\zeta_j\,
\underset{k\neq j}{\wedge}\,
d\bar\zeta_k\,
\wedge\,
d\zeta_k.
\]
This is a $(n, n-1)$-form which is 
$\mathcal{ C}^\omega$ off the diagonal $\{
\zeta = z\}$. For $n=2$, it coincides with the Cauchy kernel $\frac{
1}{ 2\pi i}\, \frac{ 1}{ \zeta - z}$. If $f$ and $\partial \Omega$ are
$\mathcal{ C}^1$, the integral formula:
\[
F(z)
:=
\int_{\partial \Omega}\,
f(\zeta)\,{\sf BM}(\zeta,z)
\]
provides the holomorphic extension $F$.

\subsection*{1.10.~Hypersurfaces of $\C^n$ and
jump theorem for CR functions} Let $M$ be a real hypersurface of
$\C^n$ without boundary. In the sequel, we shall mainly deal with
three geometric situations.

\begin{itemize}

\smallskip\item[$\bullet$]
{\bf Local:} $M$ is defined in a small open polydisc
centered at one point $p \in M$.

\smallskip\item[$\bullet$]
{\bf Global:} $M$ is a connected
orientable embedded submanifold of $\C^n$.

\smallskip\item[$\bullet$]
{\bf Boundary:} $\C^n \backslash M$ consists of
two open sets $\Omega^+$, bounded and 
$\Omega^-$, unbounded.

\end{itemize}\smallskip

Then there exists some appropriate neighborhood $\mathcal{ M}$ of $M$
in $\C^n$ in which $M$ is {\sl relatively closed}, in the sense that
$\overline{ M} \cap \mathcal{ M} = M$.

More generally, let $\mathcal{ M}$ be an arbitrary complex manifold of
dimension $n \geqslant 1$ and let $M \subset \mathcal{ M}$ be a
hypersurface of class at least $\mathcal{ C}^1$ which is 
relatively closed in $\mathcal{ M}$ and oriented.
The complement $\mathcal{ M}
\backslash M$ then consists of two connected components $\Omega^+$ and
$\Omega^-$, where $\Omega^+$ is located on the positive side to
$M$. Also, let $f : M \to \C$ be a CR function of class at least
$\mathcal{ C}^0$. By definition, $f$ is CR if the current of
integration on $M$ of bidegree $(0,1)$ defined 
by\footnote{ 
Here, $\mathcal{ D}^{ p, q}$ is the space of $\mathcal{ C}^\infty$
forms of bidegree $(p,q)$ having compact support; fundamental notions
about currents may be found in~\cite{ ch1989}.
}:
\[
f_M(\omega)
:=
\int_M\,f\,\omega,
\ \ \ \ \ \ \ \
\omega\in\mathcal{D}^{n,n-1},
\]
satisfies $\int_M \, f\,
\overline{\partial} \varpi = 0$ for every $\varpi \in \mathcal{ D}^{
n, n-2}$. Equivalently, $\overline{ \partial} f_M = 0$ in the sense of
currents, where $f_M$ is interpreted as a $(0,1)$-form having measure
coefficients.

To formulate the jump theorem in arbitrary complex manifolds, we shall
mainly assume that the Dolbeault $\overline{ \partial}$-complex on
$\mathcal{ M}$ is exact in bidegree $(0,1)$, namely $H_{\overline{
\partial}}^{ 0, 1} ( \mathcal{ M}) =0$. This assumption holds for
instance when $\mathcal{ M} = \Delta^n$, $\C^n$ or $P_n (\C)$. It
means that the equation $\overline{ \partial} u = v$, where $v$ is a
$\overline{ \partial}$-closed $(0,1)$-form on $\mathcal{ M}$ having
$\mathcal{ C}^\infty$, $L^2$ or distributional coefficients has a
$\mathcal{ C}^\infty$, $L^2$ or distributional solution $u$ on
$\mathcal{ M}$.

Consequently, there exists a distribution $F$ on $\mathcal{ M}$ with
$\overline{ \partial} F = f_M$. As ${\rm supp} \, f_M \subset M$, such
a function $F$ is holomorphic in $\mathcal{ M} \backslash M$. The
difference $F_2 - F_1$ of two solutions to $\overline{ \partial} F =
f_M$ is holomorphic in $\mathcal{ M}$. In the case where $\mathcal{ M}
= \C^n$, a solution to $\overline{ \partial} F = f_M$ may be
represented (\cite{ ch1975, ra1986}) by means of the
Bochner-Martinelli kernel as $F(z) := \int_M\,f(\zeta)\, {\sf
BM}(\zeta,z)$. In complex dimension $n=1$, such a solution coincides
with the classical Cauchy transform.
 
In 1975, after previous work of Andreotti-Hill (\cite{ ah1972b}),
Chirka obtained a several complex variables version of the
Sokhotski${\check{ \text{\i }}}$-Plemelj Theorem~2.7(IV).

\def\thetheorem{1.11}\begin{theorem} 
{\rm (\cite{ ch1975})} 
Assume that $H_{ \overline{ \partial }}^{ 0,1} ( \mathcal{ M}) = 0$
and that the hypersurface $M \subset \mathcal{ M}$ is orientable and
relatively closed, {\it i.e.} $\overline{ M} \cap \mathcal{ M} = M$.
Assume $\dim \mathcal{ M} = n\geqslant 1$ and let $(\kappa, \alpha)$
with $0 \leqslant \kappa \leqslant \infty$, $0 < \alpha < 1$. If $M$
is $\mathcal{ C}^{ \kappa+ 1, \alpha}$ and if the current $f_M$
associated to a $\mathcal{ C}^{ \kappa, \alpha}$ function $f : M \to
\C$ is CR, then every distributional solution $F \in \mathcal{ O}
(\mathcal{ M} \backslash M)$ to $\overline{ \partial } F = f_M$
extends to be $\mathcal{ C}^{ \kappa, \alpha}$ in the two closures
$\overline{ \Omega^\pm} = \Omega^\pm \cup M$, yielding two functions
$F^\pm \in \mathcal{ O} ( \Omega^\pm) \cap \mathcal{ C}^{ \kappa,
\alpha} (\Omega^\pm \cup M)$ whose jump across $M$ equals $f${\rm :}
\[
F^+(z)
-
F^-(z)
=
f(z),
\ \ \ \ \ \ \
\forall\,z\in M.
\]
A similar jump formula holds for $f \in L_{ loc, CR}^{ \sf p} (M)$,
with $M$ at least $\mathcal{ C}^1$ {\rm (}or a Lipschitz graph{\rm )}
and for $f\in \mathcal{ D}_{ CR}' (M)$, with $M \in \mathcal{ C
}^\infty$.
\end{theorem}

When $\mathcal{ M} = \C$, the conditions that $f$ is CR and that $H^{
0, 1}_{ \overline{ \partial }} (\mathcal{ M}) = 0$ are automatically
satisfied and we recover the Sokhotski${\check{ \text{\i }}}$-Plemelj
jump formula. However, we mention that in several complex variables
($n\geqslant 2$), there is no analog of the second formula $\frac{1 }{
2} \left[ F^ + (\zeta_0) + F^- (\zeta_0) \right] = {\rm p.v.}\,
\frac{1}{ 2\pi i} \int_\Gamma\, \frac{f( \zeta)}{ \zeta- \zeta_0}\,
d\zeta$. The reason is the inexistence of a universal integral formula
solving $\overline{ \partial} F = f_M$. Nevertheless, there should
exist generalized principal value integrals which depend on the
kernel.

If $M$ is only $\mathcal{ C}^1$ and $f$ is only $\mathcal{ C}^0$, it
is in general untrue that $F^-$ and $F^+$ extend continuously to
$M$. Fortunately, there is a useful substitute result, analog to
Theorem~2.9(IV). Consider a open subset $M' \subset M$ having compact
closure $\overline{ M'}$ not meeting $\partial M = \overline{ M}
\backslash M$. We may embedd $M'$ in a one-parameter family $\left(
M_\varepsilon' \right)_{ \vert \varepsilon \vert < \varepsilon_0 }$,
$\varepsilon_0 >0$, of hypersurfaces that foliates a
strip thickening of $M'$.

\def\thetheorem{1.12}\begin{theorem}
{\rm (\cite{ ch1975})} If $f$ is CR and $\mathcal{ C}^\kappa$ on
a $\mathcal{ C}^{ \kappa + 1}$ hypersurface $M$, then
\[
\lim_{\varepsilon\to 0} \big\vert \! \big\vert F\vert_{
M_\varepsilon'} - F\vert_{ M_{- \varepsilon }'} - f \big\vert \!
\big\vert_\kappa =0.
\]
\end{theorem}

\subsection*{1.13.~CR extension in the projective space}
Unlike in $\C^n$, there is no privileged ``interior''
side of an orientable connected hypersurface in the projective space
$P_n (\C)$, $n\geqslant 2$. Nevertheless, a version of the Hartogs-Bochner
theorem holds. The proof is an illustration of the use of the jump
theorem.

\def\thetheorem{1.14}\begin{theorem}
{\rm ($n\geqslant 3$: \cite{ hl1975}; $n=2$: \cite{ sa1999,
dm2002})} Let $M$ be a compact orientable connected $\mathcal{C}^2$
real hypersurface of $P_n (\mathbb{C})$ that divides the projective
space into two domains $\Omega^-$ and $\Omega^+$. Then{\rm :}
\begin{itemize}

\smallskip\item[{\bf (i)}]
there exists a side, $\Omega^-$ \text{\rm or} $\Omega^+$, to which
every function holomorphic in some neighborhood of $M$ extends
holomorphically\footnote{ 
Using propagation techniques of Section~3, the theorem holds assuming
that $M$ is globally minimal and considering continuous CR functions
on $M$.
}{\rm ;}

\smallskip\item[{\bf (ii)}] 
every function that is holomorphic in the union of the other side of
$M$ together with a neighborhood of $M$ must be constant.
\end{itemize}\smallskip
\end{theorem}

Let us summarize the proof. Let $f$ be holomorphic in some
neighborhood $\mathcal{ V}(M)$ of $M$ in $P_n (\C)$. As the Dolbeault
cohomology group $H^{0,1}(P_n(\mathbb{C}))$ vanishes for $n\geqslant 2$
(\cite{hele1984, he1985}), thanks to Theorem~1.11 above, the CR
function $f \vert_{ M}$ on $M$ can be decomposed as the jump $f = F^+
- F^-$ between two functions $F^\pm$ holomorphic in $\Omega^\pm$ which
are (at least) continuous up to $M$. It suffices then to show that
either $F^+$ or $F^-$ is constant, since clearly, if $F^+$
(resp. $F^-$) is constant equal to $c^+$ (resp. $c^-$), then $f$
extends 
holomorphically to $\Omega^-$ 
(resp. to $\Omega^+$) as $c^+ - F^-$ (resp. as $F^+ - c^-$).

By contradiction, assume that both $F^+$ and $F^-$ are nonconstant. We
choose two domains $U^+$ and $U^-$ with $\mathcal{ V} (M) \cup
\Omega^\pm \supset U^\pm \supset M \cup \Omega^\pm$. By a preliminary
(technical) deformation argument, we may assume that $F^\pm$ is
holomorphic in $U^\pm$. According to a theorem due to Takeuchi
\cite{ta1964}, holomorphic functions in an arbitrary domain of $P_n
(\C)$ ($n\geqslant 2$), either are constant or separate points. Since
$F^-$ is nonconstant, $\mathcal{ O} (U^-)$ separates
points. Conjugating with elements of the group ${\rm PGL} (n, \C)$ of
projective automorphisms of $P_n (\C)$, shrinking $\mathcal{ V} (M)$
and $U^-$ slightly if necessary, we may verify (\cite{ dm2002}) that
$\mathcal{ O} ( U^-)$ separates points and provides local system of
holomorphic coordinates at every point. By standard techniques of
Stein theory (\cite{ ho1973}), it follows that $U^-$ is embeddable in
some $\C^N$, with $N$ large. The image of $M$ under such an embedding
$\Phi$ is a compact CR submanifold of $\C^N$ that is
filled by the relatively compact complex manifold $\Sigma^- = \Phi
(U^-)$ with boundary $\Phi (M)$. Two applications of the maximum 
principle to
the nonconstant holomorphic function $F^+ \circ \Phi^{ -1}$ say that
it must decrease inside $\Sigma^-$, since $\Sigma^-$ is interior to
$\Phi (M)$ in $\C^N$, and that it must increase, since the one-sided
neighborhood $U^- \cap \mathcal{ V} (M)$is exterior to $U^+$.
This is the desired contradiction.

\subsection*{1.15.~Levi extension theorem}
A $\mathcal{ C}^2$ hypersurface $M \subset \C^n$ may always be
represented as $M = \{ z \in U : r(z) = 0\}$, where $U$ is some open
neighborhood of $M$ in $\C^n$, and where $r: U \to \R$ is a $\mathcal{
C}^2$ implicit {\sl defining function} that satisfies $d r (q ) \neq
0$ at every point $q$ of $M$. Two defining functions $r^1, r^2$ are
nonzero multiple of each other in some neighborhood $V \subset U$ of
$M$: there exists $\lambda : V \to \R$ nowhere vanishing with $r_2 =
\lambda \, r_1$. 

At a point $p\in M$, the {\sl Levi form of} $r$:
\[
\mathfrak{L}_p\,r(L_p,\overline{L}_p)
:=
\sum_{1\leqslant j,k\leqslant n}\,
\frac{\partial^2 r}{\partial z_j\partial \bar z_k}(p)\,
L_p^j\overline{L}_p^k,
\ \ \ \ \ \ \ \ \ \ 
L_p\in T_p^{1,0}M,
\]
is a Hermitian form that may be diagonalized. 
Its {\sl signature at} $p$: 
\[
(a_p,b_p)
:=
\left(
\#\,
\text{\rm
positive eigenvalues},\
\#\,
\text{\rm
negative eigenvalues}
\right)
\]
is the same for $r_1$ and $r_2$ if they are positive multiples of each
other. It is also invariant through local biholomorphic changes of
coordinates $z \mapsto \widetilde{ z} (z)$ that do not reverse the
orientation of $M$. Reversing the orientation or
taking a negative factor $\lambda$ corresponds to the
transposition $(a_p, b_p) \mapsto (b_p, a_p)$.

The Levi form may be read off a graphed equation $v=\varphi (x,y,u)$
for $M$.

\def\thelemma{1.16}\begin{lemma}
There exist local holomorphic coordinates $(z,u+iv) \in \C^{ n-1}
\times \C$ centered at $p$ in which $M$ is represented as a
graph of the form{\rm :}
\[
v 
=
\varphi(z,u)
=
\sum_{1\leqslant k\leqslant a_p+b_p}\, 
\varepsilon_k\, 
z_k\bar z_k
+
{\rm o}(\vert z\vert^2)
+
{\rm O}(\vert z\vert\,\vert u\vert)
+
{\rm O}(\vert u\vert^2),
\]
where $\varepsilon_k = +1$ for $1\leqslant k\leqslant a_p$,
$\varepsilon_k = -1$ for $a_p+1\leqslant k\leqslant a_p+b_p$. If $a_p
= n-1$, the open set $\{ v > \varphi \}$ is strongly convex {\rm (}in
the real sense{\rm )} in a neighborhood of $p$.
\end{lemma}

Assuming that $M$ is orientable, it is surrounded by two open
sides. By an {\sl open side of} $M$, we mean a connected component of
$\mathcal{ V} \backslash M$ for a (thin) neighborhood $\mathcal{ V}$
of $M$ which is divided by $M$ in two components. As germs of open
sets along $M$, there exist two open sides (if $M$ were not
orientable, there would exist only one).

Assuming that $M$ is represented either by an implicit equation $r=0$
or as a graph $-v + \varphi (x,y,u) = 0$, we adopt the convention of
denoting:
\[
\aligned
\Omega^+
&
:=
\{r<0\}
\ \ \ \ \ 
{\rm or}
\ \ \ \ \
\Omega^+
:=
\{v>\varphi(x,y,u)\},
\\
\Omega^-
&
:=
\{r>0\}
\ \ \ \ \ 
{\rm or}
\ \ \ \ \
\Omega^-
:=
\{v<\varphi(x,y,u)\}.
\endaligned
\]
Once a local side $\Omega$ of $M$ has been fixed, $M$ is oriented and
the indetermination $r \leftrightarrow -r$ disappears. By convention,
we will always represent $\Omega = \{ r < 0\}$. Then the number of
positive and of negative eigenvalues of the Levi-form of $r$ at a
point $p\in \partial \Omega$ is an invariant. By common abuse of
language, we speak of the {\sl Levi form} of $\partial \Omega$.

At one of its points $p$, a boundary $\partial \Omega$ is called {\sl
strongly pseudoconvex} (resp. {\sl strongly pseudoconcave}) if its
Levi form has all its eigenvalues $>0$ (resp. $<0$) at $p$. It is
called {\sl weakly pseudoconvex} (resp. {\sl weakly pseudoconcave}) at
$p$ if all eigenvalues are $\geqslant 0$ (resp. $\leqslant 0$).
Often, the term ``weakly'' is dropped in common use.

\def\thedefinition{1.17}\begin{definition}
{\rm
If $\Omega$ is an open side of $M$, we say that 
$\Omega$ is {\sl holomorphically extendable} at $p$ if for
every open\footnote{ Although the shape of polydiscs is not invariant
by local biholomorphisms, their topology is. To avoid dealing
implicitly with possibly
wild open sets, we prefer to speak of neighborhoods
$U_p$, $V_p$, $W_p$ of $p$ that are polydiscs. } polydisc $U_p$
centered at $p$, there exists an open polydisc $V_p$ centered at $p$
such that for every $f \in \mathcal{ O} (\Omega \cap U_p)$, there
exists $F \in \mathcal{ O} (V_p)$ with $F\vert_{ \Omega \cap V_p} = f
\vert_{ \Omega \cap V_p}$.
}\end{definition}

In 1910, Levi localized the Hartogs extension phenomenon.

\def\thetheorem{1.18}\begin{theorem}
{\rm (\cite{ bo1991, trp1996, tu1998, ber1999})}
If the Levi form of $M \subset \partial \Omega$ has one negative
eigenvalue at a point $p$, then $\Omega$
is holomorphically extendable at $p$. 
\end{theorem}

\proof
As $\Omega$ is given by $\{ v > - z_1 \bar z_1 + \cdots\}$, 
for $\varepsilon >0$ small, the
disc $A_\varepsilon (\zeta) := (\varepsilon\, \zeta, 0, \dots, 
0)$ has its boundary $A_\varepsilon (\partial \Delta)$
contained in $\Omega$ near $p$.

\def\thelemma{1.19}\begin{lemma}
Assume $M$ is $\mathcal{ C}^1$, let $p\in M$ and $\Omega$ be an
open side of $M$ at $p$. Suppose that for every open polydisc $U_p$
centered at $p$, there exists an analytic disc $A : \Delta \to U_p$
continuous in $\overline{ \Delta}$ with $A (0) = p$ and $A (\partial
\Delta) \subset \Omega$. Then $\Omega$ is
holomorphically extendable at $p$.
\end{lemma}

To draw $A (\Delta)$, decreasing by $1$ its dimension, we represent it
as a curve. The cusp illustrates the fact that $A (\Delta)$ is {\it
not}\, assumed to be embedded.

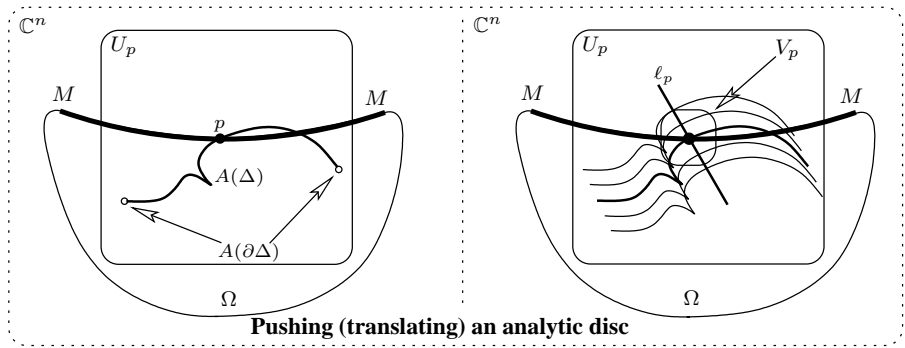
\begin{center}
\input pushing-disc.pstex_t
\end{center}

To prove the lemma, we may assume that $p=0$. Since $A (\partial
\Delta)$ is contained in $\Omega$, for $z\in \C^n$ very small, say
$\vert z \vert < \delta$, the translates $z + A (\partial \Delta)$ of
the disc boundary are also contained in $\Omega$.

Consequently, the Cauchy integral:
\[
F(z)
:= 
\frac{1}{2\pi i}\,
\int_{\partial\Delta}\,
f(z+A(\zeta))\,
\frac{d\zeta}{\zeta}
\]
is meaningful and it defines a holomorphic function of $z$ in the
polydisc $V_p := \{ \vert z \vert < \delta \}$.

Does it coincide with $f$ in $V_p \cap \Omega$\,? The assumption that
$M$ is $\mathcal{ C}^1$ yields a local real segment $\ell_p$
transversal to $M$ at $p$. If $U_p$ is sufficiently small and if $z
\in \ell_p \cap \Omega$ goes sufficiently deep in $\Omega$, the disc
$z + A (\overline{ \Delta})$ is contained in $\Omega$, so that the
Cauchy integral $F(z)$ coincides with $f (z)$ for those $z$.
\endproof

\subsection*{ 1.20.~Contact of weakly pseudoconvex domains with complex
hypersurfaces} The domain $\Omega$ is said to {\sl admit a support
complex hypersurface at} $p \in \partial \Omega$ if there
exists a local (possibly singular) complex hypersurface $\Sigma$
passing through $p$ that does not intersect $\Omega$. In this
situation, if $\Sigma = \{ h(z) = 0\}$ with $h$ holomorphic, the
function $1/h$ does not extend holomorphically at $p$, being
unbounded. With $\alpha > 0$ not integer, one may define branches
$h^\alpha$ which are uniform in $\Omega$ and continuous up to 
$\partial \Omega$, but
whose extension through $p$ would be ramified around $\Sigma$.
Consequently, the existence of a support complex hypersurface prevents
$\mathcal{ O} (\Omega)$ to be holomorphically extendable at $p$. Is it
the right obstruction\,? For instance, 
at a strongly pseudoconvex boundary point,
the complex tangent plane is support.

Nevertheless, in 1973, Kohn-Nirenberg constructed a special
pseudoconvex domain $\Omega_{\sf KN}^+$ in $\C^2$
showing that:

\begin{itemize}

\smallskip\item[$\bullet$]
not every weakly pseudoconvex smoothly bounded domain is locally
biholomorphically equivalent to a domain which is 
weakly convex in the real sense;

\smallskip\item[$\bullet$]
the holomorphic non-extendability of $\mathcal{ O} (\Omega)$ at $p$ is
totally independent from the existence of a local supporting complex
hypersurface at $p$.

\end{itemize}\smallskip 
The boundary of this domain
\[
M_{\sf KN}
:=
\big\{
(z,w)\in\C^2:
{\rm Im}\,w
=
\vert
zw
\vert^2
+
z^4\bar z^4
+
15\,
(z^7\bar z
+
\bar z^7z)/14
\big\}
\]
may be checked to be strongly pseudoconvex at every point except the
origin, where it is weakly pseudoconvex. Hence $\mathcal{ O} (
\Omega_{\sf KN}^+)$ is not holomorphically extendable at the
origin. However, $M_{\sf KN}$ has the striking property that every
local (possibly singular) complex hypersurface $\Sigma$ passing
through the origin meets both sides of $M_{\sf KN}$ in every
neighborhood of the origin. By means of a Puiseux parametrization
(\cite{ japf2000}), such a complex curve $\Sigma$ is always the image
of a certain holomorphic disc $\lambda : \Delta \to \C^2$ with
$\lambda (0) = 0$.

\def\thetheorem{1.21}\begin{theorem}
{\rm (\cite{ kn1973})} Whatever the holomorphic disc $\lambda$, for
every $\varepsilon >0$, there are $\zeta^-$ and $\zeta^+$ in $\Delta$
with $\vert \zeta^\pm \vert < \varepsilon$ such that $\lambda (
\zeta^\pm) \in \Omega_{\sf KN}^\pm$.
\end{theorem}

Clearly, $\Omega_{\sf KN}^+$ is not locally convexifiable at the
origin (otherwise, the biholomorphic image of the complex tangent line
would be support).

Often for technical reasons, certain results in several complex
variables require boundaries to be convex in the real
sense. Although this condition is {\it not} biholomorphically
invariant, it is certainly meaningful to characterize the class of
{\sl convexifiable} domains, at least locally: does there exist an
analytico-geometric criterion enabling to recognize local
convexifiability by reading a defining equation\,?


\subsection*{ 1.22.~Holomorphic extendability 
across finite type hypersurfaces} Let $M$ be a $\mathcal{ C}^\omega$
hypersurface and let $p\in M$. Then $M$ is of type $m$ at $p$, 
in the sense of Definition~4.22(III), if and
only if there exists a local graphed equation centered at $p$ of the
form:
\[
v
= 
\varphi_m(z,\bar z)
+
{\rm O}(\vert z\vert^{m+1})
+
{\rm O}(\vert z\vert\,\vert u\vert)
+
{\rm O}(\vert u\vert^2),
\]
where $\varphi_m \in \C [ z, \bar z]$ is a {\it nonzero}\, homogeneous
real-valued polynomial of degree $m$ having no pluriharmonic term,
namely $0 \equiv \varphi_m (0, \bar z) \equiv \varphi_m (z, 0)$. The
restriction $\varphi_m ( \ell (\zeta), \overline{ \ell (\zeta)})$ of
$\varphi_m$ to a complex line $\C \ni \zeta \mapsto \ell (\zeta) \in
\C^{ n-1}$, $\ell (0) = 0$, is a polynomial in $\C [ \zeta, \bar
\zeta]$. For almost every choice of $\ell$, this polynomial is
nonzero, homogeneous of the same degree $m$ and contains no harmonic
term. After a rotation, such a line is the complex
$z_1$-axis. Denoting $z' = (z_2, \dots, z_{ n-1})$, we obtain:
\def\theequation{1.23}\begin{equation}
v
=
\varphi_m(z_1,\bar z_1)
+
{\rm O}(\vert z_1\vert^{m+1})
+
O(\vert z'\vert)
+
{\rm O}(\vert z\vert\,\vert u\vert)
+
{\rm O}(\vert u\vert^2).
\end{equation}

\def\thetheorem{1.24}\begin{theorem}
{\rm (\cite{ befo1978, r1983, bt1984})}
If $m$ is even, at least one side $\Omega^+$ or $\Omega^-$
is holomorphically extendable at $p$. If $m$ is odd, 
both sides have this property.
\end{theorem}

\proof
Let $\varepsilon >0$ arbitrarily small, let $a\in \C$ with
$\vert a \vert < 1$, let $\zeta$ be in the closed unit disc
$\overline{ \Delta}$, and introduce a $\C^n$-valued analytic disc:
\[
A_\varepsilon(\zeta)
:=
\left(
\varepsilon(a+\zeta),
0,\dots,0,
\varepsilon^m\,
W(\zeta)
\right)
\]
having zero $z'$-component and $z_1$-component being a disc of radius
$\varepsilon$ centered at $-a$. We assume its $w$-component $W(
\zeta)$ be defined by requiring that the $\C^2$-valued disc
\[
B_\varepsilon (\zeta) := ( \varepsilon(a+ \zeta), \varepsilon^m\,
W(\zeta))
\]
has its boundary $B_\varepsilon (\partial \Delta)$ attached
to $v = \varphi_m (z_1, \bar z_1)$. By homogeneity, $\varepsilon^m$
drops and it is equivalent to require that $B_1$ is attached to $v =
\varphi_m (z_1, \bar z_1)$. Equivalently, the imaginary part
$V(\zeta)$ of $W = U + i \, V$ should satisfy:
\[
V(e^{i\,\theta})
=
\varphi_m(a+e^{i\,\theta},\bar a+e^{-i\,\theta}),
\]
for all $e^{i\, \theta} \in \partial \Delta$.
To obtain a harmonic extension to $\Delta$ of the function $V$ thus
defined on $\partial \Delta$, no Bishop equation is needed. It
suffices to take the harmonic prolongation by means of Poisson's
formula, as in \S2.20(IV):
\[
V(\eta)
=
{\sf P}V(\eta)
=
\frac{1}{2\pi i}\,
\int_{\partial\Delta}\,
\varphi_m(a+\zeta,\bar a+\bar\zeta)\,
\frac{1-\vert\eta\vert^2}{\vert\zeta-\eta\vert^2}\,
\frac{d\zeta}{\zeta}.
\]
Since $\varphi_m$ has no harmonic term, it may be factored as
$\varphi_m = z_1 \bar z_1 \, \psi_1 (z_1, \bar z_1)$, with $\psi_1 \in
\C[ z_1, \bar z_1]$ homogeneous of degree $(m-2)$ and nonzero. In the
integral above, we put $\eta := -a$ and we replace $\varphi_m = z_1
\bar z_1 \psi_1$ to get the value of $V$ at $-a$:
\[
\aligned
V(-a)
&
=
\frac{1}{2\pi i}\,
\int_{\partial\Delta}\,
\varphi_m(a+\zeta,\bar a+\bar\zeta)\,
\frac{1-\vert a\vert^2}{\vert\zeta+a\vert^2}\,\frac{d\zeta}{\zeta}
\\
&
=
\frac{1-\vert a\vert^2}{2\pi}\,
\int_{-\pi}^{\pi}\,
\psi_1(a+e^{i\,\theta},\bar a+e^{-i\,\theta})\,d\theta.
\endaligned
\]
As a function of $a\in \Delta$, the last integral is identically zero
if and only if the polynomial $\psi_1$ is zero. Thus, there exists $a$
such that $V( -a) \neq 0$. (However, we have no
information about the possible
signs of $V(-a)$ in terms of $\varphi_m$.) Then we define $U(\zeta)$
to be the harmonic prolongation of $-{\sf T} V$ that vanishes at $-a$
and $B_\varepsilon (\zeta) := ( \varepsilon (a + \zeta), \varepsilon^m
\, W (\zeta))$.

The positivity (resp. negativity) of the sign of $V(-a)$ means that
$B_\varepsilon (-a) = 
(0, i\, V(-a))$ is in $\Omega^+$ (resp. $\Omega^-$). Then after
translating slightly $B_\varepsilon$ in the right direction along the
$v$-axis, Lemma~1.19 applies to deduce that $\Omega_1^- \subset \C^2$
(resp. $\Omega_1^+ \subset \C^2$) is holomorphically extendable at the
origin. Since $\varphi_m (-z_1, - \bar z_1) = (-1)^m \, \varphi_m
(z_1, \bar z_1)$, in the case where $m$ is odd, the disc $-
B_\varepsilon$ will also be attached to $M_1$ and will provide
extendability of the other side.

Thanks to basic majorations of the ``${\rm O}$'' 
remainders in the equation~\thetag{ 1.23} of
$M$, if $\varepsilon >0$ is sufficiently small, then $\Omega^- \subset
\C^n$ (resp. $\Omega^+ \subset \C^n$) has the same extendability
property.
\endproof

If $M$ is a real analytic hypersurface, it is easily seen, by
inspecting the Taylor series of its graphing function, that $M$ is not
of finite type at a point $p$ if and only if it may be represented by
$v = u \, \widetilde{
\varphi } (z, u)$, with $\widetilde{
\varphi } \in \mathcal{ C}^\omega$.
Then the local complex hypersurface $\{ v = u = 0 \}$
is contained in $M$.

\def\thecorollary{1.25}\begin{corollary}
{\rm (\cite{ befo1978, r1983, bt1984})} 
If $M$ is $\mathcal{ C}^\omega$ and if $p\in M$, the
following properties are equivalent{\rm :}
\begin{itemize}

\smallskip\item[$\bullet$]
$M$ has finite type at $p${\rm ;}

\smallskip\item[$\bullet$]
$M$ does not contain any local complex analytic 
hypersurface passing through $p${\rm ;}

\smallskip\item[$\bullet$]
$\Omega^+$ or $\Omega^-$ 
is holomorphically extendable at $p$.

\end{itemize}\smallskip
\end{corollary}

\subsection*{ 1.26.~Which side is holomorphically extendable\,?} 
We claim that it suffices to study osculating domains
in $\C^2$ of the form:
\[
\Omega_{\varphi_m}^+ 
:= \{ -v + \varphi_m (z,
\bar z) < 0 \},
\ \ \ \ \ \ \
z\in\C,
\ \
w= u+iv\in\C,
\]
where $\varphi_m \neq 0$ is real, homogeneous of degree $m\geqslant 2$ and
has no harmonic term. Indeed, extendability properties of such
domains transfer to perturbations~\thetag{ 1.23}. Also, extendability
properties of $\Omega_{\varphi_m}^-$ are just the same, via $\varphi_m
\leftrightarrow -\varphi_m$. For this reason, if $m$ is odd, both
$\Omega_{ \varphi_m}^+$ and $\Omega_{ \varphi_m}^-$ are
holomorphically extendable at $p$.

The local complex line $\{ (z, 0)\}$ intersects the closure
$\overline{ \Omega_{\varphi_m}^+}$ in regions that are closed angular
sectors (cones), due to homogeneity. We call these regions {\sl
interior}. The complement $\C^2 \backslash \Omega_{ \varphi_m}^+$
intersects $\{ (z, 0)\}$ in open, {\sl exterior} sectors.

\def\thetheorem{1.27}\begin{theorem}
{\rm (\cite{ r1983, bt1984, fr1985})}
If there exists an {\rm interior} sector of angular width $>
\frac{ \pi}{ m}$, then $\Omega_{\varphi_m}^+$ 
is holomorphically extendable at $p$.
\end{theorem}

The proof consists in choosing an appropriate truncated angular sector
as the $z$-component of a disc attached to $\partial \Omega_{
\varphi_m}^+$, instead of the round disc $\zeta \mapsto \varepsilon
(a+ \zeta)$.

\def\theexample{1.28}\begin{example}
{\rm 
Every homogeneous quartic $v = \varphi_4 (z, \bar z)$ in $\C^2$
is biholomorphically equivalent to a model
\[
0
= 
r_a:=
-
v
+
z^2 \bar z^2
+ 
a\, z\bar z(z^2 + \bar z^2),
\]
for some $a\in \R$. Such a hypersurface bounds two open sides
$\Omega_a^+ = \{ r_a <0\}$ and $\Omega_a^- = \{ r_a >0\}$ which enjoy
the following properties (\cite{ r1983, bt1984}):

\begin{itemize}

\smallskip\item[$\bullet$]
$\Omega_a^-$ is holomorphically extendable at $p$, 
for every $a$;

\smallskip\item[$\bullet$]
$\vert a \vert < 2/3$ if and only if $\Omega_a^+$ is everywhere
strongly pseudoconvex;

\smallskip\item[$\bullet$]
$\vert a \vert \leqslant 1 \big/ \sqrt{ 2}$ if and only if
$\Omega_a^+$ is {\it not}\, holomorphically extendable at $p$;

\smallskip\item[$\bullet$]
$\vert a\vert > 1 \big/ \sqrt{ 2}$ if and only if $\Omega_a^+$ is
holomorphically extendable at $p$;

\smallskip\item[$\bullet$]
the above extendability property holds true for any perturbation of
$\partial \Omega$ by higher order terms.

\end{itemize}\smallskip

}\end{example}

Finer results, strictly more general than the above theorem that apply
to sixtics, were obtained in \cite{ fr1985}. If we remove all exterior
sectors of angular width $\geqslant \frac{ \pi}{ m}$, the rest of the
complex line $\{ (z, 0)\}$ is formed by disjoint closed sectors, which
are called {\sl supersectors of order $m$ of $\Omega_{ \varphi_m}^+$
at $p$}. A supersector is {\sl proper} if it contains points of
$\Omega_{ \varphi_m}^+$.

\def\thetheorem{1.29}\begin{theorem} 
{\rm (\cite{ fr1985})}

\begin{itemize}

\smallskip\item[{\bf (i)}]
If $\Omega_{ \varphi_m}^+$ has a proper supersector of angular width
$> \frac{ \pi}{ m}$, then $\Omega_{ \varphi_m}^+$ is
holomorphically extendable at $p$.

\smallskip\item[{\bf (ii)}]
If all supersectors of $\Omega_{ \varphi_m}^+$ have angular width $<
\frac{ \pi}{ m}$, then there exists $f \in \mathcal{ O} (\Omega_{
\varphi_m}^+) \cap \mathcal{ C}^0 ( \overline{ \Omega_{
\varphi_m}^+})$ that does not extend holomorphically at $p$.

\end{itemize}\smallskip

\end{theorem}

Even in the case $m=6$, some cases in this theorem are left
open. Examples may be found in~\cite{ fr1985}.

\def\theopenproblem{1.30}\begin{openproblem}
In the case where $m$ is even, find a necessary and sufficient
condition for $\Omega_{ \varphi_m}^+ = \{v > \varphi_m (z, \bar z)
\}$ to be holomorphically extendable at $p$, or show that
the problem is undecidable.
\end{openproblem}

One could generalize this (already wide open) question to a not
necessarily finite type boundary, $\mathcal{ C}^\omega$, $\mathcal{
C}^\infty$, $\mathcal{ C}^2$ or even $\mathcal{ C}^0$ graph.

\section*{ \S2.~Tr\'epreau's theorem, deformations of Bishop discs 
\\ and propagation on hypersurfaces} \

\subsection*{ 2.1.~Holomorphic extension of
CR functions via jump} Let $M$ be a hypersurface in $\C^n$ of class at
least $\mathcal{ C}^{ 1, \alpha}$ with $0 < \alpha < 1$ and let $f$ be
a continuous CR function on $M$. At each point $p$ of $M$, we may
restrict $f$ to a small open ball (or polydisc) $\Omega_p$ centered at
$p$. Applying the jump Theorem~1.11, we may represent $f = F^+ - F^-$,
with $F^\pm \in \mathcal{ O} ( \Omega_p^\pm) \cap \mathcal{ C}^0 (
\overline{ \Omega_p^\pm})$. If $\Omega_p^+$ (resp. $\Omega_p^-$) is
holomorphically extendable at $p$, then $F^+$ (resp. $F^-$) extends to
a neighborhood $\omega_p$ of $p$ in $\C^n$ as $G \in \mathcal{ O}
(\omega_p)$ (resp. $H \in \mathcal{ O} (\omega_p)$). Then $f$ extends
holomorphically to the small one-sided neighborhood $\omega_p^-$
(resp. to $\omega_p^+$) as $G - F^-$ (resp. as $F^+ - H$).

\def\thelemma{2.2}\begin{lemma}
On hypersurfaces, at a given point, local holomorphic extendability of
CR functions to one side is equivalent to holomorphic extendability to
the same side of the holomorphic functions defined in the opposite side.
\end{lemma}

Consequently, the theorems of \S1.22 yield gratuitously extension
results about CR functions. For instance:

\def\thecorollary{2.3}\begin{corollary}
{\rm (\cite{ befo1978, r1983, bt1984})} 
On a real analytic hypersurface $M$, at a given point $p$, continuous
CR functions extend holomorphically to one side if and only if $M$
does not contain any local complex hypersurface passing through $p$.
\end{corollary}

The assumption of real analyticity, or the assumption of finite
typeness in case $M$ is $\mathcal{ C }^\infty$, both consume much
smoothness. The removal of these assumptions was accomplished by
Tr\'epreau in 1986.

\def\thetheorem{2.4}\begin{theorem}
{\rm (\cite{ trp1986})} Let $M$ be a $\mathcal{ C}^2$ hypersurface of
$\C^n$, $n \geqslant 2$ and let $p \in M$. The following two
conditions are equivalent{\rm :}
\begin{itemize}

\smallskip\item[$\bullet$]
$M$ does not contain any local 
complex hypersurface passing through $p$. 

\smallskip\item[$\bullet$]
for every open subset $U_p \subset M$ containing $p$, there exists a
one-sided neighborhood $\omega_p^\pm$ of $M$ at $p$ with $\overline{
\omega_p^\pm} \cap M \Subset U_p$ such that for every $f \in \mathcal{
C}_{ CR}^0 (U_p)$, there exists $F \in \mathcal{ O} (\omega_p^\pm)
\cap \mathcal{ C}^0 ( \omega^\pm \cup U_p)$ with $F \vert_{ U_p} = f$.

\end{itemize}\smallskip

\end{theorem}

We have seen that characterizing the side of extension is an open
question, even for rigid polynomial hypersurfaces $v = \varphi_m (z,
\bar z)$ and even for $m=6$. Although the above theorem constitutes a
neat answer for holomorphic extension to some imprecise side, it does
{\it not} provide any control of the side of extension.

\smallskip

Let $M$ be a $\mathcal{ C}^2$ orientable connected hypersurface and
let $\Omega_M^+$ be an open side of $M$. One could hope to
characterize holomorphic extension to the other side at {\it every
point} of $M$, since weak pseudoconvexity characterizes holomorphic
{\it non-extendability} at {\it every point} of $M$, by Oka's theorem.

\smallskip

\def\theexample{2.5}\begin{example}{\rm
(\cite{ trp1992}) 
In $\C^3$, let $\Omega_M^+$ be $\big\{ v > \varphi_m ( z_1, \bar z_1)
- \vert z_2 \vert^2 \, \vert z_1\vert^{ 2N } \big\}$ where $\varphi_m
\not\equiv 0$, of degree $m$ with $3\leqslant m < N$ is as in Open
problem~1.30. One verifies that holomorphic extension at every point
of $M$ entails a characterization of holomorphic extension at the
origin for the domain $\big\{ v > \varphi_m (z, \bar z) - \varepsilon
\, \vert z \vert^{ 2N} \big\}$.
}\end{example}

In the sequel, we shall abandon definitely the difficult, still open
question of how to control sides of holomorphic extension.

\smallskip

Although Theorem~2.4 is well known in Several Complex Variables, there
is a more general formulation with a simpler proof than
the original one. The remainder of this section will expose such a
proof.

By a {\sl global one-sided neighborhood} of a connected (not
necessarily orientable) hypersurface $M \subset \C^n$, we mean a
domain $\Omega_M$ with $\overline{ \Omega}_M \supset M$ such that for
every point $q \in M$, at least one open side $\omega_q^\pm$ of $M$ at
$q$ is contained in $\Omega_M$. In fact, to insure connectedness,
$\Omega_M$ is the interior of the closure of the union $\cup_{ q\in
M}\, \omega_q^\pm$ of all (possibly shrunk) one-sided
neighborhoods. 

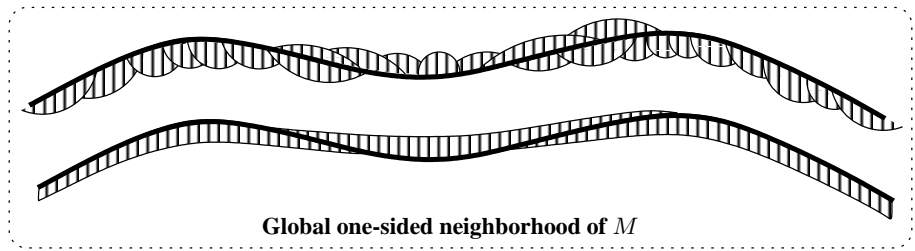
\begin{center}
\begin{picture}(0,0)%
\includegraphics{global-one-sided.pstex}%
\end{picture}%
\setlength{\unitlength}{4144sp}%
\begingroup\makeatletter\ifx\SetFigFont\undefined
\def\x#1#2#3#4#5#6#7\relax{\def\x{#1#2#3#4#5#6}}%
\expandafter\x\fmtname xxxxxx\relax \def\y{splain}%
\ifx\x\y   
\gdef\SetFigFont#1#2#3{%
  \ifnum #1<17\tiny\else \ifnum #1<20\small\else
  \ifnum #1<24\normalsize\else \ifnum #1<29\large\else
  \ifnum #1<34\Large\else \ifnum #1<41\LARGE\else
     \huge\fi\fi\fi\fi\fi\fi
  \csname #3\endcsname}%
\else
\gdef\SetFigFont#1#2#3{\begingroup
  \count@#1\relax \ifnum 25<\count@\count@25\fi
  \def\x{\endgroup\@setsize\SetFigFont{#2pt}}%
  \expandafter\x
    \csname \romannumeral\the\count@ pt\expandafter\endcsname
    \csname @\romannumeral\the\count@ pt\endcsname
  \csname #3\endcsname}%
\fi
\fi\endgroup
\begin{picture}(5424,1464)(439,-1335)
\put(1964,-1239){\makebox(0,0)[lb]{\smash{\SetFigFont{9}{10.8}{rm}{\color[rgb]{0,0,0}{\bf Global one-sided neighborhood of $M$}}%
}}}
\end{picture}

\end{center}

Then $\Omega_M$ contains a neighborhood in $\C^n$ of
every point $r\in M$ which belongs to at least two one-sided
neighborhoods that are opposite. The classical Morera theorem insures
holomorphicity in a neighborhood of such points $r$.

\smallskip

Remind that $M$ is called {\sl globally minimal} if it consists of a
single CR orbit. The assumption that $M$ does not contain any complex
hypersurface at any point means that for every $p \in M$, every open
$U_p \ni p$, the CR orbit $\mathcal{ O}_{ CR} (U_p, p)$ contains a
neighborhood of $p$ in $M$. This implies that $M$ is globally minimal
and hence, Theorem~2.4 is less general than the following.

\def\thetheorem{2.6}\begin{theorem}
{\rm (\cite{ trp1990, tu1994a})} Let $M$ be\footnote{ This theorem
also holds true with $M\in \mathcal{ C}^2$ and even with $M \in
\mathcal{ C}^{ 1, \alpha}$ ($0 < \alpha < 1$), provided one redefines
the notion of CR orbit in terms of boundaries of small attached
analytic discs.} a connected $\mathcal{ C}^{
2, \alpha}$ $(0 < \alpha < 1)$ hypersurface of $\C^n$ $(n\geqslant
2)$. If $M$ is globally minimal, then there exists a global one-sided
neighborhood $\Omega_M$ of $M$ such that for every continuous CR
function $f\in \mathcal{ C}_{ CR}^0 (M)$, there exists $F \in
\mathcal{ O} (\Omega_M) \cap \mathcal{ C}^0 ( \Omega_M \cup M)$ with
$F \vert_M = f$.
\end{theorem}

It will appear that $\Omega_M$ is contructed by gluing discs to $M$
and to subsequent open sets $\Omega ' \subset \Omega_M$ which are all
contained in the {\sl polynomial hull} of $M$:
\[
\widehat{M}
:=
\Big\{
z\in\C^n:\,
\vert P(z)\vert\leqslant
\sup_{w\in M}\,
\vert P(w)\vert, \
\forall\ P \in\C[z]
\Big\}.
\]

Let us summarize the proof. Although the assumption of global
minimality is so weak that $M$ may incorporate large open Levi-flat
regions, there exists at least one point $p \in M$ in a neighborhood
of which
\[
T_qM
=
T_q^cM
+
[T_q^cM,T_q^cM],
\ \ \ \ \ \ \ \ \
q\in U_p.
\]
Otherwise, the distribution $p \mapsto T_p^c M$ would be
Frobenius-integrable and all CR orbits would be complex
hypersurfaces\,! At such a point $p$, the classical Lewy extension
theorem (\S2.10 below) insures that $\mathcal{ C}_{ CR}^0 (M)$ extends
holomorphically to (at least) one side at $p$.

\def\thetheorem{2.7}\begin{theorem}
{\rm (\cite{ trp1990, tu1994a})}
Let $M$ be a connected $\mathcal{ C}^{ 2, \alpha}$ hypersurface, not
necessarily globally minimal. If $\mathcal{ C}_{ CR}^0 (M)$ extends
holomorphically to a one-sided neighborhood at some point $p\in M$,
then holomorphic extension to one side
$\omega_q^\pm$ holds at {\rm every} point $q
\in \mathcal{ O}_{ CR} (M,p)$.

\end{theorem}

When $\mathcal{ O}_{ CR} (M, p) = M$ as in Theorem~2.6, the global
one-sided neighborhood $\Omega_M$ will be the interior of the closure
of the union $\cup_{ q\in M}\, \omega_q^\pm$ of all (possibly
shrunk) one-sided neighborhoods.

\smallskip

The next paragraphs are devoted to expose a detailed proof of both the
Lewy theorem and of the above propagation theorem.

\subsection*{ 2.8.~Approximation theorem and maximum principle}
According to the approximation Theorem~5.2(III), for every $p \in
M$, there exist a neighborhood $U_p$ of $p$ in $M$ and a sequence
$(P_\nu (z))_{ \nu \in \N}$ of holomorphic polynomials with $\lim_{
\nu \to \infty} \, \left\vert \! \left\vert P_\nu - f \right\vert \!
\right\vert_{ \mathcal{ C}^0 (U_p)} =0$.

\def\thelemma{2.9}\begin{lemma}
For every analytic disc $A \in \mathcal{ O} (\Delta) \cap \mathcal{
C}^0 ( \overline{ \Delta})$ with $A (\partial \Delta) \subset U_p$,
the sequence $P_\nu$ also converges uniformly on the closed disc $A
(\overline{ \Delta})$, even if $A
(\overline{ \Delta})$ goes outside $U_p$.
\end{lemma}

\proof
By assumption, $\lim_{ \nu, \mu \to \infty} \, \left\vert \!
\left\vert P_\nu - P_\mu \right\vert \! \right\vert_{ \mathcal{ C}^0 
(U_p)} =0$. Let $\eta \in \overline{ \Delta}$
arbitrary. Thanks to the maximum principle and to
$A (\partial \Delta) \subset U_p$:
\[
\aligned
\left\vert\!\left\vert
P_\nu(A(\eta))
-
P_\mu(A(\eta))
\right\vert\!\right\vert
&
\leqslant
\max_{\zeta\in\partial\Delta}
\left\vert\!\left\vert
P_\nu(A(\zeta))
-
P_\mu(A(\zeta))
\right\vert\!\right\vert
\\
&
\leqslant
\sup_{z\in U_p}\,
\left\vert\!\left\vert
P_\nu(z)
-
P_\mu(z)
\right\vert\!\right\vert
\longrightarrow 0.
\endaligned
\]
The same argument shows that $P_\nu$ converges uniformly in the
polynomial hull of $U_p$ (we shall not need this).
\endproof

Next, suppose that we have some family of analytic discs $A_{ s}$,
with $s$ a small parameter, such that $\cup_s \, A_s (\Delta)$
contains an open set in $\C^n$, for instance a one-sided neighborhood
at $p \in M$. Then $(P_\nu)_{ \nu \in \N}$ converges uniformly on
$\cup_s \, A_s (\Delta)$ and a theorem due to Cauchy assures that the
limit is {\it holomorphic}\, in the interior of $\cup_s \, A_s
(\Delta)$. It then follows that $f$ extends holomorphically to the
interior of $\cup_s \, A_s (\Delta)$.

\smallskip

Remarkably, this short argument based on an application of the
approximation Theorem~5.2(III) shows that\footnote{ This is the
so-called {\sl Method of analytic discs}~; $\overline{ \partial}$
techniques are also powerful.}:

\begin{center}
\begin{minipage}[t]{10cm}
{\sf In order to establish local holomorphic extension of CR
functions, it suffices to glue appropriate families of analytic discs
to CR manifolds}.
\end{minipage}
\end{center}

\noindent
In the sequel, the geometry of glued discs will be studied for itself;
thus, it will be understood that statements about holomorphic or CR
extension follow immediately; elementary details about continuity of
the obtained extensions will be skipped.

\subsection*{ 2.10.~Lewy extension} 
Since $M$ is globally minimal, there exists a point $p$ at which $T_p
M = T_p^c M + [ T_p^cM, T_p^c M]$. Granted Lemma~2.2, holomorphic
extension to one side at such a point $p$ has already been proved in
Theorem~1.18. Nevertheless, we want to present a geometrically
different proof that will produce preliminaries and motivations for
the sequel.

Since $T^c M = {\rm Re} \, T^{ 1, 0} M = {\rm Re} \, T^{ 0, 1} M$, we
have equivalently $\left[ T^{ 1, 0} M , T^{ 0, 1} M \right] (p)
\not\subset \C \otimes T_p^c M$, namely the intrinsic Levi form of $M$
at $p$ is nonzero. In other words, there exists a local section $L$ of
$T^{ 1, 0} M$ with $L(p) \neq 0$ and $\left[ L, \overline{ L} \right]
(p) \not \in \C \otimes T_p^cM$. After a complex linear transformation
of $T_p^cM$, we may assume that $L( p ) = \frac{ \partial }{\partial
z_1}$. After removing the second order pluriharmonic terms, there
exist local coordinates $(z_1, z', w)$ vanishing at $p$ such that $M$
is represented by
\[
v
=
-
z\bar z_1
+
{\rm O}(\vert z_1\vert^{2+\alpha})
+
{\rm O}(\vert z'\vert)
+
{\rm O}(\vert z\vert\vert u\vert)
+
{\rm O}(\vert u\vert^2).
\]
The minus sign is set for clarity in the diagram of
\S2.12 below.
We denote by $\varphi (z_1, z', u)$ the right hand side. Let
$\varepsilon_1 >0$ be small. For $\varepsilon$ satisfying $0 <
\varepsilon \leqslant \varepsilon_1$, we introduce the analytic disc
\[
A_\varepsilon(\zeta)
:=
\left(
\varepsilon(1-\zeta),0',U_\varepsilon(\zeta)
+i\,V_\varepsilon(\zeta)
\right)
\]
with zero $z'$-component, with $z_1$-component equal to a (reverse)
round disc of radius $\varepsilon$ centered at $1\in \C$ and with
$u$-component $U_\varepsilon$ satisfying the Bishop-type equation:
\[
U_\varepsilon(e^{i\,\theta})
=
-
{\sf T}_1
\left[
\varphi(\varepsilon(1-\cdot),0',U_\varepsilon(\cdot))
\right]
(e^{i\,\theta}).
\]
Acoording to Theorem~3.7(IV), a unique solution $U_\varepsilon (e^{
i\, \theta})$ exists and is $\mathcal{ C}^{2, \alpha-0}$ with respect
to $(e^{i\, \theta}, \varepsilon)$. Since ${\sf T}_1 ( \psi)(1) =0$ by
definition, we have $U_\varepsilon (1) = 0$ and then $V_\varepsilon :=
{\sf T}_1 (U_\varepsilon)$ also satisfies $V_\varepsilon (1) = 0$.
Consequently, $A_\varepsilon (1) = 0$. By applying ${\sf T}_1$ to both
sides of the above equation, we see that the disc is attached to $M$:
\[
V_\varepsilon(e^{i\,\theta})
=
\varphi
\left(
\varepsilon(1-e^{i\,\theta}),0',
U_\varepsilon(e^{i\,\theta})
\right).
\]
We shall prove that for $\varepsilon_1$ sufficiently small, {\it
every}\, disc $A_\varepsilon (\Delta)$ with $0 < \varepsilon \leqslant
\varepsilon_1$ is {\it not}\, tangent to $M$ at $p$. We draw two
diagrams: a $2$-dimensional and a $3$-dimensional view. In both, the
$v$-axis is vertical, oriented down.

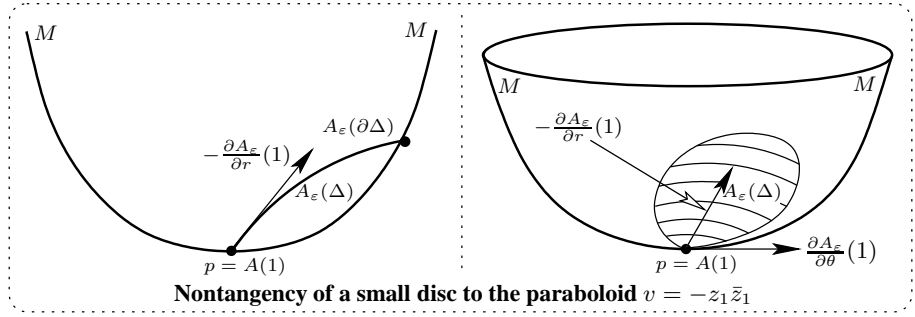
\begin{center}
\input paraboloid.pstex_t
\end{center}

Just now, we need a geometrical remark. Let $A \in \mathcal{ O}
(\Delta) \cap \mathcal{ C}^1 ( \overline{ \Delta})$ be an arbitrary
but small analytic disc attached to $M$ with $A(1)=0$. We use polar
coordinates to denote $\zeta=r\, e^{i\theta}$.

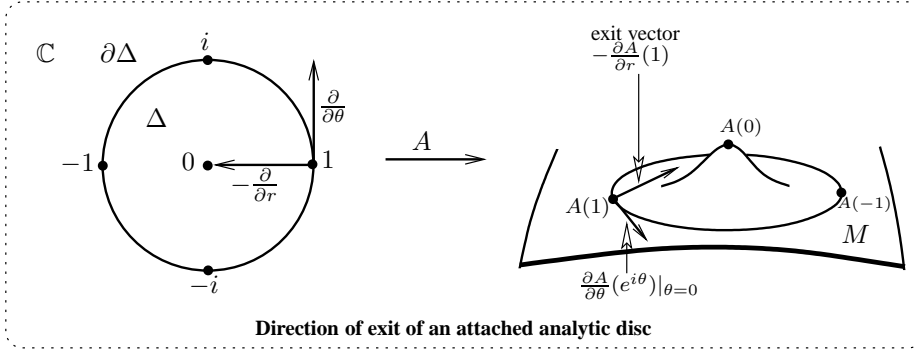
\begin{center}
\input tangent-at-one.pstex_t
\end{center}

\noindent
The holomorphicity of $A$ yields the following identities between
vectors in $T_p \C^n$:
\[
i\left.
{\partial A \over \partial \theta}
(e^{i\theta}) \right\vert_{\theta=0}=
- \left.
{\partial A\over \partial r}(r) \right\vert_{r=1}=
- \left. 
{\partial A \over \partial \zeta}(\zeta) \right\vert_{\zeta=1}.
\]
The multiplication by $i$ (or equivalently the complex structure $J$)
provides an isomorphism $T_p \C^n / T_p M \to T_p M / T_p^cM$; in
coordinates, $T_p \C^n / T_p M \simeq \R_v$, $T_p M / T_p^c M
\simeq \R_u$ and $J$ sends $\R_u$ to $\R_v$.
It follows that $\frac{ \partial A }{ \partial
r} (1)$ is not tangent to $M$ at $p$ if and only if $\frac{ \partial
A}{ \partial \theta} (1)$ is not complex tangent to $M$ at $p$.

\smallskip

Coming back to $A_\varepsilon$, we call the vector
\[
-
\frac{\partial A_\varepsilon}{
\partial r}(1)\ 
{\rm mod}\, T_pM
=
-
\frac{\partial W_\varepsilon}{
\partial r}(1)\ 
{\rm mod}\, T_pM
\]
the {\sl exit vector} of $A_\varepsilon$. By differentiating
$V_\varepsilon = \varphi$ at $\theta = 0$, taking account of
$d\varphi(0) = 0$, we get $\frac{ \partial V_\varepsilon}{ \partial
\theta} (1) = 0$. So only the real part $\frac{ \partial
U_\varepsilon}{ \partial \theta} (1)$ of $\frac{ \partial
W_\varepsilon}{ \partial \theta} (1)$ may be nonzero.

\def\thelemma{2.11}\begin{lemma}
Shrinking $\varepsilon_1$ if necessary, the exit vector of every disc
$A_\varepsilon$ with $0 < \varepsilon \leqslant \varepsilon_1$ is
nonzero{\rm :}
\[
-
\frac{\partial W_\varepsilon}{\partial r}(1)
=
i\,
\frac{\partial W_\varepsilon}{\partial \theta}(1)
=
i\,
\frac{\partial U_\varepsilon}{\partial \theta}(1)
\neq 0.
\]
\end{lemma}

\proof
The principal term of $\varphi$ is $-z_1\bar z_1$. We
compute first:
\[
\aligned
{\sf T}_1
\left[
-Z_1(\zeta)\overline{Z}_1(\zeta)
\right]
&
=
{\sf T}_1
\left[
\varepsilon^2(e^{-i\,\theta}-2+e^{i\,\theta})
\right]
\\
&
=
\frac{1}{i}\,\varepsilon^2(
-
e^{-i\,\theta}
+
e^{i\,\theta}).
\endaligned
\]
Proceeding as carefully as in Section~3(IV), we may verify that
\[
\aligned
U_\varepsilon(e^{i\,\theta})
&
=
-{\sf T}_1
\left[
-
Z_1(\zeta)\overline{Z}_1(\zeta)
+
{\sf Remainder}
\right](e^{i\,\theta})
\\
&
=
-
2\,\varepsilon^2
\sin\theta
+
\widetilde{U}_\varepsilon(e^{i\,\theta}),
\endaligned
\]
with a $\mathcal{ C}^{ 2, \alpha- 0}$ remainder satisfying $\big\vert
\! \big\vert \widetilde{ U}_\varepsilon \big\vert \! \big\vert_{ 1, 0}
\leqslant {\sf K} \, \varepsilon^{ 2+ \alpha }$, for some quantity
${\sf K} >0$. So $\frac{ \partial U_\varepsilon}{ \partial \theta}
(1) = - 2\, \varepsilon^2 + {\rm O} ( \varepsilon^{ 2+ \alpha }) \neq
0$.
\endproof

\subsection*{ 2.12.~Translations of a nontangent analytic disc}
We now fix $\varepsilon$ with $0 < \varepsilon \leqslant
\varepsilon_1$ and we denote simply by $A$ the disc
$A_\varepsilon$. So the vector
\[
\frac{\partial A}{\partial\theta}(1)
=
\left(
-
i\,\varepsilon,0',
-
2\,\varepsilon^2+{\rm O}(\varepsilon^{2+\alpha})
\right)
\]
is not tangent to $T_p^c M = \{ v= u =0\}$ at the origin. Furthermore,
it is not tangent to the $(2n-2)$-dimensional sub-plane $\{ y_1 = v =
0\}$ of $T_pM = \{ v=0\}$.

We now introduce parameters of translation $x_1^0\in \R$, $z_0' \in
\C^{ n-2}$ and $u_0 \in \R$ with $\vert x_1^0\vert, \, \vert z_0'
\vert, \, \vert u_0 \vert < \delta_1$, where $0 < \delta_1 <<
\varepsilon$. The points in $M$ of coordinates
\[
\left(
x_1^0, z_0',
u_0+i\,\varphi(x_1^0,z_0',u_0)
\right)
\]
cover a small $(2n-2)$-dimensional submanifold 
$K_p$ with $T_p K_p = \{ y_1 = v = 0 \}$
transverse to the disc
boundary $A_{ \varepsilon} (\partial \Delta)$ at $p$ that we draw
below.

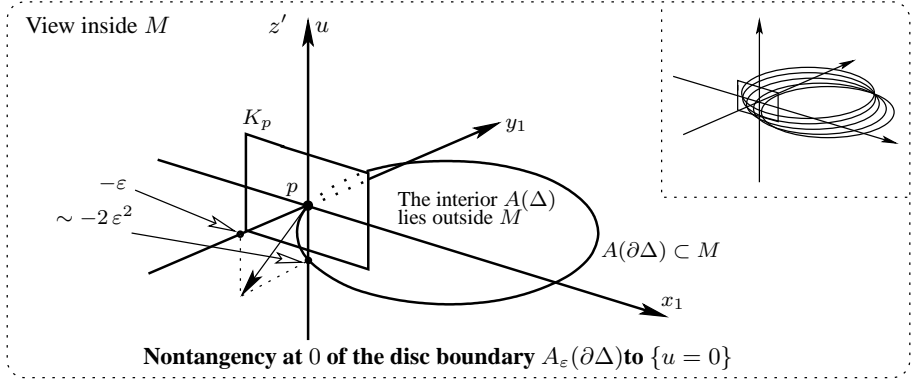
\begin{center}
\input inside-disc.pstex_t
\end{center}

To conclude the proof of one-sided holomorphic extension at the Levi
nondegenerate point $p$, it suffices to deform the disc
$A_{x_1^0,z_0',u_0}$ so that its distinguished point
$A_{x_1^0,z_0',u_0} (1)$ covers the submanifold $K_p$, namely
\def\theequation{2.13}\begin{equation}
A_{x_1^0,z_0',u_0}(1)
=
\left(
x_1^0, z_0',
u_0+i\,\varphi(x_1^0,z_0',u_0)
\right).
\end{equation}
This may be achieved easily by defining
\[
\big(
Z_{1,x_1^0}(\zeta),Z_{z_0'}'(\zeta)
\big)
:=
\big(
\varepsilon_1(1-\zeta)+x_1^0,
z_0'
\big)
\]
and by solving the Bishop-type equation:
\def\theequation{2.14}\begin{equation}
U_{x_1^0,z_0',u_0}(e^{i\,\theta})
=
u_0
-
{\sf T}_1
\left[
\varphi
\left(
Z_{1;x_1^0}(\cdot),Z_{z_0'}'(\cdot),
U_{x_1^0,z_0',u_0}(\cdot)
\right)
\right]
(e^{i\,\theta})
\end{equation}
for the $u$-component of the sought disc $A_{ x_1^0, z_0', u_0}$.
Thanks to Theorem~3.7(IV), the solution exists and is $\mathcal{C}^{
2, \alpha-0}$ with respect to all the variables. We finally define the
$v$-component of $A_{ x_1^0, z_0', u_0}$:
\def\theequation{2.15}\begin{equation}
V_{x_1^0,z_0',u_0}(e^{i\,\theta})
:=
{\sf T}_1
\left[
U_{x_1^0,z_0',u_0}(\cdot)
\right]
(e^{i\,\theta})
+
\varphi(x_1^0,z_0',u_0).
\end{equation}
Applying ${\sf T}_1$ to~\thetag{ 2.14}, we see that this disc is
attached to $M$; also, putting $e^{ i\, \theta} := 1$ in~\thetag{
2.14} and in~\thetag{ 2.15}, we see that~\thetag{ 2.13} holds.
Geometrically, the $(2 n -2)$ added parameters $(x_1^0, z_0', u_0)$
correspond to translations in $M$ of the
original disc $A_{ \varepsilon_1}$.

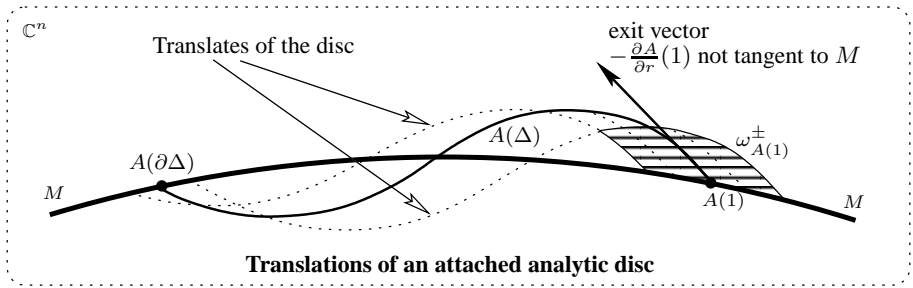
\begin{center}
\input translation-disc.pstex_t
\end{center}

Define the open circular region $\Delta_1 := \{ \zeta\in \Delta: \,
\vert \zeta -1 \vert < \delta_1\}$ around $1$ in the unit disc. Then
we claim that shrinking $\delta_1 >0$ if necessary, the set
\[
\Big\{
A_{x_1^0,z_0',u_0}(\zeta):\,
\zeta\in\Delta_1,\
\vert x_1^0\vert
<
\delta_1,\
\vert z_0'\vert
<
\delta_1,
\vert u_0\vert
<
\delta_1
\Big\}
\]
contains a one-sided neighborhood of $M$ at $p = A_{ 0, 0, 0}
(1)$. Indeed, by computation, one may check that the $2n$ vectors of
$T_p \C^n$
\[
\frac{\partial A_{0,0,0}}{\partial x_1}(1),
\ \
\frac{\partial A_{0,0,0}}{\partial\theta}(1),
\ \
\frac{\partial A_{0,0,0}}{\partial x_k'}(1),
\ \
\frac{\partial A_{0,0,0}}{\partial y_k'}(1),
\ \
\frac{\partial A_{0,0,0}}{\partial u}(1),
\ \
-\frac{\partial A_{0,0,0}}{\partial r}(1),
\]
are linearly independent; geometrically and by construction, the first
$(2n -1)$ of these vectors span $T_p M$ and the last one is linearly
independent, since by construction the exit vector of $A_{
\varepsilon_1 }$ is nontangent to $M$ at $p$.
\qed

Incidentally, we have proved an elementary but crucial statement: by
``translating'' (through a suitable Bishop-type equation) any small
attached disc whose exit vector is nonzero, we may always cover a
one-sided neighborhood.

\def\thelemma{2.16}\begin{lemma}
If a small disc $A$ attached to a hypersurface $M$ satisfies
$\frac{\partial A}{\partial \theta} (1) \not \in T_{ A(1) }^cM$, or
equivalently $- \frac{ \partial A}{ \partial r} (1) \not \in T_{ A(1)}
M$, then continuous CR functions on $M$ extend holomorphically at
$A(1)$ to the side in which points the nonzero exit vector $i\, \frac{
\partial A}{ \partial \theta} (1) = - \frac{ \partial A}{ \partial r}
(1)$.
\end{lemma}

Of course, the choice of the point $1 \in \partial \Delta$ is no
restriction at all, since after a M\"obius reparametrization, any
given point $\zeta_0 \in \partial \Delta$ becomes $1\in \partial
\Delta$.
 
\subsection*{2.17.~Propagation of holomorphic extension} 
The Levi form assumption $T_p M = T_p^c M + [ T_p^c M , T_p^c M ]$ was
strongly used to insure the existence of a disc having a nonzero exit
vector at $p$. But if a disc $A$ is attached to a highly degenerate
part of $M$, for instance to a region where the Levi form is nearly
flat, the disc $A$ might well satisfy $\frac{ \partial A }{ \partial
\theta} (\zeta_0) \in T_{ A(\zeta_0)}^c M$ (or equivalently, $- \frac{
\partial A}{ \partial r} (\zeta_0) \in T_{ A (\zeta_0)} M$), for every
$\zeta_0 \in \partial \Delta$. Then we are stuck.

\smallskip

To go through, two strategies are known in the literature.

\begin{itemize}

\smallskip\item[$\bullet$]
Devise refined pointwise ``finite type'' assumptions insuring the
existence of small discs having nonzero exit vector at a given central
point.

\smallskip\item[$\bullet$]
Devise deformation arguments that propagate holomorphic extension from
Levi nondegenerate regions up to highly degenerate regions.

\end{itemize}\smallskip

Unfortunately, the first, more developed strategy is unable to provide
any proof of Theorem~2.6. Indeed, a smooth globally minimal
hypersurface may well contain large Levi-flat regions, as for instance
$\{ (z, w) \in \C^2: \, v = \varpi (x) \}$ with a $\mathcal{
C}^\infty$ function $\varpi$ satisfying $\varpi (x) \equiv 0$ for
$x\leqslant 0$ and $\varpi (x) > 0 $ for $x > 0$ (to check global
minimality, proceed as in Example~3.10); Theorem~4.8(III) shows that a
Levi-flat portion $M_{\sf LF}$ of a hypersurface $M$ is locally
foliated by complex $(n-1)$-dimensional submanifolds; the uniqueness
in Bishop's equation\footnote{ A more general property holds true
({\it see}~\cite{ trp1990, tu1994a, mp2006b}): every small attached
disc is necessarily attached to a single (local or global) CR orbit;
here, $\Sigma$ is a local orbit, whence $A (\partial \Delta) \subset
\Sigma$.} then entails that every small analytic disc $A \in \mathcal{
O} (\Delta) \cap \mathcal{ C}^1 (\overline{ \Delta})$ with $A
(\partial \Delta) \subset M_{ \sf LF}$ must satisfy $A (\partial
\Delta) \subset \Sigma$, where $\Sigma \subset M_{ \sf LF}$ is the
unique local complex connected $(n - 1)$-dimensional submanifold of
the foliation that contains $A(1)$; then the uniqueness principle for
holomorphic maps between complex manifolds yields $A (\overline{
\Delta}) \subset \Sigma$; finally, $- \frac{ \partial A}{ \partial r}
(\zeta_0) \in T_{ A( \zeta_0)} \Sigma = T_{ A (\zeta_0)}^c M$ has exit
vector tangential to $M$ at every $\zeta_0 \in \partial \Delta$.

For this reason, we will focus our attention only on the second, most
powerful strategy, starting with a review.

After works of Sj\"ostrand (\cite{ hs1982, sj1982a, sj1982b}) on
propagation of singularities for certain classes of partial
differential operators, of Baouendi-Chang-Treves~\cite{ bct1983}, and
of Hanges-Treves~\cite{ ht1983}, Tr\'epreau~\cite{ trp1990} showed
that the hypoanalytic wave-front set of a CR function or distribution
propagates along complex-tangential curves. The microlocal technique
involves deforming $T^* M$ inside conic sets and controlling a certain
oscillatory integral called {\em Fourier-Bros-Iagolnitzer} ({\sc fbi})
transform. In 1994, Baouendi-Rothschild-Tr\'epreau~\cite{ brt1994}
showed how to deform analytic discs attached to a hypersurface in
order to get some propagation results (however, Theorem~2.7 which
appears in~\cite{ trp1990} is not formulated in~\cite{ brt1994}).
Then Tumanov~\cite{ tu1994a} showed how to deformation discs attached
to generic submanifolds of arbitrary codimension and provided
extension results that cannot be obtained by means of the usual
microlocal analysis.

Until the end of Section~4, our goal will be to describe and to
exploit this technique of propagation. The geometric idea is as
follows.

As in Theorem~2.7, assume that holomorphic extension is already known
to hold in a one-sided neighborhood $\omega_q^\pm$ at some point $q\in
M$. Referring to the diagram after the main Proposition~2.21 below, we
may pick a disc $A$ with $A(-1) = q$. Then a small part of its
boundary, namely for $e^{ i\theta}$ near $-1$, lies in $\overline{
\omega_q^\pm}$. If the vector $\frac{ \partial A}{ \partial \theta}
(1)$ is not complex tangential at the opposite point $p = A(1)$, it
suffices to apply Lemma~2.16 just above to get holomorphic extension
at $p$, almost gratuitously. On the contrary, if $\frac{\partial A}{
\partial \theta} (1)$ is complex tangential at $p$, we may well hope
that by slightly deforming $M$ as a hypersurface $M^{\sf d}$ which
goes inside $\omega_q^\pm$ a bit, there exists a deformed disc $A^{\sf
d}$ attached to $M^{\sf d}$ with again $A^{\sf d} (1) = p$ that will
be not tangential: $-\frac{ \partial A^{\sf d}}{ \partial r} (1) \not
\in T_{ A^{\sf d} (1)} M$. Then a translation of the disc $A^{ \sf d}$
as in Lemma~2.16 will provide holomorphic extension at $p$.

\subsection*{2.18.~Approximation theorem and chains of analytic discs} 
To prove Theorem~2.7, we first formulate a version of
the approximation theorem which is apppropriate for
our purposes.

\def\thelemma{2.19}\begin{lemma}
{\rm (\cite{ tu1994a})}
For every $p\in M$, there exists a neighborhood $U_p$ of $p$ in $M$
such that for every $q\in U_p$, for every one-sided neighborhood
$\Omega_q^\pm$ of $U_p$ at $q$, there exists a smaller one-sided
neighborhood $\omega_q^\pm \subset \Omega_q^\pm$ of $U_p$ at $q$ such
that the following approximation property holds{\rm :}
\begin{itemize}

\smallskip

\item[$\bullet$]
for every $F \in \mathcal{ C}^0 (M \cup \Omega_q^\pm)$ which
is CR on $M$ and holomorphic in $\Omega_q^\pm$, there exists a
sequence of holomorphic polynomials $(P_\nu(z))_{ \nu \in \N}$ such
that $0 = \lim_{ \nu \to \infty}\, \left\vert \! \left\vert P_\nu - f
\right\vert \! \right\vert_{ \mathcal{ C}^0 (U_p \cup
\omega_q^\pm)}$.

\end{itemize}\smallskip
\end{lemma}

The proof is an adaptation of Theorem~5.2(III). It suffices to allow
the maximally real submanifolds $\Lambda_u \subset M$ be slightly
deformed in $\Omega_q^\pm$. With a control of the smallness of their
$\mathcal{ C}^1$ norm, one may insure that they cover not only $U_p$
but also $\omega_q^\pm$. Further details will not be provided.

\begin{center}
\begin{minipage}[t]{10cm}
{\sf To establish local holomorphic extension of CR
functions, it is allowed to glue discs not only to 
$M$ but also to previously constructed one-sided neighborhoods}.
\end{minipage}
\end{center}

\smallskip

Pursuing, we formulate a lemma and a main proposition. 

\def\thelemma{2.20}\begin{lemma}
{\rm (\cite{ tu1994a})} Let $p\in M$ and let $U_p$ be a neighborhood
of $p$ in $M$, arbitrarily small. For every $q \in \mathcal{ O}_{ CR}
(M, p)$ and every small $\varepsilon >0$, there exist $\ell \in \N$
with $\ell = {\rm O} (1/ \varepsilon)$ and a chain of $\mathcal{ C}^{
2, \alpha - 0}$ analytic discs $A^1, A^2, \dots, A^{\ell - 1}, A^\ell$
attached to $M$ with the properties{\rm :}

\begin{itemize}

\smallskip\item[$\bullet$]
$A^1(-1) \in U_p$, i.e. this point 
is arbitrarily close to $p${\rm ;}

\smallskip\item[$\bullet$]
$A^1(1)=A^2(-1)$, $A^2(1)=A^3(-1)$, \ldots,
$A^{\ell-1}(1)=A^\ell(-1)${\rm ;}

\smallskip\item[$\bullet$]
$A^\ell(1)=q${\rm ;}

\smallskip\item[$\bullet$]
$\vert \! \vert A^k \vert \! \vert_{ \mathcal{
C}^{1, 0} (\overline{ \Delta })} \leqslant \varepsilon$, for $k=1, 2,
\dots, \ell${\rm ;}

\smallskip\item[$\bullet$]
each $A^k$ is an embedding $\overline{ \Delta} \to \C^n$.

\end{itemize}\smallskip
\end{lemma}

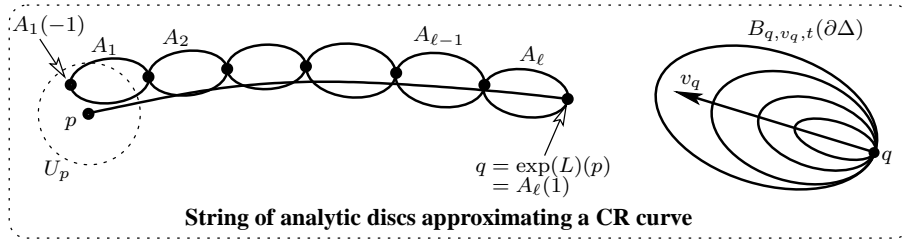
\begin{center}
\input string-discs.pstex_t
\end{center}

Such a chain of analytic discs will be constructed by approximating a
complex-tangential curve that goes from $q$ to $p$, using families of
discs $B_{ q, v_q, t} (\zeta)$ to be introduced in a while. The above
lemma is essentially obvious, whereas the next proposition constitutes
the very core of the argument.

\def\theproposition{2.21}\begin{proposition}
{\rm (\cite{ brt1994, tu1994a})} {\sf (Propagation along
a disc)} Let $A$ be a small $\mathcal{ C}^{ 2, \alpha - 0}$ analytic
disc attached to $M$ which is an embedding $\overline{ \Delta} \to
\C^n$. If $\mathcal{ C}_{ CR}^0 (M)$ extends holomorphically to a
one-sided neighborhood $\omega_{A(-1)}^\pm$ at the point $A (-1)$,
then it also extends holomorphically to a one-sided neighborhood at $A
(1)$. With more precisions{\rm :}

\begin{itemize}

\smallskip\item[$\bullet$]
if the exit vector $-\frac{ \partial A}{ \partial r} (1)$ is not
tangent to $M$ at $A(1)$, extension holds to the side in which points
$- \frac{ \partial A}{ \partial r} (1)${\rm :} this is already known,
by Lemma~2.16{\rm ;}

\smallskip\item[$\bullet$]
if the exit vector $-\frac{ \partial A}{ \partial r} (1)$ is tangent
to $M$ at $A(1)$, there exists an arbitrarily small deformation
$A^{\sf d}$ of $A$ with $A^{\sf d}(1) = A(1)$ having boundary $A^{\sf
d} (\partial \Delta)$ contained in $M \cup \omega_{ A( -1) }^\pm$ such
that the new exit vector $- \frac{\partial A^{\sf d}}{\partial r} (1)$
is not tangent to $M$ at $A^{ \sf d}(1)${\rm ;} 
then by translating $A^{ \sf d}$ as in Lemma~2.16, holomorphic 
extension holds at $A(1)$.

\end{itemize}\smallskip

\end{proposition}

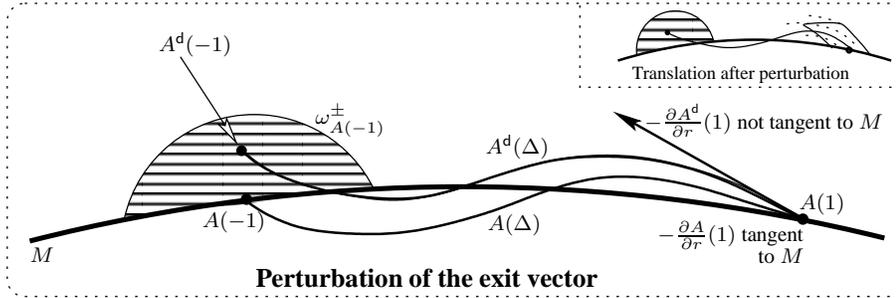
\begin{center}
\input perturbation-disc.pstex_t
\end{center}

Indeed, 
thanks to the flexibility of the solutions to the parametrized Bishop
equation provided by Theorem~3.7(IV), we can easily, 
as in Lemma~2.16, add translation
parameters $(x_1^0, z_0', u_0)$ to a slightly deformed disc $A^{\sf
d}$ attached to $M \cup \omega_{ A(-1)}^\pm$ and then $A_{ x_1^0,
z_0', u_0}^{\sf d}(\Delta_1)$ covers a small one-sided neighborhood of
$M$ at $A(1) = A^{\sf d} (1)$, thanks to the crucial condition $-
\frac{ \partial A^{ \sf d}}{ \partial r} (1) \neq 0$. We shall not
copy the details.

\smallskip
We claim that the proposition ends the proof of Theorem~2.7. By
assumption, $\mathcal{ C}_{ CR}^0 (M)$ extends holomorphically to a
one-sided neighborhood $\omega_p^\pm$ at $p$. The closure $\overline{
\omega_p^\pm}$ contains an open neighborhood $U_p$ of $p$. Let $q \in
\mathcal{ O}_{ CR} (M, p)$ and construct a chain of analytic discs
from $q$ up to a point $p' \in U_p$. The endpoint $p' = A^1 (-1)$ of the
chain of analytic discs being arbitrarily close to $p$, hence in
$U_p$, holomorphic extension holds at $A^1 (-1)$. We then apply the
proposition successively to the discs $A^1, A^2, \dots, A^\ell$ and
deduce holomorphic extension at $q$.

\smallskip

We now explain Lemma~2.20. To approximate a complex-tangential curve,
it suffices to construct families of analytic discs that are
essentially directed along given vectors $v_q \in T_q^cM$.

\def\thelemma{2.22}\begin{lemma}
For every point $q\in M$ and every nonzero complex tangent vector $v_q
\in T_q^cM \backslash \{ 0\}$, there exists a family of $\mathcal{
C}^{ 2, \alpha -0}$ analytic discs $B_{ q, v_q, t} (\zeta)$ 
parametrized by $t \in \R$ with $\vert t\vert < t_1$, for
some $t_1 >0$, that
satisfies{\rm :}

\begin{itemize}

\smallskip\item[$\bullet$]
$B_{q,v_q,t}(\partial \Delta) \subset M${\rm ;}

\smallskip\item[$\bullet$]
$q = B_{q,v_q,t} (1)${\rm ;}

\smallskip\item[$\bullet$]
$v_q = \frac{\partial B_{ q, v_q, 0}}{ \partial t} (-1)${\rm ;}

\smallskip\item[$\bullet$]
$\left\vert \! \left\vert B_{q, v_q, t} \right\vert \! \right\vert_{
\mathcal{ C}^{ 1, 0} (\overline{ \Delta})} \leqslant {\sf K} \, t$, for
some ${\sf K} >0$.

\end{itemize}\smallskip

\end{lemma}

\proof
In coordinates centered at $q$, represent $M$ by $v= \varphi (z, u)$
with $\varphi (0) = 0$ and $d\varphi (0) = 0$. The vector $v_q \in
T_p^c M = \{ w= 0\}$ has coordinates $(\dot{ z}_q, 0)$ for some 
nonzero $\dot{ z}_q \in \C^{
n-1}$. Introduce the family of analytic discs
\[
B_{q,v_q,t}(\zeta)
:=
\left(
t\,\dot{ z}_q(1-\zeta)/2,\,W_t(\zeta)
\right),
\]
where the real part $U_t$ of $W_t$ is the
unique $\mathcal{C}^{ 2, \alpha -0}$
solution of the Bishop-type equation:
\[
U_t(e^{i\,\theta})
=
-
{\sf T}_1
\big[
\varphi
\left(
t\,\dot{ z}_q(1-\cdot)/2,U_t(\cdot)
\right)
\big]
(e^{i\,\theta}).
\]
Proceeding as carefully as in Section~3(IV), we may verify that the
assumption $d\varphi (0) = 0$ implies that $\left\vert \! \left\vert
W_t \right\vert \! \right\vert_{ 1, 0} = {\rm O} (\vert t\vert^2)$.
Then it is obvious that $v_q = (\dot{ z}_q, 0) =
\frac{ \partial B_{ q, v_q, 0}}{
\partial t} (-1)$.
\endproof

We now complete the proof of Lemma~2.20. Any point $q \in \mathcal{
O}_{ CR} (M, p)$ is the endpoint of a finite concatenation of integral
curves of sections $L$ of $T^c M$. It suffices to construct the chain
of discs for a single such curve $\exp (tL) (p)$. After multiplying
$L$ by a suitable function, we may assume that $q$ is the time-one
endpoint $q = \exp ( L) (p)$.

Moving backwards, we start from $q_\ell := q$, we define $A^\ell
(\zeta) := B_{ q_\ell, -L(q_\ell), 1/\ell} (\zeta)$ and we set
$q_{\ell - 1} := B_{ q_\ell, -L(q_\ell), 1/\ell} (-1)$. Clearly,
$q_{\ell -1} = q_\ell - \frac{ 1}{\ell} \, L(q_\ell) + {\rm O} (\frac{
1}{ \ell^2})$. Starting again from $q_{ \ell - 1}$, we again move
backwards and so on, {\it i.e.} we define by descending induction:

\begin{itemize}

\smallskip\item[$\bullet$]
$A^k(\zeta)
:=
B_{q_k,-L(q_k),1/\ell}(\zeta)$;

\smallskip\item[$\bullet$]
$q_{k-1}:=B_{q_k,-L(q_k),1/\ell}(-1)$,

\end{itemize}\smallskip

\noindent
until $k = 1$. Since $q_{ k-1} = q_k - \frac{ 1}{ \ell} \, L( q_k) +
{\rm O} (\frac{ 1}{\ell^2})$ for $k=1, \dots, \ell$, the sequence of
points $q_k$ is a discrete approximation of the integral curve of $L$,
hence the endpoint $q_0 = A^1 (- 1)$ is arbitrarily close to $p$,
provided $\ell$ is large enough. Finally, by construction $\vert
\! \vert A^k \vert \! \vert_{ 1, 0} = {\rm O} (\frac{
1}{\ell})$.
\qed

\smallskip

The proof of the main Proposition~2.21 does not use special features of
hypersurfaces. For this reason, we will directly deal with generic
submanifolds of arbitrary codimension, passing to a new section.

\section*{\S3.~Tumanov's theorem, deformations of Bishop discs 
\\ 
and propagation on generic manifolds}

\subsection*{ 3.1.~Wedges and CR-wedges}
Assume now that $M$ is a connected generic submanifold in $\C^n$ of
codimension $d\geqslant 1$ and of CR dimension $m = d-n \geqslant
1$. The case $d= 1$ corresponds to a hypersurface. The notion of {\sl
local wedge} at a point $p$ generalizes to codimension $d\geqslant 2$
the notion of one-sided neighborhood at a point of a hypersurface.

More briefly that was has been done in Section~4(III), a wedge may be
defined as follows. Choose a $d$-dimensional real subspace $H_p$ of
$T_p\C^n$ satisfying $T_p \C^n = T_p H_p \oplus T_p M$ and a small
convex open salient truncated cone $C_p \subset H_p$ with vertex
$p$. Then a {\sl local wedge} of edge $M$ at $p$ is:
\[
\mathcal{W}(U_p,C_p)
:=
\{
q+{\sf c}:\,
q\in U_p,\,
{\sf c}\in C_p\}.
\]

This is not yet the most effective definition. Up to shrinking open
sets and parameter spaces, all definitions of local wedges will
coincide. Concretely, the wedges we shall construct will always been
obtained as unions of small pieces of families of analytic discs
partly attached to $M$. So we formulate all the technical conditions
that will insure that such pieces of discs cover a wedge.

\def\thedefinition{3.2}\begin{definition}{\rm
A {\sl local wedge} of edge $M$ at $p$ is 
a set of the form{\rm :}
\[
\mathcal{W}_p
:=
\Big\{
A_{t,s}\big(re^{i\,\theta}\big):\,
\vert t\vert<t_1,\,
\vert s\vert<s_1,\,
\vert \theta\vert<\theta_1,\,
r_1<r<1
\Big\},
\]
where, $t\in \R^{ d-1}$ is a rotation parameter, $t_1 >0$ is small,
$s\in \R^{ 2m+d-1}$ is a translation parameter, $s_1 >0$ is small,
$\theta_1 >0$ is small, $r_1 <1$ is close to $1$ and $A_{ t,s}
(\zeta)$, with $\zeta \in \overline{ \Delta}$, is a 
parametrized family of $\mathcal{ C}^{ 2, \alpha - 0}$
analytic discs satisfying:
\begin{itemize}

\smallskip\item[$\bullet$]
$A_{t,0} (1) = p$ for every $t$;

\smallskip\item[$\bullet$]
the boundaries $A_{t,s}( \partial \Delta)$ are partly (sometimes
completely) attached to $M$, namely $A_{ t,s} (e^{ i\, \theta}) \in
M$, at least for $\vert \theta \vert \leqslant \frac{ 3\pi }{ 2}$;

\smallskip\item[$\bullet$]
for every fixed $t$, the mapping $(s, e^{ i\, \theta}) \mapsto A_{
t,s} (e^{ i\, \theta})$ is a diffeomorphism from $\{ \vert s\vert <
s_1\} \times \{ \vert \theta \vert < \theta_1 \}$ onto a neighborhood
of $p$ in $M$;

\smallskip\item[$\bullet$]
the exit vector $- \frac{ \partial A_{0,0 }}{ \partial r} (1)$ is not
tangent to $M$ at $p$, namely it has nonzero projection ${\sf
proj\, }_{T_p \C^n/ T_pM} (- \partial A_{ t,0}/ \partial r(1))$ onto
the normal space $T_p\C^n / T_p M$ to $M$ at $p$;

\smallskip\item[$\bullet$]
choose any linear subspace $H_p$ of $T_p\C^n$ satisfying $T_p H_p
\oplus T_p M = T_p \C^n$, so that $H_p \simeq T_p \C^n / T_p M$, denote
by ${\sf proj\,}_{ H_p} : T_p \C^n \to H_p$ the projection onto $H_p$ parallel
to $T_pM$, define the associated exit vector 
\[
\text{\sf ex} (A_{t,0}) 
:=
{\sf proj\,}_{ H_p} ( - \frac{ \partial A_{t,0}}{ \partial r} (1)) \in H_p
\] 
and
the associated {\sl normalized exit vector} $\text{\sf n-ex} (A_{t,0})
:= \text{\sf ex} 
(A_{t,0}) / \vert \text{\sf ex} (A_{t,0}) \vert$; then the rank
at $t = 0$ of the mapping
\[
\R^{d-1}\ni
t
\longmapsto
\text{\sf n-ex}(A_{t,0})
\in S^{d-1}\subset
\R^d
\]
should be maximal equal to $d-1$.

\end{itemize}\smallskip
}\end{definition}

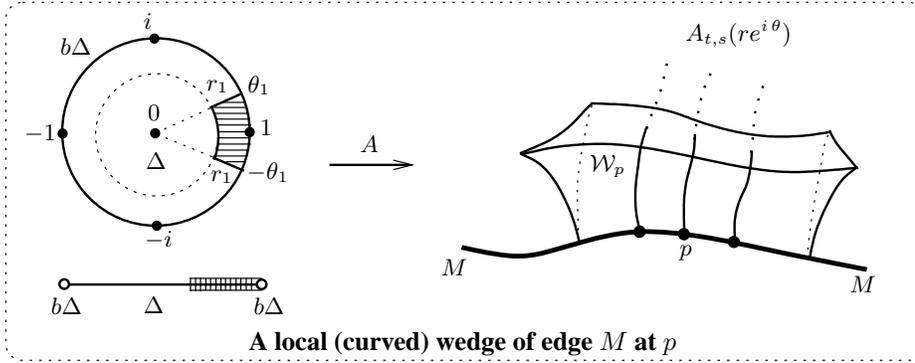
\begin{center}
\input curved-wedge.pstex_t
\end{center}

The last, most significant condition means that $\text{\sf n-ex}
(A_{t,0})$ describes an open neighborhood of $\text{\sf n-ex}
(A_{0,0})$ in the unit sphere $S^{ d-1} \subset \R^d$. This is of
course independent of the choice of $H_p$. Then, fixing $s=0$ and
$\theta = 0$, as the rotation parameter $t \in \R^{ d-1}$ varies with
$\vert t\vert < t_1$, and as the radius $r$ with $r_1 < r < 1$ varies,
the curves $A_{ t, 0} (r)$ generate an open truncated 
(curved) cone in some
$d$-dimensional local submanifold transverse to $M$ at $p$. Finally, as
the translation parameter $s$ varies, the points $A_{ t,s} (re^{i\,
\theta})$ describe a (curved) local wedge of edge $M$ at $p$.

\def\thelemma{3.3}\begin{lemma}
Shrinking $t_1>0$, $s_1 >0$, $\theta_1 >0$ and $1-r_1 >0$ if
necessary, the points of $\mathcal{ W}_p$ are covered injectively:
$A_{ t,s} (r e^{ i\, \theta}) = A_{ t',s'} ( r' e^{ i\, \theta '})$
if and only if $t=t'$, $s=s'$, $r=r'$ and $\theta = \theta'$.
\end{lemma}

This property follows directly from all the rank conditions. It will
be useful to insure uniqueness of holomorphic extension (monodromy).

\def\thedefinition{3.4}\begin{definition}{\rm
(\cite{ tu1990, trp1990}) A {\sl local CR-wedge} of edge $M$ at $p$ of
dimension $2m+ d+ e$, with $1 \leqslant e \leqslant d$, is a set
$\mathcal{ W}_p^{CR, e}$ defined similarly as a local wedge, but
assuming that the rotation parameter $t$ belongs to $\R^{ e-1}$ and
that the rank of the normalized exit vector mapping
\[
\R^{ e-1} \ni t
\longmapsto \text{\sf n-ex} (A_{ t,0}) \in
S^{d-1}\subset\R^d 
\]
is equal to $e-1$.
}\end{definition}

Then, fixing $s=0$ and $\theta = 0$, as the rotation parameter $t \in
\R^{ e-1}$ with $\vert t\vert < t_1$ varies, and as the radius $r$
with $r_1 < r < 1$ varies, the curves $A_{ t, 0} (r)$ describe an open
truncated (curved) cone in some $e$-dimensional
local submanifold transverse
to $M$ at $p$. These intermediate wedges of smaller dimension will
play a crucial technical r\^ole in the sequel.

The case $e=1$ deserves special attention. A CR-wedge is then just a
manifold with boundary $M_p^1$ with $\dim \, M_p^1 = 1+ \dim \, M$
that is attached to $M$ at $p$, namely there exists an open
neighborhood $U_p$ of $p$ in $M$ with $U_p \subset \partial M_p^1$.
If in addition $M$ has codimension $d=1$, we recover the notion of
one-sided neighborhood. It is clear that after a possible shrinking,
every $\mathcal{ C}^{ 2, \alpha - 0}$ manifold with boundary $M_p^1$
attached to $M$ at $p$ may be prolonged as a local $\mathcal{ C}^{ 2,
\alpha - 0}$ generic submanifold $\mathcal{ M}_p^1 \equiv \mathcal{
W}_p^{ CR, 1}$ containing a neighborhood of $p$ in $M$ (as shown
in the right
diagram).

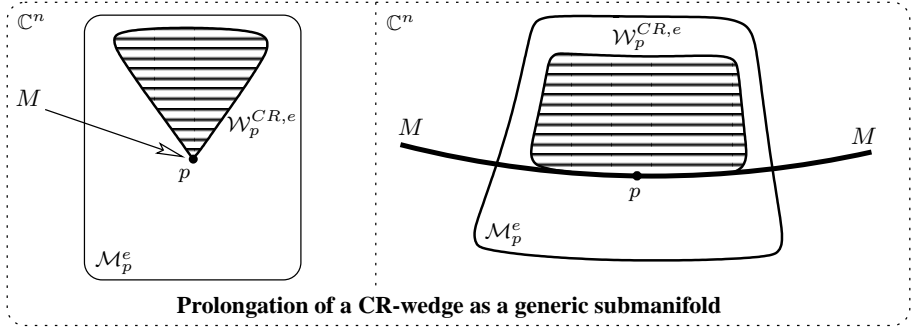
\begin{center}
\input prolonged-CR-wedge.pstex_t
\end{center}

By elementary differential geometry, for
$e \geqslant 2$, it
may be verified that
a local CR-wedge $\mathcal{ W}_p^{ CR, e}$ of
edge $M$ at $p$ defined by means of
a $\mathcal{ C}^{ 2, \alpha -
0}$ family of discs, namely
\[
\mathcal{W}_p^{CR,e}
:=
\Big\{
A_{t,s}\big(re^{i\,\theta}\big):\,
\vert t\vert<t_1,\,
\vert s\vert<s_1,\,
\vert \theta\vert<\theta_1,\,
r_1<r<1
\Big\},
\]
may also be prolonged as a local generic submanifold
$\mathcal{ M}_p^e$ of dimension $2m+ d + e$ containing a neighborhood
of $p$ in $M$. The left diagram is an illustration; in it, 
$e = d = 2$, so that $M$ of codimension $2$
is (unfortunately for intuition) collapsed to $p$.

However, the smoothness of $\mathcal{ M}_p^e$ can decrease to
$\mathcal{ C}^{ 1, \alpha - 0}$, because as in a standard local
blowing down $(z_1, z_2) \mapsto (z_1, z_1 z_2)$, the rank of the map
$(r, \theta, s, t) \longmapsto A_{t,s} \big(re^{ i\, \theta }\big)$
degenerates when $r=1$, since the discs (partial) boundaries $\big\{
A_{t,s} \big(e^{ i\, \theta }\big): \, \vert \theta \vert \leqslant
\frac{ 3 \pi}{ 2} \big\}$ are constrained to stay in $M$. For
technical reasons, we will need in the sequel the existence of a
prolongation $\mathcal{ M}_p^e$ that is $\mathcal{ C}^{ 2, \alpha -
0}$ also when $e \geqslant 2$. The following modification of the
definition of $\mathcal{ W}_p^{ CR, e}$ insures the existence of a
$\mathcal{ C}^{ 2, \alpha -0}$ prolongation $\mathcal{ M}_p^e$. It
will be applied implicitly in the sequel without further mention.

So, assume $e \geqslant 2$, let $A_{ t, s}$ be a family of discs as in
Definition~3.4 with $\text{ \sf ex} (A_{ 0, 0}) \neq 0$ in $T_p\C^n /
T_p M$ and $t\mapsto \text{ \sf n-ex} (A_{t, 0})$ of rank $e - 1$ at
$t = 0$. Fix $t := 0$ and define firstly
\[
\mathcal{W}_p^{CR,1}
:=
\big\{
A_{0,s}(re^{i\,\theta}):\,\vert s\vert<s_1,\
\vert\theta\vert<\theta_1,\
r_1<r<1
\big\}.
\]
This is a manifold with boundary attached to $M$ at $p$. So there
is a small $\mathcal{ C}^{ 2, \alpha - 0}$ prolongation $\mathcal{
M}_p^1 \supset \mathcal{ W}_p^{ CR, 1}$.

Choose $t\neq 0$ small with $A_{ t, 0}$ having exit vector nontangent
to $\mathcal{ M}_p^1$ at $p$. Introduce a one-parameter family
$M_\sigma$, $\sigma \in \R$, $\vert \sigma \vert < \sigma_1$,
$\sigma_1 >0$, of generic submanifolds obtained by deforming slightly
$M$ inside $\mathcal{ M}_p^1$ near $p$, with $M_\sigma \subset M \cup
\mathcal{ W}_p^{ CR, 1}$ for $\sigma \geqslant 0$. The $M_\sigma$ are
``translates'' of $M$ in $\mathcal{ M}_p^1$ near $p$. To understand
the process, we draw two diagrams in different
dimensions.

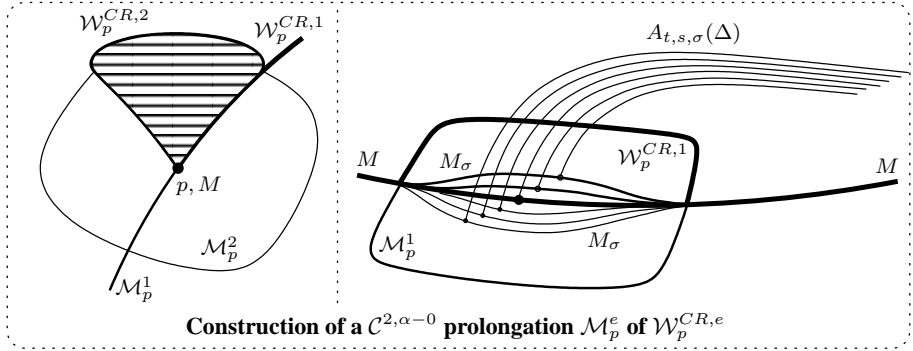
\begin{center}
\input modification-definition-CR-wedge.pstex_t
\end{center}

Thanks to the
flexibility of Bishop's equation (Theorem~3.7(IV)), the $A_{ t, s}$
may be deformed as a $\mathcal{ C}^{ 2, \alpha - 0}$ family $A_{ t, s,
\sigma}$ and we define secondly
\[
\mathcal{W}_p^{CR,2}
:=
\big\{
A_{t,s,\sigma}(r\,e^{i\,\theta}):\,
\vert s\vert<s_1,\
0<\sigma<\sigma_1,\
\vert\theta\vert<\theta_1,\
r_1<r<1
\big\}.
\]
Then this set constitutes a local CR-wedge of dimension $2m + d + 2$
with edge $M$ at $p$. Letting $\sigma$ run in $(-\sigma_1, \sigma_1)$
above, we get instead a certain manifold with boundary attached to
$\mathcal{ M}_p^1$ that may be extended as a $\mathcal{ C}^{ 2, \alpha
- 0}$ generic submanifold $\mathcal{ M}_p^2$ of dimension $2m + d +
2$. Then $\mathcal{ W}_p^{ CR, 2}$ is essentially one quarter of
$\mathcal{ M}_p^2$. We neither draw $\mathcal{ W}_p^{ CR, 2}$ nor
$\mathcal{ W}_p^2$ in the right diagram above, but the reader sees
them. By induction, using that $t\mapsto \text{ \sf n-ex} (A_{t, 0})$
has rank $e - 1$ at $t = 0$, we get the following.

\def\thelemma{3.5}\begin{lemma}
After a possible shrinking, a suitably constructed local $\mathcal{
C}^{ 2, \alpha - 0}$ CR-wedge $\mathcal{ W}_p^{ CR, e}$ of edge $M$ at
$p$ may be prolonged as a local $\mathcal{ C}^{ 2, \alpha - 0}$
generic submanifold $\mathcal{ M}_p^e$ of dimension $2m + d + e$
containing a neighborhood of $p$ in $M$.
\end{lemma}

In the sequel, similar technical constructions will be applied to
insure the existence of $\mathcal{ C}^{ 2, \alpha - 0}$ prolongations
$\mathcal{ M}_p^e \supset \mathcal{ W}_p^{ CR, e}$ without further
mention.

\subsection*{ 3.6.~Holomorphic extension of CR functions in higher
codimension} In 1988, Tumanov~\cite{ tu1988} established a theorem
that is nowadays celebrated in {\sl Several Complex Variables}.
Recall that by definition, $M$ is {\sl locally minimal at} $p$ if the
local CR orbit $\mathcal{ O}_{ CR}^{ loc} (M, p)$ contains a
neighborhood of $p$ in $M$. Equivalently, $M$ does not contain any
local submanifold $N$ passing through $p$ with ${\rm CRdim}\, N = {\rm
CRdim}\, M$ and $\dim N < \dim M$.

\def\thetheorem{3.7}\begin{theorem}
{\rm (\cite{ tu1988, brt1994, trp1996,
tu1998, ber1999})} Let $M$ be a local $\mathcal{ C}^{ 2,
\alpha}$ generic submanifold of $\C^n$ and let $p \in M$. If $M$ is
locally minimal at $p$, then there exists a local wedge
$\mathcal{ W}_p$ of edge $M$ at $p$ such that every $f \in
\mathcal{ C}_{ CR}^0 (M)$ possesses a holomorphic extension $F \in
\mathcal{ O}( \mathcal{ W}_p) \cap \mathcal{ C}^0 (M \cup \mathcal{
W}_p)$ with $F\vert_M = f$.
\end{theorem}

Conversely, recall that according to Theorem~4.41(III), if $M$ is not
locally minimal at $p$, there exists a local continuous CR function
that is not holomorphically extendable to any local wedge at $p$.

Since the literature already contains abundant
restitutions\footnote{We recommend mostly the two elegant
presentations~\cite{ trp1996} and~\cite{ tu1998}; other references
are: \cite{ brt1994, ber1999}. Excepting a conceptual abstraction
involving the implicit function theorem in Banach spaces and the
conormal bundle to $M$, the major arguments: differentiation of
Bishop's equation and a crucial correspondence between an {\sl exit
vector mapping} and an evaluation mapping defined on the space of
discs attached to $M$, the geometric structure of the proof is exactly
the same in the original article~\cite{ tu1988} as in the
restitutions.}, we will focus instead on propagation phenomena that
are less known.

\smallskip

In 1994, as an answer to a conjecture formulated by Tr\'epreau
in~\cite{ trp1990}, it was shown simultaneously by J\"oricke and by
the first author that Tumanov's theorem generalizes to globally
minimal $M$. The preceding statement is a direct corollary of the
next. Its proof given in~\cite{ me1994, jo1996} used techniques and
ideas of Tumanov~\cite{ tu1988, tu1994a} and of Tr\'epreau~\cite{
trp1990}.

\def\thetheorem{3.8}\begin{theorem}
{\rm (\cite{ me1994, jo1996})} Let $M$ be a connected
$\mathcal{ C}^{ 2, \alpha}$ generic submanifold of $\C^n$. If $M$ is
globally minimal then at every point $p\in M$, there exists a local
wedge $\mathcal{ W}_p$ of edge $M$ at $p$ such that every continuous
CR function $f \in \mathcal{ C}_{ CR}^0 (M)$ possesses a holomorphic
extension $F \in \mathcal{ O}( \mathcal{ W}_p) \cap \mathcal{ C}^0 (M
\cup \mathcal{ W}_p)$ with $F\vert_M = f$.
\end{theorem}

With this statement, the extension theorem for CR function has reached
a final, most general form. Philosophically, the main reason why it
is true lies in the propagation of holomorphic extendability along
complex-tangential curves. This was developed by Tr\'epreau in 1990,
using microlocal analysis.

\def\thetheorem{3.9}\begin{theorem}
{\rm (\cite{ trp1990})}
Let $M$ be a connected $\mathcal{ C}^\infty$ generic submanifold of
$\C^n$. If $\mathcal{ C}_{ CR}^0 (M)$ extends holomorphically to a
local wedge at some point $p\in M$, then at {\rm every} point $q \in
\mathcal{ O}_{ CR} (M, p)$, there exists a local wedge $\mathcal{
W}_q$ of edge $M$ at $q$ such that every $f \in \mathcal{ C}_{ CR}^0
(M)$ possesses a holomorphic extension $F \in \mathcal{ O}( \mathcal{
W}_q) \cap \mathcal{ C}^0 (M \cup \mathcal{ W}_q)$ with $F\vert_M =
f$.
\end{theorem}

Before surveying the original proof (\cite{ me1994, jo1996})
of this theorem in Section~5, we shall expose in length a
substantially simpler proof of Theorem~3.8 that was devised by the
second author in~\cite{ po2004}. This neat proof treats locally and
globally minimal generic submanifolds on the same
footing. It relies
partly upon a natural deformation proposition due to Tumanov in~\cite{
tu1994a}, but without any notion of defect of an analytic disc, without
any needs to control the variation of the direction of
CR-extendability, and without any partial connection, as in~\cite{
trp1990, tu1994a, me1994}. The next paragraphs and
Section~4 are devoted to the proof of this most general Theorem~3.8.

\def\theexample{3.10}\begin{example}
{\rm 
A globally minimal manifold may well be {\it
not}\, locally minimal at any point. 

Indeed,
let $\chi : \R \to \R^+$ be $\mathcal{ C}^\infty$ with $\chi = 0$ on
$(-\infty, 1]$, with $\chi > 0$ on $(1, +\infty)$ and with second
derivative $\chi_{ xx} >0$ on $( 1, + \infty)$. Consider the
generic manifold $M$ of $\C^3$ defined by the two equations 
\[
v_1
=
\chi(x),
\ \ \ \ \ \ \ \ \ \
v_2
=
\chi(-x),
\]
in coordinates $(x+i\,y, u_1+i\,v_1, u_2 + i\, v_2)$. Then
$T^{ 1, 0} M$ is generated by
\[
L
=
\frac{\partial}{\partial z}
+
i\,\chi_x(x)\,
\frac{\partial}{\partial w_1}
-
i\,\chi_x(-x)\,
\frac{\partial}{\partial w_2}.
\]
In terms of the four coordinates $(x, y, u_1, u_2)$ on $M$, the
two vector fields generating $T^cM$ are
\[
\aligned
L^1
:=
2\,{\rm Re}\,L
&
=
\frac{\partial}{\partial x},
\\
L^2
:=
2\,{\rm Im}\,L
&
=
\frac{\partial}{\partial y}
+
\chi_x(x)\,
\frac{\partial}{\partial u_1}
-
\chi_x(-x)\,
\frac{\partial}{\partial u_2}
\endaligned
\]
(we have dropped $\chi_x (x) \, \frac{ \partial}{ \partial v_1} -
\chi_x (-x)\, \frac{ \partial}{ \partial v_2}$ in $2\, {\rm Re}\, L$).
Denote by $\LL^0$ the system of these two vector fields $\{ L^1, L^2\}$ on
$\R^4 \simeq M$ and by $\LL$ the $\mathcal{ C}^\infty (\R^4)$-hull of
$\LL^0$. Observe that the Lie bracket
\[
\left[
L^1,L^2
\right]
=
\chi_{xx}(x)\,\frac{\partial}{\partial u_1}
+
\chi_{xx}(-x)\,\frac{\partial}{\partial u_2}
\]
is zero at points $p = (x_p, y_p, u_1^p, u_2^p)$ with $- 1 < x_p < 1$,
has non-zero $\frac{ \partial }{\partial u_2}$-component at points $p$
with $x_p < -1$ and has non-zero $\frac{ \partial}{\partial
u_1}$-component at points $p$ with $x_p >1$. 
It follows that the local $\LL$-orbit 
of a point $p$ with $x_p < -1$ is $\{ u_1 = u_1^p\}$, 
of a point $p$ with $- 1 < x_p < 1$ is $\{ u_1 = u_1^p, \, 
u_2 = u_2^p\}$ and of a point $x_p$ with $x_p > 1$ is 
$\{ u_2 = u_2^p\}$.
Also, observe that since
the vector field $L^1
= \frac{ \partial }{\partial x}$ belongs to $\LL$, the
local $\LL$-orbit of any point $p= (x_p, y_p, u_1^p, u_2^p)$ 
contains points of coordinates
$(x_p + t, y_p, u_1^p, u_2^p)$, with $t$ small. We deduce that the
local $\LL$-orbit of points $p$ with $x_p = - 1$ or $x_p = 1$ are
three-dimensional, hence in conclusion:
\[
\mathcal{O}_\LL^{loc}(\R^4,p)
=
\left\{
\aligned
&
U_p\cap\{u_1=u_1^p\}
\ \ \ \ \ \ \ \ \ \ 
\ \ \ \ \ \ \ \ \ \ 
\ \ \ \ \ 
{\rm if}\ x_p\leqslant -1,
\\
&
U_p\cap\{u_1=u_1^p,\,u_2=u_2^p\}
\ \ \ \ \ \ \ \ \ \
{\rm if}\ -1< x_p<1,
\\
&
U_p\cap\{u_2=u_2^p\}
\ \ \ \ \ \ \ \ \ \ 
\ \ \ \ \ \ \ \ \ \ 
\ \ \ \ \ 
{\rm if}\ x_p\geqslant 1,
\endaligned
\right.
\]
where $U_p$ is a neighborhood of $p$ in $M$.
So $\LL$ is nowhere locally minimal.

\def\thelemma{3.11}\begin{lemma}
The system $\LL$ is globally minimal.
\end{lemma}

\proof
We check that any two points $p, q\in \R^4$ are in the same
$\LL$-orbit. Using the flow of $L^1 = \frac{ \partial }{\partial x}$
and then the flow of $L^2$ on $\{ x= 0\}$, the original two points $p$
and $q$ may be joined to points, still denoted by $p = (0, 0, u_1^p,
u_2^p)$ and $q = (0, 0, u_1^q, u_2^q)$, having zero $x$-component and
zero $y$-component. 

We claim that the global $\LL$-orbit $\mathcal{ O}_{ \LL} (\R^4, p)$
of every point $p = (0, 0, u_1^p, u_2^p)$ contains a neighborhood of
$p$ in $\R^4$. Since the two-dimensional plane $\{ x = y = 0\}$ is
connected, this will assure that any two points $p = (0, 0, u_1^p,
u_2^p)$ and $q = (0, 0, u_1^q, u_2^q)$ are in the same $\LL$-orbit.

Indeed, by means of $\frac{\partial }{\partial x}$, every point $p = (
0, 0, u_1^p, u_2^p)$ is joined to the two points $p^- := (-1, 0,
u_1^p, u_2^p)$ and $p^+ := (1, 0, u_1^p, u_2^p)$. Let $U_{ p^-}$ and
$U_{ p^+}$ be small neighborhoods of $p^-$ and of $p^+$. Denote by
$H^- := \{ u_1 = u_1^p\} \cap U_{ p^-}$ and by $H^+ := \{ u_2 =
u_2^p\} \cap U_{ p^+}$ small pieces of the three-dimensional local
$\LL$-orbits of $p^-$ and of $p^+$.

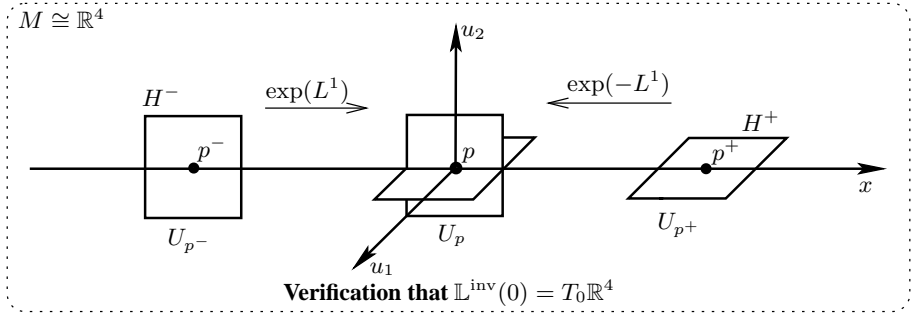
\begin{center}
\input globally-not-locally.pstex_t
\end{center}

\noindent
The flow of $L^1 = \frac{\partial }{\partial x}$ being a
translation, we deduce:
\[
\aligned
\exp(L^1)(H^-)
&
=
\{u_1=u_1^p\}
\cap 
U_p,
\\
\exp(-L^1)(H^+)
&
=
\{u_2=u_2^p\}
\cap
U_p,
\endaligned
\]
where $U_p$ is a small neighborhood of $p$ in $M\simeq \R^4$.
Observe that the
two $3$-dimensional 
planes are transversal 
in $T_p \R^4$. 
Lemma~1.28(III) yields:
\[
\aligned
\LL^{\rm inv}(p^-)
\supset
T_{p^-}\mathcal{O}_\LL^{loc}(p^-)
&
=
\{
u_1
=
u_1^p
\},
\\
\LL^{\rm inv}(p^+)
\supset
T_{p^+}\mathcal{O}_\LL^{loc}(p^+)
&
=
\{
u_2
=
u_2^p
\}.
\endaligned
\]
By the very definition of $\LL^{\rm inv}$, we necessarily
have:
\[
\aligned
\LL^{\rm inv}(p)
&
\supset
\exp(L^1)_*
\big(
\LL^{\rm inv}(p^-)
\big)
+
\exp(-L^1)_*
\big(
\LL^{\rm inv}(p^+)
\big)
\\
&
=
\{u_1=u_1^p\}
+
\{u_2=u_2\}
\\
&
=
T_p\R^4,
\endaligned
\]
so 
$\LL^{\rm inv} (p) = T_p \R^4$.
Consequently, $\mathcal{ O}_\LL (\R^4, p)$ contains
a neighborhood of $(0, 0, u_1^p, u_2^p)$ in $\R^4$.
\endproof

}\end{example}

\subsection*{ 3.12.~Setup for propagation}
Let $M$ be connected, generic and $\mathcal{ C}^{ 2, \alpha}$, let $q
\in M$ and let $\mathcal{ W}_q^{ CR, e}$ be a CR-wedge of dimension
$2m+d+e$ at $q$, with $1\leqslant e \leqslant d$. For short, we will
say that $\mathcal{ C}_{ CR}^0 (M)$ {\sl extends to be CR on}
$\mathcal{ W}_q^{ CR, e}$ if for every $f\in \mathcal{ C}_{ CR}^0$,
there exists $F \in \mathcal{ C}_{ CR}^0 (M \cup \mathcal{ W}_q^{ CR,
e})$ with $F\vert_M =f$.

\def\thetheorem{3.13}\begin{theorem}
Let $e\in \N$ with $1\leqslant e \leqslant d$. Assume that $\mathcal{
C}_{ CR}^0 (M)$ extends to be CR on a CR-wedge $\mathcal{ W}_p^{ CR,
e}$ of dimension $2m+ d+ e$ at some point $p\in M$. Then for every
$q\in \mathcal{ O}_{ CR} (M, p)$, there exists a CR-wedge $\mathcal{
W}_q^{ CR, e}$ at $q$ of the same dimension $2m + d + e$ to which
$\mathcal{ C}_{ CR}^0 (M)$ extends to be CR.
\end{theorem}

In the case $e=d$, we
recover\footnote{
Classical microlocal analysis was devised to measure the analytic wave
front set of a distribution in terms of the exponential decay ot the
Fourier transform restricted to open, conic submanifolds of the
cotangent bundle. We suspect that there might exist higher
generalizations of microlocal analysis in which one takes account of
the good decay of the Fourier transform on submanifolds of positive
codimension in the cotangent bundle.
}
Tr\'epreau's Theorem~3.9, since continuous CR functions on an open set
of $\C^n$ (here a usual wedge)
are just holomorphic. If $M$ is globally minimal, then extension
holds at every $q\in M$. Notice that this statement covers the
propagation Theorem~2.7, stated previously in the hypersurface case
$d=e=1$.

\smallskip

Let us start the proof. Through a chain of small analytic discs, every
$q\in \mathcal{ O}_{ CR} (M, p)$ is joined to a point $p'$ arbitrarily
close to $p$: indeed, Lemma~2.20 and its proof remain the same in
arbitrary codimension $d\geqslant 1$. At $p'$, CR extension holds,
because the edge of $\mathcal{ W}_p^{ CR, e}$ contains a small open
neighborhood $U_p$ of $p$ in $M$. To deduce CR extension at $q$, it
suffices therefore to propagate CR extension along a single disc, as
stated in the next main proposition.

\begin{center}
\begin{picture}(0,0)%
\includegraphics{normal-deformation.pstex}%
\end{picture}%
\setlength{\unitlength}{4144sp}%
\begingroup\makeatletter\ifx\SetFigFont\undefined
\def\x#1#2#3#4#5#6#7\relax{\def\x{#1#2#3#4#5#6}}%
\expandafter\x\fmtname xxxxxx\relax \def\y{splain}%
\ifx\x\y   
\gdef\SetFigFont#1#2#3{%
  \ifnum #1<17\tiny\else \ifnum #1<20\small\else
  \ifnum #1<24\normalsize\else \ifnum #1<29\large\else
  \ifnum #1<34\Large\else \ifnum #1<41\LARGE\else
     \huge\fi\fi\fi\fi\fi\fi
  \csname #3\endcsname}%
\else
\gdef\SetFigFont#1#2#3{\begingroup
  \count@#1\relax \ifnum 25<\count@\count@25\fi
  \def\x{\endgroup\@setsize\SetFigFont{#2pt}}%
  \expandafter\x
    \csname \romannumeral\the\count@ pt\expandafter\endcsname
    \csname @\romannumeral\the\count@ pt\endcsname
  \csname #3\endcsname}%
\fi
\fi\endgroup
\begin{picture}(5514,2049)(439,-1963)
\put(4647,-302){\makebox(0,0)[lb]{\smash{\SetFigFont{8}{9.6}{rm}{\color[rgb]{0,0,0}$H_p$}%
}}}
\put(2007,-115){\makebox(0,0)[lb]{\smash{\SetFigFont{8}{9.6}{rm}{\color[rgb]{0,0,0}$A_{t'}(\Delta)$}%
}}}
\put(972,-450){\makebox(0,0)[lb]{\smash{\SetFigFont{8}{9.6}{rm}{\color[rgb]{0,0,0}$\mathcal{W}_{A(-1)}^{CR,e}$}%
}}}
\put(5526,-1100){\makebox(0,0)[lb]{\smash{\SetFigFont{8}{9.6}{rm}{\color[rgb]{0,0,0}$C_p$}%
}}}
\put(5644,-1528){\makebox(0,0)[lb]{\smash{\SetFigFont{8}{9.6}{rm}{\color[rgb]{0,0,0}$M$}%
}}}
\put(547,-1432){\makebox(0,0)[lb]{\smash{\SetFigFont{8}{9.6}{rm}{\color[rgb]{0,0,0}$M$}%
}}}
\put(1051,-1381){\makebox(0,0)[lb]{\smash{\SetFigFont{8}{9.6}{rm}{\color[rgb]{0,0,0}$A(-1)$}%
}}}
\put(3015,-1411){\makebox(0,0)[lb]{\smash{\SetFigFont{8}{9.6}{rm}{\color[rgb]{0,0,0}$A(\Delta)$}%
}}}
\put(4806,-1473){\makebox(0,0)[lb]{\smash{\SetFigFont{8}{9.6}{rm}{\color[rgb]{0,0,0}$p$}%
}}}
\put(1870,-1876){\makebox(0,0)[lb]{\smash{\SetFigFont{9}{10.8}{rm}{\color[rgb]{0,0,0}{\bf Normal deformations of an analytic disc}}%
}}}
\end{picture}

\end{center}

\def\theproposition{3.14}\begin{proposition}
{\sf (Propagation along a disc)} {\rm (\cite{ tu1994a, mp1999},
[$*$])} Let $A$ be a small $\mathcal{ C}^{ 2, \alpha - 0}$ analytic
disc attached to $M$ which is an embedding $\overline{ \Delta} \to
\C^n$. Let $e \in \N$ with $1\leqslant e \leqslant d$. Assume that
there exists a $\mathcal{ C}^{ 2, \alpha - 0}$ CR-wedge $\mathcal{
W}_{ A (-1)}^{ CR, e}$ at $A(-1)$ of dimension $2m + d + e$ to which
$\mathcal{ C}_{ CR}^0 (M)$ extends to be CR. Then there exists a
$\mathcal{ C}^{ 2, \alpha - 0}$ CR-wedge $\mathcal{ W}_{ A (1)}^{ CR,
e}$ at $A(1)$ of the same dimension $2m + d + e$ to which $\mathcal{
C}_{ CR}^0 (M)$ extends to be CR.

With more precisions, the CR-wedge $\mathcal{ W}_{ A(1)}^{ CR, e}$ is
constructed by translating a certain family of analytic discs $A_{
t'}$ having the following properties. Setting $p:= A(1)$, there
exists a $\mathcal{ C}^{ 2, \alpha - 0}$ family $A_{ t'}$ of analytic
discs, $t'\in \R^e$, $\vert t'\vert < t_1'$, $t_1' >0$, with $A_{ t'}
\vert_{ t' = 0} = A$, with $A_{ t'} (1) = p$, satisfying $A_{ t'} (e^{
i\, \theta}) \in M$ for $\vert \theta \vert \leqslant \frac{ 3\pi }{ 2}$
and having their boundaries $A_{ t'} (\partial \Delta) \subset M \cup
\mathcal{ W}_{ A(-1)}^{ CR, e}$ for $t'$ belonging 
to some open truncated
cone ${\sf C}' \subset \R^e$, such that the exit vector mapping{\rm :}
\[
\R^e\ni t'
\longmapsto
\text{\sf ex}(A_{t'})
=
{\sf proj\,}_{H_p}
\left(
-
\frac{\partial A_{t'}}{
\partial r}(1)
\right)
\in\R^d
\]
is of maximal rank equal to $e$ at $t'=0$, where $H_p \simeq \R^d$ is
any linear subspace of $T_p\C^n$ such that $H_p \oplus T_pM = T_p
\C^n$, and where ${\sf proj\,}_{ H_p}$ is the linear projection
parallel to $T_pM$.
\end{proposition}

Geometrically, as $t'$ varies, the
exit vectors $\text{\sf ex} (A_{ t'})$ describe 
an open cone $C_p \subset H_p$, drawn in the diagram.

\smallskip

We claim that this statement covers the second, delicate case of
Proposition~2.21. Indeed assuming that $e=d=1$ and that the exit
vector $- \frac{\partial A}{\partial r} (1)$ is tangent to $M$ at
$A(1)$, the above proposition includes $A$ in a one-parameter family
$A_{ t'}$ whose direction of exit in the normal bundle has nonzero
derivative with respect to $t'$. Hence for every nonzero $t'$, the
direction of exit of $A_{ t'}$ is not tangent to $M$ at $p$. Thus, a
non-tangential deformed disc $A^{ \sf d}$ as in Proposition~2.21 may
be chosen to be any $A_{ t'}$, with $t'\neq 0$.

\proof[Proof of Proposition~3.14]
We first explain how to get CR extension at $p$ from the family $A_{
t'}$, taking for granted its existence.

\smallskip
\noindent
{\bf (I)}
Suppose firstly that the exit vector of $A = A_0$ is non-tangential to
$M$ at $p$. We have to restrict the parameter space $t' \in \R^e$ to
some parameter space $t\in \R^{ e-1}$ so as to reach Definition~3.4.

Let us take for granted the fact that the exit vector mapping has rank
$e$ at $t' = 0$. Then the normalized exit vector mapping
\[
\R^e\ni t'
\longmapsto
\text{\sf n-ex}(A_{t'})
=
\text{\sf ex}(A_{t'})/\vert \text{\sf ex}(A_{t'})\vert
\in S^{d-1}
\]
has rank $\geqslant e-1$ at $t'=0$. So there exists a small piece of an
$(e-1)$-dimensional linear subspace $\Lambda_0$ of $\R^d$,
parameterized as $t' = \phi (t)$ for some linear map $\phi$, with
$t\in \R^{ e-1}$ small, namely $\vert t \vert < t_1$, for some $t_1
>0$, such that $t\mapsto \text{\sf n-ex} (A_{\phi (t)})$ has rank
$(e-1)$ at $t=0$.

Setting $A_t := A_{ \phi (t)}$, we thus reach Definition~3.4,
without the translation parameter $s$. 

But proceeding exactly as in the hypersurface case, it is easy to
include some translation parameter getting a family $(A_{ t'})_s =
A_{ t', s}$. The proof is postponed to the end \S3.24. Then the desired
family $A_{ t, s}$ of the proposition is just $A_{ \phi(t), s}$,
shrinking $t_1 >0$ and $s_1 >0$ if necessary.

\def\thelemma{3.15}\begin{lemma}
There exists a deformation $A_{ t', s}$ of $A_{t'}$, with 
$s \in \R^{ 2m+ d- 1}$, $\vert s\vert < s_1$, such that{\rm :}
 
\begin{itemize}

\smallskip\item[$\bullet$]
the boundaries $A_{ t', s} (\partial \Delta)$ are contained in $M \cup
\mathcal{ W}_{ A(-1)}^{ CR, e}$ and $A_{ t'} (e^{ i\, \theta}) \in M$
for $\vert \theta \vert \leqslant \frac{ 3\pi }{ 2}${\rm ;}

\smallskip\item[$\bullet$]
for every fixed $t'$, the mapping $(s, e^{ i\, \theta}) \longmapsto
A_{ t', s} (e^{ i\, \theta})$ is a diffeomorphism from $\{ \vert
s\vert < s_1\} \times\{ \vert \theta \vert < \theta_1\}$ onto a
neighborhood of $p$ in $M$.

\end{itemize}\smallskip
\end{lemma}

Therefore, the final family $A_{ t, s}$ yields a CR-wedge $\mathcal{
W}_p^{ CR, e}$ at $p=A(1)$, as in Definition~3.4. The mild
generalization of the approximation Theorem~5.2(III) stated as
Lemma~2.19 above in the case $d=1$ holds in the general case
$d\geqslant 1$ without modification. Consequently, $\mathcal{
C}_{CR}^0 (M)$ extends to be CR on $\mathcal{ W}_p^{ CR, e}$.

\smallskip
\noindent
{\bf (II)}
Suppose secondly that the exit vector of $A = A_0$ is tangential to
$M$ at $p$. Thanks to the fact that the exit vector mapping has rank
$e$ at $t' = 0$, for every nonzero $t_0'$, the disc $A_{t_0'}$ is
nontangential to $M$ at $p$. In this case, we fix a small $t_0' \neq
0$ and we proceed with $A_{ t'+ t_0'}$ just as above.

\smallskip
In summary, it remains only to construct the family $A_{ t'}$ having
the crucial property that the exit vector mapping has rank $e$ at $t'
= 0$.
\endproof

\subsection*{ 3.16.~Normal deformations of analytic discs}
Thus, we now expose how to construct $A_{ t'}$. We shall introduce a
parameterized family $M_{ t'}$ of $\mathcal{ C}^{ 2, \alpha- 0}$ generic
submanifolds by pushing $M$ near $A (-1)$ inside $\mathcal{ W}_{ A
(-1)}^{ CR, e}$ in $e$ independent normal directions, $e$ being the
number of degrees of freedom offered by $\mathcal{ W}_{ A (-1)}^{ CR,
e}$. Outside a neighborhood of $A(-1)$, each $M_{ t'}$ shall coincides
with $M$ and also $M_{ t'} \vert_{ t' = 0} = M$.

We may assume that the point $p := A(1)$ is the origin in coordinates
$(z, u+ i\, v) \in \C^m \times \C^d$ in which $M$ is represented by $v
= \varphi (z, u)$, where $\varphi$ satisfies $\varphi (0) = 0$ and
$d\, \varphi (0) = 0$. Let $t' \in \R^e$ be small, namely $\vert t'
\vert < t_1'$, with $t_1 ' >0$.

In terms of graphing equations, the deformation $M_{ t'}$ may be
represented by
\[
v
=
\Phi(z,u,t'), 
\]
with $\Phi \in \mathcal{ C}^{ 2, \alpha - 0}$ defined for $\vert t' \vert
< t_1'$ satisfying $\Phi (z,u,0) \equiv \varphi (z, u)$. The point $A
(-1)$ has small coordinates $(z_{-1},u_{-1}+ i\, \varphi ( z_{ -1},
u_{ -1}))$. We require that the $e$ vectors
\[
\Phi_{t_k'}(z_{-1},u_{-1},0), 
\ \ \ \ \ \
k=1,\dots,e,
\]
are linearly independent. There exists a truncated open cone ${\sf C}'
\subset \R^e$ with the property that
\[
M_{t'}
\subset
M\cup\mathcal{W}_{A(-1)}^{CR,e},
\]
whenever $t'\in {\sf C}'$. In fact, we implicitly
assume in Proposition~3.14
that the CR-wedge based
at $A(-1)$ may be 
extended as a $\mathcal{ C}^{ 2, \alpha -0}$ generic submanifold
$\mathcal{ M}_{ A(-1)}^e$ of dimension $2m + d + e$ passing through
$A(-1)$ so that $M_{ t'}$ is contained in $M \cup \mathcal{ M}_{
A(-1)}^e$, for every $\vert t' \vert < t_1'$. The original CR-wedge
$\mathcal{ W}_{ A(-1)}^{ CR, e}$ may then be viewed as a curved real
wedge of edge $M$ which is
contained inside $\mathcal{ M}_{ A(-1)}^{ CR, e}$.

The starting $\mathcal{ C}^{ 2, \alpha}$ disc $A (\zeta) = (Z(\zeta),
W(\zeta))$ with $W(\zeta) = (U (\zeta) + i\, V(\zeta))$ is attached to
$M$ with $A (1) = 0$. Equivalently:
\[
\left\{
\aligned
V(e^{i\,\theta})
&
=
\varphi
\left(
Z(e^{i\,\theta}),U(e^{i\,\theta})
\right),
\\
U(e^{i\,\theta})
&
=
-
{\sf T}_1
\big[
\varphi\big(Z(\cdot),U(\cdot)\big)
\big]
(e^{i\,\theta}),
\endaligned\right.
\]
for every $e^{i\,\theta} \in \partial \Delta$. Thanks to the
existence Theorem~3.7(IV), there exists a $\mathcal{ C}^{ 2, \alpha
-0}$ deformation $A_{ t'}$ of $A$, where each $A_{ t'} (\zeta) :=
\left( Z (\zeta), W(\zeta, t') \right)$ with $A_{ t'} (1) = p$ has
the same $z$-component\footnote{ Since first order partial derivatives
$W_{t_k'} (\zeta, t')$, $k=1, \dots, e$, will appear in a while, we do
not write the parameter $t'$ as a lower index in $U( \zeta, t') + i\,
V (\zeta, t')$. } as $A$ and is attached to $M_{ t'}$, namely:
\def\theequation{3.17}\begin{equation}
\left\{
\aligned
V(e^{i\,\theta},t')
&
=
\Phi
\left(
Z(e^{i\,\theta}),U(e^{i\,\theta},t'),t'
\right),
\\
U(e^{i\,\theta},t')
&
=
-
{\sf T}_1
\left[
\Phi(Z(\cdot),U(\cdot,t'),t')
\right](e^{i\,\theta}),
\endaligned\right.
\end{equation}
for every $e^{i\,\theta} \in \partial \Delta$. Observe that $W( e^{
i\, \theta}, 0) \equiv W ( e^{ i\, \theta})$. We then differentiate
the first line above with respect to $t_k'$ at $t'=0$, for $k=1,
\dots, e$, which yields in matrix notation:
\def\theequation{3.18}\begin{equation}
V_{t_k'}(e^{i\,\theta},0)
=
\Phi_u\left(
Z(e^{i\,\theta}),U(e^{i\,\theta}),0
\right)
U_{t_k'}(e^{i\,\theta},0)
+
\Phi_{t_k'}
\left(
Z(e^{i\,\theta}),U(e^{i\,\theta}),0
\right).
\end{equation}
Also, the $\mathcal{ C}^{ 1, \alpha -0}$ discs $A_{ t_k'} (\zeta, 0)$
satisfy the linear Bishop-type equation
\[
U_{t_k'}(e^{i\,\theta},0)
=
-
{\sf T}_1
\left[
\Phi_u
\left(Z(\cdot),U(\cdot),0
\right)
U_{t_k'}(\cdot,0)
+
\Phi_{t_k'}
(Z(\cdot),U(\cdot),0)
\right]
(e^{i\,\theta}).
\]
As a supplementary space to $T_pM$ in $T_p\C^n$, we choose $H_p := \{
0\} \times i\,\R^d = \{ w=0, u=0\}$. Then ${\sf proj\,}_{ H_p} ( - \partial
A_{ t'} (1) / \partial r ) = - \partial V (1, t') / \partial r$, which
yields after differentiating with respect to $t_k'$ at $t'=0$:
\def\theequation{3.19}\begin{equation}
\left.
\frac{\partial}{\partial t_k'}
\right\vert_{t'=0}
{\sf proj\,}_{H_p}
\left(
-
\frac{\partial A_{t'}}{\partial r}(1)
\right)
=
-
\frac{\partial V_{t_k'}}{\partial r}(1,0),
\end{equation}
for $k=1, \dots, e$. We will establish that if the local deformations
$M_{ t'}$ of $M$ inside the CR-wedge $\mathcal{ W}_{ A(-1)}^{ CR, e}$
are concentrated in a sufficiently thin neighborhood of $A(-1)$, then
the above $e$ vectors $- \partial V_{t_k'} / \partial r (1,0)$, $k= 1,
\dots, e$, are linearly independent. This will complete the proof of
the proposition.

\smallskip
There is a singular integral operator $\mathcal{ J}$ which yields the
interior normal derivative at $1 \in \partial \Delta$ of any
$\mathcal{ C}^{ 1, \alpha -0}$ mapping $v = \overline{ \Delta} \to
\R^d$ which is harmonic in $\Delta$ and vanishes at 
$1 \in \partial \Delta$:
\def\theequation{3.20}\begin{equation}
\mathcal{J}(v)
:=
{\rm p.v.}\,
\frac{1}{\pi}\,\int_{-\pi}^\pi\,
\frac{v(e^{i\,\theta})}{\vert e^{i\,\theta}-1\vert^2}\,d\theta
=
-
\frac{\partial v}{\partial r}(1).
\end{equation}
The proof is postponed to Lemma~3.25 below. If $h: \overline{ \Delta}
\to \C^d$ is $\mathcal{ C}^{ 1, \alpha -0}$ and holomorphic in
$\Delta$, we have in addition
\[
\mathcal{J}(h)
=
-
\frac{\partial h}{\partial r}(1)
=
i\,
\frac{\partial h}{\partial \theta}(1).
\]
With the singular integral $\mathcal{ J}$, we may thus
reformulate~\thetag{ 3.19}:
\[
{\sf proj\,}_{H_p}
\left(
-
\frac{\partial^2 A_0}{\partial t_k'\partial r}(1)
\right)
=
\mathcal{J}(V_{t_k'}).
\]

\def\thelemma{3.21}\begin{lemma}
Let $u, v \in \mathcal{ C}^{ 1, \alpha-0} (\overline{ \Delta}, \R^d)$ be
harmonic in $\Delta$ and vanish at $1 \in \partial \Delta$. Then{\rm
:}
\[
0
=
\mathcal{J}
\left(
u\,v
-
{\sf T}_1u\,{\sf T}_1v
\right).
\]
In addition, $u$ {\rm (}and also $v${\rm )} 
satisfies the two equations{\rm :}
\[
\mathcal{J}(u)
=
-
\frac{\partial ({\sf T}_1u)}{\partial \theta}(1)
\ \ \ \ \
{\rm and}
\ \ \ \ \
\mathcal{J}({\sf T}_1u)
=
\frac{\partial u}{\partial\theta}(1).
\]
\end{lemma}

\proof
The holomorphic
product $w := (u+ i\,{\sf T}_1 u) (v + i\, {\sf T}_1 \, v)$
vanishes to second order at $1 \in \partial \Delta$, so $\mathcal{ J}
(w) = 0$, hence
\[
0
=
{\rm Re}\,
\mathcal{J}(w)
=
\mathcal{J}(u\,v
-
{\sf T}_1u\,{\sf T}_1v).
\]
The pair of equations satisfied by $u$ is obtained by identifying the
real and imaginary parts of $\mathcal{ J} (h) = i\, \frac{\partial
h}{\partial \theta} (1)$, where $h := u + i\, {\sf T}_1 u$.
\endproof

Following~\cite{ tu1994a},
we now introduce a $d\times d$ matrix $G$ of $\mathcal{ C}^{ 1,
\alpha}$ functions on $\partial \Delta$ defined by the functional
equation
\[
G(e^{i\,\theta})
=
I
+
{\sf T}_1
\left[
G(\cdot)\,\Phi_u
\left(
Z(\cdot),U(\cdot),0
\right)
\right](e^{i\,\theta}).
\]
Here $\Phi_u = (\Phi_{ u_l}^j )_{1\leqslant l\leqslant d}^{ 1\leqslant j\leqslant d}$ is a
$d\times d$ matrix. Since $\Phi_u (z, u, 0) \equiv \varphi_u (z, u)$
is small, the solution $G$ exists and is unique, by an application of
Proposition~3.21(IV). Notice that $G (1) = I$. Applying ${\sf T}_1$
to both sides, we get ${\sf T}_1\,G = - G\,\Phi_u + {\rm cst.}$,
without writing the arguments. In fact, the constant vanishes, since
$\Phi_u (0, 0, 0) = \varphi_u (0, 0) = 0$. So we get:
\[
{\sf T}_1\,G
=
-
G\,\Phi_u.
\]
We also notice that $V_{ t_k'} = {\sf T}_1 \, U_{ t_k'}$
and $U_{ t_k'} = - {\sf T}_1 \, V_{ t_k'}$.

Next, we rewrite~\thetag{ 3.18} without arguments: $\Phi_{ t_k'} = V_{
t_k'} - \Phi_u \, U_{ t_k'}$, $k=1, \dots, e$, we apply the matrix $G$
to both sides, we replace $G\, \Phi_u$ by $-{\sf T}_1 G$ as well as
$U_{ t_k'}$ by $- {\sf T}_1 V_{ t_k'}$ and we let appear a term $u\, v
- {\sf T}_1 u \, {\sf T}_1 v$:
\[
\aligned
G\,\Phi_{t_k'}
&
=
G\,V_{t_k'}
-
G\,\Phi_u\,U_{t_k'}
\\
&
=
G\,V_{t_k'}
-
({\sf T}_1G)({\sf T}_1V_{t_k'})
\\
&
=
V_{t_k'}
+
(G-I)V_{t_k'}
-
{\sf T}_1
(G-I)\,{\sf T}_1V_{t_k'}.
\endaligned
\]
Finally\footnote{
We can also check that $\mathcal{ J} (U_{ t_k'}) = - \mathcal{ J}
({\sf T}_1 V_{ t_k'}) = \partial V_{ t_k' }(1,0) / \partial \theta =
0$. Indeed, it suffices to differentiate~\thetag{ 3.18} with respect to
$\theta$ at $\theta = 0$, noticing that $\Phi_u (0, 0, 0) = \varphi_u
(0, 0) = 0$, that $U_{ t_k'} (1, 0) = 0$ and that $\Phi_{ t_k'} (z, u,
0) = 0$ for $(z, u)$ near $(0, 0)$.
}, 
applying the singular operator $\mathcal{ J}$ 
and remembering Lemma~3.21, we obtain:
\def\theequation{3.22}\begin{equation}
\mathcal{J}(G\,\Phi_{t_k'})
=
\mathcal{J}(V_{t_k'}).
\end{equation}
We claim that if the support of the deformation $M_{ t'}$ is
sufficiently concentrated near $A(-1)$, the $e$ vectors $\mathcal{ J}
(V_{ t_k'}) = \mathcal{ J} ( G\, \Phi_{ t_k'}) \in \R^d$ are linearly
independent.

Indeed, since the deformations $M_{ t'}$ are
localized near $A(-1)$, we have
$\Phi_{ t_k'} (Z (e^{ i\, \theta}), U( e^{i\, \theta}), 0 )
\equiv 0$, except for $\vert \theta + \pi \vert < \theta_2$, with
$\theta_2 >0$ small. We deduce:
\def\theequation{3.23}\begin{equation}
\aligned
\mathcal{J}(G\,\Phi_{t_k'})
&
=
\frac{1}{\pi}\,
\int_{\vert\theta+\pi\vert<\theta_2}
\frac{G(e^{i\,\theta})\,\Phi_{t_k'}
\big(Z(e^{i\,\theta}),
U(e^{i\,\theta}),0\big)}{
\vert e^{i\,\theta}-1\vert^2}\,
d\theta
\\
&
\approx
\frac{1}{\pi}\,
\frac{G(-1)}{4}\,
\int_{\vert\theta+\pi\vert<\theta_2}\,
\Phi_{t_k'}
\big(
Z(e^{i\,\theta}),U(e^{i\,\theta}),0
\big)\,
d\theta.
\endaligned
\end{equation}
Since, by assumption, the $e$ vectors $\Phi_{ t_k'}
(z_{ -1}, u_{ -1}, 0)$ are linearly independent, the linear independence 
of the above (concentrated) vector-valued 
integrals follows.

The proofs of Proposition~3.14 and of Theorem~3.13 are complete.
\qed

\subsection*{ 3.24.~Proofs of two lemmas}
Firstly, we check formula~\thetag{ 3.20}.

\def\thelemma{3.25}\begin{lemma}
Let $u\in \mathcal{C}^{ 1, \beta} (\overline{ \Delta})$ $(0 < \beta
< 1)$ be harmonic in $\Delta$, real-valued and satisfying $u(1) =
0$. Then the interior normal derivative of $u$ at $1\in \partial
\Delta$ is given by{\rm :}
\[
-
\frac{\partial u}{\partial r}(1)
=
{\rm p.v.}\,
\frac{1}{\pi}\,
\int_{-\pi}^\pi\,
\frac{u(e^{i\,\theta})}{
\vert e^{i\,\theta}-1\vert^2}\,d\theta
=
{\rm p.v.}\,
\frac{i}{\pi}\,
\int_{\partial\Delta}\,
\frac{u(\zeta)}{(\zeta-1)^2}\,d\zeta.
\]
\end{lemma}

\proof
The function $h := u + i\, {\sf T} u$ is holomorphic in $\Delta$ and
$\mathcal{ C}^{ 1, \beta}$ in $\overline{ \Delta}$. Since ${\sf T}
u$ is also harmonic in $\Delta$, since $\frac{ \partial h}{\partial r}
(1) = \frac{ \partial u }{\partial r} (1) + i\, \frac{ \partial {\sf
T} u}{\partial r} (1)$, and since the kernel $\vert e^{ i\, \theta} -
1 \vert^{ -2}$ is real, we may prove the lemma with $u$ replaced by $h
\in \mathcal{ O} (\Delta) \cap \mathcal{ C}^{ 1, \beta} (\overline{
\Delta})$.

Let $\zeta = r e^{i\, \theta}$ and denote $h_1 := \frac{\partial
h}{\partial \zeta} (1) = \frac{\partial h}{\partial r} (1)$, so that
$h( \zeta) = (\zeta -1) h_1 + {\rm O} (\vert \zeta - 1\vert^{ 1+
\beta})$. We remind that, for any $\zeta_0 \in \partial \Delta$, by
an elementary modification
of Cauchy's formula, we have ${\rm p.v.}\,\frac{1}{2\pi
i}\, \int_{\partial\Delta}\,\frac{d\zeta}{\zeta-\zeta_0} = \frac{
1}{2}$. We deduce that the linear term $(\zeta -1) h_1$ provides the
main contribution:
\[
{\rm p.v.}\,
\frac{i}{\pi}\,\int_{
\partial\Delta}\,
\frac{(\zeta-1)h_1}{
(\zeta-1)^2}\,d\zeta
=
-
2\,h_1\,{\rm p.v.}\,\frac{1}{2\pi i}\,
\int_{\partial\Delta}\,
\frac{d\zeta}{\zeta-1}
=
-
h_1.
\]
Thus, we have to prove that the remainder $r(\zeta) := h(\zeta) -
(\zeta -1) h_1$, which belongs to $\mathcal{ O} (\Delta) \cap
\mathcal{ C}^{ 1, \beta} (\overline{ \Delta})$, gives no contribution,
namely satisfies $\int_{ \partial \Delta} \, \frac{ r(\zeta)}{ (\zeta
- 1)^2}\, d\zeta = 0$.

Set $s(\zeta) := \frac{ r(\zeta)}{ (\zeta - 1)^2}$. Then $s\in
\mathcal{ O} (\Delta)$ is continuous on $\overline{ \Delta} \backslash
\{ 1\}$ and satisfies $\vert s(\zeta) \vert \leqslant {\sf K} \, \vert
\zeta - 1\vert^{ \beta -1}$ for some ${\sf K} >0$. We claim that by an
application of Cauchy's theorem, the integral $\int_{ \partial \Delta}
\, s(\zeta) \, d\zeta$, which exists without principal value,
vanishes.

Indeed, let $\varepsilon$ with $0 < \varepsilon << 1$ and consider the
open disc $\Delta (1, \varepsilon)$ of radius $\varepsilon$ centered
at $1$. The drawing of this disc delineates three arcs of $\overline{
\Delta}$:

\begin{itemize}

\smallskip\item[{\bf (i)}]
the open arc $\partial \Delta \backslash
\overline{ \Delta (1, \varepsilon)}$, of length 
$\approx 2\pi - 2\, \varepsilon$;

\smallskip\item[{\bf (ii)}] 
the closed arc $\partial \Delta
\cap \overline{ \Delta (1,
\varepsilon)}$, of length $\approx 2\,
\varepsilon$;

\smallskip\item[{\bf (iii)}]
the closed arc $\partial \Delta (1, \varepsilon) \cap
\overline{ \Delta}$, of length is $\approx \pi \varepsilon$.

\end{itemize}\smallskip

\vskip -2cm

\hskip 9cm
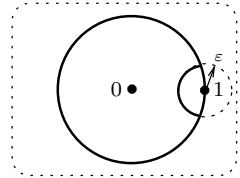
\begin{picture}(0,0)%
\includegraphics{three-arcs.pstex}%
\end{picture}%
\setlength{\unitlength}{4144sp}%
\begingroup\makeatletter\ifx\SetFigFont\undefined
\def\x#1#2#3#4#5#6#7\relax{\def\x{#1#2#3#4#5#6}}%
\expandafter\x\fmtname xxxxxx\relax \def\y{splain}%
\ifx\x\y   
\gdef\SetFigFont#1#2#3{%
  \ifnum #1<17\tiny\else \ifnum #1<20\small\else
  \ifnum #1<24\normalsize\else \ifnum #1<29\large\else
  \ifnum #1<34\Large\else \ifnum #1<41\LARGE\else
     \huge\fi\fi\fi\fi\fi\fi
  \csname #3\endcsname}%
\else
\gdef\SetFigFont#1#2#3{\begingroup
  \count@#1\relax \ifnum 25<\count@\count@25\fi
  \def\x{\endgroup\@setsize\SetFigFont{#2pt}}%
  \expandafter\x
    \csname \romannumeral\the\count@ pt\expandafter\endcsname
    \csname @\romannumeral\the\count@ pt\endcsname
  \csname #3\endcsname}%
\fi
\fi\endgroup
\begin{picture}(1367,1055)(442,-642)
\put(1044,-152){\makebox(0,0)[lb]{\smash{\SetFigFont{8}{9.6}{rm}{\color[rgb]{0,0,0}$0$}%
}}}
\put(1655,-156){\makebox(0,0)[lb]{\smash{\SetFigFont{8}{9.6}{rm}{\color[rgb]{0,0,0}$1$}%
}}}
\put(1662, 56){\makebox(0,0)[lb]{\smash{\SetFigFont{6}{7.2}{rm}{\color[rgb]{0,0,0}$\varepsilon$}%
}}}
\end{picture}

We then decompose the integral of $s$ on $\partial \Delta$ as
integrals on the first two arcs:
\[
\int_{\partial\Delta}\,s(\zeta)\,d\zeta
=
\int_{\partial\Delta\backslash\overline{\Delta(1,\varepsilon)}}\,
s(\zeta)\,d\zeta
+
\int_{\partial\Delta\cap\overline{\Delta(1,\varepsilon)}}\,
s(\zeta)\,d\zeta.
\]
The estimate $\vert s (\zeta) \vert \leqslant {\sf K} \, \vert \zeta - 1
\vert^{ \beta -1}$ insures the smallness of the second integral:
\[
\left\vert
\int_{\partial\Delta\cap\overline{\Delta(1,\varepsilon)}}\,
s(\zeta)\,d\zeta
\right\vert
\leqslant C_1\, \varepsilon^\beta.
\]
To transform the first integral, we observe that Cauchy's theorem
entails that integration of $s(\zeta) \, d\zeta$ on the closed contour
$\big[ \partial \Delta \backslash \overline{ \Delta (1, \varepsilon)}
\big] \cup \left[ \partial \Delta (1, \varepsilon) \cap \overline{
\Delta} \right]$ vanishes:
\[
0
=
\int_{\partial\Delta\backslash\overline{\Delta(1,\varepsilon)}}\,
s(\zeta)\,d\zeta
+
\int_{\partial\Delta(1,\varepsilon)\cap\overline{\Delta}}\,
s(\zeta)\,d\zeta.
\]
Hence the first integral $\int_{\partial \Delta \backslash \overline{
\Delta(1, \varepsilon) }}$ may be replaced by the integral $- \int_{
\partial \Delta (1, \varepsilon) \cap \overline{ \Delta}}$ on the
third arc. The estimate $\vert s (\zeta) \vert \leqslant {\sf K} \, \vert
\zeta - 1 \vert^{ \beta -1}$ again insures that this second integral
is bounded by $C_2 \, \varepsilon^\beta$. In conclusion $\vert
\int_{\partial \Delta} \, s(\zeta) \, d\zeta \vert \leqslant (C_1 + C_2) \,
\varepsilon^\beta$. 
\endproof

\proof[Proof of Lemma~3.15.]
Secondly, we provide the details for the translation of the family
$A_{ t'}$. Let $v = \varphi (z, u)$ represent $M$ in a neighborhood of
$p$. By assumption, $A_{ t'} (\zeta) = ( Z(\zeta), W (\zeta, t'))$ is
attached to $M_{ t'}$, with $A_{ t'} (1) = p$. Equivalently, the two
equations~\thetag{ 3.17} hold. Since $A = A_{ t'} \vert_{ t' = 0}$ is
an embedding, the vector $v_p := \frac{ \partial A}{\partial \theta}
(1) \in T_p M$ is nonzero. As in \S2.12, we choose a small $(2m+ d -
1)$-dimensional submanifold $K_p$ passing through $p$ with $\R v_p
\oplus T_p K_p = T_p M$ and we parametrize it by $s \mapsto (z(s),
u(s)+ i\, \varphi (z(s), u(s)))$, where $s\in \R^{ 2m+ d - 1}$ is
small, $\vert s\vert < s_1$, $s_1>0$. Then the translation
\[
A_{ t', s} (\zeta)
= 
\left(
Z(\zeta)+z(s), 
W(\zeta,t',s)\right)
\]
is constructed by perturbing the two equations~\thetag{ 3.17},
requiring only that
\[
A_{t',s}(1) 
= 
(z(s),u(s)+i\,\varphi(z(s),u(s))).
\]
This is easily done:
\[
\left\{
\aligned
V(e^{i\,\theta},t',s)
&
=
\Phi\big(
Z(e^{i\,\theta})+z(s),U(e^{i\,\theta},t',s),t'
\big),
\\
U(e^{i\,\theta},t',s)
&
=
u(s)
-
{\sf T}_1
\left[
\Phi
\left(
Z(\cdot)+z(s),
U(\cdot,t',s),t'
\right)
\right](e^{i\,\theta}).
\endaligned\right.
\]
The non-tangency of $v_p$ with $K_p$ at $p$ then insures that for every
small fixed $t'$, the mapping $(\theta, s) \mapsto A_{ t', s} (e^{ i\,
\theta})$ is a diffeomorphism onto a neighborhood of 
$p$ in $M$.
\endproof

\section*{\S4.~Holomorphic extension \\ 
on globally minimal generic submanifolds}

\subsection*{4.1.~Structure of the proof of Theorem~3.8} 
Let $M$ be a $\mathcal{ C}^{ 2, \alpha}$ globally minimal generic
submanifold of $\C^n$. For clarity, we begin by a summary of the main
steps of the proof of Theorem~3.8.

\begin{itemize}

\smallskip\item[{\bf (a)}]
Since $M$ is globally minimal, the distribution $q \mapsto T_q^cM$
must be somewhere not involutive, namely there must 
exist a point $p\in M$
and a section $L$ of $T^{ 1, 0} M$ defined in an open neighborhood
$U_p$ of $p$ in $M$ with $L(p) \neq 0$ such that $\left[ L, \overline{
L} \right] (p) \not \in T_p^{ 1, 0} M \oplus T_p^{ 0, 1} M$.

\smallskip\item[{\bf (b)}]
Thanks to an easy generalization of the Lewy extension theorem
(\S2.10), there exists a manifold $M_p^1$ attached to $M$ at $p$ with
$\dim M_p^1 = 1 + \dim M$ to which $\mathcal{ C}_{ CR}^0 (M)$ extends
to be CR.

\smallskip\item[{\bf (c)}]
Thanks to the main propagation Proposition~3.14, CR extension to a
similar manifold $M_q^1$ attached to $M$ holds at every point $q\in M =
\mathcal{ O}_{ CR } (M, p)$.

\smallskip\item[{\bf (d)}]
Since there are as many manifolds with boundary as points in $M$, it
may well happen that at some point $p \in M$ which belongs to the edge
of {\it two different}\, manifolds $M_{ p '}^1$ and $M_{ p''}^1$, the
tangent spaces $T_p M_{ p'}^1$ and $T_p M_{ p''}^1$ are distinct.
Refering to the diagram of \S4.5 below, we may then immediately profit
of such a situation, if it occurs.

\smallskip\item[{\bf (e)}]
Indeed, in this case, an appropriate version of the edge-of-the-wedge
theorem guarantees that $\mathcal{ C}_{ CR}^0 (M)$ extends to be CR on
a $\mathcal{ C}^{ 2, \alpha -0}$ CR-wedge $\mathcal{ W}_p^{ CR, e}$ at
$p$ whose dimension $e$ is $\geqslant 1+1 = 2$.

\smallskip\item[{\bf (f)}]
To reason abstractly, let $e_{ \rm max}$ be the maximal integer $e$
with $1\leqslant e \leqslant d$ such that there exists a point $p\in M$ and a
$\mathcal{ C}^{ 2, \alpha -0}$ CR-wedge $\mathcal{ W}_p^{ CR, e
}$ at $p$ of dimension $2m+ d+ e$ to which $\mathcal{
C}_{ CR}^0 (M)$ extends to be CR. Thanks to the main propagation
Proposition~3.14, CR extension to a $\mathcal{ C}^{ 2, \alpha -0}$
CR-wedge $\mathcal{ W}_q^{ CR, e_{\rm max}}$ holds at every point $q
\in M = \mathcal{ O}_{ CR} (M, p)$.

\smallskip\item[{\bf (g)}]
If $e_{ \rm max} = d$, we are done, Theorem~3.8 is proved. Assuming
$e_{ \rm max} \leqslant d-1$, we must construct a contradiction
in order to complete the proof.

\smallskip\item[{\bf (h)}]
Since $e_{ \rm max}$ is maximal, again because of the
edge-of-the-wedge theorem, the transversal situation {\bf (d)} cannot
occur; in other words, every point $p\in M$ that belongs to the edges
of two different CR-wedges $\mathcal{ W}_{p'}^{ CR, e_{\rm max}}$ and
$\mathcal{ W}_{ p''}^{ CR, e_{\rm max}}$ has the property that $T_p
\mathcal{ W}_{p'}^{ CR, e_{\rm max}} = T_p \mathcal{ W}_{ p''}^{ CR,
e_{\rm max}}$.

\smallskip\item[{\bf (i)}]
It follows that, as $p$ runs in $M$, the
$(2m + d+ e)$-dimensional tangent planes $T_p \mathcal{
W}_p^{ CR, e_{ \rm max}}\cap T_pM$ glue together and they define a
$\mathcal{ C}^{ 1, \alpha- 0}$ sub-distribution $KM$ of the tangent
bundle $TM$, of dimension $2 m + e_{ \rm max}$, which contains $T^cM$.

\smallskip\item[{\bf (j)}]
Since $M$ is globally minimal, such a distribution $p \mapsto KM(p)$ must
be somewhere not involutive, namely there must 
exist a point $p\in M$ such
that $\left[ KM, KM\right] (p ) \not \subset KM (p)$.

\smallskip\item[{\bf (k)}]
The $\mathcal{ C}^{ 2, \alpha-0}$ CR-wedge $\mathcal{ W}_p^{ CR,
e_{\rm max}}$ may be included in some $\mathcal{ C}^{ 2, \alpha - 0}$
local generic submanifold $\mathcal{ M}_p^{ e_{ \rm max}}$ passing
through $p$ and containing $M$ in a neighborhood of $p$.

\smallskip\item[{\bf (l)}]
Multiplication by $i$ gives $T_p^c \mathcal{ M}_p^{ e_{\rm max}} =
KM(p) + i\, KM(p)$ and the nondegeneracy $\left[ KM, KM\right] (p )
\not \subset KM (p)$ implies that the Levi-form of $\mathcal{ M}_p^{
e_{ \rm max}}$ is not identically zero at $p$, namely $\left[T_p^c
\mathcal{ M}_p^{ e_{\rm max}}, T_p^c \mathcal{ M}_p^{ e_{\rm
max}} \right] (p) \not \subset T_p^c \mathcal{ M}_p^{ e_{\rm
max}}$.

\smallskip\item[{\bf (m)}]
Then a version of the Lewy-extension theorem on conic generic
manifolds having a generic edge guarantees that $\mathcal{ C}_{ CR}^0
(M)$ extends to be CR on a CR-wedge $\widetilde{ \mathcal{ W}}_p^{ CR,
1+ e_{ \rm max} }$ of dimension $2m + d + 1 + e_{ \rm max}$ at
$p$. This new CR-wedge is constructed by means of discs attached to $M
\cup \mathcal{ W}_p^{ CR, e_{\rm max}}$, exploiting the nondegeneracy
of the Levi form of $\mathcal{ M}_p^{ e_{\rm max}}$. This
contradicts the assumption that $e_{ \rm max} \leqslant d-1$ was maximal,
hence completes the proof of Theorem~3.8.

\end{itemize}\smallskip

The remainder of Section~4 is devoted to provide all the details of
the proof.

\subsection*{ 4.2.~Lewy extension in arbitrary codimension}
As observed in {\bf (a)} above, there exists a point $p\in M$ and a
local section $L$ of $T^{ 1, 0} M$ with $L (p) \neq 0$ such that
$\left[ L, \overline{ L} \right] (p) \not \subset \C\otimes T_p^cM$.

\def\thelemma{4.3}\begin{lemma}
{\rm (\cite{ we1982, bpo1982})}
There exists a manifold with boundary $M_p^1$ attached to a
neighborhood of $p$ in $M$ with $\dim \, M_p^1 = 1+ \dim \, M$ to
which $\mathcal{ C}_{ CR}^0 (M)$ extends to be CR.
\end{lemma}

We shall content ourselves with only one direction of extension, since
this will be sufficient for the sequel. Nevertheless, we mention that
finer results expressed in terms of the Levi-cone of $M$ at $p$ may be
found in~\cite{ bpo1982, bo1991}. 
Anyway, all the extension results that are
based on pointwise nondegeneracy conditions as the openness of
Levi-cone or the finite typeness of $M$ at a point are by far less
general than Theorem~3.8, in which propagational aspects are
involved.

\proof
The arguments are an almost straightforward generalization of the
proof of the Lewy extension theorem (hypersurface case), already
exposed in \S2.10 above. Here is a summary.

By linear algebra reasonings, we may find local coordinates $(z, w)
\in \C^m \times \C^d$ vanishing at $p$ with $L (p) =
\frac{\partial }{\partial z_1} \big\vert_p$, with $M$ given by
$v = \varphi (z, u)$, where $\varphi (0) = 0$, $d\varphi ( 0) =0$, and
with first equation given by
\[
v_1 
= 
\varphi_1
=
z_1\bar z_1
+
{\rm O}(\vert z_1\vert^{2+\alpha})
+
{\rm O}(\vert \widetilde{ z}\vert)
+
{\rm O}(\vert z\vert\,\vert u\vert)
+
{\rm O}(\vert u\vert^2),
\]
where we have split further the coordinates as $(z_1, \widetilde{
z}, w_1, \widetilde{ w})$, with $\widetilde{ z} \in \C^{ m-1}$ and
$\widetilde{ w} \in \C^{ d-1}$. For $\varepsilon >0$ small, we
introduce the disc defined by
\[
A_\varepsilon(\zeta)
:=
\left(
\varepsilon(1-\zeta),
\widetilde{0},
W_\varepsilon^1(\zeta),
\widetilde{W}_\varepsilon(\zeta)
\right),
\]
where $W_\varepsilon (\zeta) = U_\varepsilon (\zeta) + i\,
V_\varepsilon (\zeta)$ is uniquely defined by requiring that
$A_\varepsilon$ is attached to $M$ and satisfies $A_\varepsilon (1) =
p$. As in \S2.10, one verifies that
\[
-
\frac{\partial V_\varepsilon^1}{\partial r}(1)
=
2\,\varepsilon^2
+
{\rm O}(\varepsilon^{2+\alpha}).
\]
Hence the exit vector of $A_\varepsilon$ at $1 \in \partial \Delta$
is nontangential to $M$ at $p$, provided 
$\varepsilon >0$ is small enough and fixed. By
translating $A_\varepsilon$, we construct the desired manifold with
boundary $M_p^1$.
\endproof

\subsection*{ 4.4.~Maximal dimension for CR extension}
As in \S4.1{\bf (f)}, let $e_{ \rm max}$ be the maximal integer $e
\leqslant d$ such that there exists a point $p\in M$ and a $\mathcal{ C}^{
2, \alpha - 0}$ CR-wedge $\mathcal{ W}_p^{ CR, e}$ at $p$ of dimension
$2m + d + e$ to which $\mathcal{ C}_{ CR}^0 (M)$ extends to be CR. By
the above Lewy extension, we have $e_{ \rm max} \geqslant 1$. Thanks to
the main propagation Proposition~3.14, it immediately follows that CR
extension to a $\mathcal{ C}^{ 2, \alpha -0}$ CR-wedge $\mathcal{
W}_q^{ CR, e_{\rm max}}$ holds at every point $q \in M = \mathcal{
O}_{ CR} (M, p)$. If $e_{\rm max} = d$, Theorem~3.8 is proved,
gratuitously.

Assuming that $1\leqslant e_{ \rm max} \leqslant d-1$, in order to establish
Theorem~3.8, we must construct a contradiction. In the sequel, we
shall simply denote $e_{ \rm max}$ by $e$.

To proceed further, we must reformulate with high precision how were
constructed all the CR-wedges obtained by the propagation
Proposition~3.14.

For every point $p\in M$, there exists a local CR-wedge $\mathcal{
W}_p^{ CR, e}$ attached to a neighborhood of $p$ in $M$ which is
described by means of a family of analytic discs $A_{ p, t, s}
(\zeta)$, where $t$ and $s$ are parameters. Here, the subscript $p$ is
not a parameter, it indicates only that $p$ is the base point of $A_{
p, t, s}$, namely $A_{ p, t, 0} (1) = p$. The family $A_{ p, t, s}$
enjoys properties that are listed below. In this list, the conditions
are more uniform than those formulated in Definition~3.4, but one
immediately verifies that both formulations are equivalent, up to a
shrinking of $t_1 (p) >0$, of $s_1 (p) >0$, of $\theta_1(p) >0$ and of
$1 - r_1(p) >0$.

\begin{itemize}

\smallskip\item[$\bullet$]
The rotation parameter $t \in \R^{ e-1}$ runs in $\{ \vert t\vert <
t_1 (p) \}$, for some small $t_1 (p) >0$.

\smallskip\item[$\bullet$]
The translation parameter $s \in \R^{ 2m + d - 1}$ runs in 
$\{ \vert s\vert < s_1(p)\}$, for some small $s_1 (p) >0$.

\smallskip\item[$\bullet$]
The point $q(p) := A_{ p, 0, 0} (-1) \in M$ is close to $p$.

\smallskip\item[$\bullet$]
At $q(p)$, there is a CR-wedge $\mathcal{ W}_{ q(p)}^{ CR, e}$.

\smallskip\item[$\bullet$]
The family $A_{ p, t, s}$ satisfies $A_{ p, t, s} (\partial \Delta)
\subset M \cup \mathcal{ W}_{ q(p)}^{ CR, e}$.

\smallskip\item[$\bullet$]
A small angle $\theta_1 (p) > 0$ and a radius $r_1 (p) >0$ 
close to $1$ are chosen.

\smallskip\item[$\bullet$] 
A family $H_{ p'}$ of linear subspaces of $T_{ p'}\C^n$
satisfying $T_{ p'} H_{ p'} \oplus T_{ p'} M = T_{ p'} \C^n$ for all
$p'\in M$ in a neighborhood of $p$ is chosen.

\smallskip\item[$\bullet$] 
For every $t$ with $\vert t\vert < t_1 (p)$, every $s$ with $\vert
s\vert < s_1 (p)$ and every $\theta$ with $\vert \theta \vert <
\theta_1 (p)$, the exit vector of $A_{ p, t, s} (e^{ i\, \theta})$ at
$e^{i\,\theta}$ is not tangent to $M$:
\[
\text{\sf ex}(A_{p,t,s})(e^{i\,\theta})
:=
{\sf proj\,}_{H_{A_{p,t,s}(e^{i\,\theta})}}
\left(
i\,
\frac{\partial A_{p,t,s}}{\partial \theta}(e^{i\,\theta})
\right)
\neq 
0.
\]

\smallskip\item[$\bullet$] 
For every fixed $s$ with $\vert s\vert < s_1 (p)$ and every fixed
$\theta$ with $\vert \theta\vert < \theta_1 (p)$, the normalized exit
vector mapping
\[
\R^{e-1}
\ni\,t\
\longmapsto\
\text{\rm n-ex}
(A_{p,t,s})(e^{i\,\theta})
\in S^{d-1}
\]
is of rank $(e-1)$ at every $t \in 
\{ \vert t\vert < t_1 (p)\}$.

\smallskip\item[$\bullet$]
For some $t_2 (p)$, $s_2 (p)$, $\theta_2(p)$ and $r_2(p)$ satisfying
$0 < t_2 (p) < t_1(p)$, $0 < s_2 (p) < s_1 (p)$, $0 < \theta_2 (p) <
\theta_1 (p)$ and $0 < 1 - r_2 (p) < 1 - r_1 (p) < 1$, the CR-wedge is
precisely defined as:
\[
\mathcal{W}_p^{CR,e}
:=
\left\{
A_{p,t,s}(re^{i\,\theta}):\
\vert t\vert<t_2(p),\
\vert s\vert<s_2(p),\
\vert \theta\vert<\theta_2(p),\
r_2(p)<r<1
\right\}.
\]

\smallskip\item[$\bullet$]
Finally, the CR-wedge $\mathcal{ W}_p^{ CR, e}$ is contained in a
$\mathcal{ C}^{ 2, \alpha - 0}$ local generic submanifold $\mathcal{
M}_p^e$ of the same dimension $2m + d +e$ that contains a neighborhood
of $p$ in $M$. At a point $p ' = A_{ p, t', s'} (e^{ i\, \theta'})
\in M$ of the edge of $\mathcal{ W}_p^{ CR, e}$, the tangent space of
$\mathcal{M }_p^e$
is:
\[
T_{p'}\mathcal{M}_p^e
=
T_pM
\oplus\,
\R\left(
i\,
\frac{\partial A_{p,t',s'}}{\partial\theta}
(e^{i\,\theta'})
\right)
\bigoplus_{1\leqslant k\leqslant e-1}
\R
\left(
i\,
\frac{\partial^2A_{p,t',s'}}{\partial\theta\partial t_k}
(e^{i\,\theta'})
\right).
\]

\end{itemize}\smallskip

\subsection*{4.5.~An edge-of-the-wedge theorem}
There are as many generic submanifolds $\mathcal{ M}_{ p'}^{ CR, e}$
of codimension $d-e$ as points $p' \in M$. At a point $p = A_{ p', t',
s'} (e^{ i\, \theta'})$ that belongs to the edge of such an $\mathcal{
M}_{ p'}^{ CR, e}$, we may define a linear subspace of $T_p M$ by
\[
KM_{p'}(p)
:=
T_p^c\mathcal{M}_{p'}^e\cap
T_pM.
\]
Since $\mathcal{ M}_{p'}^e$ is generic and contains $M$ in a
neighborhood of $p$, this space $KM_{ p'}(p)$ contains $T_p^cM$ and is
$(2m+ e)$-dimensional. Also, multiplication by $i$ induces an
isomorphism $KM_{ p'}( p)/ T_p^c M \simeq T_p \mathcal{ M}_{ p'}^e
/T_pM $.

In general, two different $KM_{ p'} (p)$ and $KM_{ p''} (p)$ need 
not coincide, or equivalently, two different
tangent spaces $T_p \mathcal{ M}_{ p'}^e$
and $T_p \mathcal{ M}_{ p''}^e$ are unequal.

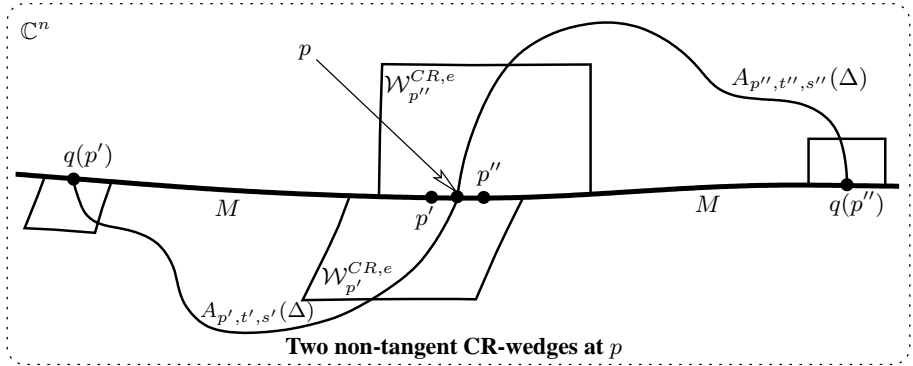
\begin{center}
\input transversal-CR-wedge.pstex_t
\end{center}

\noindent
More precisely, there is a dichotomy.

\begin{itemize}\smallskip\item[{\bf (I)}]
Either for every two points $p', p''\in M$ such that there exists a
point $p$ belonging to the intersection of the edges of the two
CR-wedges $\mathcal{ W}_{ p'}^{ CR, e}$ and $\mathcal{ W}_{ p''}^{ CR,
e}$, namely of the form:
\[
p
=
A_{p',t',s'}(e^{i\,\theta'})
=
A_{p'',t'',s''}(e^{i\,\theta''}),
\]
for some values 
\[
\aligned
\vert t'\vert
&
<
t_2(p'),
\ \ \ \ \ \
\vert s'\vert
<
s_2(p'),
\ \ \ \ \ \
\vert\theta'\vert
<
\theta_2(p'),
\\
\vert t''\vert
&
<
t_2(p''),
\ \ \ \ \
\vert s''\vert
<
s_2(p''),
\ \ \ \ \ 
\vert\theta''\vert
<
\theta_2(p''),
\endaligned
\]
the two spaces $T_p \mathcal{ M}_{ p'}^e$ and $T_p \mathcal{ M}_{
p''}^e$ coincide. Equivalently, $KM_{ p'} (p) = KM_{ p''} (p)$.

\smallskip\item[{\bf (II)}]
Or there exist two points $p', p'' \in M$ and a point $p = A_{p',t',
s'}(e^{ i\,\theta'}) = A_{p'', t'',s''}(e^{i \,\theta'' })$ in the
intersection of the edges of the two CR-wedges $\mathcal{ W}_{ p'}^{
CR, e}$ and $\mathcal{ W}_{ p''}^{ CR, e}$ such that
\[
T_p\mathcal{ M}_{p'}^e
\neq
T_p\mathcal{ M}_{p''}^e.
\]
\end{itemize}\smallskip

\def\thelemma{4.6}\begin{lemma}
The case $T_p \mathcal{ M}_{ p'}^e \neq T_p\mathcal{ M}_{ p''}^e$
implies that $\mathcal{ C}_{ CR}^0 (M)$ extends to be CR on a CR-wedge
$\widetilde{ \mathcal{ W}}_p^{CR, 1+ e}$ at $p$
whose dimension equals $2m + d + 1+ e$,
contradicting the maximality of $e = e_{ \rm max}$.
\end{lemma}

Of course, this lemma follows by a known CR version of the
edge-of-the-wedge theorem (\cite{ ai1989}), but for completeness, we
summarize a shorter proof that exploits the existence of the discs
$A_{ p', t, s}$, as in~\cite{ po2004}.

\proof
By construction, the family $A_{ p', t, s} (\zeta)$ covers the
CR-wedge $\mathcal{ W}_{ p'}^{ CR, e}$. The point $p$ belongs to the
edge of $\mathcal{ W}_{ p'}^{ CR, e}$. 

Since $T_p \mathcal{ M}_{ p'}^e \neq T_p\mathcal{ M}_{ p''}^e$, there
exists a manifold $M_p^1 \subset \mathcal{ W}_{ p''}^{ CR, e}$
attached to $M$ at $p$ with $\dim \, M_p^1 = 1 + \dim \, M$ such that
\[
1+e
=
\dim
\left(
\left[
T_pM_p^1
+
T_p\mathcal{W}_{p'}^{CR,e}
\right]\Big/
T_pM
\right).
\]

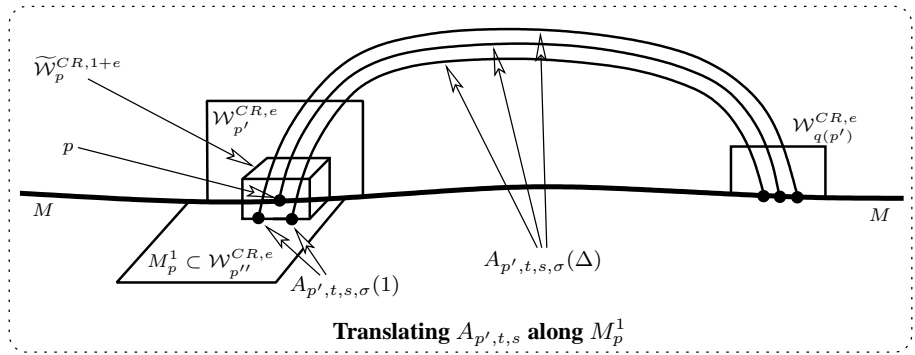
\begin{center}
\input CR-edge-of-the-wedge.pstex_t
\end{center}

We may deform the family $A_{ p', t, s}$ by translating it along
$M_p^1$, as in the diagram. So we introduce a supplementary parameter
$\sigma >0$ and we require that the point $A_{ p', t, s, \sigma} (1)$
should cover a one-sided neighborhood of $p$ in $M_p^1$ as $\sigma$
runs in $(0, \sigma_1)$, for some small $\sigma_1 >0$, and as the
previous translation parameter $s \in \R^{ 2m+ d -1}$ runs in
$\{ \vert s\vert < s_2 (p') \}$. Thanks to Theorem~3.7(IV), the
corresponding Bishop-type equation has $\mathcal{ C}^{ 2, \alpha - 0}$
solutions.

If we choose $t_3 >0$ with $\vert t' \vert + t_3 < t_2 (p')$, $s_3 >0$
with $\vert s' \vert + s_3 < s_2 (p')$, $\theta_3 >0$ with $\vert
\theta '\vert + \theta_3 < \theta_2 (p')$, $\sigma_3 >0$ with
$\sigma_3 < \sigma_1$ and $r_3< 1$ with $r_2 (p') < r_3 <1$, the set:
\[
\aligned
\widetilde{W}_p^{CR,1+e}
&
:=
\big\{
A_{p',t,s,\sigma}(re^{i\,\theta}):\
\vert t-t'\vert<t_3,\
\vert s-s'\vert<s_3,\
\\
&
\ \ \ \ \ \ \ \ \ \ \ \ \ \ \ \ \ \
\ \ \ \ \ \ \ \ \ \ \ \ \ \ \ \ \ \
\vert\theta-\theta'\vert<\theta_3,\
r_3<r<1,\
0<\sigma<\sigma_3
\big\}
\endaligned
\]
will constitute a CR-wedge of dimension $2m + d + 1+ e$ at $p$. By a
technical adaptation of the approximation Theorem~5.2(III) ({\it cf.}
Lemma~2.19), $\mathcal{ C}_{ CR }^0 (M)$ extends to be CR on
$\widetilde{ W}_p^{ CR,1 +e}$.
\endproof

\subsection*{ 4.7.~Definition of the (non-integrable) subbundle 
$KM \subset TM$} Consequently, case {\bf (II)} cannot occur, because
of the definition of $e = e_{ \rm max}$. Thus, case {\bf (I)} holds.
In other words, as $p'$ runs in $M$, the $\mathcal{ C}^{ 1, \alpha -
0}$ distributions $p \mapsto KM_{ p'} (p)$ defined for $p$ in the edge
of $\mathcal{ W }_{ p'}^{ CR, e}$ (a neighborhood of $p'$ in $M$) glue
together in a well-defined $\mathcal{ C }^{ 1, \alpha - 0}$ vector
subbundle of $TM$. Observe that $T^cM$ is a subbundle of $KM$ of
codimension $e$. For every point $p\in M$, we have:
\[
T_p^cM
\subset
KM(p)
=
T_p^c\mathcal{M}_p^e
\cap
T_pM.
\]

As in \S4.1{\bf (j)}, since $M$ is globally minimal and since $KM$ is
of codimension $d-e \geqslant 1$ in $TM$, there must exist a point $p\in M$
such that $\left[ KM, KM \right] (p) \not\subset KM (p)$.

\def\thelemma{4.8}\begin{lemma}
At such a point $p$, the Levi form of $\mathcal{ M}_p^e$
does not vanish identically{\rm :}
\[
\left[
T^c\mathcal{M}_p^e, 
T^c\mathcal{M}_p^e
\right](p)
\not\subset
T_p^c\mathcal{M}_p^e.
\]
\end{lemma}

\proof
We reason by contradiction, assuming that $\left[ T^c\mathcal{ M}_p^e,
T^c\mathcal{ M}_p^e \right](p) \subset T_p^c\mathcal{ M}_p^e$.
Let $K^1$ and $K^2$ be two arbitrary $\mathcal{ C}^{ 1, \alpha - 0}$
sections of $KM$ defined in a small neighborhood $U_p$ of $p$ in $M$.
Since $KM\vert_{ U_p}$ is a subbundle of $TM\vert_{ U_p}$, we have
\[
\left[K^1, 
K^2\right](p)\in 
T_pM.
\]
We may extend $K^1$ and $K^2$ to a neighborhood $\mathcal{ U}_p$ of
$p$ in $\mathcal{ M}_p^e$ that contains $U_p$ as sections $\mathcal{
K}^1$ and $\mathcal{ K}^2$ of $T^c \mathcal{ M}_p^e \vert_{ \mathcal{
U}_p}$. Since $K^1$ and $K^2$ are tangent to $M\cap U_p$, one verifies
that, independently of the extension:
\[
\left[
K^1,K^2
\right](p)
=
\left[
\mathcal{K}^1,\mathcal{K}^2
\right](p)
\in
T_p^c\mathcal{M}_p^e,
\]
where the second Lie bracket belongs to $T_p^c\mathcal{ M}_p^e$,
because we assumed that the Levi form of $\mathcal{ M}_p^e$ vanishes
at $p$. We deduce
\[
\left[
K^1,K^2
\right](p)
\in
T_p^c\mathcal{M}_p^e\cap
T_pM
=
KM(p).
\]
This contradicts $\left[ KM, KM\right] (p ) \not\subset
KM(p)$.
\endproof

\subsection*{4.9.~Lewy extension on CR-wedges}
To contradict the maximality of $e= e_{\rm max}$ at a point $p$ at
which $\left[ KM, KM\right] (p) \not \subset KM (p)$, we formulate a
Lewy extension theorem on the conic manifold with edge $\mathcal{
W}_p^{ CR, e}$.

\def\theproposition{4.10}\begin{proposition}
Let $p\in M$ and assume that $\left[ T_p^c \mathcal{ M}_p^e, T_p^c
\mathcal{ M}_p^e \right] (p) \not \subset T_p^c \mathcal{ M}_p^e$.
Then there exists a $(2m + d + 1 + e)$-dimensional local CR-wedge
$\widetilde{ W}_p^{ CR, 1+ e}$ of edge $M$ at $p$ to which $\mathcal{
C}_{ CR}^0 (M \cup \mathcal{ W}_p^{ CR, e})$ extends to be CR.
\end{proposition}

Thus, this proposition concludes the proof of Theorem~3.8.

\proof
There exists a local section $L$ of $T^{ 1, 0} \mathcal{ M }_p^e$ with
$L (p ) \neq 0$ such that $\left[ L, \overline{ L} \right] (p) \not\in
T_p^{ 1, 0} \mathcal{ M}_p^e \oplus T_p^{ 0, 1} \mathcal{ M}_p^e$. It
is appropriate to distinguish two cases.

\smallskip

Firstly, assume that $L( p ) \in T^{ 1, 0} M$.
Then as in \S4.2, we
may construct a small analytic disc $A_\varepsilon$ attached to $M$ in
a neighborhood of $p$ having exit vector $- \frac{ \partial
A_\varepsilon}{ \partial r} (1)$ 
approximately directed by $\left[ L, \overline{ L} \right] (p)
\not \in \C \otimes T_p \mathcal{ M}_p^e$. So this disc has exit vector
nontangential to $\mathcal{ M}_p^e$ at $p$.
By
translating it along $M$ and along the $e$ supplementary directions
offered by $\mathcal{ W}_p^{CR, e}$, we deduce CR extension to a $(2m+
d + 1 + e)$-dimensional CR wedge $\widetilde{ W}_p^{ CR, 1+e}$.

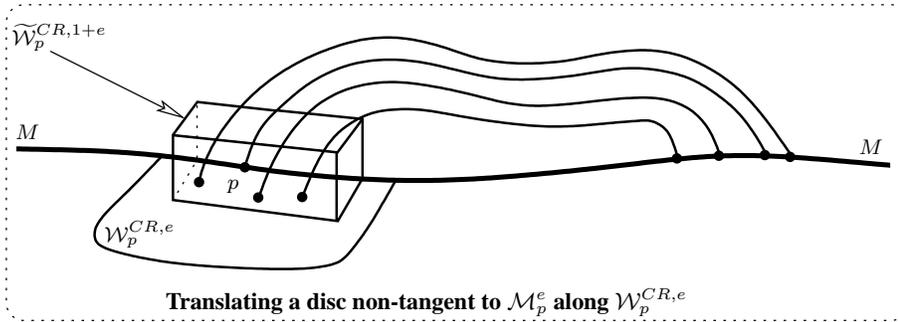
\begin{center}
\begin{picture}(0,0)%
\includegraphics{1-plus-e.pstex}%
\end{picture}%
\setlength{\unitlength}{4144sp}%
\begingroup\makeatletter\ifx\SetFigFont\undefined
\def\x#1#2#3#4#5#6#7\relax{\def\x{#1#2#3#4#5#6}}%
\expandafter\x\fmtname xxxxxx\relax \def\y{splain}%
\ifx\x\y   
\gdef\SetFigFont#1#2#3{%
  \ifnum #1<17\tiny\else \ifnum #1<20\small\else
  \ifnum #1<24\normalsize\else \ifnum #1<29\large\else
  \ifnum #1<34\Large\else \ifnum #1<41\LARGE\else
     \huge\fi\fi\fi\fi\fi\fi
  \csname #3\endcsname}%
\else
\gdef\SetFigFont#1#2#3{\begingroup
  \count@#1\relax \ifnum 25<\count@\count@25\fi
  \def\x{\endgroup\@setsize\SetFigFont{#2pt}}%
  \expandafter\x
    \csname \romannumeral\the\count@ pt\expandafter\endcsname
    \csname @\romannumeral\the\count@ pt\endcsname
  \csname #3\endcsname}%
\fi
\fi\endgroup
\begin{picture}(5424,1914)(446,-1910)
\put(514,-812){\makebox(0,0)[lb]{\smash{\SetFigFont{8}{9.6}{rm}{\color[rgb]{0,0,0}$M$}%
}}}
\put(1781,-1112){\makebox(0,0)[lb]{\smash{\SetFigFont{8}{9.6}{rm}{\color[rgb]{0,0,0}$p$}%
}}}
\put(1045,-1424){\makebox(0,0)[lb]{\smash{\SetFigFont{8}{9.6}{rm}{\color[rgb]{0,0,0}$\mathcal{W}_p^{CR,e}$}%
}}}
\put(5558,-898){\makebox(0,0)[lb]{\smash{\SetFigFont{8}{9.6}{rm}{\color[rgb]{0,0,0}$M$}%
}}}
\put(491,-246){\makebox(0,0)[lb]{\smash{\SetFigFont{8}{9.6}{rm}{\color[rgb]{0,0,0}$\widetilde{\mathcal W}_p^{CR,1+e}$}%
}}}
\put(1411,-1831){\makebox(0,0)[lb]{\smash{\SetFigFont{9}{10.8}{rm}{\color[rgb]{0,0,0}{\bf Translating a disc non-tangent to $\mathcal{M}_p^e$ along $\mathcal{ W}_p^{CR,e}$}}%
}}}
\end{picture}

\end{center}

Secondly, assume that $L (p) \not \in T_p^{ 1, 0} M$ for every local
section $L$ of $T^{ 1, 0} M$ such that $\left[ L, \overline{ L}
\right] (p) \not\in T_p^{ 1, 0} \mathcal{ M}_p^e \oplus T_p^{ 0, 1}
\mathcal{ M}_p^e$.

We explain the case $d = 2$, $e=1$ first, since this case is easier to
understand. Under this assumption, $\mathcal{ M}_p^1$ is a
hypersurface of $\C^n$ divided in two parts by $M$, one part being
$\mathcal{ W}_p^{ CR, 1}$. We draw a diagram.

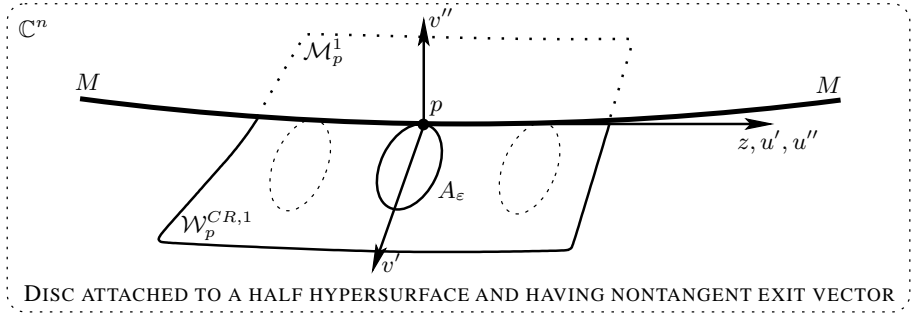
\begin{center}
\input attached-to-half.pstex_t
\end{center}

There exist coordinates $(z, w', w'') \in \C^{ n-2} \times \C \times
\C$ centered at $p$ in which $\mathcal{ M}_p^1$ is given by $v'' =
\psi (z, w', u'')$, with $\psi (0) = 0$ and $d\psi (0) = 0$ and in
which $M$ is given by a supplementary equation $v ' = \varphi' (z, u',
u'')$ with $\varphi' (0) = 0$ and $d\varphi ' ( 0) = 0$. Changing the
orientation of the $v'$-axis if necessary, it follows $\mathcal{
W}_p^{ CR, 1}$ is given by the equation $v'' = \psi (z, w', u'')$ and
the inequation $v' > \varphi '( z, u', u'')$, with $\varphi' ( 0) = 0$
and $d\varphi ' (0) = 0$. In the diagram, $T_p^c \mathcal{ M}_p^1$ is
the direct sum of the $z$-coordinate space with the $u'+ i\,
v'$-coordinate axis.

The Levi form of $\mathcal{ M}_p^1$ is represented by a scalar
Hermitian form $H (z, w', \bar z, \bar w')$. By assumption, its
restriction to $T_p^c M$ vanishes (otherwise, the first case holds),
so $H (z, 0, \bar z, 0 ) \equiv 0$. The assumption that the Levi form
of $\mathcal{ M}_p^1$ does not vanish identically insures that $H$ is
nonzero. To proceed further, we need $H (0, w', 0, \bar w') \not
\equiv 0$. If $H (0, w', 0, \bar w') \equiv 0$, since $H$ is nonzero,
by a linear coordinate change of the form $\widetilde{ w}' = w'$,
$\widetilde{ z}_k = z_k + a_k\, w'$, $k=1, \dots, n-2$, $\widetilde{
w}'' = w''$, we may insure that $H (0, w', 0, \overline{w}') \not
\equiv 0$. Observe that such a change of coordinates stabilizes both
$T_pM$ and $T_p \mathcal{ M}_p^1$. After a real dilation, we can
assume that the equation of $\mathcal{ M}_p^1$ is of the form:
\[
v''
=
w'\bar w'
+
{\rm O}(\vert w'\vert^{2+\alpha-0})
+
{\rm O}(\vert z\vert \vert (z,w')\vert)
+
{\rm O}(\vert u''\vert \vert (z,w')\vert)
+
{\rm O}(\vert u''\vert^2).
\]
To the hypersurface $\mathcal{ M}_p^1$, we attach a disc
$A_\varepsilon (\zeta)$ with zero $z$-component, with $w'$-component
equal to $i\, \varepsilon\, (1 - \zeta)$ and with $w''$-component
$(U_\varepsilon'' (\zeta) + i\, V_\varepsilon'' (\zeta))$ of class
$\mathcal{ C}^{ 2, \alpha - 0}$ satisfying the corresponding
Bishop-type equation. Exactly as in the Lewy extension theorem
(\S2.10), for $\varepsilon >0$ small enough and fixed, the exit vector
of $A_\varepsilon$ at $p$ is nontangent to $\mathcal{ M}_p^1$ (this is
uneasy to draw in the diagram above, but imagine that the disc drawn
in \S2.10 is attached to a half-paraboloid). Furthermore, using the
inequality $v' (e^{i\, \theta}) = \varepsilon (1 - \cos \theta)
\geqslant \varepsilon \, \frac{ \theta^2}{ \pi}$ for $\vert \theta
\vert \leqslant \pi$ together with the property $d\varphi' ( 0) = 0$,
it is elementary to verify that $A_\varepsilon (\partial \Delta
\backslash \{ 1\})$ is contained in the open half-hypersurface $\{ v'
> \varphi'\}$, as shown in the diagram.

Since the exit vector of $A_\varepsilon$ is nontangent to $\mathcal{
M}_p^1$, in order to get holomorphic extension to a wedge at $p$, it
suffices to translate the disc $A_\varepsilon$ in the
half-hypersurface $\mathcal{ W}_p^{ CR, 1}$.

However, if we translate $A_\varepsilon$ as usual by requiring that
the base point $A_{ \varepsilon, s} (1) = p_s$, with $s\in \R^{ 2n -
2}$ small, covers a neighborhood of $p$ in $M$, it may well happen
that, due to the curvature of $M$ in a neighborhood of $p$, the
boundary of the translated disc enters slightly the other side of
$\mathcal{ M}_p^1$, which is forbidden.

To remedy this imperfection, two equally good options present
themselves. The first option would be to rotate slightly the
translated disc $A_{ \varepsilon, s}$ in order that it becomes tangent
to $M$ at the point $p_s = A_{ \varepsilon, s} (1)$. Then adding a
small parameter $\sigma >0$, we would
translate it slightly in the positive
direction of $\mathcal{ W}_p^{ CR, 1}$, essentially along the positive
$v'$-direction.

The second option is to introduce a family of complex affine
biholomorphisms $\Psi_s$ that transfer $p_s \in M$ to the origin and
transfer the tangent spaces at $p_s$ 
of $\mathcal{ M}_p^1$ and of $M$ to $\{
v'' = 0\}$ and to $\{ v'' = v' = 0\}$. So $\Psi_s ( \mathcal{
M}_p^1)$ is given by $v'' = \psi'' (z, w', u'': s)$ with $\psi$ of
class $\mathcal{ C}^{ 2, \alpha - 0}$ with respect to all variables
and with the map $(z,w',u'') \mapsto \psi'' (z, w', u'' : s)$
vanishing to second order at the origin for every $s \in \R^{ 2n - 2}$
small. Also, $\Psi_s \big(
\mathcal{ W}_p^{ CR, 1} \big)$ is given by a
supplementary inequation $v' > \varphi' ( z, u', u'' : s)$, with
$\varphi'$ of class $\mathcal{ C}^{ 2, \alpha}$ (the smoothness of
$M$) with respect to all variables and with $(z, u', u'') \mapsto
\varphi (z, u', u'' : s)$ vanishing to second order at the origin.

To the hypersurface $\Psi_s (\mathcal{ M}_p^1)$, we attach the family
of discs
\[
\widetilde{A}_{ \varepsilon, s, \sigma} (\zeta)
=
\left(
0,\,
i\,\sigma
+
i\,\varepsilon\,(1-\zeta),\,
\widetilde{U}_{\varepsilon,s,\sigma}''(\zeta)
+
i\,\widetilde{V}_{\varepsilon,s,\sigma}''(\zeta)
\right)
\]
having zero $z$-component and $w'$-component equal to
$i\,\sigma+i\,\varepsilon\,(1 -\zeta)$, where $\sigma \in \R$ with
$\vert \sigma \vert < \sigma_1$, $\sigma_1 >0$, 
is a small parameter of translation
along the $v'$-axis. Of course:
\[
\left\{
\aligned
\widetilde{U}_{\varepsilon,s,\sigma}''(e^{i\,\theta})
&
=
-
{\sf T}_1
\big[
\psi
\big(
0,i\,\sigma
+
i\,\varepsilon(1-\cdot),
\widetilde{U}_{\varepsilon,s,\sigma}''(\cdot):s
\big)
\big](e^{i\,\theta}),
\\
\widetilde{V}_{\varepsilon,s,\sigma}''(e^{i\,\theta})
&
=
{\sf T}_1
\big[
\widetilde{U}_{\varepsilon,s,\sigma}''
\big](e^{i\,\theta}).
\endaligned\right.
\]
By means of elementary computations involving Taylor's formula, we
verify two facts.

\begin{itemize}

\smallskip\item[$\bullet$]
If $\varepsilon >0$ is sufficiently small and fixed, $\widetilde{
A}_{\varepsilon, s, 0} (\partial \Delta \backslash \{ 1\})$ is
contained in the open half-hypersurface $\{ v' > \varphi' (z, u', u''
: s)\}$, for all $s\in \R^{ 2n - 2}$ with $\vert s \vert < s_1$, $s_1
>0$ small.

\smallskip\item[$\bullet$]
Furthermore, for all $\sigma$ with $0 < \sigma \leqslant \sigma_1$, and 
all
$s$ with
$\vert s\vert < s_1$, the disc boundary $\widetilde{ A}_{ \varepsilon,
s, \sigma} (\partial \Delta)$ is contained in the open
half-hypersurface $\{ v' > \varphi' (z, u', u'' : s)\}$.

\end{itemize}\smallskip

Coming back to the old system of coordinates, it follows that the
family of discs $A_{ \varepsilon, s , \sigma} := \Psi_s^{ -1} \circ
\widetilde{ A}_{ \varepsilon, s, \sigma}$ has base point $A_{
\varepsilon, s, \sigma} (1)$ covering a neighborhood of $p$ in the
half-hypersurface $\mathcal{ W}_p^{CR, 1}$, 
as $s$ and $\sigma$ vary. Since the exit vector of $A_\varepsilon$
is not tangent to $\mathcal{ M}_p^1$ at $p$, this
family of discs covers a $2n$-dimensional wedge $\widetilde{ W}_p^{
CR, 2n}$ of edge $M$ at $p$. This completes the proof of the second
case of the proposition when $e=1$ and $d=2$.

\smallskip

Based on these explanations, we may now summarize the general case.
There exist coordinates $(z, w', w'') \in \C^m \times \C^e \times \C^{
d-e}$ vanishing at $p$ in which the
$\mathcal{ C}^{ 2, \alpha - 0}$
generic submanifold $\mathcal{ M}_p^e$ is represented by
$v'' =\psi (z, w', u'')$, with $\psi (0) = 0$ and $d\psi ( 0) = 0$.
After killing the second order pluriharmonic quadratic
terms in every right hand side $\psi_{ j''} (z,w',0)$, $j'' = 1,
\dots, d-e$, we may assume that the quadratic terms are Hermitian
forms $H_{ j''} (z, w', \bar z, \bar w')$.

After a linear 
change of coordinates in the $w'$-space, $T_p M = \{ v' = v''
= 0\}$, the $\mathcal{ C}^{ 2, \alpha}$ generic edge $M$ is defined by
$v = \varphi (z, u)$ with $\varphi ( 0) = 0$, $d\varphi (0) = 0$ and
the conic open submanifold $\mathcal{ W}_p^{ CR, e}$ of $\mathcal{
M}_p^e$ is defined by $v'' = \psi(z,w',u'')$ together with the
inequations
\[
v_{j'}'
>
\varphi_{j'}'(z,u',u''),
\ \ \ \ \ \ \ \ \ \ \ 
j'=1,\dots,e,
\]
where $\varphi = (\varphi', \varphi'')$. In fact, we may assume that
the cone defining the CR-wedge on the tangent space is slightly larger
than the salient cone $v_{ j'} ' >0$, $j'=1, \dots, e$.

The nonvanishing of the Levi form of $\mathcal{ M}_p^e$ at $p$ entails
that at least one Hermitian form $H_{ j''} (z, w', \bar z, \bar w')$
is nonzero. After renumbering, $H_1$ is nonzero. Also, since $T_p^c
M$ is the $z$-coordinate space, we have $H_1 (z, 0, \bar z, 0) \equiv
0$ (otherwise, the first case holds). After a complex linear
coordinate change of the form $\widetilde{ w}' = w'$, $\widetilde{
z}_k = z_k + \sum_{ j'=1}^e\, a_k^{j'}\, w_{j'}'$, $\widetilde{ w}''=
w''$, we may insure that $H_1 (0, w', 0, \overline{ w}') \not \equiv
0$. Then the set of vectors $(0,w')$ on which $H_1$ vanishes is a
proper real quadratic cone of $\C^e$. Consequently, for almost every real
vector $(0, i\, v')$, the quadratic form $H_1$ is nonzero on the
complex line $\C (0, i\, v')$. Since the cone defining $\mathcal{
W}_p^{ CR, e}$ is open and may be slightly shrunk, we can assume that
$H_1$ does not vanish on $\C (0, i\, v_1')$, with $v_1 ' = (1, \dots,
1) \in \R^e$. It follows that the disc $A_{ \varepsilon}$ attached to
$\mathcal{ M}_p^e$ having zero $z$-component and $w'$-component equal
to $\left( i\,\varepsilon\, (1 - \zeta), \dots, i\,\varepsilon\, (1 -
\zeta)\right)$ is nontangent to $\mathcal{ M}_p^e$ 
at $p$.

Furthermore, letting a point $p_s\in M$ of coordinates $s := (z, u)$
vary in a small neighborhood of $p$ in $M$, we may construct a family
of biholomorphisms $\Psi_s$ sending
$p_s$ to the origin and normalizing the equations of $M$, of $\mathcal{
M}_p^e$ and of $\mathcal{ W}_p^{ CR, e}$ under the form $v = \varphi(
z, u: s)$, $v'' =\psi (z, w', u'' : s)$ and $v_{ j'} ' > \varphi_{ j'}
' (z, u', u'': s)$, with $\varphi$ being
$\mathcal{ C}^{ 2, \alpha}$ and with
$\psi$ being $\mathcal{ C}^{ 2, \alpha - 0}$ with respect to all
variables and both vanishing to second order at the origin. 

Let $\sigma \in \R^e$, $\vert \sigma \vert < \sigma_1$, be a small
parameter of translation along the $v'$-coordinate space.
To the generic submanifold $\Psi_s (\mathcal{ M}_p^e)$, we attach the
family of discs
\[
\widetilde{A}_{\varepsilon,s,\sigma}(\zeta)
=
\left(
0,
W_{\varepsilon,\sigma}'(\zeta),
U_{\varepsilon,s,\sigma}''(\zeta)
+i\,
V_{\varepsilon,s,\sigma}''(\zeta)
\right),
\]
where 
\[
W_{\varepsilon,\sigma}'(\zeta)
=
\big(
i\,\sigma_1
+
i\,\varepsilon\,(1-\zeta),
\dots,
i\,\sigma_e
+
i\,\varepsilon\,(1-\zeta)
\big),
\]
and where
\[
\left\{
\aligned
\widetilde{U}_{\varepsilon,s,\sigma}''(e^{i\,\theta})
&
=
-
{\sf T}_1
\big[
\psi
\big(
0,i\,\sigma
+
i\,\varepsilon(1-\cdot),
\widetilde{U}_{\varepsilon,s,\sigma}''(\cdot):s
\big)
\big](e^{i\,\theta}),
\\
\widetilde{V}_{\varepsilon,s,\sigma}''(e^{i\,\theta})
&
=
{\sf T}_1
\big[
\widetilde{U}_{\varepsilon,s,\sigma}''
\big](e^{i\,\theta}).
\endaligned\right.
\]
By means of elementary computations involving Taylor's formula, we may
verify that for all $\sigma\in \R^e$ with $0 < \sigma_{ j'} \leqslant
\sigma_1$, $j'= 1, \dots, e$, and all $s\in \R^{ 2m + d}$, $\vert
s\vert < s_1$, the disc boundary $\widetilde{ A}_{ \varepsilon, s,
\sigma} (\partial \Delta)$ is contained in $\{ v_{ j'}' > \varphi_{
j'}' (z, u', u'' : s), \, j' = 1, \dots, e\}$.

Coming back to the old system of coordinates, it follows that the
family of discs $A_{ \varepsilon, s , \sigma} := \Psi_s^{ -1} \circ
\widetilde{ A}_{ \varepsilon, s, \sigma}$ has base point $A_{
\varepsilon, s, \sigma} (1)$ covering a neighborhood of $p$ in the
CR-wedge $\mathcal{ W}_p^{ CR, e}$, as $s$ and $\sigma$ vary. Since
its exit vector is not tangent to $\mathcal{ M}_p^1$ at $p$, this
family of discs covers a $(2m+ d+ 1 + e)$-dimensional CR-wedge
$\widetilde{ W}_p^{ CR, 1+ e}$ of edge $M$ at $p$.

The proofs of the proposition and of Theorem~3.8 are complete.
\endproof

\subsection*{ 4.11.~Wedgelike domains}
On a globally minimal $M$, at every point $p\in M$, we have
constructed a local wedge $\mathcal{ W}_p$ by gluing deformations of
discs. It may well happen that at a point $p$ that belongs to the
edges of two different wedges $\mathcal{ W}_{ q'}$ and $\mathcal{ W}_{
q''}$, the wedges have empty intersection in $\C^n$ (imagine two thin
opposite 
cones having vertex at $0 \in \R^2$). Fortunately, by means of the
translation trick presented in \S4.5 ({\it cf.} the diagram), we can
fill in the space in between. Achieving this systematically, by a sort
of gluing-shrinking processus, we obtain some connected open set
$\mathcal{ W}$ attached to $M$ containing possibly smaller wedges
$\mathcal{ W}_p ' \subset \mathcal{ W}_p$ at every point.

To set-up a useful definition, by a {\sl wedgelike domain} $\mathcal{
W}$ attached to $M$ we mean a {\it connected}\, open set that contains
a local wedge of edge $M$ at every point. Geometrically speaking, the
requirement of connectedness prevents jumps of the directions of local
wedges in the normal bundle $T\C^n\vert_M / TM$.

We may finally conclude this section by the formulation of a statement
that is the very starting point of the study of removable
singularities for CR functions (\cite{ mp1998, mp1999, mp2000, 
mp2002, mp2006a}).

\def\thetheorem{4.12}\begin{theorem}
{\rm (\cite{ me1997, mp1999})}
If $M$ is a globally minimal $\mathcal{ C}^{ 2, \alpha}$ generic
submanifold of $\C^n$, there exists a wedgelike domain $\mathcal{ W}$
attached to $M$ such that every continuous CR function $f\in \mathcal{
C}_{ CR}^0 (M)$ possesses a holomorphic extension $F \in \mathcal{ O}
(\mathcal{ W}) \cap \mathcal{ C}^0 (M\cup \mathcal{ W}) $ with
$F\vert_M = f$.
\end{theorem}

Its $L^{ \sf p}$ version deserves special attention. Let $\mathcal{
W}$ be a wedgelike domain attached to $M$. A holomorphic function $F
\in \mathcal{ O} (\mathcal{ W})$ is said to {\sl belong to the Hardy
space} $H_{ loc}^{ \sf p} (\mathcal{ W})$ if, for every $p\in M$, for
every local coordinate system centered at $p$ in which $M$ is given by
$v = \varphi (x, y, u)$, for every local wedge of edge $M$ at $p$
contained in $\mathcal{ W}$ of the form
\[
\aligned
\mathcal{W}
=
\mathcal{W}(\rho,\sigma,C)
:=
\big\{
&
(x+iy,u+iv)
\in\Delta_\rho^m
\times
\square_\rho^d
\times
i\square_\sigma^d:
\\
& \
v-\varphi (x,y,u)\in C\big\},
\endaligned
\]
as defined in \S4.29(III), for every cone $C' \subset \R^d$ with $C '
\cap S^{ d-1} \subset \subset C \cap S^{ d-1}$ and for every $\rho ' <
\rho$, the supremum:
\[
\sup_{\theta'\in C'}\,
\int_{\Delta_{\rho'}^m
\times\square_{\rho'}^d}\,
\left\vert
F\left(x+iy,u+i\varphi(x,y,u)+i\theta'\right)
\right\vert^{\sf p}\,
dx\wedge dy \wedge du \ 
< 
\ \infty
\]
is finite. An adaptation of the proof of the preceding theorem yields
its $L^{ \sf p}$ version.

\def\thetheorem{4.13}\begin{theorem}
{\rm (\cite{ po1997, po2000})} If $M$ is a globally minimal $\mathcal{
C}^{ 2, \alpha}$ generic submanifold of $\C^n$, there exists a
wedgelike domain $\mathcal{ W}$ attached to $M$ such that every
function $f\in L_{ loc, CR}^{ \sf p}$, $1\leqslant {\sf p} \leqslant
\infty$, possesses a Hardy space holomorphic extension $F \in {\sf
H}_{ loc}^{ \sf p} (\mathcal{ W})$.
\end{theorem}

To conclude, we would like to mention that arguments similar to those
of Theorem~4.12 yield a mild generalization, worth to be mentioned:
the CR extension theory is valid for $\mathcal{ C}^{2, \alpha}$ CR
manifolds that are only locally embeddable.

However, for concreteness reasons, we preferred to set up the theory in
a globally embedded context. In the remainder of the memoir, not to
enter superficial corollaries, we will formulate all our results under
the paradigmatic assumption of global minimality. Thus, {\em
Theorems~4.12 and~4.13 will be our basic main starting point}.

The two monographs~\cite{ trv1992, bch2005} deal not only with
embedded structures but also with locally integrable
structures. Nevertheless, most topics exposed here are not yet
embraced in a comprehensive theory ({\it cf.} \S3.29(III)). So it is
an open direction of research to transfer the theory of holomorphic
extension of CR functions (including removable singularities) to
locally integrable structures.

\newpage

\begin{center}
{\Large\bf VI:~Removable singularities}
\end{center}

\bigskip\bigskip\bigskip

\begin{center}
\begin{minipage}[t]{11cm}
\baselineskip =0.35cm
{\scriptsize

\centerline{\bf Table of contents}

\smallskip

{\bf 1.~Removable singularities for linear partial differential 
operators \dotfill 206.}

{\bf 2.~Removable singularities for holomorphic functions 
of one or several complex variables \dotfill 210.}

{\bf 3.~Hulls and removable singularities at the boundary \dotfill 227.}

{\bf 4.~Smooth and metrically thin removable singularities for
CR functions \dotfill 241.}

{\bf 5.~Removable singularities in CR dimension $1$ \dotfill 252.}

\smallskip

{\footnotesize\tt \hfill [7 diagrams]}

}\end{minipage}
\end{center}

\bigskip


{\small

Removable singularities for general linear partial differential
operators $P = \sum_{ \beta \in \N^m}\, a_\beta (x) \,
\partial_x^\beta$ on domains $\Omega \subset \R^n$ having order
$m\geqslant 1$ and $\mathcal{ C}^\infty$ coefficients have been
studied by Harvey and Polking (1970) in a general setting.
Assumptions of metrical thinness of singularities, in the sense of
Minkowski content or of Hausdorff measure, insure automatic
removability. For instance, relatively closed sets $C \subset \Omega$
whose $(n - m)$-dimensional Hausdorff measure is null are $(P, L_{ loc
}^\infty)$-removable. For structural reasons, these general results
(valid whatever the structure of the operator) necessitate a control
of growth when dealing with $L_{ loc }^1$-removability. In addition,
when $P$ is an embedded complex-tangential operator, this approach
does not convey to the adequate results, because removable
singularities for holomorphic or CR functions must take advantage of
automatic extension to larger sets.

Since almost two decades, thanks to the impulse of Stout, removable
singularities have attracted much attention in several complex
variables. A natural question is whether the Hartogs-Bochner extension
Theorem~1.9(V) holds when considering CR functions that are defined
only in the complement $\partial \Omega \backslash K$ of some compact
set $K \subset \partial \Omega$ of a connected smooth boundary
$\partial \Omega \subset \C^n$. In complex dimension $n=2$, Stout
showed that the answer is positive {\it if and only if}\, $K$ is
convex with respect to the space of functions that are holomorphic in
a neighborhood of $\overline{ \Omega}$. In complex dimension
$n\geqslant 3$, a complete cohomological characterization of different
nature was obtained by Lupacciolu (1994).

In another direction, by means of the above-cited global continuity
principle, J\"oricke (1995) generalized Stout's theorem to weakly
pseudoconvex domains. Recently, J\"oricke and the second author were
able to remove the pseudoconvexity assumption by applying purely
geometrical constructions without integral formulas, controlling
uniqueness of the extension (monodromy) by fine arguments.

Within the general framework of CR extension theory (exposed in
Part~V), the study of removable singularities has been endeavoured by
J\"oricke in the hypersurface case since 1988, and after by the two
authors in arbitrary codimension since 1995. The notions of CR-,
$\mathcal{ W }$- and $L^{ \sf p}$-removability, although different,
may be shown to be essentially equivalent, thanks to technical
deformation arguments. All the surveyed results hold in $L_{ loc}^{
\sf p}$ with $1 \leqslant {\sf p} \leqslant \infty$, including ${\sf
p} = 1$ and without any growth assumption near the singularity. On a
generic globally minimal $\mathcal{ C }^{ 2, \alpha }$ generic
submanifold $M$ of $\C^n$, closed sets $C \subset M$ having vanishing
$(\dim M - 2)$-dimensional Hausdorff measure are CR-, $\mathcal{ W}$-
and $L^{ \sf p}$-removable. As an application, CR meromorphic
functions defined on an everywhere locally minimal $M$ do extend
meromorphically to a wedgelike domain attached to $M$.

In conjunction with the Harvey-Lawson complex Plateau theorem,
singularities $C$ that are {\it a priori}\, contained in a
$2$-codimensional $\mathcal{ C }^{ 2, \alpha}$ submanifold $N$ of a
strongly pseudoconvex $\mathcal{ C}^{ 2, \alpha }$ boundary $\partial
\Omega \subset \C^n$ ($n \geqslant 3$) are shown by J\"oricke to be
{\it not}\, removable {\it if and only if}\, $N$ is a maximally
complex cycle. The condition that $N$ be somewhere generic was shown
by the two authors to be sufficient for its removability in arbitrary
codimension.

Concerning more massive singularities, a compact subset $K$ of a
one-codimensional submanifold $M^1 \subset\partial \Omega \subset
\C^n$ is CR-, $\mathcal{ W }$- and $L^{ \sf p}$-removable provided the
CR dimension of $\partial \Omega$ is $\geqslant 2$ (viz. $n\geqslant
3$) and provided $K$ does not contain any CR orbit of $M^1$
(J\"oricke, 1999). The second author generalized this theorem to
higher codimension, assuming that $M$ is globally minimal of CR
dimension $m \geqslant 2$. The main geometric argument (called {\sl
sweeping out by wedges}) being available only in CR dimension $m
\geqslant 2$, the more delicate case of CR dimension $m = 1$ is
studied extensively in the research article~\cite{ mp2006a}, 
placed in direct continuation to this survey.

}

\section*{ \S1.~Removable singularities for 
\\
linear partial differential operators}

\subsection*{ 1.1.~Hausdorff measure}
Let $M$ be a $\mathcal{ C}^1$ abstract manifold of dimension
$n\geqslant 1$ equipped with some Riemannian metric. For $\ell \in \R$
with $0 \leqslant \ell \leqslant n$, we remind (\cite{ ch1989}) the
definition of the notion of $\ell$-dimensional Hausdorff measure ${\sf
H}^\ell$ on $M$, that generalizes the notion of integer dimension of
submanifolds.

If $C \subset M$ is an arbitrary subset and if $\delta >0$ is small,
we define
\[
{\sf H}_{\delta}^\ell(C)
=
\inf
\left\{
\sum_{j=1}^{\infty} 
r_j^\ell:
C \
\text{\rm is covered by geodesic balls} \ B_j \
\text{\rm of radius} \ r_j \leqslant \delta
\right\}.
\]
Clearly, ${\sf H }_{ \delta }^\ell (C) \leqslant {\sf H}_{ \delta'
}^{\ell } (C)$, for $\delta ' \leqslant \delta$, so the limit ${\sf H
}^\ell (C) = \underset{ \delta \to 0^+}{ \lim} \, {\sf H}_{
\delta}^\ell (C)$ exists in $[0, \infty]$. This limit is called the
{\sl $\ell$-dimensional Hausdorff measure} of $C$. The value of ${\sf
H}^\ell (C)$ depends on the choice of a metric, but the two properties
${\sf H}^\ell (C) = 0$ and ${\sf H}^\ell (C) = \infty$ are
independent. The most significant property is that there exists a
critical exponent $\ell_C \geqslant 0$, called the {\sl Hausdorff
dimension} of $C$, such that ${\sf H}^\ell (C) = \infty$ for all $\ell
< \ell_C$ and such that ${\sf H}^{ \ell} (C) = 0$ for all $\ell>
\ell_C$. Then the value ${\sf H}^{ \ell_C } (C)$ may be arbitrary in
$[ 0, \infty]$.

\def\theproposition{1.2}\begin{proposition} 
{\rm (\cite{ fe1969, ch1989})} The following properties hold
true{\rm :}

\begin{itemize}

\smallskip\item[{\bf (1)}] 
${\sf H}^0(C) = {\rm Card} (C)${\rm ;}

\smallskip\item[{\bf (2)}] 
${\sf H}^n (C)$ coincides with the outer Lebesgue measure of $C
\subset M${\rm ;}

\smallskip\item[{\bf (3)}]
a $\mathcal{ C}^1$ submanifold $N \subset M$ has Hausdorff dimension
$\ell_N = \dim N${\rm ;}

\smallskip\item[{\bf (4)}] 
if ${\sf H}^{n -1} (C) = 0$, then $M\backslash C$ is locally
connected{\rm ;}

\smallskip\item[{\bf (5)}] 
if $f: M \to N$ is a $\mathcal{ C }^1$ map and if $C \subset M$
satisfies ${\sf H}^\ell (C) = 0$ for some $\ell \geqslant \dim N$,
then for almost every $q\in N$, it holds that ${\sf H}^{ \ell - \dim
N} (C \cap f^{ -1} (q)) = 0$.

\smallskip\item[{\bf (6)}] 
${\sf H}^\ell (C) = 0$ if and only if ${\sf H}^\ell (K) = 0$ for each
compact set $K \subset C$.

\end{itemize}\smallskip

\end{proposition}

\subsection*{1.3.~Metrically thin
singularities of linear partial differential operators}
Let $\Omega$ be a domain in $\R^n$, where $n\geqslant 1$. We shall
denote the Lebesgue measure by ${\sf H}^n$. Consider a class of
$\mathcal{ F} ( \Omega)$ of distributions defined on $\Omega$, for
instance $L_{ \rm loc}^{\rm p} ( \Omega)$, $\mathcal{ C}^{ \kappa,
\alpha} (\Omega)$ ($\kappa \in \N$, $0 \leqslant \alpha \leqslant 1$)
or $\mathcal{ C }^\infty ( \Omega)$. Consider a linear partial
differential operator
\[
P
= 
P(x,\partial_x) 
= 
\sum_{\beta \in \N^n,\ 
\vert\beta\vert\leqslant m}\,a_\beta(x)\,\partial_x^\beta
\]
of order $m\geqslant 1$, defined in $\Omega$ and having $\mathcal{
C}^\infty$ coefficients $a_\beta (x )$.

\def\thedefinition{1.4}\begin{definition}{\rm
A relatively closed subset $C$ of $\Omega$ is called $(P, \mathcal{ F}
)$-{\sl removable} if every $f\in \mathcal{ F} (\Omega)$ satisfying $P
f =0$ in $\Omega \backslash C$ does satisfy $P f= 0$ in all of
$\Omega$, in the sense of distributions.
}\end{definition}

For instance, according to the classical {\sl Riemann removability
theorem}, discrete subsets $\{ p_k \}_{ k\in \N}$ of a domain $\Omega$
in $\C$ are $(\overline{ \partial}, L^\infty )$-removable. In fact,
since every distribution solution of $\overline{ \partial}$ is
holomorphic (hypoellipticity), functions extend to be true holomorphic
functions in a neighborhood of each $p_k$. The Riemann removability
theorem also holds under the weaker assumption that $f\in \mathcal{ O}
(\Omega \backslash \{ p_k\}_{ k\in \N})$ satisfies $f (z - p_k) = {\rm
o} (\vert z - p_k \vert^{ -1})$ as $z$ approaches $p_k$.

In several complex variables, the classical {\sl Riemann removability
theorem} may be stated as follows.

\def\thetheorem{1.5}\begin{theorem}
{\rm (\cite{ ch1989})} Let $\Sigma$ be a complex analytic subset of
$\Omega$. Holomorphic functions in $\Omega \backslash \Sigma$ extend
uniquely through $\Sigma$ either if $\dim_\C \Sigma \leqslant n-2$ or
if $\dim_\C \Sigma = n-1$ and they belong to $L_{ loc}^\infty
(\Omega)$.
\end{theorem}

The second case also holds true for functions that belong to $L_{
loc}^2 (\Omega)$. The proofs are elementary and short: in one or
several complex variables, everything comes down to observing that
$\frac{ 1}{ z}$ is a true ${\rm O} (\frac{ 1}{ \vert z \vert})$ near
$z = 0$ and does not belong to $L_{ loc}^2$.

These preliminary statements are superseded by more general
removability theorems, exposed in~\cite{ hp1970}, that we shall now
restitute. Some of the (elementary) proofs will be surveyed to give
the flavour of the arguments.

\smallskip

In 1956, S.~Bochner (\cite{ bo1956, hp1970}) established
remarkable removability theorems, valid for general linear differential
operators $P$, in which the metrical conditions on the size of the
singularity $C$ depend only on the order $m$ of $P$. Some preliminary
material is needed.

\def\thelemma{1.6}\begin{lemma}
Let $K \subset \R^n$ be a compact set. For every $\varepsilon >0$,
there exists a function $\varphi_\varepsilon \in \mathcal{ C}_c^\infty
(\R^n)$ with $\varphi_\varepsilon \equiv 1$ in a neighborhood of $K$
and with ${\rm supp}\, \varphi_\varepsilon \subset K_\varepsilon$ such
that $\vert \partial_x^\beta \varphi_\varepsilon (x) \vert \leqslant
C_\beta \, \varepsilon^{ - \vert \beta \vert}$ for all $x \in \R^n$
and all $\beta \in \N^n$.
\end{lemma}

\proof
Denote by ${\bf 1}_B (\cdot)$ the characteristic function of a set $B
\subset \R^n$. It suffices to define the (rescaled) convolution
integral $\varphi_\varepsilon (x) := \varepsilon^{ -n}\, \int_{
\R^n}\, {\bf 1}_{ K_{ \varepsilon /2}} ( y) \, \psi( (x-y)/\varepsilon)\,
dy$, where $\psi \in \mathcal{ C}_c^\infty (\R^n)$ has support
contained in $\{ \vert x\vert \leqslant 1/3 \}$ and satisfies $\int \, \psi
(y) \, dy = 1$.
\endproof

It may happen that $C$ is not $( P, \mathcal{ F})$-removable, whereas
$C$ is removable for some individual function $f\in \mathcal{ F}
(\Omega)$ satisfying certain supplementary conditions. In this case,
we shall say that $C$ is an {\sl illusory singularity} of $f$.

\def\thetheorem{1.7}\begin{theorem}
{\rm (\cite{ bo1956, hp1970})}
Let $f\in L_{ loc}^1 (\Omega)$. If, for each compact set $K \subset
C$, we have 
\[
\underset{\varepsilon\to0^+}{\liminf}\
\left[
\varepsilon^{-m}\,
\left\vert\!\left\vert 
f\,{\bf 1}_{K_\varepsilon}
\right\vert\!\right\vert_{L^1}
\right]
=
0,
\]
then $C$ is an
illusory singularity of $f$.
\end{theorem}

Whenever the integral is meaningful, for instance if $f\in L_{ loc}^1
(\Omega)$ and $\varphi \in \mathcal{ C}_c^\infty (\Omega)$, we define
$(f, \varphi) := \int_\Omega \, f\, \varphi$, where the integral is
computed with respect to the Lebesgue measure. The {\sl formal
adjoint} of $P$, denoted by ${}^t\!P$, satisfies the relations 
$( P
\varphi, \psi ) = \big( \varphi, {}^t\!P \psi
\big)$ for all $\varphi, \, \psi
\in \mathcal{ C}_c^\infty (\Omega)$, and these relations define it
uniquely as
\[
{}^t\!P(\varphi) 
:= 
\sum_{\vert\beta\vert\leqslant m}\,(-1)^{\vert\beta\vert}\, 
\partial_x^\beta(a_\beta\,\varphi).
\]

\proof[Proof of Theorem~1.7.] Let $K := ({\rm supp}\, \varphi) \cap C$
and let $\varphi_\varepsilon$ be the family of functions constructed
in Lemma~1.6. Since ${\rm supp}\, P f \subset C$, we have $( P f,
\varphi) = ( P f, \varphi_\varepsilon\, \varphi) = (f, {}^t\! P(
\varphi_\varepsilon \, \varphi))$. Lemma~1.6 entails that $\left\vert
\! \left\vert {}^t \! P (\varphi_\varepsilon \, \varphi) \right\vert
\! \right\vert_{ L^\infty} \leqslant C \, \varepsilon^{ - m}$, for
some quantity $C >0$ that is independent of $\varepsilon$. We deduce
that $\left\vert (P f, \varphi ) \right\vert \leqslant C \,
\varepsilon^{ - m}\, \left\vert \! \left\vert f\, {\bf 1}_{
K_\varepsilon} \right\vert \! \right\vert_{ L^1}$ for all
$\varepsilon >0$. Thanks to the main assumption, this implies that $(
P f, \varphi) = 0$.
\endproof

If ${\sf p} \in \R$ with $1\leqslant {\sf p}\leqslant \infty$ is the
exponent of an $L^{ \sf p}$-space, we denote by ${\sf p}' := \frac{
{\sf p}}{{\sf p}- 1} \in [ 1, \infty]$ the {\sl conjugate exponent},
also defined by the relation $1 = \frac{ 1}{ {\sf p}} + \frac{ 1}{
{\sf p}'}$. By H\"older's inequality, we have $\left\vert \!
\left\vert f\, {\bf 1}_{ K_\varepsilon} \right\vert \! \right\vert_{
L^1} \leqslant \big( {\sf H}^n ( K_\varepsilon )\big)^{ 1/{\sf p}'}\,
\left\vert \! \left\vert f\, {\bf 1}_{ K_\varepsilon} \right\vert \!
\right\vert_{ L^{\sf p}}$.

\def\thecorollary{1.8}\begin{corollary}
Let $f\in L_{ loc}^{ \sf p} ( \Omega)$, where $1\leqslant {\sf p}
\leqslant \infty$. If, for each compact set $K \subset C$,
\[
\underset{\varepsilon\to 0^+}{\liminf}\,
\left[
\big(
\varepsilon^{-m{\sf p}'}
{\sf H}^n(K_\varepsilon)
\big)^{1/{\sf p}'}\,
\left\vert\!\left\vert
f\,{\bf 1}_{K_\varepsilon}
\right\vert\!\right\vert_{L^{\sf p}}
\right]
=
0,
\]
then $C$ is an illusory singularity of $f$.
\end{corollary}

The next theorem translates Corollary~1.9 in terms of Hausdorff
measures, a finer concept than the Minkowski content.

\def\thetheorem{1.9}\begin{theorem}
{\rm (\cite{ hp1970})}
{\bf (i)} Let $1 < {\sf p} < \infty$ and assume that
$n - m {\sf p}' \geqslant 0$. If ${\sf H}^{ n - m {\sf p}'}
(K) < \infty$ for every compact set $K \subset C$, then $C$ is $( P,
L_{ loc}^{\sf p})$-removable.

\smallskip
{\bf (ii)}
Let ${\sf p} = \infty$ and assume that
$n - m \geqslant 0$. If ${\sf H}^{ n- m} (C) = 0$, then $C$ is $(P,
L_{ loc}^\infty)$-removable.

\smallskip
{\bf (iii)}
Let ${\sf p} = \infty$ and assume that
$n - m \geqslant 0$. If, ${\sf H}^{ n- m} (K) < \infty$ for each
compact set $K \subset C$, then $P f$ is a measure supported on $C$,
for every $f \in L_{ loc}^\infty$ satisfying $P f= 0$ on $\Omega
\backslash C$.

\end{theorem}

An application of {\bf (ii)} to $P = \overline{ \partial}$ in one or
several complex variables yields the Riemann removability
Theorem~2.31 below.

\proof
We survey only the proof of {\bf (i)}. Let $\varphi \in \mathcal{
C}_c^\infty (\Omega)$ and set $K := C \cap {\rm supp}\, \varphi$.

\def\thelemma{1.10}\begin{lemma}
{\rm (\cite{ hp1970})}
Let $K \subset \R^n$ be a compact set. Let ${\sf p} '$ with $1\leqslant
{\sf p} ' < \infty$ and assume $n - m {\sf p}' \geqslant 0$. For every
$\varepsilon >0$, there exists $\varphi_\varepsilon \in \mathcal{
C}_c^\infty (\R^n)$ with $\varphi_\varepsilon \equiv 1$ in a
neighborhood of $K$ and with ${\rm supp} \, \varphi_\varepsilon
\subset K_\varepsilon$ such that for all $\beta \in \N^n$ with $\vert
\beta \vert \leqslant m$, we have
\[
\big\vert\!\big\vert
\partial_x^\beta\,\varphi_\varepsilon
\big\vert\!\big\vert_{L^{{\sf p}'}}
\leqslant
C\,\varepsilon^{m-\vert\beta\vert}
\big(
{\sf H}^{n-m{\sf p}'}(K)
+
\varepsilon
\big)^{1/{\sf p}'},
\]
where $C > 0$ is independent of $\varepsilon$.
\end{lemma}

With such cut-off functions $\varphi_\varepsilon$, since ${\rm supp}\,
Pf \subset C$, we have $(Pf, \varphi) = (Pf, \varphi_\varepsilon \,
\varphi) = (f, {}^t\! P (\varphi_\varepsilon \, \varphi))$. By
H\"older's inequality and the preceding lemma:
\[
\vert
(Pf,\varphi)
\vert
\leqslant
\left\vert\!\left\vert
f\,{\bf 1}_{K_\varepsilon}
\right\vert\!\right\vert_{L^{\sf p}}\,
\left\vert\!\left\vert
{}^t\!P(\varphi_\varepsilon\,\varphi)
\right\vert\!\right\vert_{L^{{\sf p}'}}
\leqslant
C\,
\left\vert\!\left\vert
f\,{\bf 1}_{K_\varepsilon}
\right\vert\!\right\vert_{L^{\sf p}}\,
\big(
{\sf H}^{n-m{\sf p}'}(K)
+
\varepsilon
\big)^{1/{\sf p}'}.
\]
The theorem follows from 
\[
\underset{\varepsilon\to 0^+}{\lim}
\left\vert\!\left\vert
f\,{\bf 1}_{K_\varepsilon}
\right\vert\!\right\vert_{L^{\sf p}}
=
0,
\]
since ${\sf H}^n ( K_\varepsilon) \to 0$ (remind ${\sf H}^{ n - m {\sf
p}'} (K ) < \infty$).
\endproof

It seems impossible to get $L^1$ removability without an assumption of
growth. At the opposite, in a CR context, the techniques introduced
in~\cite{ jo1999b, mp1999} that are developed in Section~5 and
in~\cite{ mp2006a} will exhibit $L^1$-removability of certain closed
subsets of generic submanifolds with only metrico-geometric
assumptions.






\section*{ \S2.~Removable singularities for holomorphic functions
of one or several complex variables}

\subsection*{ 2.1.~Painlev\'e problem, zero length and analytic 
capacity} The classical {\sl Painlev\'e problem} (\cite{ pa1888,
ah1947}) is to find metric or geometric characterizations of compact
sets $K \subset \C$ that are $(\overline{ \partial},
L^\infty)$-removable, {\it i.e.} such that every $f\in \mathcal{ O} (
\C \backslash K) \cap L^\infty ( \C \backslash K)$ extends
holomorphically through $K$.

Theorem~1.9{\bf (ii)} says that ${\sf H}^1 (K) =0$ suffices. It is
also known (\cite{ ma1984, pa2005}) that if ${\sf H}^{ 1 +
\varepsilon} (K) >0$ for some $\varepsilon >0$, then $K$
has positive analytic capacity (definition below) and is never
$(\overline{ \partial}, L^\infty)$-removable. Furthermore, Garnett
(\cite{ gar1970}) constructed a self-similar Cantor compact set $K
\subset \C$ with $0 < {\sf H}^1 (K ) < + \infty$ which is $(
\overline{ \partial}, L^\infty)$-removable. Consequently, Hausdorff
measure is not fine enough.

Under a geometric tameness assumption a converse to the sufficiency of
$H^1 (K) = 0$ holds and is usually called the solution to {\sl
Denjoy's conjecture}.

\def\thetheorem{2.2}\begin{theorem}
{\rm (\cite{ cal1977, cmm1982})} A compact set $K \subset \C$ that is
{\rm a priori} contained in a Lipschitz curve is $(\overline{
\partial}, L^\infty )$-removable {\rm if and only if} it has zero
$1$-dimensional Hausdorff measure.
\end{theorem} 

Classically, this statement is an application of the celebrated result
of Calder\'on, Coifman, McIntosh and Meyer about the $L^2$-boundedness
of the Cauchy integral on Lipschitz curves. Let us survey one of
the simplified proofs (\cite{ mv1995}) which involves {\sl Menger
curvature}, a concept useful in a recent answer to Painlev\'e's
problem obtained in~\cite{ to2003}.

\smallskip

Let $\Gamma := \big\{ (x, y) \in \R^2 : \, y = \varphi (x) \big\}$ be
a (global) Lipschitz graph; here $\varphi \in \mathcal{ C}^{ 0, 1}$ is
locally absolutely continuous and $\varphi'$ exists almost everywhere
(a.e.) with $\vert \! \vert \varphi ' \vert \! \vert_{ L^\infty} < +
\infty$.

\def\thetheorem{2.3}\begin{theorem}
{\rm (\cite{ cal1977, cmm1982, mv1995})}
If $f \in L^2 (\Gamma)$, the Cauchy 
principal value integral 
\[
{\sf C}^0f(z)
:=
\lim_{\varepsilon\to 0}\,
\frac{1}{2\pi i}\,
\int_{\vert\zeta-z\vert>\varepsilon}\,
\frac{f(\zeta)}{\zeta-z}\,d\zeta
\]
exists for almost every $z\in \Gamma$ and defines a
function ${\sf C}^0 f (z)$ on $\Gamma$, the
{\rm Cauchy transform} of $f$, which 
belongs to $L^2 (\Gamma)$ and satisfies
in addition
\[
\vert\!\vert
{\sf C}^0f
\vert\!\vert_{L^2(\Gamma)}
\leqslant
C_1\,
\vert\!\vert
f
\vert\!\vert_{L^2(\Gamma)},
\]
for some positive constant $C_1 = C_1 \big( 
\vert \! \vert \varphi ' \vert \! \vert_{ L^\infty} \big)$.
\end{theorem}

Parametrizing $\Gamma$ by $\zeta (t) = t + i \, \varphi (t)$, dropping
the innocuous factor $1 + i\, \varphi ' (t)$ and setting $z := x + i
\, \varphi (x)$, one has to estimate the $L^2$-norm of the truncated
integral
\[
{\sf C}_\varepsilon'(f)(x)
:=
\int_{\vert t-x\vert>\varepsilon}\,
\frac{f(t)}{\zeta(t)-\zeta(x)}\,dt,
\]
with a constant independent of $\varepsilon$. Even more,
interpolation arguments reduce the task to a
single estimate of the form
\[
\int_\R\,
\big\vert
{\sf C}_\varepsilon'(\chi_I)
\big\vert^2
\leqslant
C_1\,\vert I\vert,
\]
where $C_1 = C_1 \big( \vert \! \vert \varphi ' \vert \! \vert_{
L^\infty} \big)$ and where $\chi_I$ is the characteristic function of
an interval $I \subset \R$ of length $\vert I \vert$. Following
\cite{ mv1995}, a symmetrization 
of the (implicitely triple) integral 
$\int_I \, {\sf C}_\varepsilon ' (\chi_I) \, 
\overline{ {\sf C}_\varepsilon ' (\chi_I)}$
provides
\[
\aligned
6\,\int_I\,
\big\vert
{\sf C}_\varepsilon'(\chi_I)
\big\vert^2
=
&
\int\int\int_{S_\varepsilon}
\left(
\sum_{\sigma\in\mathfrak{S}_3}\,
\frac{1}{
\zeta(x_{\sigma(2)})
-
\zeta(x_{\sigma(1)})}\,
\frac{1}{
\overline{
\zeta(x_{\sigma(3)})
-
\zeta(x_{\sigma(1)})}}
\right)
\cdot
\\
&
\ \ \ \ \ \ \ \ \ \ \ \ \ \ \ \ \ \ \ \ \
\ \ \ \ \ \ \ \ \ \ \ \ \ \ \ \ \ \ \ \ \
\ \ \ \ \ \ \ \ \ \ \ \ \ \ \ \ \ \ \ \ \
\cdot
dx_1dx_2dx_3
+
{\rm O}(\vert I\vert),
\endaligned
\]

\noindent
where $S_\varepsilon := \big\{ (x,y,t) \in I^3 : \, \vert y - x \vert
> \varepsilon, \ \vert t - x \vert > \varepsilon, \ \vert t - y \vert
> \varepsilon \big\}$ and where $\mathfrak{ S}_3$ is the permutation
group of $\{ 1, 2, 3 \}$.

Then a ``magic'' (\cite{ po2005}) formula enters the scene:
\[
\left(
\frac{4\,S( z_1, z_2, z_3)}{
\vert z_1-z_2\vert\,
\vert z_1-z_3\vert\,
\vert z_2-z_3\vert
}
\right)^2
=
\sum_{\sigma\in\mathfrak{S}_3}\,
\frac{1}{
\zeta(x_{\sigma(2)})
-
\zeta(x_{\sigma(1)})}\,
\frac{1}{
\overline{
\zeta(x_{\sigma(3)})
-
\zeta(x_{\sigma(1)})}},
\]
where $S (z_1, z_2, z_3)$ denotes the enclosed area; the left hand
side measures the ``flatness'' of the triangle. This crucial formula
enables one to link rectifiability properties to the Cauchy kernel.

\def\thedefinition{2.4}\begin{definition}{\rm
The {\sl Menger curvature}\, of the triple $\{ z_1, z_2, z_3 \}$
is the square root of the above
\[
c(z_1,z_2,z_3)
:=
\frac{4\,S(z_1,z_2,z_3)}{
\vert z_1-z_2\vert\,
\vert z_1-z_3\vert\,
\vert z_2-z_3\vert
};
\]
one sets $c := 0$ if the points are aligned. One also verifies that
$c (z_1, z_2, z_3) = 1 / R (z_1, z_2, z_3)$, where $R$ is the radius
of the circumbscribed circle.
}\end{definition}

Thanks to the nice formula and to the basic inequality
\[
c\big(\zeta(x),\zeta(y),\zeta(t)\big)
\leqslant
2\,
\left\vert
\frac{
\frac{\varphi(y)-\varphi(x)}{y-x}
-
\frac{\varphi(t)-\varphi(x)}{t-x}
}{
\vert t-y\vert}
\right\vert
\]
the previous symmetric Cauchy triple integral is transformed to an
integral involving geometric Lipschitz
properties of $\Gamma$. After some
computations (\cite{ mv1995}), one gets $\int_I \, \big\vert {\sf
C}_\varepsilon ' (\chi_I) \big\vert^2 \leqslant C_1 \big( \vert \!
\vert \varphi '\vert \! \vert_{ L^\infty} \big) \cdot \vert I \vert$.

\smallskip

Menger curvature also appears in a recent result, considered to 
be an answer to Painlev\'e's problem.

\def\thetheorem{2.5}\begin{theorem}
{\rm (\cite{ to2003, pa2005})}
A compact set $K \subset \C$ is
{\rm not} removable for 
$\mathcal{ O} (\C \backslash K) \cap L^\infty
(\C \backslash K)$ {\rm if and only if}
there exists a {\rm nonzero} positive
Radon measure $\mu$ with ${\rm supp} \, \mu \subset K$
such that

\begin{itemize}

\smallskip\item[$\bullet$]
there exists $C_1 >0$ with $\mu \big( \Delta (z, \rho) \big) \leqslant
C_1 \, \rho$ for every $z\in \C$ and $\rho >0${\rm ;}

\smallskip\item[$\bullet$]
$\int \int \int \, \big[ c ( x, y, z) \big]^2 \, 
d\mu(x) d\mu(y) d\mu(z) < + \infty$.

\end{itemize}\smallskip
\end{theorem}

\bigskip

The first condition concerns the size of $K$; the second
one is of quantitative-geometric nature.

\smallskip

We conclude by mentioning a classical functional characterization due
to Ahlfors, usually considered to be only a reformulation of
Painlev\'e's problem, but which has already found generalizations in
locally integrable structures (\S2.16 below). The {\sl analytic
capacity} of a compact set $K \subset \C$ is\footnote{ If $\vert \!
\vert f \vert \! \vert_{ L^\infty} \leqslant 1$, setting $g(z) :=
\big[ f(\infty) - f(z) \big] \big/ \big[ 1 - \overline{ f(\infty)} \,
f(z) \big]$, we have $g (\infty) = 0$, $\vert \! \vert g \vert \!
\vert_{ L^\infty} \leqslant 1$ and $\vert g' (\infty) \vert = \vert
f'(\infty) \vert \big/ \big( 1 - \vert f (\infty) \vert^2 \big)
\geqslant \vert f'( \infty) \vert$, so that in the definition of
analytic capacity, we may restrict to take the supremum over functions
$g\in H^\infty (\C \backslash K)$ with $\vert \! \vert g \vert \!
\vert_{ L^\infty} \leqslant 1$ and $g (\infty) = 0$.
}
\[
\text{\sf an-cap}(K)
:=
\sup\big\{
\vert f'(\infty)\vert:\
f\in H^\infty(\C \backslash K),\,
\vert\!\vert f\vert\!\vert_{L^\infty}\leqslant1
\big\},
\]
where $H^\infty (\C \backslash K)$ denotes the space of bounded
holomorphic functions defined in $\C \backslash K$ (or defined in the
complement of $K$ in the Riemann sphere $\C \cup \{ \infty \}$,
because $\{ \infty \}$ is removable).

\def\thetheorem{2.6}\begin{theorem}
{\rm (\cite{ ah1947, ma1984, ht1997})} A compact set $K \subset \C$ is
removable for $\mathcal{ O} (\C \backslash K) \cap L^\infty (\C
\backslash K)$ if and only if $\text{\sf an-cap} (K) = 0$.
\end{theorem}

\subsection*{ 2.7.~Rad\'o-type theorems}
A classical theorem due to Rad\'o (\cite{ ra1924, stu1968, rs1989,
ch1994}) asserts that a continuous function $f$ defined in a domain
$\Omega \subset \C$ that is holomorphic outside its zero-set $f^{ -1}
(0)$ is in fact holomorphic everywhere. By a separate holomorphicity
argument, this statement extends directly to several complex
variables. In~\cite{ stu1993}, it is shown that $f^{ -1} (0)$ may be
replaced by $f^{ -1} (E)$, where $E \subset \C$ is compact and has
null analytic capacity. In~\cite{ rs1989}, it is shown that a
continuous function defined in a strongly pseudoconvex $\mathcal{
C}^2$ hypersurface $M \subset \C^n$ ($n\geqslant 2$) that is CR
outside its zero-set is CR everywhere; a thin subset of weakly
pseudoconvex points is allowed, but the case of general hypersurfaces
is not covered. In~\cite{ al1993}, it is shown that closed sets $f^{
-1} (E)$ are removable in the same situation, wehere $E \subset \C$ is
a closed polar set, viz. $E \subset \{ u = - \infty \}$ for some
subharmonic function $u \not \equiv - \infty$. Chirka strengthens
these results in the following theorem, where no assumption is made on
the geometry of the hypersurface.

Remind (\cite{ ch1989, de1997}) 
that $E \subset \C^m$ is called {\sl complete
pluripolar} if $E = \{ \varphi = -\infty \}$ for some
plurisubharmonic function $\varphi \not\equiv - \infty$ on $\C^m$.

\def\thetheorem{2.8}\begin{theorem}
{\rm (\cite{ ch1994})}
Let $M \subset \C^n$ {\rm (}$n\geqslant 2${\rm )} be hypersurface that
is a local Lipschitz {\rm (}$\mathcal{ C}^{ 0, 1}${\rm )} graph at
every point, let $C$ be a closed subset of $M$ and let $f : M
\backslash C \to \C^m \backslash E$ be a continuous mapping satisfying
$\vert \! \vert f \vert \! \vert_{ \mathcal{ C}^0 (M\backslash C)} <
\infty$ such that the set of limit values of $f$ from $M\backslash C$
up to $C$ is contained in a closed complete pluripolar set $E \subset
\C^m$ {\rm (}$m\geqslant 1${\rm )}. Then the trivial extension
$\widetilde{ f}$ of $f$ to $C$ defined by $\widetilde{ f} := 0$ on $C$
is a CR mapping of class $L^\infty$ on the whole of $M$.
\end{theorem}

In higher codimension, nothing is known.

\def\theopenquestion{2.9}\begin{openquestion}
Let $M \subset \C^n$ {\rm (}$n\geqslant 3${\rm )} be a generic
submanifold of codimension $d\geqslant 2$ and of CR dimension
$m\geqslant 1$ that is at least $\mathcal{ C}^1$. Let $f\in \mathcal{
C}^0 (M)$ that is CR outside its zero-set $f^{ -1} (0)$. Is $f$ CR
everywhere\,?
\end{openquestion}

Remind that condition {\bf (P)} (Definition~3.5(III)) for a linear
partial differential operator $P$ of principal type assures local
solvability of the equations $P f = g$. Remind also that nowhere
vanishing vector fields are of principal type.

\def\thetheorem{2.10}\begin{theorem}
{\rm (\cite{ ht1993})} Let $\Omega \subset \R^n$ {\rm (}$n\geqslant
2${\rm )} be a domain and let $L$ be a nowhere vanishing vector field
on $\Omega$ having $\mathcal{ C}^\infty$ complex-valued coefficients
and satisfying condition {\bf (P)} of Nirenberg-Treves. If $f \in
\mathcal{ C}^0 (\Omega)$ satisfies $L f = 0$ in $\Omega \backslash f^{
-1} (0)$ in the sense of distributions, then $f$ is a weak solution of
$L f= 0$ all over $\Omega$.
\end{theorem}

\subsection*{ 2.11.~Capacity and partial differential operators having 
constant coefficients} The preceding results admit partial
generalizations to vector field systems. Let $\Omega \subset \R^n$ be
an open set and let $P = P (\partial_x) = \sum_{ \beta\in \N^n}\,
a_\beta \, \partial_x^\beta$ be a linear partial differential operator
having {\it constant coefficients}\, $a_\beta \in \C$. By a theorem
due to Malgrange, Ehrenpreis and Palamodov (\cite{ ho1963}), such a
$P$ always admits a {\sl fundamental solution}, namely there exists a
distribution $E\in \mathcal{ D} ' (\R^n)$ such that $P (\partial_x) E
= \delta_0$ is the Dirac measure at the origin.

Let $\mathcal{ F} \subset \mathcal{ D}' (\Omega)$ be a Banach space,
{\it e.g.} $\mathcal{ F} = L^{\sf p} (\Omega)$ with $1 \leqslant {\sf
p} \leqslant \infty$, or $\mathcal{ F} = L^\infty (\Omega) \cap
\mathcal{ C}^0 (\Omega)$, or $\mathcal{ F} = \mathcal{ C}^{ 0, \alpha}
(\Omega)$ with $0 < \alpha \leqslant 1$.

\def\thedefinition{2.12}\begin{definition}{\rm
For each relatively closed set $C \subset \Omega$, the {\sl $\mathcal{
F}$-capacity of $C$ with respect to $P$} is
\[
\mathcal{ F}\text{\sf -cap}_P
(C,\partial \Omega)
:=
\sup
\Big\{
\big\vert(Pf,{\bf 1}_\Omega)\big\vert:\,
f\in\mathcal{F},\
\vert\!\vert f\vert\!\vert_{\mathcal{F}}\leqslant 1,\
{\rm supp}\, 
(Pf)\Subset C
\Big\}.
\]
}\end{definition}

If a closed set $C \subset \Omega$ is $(P, \mathcal{ F}$)-removable,
by definition $P f = 0$ everywhere, hence $\mathcal{ F}$-$\text{\sf
Cap}_P (C, \Omega) = 0$. The following theorem establishes the
converse for a wide class of differential operators having constant
coefficients. For $\beta\in \N^n$, denote by $Q^{(\beta)} (x) :=
\partial_x^\beta \, Q (x)$ the $\beta$-th partial derivative of a
polynomial $Q (x) \in \R [ x]$.

\def\thetheorem{2.13}\begin{theorem}
{\rm (\cite{ hp1972})} Assume that $P$ possesses a fundamental
solution $E\in \mathcal{ D}' (\R^n)$ such that $P^{ (\beta)}
(\partial_x) E$ is a regular Borel measure on $\R^n$ for every $\beta
\in \N^n$. Let $\Omega \Subset \R^n$ be a bounded domain and let $C
\subset \Omega$ be a relatively closed subset. Then

\begin{itemize}

\smallskip\item[$\bullet$]
$C$ is $(P, L_{ loc}^{\sf p})$-removable, $1 < {\sf p}
\leqslant \infty$, if and only if $L^{\sf p}$-$\text{\sf cap}_P (C,
\Omega) = 0${\rm ;}

\smallskip\item[$\bullet$]
$C$ is $(P, L^\infty \mathcal{ C}^0)$-removable if and only if
$L^\infty \mathcal{ C}^0$-$\text{\sf cap}_P (C, \Omega ) =0${\rm ;}

\smallskip\item[$\bullet$]
$C$ is $(P, \mathcal{ C}^{ 0, \alpha})$-removable, $0 < \alpha
\leqslant 1$, if and only if $\mathcal{ C}^{ 0, \alpha}$-$\text{\sf
cap}_P (C, \Omega) = 0$.

\end{itemize}\smallskip
\end{theorem}

This hypothesis about $P$ is satisfied by elliptic, semi-elliptic,
and parabolic operators and also by the wave operator in $\R^2$ (\cite{
hp1972}). The theorem (whose proof is rather short) also holds true
if $P = P (x, \partial_x)$ has real analytic coefficients and admits a
fundamental solution $E$ such that $P^{ (\beta)} E$ is a regular Borel
measure for every $\beta \in \N^n$. But it is
void in $L^1$.

\def\thetheorem{2.14}\begin{theorem}
{\rm (\cite{ hp1972})} There is a {\rm unique} function, called a {\rm
capacitary extremal}, $f^{\sf cap} \in L^{ \sf p} (\Omega)$ with
$\vert \! \vert f^{\sf cap} \vert \! \vert_{ L^{\sf p}} \leqslant 1$
and $P f^{ \sf cap} = 0$ in $\Omega \backslash K$ such that $\big( P
f^{\sf cap}, {\bf 1}_\Omega \big) = L^{\sf p}$-$\text{\sf cap}_P (K,
\Omega)$.
\end{theorem}

We observe that the definition of $L^{\sf p}$-$\text{\sf cap}_P (K,
\Omega)$ is inspired from Ahlfors' notion of analytic capacity and we
mention that the capacitary extremal $f^{\sf cap}$ is linked to the
Riemann uniformization theorem.

\def\theexample{2.15}\begin{example}{\rm
In fact, with $\Omega = \C$ and $P = \partial / \partial \bar z =:
\overline{ \partial}$, the $L^\infty$-capacity of a compact set $K
\subset \C$ with respect to $\overline{ \partial}$ may be shown to be
equal, up to the constant $\pi$, to the analytic capacity of $K$,
namely
\[
L^\infty\text{\sf -cap}_{\overline{\partial}}
(K,\C)
=
\pi\,
\text{\sf an-cap}(K).
\]
Indeed, letting $f\in L^\infty (\C)$, assuming that $\overline{
\partial }f$ is supported by $K$, choosing a big open disc $D \Subset \C$
containing $K$, integrating by parts (Riemann-Green) and performing
the change of variables $w := 1 / z$, we may compute
\[
\big(
\overline{\partial}f,{\bf 1}_\C
\big)
=
\big(
\overline{\partial}f,{\bf 1}_D
\big)
=
\frac{1}{2i}\,
\int\!\!\int_D\,
\frac{\partial f}{\partial\bar z}\,
d\bar z\wedge dz
=
\frac{1}{2i}\,
\int_{\partial D}\,
f(z)\,dz
=
\pi\,f'(\infty).
\]
Remind (footnote) that in the definition of $\text{\sf an-cap}
(K)$ given in \S2.1, one may assume that $f( \infty) = 0$. If in
addition the complement of $K$ in the Riemann sphere $\C \cup \{
\infty \}$ is simply connected, the unique solution $f^{\sf cap}$ of
$\big( \overline{ \partial} f^{\sf cap}, {\bf 1}_\C \big) =
L^\infty\text{\sf -cap}_{ \overline{ \partial }} (K, \C)$ asserted by
Theorem~2.14, viz. the unique solution of the extremal problem
\[
\sup
\Big\{
\vert f'(\infty)\vert:\,
f\in L^\infty(\C),\
\partial f/\partial\bar z=0\
\text{\rm in}\
\C\backslash K,\
f(\infty)=0\
\text{\rm and}\
\vert\!\vert f\vert\!\vert_{L^\infty}\leqslant 1
\Big\}
\]
is the (unique) Riemann uniformization map $f^{\sf cap} : \big( \C
\cup \{ \infty \} \big) \backslash K \to \Delta$ satisfying $f^{\sf
cap} (\infty) = 0$ and $\partial_z f^{\sf cap} (\infty) > 0$.
}\end{example}

\subsection*{ 2.16.~Removable singularities of locally solvable
vector fields} Let
\[
\mathcal{S}(\R^n)
:=
\Big\{
f\in\mathcal{C}^\infty(\R^n):\,
\lim_{\vert x\vert\to\infty}\,
\big\vert
x^\alpha\partial_x^\beta f(x)
\big\vert
=0,\
\forall\,\alpha,\,\beta\in\N^n
\Big\}
\]
be the space of $\mathcal{ C}^\infty$ functions defined in $\R^n$ and
having tempered growth. As is known (\cite{ ho1963}), the Fourier
transform
\[
{\sf F}f(\xi)
:=
\int_{\R^n}\,
e^{-2\pi\,i\left<x,\xi\right>}\,f(x)\,dx,
\ \ \ \ \ \ \ \ \ \ \ \ \ \ \ 
f\in\mathcal{S}(\R^n),
\]
$\left< x, \xi\right> := \sum_{ k=1}^n\, x_k \, \xi_k$, is an
automorphism of $\mathcal{ S} (\R^n)$ having as inverse
\[
{\sf F}^{-1}f(\xi)
:=
\int_{\R^n}\,
e^{2\pi\,i\left<x,\xi\right>}\,f(x)\,dx
=
{\sf F}f(-\xi).
\]
Equipping $\mathcal{ S} (\R^n)$ with the countable family of
semi-norms $p_{ \alpha, \beta} (f) := \sup_{ x\in \R^n}\, \big\vert
x^\alpha \partial_x^\beta f(x) \big\vert$, the space $\mathcal{ S}'
(\R^n)$ of {\sl tempered distributions} consists of linear functionals
${\sf T}$ on $\mathcal{ S} (\R^n)$ that are continuous, viz. there
exists $C >0$ and $\alpha, \, \beta \in \N^n$ such that $\big\vert \!
\left< {\sf T}, f \right> \! \big\vert \leqslant C \, p_{\alpha, \beta}
(f)$ for every $f\in \mathcal{ S} (\R^n)$.

For ${\sf p} \in \R$ with $1\leqslant {\sf p} \leqslant \infty$ and for
$\sigma \in \R$, we remind the definition of the {\sl Sobolev space}
\[
L_\sigma^{\sf p}(\R^n)
:=
\Big\{
{\sf T}\in\mathcal{S}'(\R^n):\
\vert\!\vert{\sf T}\vert\!\vert_{L_\sigma^{\sf p}}
:=
\big\vert\!\big\vert\Lambda^{-\sigma}\,{\sf T}
\big\vert\!\big\vert_{L^{\sf p}}
<
\infty
\Big\},
\]
where $\Lambda^{ -\sigma}\, {\sf T} (x) := {\sf F}^{ -1} \big[ (1 +
\vert \xi \vert^2)^{ -\sigma /2}\, {\sf F}\, {\sf T} (\xi) \big](x)$.
For $\sigma = \kappa \in \N$ and ${\sf p}$ in the range $1 < {\sf p} <
\infty$, the space $L_\kappa^{ \sf p} (\R^n)$ is exactly the subspace
of functions $u\in L^{ \sf p} (\R^n)$ whose partial derivatives of
order $\leqslant \kappa$ (in the distributional sense) belong to
$L^{\sf p} (\R^n)$. This space is equivalently normed by $\vert \!
\vert u \vert \! \vert_{ L_\kappa^{ \sf p }} := \sum_{ \vert \beta
\vert \leqslant \kappa}\, \big\vert \! \big\vert \partial_x^\beta u
\big\vert \! \big\vert_{ L^{\sf p}}$.

\smallskip

Let $\Omega \subset \R^n$ be a domain and let $P = P (x, \partial_x) =
\sum_{ \vert \beta \vert \leqslant m}\, a_\beta (x) \,
\partial_x^\beta$ be a linear partial differential operator of order
$m \geqslant 1$ defined in $\Omega$ and having $\mathcal{ C }^\infty$
coefficients $a_\beta (x)$.

\def\thedefinition{2.17}\begin{definition}{\rm
We say that $P$ is {\sl locally solvable in $L^{\sf p}$ with one loss
of derivative}\, if every point $p \in \Omega$ has an open
neighborhood $U_p \subset \Omega$ such that for every compactly
supported ${\sf T} \in L_\sigma^{\sf p} (U_p)$, the equation ${\sf
P}\, {\sf S} = {\sf T}$ has a solution ${\sf S} \in L_{\sigma + m -
1}^{\sf p} (U_p)$.
}\end{definition}

\def\thetheorem{2.18}\begin{theorem}
{\rm (\cite{ befe1973, hp1996})} Let $L$ be a nowhere vanishing vector
field having $\mathcal{ C}^\infty$ coefficients in a domain $\Omega
\subset \R^n$ {\rm (}$n\geqslant 2${\rm )} and assume that $L$ satisfies
condition {\bf (P)} of Nirenberg-Treves. Then for every ${\sf p} \in \R$
with $1 < {\sf p} < \infty$, the operator $L$ is locally solvable in
$L^{\sf p}$ with loss of one derivative.
\end{theorem}

\def\theexample{2.19}\begin{example}{\rm
As discovered in~\cite{ ht1996}, local solvability fails to hold in
$L^\infty$ for the (locally solvable) vector field $\frac{ \partial
}{\partial x} - \frac{
i}{ x^2} \, e^{ -1 / \vert x \vert} \, \frac{
\partial }{ \partial y}$ satisfying {\bf (P)} on $\R^2$.
}\end{example}

Removable singularities for vector fields in $L^{\sf p}$
have been studied in~\cite{ ht1996, ht1997}.
Because of the example, results in $L^{\sf p}$
with $1 < {\sf p} < \infty$ differ from
results in $L^\infty$.

\def\thedefinition{2.20}\begin{definition}{\rm
A relatively closed set $C \subset \Omega$ of an open set $\Omega
\subset \R^n$ is {\sl everywhere $(P, L^{\sf p})$-removable}\, if for
every open subset $U \subset \Omega$ and for every $f\in L^{ \sf p}
(U)$ satisfying $P f = 0$ in $U \backslash C$, then $f$ also satisfies
$P f = 0$ in all of $U$.
}\end{definition}

\def\thetheorem{2.21}\begin{theorem}
{\rm (\cite{ ht1996, ht1997})} Let $L$ be a nowhere vanishing vector
field having $\mathcal{ C}^\infty$ coefficients in an open subset
$\Omega \subset \R^n$ {\rm (}$n\geqslant 2${\rm )} and assume that $L$
satisfies condition {\bf (P)}. Let ${\sf p} \in \R$ with $1 < {\sf p}
< \infty$. Then a relatively closed set $C \subset \Omega$ is
everywhere $(L, L^{\sf p})$-removable {\rm if and only if} there is an
open covering of $C$ by open sets $\Omega_j$, $j\in J$, such that
\[
L^{\sf p}\text{\sf -cap}_L
\big(C\cap\Omega_j,\Omega_j\big)
=0,
\]
for every $j\in J$.
\end{theorem}

In $L^\infty$, when trying to perform the proof of this theorem,
local solvability of positive multiples of $L$ is technically
needed. Observing that $P f = 0$ is equivalent to $e^\psi P f = 0$,
the following notion appeared to be appropriate to deal with $(L,
L^\infty )$-removability.

\def\thedefinition{2.22}\begin{definition}{\rm
{\rm (\cite{ ht1996, ht1997})}
The {\sl full $L^\infty$-capacity}\,
of a relatively closed set $C \subset \Omega$ with respect
to $L$ is
\[
\text{\sf full-}L^\infty\text{\sf -cap}_L
(C,\Omega)
:=
\sup_{\widetilde{L}}\,
\Big\{
L^\infty\text{\sf -cap}_{\widetilde{L}}(C,\Omega)
\Big\},
\]
where the supremum is taken over all vector fields $\widetilde{ L} =
e^\psi \, L$ with $\psi \in \mathcal{ C}^\infty (\Omega)$ satisfying
$\sup_\Omega \, \big\vert \partial_x^\alpha \psi (x) \big\vert
\leqslant 1$ for $\vert \alpha \vert \leqslant 1$.
}\end{definition}

By a fine analysis of the degeneracies of $L$ and of the structure of
the Sussmann orbits of $\big\{ {\rm Re}\, L, \, {\rm Im}\, L \big\}$,
Hounie-Tavares were able to substantially generalize Ahlfors'
characterization.

\def\thetheorem{2.23}\begin{theorem}
{\rm (\cite{ ht1996, ht1997})} A relatively closed set $C \subset
\Omega$ is everywhere $(L, L^\infty )$-removable {\rm if and only if}
there is an open covering of $C$ by open sets $\Omega_j$, $j\in J$,
such that
\[
\text{\sf full-}L^\infty\text{\sf -cap}_L
(C\cap\Omega_j,\Omega_j)
=
0,
\]
for every $j\in J$.
\end{theorem}

On orbits of dimension one, $L$ behaves as a multiple of a real vector
field (one-dimensional behavior); on orbits of dimension two, $L$ has
the behavior of $\overline{ \partial}$ on a Riemann surface $\Sigma
\subset \mathcal{ O}$, but on the complement $\mathcal{ O} \backslash
\Sigma$ which is a union of curves with different endpoints along
which ${\rm Re}\, L$ and ${\rm Im}\, L$ are both tangent (degeneracy),
$L$ behaves again as a multiple of a real vector field
(one-dimensional behavior). As shown in~\cite{ ht1996} (main
Theorem~7.3 there) a relatively closed set $C \subset \Omega$ is
everywhere removable if and only if $C$ does not disconnect almost
every curve on which $L$ has one-dimensional behavior and furthermore,
the intersection of $C$ with almost every (reduced) orbit of dimension
two has zero analytic capacity for its natural holomorphic structure.

\def\theopenproblem{2.24}\begin{openproblem}
Study removability of a $\mathcal{ C }^\infty$ locally integrable
involutive structure of rank $\lambda \geqslant 2$ in terms of analytic
capacity.
\end{openproblem}

\subsection*{ 2.25.~Cartan-Thullen argument and a local continuity 
principle} 
The Behnke-Sommer {\sl Kontinuit\"atssatz}, alias {\sl Continuity
Principle}, states informally as follows (\cite{ sh1990}). Let
$(\Sigma_\nu )_{ \nu \in \N}$ be a sequence of complex manifolds with
boundary $\partial \Sigma_\nu$ contained in a domain $\Omega$ of
$\C^n$. If $\Sigma_\nu$ converges to a set $\Sigma_\infty \subset
\overline{ \Omega}$ and if $\partial \Sigma_\infty$ is contained
in $\Omega$, then every holomorphic function $f\in \mathcal{ O}
(\Omega)$ extends holomorphically to a neighborhood of the set
$\Sigma_\infty$ in $\C^n$. The geometries of the $\Sigma_\nu$ and of
$\Sigma_\infty$ have to satisfy certain assumptions in order that the
statement be correct; in addition, monodromy questions have to be
considered carefully. For applications to removable singularities 
in~\cite{ mp2006a}, the rigorous Theorem~2.27 below
is formulated, with the $\Sigma_\nu$ being embedded
analytic discs.

We denote by $z = (z_1, \dots, z_n)$ the complex
coordinates on $\C^n$ and by $\vert z\vert = \max_{ 1\leqslant
i\leqslant n} \, \vert z_i\vert$ the polydisc norm. If $E\subset \C^n$
is an arbitrary subset, for $\rho >0$, we denote by
\[
\mathcal{V}_\rho(E)
:=
\bigcup_{p\in E}\, 
\{z\in\C^n:
\vert z-p\vert 
< 
\rho\}
\]
the union of all open polydiscs of radius $\rho$ centered at points
of $E$.

\def\thelemma{2.26}\begin{lemma}
{\rm (\cite{ me1997})}
Let $\Omega$ be a nonempty domain of $\C^n$ and let $A: \overline{
\Delta} \to \Omega$, $A\in \mathcal{ O} (\Delta) \cap \mathcal{ C}^1 
(\overline{ \Delta})$, be an analytic disc contained in $\Omega$
having the property that there exist two constants $c$ and $C$ with $0
< c < C$ such that
\[
c\,
\vert\zeta_2 
-
\zeta_1\vert 
<
\vert A(\zeta_2) 
-
A(\zeta_1)\vert 
<
C\vert\zeta_2
-
\zeta_1\vert,
\]
for all distinct points $\zeta_1, \zeta_2 \in \overline{ \Delta}$.
Set 
\[
\rho 
:= 
\inf\big\{\vert z 
- 
A(\zeta)\vert:\, 
z\in\partial\Omega,\,
\zeta\in\partial\Delta\big\},
\]
namely $\rho$ is the polydisc distance between $A( \partial \Delta)$
and $\partial \Omega$, and set $\sigma := \rho \, c / 2C$. Then for
every holomorphic function $f\in \mathcal{ O }( \Omega )$, there
exists a holomorphic function $F\in \mathcal{ O } \big( \mathcal{
V}_\sigma ( A( \overline{ \Delta })) \big)$ such that $F = f$ on
$\mathcal{ V }_\sigma \big( A( \partial \Delta ) \big)$.
\end{lemma}

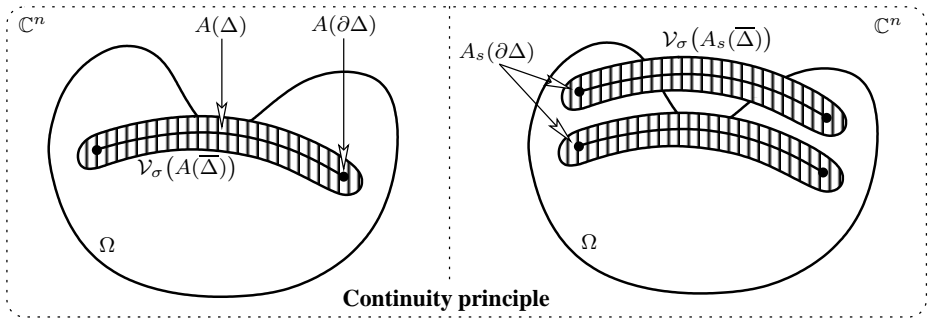
\begin{center}
\input continuity-principle.pstex_t
\end{center}

The inequalities involving $c$ and $C$ are satisfied for instance when
$A$ is $\mathcal{ C}^1$ embedding of $\overline{ \Delta}$ into
$\Omega$. Whereas $A ( \overline{ \Delta})$ is contained in $\Omega$,
the neighborhood $\mathcal{ V}_\sigma ( A( \overline{ \Delta}))$ is
allowed to go beyond. We do not claim that the two functions $f \in
\mathcal{ O} (\Omega)$ and $F \in \mathcal{ O} \big( \mathcal{
V}_\sigma ( A (\overline{ \Delta })) \big)$ stick together as a
holomorphic function globally defined in $\Omega \cup \mathcal{
V}_\sigma ( A (\overline{ \Delta }) )$. In fact, $\Omega \cap
\mathcal{ V }_\sigma ( A (\overline{ \Delta}) )$ may have several
connected components.

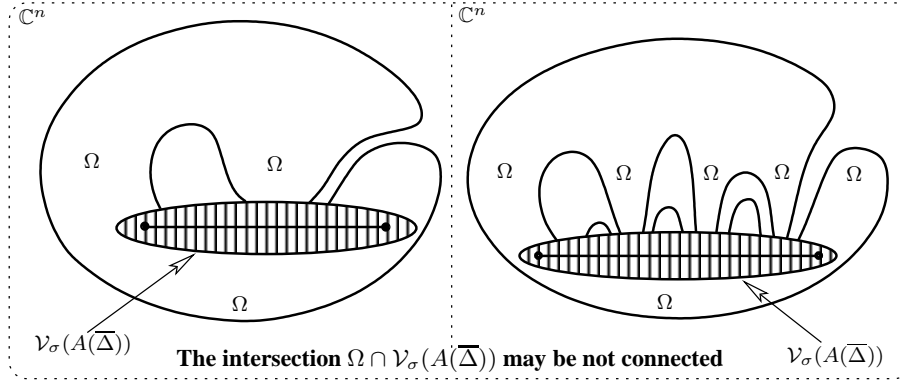
\begin{center}
\input G-teeth.pstex_t
\end{center}

In the geometric situations we encounter in~\cite{ me1997, mp1999,
mp2002, mp2006a}, after shrinking $\Omega$ somehow slightly to some
subdomain $\Omega'$, we shall be able to insure that the intersection
$\Omega' \cap \mathcal{ V }_\sigma ( A (\overline{ \Delta }))$ is
connected and that the union $\Omega ' \cup \mathcal{ V }_\sigma ( A
(\overline{ \Delta }))$ is still significantly ``bigger'' than
$\Omega$.

\proof[Proof of Lemma~2.26.]
Let $f\in \mathcal{ O} (\Omega)$. For $\zeta \in \overline{ \Delta}$
arbitrary, we consider the locally converging Taylor series
$\sum_{\alpha \in \N^n}\, f_\alpha\, (z - A(\zeta) )^\alpha$ of $f$ at
$A(\zeta)$. For $\rho'$ with $0 < \rho' < \rho$ arbitrarily close to
$\rho$, since $\mathcal{V }_{ \rho' }( A (\partial \Delta )) \Subset
\Omega$, the quantity
\[
M_{\rho'}(f)
:=
\sup\big\{
\vert f(z) \vert:\, 
z\in \mathcal{V}_{\rho'}(
A(\partial \Delta))
\big\}<
\infty,
\]
is finite (it may explode as $\rho' \to \rho$). Thus, Cauchy's
inequality on a polydisc of radius $\rho '$ centered at an arbitrary
point $A ( e^{ i \, \theta})$ of $\partial \Delta$ yields
\[
\frac{1}{\alpha!}\,
\left\vert 
\frac{\partial^{\vert\alpha\vert}f}{
\partial z^\alpha}\big(A(e^{i\,\theta})\big)
\right\vert
\leqslant
\frac{M_{\rho'}(f)}{{\rho'}^\alpha},
\]
uniformly for all $e^{ i\, \theta} \in \partial \Delta$. Then the
maximum principle applied to the function $\zeta \mapsto \frac{
\partial^{ \vert \alpha \vert} f}{ \partial z^\alpha} (A (\zeta))$
holomorphic in $\Delta$ provides the crucial inequalities
(Cartan-Thullen argument):
\[
\aligned
\vert f_\alpha \vert 
=
\frac{1}{\alpha!}\,
\left\vert
{\partial^{\vert \alpha \vert} f\over\partial z^\alpha}\, 
\big(A(\zeta)\big)\right\vert 
&
\leqslant 
{1\over \alpha!}\, \sup_{e^{i\,\theta}\in\partial\Delta}\, 
\left\vert
{\partial^{\vert \alpha \vert} f\over\partial z^\alpha}\,
\big(A(e^{i\,\theta})\big)
\right\vert
\\
&
\leqslant
\frac{M_{\rho'}(f)}{
{\rho'}^\alpha}.
\endaligned
\]
Consequently, the Taylor series of $f$ converges normally in the
polydisc $\Delta_\rho^n (A( \zeta ))$ of center $A( \zeta)$ and of
radius $\rho$, this being true for every $A ( \zeta) \in A( \overline{
\Delta})$. These series define a collection of holomorphic functions
$F_{ A(\zeta), \rho} \in \mathcal{ O}( \Delta_\rho^n ( A (\zeta)))$
parametrized by $\zeta \in \overline{ \Delta}$. We claim that the
restrictions of all these functions to the smaller polydiscs
$\Delta_\sigma^n (A (\zeta))$ stick together in a well defined
holomorphic function $F \in \mathcal{O } \big( \mathcal{ V}_\sigma ( A
( \overline{ \Delta})) \big)$.

Indeed, assume that two distinct points $\zeta_1 , \zeta_2 \in
\overline{ \Delta}$ are such that the intersection of the two small
polydiscs $\Delta_\sigma^n (A( \zeta_1)) \cap \Delta_\sigma^n (A
(\zeta_2))$ is nonempty, so $\vert A (\zeta_2) - A (\zeta_1) \vert <
2\, \sigma$. It follows that for every $\zeta$ belonging to the
segment $[\zeta_1, \zeta_2]$, we have:
\[
\vert\zeta 
-
\zeta_1\vert 
\leqslant 
\vert\zeta_2
- 
\zeta_1\vert 
<
\vert A(\zeta_2)
-
A(\zeta_1)\vert/c
<
2\,\sigma/c
\] 
whence
\[
\vert A(\zeta) 
- 
A(\zeta_1)\vert
< 
C\vert\zeta 
-
\zeta_1\vert 
< 
2\,C\sigma/c
= 
\rho.
\]
This means that the curved segment $A( [ \zeta_1, \zeta_2 ])$ is
contained in the connected intersection of the two large polydiscs
$\Delta_\rho^n (A( \zeta_1)) \cap \Delta_\rho^n (A( \zeta_2))$. In a
small neighborhood of $A (\zeta_1)$ and of $A (\zeta_2)$, the two
holomorphic functions $F_{ A(\zeta_1) ,\rho}$ and $F_{ A (\zeta_2),
\rho}$ coincide with $f$ by construction. Thanks to the principle of
analytic continuation, it follows that they even coincide with $f$ in
a thin connected neighborhood of the segment $A( [ \zeta_1, \zeta_2
])$. Again thanks to the principle of analytic continuation, $F_{
A(\zeta_1) ,\rho}$ and $F_{ A (\zeta_2), \rho}$ coincide in the
connected intersection $\Delta_\rho^n (A( \zeta_1)) \cap \Delta_\rho^n
(A( \zeta_2))$. It follows that they stick together to provide a well
defined function $F_{ A(\zeta_1), A(\zeta_2), \rho}$ that is
holomorphic in $\Delta_\rho^n (A( \zeta_1)) \cup \Delta_\rho^n (A(
\zeta_2))$. In conclusion, the restriction $F_{ A(\zeta_1),
A(\zeta_2), \rho} \big\vert_{ \Delta_\sigma (A (\zeta_1)) \cup
\Delta_\sigma (A (\zeta_2))}$ is holomorphic in the union of the two
small polydiscs $\Delta_\sigma (A (\zeta_1)) \cup \Delta_\sigma (A
(\zeta_2))$, whenever the intersection $\Delta_\sigma (A (\zeta_1))
\cap \Delta_\sigma (A (\zeta_2))$ is nonempty. This proves that all
the restrictions $F_{ A( \zeta), \rho} \big\vert_{ \Delta_\sigma^n
(A(\zeta ))}$ stick together in a well defined holomorphic function $F
\in \mathcal{ O} \big( \mathcal{ V}_\sigma ( A ( \overline{ \Delta }))
\big)$.
\endproof

In the next theorem (a local {\sl continuity principle} often used
in~\cite{ me1997, mp1999, mp2002, mp2006a}), $A_1 (\overline{ \Delta})
\subset \Omega$, but contrary to Lemma~2.26, $A_s( \overline{
\Delta})$ may well be not contained in $\Omega$ for $s < 1$;
nevertheless, the boundaries $A_s ( \partial \Delta)$ must always stay
in $\Omega$.

\def\thetheorem{2.27}\begin{theorem}
(\cite{ me1997})
Let $\Omega$ be a nonempty domain in $\C^n$ and let $A_s: \overline{
\Delta} \to \C^n$, $A_s \in \mathcal{ O} (\Delta) \cap \mathcal{ C}^1 
(\overline{ \Delta})$, be a one-parameter family of analytic
discs, where $s \in [0, 1]$. Assume that there
exist two constants $c_s$ and $C_s$ with $0 < c_s < C_s$ such that
\[
c_s\,
\vert\zeta_2 
-
\zeta_1\vert 
< 
\vert A(\zeta_2) 
- 
A(\zeta_1)\vert 
< 
C_s\,
\vert\zeta_2 
-
\zeta_1\vert,
\]
for all distinct points $\zeta_1, \zeta_2 \in \overline{
\Delta}$. Assume that $A_1 (\overline{ \Delta }) \subset \Omega$, set
\[
\rho_s
:=
\inf\big\{\vert z 
- 
A_s(\zeta)\vert:\, 
z\in\partial\Omega,\,
\zeta\in\partial\Delta\big\}, 
\]
namely $\rho_s$ is the polydisc distance between $A_s(\partial
\Delta)$ and $\partial \Omega$, and set $\sigma_s := \rho_s c_s /2C_s$.
Then for every holomorphic functions $f\in \mathcal{O }( \Omega)$, and
for all $s\in [0,1]$, there exist holomorphic functions $F_s\in
\mathcal{ O} \big( 
\mathcal{ V}_{ \sigma_s} (A_s (\overline{ \Delta})) \big)$ such
that $F_s =f$ in $\mathcal{ V}_{ \sigma_s} ( A_s( \partial \Delta))$.
\end{theorem}

\proof 
Let $\mathcal{ I} \subset [0, 1]$ be the connected set of real $s_0$
such that the statement is true for all $s$ with $s_0 \leqslant s
\leqslant 1$. By Lemma~2.26, we already know that $1\in \mathcal{
I}$. We want to prove that $\mathcal{ I} =[ 0, 1]$. It suffices to
prove that $\mathcal{ I}$ is both open and closed.

The fact that $\mathcal{ I}$ is closed follows by ``abstract
nonsense''. We claim that $\mathcal{ I}$ is also open. Indeed, let
$s_0 \in \mathcal{ I}$ and let $s_1< s_0$ be such that $A_{s_1 }(
\overline{ \Delta})$ is contained in $\mathcal{ V }_{ \sigma_{ s_0}}
(A_{ s_0}( \overline{ \Delta} ))$. Since $F_{ s_0}=f$ in $\mathcal{
V}_{ \sigma_{ s_0}} (A_{ s_0}( \partial \Delta ))$ and since the
polydisc distance between $A_{ s_1} (\partial \Delta)$ and $\partial
\Omega$ is equal to $\rho_{s_1}$, it follows as in the first part of
the proof of Lemma~2.26, that the Taylor series of $F_{ s_0}$ at
arbitrary points of the form $A_{ s_1}( \zeta)$, $\zeta \in \overline{
\Delta}$, converges in the polydisc $\Delta_{ \rho_{s_1 } }^n (A_{
s_1} (\zeta))$. This gives holomorphic functions $F_{ A(\zeta),
\rho_{ s_1}} \in \mathcal{ O} \big( \Delta_{ \rho_{ s_1} }^n (A_{ s_1}
(\zeta )) \big)$, for every $\zeta \in \overline{ \Delta}$. Reasoning
as in the second part of the proof of Lemma~2.26, we obtain a function
$F_{ s_1} \in \mathcal{ O} \big( \mathcal{ V}_{ \sigma_{ s_1}} ( A_{
s_1}( \overline{ \Delta })) \big)$ with $F_{s_1 } = f$ in $\mathcal{
V}_{ \sigma_{ s_1 }}( A_{s_1 }( \partial \Delta))$. This shows that
$\mathcal{ I}$ is open, as claimed and completes the proof.
\endproof

\subsection*{ 2.28.~Singularities as complex hypersurfaces}
Let $\Omega \subset \C^n$ ($n\geqslant 2$) be a domain. A typical
elementary singularity in $\Omega$ is just the zero set $Z_f := \{ f =
0 \}$ of a holomorphic function $f \in \mathcal{ O} (\Omega)$ since
for instance, the functions $1/ f^k$, $k\geqslant 1$, and $e^{ 1/f}$
are holomorphic in $\Omega \backslash Z_f$ and singular along
$Z_f$. Because $\C$ is algebraically closed, the closure in
$\overline{ \Omega}$ of such $Z_f$'s necessarily intersects $\partial
\Omega$. Early in the twentieth century, the italian mathematicians
Levi, Severi and B.~Segre (\cite{ se1932}) interpreted Hartogs'
extension theorem as saying that compact sets $K \subset \Omega$ are
removable, confirming the observation $\overline{ Z_f} \cap \partial
\Omega \neq \emptyset$.

\def\thedefinition{2.29}\begin{definition}{\rm
{\bf (i)} A relatively closed subset $C$ of a domain $\Omega \subset
\C^n$ is called {\sl removable} if the restriction map $\mathcal{ O}
(\Omega) \to \mathcal{ O} (\Omega \backslash C)$ is surjective.

\smallskip
{\bf (ii)} Such a set $C$ is called {\sl locally removable} if for
every $p\in C$, there exists an open neighborhood $U_p$ of $p$ in
$\Omega$ such that the restriction map $\mathcal{ O} \big( (U_p \cup
\Omega) \backslash C \big) \to \mathcal{ O} (\Omega \backslash C)$ is
surjective.

}\end{definition}

Under the assumption that $C$ is contained in a real submanifold of
$\Omega$, the general philosophy of removable singularities is that a
set too small to be a $Z_f$ (viz. a complex $(n- 1 )$-dimensional
variety) is removable. The following theorem collects five statements
saying that $C$ is removable provided it cannot contain any complex
hypersurface of $\Omega$. Importantly, our submanifolds $N$ of
$\Omega$ will always be assumed to be {\sl embedded}, namely for every
$p\in N$, there exist an open neighborhood $U_p$ of $p$ in $\Omega$
and a diffeomorphism $\psi_p : U_p \to \R^{ 2n}$ such that $\psi_p (N
\cap U_p) = \R^{ \dim N} \times \{ 0\}$.

\def\thetheorem{2.30}\begin{theorem}
Let $\Omega$ be a domain of $\C^n$ {\rm (}$n
\geqslant 2${\rm )} and let $C
\subset \Omega$ be a relatively closed subset. The restriction map
$\mathcal{ O} (\Omega) \to \mathcal{ O} (\Omega \backslash C)$ is
surjective, namely $C$ is removable, under each one of the following
five circumstances.
\begin{itemize}

\smallskip\item[{\bf (rm1)}]
$C$ is contained in a connected
submanifold $N \subset \Omega$ of
codimension $\geqslant 3$.

\smallskip\item[{\bf (rm2)}]
${\sf H}^{ 2n - 2} (C) = 0$.

\smallskip\item[{\bf (rm3)}]
$C$ is a relatively closed proper subset of a connected $\mathcal{
C}^2$ submanifold $N \subset \Omega$ of codimension $2$.

\smallskip\item[{\bf (rm4)}]
$C = N$ is a connected $\mathcal{ C}^2$ submanifold
$N \subset \Omega$ that is not a complex hypersurface
of $\Omega$.

\smallskip\item[{\bf (rm5)}]
$C$ is a closed subset of a connected
$\mathcal{ C}^2$ real hypersurface $M^1
\subset \Omega$ that does not contain any CR orbit of $M^1$.

\end{itemize}\smallskip
\end{theorem}

In {\bf (rm1)} and in {\bf (rm2)}, $C$ is in fact locally removable.
In {\bf (rm5)}, two kinds of CR orbits coexist: those of real
dimension $(2n -2)$, that are necessarily complex hypersurfaces, and
those of real dimension $(2n - 1)$, that are open subsets of $M$. It
is necessary to exclude them also. Indeed, if for instance $\Omega$ is
divided in two connected components $\Omega^\pm$ by a globally minimal
$M^1$, taking $C = M^1$, any locally constant function on $\Omega
\backslash M^1$ equal to two {\it distinct}\, constants $c^\pm$ in
$\Omega^\pm$ does not extend holomorphically through $M^1$. The proof
of the theorem is elementary\footnote{ Some refinements of the
statements may be formulated, for instance assuming in {\bf (rm3)} and
{\bf (rm4)} that $N$ is $\mathcal{ C}^1$ and has some metrically thin
singularities (\cite{ cst1994}).} and we
will present it in \S2.34 below, as a relevant preliminary to
Theorems~4.9, 4.10, 4.31 and~4.32, and to the main Theorem~1.7 
of~\cite{ mp2006a}).

We mention that under the assumption of local boundedness, more
massive singularities may be removed. An application of
Theorem~1.9{\bf (ii)} to several complex variables deserves to be
emphasized.

\def\thetheorem{2.31}\begin{theorem}
{\rm (\cite{ hp1970})}
If $\Omega \subset \C^n$ is a domain and if $C \subset \Omega$ is a
relatively closed subset satisfying ${\sf H}^{ 2n - 1} (C) = 0$, then
$C$ is removable for functions holomorphic in $\Omega \backslash C$
that are locally bounded in $\Omega$.
\end{theorem}

Following~\cite{ stu1989} and~\cite{ lu1990}, we now provide
variations on {\bf (rm2)}. Any global complex variety of codimension
one in $\C^n$ is certainly of infinite $(2n - 2)$-dimensional area.

\def\thetheorem{2.32}\begin{theorem}
{\rm (\cite{ stu1989})} Every closed set $C \subset \C^n$ satisfying
${\sf H}^{ 2n - 2} (C ) < \infty$ is removable for $\mathcal{ O}
(\C^n)$.
\end{theorem}

The result also holds true in the unit ball $\B_n$ of $\C^n$, provided
one computes the $(2n-2)$-dimensional Hausdorff measure with respect
to the distance function derived from the Bergman metric (\cite{
stu1989}). Also, if $\Sigma$ is an arbitrary complex $k$-dimensional
closed submanifold of $\C^n$, every closed subset $C \subset \Sigma$
with ${\sf H}^{ 2k-2} (C) < \infty$ is removable for $\mathcal{ O}
(\Sigma \backslash C)$ (\cite{ stu1989}).

A finer variation on the theme requires that ${\sf H}^{ 2n -2} (C
\cap R\, \B_n) $ does not grow too rapidly as a function of the radius
$R \to \infty$. For instance (\cite{ stu1989}), a closed subset $C
\subset \C^2$ that satisfies ${\sf H}^2 ( C \cap R\, \B_n) < c\, R^2$
for all large $R$ is removable, provided $c < \frac{ \pi^2}{ 4\,
\sqrt{ 2}}$. It is expected that $c < \pi$ is optimal, since the line
$L := \{ (z_1, 0)\}$ satisfies ${\sf H}^2 ( L \cap R\, B_2) = \pi \,
R^2$.

Yet another variation, raised in~\cite{ stu1989}, is as
follows. Consider a closed set $C$ in the complex projective space
$P_n (\C)$ ($n\geqslant 2$) such that the Hausdorff $(2n-2)$-dimensional
measure (with respect to the Fubini-Study metric) of $C$ is strictly
less than that of any complex algebraic hypersurface of $P_n (\C)$.
Is it true that $C$ is a removable singularity for meromorphic
functions, in the sense that every meromorphic function on $P_n (\C)
\backslash C$ extends meromorphically through $C$\,? This question was
answered by Lupacciolu.

Let $d_{\rm FS} (z, w)$ denote the geodesic distance between two
points $z, w \in P_n (\C)$ relative to the Fubini-Study metric and let
${\sf H}_{ \rm FS}^\ell$ denote the $\ell$-dimensional Hausdorff
measure in $P_n (\C)$ computed with $d_{ \rm FS}$. Given a nonempty
closed subset $C \subset P_n (\C)$, define:
\[
\rho(C)
:=
\frac{\max_{z\in P_n(\C)}d_{\rm FS}(z,C)}{
\max_{z,w\in P_n(\C)}d_{\rm FS}(z,w)}
=
\frac{\max_{z\in P_n(\C)}d_{\rm FS}(z,C)}{
{\rm diam}_{\rm FS}(P_n(\C))}
\leqslant
1.
\]
If the Fubini-Study metric is normalized so that ${\rm vol}\, P_n (\C)
= 1$, the $(2n-2)$-dimensional volume of an irreducible complex
algebraic hypersurface $\Sigma \subset P_n (\C)$ is equal to $\deg V$
and ${\sf H}_{ \rm FS}^{ 2n-2} (C) = (4/\pi)^{ n-1}\, (n-1)!\, \deg
V$. It follows that the minimum value of ${\sf H}^{ 2n -2} (\Sigma)$
is equal to $(4/\pi)^{ n-1}\, (n-1)!$ and is attained for $V$ equal to
any hyperplane of $P_n (\C)$. Let $\mathcal{ M}$ denote the sheaf of
meromorphic functions on $P_n (\C)$.

\def\thetheorem{2.33}\begin{theorem}
{\rm (\cite{ lu1990})} Let $C \subset P_n (\C)$ be a closed subset
such that
\[
{\sf H}_{\rm FS}^{2n-2}(C)
<
\big[
\rho(C)
\big]^{4n-4}\
(4/\pi)^{ n-1}\, (n-1)!
\]
Then the restriction map $\mathcal{ M}( P_n(\C)) \longrightarrow
\mathcal{ M} ( P_n (\C) \backslash C)$ is onto.
\end{theorem}

\subsection*{ 2.34.~Proof of Theorem~2.30} 
We claim that we may focus our attention only on {\bf (rm2)} and on
{\bf (rm5)}. Indeed, since a submanifold $N \subset \Omega$ of
codimension $\geqslant 3$ satisfies ${\sf H}^{ 2n - 2} (N) = 0$, {\bf
(rm1)} is a corollary of {\bf (rm2)}.

In both {\bf (rm3)} and {\bf (rm4)}, we may include $N$ in some
$\mathcal{ C}^2$ hypersurface $M^1$ of $\C^n$, looking like a thin
strip elongated along $N$. We claim that $C$ then does not contain
any CR orbit of any such $M^1$, so that {\bf (rm5)} applies. Indeed,
CR orbits of $M^1$ being of dimension $(2n-2)$ or $(2n-1)$ and $C$
being already contained in the $(2n-2)$-dimensional $N \subset M^1$,
it could only happen that $C = N = \Sigma$ identifies as a whole to a
connected (CR orbit) complex hypersurface $\Sigma \subset M^1$. But
this is excluded by the assumption that $C\neq N$ in {\bf (rm3)} and
by the existence of generic points in {\bf (rm4)}.

\smallskip

Firstly, we prove {\bf (rm2)}. We show that $C$ is locally
removable. Let $p\in C$ and let $B_p \subset \Omega$ be a small open
ball centered at $p$. By a relevant 
application of Proposition~1.2{\bf (5)},
one may verify that for almost every complex line $\ell$ passing
through $p$, the intersection $\ell \cap B_p \cap C$ is reduced to $\{
p\}$. Choose such a line $\ell_1$. Centering coordinates at $p$ and
rotating them if necessary, we may assume that $\ell_1 = \{ (z_1, 0,
\dots, 0)\}$, whence for $\varepsilon >0$ small and fixed, the disc
$A_\varepsilon (\zeta) := (\varepsilon \, \zeta, 0, \dots, 0)$
satisfies $A_\varepsilon (\partial \Delta) \cap C = \emptyset$.

Fix such a small $\varepsilon_0 >0$ and set $\rho_0 := {\rm dist}\,
\big( A_{ \varepsilon_0} (\partial \Delta), C \big) >0$. For $\tau =
(\tau_2, \dots, \tau_n) \in \C^{ n-1}$ satisfying $\vert \tau \vert <
\frac{ 1}{ 2} \, \rho_0$, set $A_{ \varepsilon_0, \tau} (\zeta) :=
(\varepsilon_0 \zeta, \tau_2, \dots, \tau_n)$. Letting $s\in [0, 1]$,
we interpolate between $A_{ \varepsilon_0, 0}$ and $A_{ \varepsilon_0,
\tau}$ by defining
\[
A_{\varepsilon_0,\tau,s}(\zeta)
:=
\big(
\varepsilon_0\zeta,s\,\tau_2,\dots,s\,\tau_n
\big).
\]
Since $\vert \tau\vert < \frac{ 1}{ 2}\, \rho_0$, these discs all
satisfy ${\rm dist}\, \big( A_{ \varepsilon_0, \tau, s} (\partial
\Delta), C \big) \geqslant \frac{ 1}{ 2} \, \rho_0$. Since the embedded disc
$A_{ \varepsilon_0, \tau, 1} (\overline{ \Delta})$ is $2$-dimensional
and since ${\sf H}^{ 2n-2} (C) = 0$, Proposition~1.2{\bf (5)} assures
that for almost every $\tau$ with $\vert \tau \vert < \frac{ 1}{ 2}\,
\rho_0$, its intersection with $C$ is empty. Thus, we may apply the
continuity principle Theorem~2.27, setting $c_s = \frac{ 1}{ 2}\,
\varepsilon_0$, $C_s = 2\, \varepsilon_0$ and $\rho_s := \frac{ 1}{
2}\, \rho_0$ uniformly for every $s \in [ 0, 1]$, whence $\sigma_s =
\frac{ 1}{ 8} \, \rho_0$ independently of the smallness of $\tau$: for
every $f\in \mathcal{ O} (\Omega \backslash C)$, there exists $F_0 \in
\mathcal{ O} \big( \mathcal{ V}_{ \frac{ \rho_0}{ 8} } (A_{
\varepsilon_0, \tau, 0} ( \overline{ \Delta })) \big)$ with $F_0 = f$
in $\mathcal{ V}_{ \frac{ \rho_0 }{ 8} } (A_{ \varepsilon_0, \tau, 0}
( \partial \Delta))$. But since $H^{ 2n - 2} (C) = 0$, for every
connected open set $\mathcal{ V} \subset \Omega$, the intersection
$\mathcal{ V} \cap (\Omega \backslash C)$ is also connected
(Proposition~1.2{\bf (4)}), so $F_0$ and $f$ stick together as a well
defined function holomorphic in $\Omega \cup \mathcal{ V}_{ \frac{
\rho_0}{ 8} } (A_{ \varepsilon_0, \tau, 0} ( \overline{ \Delta
}))$. If $\tau$ was chosen sufficiently small, it is clear that $p = 0
\in \C^n$ is absorbed in $\mathcal{ V}_{ \frac{ \rho_0}{ 8} } (A_{
\varepsilon_0, \tau, 0} ( \overline{ \Delta }))$, hence removable.

\smallskip

Secondly, we prove {\bf (rm5)}. Let $C \subset M^1$ containing no CR
orbit and define
\[
\mathcal{C}'
:=
\Big\{
C'\subset C\
\text{\rm closed},\
\forall\,f\in\mathcal{O}(\Omega\backslash C), 
\exists\,f'\in\mathcal{O}(\Omega\backslash C')\
\text{\rm with}\
f'\big\vert_{\Omega\backslash C}
=
f
\Big\}.
\]

\def\thelemma{2.35}\begin{lemma}
If $C_1', C_2' \in \mathcal{ C}'$, then $C_1 ' \cap C_2 ' \in
\mathcal{ C}'$.
\end{lemma}

\proof
Let $f_j' \in \mathcal{ O} (\Omega \backslash C_j')$, $j=1, 2$, with
$f_j' \big\vert_{ \Omega \backslash C} = f$. We claim that $f_1'$ and
$f_2'$ match up on $C \backslash (C_1' \cup C_2')$, hence define
together a holomorphic function $f_{ 12}' \in \mathcal{ O} \big(
\Omega \backslash (C_1' \cap C_2') \big)$ with $f_{ 12}' \big\vert_{
\Omega \backslash C} = f$. Indeed, choose an arbitrary point $p\in C
\backslash (C_1 ' \cup C_2')$. There exists a small open ball $B_p$
centered at $p$ with $B_p \cap (C_1 ' \cup C_2 ') = \emptyset$. Since
$f_1 ' = f_2 ' = f$ at least in the dense subset $B_p \backslash M^1$
of $M^1 \cap B_p$, by continuity of $f_1 '$ and of $f_2'$ in $B_p$,
necessarily $f_1 ' (p) = f_2 ' (p)$.
\endproof

Next, we define
\[
\widetilde{C}
:=
\bigcap_{C'\in\mathcal{C}'}\,C'.
\]
Intuitively, $\widetilde{ C}$ is the ``nonremovable core'' of $C$. By
the lemma, for every $f\in \mathcal{ O} ( \Omega \backslash C)$, there
exists $\widetilde{ f} \in \mathcal{ O} ( \Omega \backslash
\widetilde{ C})$ with $\widetilde{ f} \big\vert_{ \Omega \backslash C}
= f$. To prove {\bf (rm5)}, we must establish that $\widetilde{ C} =
\emptyset$. Reasoning by contradiction, we assume that $\widetilde{
C} \neq \emptyset$ and we apply to $C := \widetilde{ C}$ the lemma
below, which is in fact a corollary of Tr\'epreau's Theorem~2.4(V).
Of course, $\widetilde{ C}$ cannot contain any CR orbit of $M^1$.
Remind (\S1.27(III), \S4.9(III)) that for us, local CR orbits are {\it
not}\, germs, but local CR submanifolds of a certain small size.

\def\thelemma{2.36}\begin{lemma}
Let $M^1 \subset \C^n$ {\rm (}$n\geqslant 2${\rm )} be a $\mathcal{
C}^2$ hypersurface and let $C \subset M^1$ be a closed subset
containing no CR orbit of $M^1$. Then there exists at least one point
$p\in C$ such that for every neighborhood $V_p^1$ of $p$ in $M^1$, we
have{\rm :}
\[
V_p^1
\cap
\mathcal{O}_{CR}^{loc}(M^1,p)
\not\subset
C,
\]
and in addition, all such points $p$ are
locally removable.
\end{lemma}

Then $\widetilde{ f}$ extends holomorphically to a neighborhood of all
such points $p \in \widetilde{ C}$, contradicting the
definition of $\widetilde{ C}$, hence completing
the proof of {\bf (rm5)}.

\proof
If $\mathcal{ O}_{ CR}^{ loc} (M^1, q) \subset C$ for every $q\in C$,
then small complex-tangential curves issued from $q$ necessarily
remain in $\mathcal{ O}_{ CR}^{ loc} (M^1, q)$, hence in $C$, and
pursuing from point to point, global complex-tangential curves issued
from $q$ remain in $C$, whence $\mathcal{ O}_{ CR} (M^1, q) \subset C$,
contrary to the assumption.

So, let $p\in C$ with $\mathcal{ O}_{ CR}^{ loc} (M^1, p) \not\subset
C$. To pursue, we need that $M^1$ is minimal at $p$. Since $M^1$ is
a hypersurface, it might only happen that $\mathcal{ O}_{ CR}^{ loc}
(M^1, p)$ is a local complex hypersurface, a bad situation that has to
be changed in advance.

Fortunately, without altering the conclusion of the lemma (and of {\bf
(rm5)}), we have the freedom of perturbing the auxiliary hypersurface
$M^1$, leaving $C$ fixed of course. Thus, assuming that $\mathcal{
O}_{ CR} (M^1, p)$ is a complex hypersurface, we claim that there
exists a small (in $\mathcal{ C}^2$ norm) deformation $M_d^1$ of $M^1$
supported in a neighborhood of $p$ with $M_d^1 \supset C$ such that
$M_d^1$ is minimal at $p$.

Indeed, let $(q_k)_{ k\in \N}$ be a sequence of points tending to $p$
in $\mathcal{ O}_{ CR}^{ loc} (M^1, p)$ with $q_k \not\in C$. To
destroy the local complex hypersurface $\mathcal{ O }_{ CR}^{ loc}
(M^1, p)$, it suffices to achieve, by means of cut-off functions,
small bump-deformations of $M$ centered at all the $q_k$; it is easy
to write the technical details in terms of a local graphed
representation $v = \varphi^1 (z, u)$ for $M^1$. Outside a small
neighborhood of the union of the $q_k$, $M_d^1$ coincides with $M^1$.
Then the resulting $M_d^1$ is necessarily minimal at $p$, since if it
where not, the uniqueness principle\footnote{ Minimalization at a
point takes strong advantage of the rigidity of complex hypersurfaces
in this argument. } for complex
manifolds would force $\mathcal{ O}_{ CR}^{ loc} (M^1, p) = \mathcal{
O}_{ CR}^{ loc} (M_d^1, p)$, but $M_d^1$ does not contains the $q_k$.

So we can assume that $M^1$ is minimal at every point $p\in C$ at
which $\mathcal{ O}_{ CR}^{ loc} (M^1, p) \not \subset C$. Let $B_p
\subset \C^n$ be a small open ball centered at $p$ with $B_p \cap M^1
\subset \mathcal{ O}_{ CR}^{ loc} (M^1, p)$. We will show that
$\mathcal{ O} (\Omega \backslash C)$ extends holomorphically to
$B_p$. By assumption, $B_p \cap M^1 \not\subset C$, hence $B_p
\backslash C$ is connected, a fact that will insure monodromy.

Fixing $v_p \in T_p\C^n \backslash \{ 0\}$ with $v_p \not\in T_p M^1$,
we consider the global translations
\[
M_s^1
:=
M^1
+
s\,v_p,
\ \ \ \ \ \ \ \ \ 
s\in\R,
\]
of $M^1$. Let $f\in \mathcal{ O} (\Omega \backslash C)$ be
arbitrary. For small $s\neq 0$, $M_s^1 \cap B_p$ does not intersect
$M^1$, hence the restriction $f\vert_{ M_s^1 \cap B_p}$ is a
$\mathcal{ C}^2$ CR function on $M_s^1 \cap B_p$ (but $f\vert_{ M_0^1
\cap B_p}$ has possible singularities at points of $C\cap B_p$).

With $U_p := M^1 \cap B_p$, Theorem~2.4(V) says that $\mathcal{ C}_{
CR}^0 (U_p)$ extends holomorphically to some one-sided neighborhood
$\omega_p^\pm$ of $M^1$ at $p$. Reorienting if necessary, we may
assume that the extension side is $\omega_p^-$ and that $p + s\, v_p
\in B_p^+$ for $s > 0$ small. The statement and the proof of
Theorem~2.4(V) are of course invariant by translation. Hence $\mathcal{
C}_{ CR}^0 (U_p + s\, v_p)$ extends holomorphically to $\omega_p^- +
s\, v_p$, for every $s >0$. It is geometrically clear that for $s >0$
small enough, $\omega_p^- + s\, v_p$ contains $p$. Thus $f
\big\vert_{ M_s^1 \cap B_p}$ extends holomorphically to a neighborhood
of $p$ for such $s$. Monodromy of the extension follows from the fact
that $B_p' \backslash C$ is connected for every open ball $B_p'$
centered at $p$. This completes the proof of the lemma.
\endproof

\subsection*{ 2.37.~Removability and extension of complex hypersurfaces}
Let $\Omega \subset \C^n$ be a domain. Theorem~2.30
{\bf (rm4)} shows that a
connected $2$-codimensional submanifold $N \subset \Omega$ is
removable provided it is not a complex hypersurface, or equivalently,
is generic somewhere. Conversely, assume that $\Omega$ is
pseudoconvex and let $H \subset \Omega$ be a (not necessarily
connected) closed complex hypersurface. Then $\Omega \backslash H$ is
(obviously) locally pseudoconvex at every point, hence the
characterization of domains of holomorphy yields a
function $f \in \mathcal{ O} (\Omega \backslash H)$ whose domain of
existence is exactly $\Omega \backslash H$. Thus, $H$ is
nonremovable. But in a nonpseudoconvex domain, closed complex
hypersurfaces may be removable.

\def\theexample{2.38}\begin{example}{\rm
For $\varepsilon >0$ small, consider the following nonpseudoconvex
subdomain of $\B_2$, defined as the union of a spherical shell
together with a thin rod of radius $\varepsilon$
directed by the $y_2$-axis:
\[
\Omega_\varepsilon
:=
\big\{1/2<
\vert z_1\vert^2
+
\vert z_2\vert^2<1
\big\}
\bigcup
\left(
\B_2 \cap
\big\{
x_1^2+y_1^2+x_2^2<\varepsilon^2
\big\}
\right).
\] 
Then the intersection of the $z_1$-axis with the spherical shell, 
namely
\[
H 
:= 
\big\{ 
(z_1,0):\,
1/2 <\vert z_1\vert<1
\big\},
\] 
is a relatively closed complex hypersurface of $\Omega_\varepsilon$
homeomorphic to an open annulus. We claim that $H$ is removable.

Indeed, applying the continuity principle along discs parallel to the
$z_1$-axis, $\mathcal{ O} (\Omega_\varepsilon)$ extends
holomorphically to $\B_2 \backslash \{ z_2 = 0 \}$. Since the open
small disc $\{ (z_1 , 0 ) : \, \vert z_1 \vert < \varepsilon \}$,
considered as a subset of the closed complex hypersurface $\widetilde{
H}$ of $\B_2$ defined by
\[
\widetilde{H}
:=
\{(z_1, 0):\, 
\vert z_1 \vert < 1\} 
\]
is contained in the thin rod, hence
in $\Omega_\varepsilon$, Theorem~2.30 {\bf (rm3)} finishes to show
that 
\[
{\sf E}( \Omega_\varepsilon \backslash H)
= 
{\sf E}(\Omega_\varepsilon)
=
\B_2.
\] 
In such an example, we point out that the closed complex hypersurface
$H \subset \Omega_\varepsilon$ extends as the closed complex hypersurface
$\widetilde{ H } \subset {\sf E} (\Omega_\varepsilon \backslash H)$
but that the intersection
\[
\widetilde{ H } \cap \Omega_\varepsilon
=
\big\{(z_1,0):\, 
\vert z_1\vert 
< 
\varepsilon\big\}
\bigcup
\big\{(z_1,0):\,1/2<\vert z_1\vert<1\big\}
\]
is strictly bigger than $H$.

}\end{example}

\def\theproblem{2.39}\begin{problem}
Understand which relatively closed complex hypersurfaces of a general
domain $\Omega \subset \C^n$ are removable.
\end{problem}

We thus consider a (possibly singular and reducible) closed complex
hypersurface $H$ of $\Omega$. Basic properties of complex analytic
sets (\cite{ ch1991}) insure that $H = \bigcup_{ j\in J} \, H_j$
decomposes into at most countably many closed complex hypersurfaces
$H_j \subset \Omega$ that are irreducible.

\def\thedefinition{2.40}\begin{definition}{\rm
We say that $H_j$ allows an {\sl $H$-compatible extension} to ${\sf
E}( \Omega)$ if there exists an irreducible closed complex hypersurface
$\widetilde{ H }_j$ of ${\sf E}( \Omega)$ extending $H_j$ in the sense
that $H_j \subset \widetilde{ H }_j \cap \Omega$ whose intersection
with $\Omega$ remains contained in $H$:
\[
\widetilde{ H }_j \cap
\Omega
\subset 
\bigcup_{j'\in J}\,H_{j'}.
\]
}\end{definition}

The principle of analytic continuation for irreducible complex
analytic sets (\cite{ ch1991}) assures that $H_j$ has at most one
$H$-compatible extension. On the other hand,
$\widetilde{H}_j$ may be an $H$-compatible extension of several $H_{
j'}$. In the above example, the removable
annulus $H$ had no $H$-compatible extension to ${\sf E}(
\Omega)$.

\def\thetheorem{2.41}\begin{theorem}
{\rm (\cite{ dl1977, japf2000})} Let $\Omega \subset \C^n$ {\rm
(}$n\geqslant 2${\rm )} be a domain and let $H = \bigcup_{ j\in J}\,
H_j$ be a closed complex hypersurface of $\Omega$, decomposed into
irreducible components $H_j$. Set
\[
J_{\rm comp}:=
\big\{j\in J:\,H_j\
\text{\rm allows an}\
H\text{\rm -compatible extension}\
\widetilde{H}_j\
\text{\rm to}\
{\sf E}(\Omega)
\big\}.
\]
Then 
\[
{\sf E}(\Omega\backslash H)
=
{\sf E}(\Omega)\Big\backslash
\bigcup_{j\in J_{\rm comp}}\,
\widetilde{H}_j.
\]
In particular, $H$ is removable {\rm (}resp. nonremovable{\rm )} if
and only if $J_{\rm comp} = \emptyset$ {\rm (}resp. $J_{ \rm comp}
\neq \emptyset${\rm )}.
\end{theorem}

This statement was obtained after a chain of generalizations
originating from the classical results of Hartogs \cite{ ha1909} and of
Oka \cite{ ok1934}. In \cite{ gr1956} it was proved for the case that
$H$ is of the form $\Omega \cap \widetilde{ H}$, $\widetilde{ H}
\subset{ \sf E}( \Omega)$, and in \cite{ nis1962} under the additional
assumption that ${\sf E}( \Omega \backslash H)$ is a subset of ${\sf
E}( \Omega)$ ({\it a priori}, it is only a set {\it over}\, ${\sf E}(
\Omega)$). Actually Theorem 2.41 was stated in~\cite{ dl1977} even for
Riemann domains $\Omega$. But it was remarked in~\cite{ japf2000}
(p.~306) that the proof in~\cite{ dl1977} is complete only if the
functions of $\mathcal{ O}( \Omega)$ separate the points of $\Omega$,
{\it i.e.}~if $\Omega$ can be regarded as a subdomain of ${\sf
E}( \Omega)$. Actually the proof in~\cite{ dl1977} starts from the
special case where $\Omega$ is a Hartogs figure, which can be treated
by a subtle geometric examination. Then a localization argument shows
that extension of hypersurfaces which are singularity loci of
holomorphic functions cannot stop when passing from $\Omega$ to ${\sf
E}(\Omega)$. But in the nonseparated case, the global effect of
identifying points of $\Omega$ interferes nastily, and it is unclear
how to justify the localization argument. The final step for general
Riemann domains was achieved by the second author by completely
different methods.

\def\thetheorem{2.42}\begin{theorem}
{\rm (\cite{ po2002})} Let $\pi: X \rightarrow \C^n$ be an arbitrary
Riemann domain, and let $H \subset X$ be a closed complex
hypersurface. Denote by $\alpha : X \rightarrow {\sf E}( X )$ the
canonical immersion of $X$ into ${\sf E}( X)$. Then there is a closed
complex hypersurface $\widetilde{ H}$ of ${\sf E} (X)$ with $\alpha^{
- 1}( \widetilde{ H } ) \subset H$ such that
\[
{\sf E}(X\backslash H)
=
{\sf E }(X) 
\big\backslash 
\widetilde{ H}.
\]
\end{theorem}

Let us briefly sketch the main idea of the proof (\cite{ po2002}).
The essence of the argument is to reduce extension of
hypersurfaces to that of meromorphic functions.
For every pseudoconvex Riemann domain $\pi: X \rightarrow\C^n$, there
exists $f \in \mathcal{ O}( X) \cap L^2(X)$ having $X$ as domain of
existence whose growth is controlled by some power of the polydisc
distance to the abstract boundary $\breve{\partial} X$. At boundary
points where $\breve{ \partial} X$ can be locally identified with a
complex hypersurface, $f$ has just a pole of positive order. One can
deduce that those hypersurfaces $H$ of $X$ along which some
holomorphic function on $X \backslash H$ becomes singular can be
represented as the polar locus of some meromorphic function $g$
defined in $X$. But $g$ extends meromorphically to ${\sf E} (X)$, and
the polar locus of the extension yields the desired extension of
$H$. 
\section*{ \S3.~Hulls and removable singularities at the boundary}

\subsection*{ 3.1.~Motivations for removable singularities at the
boundary} As already observed in Section~1, 
beyond the harmonious realm of pseudoconvexity, the general
problem of understanding compulsory holomorphic (or CR) extension is
intrinsically rich and open. Some
elementary Baire category arguments show that most
domains are not pseudoconvex, most CR manifolds have nontrivial
disc-envelope, and most compact sets have nonempty essential
polynomial hull. Hence, the Grail for the theory of holomorphic
extension would comprise:

\begin{itemize}

\smallskip\item[$\bullet$]
a geometric and constructive view of the envelope of holomorphy of
most domains, following the Behnke-Sommer Kontinuit\"atssatz and
Bishop's philosophy;

\smallskip\item[$\bullet$]
a clear correspondence between function-theoretic techniques, for
instance those involving $\overline{ \partial}$ arguments, and
geometric techniques, for instance those involving
families of complex analytic varieties.

\end{itemize}\smallskip

Several applications of the study of envelopes of holomorphy appear,
for instance in the study of boundary regularity of solutions of the
$\overline{ \partial}$-complex, in the complex Plateau problem, in the
study of CR mappings, in the computation of polynomial hulls and in
removable singularities, the topics of this Part~VI and of~\cite{
mp2006a}.

In the 1980's, rapid progress in the understanding of the boundary
behavior of holomorphic functions led many authors to study the
structure of singularities up to the boundary. In \S2.28, we discussed
removability of relatively closed subsets $C$ of domains $\Omega
\subset \C^n$, {\it i.e.}~the problem whether $\mathcal{ O} (\Omega)
\to \mathcal{ O} (\Omega \backslash C)$ is surjective. Typically $C$
was supposed to be lower-dimensional and its geometry near $\partial
\Omega$ was irrelevant. Now we assume $\Omega$ to be bounded in $\C^n$
($n \geqslant 2$) and we consider compact subsets $K$ of $\overline{
\Omega}$, possibly meeting $\partial \Omega$.

\def\theproblem{3.2}\begin{problem}
Find criteria of geometric, or of function-theoretic nature, assuring
that the restriction map $\mathcal{ O} (\Omega) \to \mathcal{ O}
(\Omega \backslash K)$ is surjective.
\end{problem}

If $K \subset \overline{ \Omega} \cap\partial \Omega=\emptyset$ and
$\Omega\backslash K$ is connected, surjectivity follows from the
Hartogs-Bochner extension Theorem~1.9(V). Since this
theorem even gives extension of CR functions on $\partial\Omega$, it
seems reasonable to ask for holomorphic extension of CR functions on
$\partial \Omega \backslash K$, and then it is natural to assume that
$K$ is contained in $\partial \Omega$.
Hence the formulation of a second trend of questions\footnote{ CR
distributions may also be considered, but in the sequel, we shall
restrict considerations to continuous and integrable CR functions.}.

\def\theproblem{3.3}\begin{problem}
Let $K$ is a compact subset of $\partial \Omega$ such that
$\partial \Omega \backslash K$ is a hypersurface of class at least
$\mathcal{ C}^1$. Understand under which circumstances CR functions 
of class $\mathcal{ C}^0$ or $L_{ loc}^{ \sf p}$ on $\partial
\Omega \backslash K$ extend holomorphically to $\Omega$.
\end{problem}

A variant of these two problems consists in assuming that functions
are holomorphic in some thin (one-sided) neighborhood of $\partial
\Omega \backslash K$. In all the theorems that will be surveyed 
below, it appears that the thinness of the (one-sided) neighborhood of
$\partial \Omega \backslash K$ has no influence on extension, as in
the original Hartogs theorem. In this respect, it is of interest to
immediately indicate the connection of these two problems with the
problem of determining certain envelopes of holomorphy.

In the second problem, the hypersurface $\partial \Omega \backslash K$
is often globally minimal, a fact that has to be verified or might be
one of the assumptions of a theorem. For instance, several
contributions deal with the paradigmatic case where $\partial \Omega$
is at least $\mathcal{ C}^2$ and strongly pseudoconvex (hence
obviously globally minimal). Then thanks to the elementary Levi-Lewy
extension theorem (Theorem~1.18(V), lemma~2.2(V) and \S2.10(V)), there
exists a one-sided neighborhood $\mathcal{ V} (\partial \Omega
\backslash K)$ of $\partial \Omega \backslash K$ contained in $\Omega$
to which both $\mathcal{ C}_{ CR}^0 (\partial \Omega \backslash K)$
and $L_{ loc, CR}^{ \sf p} (\partial \Omega \backslash K)$ extend
holomorphically. The size of $\mathcal{ V} (\partial \Omega \backslash
K)$ depends only on the local geometry of $\partial \Omega$, because
$\mathcal{ V} (\partial \Omega \backslash K)$ is obtained by gluing
small discs (Part~V). In fact, an inspection of the proof of the
Levi-Lewy extension theorem together with an application of the
continuity principle shows also that the envelope of holomorphy of any
thin one-sided neighborhood $\mathcal{ V}' (\partial \Omega \backslash
K)$ (not necessarily contained in $\Omega$\,!) contains a one-sided
neighborhood $\mathcal{ V} (\partial \Omega \backslash K)$ of
$\partial \Omega \backslash K$ contained in the pseudoconvex domain
$\Omega$ that has a fixed, incompressible
size\footnote{ 
To be rigorous: for every holomorphic function $f\in \mathcal{ O}
\big( \mathcal{ V}' (\partial \Omega \backslash K) \big)$, there
exists a holomorphic function $F \in \mathcal{ O} \big( \mathcal{ V}
(\partial \Omega \backslash K) \big)$ that coincides with $f$ in a
possibly smaller one-sided neighborhood $\mathcal{ V}'' (\partial
\Omega \backslash K) \subset \mathcal{ V}' (\partial \Omega \backslash
K)$. Details of the proof (involving
a deformation argument) will not be provided here
({\it see}~\cite{ me1997, jo1999a}).
}.

As they are formulated, the above two problems turn out to be slightly
too restrictive. In fact, the final goal is to understand the
envelope ${\sf E} \big( \mathcal{ V} ( \partial \Omega \backslash K )
\big)$, or at least to describe some significant part of ${\sf E}
\big( \mathcal{ V} ( \partial \Omega \backslash K ) \big)$ lying above
$\Omega$. Of course, the question to which extent is the geometry of
${\sf E} \big( \mathcal{ V} (\partial \Omega \backslash K) \big)$
accessible (constructively speaking) depends sensitively on the shape
of $\Omega$. Surely, the strictly pseudoconvex case is the easiest and
the best understood up to now. In what follows we will encounter
situations where ${\sf E} \big( \mathcal{ V} ( \partial \Omega
\backslash K ) \big)$ contains $\Omega \backslash \widehat{ K}$, for
some subset $\widehat{ K} \subset \overline{ \Omega}$ defined in
function-theoretic terms and depending on $K \subset \partial
\overline{ \Omega}$. We will also encounter situations where ${\sf E}
\big( \mathcal{ V} ( \Omega \backslash K) \big)$ is necessarily
multisheeted over $\C^n$. In this concern, we will see a very striking
difference between the complex dimensions $n = 2$ and $n \geqslant 3$.

\smallskip

In the last two decades, a considerable interest has been devoted to a
subproblem of these two problems, especially with the objective of
characterizing the singularities at the boundary that are removable.

\def\thedefinition{3.4}\begin{definition}{\rm
In the second Problem~3.3, the compact subset $K \subset \partial
\Omega$ is called {\sl CR-removable} if for every CR function $f\in
\mathcal{ C}_{ CR}^0 (\partial \Omega \backslash K)$ (resp. $f\in L_{
loc, CR}^{ \sf p} (\partial \Omega \backslash K)$), there exists $F \in
\mathcal{ O} (\Omega) \cap \mathcal{ C}^0 (\overline{ \Omega}
\backslash K)$ (resp. $F \in \mathcal{ O} (\Omega) \cap H_{ loc}^{ \sf
p} (\overline{ \Omega} \backslash K)$) with $F \vert_{ \partial \Omega
\backslash K} = f$ (resp. locally at every point $p \in \partial
\Omega \backslash K$, the $L_{ loc, CR}^{ \sf p}$ boundary value of
$F$ equals $f$).
}\end{definition}

Before exposing and surveying some major results we would like to
mention that a complement of information and different
approaches may be found in the two surveys \cite{
stu1993, cst1994}, in the two monographs \cite{ ky1995, lt1997} and in
the articles~\cite{ stu1981, lt1984, lu1986, lu1987, jo1988, lt1988,
stu1989, ky1990, ky1991, stu1991, fs1991, jo1992, kn1993, du1993,
ls1993, lu1994, ac1994, jo1995, kr1995, jo1999a, jo1999b,
js2000, jp2002, js2004}.

\subsection*{ 3.5.~Characterization of removable sets contained
in strongly pseudoconvex boundaries} Taking inspiration from the
pivotal Oka theorem, one of the goals of the study of removable
singularities (\cite{ stu1993}) is to characterize removability in
function-theoretically significant terms, especially in terms of
convexity with respect to certain spaces of functions. In the very
beginnings of {\sl Several Complex Variables}, polynomial convexity
appeared in connexion with holomorphic approximation. According to the
Oka-Weil theorem (\cite{ aw1998}), functions that are holomorphic in
some neighborhood of a polynomially convex compact set $K \subset
\C^n$ may be approximated uniformly by polynomials. Later on,
holomorphic convexity appeared to be central in Stein theory (\cite{
ho1973}), one of the seminal frequently used idea being to encircle
convex compact sets by convenient analytic polyhedra.

The notion of convexity adapted to our pruposes is the following. By
$\mathcal{ O} (\overline{ \Omega})$, we denote the ring of functions
that are holomorphic in some neighborhood of the closure $\overline{
\Omega}$ of a domain $\Omega \subset \C^n$. As in the concept of
germs, the neighborhood may depend on the function.

\def\thedefinition{3.6}\begin{definition}{\rm
Let $\Omega \Subset \C^n$ be a bounded domain and let $K \subset
\overline{ \Omega}$ be a compact set. The {\sl $\mathcal{ O} (
\overline{ \Omega})$-convex hull} of $K$ is
\[
\widehat{K}_{\mathcal{O}(\overline{\Omega})}
:=
\Big\{
z\in\overline{\Omega}:\
\vert g(z)\vert\leqslant\max_{w\in K}\,
\vert g(w)\vert\
\text{\rm for all}\
g\in\mathcal{O}(\overline{\Omega})
\Big\}.
\]
If $K = \widehat{ K}_{ \mathcal{ O} (\overline{ \Omega })}$, then $K$
is called {\sl $\mathcal{ O} (\overline{ \Omega })$-convex}.
}\end{definition}

If $\Omega$ is strongly pseudoconvex, a generalization of the Oka-Weil
theorem shows that every function which is holomorphic in a
neighborhood of some $\mathcal{ O}( \overline{ \Omega})$-convex
compact set $K \subset \overline{ \Omega }$ may be approximated
uniformly on $K$ by functions of $\mathcal{ O}(\overline{ \Omega})$
(nevertheless, for nonpseudoconvex domains, this approximation
property fails\footnote{
Indeed, consider for instance the Hartogs figure $\Omega := \big\{
\vert z_1\vert <1, \, \vert z_2\vert < 2 \big\} \cup \big\{ 1
\leqslant \vert z_1 \vert < 2,1 < \vert z_2 \vert <2 \big\}$ in
$\C^2$. Then the annulus $K = \{z_1 =1, 1 \leqslant \vert z_2\vert
\leqslant 2)\} \subset \partial \Omega$ is $\mathcal{ O}( \overline{
\Omega })$-convex. We claim that the function $g := 1/ z_2$,
holomorphic in a neighborhood of $K$, cannot be approximated on $K$ by
functions $f\in\mathcal{ O}(\overline{ \Omega})$. Indeed, by Hartogs
extension $\mathcal{ O}(\overline{ \Omega})=\mathcal{ O} \big( \{\vert
z_1 \vert \leqslant 2, \, \vert z_2 \vert \leqslant 2\} \big)$, which
implies that every $f \in \mathcal{ O} ( \overline{ \Omega} )$ has to
satisfy the maximum principle on the disc $\{z_1 = 1, \vert z_2 \vert
\leqslant 2) \}\supset K$. Rounding off the corners, we get an example
with $\partial \Omega \in \mathcal{ C}^\infty$.}).

\smallskip

We may now begin with the formulation of a seminal theorem due to
Stout that inspired several authors. We state the CR version, due to
Lupacciolu\footnote{Said differently, the envelope of holomorphy of an
arbitrarily thin (interior) one-sided neighborhood of $\partial \Omega
\backslash K$ is one-sheeted and identifies with $\Omega$.}.

\def\thetheorem{3.7}\begin{theorem}
{\rm (\cite{ stu1981, lu1986, stu1993})} In complex dimension
$n=2$, a compact subset $K$ of a $\mathcal{ C }^2$ strongly
pseudoconvex boundary $\partial \Omega \Subset \C^2$ is CR-removable
{\rm if and only if} it is $\mathcal{ O } (\overline{ \Omega }
)$-convex.
\end{theorem}

The ``only if'' part is the easiest, relies on a lemma due to
S{\l}odkowski (\cite{ rs1989, stu1993}) and will be presented
after Lemma~3.11. Let us sketch the beautiful key idea of the ``if'' part
(\cite{ stu1981, lu1986, stu1993, po1997}).

From \S1.7(V), remind the expression of the Bochner-Martinelli kernel:
\[
{\sf BM}
(\zeta,z)
=
\frac{1}{
(2\pi i)^2\vert\zeta-z\vert^4}
\left[
\overline{(\zeta_2-z_2)}\,d\bar\zeta_1
-
\overline{(\zeta_1-z_1)}\,d\bar\zeta_2
\right] 
\wedge d\zeta_1 \wedge d\zeta_2.
\]
Let $\mathcal{ M} \subset \C^2$ be a thin strongly pseudoconvex
neighborhood of $\overline{ \Omega}$. By means of a fixed function $g
\in \mathcal{ O} ( \mathcal{ M})$, it is possible to construct some
explicit primitive of ${\sf BM}$ as follows. This idea goes back to
Martinelli and has been exploited by Stout, Lupacciolu, Leiterer,
Laurent-Thi\'ebaut, Kytmanov and others. By a classical result (\cite{
hele1984}), $g$ admits a Hefer decomposition
\[
g(\zeta)-g(z)
=
g_1(\zeta,z)[\zeta_1-z_1]
+
g_2(\zeta,z)[\zeta_2-z_2],
\] 
with $g_1, g_2 \in \mathcal{ O} ( \mathcal{ M } \times \mathcal{ M }
)$. Then a direct calculation shows that for $z\in \mathcal{ M}$
fixed, the $(0, 2)$-form
\[
\Theta_{g,z}(\zeta)
=
\frac{g_2(\zeta,z)\overline{
(\zeta_1-z_1)}-g_1(\zeta,z)\overline{
(\zeta_2-z_2)}}
{(2\pi i)^2\,
\vert\zeta-z\vert^2\,
\big[g(\zeta)-g(z)\big]}
\ d\zeta_1\wedge d\zeta_2
\]
satisfies 
\[
\overline{ \partial}_\zeta \Theta_{g,z} (\zeta) = 
d_\zeta\Theta_{g,z} (\zeta)
=
{\sf BM}
(\zeta, z),
\]
on $\{\zeta \in \mathcal{ M}:\, g (\zeta) \neq g (z) \}$, {\it i.e.}
provides a primitive of ${\sf BM}$ outside some thin set. In $\C^n$
for $n\geqslant 3$, there is also a similar explicit primitive.

Let $K \subset \partial \Omega$ be as in Theorem~3.7 and fix $z \in
\overline{ \Omega} \big\backslash K$. By $\mathcal{ O} (\overline{
\Omega })$-convexity of $K$, there exists $g \in \mathcal{ O}(
\overline{ \Omega})$ with $g(z) = 1$ and $\max_{ w\in K} \vert g (w)
\vert < 1$. After a slight elementary modification of $g$ (\cite{
jo1995, po1997}), one can insure that the set $\big\{ w\in \mathcal{
M}: \, \vert g (w) \vert = 1 \big\}$ is a geometrically smooth
$\mathcal{ C}^\omega$ Levi-flat hypersurface of $\mathcal{ M}$
transverse to $\partial \Omega$. Then the region $\Omega_g := \Omega
\cap \big\{ \vert g\vert >1 \big\}$ has piecewise smooth connected
boundary
\[
\partial \Omega_g = \big( \partial \Omega \cap \{ \vert g
\vert > 1\} \big) \bigcup \big( \Omega \cap \{ \vert g
\vert = 1\} \big)
\]
and its closure $\overline{ \Omega}_g$ in $\C^2$ does not intersect
$K$. 

Let $f$ be an arbitrary continuous CR function on $\partial \Omega
\backslash K$. Supposing for a while that $f$ already enjoys a
holomorphic extension $F \in \mathcal{ O}( \Omega) \cap \mathcal{ C}^0
\big( \overline{ \Omega } \backslash K \big)$, the Bochner-Martinelli
representation formula then provides for every $z \in \Omega_g$
the value
\[
F(z)
=
\int_{\partial\Omega_g}\,
f(\zeta)\,
{\sf BM}(\zeta,z).
\]
Decomposing $\partial \Omega_g$ as above and 
using the primitive $\Theta_{ g,z}$,
we may write
\[
F(z)
=
\int_{\partial\Omega\cap\{\vert g\vert >1\}}\,f(\zeta)\,
{\sf BM}(\zeta,z)
+
\int_{\Omega\cap\{\vert g\vert =1\}}\,f(\zeta)\,
d_\zeta\Theta_{g,z}(\zeta).
\] 
Supposing $f\in \mathcal{ C}^1$ and applying Stokes' theorem\footnote{
In the general case $f \in \mathcal{ C}^0$, one shrinks slightly
$\Omega_g$ inside $\Omega$, rounds off its corners and passes to the
limit.} to the (Levi-flat) hypersurface $\Omega \cap \{
\vert g \vert = 1\}$ with boundary equal to $\partial \Omega \cap \{
\vert g \vert = 1\}$, we get
\[
F(z)
=
\int_{\partial\Omega\cap\{\vert g\vert >1\}}\,f(\zeta)\,
{\sf BM}(\zeta,z)
+
\int_{\partial\Omega\cap\{\vert g\vert =1\}}\,
f(\zeta)\,\Theta_{g,z}(\zeta).
\]

But the holomorphic extension $F$ of an arbitrary $f \in \mathcal{
C}_{ CR}^0 \big( \partial \Omega \backslash K \big)$ is still unknown
and in fact has to be constructed\,! Since the two integrations in the
above formula are performed on parts of $\partial \Omega \backslash K$
where $f$ is defined, we are led to set:
\[
F_g(z)
=
\int_{\partial\Omega\cap\{\vert g\vert >1\}}\, 
f(\zeta)\,{\sf BM}(\zeta,z)
+
\int_{\partial\Omega\cap\{\vert g\vert =1\}} 
f(\zeta)\,\Theta_{g,z}(\zeta),
\]
as a candidate extension of $f$ at every $z \in \Omega_g$. Since $K$
is $\mathcal{ O}( \overline{ \Omega })$-convex, $\Omega$ is the union
of the regions $\Omega_g$ for $g$ running in $\mathcal{ O} (
\overline{ \Omega} )$, but these extensions $F_g(z)$ do depend on $g$,
because $\Theta_{g, z}$ does. The remainder of the proof (\cite{
stu1993, po1997}) then consists in:

\begin{itemize}

\smallskip\item[\bf (a)]
verifying that $F_g$ is holomorphic (the kernels are
not holomorphic with respect to $z$);

\smallskip\item[\bf (a)] 
showing that two differente candidates $F_{ g_1}$ and $F_{ g_2}$
coincide in fact on $\Omega_{ g_1} \cap \Omega_{ g_2}$;
 
\smallskip\item[\bf (b)] 
verifying that at least one candidate $F_g$ has boundary value equal
to $f$ on some controlled piece of $\partial \Omega \backslash K$.

\end{itemize}\smallskip

The reader is referred to \cite{ stu1981, lu1986, stu1993}
for complete arguments.

\smallskip

In the above construction, the strict pseudoconvexity of $\Omega$
insured the existence of a Stein ({\it i.e.} pseudoconvex)
neighborhood basis $\big( \mathcal{ M}_j \big)_{ j\in J}$ of
$\overline{ \Omega}$ which guaranteed in turn the existence of a Hefer
decomposition. It was pointed out by Ortega that Hefer decomposition
(called {\sl Gleason decomposition} in~\cite{ or1987}) holds on
$\mathcal{ C}^\infty$ pseudoconvex boundaries $\partial \Omega \Subset
\C^n$, but may fail in the nonpseudoconvex realm. So, let $\Omega
\Subset \C^n$ be a bounded domain having $\mathcal{ C}^\infty$
boundary. Denote by $A^\infty (\Omega) := \mathcal{ O} (\Omega) \cap
\mathcal{ C}^\infty (\overline{ \Omega})$ the ring of holomorphic
functions in $\Omega$ that are $\mathcal{ C}^\infty$ up to the
boundary.

\def\thetheorem{3.8}\begin{theorem}
{\rm (\cite{ or1987})} If $\partial \Omega \Subset \C^n$ {\rm
(}$n\geqslant 1${\rm )} is $\mathcal{ C}^\infty$ and pseudoconvex,
then every $g\in A^\infty ( \Omega )$ has a decomposition
\[
g(z)
-
g(w)
=
\sum_{k=1}^n\,g_k(z,w)
[z_k-w_k],
\]
with the $g_k \in A^\infty (\Omega \times \Omega)$.
\end{theorem}

This decomposition formula also holds under the assumption that
$\Omega$ is a domain of holomorphy (having possibly nonsmooth
boundary), but provided that $\overline{ \Omega}$ has a basis of
neighborhoods consisting of Stein domains. However, not every
$\mathcal{ C}^\infty$ weakly pseudoconvex boundary admits a Stein
neighborhood basis, as is shown by the so-called {\sl worm domains}
(\cite{ df1977, fs1987}).

\def\theexample{3.9}\begin{example}{\rm 
Furthermore, the above decomposition theorem fails to hold on general
domains. Following~\cite{ or1987}, consider the union $\Omega_1 \cup
\Omega_2$ in $\C^2$ of the two sets
\[
\aligned
\Omega_1
&
:=
\{
-4<x_1<0,\
\vert z_2\vert<e^{x_1}
\}
\ \ \ \ \ 
\text{\rm and}
\\
\Omega_2
&
:=
\{
0\leqslant x_1<4,\
e^{-1/x_1}<\vert z_2\vert<1
\}.
\endaligned
\]
The continuity principle along families of analytic discs parallel to
the $z_2$-axis shows that the envelope of holomorphy of $\Omega_1
\cup \Omega_2$ contains $\Omega_1 \cup \Omega_3$,
where $\Omega_3 := \{ 0 \leqslant x_1 < 4, \ \vert z_2\vert < 1\}$.

The holomorphic mapping $R (z_1, z_2) := (e^{ i\, z_1}, z_2)$ is
one-to-one from $\Omega_1 \cup \Omega_2$ onto its image $R (\Omega_1
\cup \Omega_2)$. However, the extension of $R$ to $\Omega_1 \cup
\Omega_3$ is not injective, because $R$ takes the same value at the
two points $(\pm \pi, e^{ -2\,\pi} ) \in \Omega_1 \cup \Omega_3$. If
Theorem~3.8 were true on the domain $R (\Omega_1 \cup \Omega_2)$,
pulling the decomposition formula back to $\Omega_1 \cup \Omega_2$, it
would follow that every $g\in \mathcal{ O} (\Omega_1 \cup \Omega_2)$
has a decomposition
\[
\aligned
g(z)
-
g(w)
&
=
\widetilde{g}_1(e^{i\,z},w)
\big[
e^{i\,z_1}
-
e^{i\,z_2}
\big]
+
\widetilde{g}_2(e^{i\,z},w)
[z_2
-
w_2]
\\
&
=
g_1(z,w)
\big[
e^{i\,z_1}
-
e^{i\,z_2}
\big]
+
g_2(z,w)
[z_2
-
w_2].
\endaligned
\]
Then the same decomposition would hold for every $g\in \mathcal{ O}
(\Omega_1 \cup \Omega_3)$, by automatic holomorphic extension of $g,
g_1, g_2$. Choosing $z = (- \pi, e^{ - 2\, \pi})$, $w = (\pi, e^{
-2\, \pi})$ and $g$ such that $g(z) \neq g(w)$ ($g := z_1$ will do\,!), 
we reach a contradiction.
}\end{example}

\def\thecorollary{3.10}\begin{corollary}
{\rm (\cite{ or1987, lp2003})}
Every function holomorphic in a domain $\Omega \subset \C^n$ enjoys
the Hefer division property precisely when the envelope of holomorphy
of $\Omega$ is schlicht.
\end{corollary}

The above results mean that a direct application of the integral
formula approach sketched above becomes impossible for domains having
nonschlicht envelope. Nevertheless, in~\cite{ lt1988}, using more
general divison methods (\cite{hele1984}), a Bochner-Martinelli kernel
on an arbitrary Stein manifold was constructed that enabled to obtain
Theorem~3.28 below, valid for nonpseudoconvex domains.

\smallskip

We conclude our presentation of Theorem~3.7 by exposing the ``only
if'' of Theorem~3.7.

\def\thelemma{3.11}\begin{lemma}
{\rm (\cite{ rs1989, stu1993})} Let $\partial
\Omega \Subset \C^2$ be a $\mathcal{ C}^2$
strongly pseudoconvex boundary and let $K \subset \partial \Omega$
be a compact set. Then $\Omega \big \backslash \widehat{ K}_{
\mathcal{ O} (\overline{ \Omega})}$ is pseudoconvex.
\end{lemma}

Taking for granted the lemma, by contraposition, suppose that $K
\subset \partial \Omega$ is not $\mathcal{ O} (\overline{
\Omega})$-convex, viz. $K \subsetneqq \widehat{ K}_{ \mathcal{ O}
(\overline{ \Omega})}$ and show that $K$ is not removable. It follows
from strict pseudoconvexity of $\partial \Omega$ that $\Omega \cap
\widehat{ K}_{ \mathcal{ O} (\overline{ \Omega})}$ is nonempty (\cite{
stu1993}). Leaving $K$ fixed, by deforming
$\partial \Omega$ away from $\Omega$, we may enlarge slightly $\Omega$
as a domain $\Omega' \supset \Omega$ with $\partial \Omega ' \supset
K$ and $\Omega ' \supset \partial \Omega \backslash K$, having
$\mathcal{ C}^2$ boundary $\partial \Omega '$ close to $\partial
\Omega$, as illustrated. Since by the lemma, $\Omega \big \backslash
\widehat{ K}_{ \mathcal{ O} ( \overline{ \Omega})}$ is pseudoconvex,
it follows easily (\cite{ stu1993}) that $\Omega ' \big \backslash
\widehat{ K}_{ \mathcal{ O} ( \overline{ \Omega})}$ is also
pseudoconvex. Consequently (\cite{ ho1973}), 
there exists a holomorphic
function $F' \in \mathcal{ O} \big( \Omega' \big \backslash \widehat{
K}_{ \mathcal{ O} ( \overline{ \Omega })} \big)$ that does not extend
holomorphically at any point of the boundary of $\Omega' \big
\backslash \widehat{ K}_{ \mathcal{ O} ( \overline{ \Omega })}$. The
restriction of $F'$ to $\partial \Omega \backslash K$ is a CR function
on $\partial \Omega \backslash K$ for which $K$ is {\it not}\,
removable, since $\Omega \cap \widehat{ K}_{ \mathcal{ O} (\overline{
\Omega })} \neq \emptyset$.

\subsection*{ 3.12.~Removability, polynomial hulls and Cantor
sets} A generalization of Theorem~3.7, essentially with the same proof
(excepting notational complications) holds in arbitrary complex
dimension $n\geqslant 2$.

\def\thetheorem{3.13}\begin{theorem}
{\rm (\cite{ lu1986, stu1993})} Let $\Omega \Subset \C^n$, $n\geqslant
2$, be a bounded pseudoconvex domain such that $\overline{ \Omega}$
has a Stein neighborhood basis. If $K \subset \partial \Omega$ is
compact and $\mathcal{ O} (\overline{ \Omega})$-convex, and if
$\partial \Omega = K \cup M$, where $M$ is a connected $\mathcal{
C}^1$ hypersurface of $\C^n \backslash K$, then $K$ is CR-removable.
\end{theorem}

\def\theexample{3.14}\begin{example}{\rm
{\rm (\cite{ cst1994, jo1999a})} Let $M$ be a connected compact
orientable $(2n-3)$-dimensional maximally complex (Definition~4.7
below) CR manifold of class $\mathcal{ C}^1$ contained in the unit
sphere $\partial \B_n$ ($n \geqslant 2$) with empty boundary in the
sense of currents. Such an $M$ is called a {\sl maximally complex
cycle}. By a theorem due to Harvey-Lawson (reviewed as Theorem~4.16
below), 
if $M$ satisfies the moments' condition, then
$M$ is the boundary of a unique complex $(n-1)$-dimensional
complex subvariety $\Sigma
\subset \B_n$. Since the cohomology group $H^2 (\B_n , \Z)$ vanishes,
by a standard Cousin problem, $\Sigma$ may be defined as the zero-set
of some global holomorphic function $f \in \mathcal{ O} (\B_n) \cap
\mathcal{ C}^0 (\overline{ \B}_n)$. The maximum principle yields that
the compact set $K := \Sigma \cup M = \overline{ \Sigma}$ is
$\mathcal{ O} (\overline{ \B}_n)$-convex. Consequently, the envelope
of holomorphy of an arbitrarily thin one-sided neighborhood of
$\partial \B_n \backslash M$ is equal to the pseudoconvex domain $\B_n
\backslash \big( M \cup \Sigma \big)$.

}\end{example}

If in addition $\Omega$ is Runge (\cite{ ho1973}) or if $\overline{
\Omega}$ is polynomially convex, then every $f \in \mathcal{ O}
(\overline{ \Omega})$ may be approximated uniformly by polynomials on
some sufficiently small neighborhood of $\overline{ \Omega}$ (whose
size depends on $f$). It then follows that polynomial convexity and
$\mathcal{ O} (\overline{ \Omega })$-convexity are equivalent. As a
paradigmatic example, this holds when $\Omega = \B_n$ is the unit
ball.

\def\thecorollary{3.15}\begin{corollary}
{\rm (\cite{ stu1993})} Let $\Omega$ and $K \subset \Omega$ be as in
Theorem~3.13 and assume that $\Omega$ is Runge in $\C^n$, for instance
$\Omega = \B_n$. If $K$ is polynomially convex, then $K$ is
CR-removable. If $n=2$, the CR-removability of $K$ is equivalent to
its polynomial convexity.
\end{corollary}

Although the last necessary and sufficient condition seems to be
satisfactory, we must point out that concrete geometric
characterizations of polynomial convexity usually are hard to provide.
In \S5.14 below, we shall describe a class of removable compact
sets whose polynomial convexity may be established directly.

\smallskip

A compact subset $K$ of $\R^n$ ($n \geqslant 1$) is a {\sl Cantor
set}\, if it is perfect, viz. coincides with its first derived set
$K'$. It is called {\sl tame}\, if there is a homeomorphism of $\R^n$
onto itself that carries $K$ onto the standard middle-third Cantor set
contained in the coordinate line $\R_{ x_1}$.

Tame Cantor sets $K$ in a $\mathcal{ C}^2$ strongly pseudoconvex
boundary $\partial \Omega \Subset \C^2$ were shown to be CR-removable
in~\cite{ fs1991}, provided there exists a Stein neighborhood
$\mathcal{ D}$ of $K$ in $\C^2$ such that $K$ is $\mathcal{ O}
(\overline{ \mathcal{ D } })$-convex. By further analysis, this last
assumption was shown later to be redundant and in general, tame Cantor
sets are CR-removable. It was then suggested in~\cite{ stu1993} that
all Cantor subsets of $\partial \B_n$ ($n \geqslant 2$) are removable,
or equivalently polynomially convex. Nevertheless, Rudin and then
Vitushkin, Henkin and others had constructed Cantor sets $K \subset
\C^2$ having large polynomial hull $\widehat{ K}$, {\it e.g.} so that
$\widehat{ K}$ contains a complex curve, or even contains interior
points. Recently, in a beautiful paper, J\"oricke showed how to put
such sets in the $3$-sphere $\partial \B_2$, thus solving the question
in the negative.

\def\thetheorem{3.16}\begin{theorem}
{\rm (\cite{ jo2005})} For every positive number $r < 1$, there exists
a Cantor set $K \subset \partial \B_2$ whose polynomial hull
$\widehat{ K}$ contains the closed ball $r\, \overline{ \B}_2$.
\end{theorem}

\subsection*{ 3.17.~$L^{\sf p}$-removability
and further results} In the definition of CR-removability, nothing is
assumed about the behavior from $\overline{ \Omega} \big \backslash K$
up to $K$: the rate of growth may be arbitrarily high. If,
differently, functions are assumed to be tame on $\partial \Omega$
(including $K$), better removability assertions hold.

\def\thedefinition{3.18}\begin{definition}{\rm
A compact subset $K$ of a $\mathcal{ C}^1$ boundary $\partial \Omega
\Subset \C^n$ ($n\geqslant 2$) is called {\sl $L^{\sf p}$-removable} ($1
\leqslant {\sf p} \leqslant \infty$) if every function $f\in L^{
\sf p} (\partial \Omega)$ which is CR on $\partial \Omega \backslash
K$ is in fact CR on the whole boundary $\partial \Omega$.
}\end{definition}

Then by the Hartogs-Bochner theorem, $f$ admits a holomorphic
extension to $\Omega$ that may be checked to belong to $H^{ \sf p}
(\Omega)$.

\def\thetheorem{3.19}\begin{theorem}
{\rm (\cite{ ac1994})} Let $\Omega \Subset \C^n$ {\rm (}$n \geqslant
2${\rm )} be a bounded domain having $\mathcal{ C}^2$ boundary
$\partial \Omega$ and let $M$ be a $\mathcal{ C}^2$ totally real
embedded submanifold of $\partial \Omega$. If $K \subset M$ is a
polynomially convex compact subset, then $K$ is $L^{\sf p}$-removable.
\end{theorem}

\smallskip

In complex dimension $n \geqslant 3$, the two extension Theorems~3.13
and~3.19 are not optimal. In general, additional extension phenomena
occur, which are principally overlooked by assumptions on the hull of
the singularity. A more geometric point of view (\S3.23
below) shows that these theorems may be established by means
of holomorphic extension along one-parameter families of complex
analytic hypersurfaces, whereas the (finer) {\sl
Kontinuit\"atssatz} holds along families of analytic {\sl discs},
whose thinness offers more freedom to fill in maximal domains of
extension.

\def\theexample{3.20}\begin{example}{\rm
Let $\Omega := \B_3$ be the unit ball in $\C^3$, and let 
\[
K=
\big\{
(z_1,z_2,0)\in\partial\B_3:\,
\vert z_1\vert \geqslant 1/2
\big\}
\]
be a 3-dimensional ring in the intersection of $\partial \B_3$ with
the $(z_1, z_2)$-plane. The maximum principle along discs parallel to
the $z_2$-axis yields:
\[
\widehat{K}_{\mathcal{ O}
(\overline{\B}_3)}
=
\big\{
(z_1,z_2,0)\in\overline{\B}_3:
\vert z_1\vert \geqslant 1/2
\big\}
\neq
K,
\]
so $K$ is not $\mathcal{ O} ( \overline{ \B}_3 )$-convex.
Nevertheless, this $K$ is removable. Indeed, applying the continuity
principle, we may first fill in $\B_3 \big \backslash \widehat{ K}$ by
means of discs parallel to the $z_2$-axis and then fill in the
complete ball $\B_3$, by means of discs parallel to the $z_3$-axis.
}\end{example}

In higher dimensions $n \geqslant 3$, the relevant characterizations of
CR-removable compact sets contained in strongly pseudoconvex frontiers
are of cohomological nature (\S3.33 below). In another vein, the
assumption that $\overline{ \Omega }$ possesses a Stein neighborhood
basis in Theorem~3.13 above inspired some authors to generalize
Stout's theorem as follows.

\def\thedefinition{3.21}\begin{definition}{\rm
Let $\Omega$ be a relatively compact domain of a Stein
manifold $\mathcal{ M }$ and let $K \subset \overline{ \Omega}$ be a
compact set. The {\sl $\mathcal{ O} ( \mathcal{ M })$-convex hull} of
$K$ is
\[
\widehat{K}_{\mathcal{O}(\mathcal{M})}
:=
\Big\{
z\in\mathcal{M}:\
\vert g(z)\vert\leqslant\max_{w\in K}\,
\vert g(w)\vert\
\text{\rm for all}\
g\in\mathcal{O}(\mathcal{M})
\Big\}.
\]
If $K = \widehat{ K}_{ \mathcal{ O} (\mathcal{ M })}$, then $K$
is called {\sl $\mathcal{ O} (\mathcal{ M })$-convex}.
}\end{definition}

In $\C^n$, the $\mathcal{ O} (\mathcal{ M })$-convex hull coincides
with the polynomial hull. Notice that the next theorem is valid
without pseudoconvexity assumption on $\Omega$.

\def\thetheorem{3.22}\begin{theorem}
{\rm (\cite{ stu1981, lt1988, ky1991, stu1993, jo1995})} Let
$\mathcal{ M}$ be a Stein manifold of dimension $n\geqslant 2$, let
$\Omega \Subset \mathcal{ M}$ be a relatively compact domain such that
$\mathcal{ M} \backslash \overline{ \Omega}$ is connected and let $K
\subset \overline{ \Omega}$ be a compact set with $K = \widehat{ K}_{
\mathcal{ O} (\mathcal{ M })} \cap \partial \Omega$. Then every CR
function $f$ defined on $\partial \Omega \backslash K$ extends
holomorphically to $\Omega \backslash \widehat{ K}_{ \mathcal{ O}
(\mathcal{ M })}$, {\it i.e.}{\rm :}

\begin{itemize}

\smallskip\item[$\bullet$]
if $\partial \Omega \backslash K$ is a $\mathcal{ C}^{ \kappa,
\alpha}$ hypersurface, with $\kappa \geqslant 1$ and $0\leqslant
\alpha \leqslant 1$, and if $f \in \mathcal{ C}_{ CR}^{ \kappa,
\alpha} (\partial \Omega \backslash K)$, then the holomorphic
extension $F \in \mathcal{ O} ( \Omega \backslash K)$ belongs to the
class $\mathcal{ C}^{ \kappa, \alpha} \big( \overline{ \Omega}
\backslash \widehat{ K}_{ \mathcal{ O} (\mathcal{ M })} \big)${\rm ;}

\smallskip\item[$\bullet$]
if $\partial \Omega \backslash K$ is a $\mathcal{ C}^1$ hypersurface
and if $f \in L_{ loc}^{\sf p} (\partial \Omega \backslash K)$ with
$1\leqslant {\sf p}\leqslant \infty$, then at every point $p\in
\partial \Omega \backslash K$, the holomorphic extension $F \in
\mathcal{ O}\big( \overline{ \Omega} \backslash \widehat{ K}_{
\mathcal{ O} (\mathcal{ M })} \big)$ belongs to the Hardy space $H_{
loc}^{ \sf p} ( U_p \cap \Omega)$, for some small neighborhood $U_p$
of $p$ in $\mathcal{ M}$.

\end{itemize}\smallskip

\end{theorem}

\subsection*{ 3.23.~$A (\Omega)$-hull and 
removal of singularities on pseudoconvex boundaries} Following~\cite{
jo1995, po1997, jp2002}, we now expose a geometric aspect of some of
the preceding removability theorems. Let $\Omega \Subset \C^n$ with
$n\geqslant 2$ be a bounded domain having frontier of class at least
$\mathcal{ C}^1$. By $A (\Omega) =\mathcal{ O} (\Omega) \cap \mathcal{
C}^0 ( \overline{ \Omega})$, we denote the ring of holomorphic
functions in $\Omega$ that are continuous up to the boundary. Let $K
\subset \overline{ \Omega}$ be a compact set.

\def\thedefinition{3.24}\begin{definition}{\rm
The $A (\Omega )$-hull of $K$ is
\[
\widehat{K}_{A(\Omega)}
:=
\Big\{
z\in\overline{\Omega}:\
\vert g(z)\vert\leqslant\max_{w\in\overline{\Omega}}\,
\vert g(w)\vert\
\text{\rm for all}\
g\in\mathcal{A}(\overline{\Omega})
\Big\}.
\]
If $K = \widehat{ K}_{ A (\Omega)}$, then $K$ is called $A ( \Omega
)$-{\sl convex}. If $K = \partial \Omega \cap \widehat{ K}_{ A (
\Omega )}$, then $K$ is called {\sl CR-convex}.
}\end{definition}

The next theorem is stronger than Theorem~3.13 in two aspects: 

\begin{itemize}

\smallskip\item[$\bullet$]
the inclusion $\widehat{ K}_{ A (\Omega)} \subset \widehat{ K}_{
\mathcal{ O} (\overline{ \Omega})}$ holds in general and may be
strict;

\smallskip\item[$\bullet$]
it is not assumed that the pseudoconvex domain $\Omega$ has a Stein
neighborhood basis.

\end{itemize}\smallskip

\def\thetheorem{3.25}\begin{theorem}
{\rm (\cite{ jo1995})}
Let $\Omega$ be a bounded weakly pseudoconvex domain in
$\C^2$ having frontier of class $\mathcal{ C }^2$ and let $K$ be a
compact subset of $\partial \Omega$ with $K \neq \partial \Omega$ such
that $K$ is CR-convex, namely $K = \partial \Omega \cap \widehat{ K}_{
A( \Omega )}$. Then the following are true.

\begin{itemize}

\smallskip\item[{\bf 1)}]
Let $\mathcal{ V} (\partial \Omega \backslash K)$ be an interior
one-sided neighborhood of $\partial \Omega \backslash K$ with the
property that each connected component of $\mathcal{ V} (\partial
\Omega \backslash K)$ contains in its boundary exactly one component
of $\partial \Omega \backslash K$ and no other point of $\partial
\Omega \backslash K$. Then for every holomorphic function $f \in
\mathcal{ O } \big(\mathcal{ V} (\partial \Omega \backslash K ) \big)$,
there exists a holomorphic function $F \in \mathcal{ O} ( \Omega
\backslash \widehat{ K}_{ A (\Omega)})$ with $F = f$ in $\mathcal{ V}
(\partial \Omega \backslash K)$.

\smallskip\item[{\bf 2)}] 
{\rm (\cite{ as1990})}
There is a one-to-one correspondence between connected components of
$\partial \Omega \backslash K$ and connected components of $\Omega
\backslash \widehat{ K}_{ A (\Omega )}$, namely the boundary of each
component of $\Omega \backslash \widehat{ K}_{ A( \Omega )}$ contains
exactly one connected component of $\partial \Omega \backslash K$ and
does not intersect any other component.

\smallskip\item[{\bf 3)}]
If the boundary $\partial \Omega$ is of class $\mathcal{ C}^\infty$,
then $\Omega \backslash \widehat{ K}_{ A (\Omega)}$ is pseudoconvex,
hence it is the envelope of holomorphy of $\mathcal{ V} ( \partial
\Omega \backslash K)$.

\end{itemize}\smallskip

\end{theorem}

If $K$ is not CR-convex, the one-to-one correspondence between the
connected components of $\partial \Omega \backslash K$ and
those of $\partial \Omega \backslash \widehat{ K}_{ A (\Omega)}$
may fail.

\def\theexample{3.26}\begin{example}{\rm 
Indeed, let $\Omega := \B_2 \cap \big\{ x_1 < \frac{ 1}{ 2} \big\}$ be
a truncation of the unit ball and let $K : = \partial \B_2 \cap \big\{
x_1 = \frac{ 1}{ 2} \big\}$ be the intersection of the three-sphere
$\partial \B_2$ with the real hyperplane $\big\{ x_1 = \frac{ 1}{ 2}
\big\}$ ({\it see} only the left hand side of the diagram).

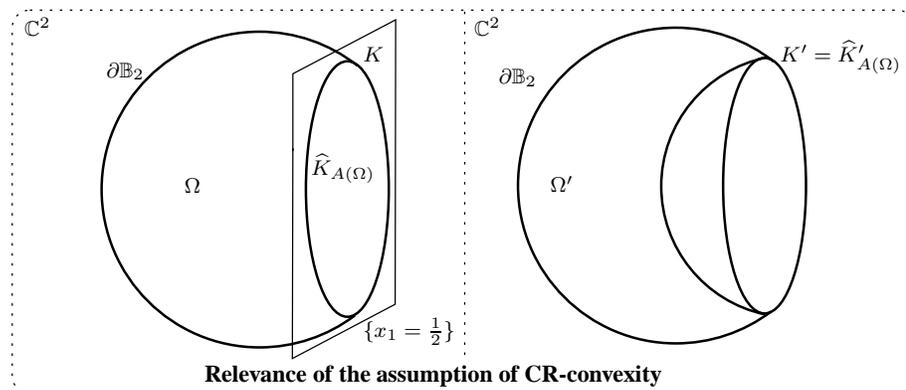
\begin{center}
\input 2-ball-truncated.pstex_t
\end{center}

The Levi-flat $3$-ball $\B_2 \cap \big\{ x_1 = \frac{ 1}{ 2} \big\}$
being foliated by complex discs, the maximum principle entails that
$\widehat{ K}_{ A (\Omega)} = \overline{ \B}_2 \cap \big\{ x_1 =
\frac{ 1}{ 2} \big\} = \widehat{ K}_{ A (\Omega)} \cap \partial \Omega
\neq K$, hence $K$ is not CR-convex. Also, $\partial \Omega \backslash
K$ has two connected components $\partial \B_2 \cap \big\{ x_1 <
\frac{ 1}{ 2} \big\}$ and $\B_2 \cap \big\{ x_1 = \frac{ 1}{ 2}
\big\}$, whereas $\partial \Omega \backslash \widehat{ K}_{ A
(\Omega)} = \partial \B_2 \cap \big\{ x_1 < \frac{ 1}{ 2} \big\}$ is
connected. Any function on $\partial \Omega \backslash K$ equal to two
distinct constants on the two connected components of $\partial \Omega
\backslash K$ is CR and not holomorphically extendable to $\Omega =
\Omega \backslash \widehat{ K}_{ A (\Omega)}$. Finally, by smoothing
out $\partial \Omega$ near the two-sphere $\partial \B_2 \cap \big\{
x_1 = \frac{ 1}{ 2} \big\}$, we obtain an example with $\mathcal{
C}^\infty$ boundary. 
}\end{example}

\subsection*{ 3.27.~Hulls and holomorphic extension from nonpseudoconvex
boundaries} Since the work~\cite{ lu1986} of Lupacciolu, the extension
of Theorem~3.7 to nonpseudoconvex boundaries was a daring open problem
(\cite{ stu1993}).

\def\thetheorem{3.28}\begin{theorem}
{\rm (\cite{ po1997, jp2002, lp2003})} Let $\Omega$ be a {\rm not
necessarily pseudoconvex} bounded domain in $\C^n$ {\rm (}$n\geqslant
2${\rm )} having connected $\mathcal{ C}^2$ frontier and let $K
\subset \partial \Omega$ be a compact set with $\partial \Omega
\backslash K$ connected such that $K = \partial \Omega \cap \widehat{
K}_{ A (\Omega)}$. Then for every continuous CR function $f\in
\mathcal{ C}_{ CR}^0 (\partial \Omega \backslash K)$, there exists a
holomorphic function $F \in \mathcal{ O} \big( \Omega \backslash
\widehat{ K}_{ A (\Omega)} \big) \cap \mathcal{ C}^0 \big( [\Omega
\backslash \widehat{ K}_{ A (\Omega )}] \cup [\partial \Omega
\backslash K] \big)$ such that $F \vert_{ \partial \Omega \backslash
K} = f$.
\end{theorem}

A purely geometrical proof applying a global continuity principle
together with a fine control of monodromy may be found in~\cite{
po1997, jp2002}; {\it cf.} also~\cite{ mp2007}. By a topological
device, a second proof (\cite{ lp2003}) derives the theorem from the
following statement, established by means of $\overline{ \partial}$
techniques.

\def\thetheorem{3.29}\begin{theorem}
{\rm (\cite{ lt1988})} Let $\mathcal{ M}$ be a Stein manifold of
complex dimension $n\geqslant 2$, let $K \subset \mathcal{ M}$ be a
compact set that is $\mathcal{ O} (\mathcal{ M})$-convex and let
$\Omega \subset \mathcal{ M}$ be a relatively compact {\rm not
necessarily pseudoconvex} domain such that $\partial \Omega \backslash
K$ is a connected $\mathcal{ C}^1$ hypersurface of $\mathcal{ M}
\backslash K$. Then for every continuous CR function $f$ on $\partial
\Omega \backslash K$, there exists a holomorphic function $F \in
\mathcal{ O} (\Omega \backslash K) \cap \mathcal{ C}^0 ( \overline{
\Omega} \backslash K)$ with $F \vert_{ \partial \Omega \backslash K} =
f$.
\end{theorem}

Contrary to the case where $\partial \Omega$ is pseudoconvex (as in
Theorem~3.25), even if $K$ is CR-convex, the one-to-one correspondence
between the connected components of $\partial \Omega \backslash K$ and
those of $\Omega \backslash \widehat{ K}_{ A (\Omega)}$ may fail to
hold. For this reason, $\partial \Omega \backslash K$ is assumed to be
connected in Theorem~3.29.

\def\theexample{3.30}\begin{example}{\rm (\cite{ lp2003})
We modify Example~3.26 so as to get a nonpseudoconvex boundary as
follows ({\it see} the right hand side of the diagram above). Let
$\Omega'$ be the unit ball $\B_2$ from which we substract the
closed ball
$\overline{ B 
( q, 1 )}$ of radius $1$ 
centered at the point $q$ of coordinates $(1,0)$. 
A computation with defining (in)equations shows
that $\Omega '$ is contained in $\big\{ x_1 < \frac{ 1}{ 2}
\big\}$. Notice that $\Omega '$ is not pseudoconvex and in fact, its
envelope of holomorphy is single-sheeted and equal to the domain
$\Omega = \B_2 \cap \big\{ x_1 < \frac{ 1}{ 2} \big\}$ 
drawn in the left hand side. 
Let $K' := \B_2 \cap \big\{ x_1 = \frac{ 1}{ 2}\big\}
\subset \partial \Omega'$ (this set is the same 
$2$-sphere as the set $K$ of the preceding example). 
Then $K'$ is CR-convex, since the candidate for its $A
(\Omega')$-hull is the three-sphere $\overline{ \B_2} \cap \big\{ x_1
= \frac{ 1}{ 2} \big\}$ that lies outside $\Omega'$. However,
$\partial \Omega' \backslash K'$ has two connected components, namely
$\partial \B_2 \cap \big\{ x_1 < \frac{ 1}{ 2}\big\}$ and $\partial B(
q, 1) \cap \big\{ x_1 < \frac{ 1}{ 2}\big\}$,
whereas $\Omega' \backslash \widehat{ K}_{ A (\Omega)}' = 
\Omega ' \backslash K' = \Omega'$ is
connected. Hence any CR function equal to two distinct constants on
these two components fails to extend holomorphically to $\Omega'$.
Finally, by smoothing out $\partial \Omega$ near the two-sphere
$\partial \B_2 \cap \big\{ x_1 = \frac{ 1}{ 2}\big\}$, we obtain an
example with $\mathcal{ C }^\infty$ boundary.
}\end{example}

If we drop CR-convexity of $K$, viz. if $K \neq \widehat{ K}_{ A
(\Omega)} \cap \partial \Omega$, then monodromy problems come on
scene: the natural embedding of $\Omega \backslash \widehat{ K}_{ A
(\Omega)}$ into the envelope of holomorphy of a one-sided neighborhood
of $\partial \Omega \backslash K$ may fail to be one-to-one.

\def\theexample{3.31}\begin{example}{\rm 
(\cite{ lp2003}) Consider the real four-dimensional open cube $C :=
(-1, 1) \times i\, (-1, 1) \times (-1, 1) \times i (-1, 1)$ in $\C^2
\simeq \R^4$.

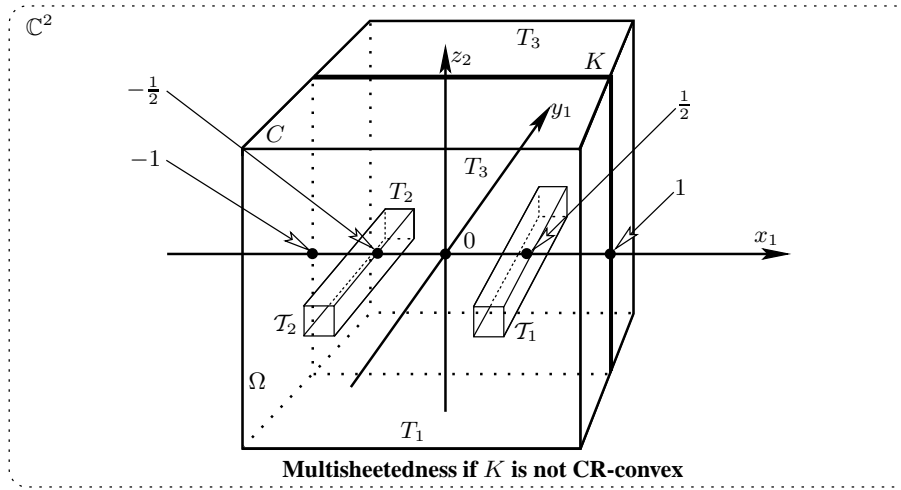
\begin{center}
\input widehatK_AD.pstex_t
\end{center}

Choose $\varepsilon >0$ small and remove from this cube $C$ firstly
the narrow tunnel $\mathcal{ T}_1 := \{ \vert z_2 \vert \leqslant
\varepsilon, \, \vert x_1 - 1/2 \vert \leqslant \varepsilon \}$ having
an entrance and an exit and secondly the (incomplete) narrow tunnel
$\mathcal{ T}_2 := \{ \vert z_2 \vert \leqslant \varepsilon, \, \vert
x_1 + 1/2 \vert \leqslant \varepsilon, \, - 1 < y_1 \leqslant 1/2 \}$
having only an entrance, and call $\Omega$ the obtained domain. Let $K
:= \partial C \cap \{ y_1 = 0\}$. The complete tunnel insures that
$\partial \Omega \backslash K$ is connected. Moreover, the maximum
principle along families of analytic discs parallel to the complex
$z_2$-axis enables to verify that
\[
\widehat{K}_{A(\Omega)} 
=
\big(
\Omega\cap\{y_1=0\}
\big)
\bigcup
K
\bigcup
\big(
\partial\mathcal{T}_1\cap\{y_1=0\}
\big)
\bigcup
\big(
\partial\mathcal{T}_2\cap\{y_1=0\}
\big).
\]
It follows that $\partial \Omega \backslash \widehat{ K}_{ A
(\Omega)}$ has three connected components, firstly the part $T_1$ of
$\partial \Omega$ that lies in the half-space $\{ y_1<0\}$; secondly
the dead-lock part $T_2$ of the second tunnel that lies in $\{ y_1
>0\}$; and thirdly, the remainder $T_3$ of the boundary, that lies in
$\{ y_1 > 0\}$. 

The branch of ${\rm log} \, z_1$ satisfying ${\rm log}\, 1 = 0$ is
uniquely defined in $\C^2 \backslash \{ (x_1, z_2) : \, x_1 \leqslant
0\}$, hence ${\rm log}\, z_1$ is holomorphic in a neighborhood of
$\partial \Omega \backslash \overline{ T}_2$, where $\overline{ T}_2
:= \partial \mathcal{ T}_2 \cap \{ y_1 \geqslant 0\}$. In addition,
${\rm log}\, z_1$ extends from points near $\mathcal{ T}_2$ in $\{ y_1
< 0\}$ to a neighborhood of $\overline{ T}_2$. In sum, it defines a
single-valued function that is holomorphic in a neighborhood of
$\partial \Omega$.

Observe that $\big( - \frac{ 1}{ 2} + \frac{ i}{ 2}, 0 \big) \in T_2
\subset \partial \Omega$. The value of ${\rm log}\, z_1$ thus defined
at this point is ${\rm log} \big( \frac{ 1}{ \sqrt{ 2}} \, e^{ -i
5\pi/4}\big) = {\rm log}\, \frac{ 1}{ \sqrt{ 2}} - i\, \frac{ 5\pi}{
4}$. On the other hand, ${\rm log}\, z_1$ restricted to a
neighborhood of $\partial C \cap \{ y_1 >0\} \subset \partial \Omega$
extends holomorphically to $C \cap \{ y_1 >0 \}$ (by means of unit
discs parallel to the $z_2$-axis) as ${\rm log}\, z_1$ itself! But
the value of this extension at $\big( -\frac{ 1}{ 2} + \frac{ i}{ 2},
0 \big)$ is different: ${\rm log} \big( \frac{ 1}{ \sqrt{ 2}} \, e^{ i
3\pi/4}\big) = {\rm log}\, \frac{ 1}{ \sqrt{ 2}} + i\, \frac{ 3\pi}{
4}$.
}\end{example}

To conclude this paragraph, before surveying the cohomological
characterizations of removable singularities in dimension $n\geqslant
3$, we reformulate the obtained characterization in complex dimension
$n = 2$. It is known that a compact set $K \subset \C^n$ is
polynomially convex if and only if the $\overline{
\partial}$-cohomology group $H_{ \overline{ \partial}}^{ 0, 1} (K)$ is
trivial and holomorphic functions in a neighborhood of $K$ can be
approximated by polynomials uniformly on $K$. Thus, we can state a
complete formulation of Theorem~3.7, with the supplementary assumption
that $\mathcal{ O} (\overline{ \Omega})$ may be approximated uniformly
by polynomials. This insures that polynomial convexity coincides with
$\mathcal{ O} (\overline{ \Omega})$-convexity. As a major example, the
theorem holds for $\Omega$ equal to the unit ball $\B_2$
(Corollary~3.15).

\def\thetheorem{3.32}\begin{theorem}
{\rm (\cite{ stu1989, stu1993, lu1994, cst1994})} The following four
conditions for a compact subset $K$ of a $\mathcal{ C}^2$ strongly
pseudoconvex compact boundary $\partial \Omega \Subset \C^2$ with
$\Omega$ Runge or $\overline{ \Omega}$ polynomially convex are
equivalent{\rm :}

\begin{itemize}

\smallskip\item[$\bullet$] 
$K$ is $\mathcal{ O}( \overline{ \Omega})$-convex.

\smallskip\item[$\bullet$]
$K$ is polynomially convex.

\smallskip\item[$\bullet$] 
$H_{ \overline{ \partial}}^{ 0, 1} (K) = 0$ and 
holomorphic functions in a neighborhood of $K$ can be
approximated by polynomials uniformly on $K$

\smallskip\item[$\bullet$] 
$K$ is removable.

\end{itemize}\smallskip

\end{theorem}

Thus, in this situation, removability amounts to polynomial convexity.
Nevertheless, the problem of characterizing geometrically the
polynomial convexity of compact sets hides several fine questions. We
shall come back to this topic in Section~5.

\subsection*{ 3.33.~Luppaciolu's characterizations}
An outstanding theorem due to Lupacciolu provides complete
cohomological characterizations of removable sets that are contained
in strongly pseudoconvex boundaries, for general $n \geqslant 2$.

Let $\mathcal{ M}$ be a Stein manifold of dimension $n \geqslant 2$ and
let $\Omega \Subset \mathcal{ M}$ be a relatively compact strongly
pseudoconvex domain having $\mathcal{ C}^2$ boundary.

Let $H_{ \overline \partial }^{ p, q} := \mathcal Z_{ \overline
\partial }^{p, q} / \overline \partial \mathcal E^{ p, q-1 }$ denote
the usual $(p, q)$-th Dolbeault cohomology group\footnote{ Appropriate
background, further survey of Lupacciolu's results and additional
material may be found in~\cite{ cst1994}.}. We endow the space
$\mathcal Z_{ \overline \partial }^{ n, n-2} (K)$ of $\overline{
\partial }$-closed $(n, n- 2)$-forms defined in a neighborhood of a
compact set $K \subset \mathcal{ M}$ with the standard locally convex
inductive limit topology derived from the inductive system of the
Fr\'echet-Schwartz spaces $\mathcal Z_{ \overline \partial }^{ n, n-2
} (U)$, as $U$ ranges through a fundamental system of open
neighborhoods of $K$ in $\mathcal{ M }$.

\def\thetheorem{3.34}\begin{theorem}
{\rm (\cite{ lu1994, cst1994})}
Assume that $\overline{ \Omega}$ is $\mathcal{ O}( \mathcal{
M})$-convex. A proper closed subset $K$ of $\partial \Omega$ is
removable {\rm if and only if} $H_{ \overline \partial }^{ n,n-1} (K)
=0$ and the restriction map $\mathcal Z_{\overline \partial}^{n, n-2}
( \mathcal{ M}) \to \mathcal Z_{ \overline \partial }^{n, n-2} (K)$
has dense image.
\end{theorem}

For $n=2$, the two conditions of the theorem reduce to the $\mathcal
O( \mathcal{ M })$-convexity of $K$ (\cite{ lu1994, cst1994}).
For $n\geqslant 3$, the following improvement is valid. By $\sp{\sigma}E$,
we denote the {\sl separated space associated to} a given topological
vector space $E$, namely the quotient $E/\overline 0$ of $E$ by the
closure of $0$.

\def\thetheorem{3.35}\begin{theorem}
{\rm (\cite{ lu1994, cst1994})} Assume that $n\geqslant 3$. Without the
assumption that $\overline{ \Omega}$ is $\mathcal O( \mathcal{
M})$-convex, the compact set $K \subset \partial \Omega$ is removable
{\rm if and only if} $H_{\overline \partial }^{ n, n-1} (K)=0$ and
$\sp{ \sigma} H_{ \overline \partial}^{ n, n-2} (K)=0$.
\end{theorem}

Lupacciolu also obtains an extrinsic characterization as follows. Let
$\Phi$ be the paracompactifying family of all closed subsets of
$\mathcal{ M} \backslash K$ that have compact closure in $\mathcal{
M}$. Let $H_{ \Phi }^{ p, q}$ the Dolbeault cohomology groups with
support in $\Phi$.

\def\thetheorem{3.36}\begin{theorem}
{\rm (\cite{ lu1994, cst1994})} For $n \geqslant 3$, a compact
subset $K$ of the boundary $\partial
\Omega$ of a $\mathcal C^2$-bounded strongly
pseudoconvex domain $\Omega \Subset \mathcal{ M}$ is removable {\rm if
and only if} $H_{ \Phi }^{ 0, 1} ( \mathcal{ M} \backslash K) = 0$.
\end{theorem}

Notice that, for $n\geqslant 3$, this theorem has the striking
consequence that the condition that $K$ be removable in a strongly
pseudoconvex boundary does not depend on the domain in question, but
rather on the situation of $K$ itself in the ambient manifold. Also,
Lupacciolu provides analogous characterizations for weak removability
(\cite{ lu1994, cst1994}).

\section*{\S4.~Smooth and metrically thin removable singularities 
\\
for CR functions}

\subsection*{ 4.1.~Three notions of removability}
We formulate the concerned notions of removability directly in
arbitrary codimension. Let $M \subset \C^n$ be a $\mathcal{ C}^{ 2,
\alpha}$ generic submanifold of positive codimension $d\geqslant 1$
and of positive CR dimension $m \geqslant 1$. Such $M$ will always be
supposed {\it connected}. In the sequel, not to mention superficial
corollaries, we will systematically assume that $M$ is globally
minimal.

\def\thedefinition{4.2}\begin{definition}{\rm
{\rm (\cite{ me1997, mp1998, jo1999a, jo1999b, mp1999, mp2002})}
A closed subset $C$ of $M$ is
said to be{\rm :}

\begin{itemize}

\smallskip\item[$\bullet$] {\sl CR-removable} if there exists a
wedgelike domain $\mathcal{ W}$ attached to $M$ to which every
continuous CR function $f\in\mathcal{ C}_{CR}^0(M\backslash C)$
extends holomorphically{\rm ;}

\smallskip\item[$\bullet$]
{\sl $\mathcal{ W }$-removable} if for every wedgelike domain
$\mathcal{ W }_1$ attached to $M\backslash C$, there is a wedgelike
domain $\mathcal{ W }_2$ attached to $M$ and a wedgelike domain
$\mathcal{ W }_3 \subset \mathcal{ W }_1 \cap \mathcal{ W}_2$ attached
to $M \backslash C$ such that for every holomorphic function $f \in
\mathcal{ O}(\mathcal{ W}_1)$, there exists a holomorphic function $F
\in \mathcal{ O}(\mathcal{ W }_2)$ which coincides with $f$ in
$\mathcal{ W }_3${\rm ;}

\smallskip\item[$\bullet$] 
{\it $L^{\sf p}$-removable}, where $1\leqslant {\sf p} \leqslant
\infty$, if every locally integrable function $f\in L^{\sf
p}_{loc}(M)$ which is CR in the distributional sense on $M \backslash
C$ is in fact CR on all of $M$.

\end{itemize}

}\end{definition}

A few comments are welcome. CR-removability requires at least $M
\backslash C$ to be globally minimal, in order that the main
Theorem~4.12(V) applies, yielding a wedgelike domain $\mathcal{ W}_1$
attached to $M\backslash C$. Then $\mathcal{ W}$-removability of $C$
implies its CR-removability. In both CR- and $\mathcal{
W}$-removabililty, after the removal of $C$, nothing is demanded about
the growth of the holomorphic extension to a global wedgelike domain
$\mathcal{ W}_2$ attached to $M$. Such extensions might well have
essential singularities at some points of $C$, although they are
holomorphic in $\mathcal{ W}_2$. On the contrary, for $L^{ \sf
p}$-removability of $C$, CR functions on $M \backslash C$ should
really extend to be CR through $C$.

Notwithstanding this difference, the sequel will reveal that $L^{ \sf
p}$-removability is also a consequence of $\mathcal{ W}$-removability,
thanks to some Hardy-space control of the holomorphic extension
$F \in \mathcal{ O} (\mathcal{ W}_2)$. In
fact, functions are assumed to be $L_{ loc}^{ \sf p}$ (a variant is to
assume continuity on $M$ instead of integrability) even near points of
$C$. This strong assumption enables to get a control of the growth of
the wedge extension. Before providing more explanations, we assert in
advance that $\mathcal{ W}$-removability is the most general notion of
removability, focusing the question on envelopes of
holomorphy. 

In codimension $d=1$, wedgelike domains identify to one-sided
neighborhoods. Then $\mathcal{ W}$-removability of $C$ means that the
envelope of holomorphy of every (arbitrarily thin) one-sided
neighborhood of $M\backslash C$ contains a complete one-sided
neighborhood of the hypersurface $M$ in $\C^n$. If $M = \partial
\Omega$ is the boundary of a bounded domain $\Omega \subset \C^n$
(having connected boundary), then $\mathcal{ W}$-removability of a
compact set $K \subset \partial \Omega$ entails its removability in the
sense of Problem~3.2, thanks to Hartogs Theorem~1.8(V).

As in~\cite{ jo1999b, mp1999}, we would like to emphasize that all
the general theorems presented in Sections~3 and~4 are void for $L_{
loc}^1$ functions, or require a strong assumption of growth. On the
contrary, the results that will be presented below hold in all spaces
$L_{ loc}^{\sf p}$ with $1 \leqslant {\sf p} \leqslant \infty$,
without any assumption of growth. The concept of $\mathcal{
W}$-removability, interpreted as a result about envelopes of
holomorphy, yields a (crucial) external drawing near the illusory
singularity, an opportunity that is intrinsically attached to locally
embeddable Cauchy-Riemann structures, but is of course absent for
general linear partial differential operators.

\subsection*{ 4.3.~Removable singularities on hypersurfaces}
In~\cite{ ls1993}, it is shown that if $\Omega \subset \C^n$ is a
pseudoconvex bounded domain having $\mathcal{ C}^2$ boundary, then
every compact subset $K \subset \partial \Omega$ with ${\sf H}^{ 2n -
3} (K) = 0$ is removable in the sense of Definition~3.4. In fact,
Lemma~4.18(III) shows that $\partial \Omega$ is globally minimal
and the next lemma shows that in codimension $d = 1$, metrically thin
singularities do not perturb global minimality.

\def\thelemma{4.4}\begin{lemma}
{\rm (\cite{ mp2002})} If $M\subset \C^n$ is a globally minimal
$\mathcal{ C}^2$ hypersurface, then for every closed set $C \subset M$
with ${\sf H}^{ 2n - 3} (C) = 0$, the complement $M \backslash C$ is
also globally minimal.
\end{lemma}

\def\theexample{4.5}\begin{example}{\rm
However, this is untrue if ${\sf H}^{ 2n -3} (C) > 0$. Let $n \geqslant 2$
and $\varphi (z, u)$ be $\mathcal{ C}^2$ defined for $\vert z\vert,
\vert u \vert < 1$ and satisfying $\varphi(z,0) \equiv 0$ for ${\rm
Re}\, z_1 \leqslant 0$. Let $M \subset \C^n$ be the graph $v = \varphi (z,
u)$ and define $C := \{ (i\, y_1, z_2, \dots, z_{ n-1}, 0)
\}$. Clearly $\dim C = 2n-3$, ${\sf H}^{ 2n-3} (C) >0$ and $\{ (z, 0):
{\rm Re}\, z_1 < 0\}$ is a single CR orbit $\mathcal{ O}_-$ of $M
\backslash C$. Also, the function $\varphi$ may be chosen so that $M$
is of finite type at every point of $M \backslash \mathcal{ O}_-$,
whence $M\backslash C$ consists of exactly two CR orbits, namely
$\mathcal{ O}_-$ and $M\backslash (\mathcal{ O}_- \cup C)$. It follows
that $M$ is globally minimal.
}\end{example}

\def\thetheorem{4.6}\begin{theorem}
{\rm (\cite{ ls1993, cst1994, mp1998, mp2002})} If $M \subset \C^n$ is
a globally minimal $\mathcal{ C}^{ 2, \alpha}$ {\rm (}$0 < \alpha <
1${\rm )} hypersurface, then every closed set $C \subset M$ with ${\sf
H}^{ 2n -3} (C) = {\sf H}^{ \dim M - 2} ( C) = 0$ is {\rm locally}
CR-, $\mathcal{ W}$- and $L^{ \sf p}$-removable.
\end{theorem}

Sometimes, we shall say that $C$ is of codimension $2^{+0}$ in $M$.
This is a version of {\bf (rm1)} and of {\bf (rm2)} of Theorem~2.30
for CR functions on general hypersurfaces. Except for $L^{ \sf
p}$-removability, refinements about smoothness assumptions may be
found in~\cite{ cst1994}. 

\smallskip

The smallest (Hausdorff) dimension of $C \subset M
\subset \C^n$ for which its removability may fail is equal to
$2n-3$. Indeed, if $C = M\cap \Sigma$ is equal to the intersection of
$M$ with some local complex hypersurface $\Sigma = \{ f = 0\}$, the
functions $1 / f^k$, $k\geqslant 1$ and $e^{ 1/f}$
restrict to be CR on $M \backslash C$, but not holomorphically
extendable to a one-sided neighborhood at points of $C$, since
$\Sigma$ visits both sides of $M$. In such a situation, the real
hypersurface $M \cap \Sigma$ of the complex hypersurface $\Sigma$ has
dimension $(2n-3)$ and CR dimension $(n-2)$.

\def\thedefinition{4.7}\begin{definition}{\rm
A CR submanifold $N \subset \C^n$ is called {\sl maximally complex} if
it is of odd dimension satisfying $\dim N = 1 + 2\, {\rm CRdim}\, N$.
}\end{definition}

Every real hypersurface of a complex manifold is maximally complex. The
next step in to study singularities $C$ contained in $( 2n-3
)$-dimensional submanifolds $N \subset M$. 

\def\theexample{4.8}\begin{example}{\rm
We show the necessity of assuming that $M\backslash C$ is also
globally minimal (\cite{ mp1999}). Take the complex hypersurface
$\mathcal{ O}_-$ of the preceding example having boundary $\partial
\mathcal{ O}_- = C = N$. Applying Proposition~4.38(III) to $S :=
\mathcal{ O}_-$, we may construct a measure on $M\backslash C$
supported by $\mathcal{ O}_-$ that is CR on $M\backslash C$ but does
not extend holomorphically to a wedge at any point of $\overline{
\mathcal{ O}}_- = \mathcal{ O}_- \cup C$, for the same reason as in
Corollary~4.39(III).
}\end{example}

Because of this example, we shall systematically assume that
$M \backslash C$ is also globally minimal, if this is not a consequence
of other hypotheses. Here is a CR version of {\bf (rm3)} and of {\bf
(rm4)} of Theorem~2.30. It says that true singularities should be
maximally complex. Before stating it, we point out that all
submanifolds of given manifolds will constantly be assumed to be {\it
embedded}\, submanifolds. Also, all subsets $C$ of a submanifold $N$
of manifold $M$ that are called {\sl closed}\,
are assumed to be closed both in $M$ and in $N$.

\def\thetheorem{4.9}\begin{theorem}
{\rm (\cite{ jo1992, me1997, jo1999a, jo1999b})} Let $M \subset \C^n$
be a $\mathcal{ C}^{ 2, \alpha}$ {\rm (}$0 < \alpha < 1${\rm )}
globally minimal hypersurface and let $N \subset M$ be a connected
$\mathcal{ C}^{2, \alpha}$ embedded submanifold of dimension $(2n-3)$,
viz. of codimension $2$ in $M$. A closed set $C \subset N$ such that
$M\backslash C$ is also globally minimal is CR-, $\mathcal{ W}$- and
$L^{\sf p}$-removable under each one of the following two
circumstances{\rm :}

\begin{itemize}

\smallskip\item[{\bf (i)}]
$n\geqslant 2$ and $C \neq N${\rm ;}

\smallskip\item[{\bf (ii)}]
$n\geqslant 3$ and $C=N$ is not maximally complex, viz. there exists at
least one point $p \in N$ at which $N$ is generic.

\end{itemize}\smallskip
\end{theorem}

One may verify (\cite{ jo1999a, mp1999}) that generic points of $N$
are locally removable and then after erasing them by deforming
slightly $M$ inside the extensional wedge existing above, {\bf (ii)}
is seen to be a consequence of {\bf (i)}. For various smoothness
refinements, the reader is referred to \cite{ jo1992, cst1994, mp1998,
jo1999a, jo1999b, mp1999}. One may also combine Theorem~4.6 and~4.9,
assuming that the submanifold $N$ is smooth, except perhaps at all
points of some metrically thin closed subset. The proof will not be
restituted.

The study of more massive singularities contained in
$(2n-2)$-dimensional submanifolds has been initiated by J\"oricke
(\cite{ jo1988}), having in mind some generalization of Denjoy's
approach to Painlev\'e's problem.

\def\thetheorem{4.10}\begin{theorem}
{\rm (\cite{ jo1999a, jo1999b})} Let $M \subset \C^n$ be a $\mathcal{
C}^{ 2, \alpha}$ {\rm (}$0 < \alpha < 1${\rm )} globally minimal
hypersurface and let $M^1 \subset M$ be a connected $\mathcal{ C}^{ 2,
\alpha}$ embedded submanifold\footnote{ We believe that $\mathcal{
C}^{ 2, \alpha}$-smoothness of $M^1$ is required in the proof built
there, since the map $w \mapsto \widehat{ h} (w)$ appearing in
equation ~\thetag{ 3.12} of~\cite{ jo1999a} \big(that corresponds
essentially to the singular integral $\mathcal{ J} (v)$ defined
in~\thetag{ 3.20}(V)\big) already requires $M^1$ to be $\mathcal{ C}^{
1, \alpha}$ with $0 < \alpha < 1$ to exist; then to compute the
differential of $w \mapsto \widehat{ h} (w)$, one must require $M^1$
to be at least $\mathcal{ C}^{ 2, \alpha}$.} of dimension $(2n-2)$,
viz. of codimension $1$ in $M$, that is {\rm generic} in $\C^n$. If
$n\geqslant 3$, a closed set $C \subset M^1$ is CR-, $\mathcal{ W}$-
and $L^{\sf p}$-removable provided it does not contain any CR orbit of
$M^1$.
\end{theorem}

It may be established \big({\it see} {\it e.g.}
Lemma~3.3 in~\cite{ mp2006a}\big) 
that $M^1 \backslash C'$ is also globally minimal
for every closed $C' \subset M^1$ containing no CR orbit of $M^1$.

\smallskip

We would like to mention that the removability of two-codimensional
singularities (Theorem~4.9) is {\it not}\, a consequence of the
removability of the bigger one-codimensional singularities
(Theorem~4.10). Indeed, it may happen that $T_p N$ contains $T_p^c M$
at several points $p\in N$ in Theorem~4.9, preventing the existence of
a generic $M^1 \subset M$ containing $N$. In addition, even if $T_p N
\not\supset T_p^c M$ for every $p\in N$, Theorem~4.9 is {\it not
anymore}\, a corollary of Theorem~4.10. Indeed, with $m=2$ and $d=1$,
choosing a local hypersurface $M \subset \C^3$ containing a complex
curve $\Sigma$, choosing $N \subset M$ of dimension $3$ containing
$\Sigma$ and being maximally real outside $\Sigma$, and choosing an
arbitrary generic $M^1 \subset M$ containing $N$ (some explicit local
defining equations may easily be written), then $\Sigma$ is a CR orbit
of $M^1$, so $N \supset \Sigma$ is not considered to be removable by
Theorem~4.10, whereas Theorem~4.9{\bf (ii)} asserts that $N$ is
removable.

Although singularities are more massive in Theorem~4.10, the
assumption $n \geqslant 3$ in it entails that the CR dimension $(n-1)$
of $M$ is $\geqslant 2$, whence $M^1$ has {\it positive}\, CR
dimension $\geqslant 1$. This insures the existence of small analytic
discs with boundary in $M^1$. Section~5 below and~\cite{ mp2006a} as a
whole are devoted to the more delicate case where $M^1$ has null CR
dimension.

\def\theexample{4.11}\begin{example}{\rm (\cite{ jo1999a}) In $\C^3$,
let $M = \partial \B_3$ and let $M^1 := \big\{ (z_1, z_2, z_3): \ 0 <
x_1 < 1/2, \ y_1 = 0 \big\}$. Clearly, $M^1$ is foliated by the
$3$-spheres 
\[
S_{x_1^*}^3 
:= 
\big\{ 
z_1=x_1^*,\
\vert z_2\vert^2 
+
\vert z_3\vert^2
= 
1-\vert x_1^*\vert^2
\big\}, 
\]
$x_1^* \in (0, 1/2)$, that are globally minimal compact
$3$-dimensional strongly pseudoconvex maximally complex CR
submanifolds of CR dimension $1$ bounding the $2$-dimensional complex
balls
\[
\B_{2,x_1^*}
:= 
\big\{ 
z_1=x_1^*,\
\vert z_2\vert^2 
+
\vert z_3\vert^2
<
1-\vert x_1^*\vert^2
\big\}.
\]
Theorem~4.10 asserts that a compact set $K \subset M^1$ is removable
if and only if it does not contain a whole sphere $S_{ x_1^*}$, for
some $x_1^* \in ( 0, 1/ 2)$. If $K$ contains such a sphere $S_{
x_1^*}^3$, the complex $2$-ball $\B_{2, x_1^*}$ coincides with the $A
(\Omega)$-hull of $S_{ x_1^*}$ and is nonremovable. More generally, an
application of both Theorems~4.10 and~3.25 yields the following.
}\end{example}

\def\thecorollary{4.12}\begin{corollary}
Let $K$ be a compact subset of $M^1$. For every 
{\rm (}interior{\rm )} one-sided
neighborhood $\mathcal{ V}^- \big( \partial \B_3 \backslash K \big)$
that is contained in $\B_3$ and every function $f$ holomorphic in
$\mathcal{ V}^- \big( \partial \B_3 \backslash K \big)$, there exists
a function $F$ holomorphic in $\B_3 \big\backslash \bigcup_{ x_1^*: \
S_{ x_1^*}^3 \subset K}\, \B_{ 2, x_1^*}$ with $F = f$ in $\mathcal{
V}^- \big( \partial \B_3 \backslash K \big)$.
\end{corollary}

By means of the complex Plateau problem, the next paragraph discusses
the necessity for $N$ not to be maximally complex in Theorem~4.9 and
for $M^1$ not to contain any CR orbit in Theorem~4.10, in a
more general context than $M = \partial \B_n$.

\subsection*{ 4.13.~Complex Plateau problem and nonremovable
singularities contained in strongly pseudoconvex boundaries} Let
$\mathcal{ M}$ be a complex manifold of dimension $n\geqslant 2$. If
$\Sigma \subset \mathcal{ M}$ is a closed pure $k$-dimensional complex
subvariety, we denote by $[ \Sigma ]$ the current of integration on
$\Sigma$, whose existence was established by Lelong in 1957 (\cite{
ch1989, de1997}).

\def\thedefinition{4.14}\begin{definition}{\rm
(\cite{ hl1975, ha1977}) A current ${\sf T}$ on $\mathcal{ M}$ is
called a {\sl holomorphic $k$-chain} if it is of the form
\[
{\sf T}
=
\sum_{\rm finite}\,
n_j[\Sigma_j],
\]
where the $\Sigma_j$ denote the irreducible components of a pure
$k$-dimensional complex subvariety $\Sigma$ of $\mathcal{ M}$ and
where the {\sl multiplicity} $n_j$ of each $\Sigma_j$ is an integer.
}\end{definition}

The complex Plateau problem consists in filling boundaries $N$ by
complex subvarieties $\Sigma$, or more generally by holomorphic chains
${\sf T}$. Maximal complexity of the boundary $N$ is naturally
required and since $N$ might encounter singular points of $\Sigma$, it
should be allowed in advance to be ``scarred'' somehow. Also, the
boundary $N$ inherits an orientation from $\Sigma$ and as the boundary
of $\Sigma$, it should have empty boundary.

\def\thedefinition{4.15}\begin{definition}{\rm
A {\sl scarred $\mathcal{ C}^\kappa$ {\rm (}$1\leqslant \kappa
\leqslant \infty${\rm )} maximally complex cycle} of dimension $(2m +
1)$, $m \geqslant 0$, is a compact subset $N \subset \mathcal{ M}$
together with a thin compact {\sl scar set} $\text{\sf sc}_N \subset
N$ such that

\begin{itemize}

\smallskip\item[$\bullet$]
${\sf H}^{ 2m+1} (\text{\sf sc}_N) = 0${\rm ;}

\smallskip\item[$\bullet$]
$N \backslash \text{ \sf sc}_N$ is an {\it oriented}\, $( 2m + 1
)$-dimensional embedded maximally complex $\mathcal{ C }^\kappa$
submanifold of $\mathcal{ M} \backslash \text{ \sf sc}_N$ having
finite $(2m+1)$-dimensional Hausdorff measure{\rm ;}

\smallskip\item[$\bullet$]
the current of integration over $N \backslash \text{\sf sc}_N$,
denoted by $[ N]$, has no boundary: $d[ N] = 0$.

\end{itemize}\smallskip
}\end{definition}

This definition was essentially devised by Harvey-Lawson and
appears to be adequately large, but sufficiently stringent to maintain
the possibility of filling a maximally complex cycle by a complex
analytic set.

\def\thetheorem{4.16}\begin{theorem}
{\rm (\cite{ hl1975, ha1977})} Suppose $N$ is a scarred $\mathcal{
C}^\kappa$ {\rm (}$1\leqslant \kappa \leqslant \infty${\rm )}
maximally complex cycle of dimension $(2m+1)$, $m\geqslant 0$, in a
Stein manifold $\mathcal{ M}$.

\begin{itemize}

\smallskip\item[$\bullet$]
If $m=0$, assume that $N$ satisfies the {\rm moment condition},
viz. $\int_N \, \omega = 0$ for every holomorphic $1$-form $\omega =
\sum_{ k = 1}^n \, \omega_k (z) \, dz_k$ having entire coefficients
$\omega_k \in \mathcal{ O} (\C^n)$.

\smallskip\item[$\bullet$]
If $m\geqslant 1$, assume nothing, since the corresponding appropriate
moment condition follows automatically from the assumption of maximal
complexity {\rm (\cite{ hl1975})}.

\end{itemize}\smallskip

Then there exists a unique holomorphic $(m+1)$-chain ${\sf T}$ in
$\mathcal{ M} \backslash N$ having compact support and finite mass in
$\mathcal{ M}$ such that
\[
d{\sf T}
=
[N]
\]
in the sense of currents in $\mathcal{ M}$. Furthermore, there is a
compact subset $K$ of $N$ with ${\sf H}^{ 2m+1} (K) = 0$ such that
every point of $N \backslash \big( K \cup \text{\sf sc}_N \big)$
possesses a neighborhood in which $({\rm supp}\, {\sf T} ) \cup N$ is
a regular $\mathcal{ C}^\kappa$ complex manifold with boundary.
\end{theorem}

A paradigmatic example, much considered since Milnor studied it,
consists in intersecting a complex algebraic subvariety of $\C^n$
passing through the origin with a spere centered at $0$; topologists
usually require that $0$ is an isolated singularity and that the
sphere is small or that the defining polynomial is homogeneous.

\smallskip

We apply this filling theorem in a specific situation. Let
$\partial \Omega \Subset \C^n$ {\rm (}$n\geqslant 3${\rm )} be a
strongly pseudoconvex $\mathcal{ C}^2$ boundary and let $M^1 \subset
\partial \Omega$ be an embedded $\mathcal{ C}^2$ one-codimensional
submanifold that is generic in $\C^n$. We assume that $M^1$ has no
boundary and is closed, viz. is a compact submanifold. Since $M^1$ has
CR dimension $(n-2)$, its CR orbits have dimension equal to either
$(2n-4)$, or to $(2n-3)$ or to $(2n-2)$. Because of
Corollary~4.19(III), no CR orbit of $M^1$ can be an immersed complex
$2$-codimensional submanifold, of real dimension $(2n - 4)$, since its
closure in $M^1$ would be a compact set laminated by complex
manifolds. Nevertheless, there may exist $(2n-3)$-dimensional CR
orbits.

\def\theproposition{4.17}\begin{proposition}
{\rm (\cite{ jo1999a})} Every CR orbit $\mathcal{ O}_{ CR}^1$ of a
connected $\mathcal{ C }^2$ hypersurface $M^1 \subset \partial \Omega$
of a $\mathcal{ C }^2$ strongly pseudoconvex boundary $\partial \Omega
\Subset \C^n$ is of the following types{\rm :}

\begin{itemize}

\smallskip\item[{\bf (i)}]
$\mathcal{ O}_{ CR}^1$ is an open subset of $M${\rm ;}

\smallskip\item[{\bf (ii)}]
$\mathcal{ O}_{ CR}^1$ is a closed maximally complex $\mathcal{ C}^1$
cycle embedded in $M^1${\rm ;}

\smallskip\item[{\bf (iii)}]
$\mathcal{ O}_{ CR}^1$ is a maximally complex $\mathcal{ C}^1$
submanifold injectively immersed in $M^1$ whose closure $C$ consists
of an uncountable union of similar CR orbits.

\end{itemize}\smallskip

\end{proposition}

In the last situation, $C$ will be called a {\sl maximally complex
exceptional minimal compact CR-invariant set}. The intersection of
$C$ with a local curve transversal to a piece CR orbit in $M^1$ may
consist of either an open segment or of a Cantor (perfect) subset.

Here is the desired converse to both Theorems~4.9 and~4.10 in a
situation where the Plateau complex filling works.

\def\thecorollary{4.18}\begin{corollary}
{\rm (\cite{ jo1999a})} Suppose that $\partial \Omega \in \mathcal{
C}^{ 2, \alpha}$ contains a compact embedded $(2n-3)$-dimensional
maximally complex submanifold $N$ {\rm (}without boundary{\rm )}.
Then $N$ is {\rm not} removable.
\end{corollary}

\proof
Indeed, the scar set of $N$ is empty and the filling of $N$ by a
holomorphic chain consists of an irreducible complex subvariety
$\Sigma$ that is necessarily contained in $\Omega$, since $\partial
\Omega$ is strongly pseudoconvex. Then the domain $\Omega \backslash
\Sigma$ is seen to be pseudoconvex and $\widehat{ N}_{ A (\Omega)} = N
\cup \Sigma$. Theorem~3.25 entails that CR functions on $\partial
\Omega \backslash N$ extend holomorphicaly to $\Omega \backslash
\Sigma$.
\endproof

A very natural problem, raised in~\cite{ jo1999a} and inspired by a
perturbation of Example~4.11, is to determine for which compact
CR-invariant subsets $K$ of a strongly pseudoconvex boundary $\partial
\Omega \subset \C^n$ the envelope of holomorphy of $\partial \Omega
\backslash K$ is multi-sheeted.

\def\thetheorem{4.19}\begin{theorem}
{\rm (\cite{ js2004})} Let $M^1 \subset \partial \B_n$ be an
orientable $(2n-2)$-dimensional generic $\mathcal{ C}^{ 2, \alpha }$
submanifold of $\partial \B_n$ {\rm (}$n\geqslant 3${\rm )} and let $K
\subset M^1$ be a compact CR-invariant subset of $M^1$ such that

\begin{itemize}

\smallskip\item[$\bullet$]
the boundary of $K$ in $M^1$ is the disjoint union of finitely many
connected compact maximally complex CR manifolds $N_1, \dots, N_\ell$
of dimension $(2n-3)$ that are $\mathcal{ C}^{ 2, \alpha -0}$ CR
orbits of $M^1${\rm ;}

\smallskip\item[$\bullet$]
the interior of $K$ with respect to $M^1$ 
is globally minimal.

\end{itemize}\smallskip

Then the envelope of holomorphy ${\sf E} 
\big ( \mathcal{ V} (\partial \B_n \backslash K) \big)$
is multi-sheeted in every neighborhood $\overline{ U}_p \subset
\overline{ \B_n}$ of every point $p \in {\rm Int} \, K$.
\end{theorem}

We conclude these considerations by formulating a deeply open problem
raised by J\"oricke. The complex Plateau problem for laminated
boundaries is a virgin mathematical landscape.

\def\theopenquestion{4.20}\begin{openquestion}
{\rm (\cite{ jo1999a})} Let $\partial \Omega \Subset \C^n$,
$n\geqslant 3$, be a strongly pseudoconvex boundary of class at least
$\mathcal{ C}^2$. Suppose that $\partial \Omega$ contains a maximally
complex exceptional minimal compact CR-invariant set $C$. Does $C$
bound a relatively compact subset $\Sigma \subset \Omega$ laminated by
complex manifolds\,?
\end{openquestion}

As observed in~\cite{ dh1997, mp1998, sa1999, ds2001}, removable
singularities have an unexpected interesting application to wedge
extension of CR-meromorphic functions.

\subsection*{ 4.21.~CR-meromorphic functions
and metrically thin singularities} For $n\geqslant 2$, a local
meromorphic map $f$ from a domain $\Omega \subset \C^n$ to the Riemann
sphere $P_1 (\C)$ has an exceptional locus $I_f \subset \Omega$, at
every point $p$ of which the value $f(p)$ is undefined. For instance
the origin $(0, 0) \in \C^2$ with $f = \frac{ z_1}{ z_2}$ (notice that
every complex number in $\C \cup \{ \infty\}$ is a limit of $\frac{
z_1}{ z_2}$). This exceptional set $I_f$ is a complex analytic subset
of $\Omega$ having codimension $\geqslant 2$ (\cite{ de1997}). It is
called the {\sl indeterminacy set} of $f$.

A meromorphic function may be more conveniently defined as a
$n$-dimensional irreducibe complex analytic subset $\Gamma_f$ of
$\Omega \times P_1 (\C)$ having surjective projection onto $\Omega$,
viz. $\pi_\Omega (\Gamma_f) = \Omega$. Here, $\Omega$ might be any
complex manifold. Indeterminacy points correspond precisely to points
$p\in \Omega$ satisfying $\pi_\Omega^{ - 1} (p) \cap \Gamma_f = \{ p\}
\times P_1 (\C)$. So, the generalization of meromorphy to the CR
category incorporates indeterminacy points.

\def\thedefinition{4.22}\begin{definition}{\rm 
(\cite{hl1975, dh1997, mp1998, sa1999}) Let $M \subset \C^n$ be a
scarred $\mathcal{ C}^1$ generic submanifold of codimension
$d\geqslant 1$ and of CR dimension $m= n -d \geqslant 1$. Then a {\sl
CR meromorphic function} on $M$ with values in $P_1 (\C)$ consists of
a triple $(f, \mathcal{ D}_f, \Gamma_f )$ such that:

\begin{itemize}

\smallskip\item[{\bf 1)}]
$\mathcal{ D}_f \subset M$ is a dense open subset of $M$ and $f:
\mathcal{ D}_f \to P_1(\C)$ is a $\mathcal{ C}^1$ map;

\smallskip\item[{\bf 2)}] 
the closure $\Gamma_f$ in $\C^n \times P_1(\C)$ of the graph $\{ (p,
f(p)): \, p\in \mathcal{ D}_f\}$ defines an oriented scarred
$\mathcal{ C}^1$ CR submanifold of $\C^n \times P_1 (\C)$ of the
same CR dimension as $M$ having empty boundary
in the sense of currents.

\end{itemize}\smallskip

\noindent
The {\sl indeterminacy locus} of $f$ is denoted by
\[
I_f 
:= 
\big\{
p\in M:\ 
\{p\}\times
P_1(\C)
\subset\Gamma_f 
\big\}.
\]
}\end{definition}

\noindent
In the CR category, $I_f$ is not as thin as in the holomorphic
category (where it has real codimension $\geqslant 4$), but it is
nevertheless thin enough for future purposes, as we shall see. A
standard argument from geometric measure theory yields almost
everywhere smoothness of almost every level set.

\def\thelemma{4.23}\begin{lemma}
{\rm (\cite{ fe1969, hl1975, ha1977})} Let $M \subset \C^n$ be a scarred
$\mathcal{ C}^1$ generic submanifold. Let $(f, \mathcal{ D}_f,
\Gamma_f)$ be a CR meromorphic function on $M$. Then for almost every
$w \in P_1(\C)$, the level set
\[
N_f(w) 
:=
\big\{
p\in M:\,
(p,w)\in\Gamma_f
\big\}
\]
is a scarred $2$-codimensional $\mathcal{ C}^1$ submanifold of $M$.
\end{lemma}

Let $p\in I_f$. Since $(p, w) \in \Gamma_f$ for every $w \in P_1
(\C)$, it follows that $I_f \subset N_f (w)$ for every $w$. Fixing
such a $w \in P_1 (\C)$, we simply denote $N_f := N_f (w)$. In
particular, the scar set $\text{\sf sc}_{N_f}$ of $N_f$ is always of
codimension $2^{ +0}$ in $M$, namely
${\sf H}^{ \dim M - 2} ({\sf sc}_{ N_f}) = 0$.

So $I_f \subset N_f$ and by definition $I_f \times P_1 (\C) \subset
\Gamma_f$. We claim that, in addition, $I_f$ has empty interior in
$N_f \backslash \text{\sf sc}_{ N_f}$. Otherwise, there exist a point
$p\in N_f \backslash \text{\sf sc}_{ N_f}$ and a neighborhood $U_p$ of
$p$ in $M$ with $U_p \cap \text{\sf sc}_{ N_f} = \emptyset$ such that
$I_f$ contains $U_p \cap N_f$, whence 
\[
(U_p \cap N_f)\times P_1 (\C)
\subset
\Gamma_f.
\]
Since $(U_p \cap N_f) \times P_1 (\C)$
has dimension equal to $\dim M = \dim \Gamma_f$, it follows that
\[
\Gamma_f\cap\big(U_p\times P_1(\C)\big)
\equiv
(U_p\cap N_f) 
\times P_1(\C).
\]
But $U_p \cap N_f$ having codimension two in $U_p$, this contradicts the
assumption that $\Gamma_f$ is a (nonempty!) graph above the dense open
subset $U_p \cap \mathcal{ D}_f$ of $U_p$.

\def\thelemma{4.24}\begin{lemma}
{\rm (\cite{ mp1998, sa1999})}
The indeterminacy set $I_f$ of $f$ is a closed set of empty interior
contained in some $2$-codimensional scarred $\mathcal{ C }^1$
submanifold $N_f$ of $M$. Moreover, the scar set $\text{\sf sc}_{N_f}$
of $N_f$ is always of codimension $2^{ +0}$ in $M$, viz. ${\sf H}^{ 2m
+ d - 2} (\text{\sf sc}_{N_f}) = 0$.
\end{lemma}

The statement below and its proof are clear if $\mathcal{ D}_f = M$;
in it, the condition $d[ \Gamma_f ] = 0$ helps in an essential way to
keep it true when the closure of $\Gamma_f$ possesses a nonempty scar
set.

\def\theproposition{4.25}\begin{proposition}
{\rm (\cite{ mp1998, sa1999})} There exists a unique CR measure ${\sf
T}_f$ on $M \backslash I_f$ with ${\sf T}_f \vert_{\mathcal{
D}_f}$ coinciding with the $\mathcal{ C}^1$ CR function $f : \mathcal{
D}_f \to P_1 (\C)$.
\end{proposition}

It is defined locally as follows. Let $p\in M \backslash I_f$
and let $U_p$ be an open neighborhood of $p$ in $M$. Since $p \not\in
I_f$, there exists $w_p \in P_1 (\C)$ with $(p, w_p) \not\in
\Gamma_f$. Composing with an automorphism of $P_1 (\C)$ and shrinking
$U_p$, we may assume that $w_p = \infty$ and that $\big( U_p \times \{
\infty\} \big) \cap \Gamma_f = \emptyset$. Letting $d {\rm Vol}_{
U_p}$ be some $(2m+d)$-dimensional volume form on $U_p$, letting
$\pi_{ \Gamma_f} : \Gamma_f \to M$ denote the natural projection, the
CR measure ${\sf T}_f \big\vert_{ U_p}$ is defined by
\[
\left<
{\sf T}_f,\varphi
\right>
:=
\int_{\Gamma_f}\,
w\cdot\pi_{\Gamma_f}^* (\varphi \ d{\rm Vol}_{ U_p}),
\]
for every $\varphi \in \mathcal{ C}_c^1 (U_p)$.

\smallskip

Thus, on $M \backslash I_f$, the CR-meromorphic function $(f,
\mathcal{ D}_f, \Gamma_f )$ behaves like an order zero CR
distribution. With $\mathcal{ C}_{ CR }^0$, $L_{ CR, loc }^{ \sf p}$,
it therefore enjoys the extendability properties of Part~V on $M
\backslash I_f$, provided that $M$ is $\mathcal{ C}^{ 2, \alpha
}$. The next theorem should be applied to $C := I_f$. Its final proof
(\cite{ mp2002}) under the most general assumptions combines both the
CR extension theory and the application of the Riemann-Hilbert problem
to global discs attached to maximally real submanifolds (\cite{
gl1994, gl1996}). We cannot restitute the proof here.

\def\thetheorem{4.26}\begin{theorem} 
{\rm (\cite{ mp1998, ds2001, mp2002})} Suppose $M \subset \C^n$ is
$\mathcal{ C}^{ 2, \alpha }$ {\rm (}$0 < \alpha < 1${\rm )} of
codimension $d \geqslant 1$ and of CR dimension $m \geqslant 1$. Then every
closed subset $C$ of $M$ such that $M$ and $M \backslash C$ are
globally minimal and such that ${\sf H }^{ 2m + d - 2} (C) = 0$ is CR-,
$\mathcal{ W }$- and $L^{ \sf p}$-removable.
\end{theorem}

However, if $f$ is a CR-meromorphic function defined on such a $M$,
with $\mathcal{ C}^1$ replaced by $\mathcal{ C}^{ 2, \alpha}$ in
Definition~4.22, the complement $M \backslash I_f$ need not be globally
minimal if $M$ is, and it is easy to construct manifolds $M$ and
closed sets $C \subset M$ with ${\sf H}^{ 2m-1} (C) < \infty$ which
perturb global minimality, {\it cf.} Example~4.8. It is therefore
natural to make the additional assumption that $M$ is locally minimal
at {\it every} point. This assumption is the weakest one that insures
that $M \backslash C$ is globally minimal, for arbitrary closed sets
$C \subset M$.

\def\thecorollary{4.27}\begin{corollary}
Assume that $M \in \mathcal{ C}^{ 2, \alpha}$ is locally minimal at
every point and let $f$ be a CR-meromorphic function. Then $I_f$ is
CR-, $\mathcal{ W }$- and $L^{\sf p}$-removable.
\end{corollary}

\proof
Lemma~4.24 holds with $\mathcal{ C}^1$ replaced by $\mathcal{ C}^{ 2,
\alpha}$. It says that $I_f$ is a closed subset with empty interior of
some scarred $\mathcal{ C}^{ 2, \alpha}$ submanifold $N_f$ of $M$. The
removability of the portion of $I_f$ that is contained in the regular
part of $N_f$ follows from Theorem~4.9{\bf (i)}. The removability of
the remaining scar set $\text{\sf sc}_{ N_f}$ follows from
Theorem~4.26 above.
\endproof

Thus the CR measure ${\sf T}_f$ on $M \backslash I_f$
(Proposition~4.25) extends holomorphically to some wedgelike domain
$\mathcal{ W}_1$ attached to $M \backslash I_f$. The $\mathcal{
W}$-removability of $I_f$ entails that the envelope of holomorphy of
$\mathcal{ W}_1$ contains a wedgelike domain $\mathcal{ W}_2$ attached
to $M$. Performing supplementary gluing of discs, the CR extension
theory (Part~V) insures that such a $\mathcal{ W}_2$ depends only on
$M$, not on $f$. As envelopes of meromorphy and envelopes of
holomorphy of domains in $\C^n$ coincide by a theorem going back to
Levi (\cite{ ks1967, iv1992}), we may conclude.

\def\thetheorem{4.28}\begin{theorem}
{\rm (\cite{ mp2002})} Suppose $M \subset \C^n$ is $\mathcal{ C}^{ 2,
\alpha }$ and locally minimal at every point. Then there exists a
wedgelike domain $\mathcal{ W}$ attached to $M$ to which every
CR-meromorphic function on $M$ extends meromorphically.
\end{theorem}

\subsection*{ 4.29.~Peak and smooth removable singularities in 
arbitrary codimension} A closed set $C \subset M$ is called a
$\mathcal{ C}^{ 0, \beta}$ peak set, $0< \beta <1$, if there exists a
{\it nonconstant} function $\varpi \in \mathcal{ C}_{ C R}^{0, \beta}
(M )$ such that $C = \{\varpi = 1\}$ and $\max_{ p\in M}\, \vert
\varpi (p) \vert \leqslant 1$.

\def\thetheorem{4.30}\begin{theorem}
{\rm (\cite{ kr1995, mp1999})} Let $M$ be $\mathcal{ C }^{ 2, \alpha
}$ {\rm (}$0 < \alpha < 1${\rm )} globally minimal. Then every
$\mathcal{ C}^{ 0, \beta}$ peak set $C$ satisfies ${\sf H}^{ \dim M} (
C) = 0$ and is $L^{\sf p}$-removable.
\end{theorem}

To conclude, we mention two precise generalizations of
Theorems~4.9 and~4.10 to higher codimension. If $\Sigma = \{ z: \,
g(z) = 0\}$ is a local complex hypersurface passing through a point
$p$ of a generic submanifold $M \subset \C^n$ that is transverse to
$M$ at $p$, viz. $T_p \Sigma + T_p M = T_p \C^n$, the intersection
$\Sigma \cap M$ is a two-codimensional submanifold of $M$ that is
nowhere generic in a neighborhood of $p$ and certainly not 
(locally) removable,
since the CR function $\frac{ 1}{ g (z)} \big \vert_{ M \backslash
(\Sigma \cap M)}$ is not extendable to any local wedge at $p$.

\def\thetheorem{4.31}\begin{theorem}
{\rm (\cite{ me1997, mp1999})} Let $M \subset \C^n$ be a $\mathcal{
C}^{ 2, \alpha}$ {\rm (}$0 < \alpha < 1${\rm )} globally minimal
generic submanifold of positive codimension $d\geqslant 1$ and of
positive CR dimension $m = n - d \geqslant 1$. Let $N \subset M$ be a
connected two-codimensional $\mathcal{ C }^{ 2, \alpha}$ submanifold
and assume that $M\backslash N$ is also globally minimal. A closed
set $C \subset N$ is CR-, $\mathcal{ W }$- and $L^{ \sf p}$-removable
under each one of the following two circumstances{\rm :}

\begin{itemize}

\smallskip\item[{\bf (i)}]
$m\geqslant 1$ and $C \neq N${\rm ;}

\smallskip\item[{\bf (ii)}]
$m\geqslant 2$ and there exists at least one point $p\in N$
at which $N$ is generic.

\end{itemize}\smallskip
\end{theorem}

In {\bf (ii)}, the assumption that $m \geqslant 2$ is essential.
Generally, if $m = 1$, whence $d = n - 1$ and $\dim M = n + 1$, a
local transverse intersection $C = \Sigma \cap M$ has dimension $n-1$,
hence cannot be generic, and is not (locally) removable by
construction. In the next statement, the similar assumption that
$m\geqslant 2$ is strongly used in the proof: the one-codimensional
submanifold $M^1 \subset M$ has then CR dimension $m- 1 \geqslant 1$,
hence there exist small Bishop discs attached to $M^1$.

\def\thetheorem{4.32}\begin{theorem}
{\rm (\cite{ po1997, me1997, po2000})} Let $M \subset \C^n$ be a
$\mathcal{ C}^{ 2, \alpha}$ {\rm (}$0 < \alpha < 1${\rm )} globally
minimal generic submanifold of positive codimension $d\geqslant
1$. Assume that the CR dimension $m = n - d$ of $M$ satisfies
$m\geqslant 2$. Let $M^1 \subset M$ be a connected $\mathcal{ C}^{ 2,
\alpha}$ one-codimensional submanifold that is generic in $\C^n$. A
closed set $C \subset M^1$ is CR-, $\mathcal{ W}$- and $L^{ \sf
p}$-removable provided it does not contain any CR orbit of $M^1$.
\end{theorem}

Three geometrically different proofs of this theorem will be
restituted in Section~10 of~\cite{ mp2006a}. The next Section~5
and~\cite{ mp2006a} are devoted to the
study of the more delicate case where $m=1$ and where $C$ is contained
in some one-codimensional submanifold $M^1 \subset M$.

\section*{ \S5.~Removable singularities in CR dimension $1$}

\subsection*{ 5.1.~Removability of totally real 
discs in strongly pseudoconvex boundaries} In 1988, applying a global
version of the Kontinuit\"atssatz, J\"oricke~\cite{ jo1988}
established a remarkable theorem, opening the way to a purely
geometric study of removable singularities.

\def\thetheorem{5.2}\begin{theorem}
{\rm (\cite{ jo1988})} Let $\partial \Omega \Subset \C^2$ be a
strongly pseudoconvex $\mathcal{ C}^2$ boundary and let $D \subset
\partial \Omega$ be a $\mathcal{ C}^2$ one-codimensional submanifold
that is diffeomorphic to the unit open $2$-disc of $\R^2$ and
maximally real at every point. Then every compact subset $K$ of $D$
is CR-, $L^\infty$ and $\mathcal{ W}$-removable.
\end{theorem}

By maximal reality of $D$, the line distribution $D \ni p \mapsto
\ell_p := T_p D \cap T_p^c M$ is nowhere vanishing and may be
integrated. This yields the {\sl characteristic foliation} $\mathcal{
F}_D^c$ on $D$. The compact set $K$ is contained in a slightly
smaller disc $D' \Subset D$ having $\mathcal{ C}^2$ boundary $\partial
D'$. Poincar\'e-Bendixson's theorem on such a disc $D'$ together with
the inexistence of singularities of $\mathcal{ F}_D^c$ entail that
every characteristic curve that enters into $D'$ must exit from $D'$.
Orienting then the real $2$-disc $D$ and its characteristic foliation,
we have the following topological observation (at the very core of the
theorem) saying that there always exists a characteristic leaf that is
not crossed by the removable compact set.

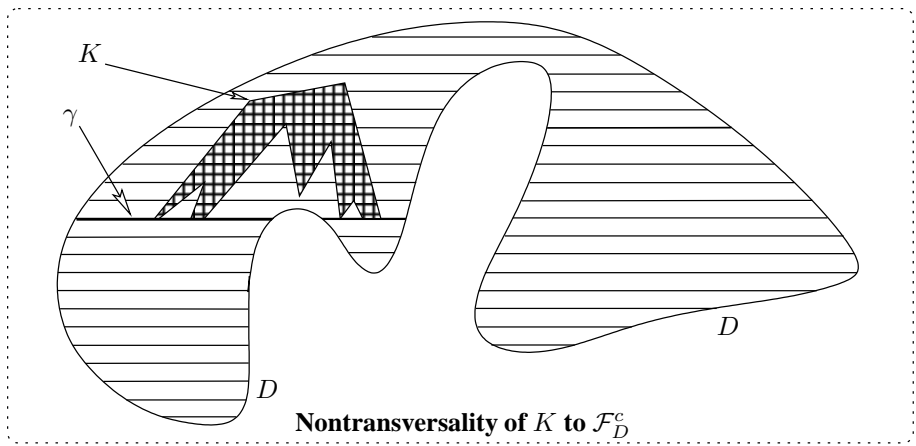
\begin{center}
\begin{picture}(0,0)%
\includegraphics{F_DcK.pstex}%
\end{picture}%
\setlength{\unitlength}{4144sp}%
\begingroup\makeatletter\ifx\SetFigFont\undefined
\def\x#1#2#3#4#5#6#7\relax{\def\x{#1#2#3#4#5#6}}%
\expandafter\x\fmtname xxxxxx\relax \def\y{splain}%
\ifx\x\y   
\gdef\SetFigFont#1#2#3{%
  \ifnum #1<17\tiny\else \ifnum #1<20\small\else
  \ifnum #1<24\normalsize\else \ifnum #1<29\large\else
  \ifnum #1<34\Large\else \ifnum #1<41\LARGE\else
     \huge\fi\fi\fi\fi\fi\fi
  \csname #3\endcsname}%
\else
\gdef\SetFigFont#1#2#3{\begingroup
  \count@#1\relax \ifnum 25<\count@\count@25\fi
  \def\x{\endgroup\@setsize\SetFigFont{#2pt}}%
  \expandafter\x
    \csname \romannumeral\the\count@ pt\expandafter\endcsname
    \csname @\romannumeral\the\count@ pt\endcsname
  \csname #3\endcsname}%
\fi
\fi\endgroup
\begin{picture}(5442,2624)(412,-2061)
\put(1909,-1780){\makebox(0,0)[lb]{\smash{\SetFigFont{10}{12.0}{rm}{\color[rgb]{0,0,0}$D$}%
}}}
\put(4664,-1400){\makebox(0,0)[lb]{\smash{\SetFigFont{10}{12.0}{rm}{\color[rgb]{0,0,0}$D$}%
}}}
\put(2144,-1975){\makebox(0,0)[lb]{\smash{\SetFigFont{10}{12.0}{rm}{\color[rgb]{0,0,0}{\bf Nontransversality of $K$ to $\mathcal{F}^c_D$}}%
}}}
\put(754,-125){\makebox(0,0)[lb]{\smash{\SetFigFont{10}{12.0}{rm}{\color[rgb]{0,0,0}$\gamma$}%
}}}
\put(844,230){\makebox(0,0)[lb]{\smash{\SetFigFont{10}{12.0}{rm}{\color[rgb]{0,0,0}$K$}%
}}}
\end{picture}

\end{center}

\begin{itemize}

\smallskip\item[$\mathcal{ F }_D^c \{K \}:$]
{\it For every compact subset $K' \subset K$, there exists a Jordan
curve $\gamma : [ -1, 1] \to D$, whose range is contained in a single
leaf of the characteristic foliation $\mathcal{ F }_D^c$, with $\gamma
(-1) \not \in K'$, $\gamma (0) \in K'$ and $\gamma (1) \not \in K'$,
such that $K'$ lies completely in one closed side of $\gamma [-1, 1]$
with respect to the topology of $D$ in a neighborhood of $\gamma[ -1,
1]$.}

\end{itemize}\smallskip

In the more general context of~\cite{ mp2006a}, we will argue that
$\mathcal{ F}_D^c \{ K\}$ is the very reason why $K$ is removable. We
will then remove locally a well chosen special point $p_{ {\rm sp}} '
\in K' \cap \gamma [ - 1, 1]$. In fact, we shall establish
removability of compact subsets $K$ of general surfaces $S$ that are
not necessarily diffeomorphic to the unit $2$-disc, provided that an
analogous topological condition holds. Also, getting rid of strong
pseudoconvexity, we shall work with a globally minimal $\mathcal{ C}^{
2, \alpha}$ hypersurface of $\C^2$. Finally, we shall relax slightly
the assumption of total reality, admitting some complex tangencies.

\def\theexample{5.3}\begin{example}{\rm
Let $\Omega = \B_2$ and let $P(z) \in \C [ z]$ be a homogeneous
polynomial of degree $\geqslant 2$ having $0$ has its only
singularity. The intersection $K := \partial \B_2 \cap \{ P = 0 \}$
is a finite union of closed real algebraic curves $\simeq S^1$ that
are everywhere transverse to $T^c \partial \B_2$. We may enlarge each
curve of $K$ as a thin $\mathcal{ C}^\omega$ annulus. There is much
freedom, but every such annulus is necessarily totally real. Denote by
$S$ the union of all annuli, a surface in $\partial \B_2$. Clearly, no
component of $K$ is removable. But the theorem does not apply: on each
annulus, the characteristic foliation $\mathcal{ F}_S^c$ is radial and
$K$ crosses each characteristic leaf.

}\end{example}

\def\theexample{5.4}\begin{example}{\rm
The theorem may fail with the disc $D$ replaced by a surface $S$
having nontrivial fundamental group, even with $S$ compact without
boundary. For instance, in $\partial \B_2 = \{ \vert z_1 \vert^2 +
\vert z_2 \vert^2 = 1 \}$, the two-dimensional torus $T_2 := \big\{
\big( \frac{ 1}{ \sqrt{ 2}} \, e^{ i\, \theta_1}, \, \frac{ 1}{ \sqrt{
2}} \, e^{ i\, \theta_2 } \big) : \, \theta_1, \, \theta_2 \in \R
\big\}$ is compact and $K := T_2$ is not removable, since $\partial
\B_2 \backslash T_2$ has exactly two connected components.

}\end{example}

\def\theexample{5.5}\begin{example}{\rm
{\rm (\cite{ jo1988})}
In the same torus $T_2$, consider instead the proper compact subset $K
:= \big\{ \big( \frac{ 1}{ \sqrt{ 2}} \, e^{ i\, \theta_1}, \, \frac{
1}{ \sqrt{ 2}} \, e^{ i\, \theta_2} \big) : \vert \theta_1 \vert
\leqslant \frac{ 3\pi}{ 2}, \, \theta_2 \in \R \big\}$, diffeomorphic
to a closed annulus. It is a set fibered by circles (contained in
$\C_{ z_2}$) over the curve $\widehat{ \gamma} := \big\{ \frac{ 1}{
\sqrt{ 2}}\, e^{ i\, \theta_1} : \, \vert \theta_1 \vert \leqslant
\frac{ 3 \pi}{ 2} \big\}$ that is contained in $\C_{
z_1}$. One may verify that the condition $\mathcal{ F}_{T_2}^c \{ K
\}$ insuring removability does not hold. In fact, applying Theorem~2.2
(in the much simpler version due to Denjoy where the curve is real
analytic), the curve $\widehat{ \gamma}$ is not $(\overline{
\partial}, L^\infty )$-removable in $\C_{ z_1}$. So we may pick a
holomorphic function $\widehat{ f} (z_1) \in \mathcal{ O} \big( \C
\backslash \widehat{ \gamma} \big)$ that is bounded in $\C \cup \{
\infty \}$ but does not extend holomorphically through $\widehat{
\gamma}$. The restriction $\widehat{ f} \big\vert_{ \partial \B_2
\backslash K}$ belongs to $L^\infty (\partial \B_2)$, is CR on
$\partial \B_2 \backslash K$ but does not extend holomorphically to
$\B_2$.

}\end{example}

\smallskip

Before pursuing, we compare Theorem~5.2 and Theorem~4.10. 

In codimension $\geqslant 2$ ({\it e.g.} for curves in $\R^3$), no
satisfactory generalization of the Poincar\'e-Bendixson theory is
known and perhaps is out of reach. This gap is caused by the
complexity of the topology of phase diagrams, by the freedom that
curves have to wind wildly around limit cycles, and by the intricate
structure of singular points.

Nevertheless, in higher complex dimension $n \geqslant 3$, CR orbits
are thicker than curves and often of codimension $\leqslant 1$. For
triples $(M, M^1, C)$ as in Theorem~4.10 with $M = \partial \Omega$
being strongly pseudoconvex, one could expect that a statement
analogous to Theorem~5.2 holds true, in which the assumption that
$M^1$ has simple topology would imply automatic removability of every
compact subset $K \subset M^1$.

To be precise, let $\partial \Omega \Subset \C^n$ ($n\geqslant 3$) be
a $\mathcal{ C}^{ 2, \alpha}$ strongly pseudoconvex boundary and let
$M^1 \subset \partial \Omega$ be a $\mathcal{ C}^{ 2, \alpha}$
one-codimensional submanifold that is generic in $\C^n$. Strong
pseudoconvexity of $\partial \Omega$ entails that CR orbits of $M^1$
are necessarily of codimension $\leqslant 1$ in $M^1$. Remind that
Theorem~4.10 says that a compact subset $K$ of $M^1$ is removable
provided it does not contain any CR orbit of $M^1$. Conversely, in
the case where $M^1$ has no exceptional CR orbit, if $K$ contains a
(then necessarily compact and maximally complex) CR orbit $N$ of
$M^1$, then $K$ is not removable, since $N$ is fillable by some
$(n-1)$-dimensional complex subvariety $\Sigma \subset \Omega$ with
$\partial \Sigma = N$. Thus, while comparing the two Theorems~4.10
and~5.2, the true question is whether the assumption that $M^1 \subset
\partial \Omega = M$ be diffeomorphic to the real $(2n-2)$-dimensional
real ball $B^{ 2n -2} \subset \R^{ 2n -2}$ prevents the existence of
compact $(2n-3)$-dimensional CR orbits of $M^1$. This would yield a
neat statement, valid in arbitrary complex dimension.

For instance, let $N := \partial \B_n \cap H$ be the intersection of
the sphere $\partial \B_n \simeq S^{ 2n -1}$ with a complex linear
hyperplane $H \subset \C^n$. With such a simple $N$ homeomorphic to a
$(2n-3)$-dimensional sphere, one may verify that every $\mathcal{
C}^\infty$ submanifold $M^1 \subset \partial \Omega$ containing $N$
which is diffeomorphic to $B^{ 2n -2}$ must contain at least one
nongeneric point. Nevertheless, admitting that $N$ has slightly more
complicated topology, the expected generalization of Theorem~5.2
appears to fail, according to a discovery of
J\"oricke-Shcherbina. This confirms the strong differences between CR
dimension $m = 1$ and CR dimension $m\geqslant 2$.

\def\thetheorem{5.6}\begin{theorem}
{\rm (\cite{ js2000})} For $\varepsilon \in \R$ with $0 < \varepsilon
< 1$ close to $1$, consider the
intersection 
\[
N_\varepsilon
:=
\big\{
z_1z_2z_3
=
\varepsilon
\big\}
\cap 
\sqrt{3}\,\partial\B_3
\]
of the complex cubic $\{ z_1 z_2 z_3 = \varepsilon \}$ with the sphere
$\sqrt{ 3} \, \partial \B_3 = \{ \vert z_1\vert^2 + \vert z_2 \vert^2
+ \vert z_3 \vert^2 = 3 \}$. Then $N_\varepsilon$ is a maximally
complex cycle diffeomorphic to $S^1 \times S^1 \times S^1$ bounding
the {\rm (}nonempty{\rm )} complex surface $\Sigma_\varepsilon := \{
z_1 z_2 z_3 = \varepsilon \} \cap \B_3$. Furthermore, there exists a
suitably constructed $\mathcal{ C}^\infty$ generic one-codimensional
submanifold $M^1 \subset \partial \B_3$ diffeomorphic to the real
$(2n-2)$-dimensional unit ball $B^{ 2n -2}$ containing
$N_\varepsilon$. Finally, since $N_\varepsilon$ bounds
$\Sigma_\varepsilon$, every compact subset $K \subset M^1$ containing
$N_\varepsilon$ is nonremovable.
\end{theorem}

\subsection*{ 5.7.~Elliptic isolated complex tangencies and Bishop 
discs} Coming back to complex dimension $n=2$, we survey known
properties of isolated CR singularities of surfaces. So, let $S$ be a
two-dimensional surface $S$ in $\C^2$ of class at least $\mathcal{
C}^2$. At a point $p \in S$, the complex tangent plane $T_p S$ is
either totally (and in fact maximally) real, viz. $T_p S \cap J T_p S
= \{ 0 \}$ or it is a complex line, viz. $T_p S = J T_p S = T_p^c S$.
An appropriate application of the jet transversality theorem shows
that after an arbitrarily small perturbation, the number of complex
tangencies of $S$ is locally finite.

If $S$ has an isolated complex tangency at one of its points $p$,
Bishop (\cite{ bi1965}) showed that
there exist local coordinates $(z, w)$ centered at $p$ in which $S$
may be represented by $w = z \bar z+ \lambda (z^2+ \bar z^2) + {\rm
o}(\vert z \vert^2)$, where the real parameter $\lambda \in [0,
\infty]$ is a biholomorphic invariant of $S$. The point $p$ is said to
be {\sl elliptic} if $\lambda \in [0, \frac{1 }{ 2})$, {\sl parabolic}
if $\lambda = \frac{1 }{ 2}$ and {\sl hyperbolic} if $\lambda \in
(\frac{ 1}{ 2}, \infty ]$. The case $\lambda = \infty$ should be
understood as the surface $w = z^2 + \bar z^2 + {\rm o} (\vert
z\vert^2)$. The shape of the projection of such a surface onto the
real hyperplane $\{ {\rm Im}\, w = 0 \} \simeq \R^3$ is essentially
ellipsoid-like for $0 < \lambda < 1/2$ and essentially saddle-like for
$\lambda > 1/2$.

In the seminal article~\cite{ bi1965}, Bishop introduced this
terminology and showed that at an elliptic point, $S$ has a nontrivial
polynomial hull $\widehat{ S}$, foliated by a continuous one-parameter
family of analytic discs attached to $M$. The geometric structure of
this family has been explored further by Kenig and Webster.

\def\thetheorem{5.8}\begin{theorem}
{\rm (\cite{ kw1982, bg1983, kw1984, hu1998})} Let $S \subset \C^2$ be
a $\mathcal{ C }^\kappa$ {\rm (}$\kappa \geqslant 7${\rm )} surface
having an elliptic complex tangency at one of its points $p$. Then
there exists a $\mathcal{ C}^{ (\kappa - 7) / 3}$ one-parameter family
of disjoint regularly embdedded analytic discs attached to $S$ and
converging to $p$. If $S$ is $\mathcal{ C}^5$, then $\widehat{ S}$ is
$\mathcal{ C}^{ 0, 1}$. Furthermore, every small analytic disc
attached to $M$ near $p$ is a reparametrization of one of the discs of
the family.

For $\kappa = \infty$, the union of these discs form a $\mathcal{ C
}^\infty$ hypersurface $\widehat{ S }$ with boundary $\partial
\widehat{ S} = S$ in a neighborhood of $p$. Furthermore, $\widehat{
S}$ is the local hull of holomorphy of $S$ at $p$.
\end{theorem}

In the case where $S$ is real analytic, local normal forms may be
found that provide a classification up to biholomorphic
changes of coordinates.

\def\thetheorem{5.9}\begin{theorem}
Let $S: w = z\bar z + \lambda (z^2 + \bar z^2) + {\rm O} (\vert
z\vert^3)$ be a local {\rm real analytic} surface in $\C^2$ passing
through the origin and having an elliptic complex tangency there.

\begin{itemize}

\smallskip\item[$\bullet$]
{\rm (\cite{ mw1983})} For every $\lambda$ satisfying $0 < \lambda <
1/2$, either $S$ is locally biholomorphic to the quadric $w = z \bar z
+ \lambda (z^2+ \bar z^2)$ or there exists an integer $s\in \N$,
$s\geqslant 1$, such that $S$ is locally biholomorphic to $w = z\bar z
+ [ \lambda + \delta u^s ] (z^2 + \bar z^2)$, where $u = {\rm Re}\, w$
and $\delta = \pm 1$.

\smallskip\item[$\bullet$]
{\rm (\cite{ mo1985})} For $\lambda = 0$, either $S$ is locally
biholomorphic to $w = z\bar z + z^s + \bar z^s + {\rm O} (\vert
z\vert^{ s+1})$ for some integer $s\geqslant 3$ or 
$S$ is locally biholomorphic to $w = z\bar z$.

\smallskip\item[$\bullet$]
{\rm (\cite{ hukr1995})} For $\lambda = 0$ and $s < \infty$, the
surface $S$ is locally biholomorphic to the surface $w = z\bar z + z^s
+ \bar z^s + \sum_{ j+k > s} \, a_{ jk} \, z^j \bar z^k$, with $a_{
jk} = \overline{ a}_{ kj}$.

\end{itemize}\smallskip

In all cases, after the straightening, $S$ is contained in the real
hyperplane $\{ {\rm Im}\, w = 0 \}$.

\end{theorem}

In the third case $\lambda = 0$, $s < \infty$, it 
is still unknown how many biholomorphic invariants
$S$ can have.

\subsection*{ 5.10.~Hyperbolic isolated complex tangencies}
The existence of small Bishop discs attached to $S$ and growing at an
elliptic complex tangency impedes local polynomial convexity. At the
opposite, if $S$ is hyperbolic, Bishop's construction fails, discs are
inexistent, and in fact $S$ is locally polynomially convex.

\def\thetheorem{5.11}\begin{theorem}
{\rm (\cite{ fs1991})} Let $S \subset \C^2$ be a $\mathcal{ C}^2$
surface represented by $w = z \bar z + \lambda (z^2 + \bar z^2) + r
(z, \bar z)$, with a $\mathcal{ C }^2$ remainder $r = {\rm o} (\vert z
\vert^2)$. If $\lambda > 1 / 2$, viz. if $S$ is hyperbolic at the
origin, then for every $\rho_1 > 0$ sufficiently small, $S \cap \big(
\rho_1 \overline{ \B_2 } \big)$ is polynomially convex.
\end{theorem}

The Oka-Weil approximation theorem then assures that continuous
functions in $S \cap \big( \rho_1 \overline{ \B }_2 \big)$ are
uniformly approximable by polynomials.

A local Bishop surface $S$ is called {\sl quadratic}\, if it is
locally biholomorphic to the quadric $w = z \bar z + \lambda (z^2 +
\bar z^2)$. An isolated complex point $p$ of $S$ is called {\sl
holomorphically flat}\, if there exist local coordinates centered at
$p$ in which $S$ is locally contained in $\{ {\rm Im}\, w = 0
\}$. Unlike elliptic points of $\mathcal{ C}^\omega$ surfaces that are
always flat, hyperbolic complex points of $\mathcal{ C }^\omega$
surfaces may fail to be flat.

\def\theexample{5.12}\begin{example}{\rm 
(\cite{ mw1983})
The algebraic hyperbolic surface ($\lambda > 1/2$)
\[
w 
= 
z\bar z
+
\lambda(
z^2
+
\bar z^2)
+
\lambda
z^3 \bar z
\] 
cannot be biholomorphically transformed into a real hyperplane.

}\end{example}

Theorem~5.11 establishes local polynomial pseudoconvexity of surfaces
at hyperbolic complex tangencies. By patching together local
plurisubharmonic defining functions, one may easily construct a Stein
neighborhood basis of every surface having only finitely many
hyperbolic complex tangencies. Unfortunately, in this way one does
not control well the topology of such neighborhoods. A finer result
answering a question of Forstneri$\check{\text{\rm c }}$ is as
follows.

\def\thetheorem{5.13}\begin{theorem}
{\rm (\cite{ sl2004})} Let $S$ be a compact real $\mathcal{ C}^\infty$
surface embedded in a complex surface $\mathcal{ X}$ having only
finitely many complex points that are all hyperbolic and 
holomorphically flat. Then
$S$ possesses a basis of open neighborhoods $\big( \mathcal{
V}_\varepsilon \big)_{0 < \varepsilon < \varepsilon_1}$,
$\varepsilon_1 >0$, such that{\rm :}

\begin{itemize}

\smallskip\item[$\bullet$]
$S = \bigcap_{ \varepsilon >0} \, 
\mathcal{ V}_\varepsilon${\rm ;}

\smallskip\item[$\bullet$]
$\mathcal{ V}_\varepsilon = \bigcup_{ \varepsilon ' < \varepsilon}\, 
\mathcal{ V}_{ \varepsilon'}${\rm ;}

\smallskip\item[$\bullet$]
$\overline{ \mathcal{ V}}_\varepsilon = \bigcap_{ \varepsilon ' >
\varepsilon} \, \mathcal{ V}_{ \varepsilon '}${\rm ;}

\smallskip\item[$\bullet$]
each $\mathcal{ V}_\varepsilon$ has a $\mathcal{ C}^\infty$ strongly
pseudoconvex boundary $\partial \mathcal{ V}_\varepsilon${\rm ;}

\smallskip\item[$\bullet$]
for every $\varepsilon$ with $0 < \varepsilon < \varepsilon_1$, the
surface $S$ is a strong deformation retract of
$\mathcal{ V}_\varepsilon$.

\end{itemize}\smallskip

\end{theorem}

It is expected that the same statement remains
true without the flatness assumption.

\subsection*{ 5.14.~Real surfaces in strongly pseudoconvex
boundaries} Coming back to removable singularities, let $\partial
\Omega \Subset \C^2$ be a $\mathcal{ C}^2$ strongly pseudoconvex
boundary and let $S \subset \partial \Omega$ be a compact surface,
with or without boundary. It will be no restriction to assume that $S$
is connected. Suppose that $S$ has a finite (possibly null) number of
complex tangencies. These points then constitute the only singular
points of the characteristic foliation of $S$. At an elliptic
(resp. hyperbolic) complex tangency, the phase diagram simply looks
like a focus (resp. saddle).

\def\thetheorem{5.15}\begin{theorem}
{\rm (\cite{ fs1991})} Let $\mathcal{ M}$ be a two-dimensional Stein
manifold, let $\partial \Omega \Subset \mathcal{ M}$ be a strongly
pseudoconvex $\mathcal{ C}^2$ boundary and let $D$ be a
$\mathcal{ C}^2$ one-codimensional
submanifold that is diffeomorphic to the unit open $2$-disc of $\R^2$
and is maximally complex, except at a finite number of hyperbolic
complex tangencies. Then every compact subset $K$ of $D$ is CR- and
$\mathcal{ W}$-removable.
\end{theorem} 

Indirectly, the characterizing Theorem~3.7 of Stout yields the
following.

\def\thecorollary{5.16}\begin{corollary} 
Every compact subset $K \subset D \subset \partial \Omega$ is
$\mathcal{ O }(\overline{ \Omega})$-convex. In particular, such a $K$
is polynomially convex if $\mathcal{ M} = \C^2$ and if $\overline{
\Omega }$ is Runge or polynomially convex, {\rm e.g.} if $\Omega =
\B_2$.
\end{corollary}

The (short) 
proof mainly relies upon the (very recent in 1991 and since then
famous) works~\cite{ bk1991} and \cite{ kr1991} by Bedford-Klingenberg
and by Kruzhilin about the hulls of two-dimensional spheres contained
in such strictly pseudoconvex boundaries $\Omega \subset \mathcal{
M}$, which may be filled by Levi-flat three-dimensional spheres after
an arbitrarily small perturbation.

\def\thetheorem{5.17}\begin{theorem} {\rm (\cite{ bk1991, kr1991})} Let
$\Omega \Subset \C^2$ be a $\mathcal{ C}^6$ strongly pseudoconvex
domain and let $S \subset \partial \Omega$ be a two-dimensional sphere
of class $\mathcal{ C}^6$ embdedded into $\partial \Omega$ that is
totally real outside a finite subset consisting of $k$ hyperbolic and
$k+2$ elliptic points. Then there exist{\rm :}

\begin{itemize}

\smallskip\item[{\bf 1)}]
a smooth domain $B \subset \R^3 (x_1, x_2, x_3)$ with boundary
$\partial B$ diffeomorphic to $S$ such that $x_3 : B \to \R$ is a
Morse function on $\partial B$ having $k+2$ extreme points and $k$
saddle points, whose level sets $\{ x_3 = {\rm cst.} \} \cap B$ are
unions of finite numbers of topological discs{\rm ;} and{\rm :}

\smallskip\item[{\bf 1)}]
a continuous injective map $\Phi : B \to \Omega$ sending $\partial B$
to $S$, the extreme and saddle points of $x_3$ on $\partial B$ to the
elliptic and hyperbolic points of $S$ and the connected components of
$\{ x_3 = {\rm cst.} \} \cap B$ to geometrically
smooth holomorphic discs.

\end{itemize}\smallskip

The set $\Phi (B)$ is the envelope of holomorphy of $S$ as well as its
$\mathcal{ O} (\overline{ \Omega})$-hull, {\rm i.e.} its polynomial
hull in case $\overline{ \Omega}$ is polynomially convex.

\end{theorem}

In~\cite{ du1993}, motivated by the problem of understanding
polynomial convexity in geometric terms, the question of $\mathcal{
O}( \overline{ \Omega })$-convexity (instead of removability) of
compact subsets of arbitrary surfaces $S \subset \partial \Omega$ (not
necessarily diffeomorphic to a $2$-disc) is dealt with directly. If
$K$ is a compact subset of a totally real surface $S \subset \partial
\Omega$, denote by $\widehat{ K }_{ \rm ess } := \overline{ \widehat{
K}_{ \mathcal{ O }( \overline{ \Omega} )} \backslash K }$ the {\sl
essential $\mathcal{ O} ( \overline{ \Omega})$-hull}\, of $K$. An
application of Hopf's lemma shows that if $K = A (\partial \Delta)$ is
the boundary of a $\mathcal{ C}^1$ analytic disc $A \in \mathcal{ O}
(\Delta) \cap \mathcal{ C}^1 (\overline{ \Delta})$ attached to the
surface $S$, necessarily $K = \widehat{ K }_{ \rm ess }$ is an
immersed $\mathcal{ C}^1$ curve that is everywhere transversal to the
characteristic foliation of $S$. If $S$ has a
hyperbolic complex tangency at one of its points $p$ and if $A (1) =
p$, then $A (\partial \Delta)$ must cross at least one separatrix in
every neighborhood of $p$. When $\widehat{ K }_{ \rm ess }$ contains
no analytic disc, similar transversality properties hold.

\def\thetheorem{5.18}\begin{theorem}
{\rm (\cite{ du1993})} Let $K \Subset S \subset \partial \Omega
\Subset \C^2$ be as above, with $\partial \Omega \in \mathcal{ C}^2$
strongly pseudoconvex and $S \in \mathcal{ C}^2$ having finitely many
hyperbolic complex tangencies. In the totally real part of $S$, the
essential $\mathcal{ O} (\overline{ \Omega})$-hull $\widehat{ K }_{
\rm ess }$ of $K$ crosses every characteristic curve that it meets. If
$\widehat{ K }_{ \rm ess }$ meets a hyperbolic complex tangency, then
it meets at least two hyperbolic sectors in every neighborhood of $p$.
\end{theorem}

As a consequence (\cite{ du1993}), every compact subset $K$ of a
two-dimensional disc $D \subset \partial \Omega$ that has only
finitely many hyperbolic complex tangencies is $\mathcal{ O
}(\overline{ \Omega })$-convex.

\subsection*{ 5.19.~Totally real discs in nonpseudoconvex boundaries}
All the above results heavily relied on strong pseudoconvexity, in
contrast to the removability theorems presented in Section~6, where
the adequate statements, based on general CR extension theory, are
formulated in terms of CR orbits rather than in terms of Levi
curvature. The first theorem for the non-pseudoconvex situation was
established by the second author.

\def\thetheorem{5.20}\begin{theorem}
{\rm (\cite{ po2003})} Let $M$ be a $\mathcal{ C }^\infty$ globally
minimal hypersurface of $\C^2$ and let $D \subset M$ be a $\mathcal{ C
}^\infty$ one-codimensional submanifold that is diffeomorphic to the
unit open $2$-disc of $\R^2$ and maximally real at every point. Then
every compact subset $K$ of $D$ is CR-, $L^{ \sf p }$- and $\mathcal{
W}$-removable.
\end{theorem} 

We would like to point out that, seeking theorems without any
assumption of pseudoconvexity leads to substantial open problems,
because one loses almost all of the strong interweavings between
function-theoretic tools and geometric arguments which are valid in
the pseudoconvex realm, for instance: Hopf Lemma, plurisubharmonic
exhaustions, envelopes of function spaces, local maximum modulus
principle, Stein neighborhood basis, {\it etc.}

\smallskip

We sketch the proof of the theorem. We first claim that $M \backslash
K$ is (also) globally minimal. Indeed, if there were a
lower-dimensional orbit $\mathcal{ O}$ of $M \backslash K$, we would
obtain a lower-dimensional orbit of $M$ by adding all characteristic
arcs intersecting $\overline{ \mathcal{ O }}$ (\cite{po2003}, Lemma 1;
\cite{ mp2006a}, Lemma~3.5). Then by Theorem~4.12(V), continuous CR
functions on $M \backslash K$ extend holomorphically to a one-sided
neighborhood $\mathcal{ V}^b ( M \backslash K)$.

For later application of the continuity principle, similarly as
in~\cite{mp2002, po2003, mp2006a}, we deform $M \backslash K$ in
$\mathcal{ V}^b ( M \backslash K)$, so that the functions are
holomorphic in some ambient neighborhood $\mathcal{ U}$ of $M
\backslash K$ in $\C^2$.

The first key idea is to construct an embedded 2-sphere containing a
neighborhood of $K$ in $D$ and to apply the filling Theorem~5.17. This
will give us a Levi flat $3$-ball foliated by analytic discs, which by
translations, will enable us to fill in a one-sided neighborhood of
$K$.

In the case where $M = \partial \Omega$ is a strictly pseudoconvex
boundary, the construction of the $2$-sphere is quite direct: we pick
an open $2$-disc $D'$ having $\mathcal{ C}^\infty$ boundary $\partial
D' \simeq S^1$ with $K \subset D' \Subset D$; translating it slightly
and smoothly within $\partial \Omega$, we obtain an almost parallel
copy $D'' \subset \partial \Omega$; then we construct the $2$-sphere
$S'$ by gluing (inside $\partial \Omega$) a thin closed strip $\simeq
[ - \varepsilon, \varepsilon] \times S^1$ to $\partial D' \simeq S^1$
and to $\partial D'' \simeq S^1$; finally, we perturb the strip part
of $S'$ in a generic way to assure that $S'$ has only (a finite number
of) isolated complex tangencies of elliptic or of 
hyperbolic type\footnote{
Observe that since $D$ is totally real, the last step can be done
without changing $S'$ along $D'$.}. Then Theorem~5.17 yields a
Levi-flat 3-ball $B' \subset \Omega$ with $\partial B' = S'$.

If $M$ is not strongly pseudoconvex, the filling of $S'$ by a
Levi-flat ball $B'$ may fail, because of a known counter-example
\cite{ fm1995}. As a trick, we modify the construction. Using the
fact that the squared distance function ${\rm dist} \, ( \cdot, D'
)^2$ is strictly plurisubharmonic in a neighborhood of $\overline{
D}'$ (by total reality), for $\varepsilon >0$ small, the sublevel sets
\[
\Omega_\varepsilon' 
:= 
\big\{q\in\C^2:\, 
{\rm dist}\,
\big(
q,\overline{D'}
\big)
<
\varepsilon\big\}
\]
are strongly pseudoconvex neighborhoods of $\overline{ D' }$
intersecting $M$ transversally along the 2-spheres $\partial
\Omega_\epsilon' \cap M$.

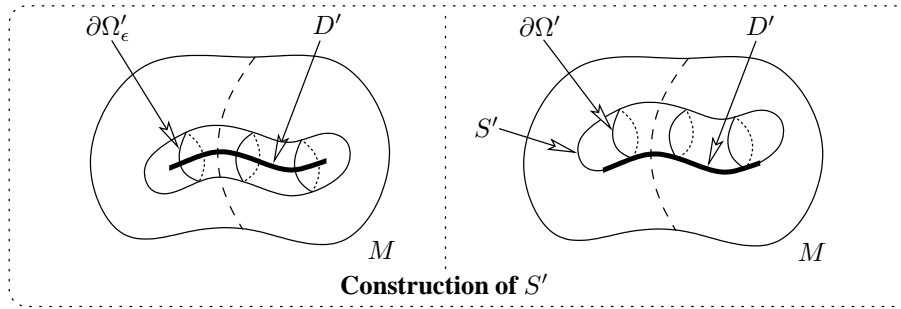
\begin{center}
\begin{picture}(0,0)%
\includegraphics{wurst.pstex}%
\end{picture}%
\setlength{\unitlength}{4144sp}%
\begingroup\makeatletter\ifx\SetFigFont\undefined
\def\x#1#2#3#4#5#6#7\relax{\def\x{#1#2#3#4#5#6}}%
\expandafter\x\fmtname xxxxxx\relax \def\y{splain}%
\ifx\x\y   
\gdef\SetFigFont#1#2#3{%
  \ifnum #1<17\tiny\else \ifnum #1<20\small\else
  \ifnum #1<24\normalsize\else \ifnum #1<29\large\else
  \ifnum #1<34\Large\else \ifnum #1<41\LARGE\else
     \huge\fi\fi\fi\fi\fi\fi
  \csname #3\endcsname}%
\else
\gdef\SetFigFont#1#2#3{\begingroup
  \count@#1\relax \ifnum 25<\count@\count@25\fi
  \def\x{\endgroup\@setsize\SetFigFont{#2pt}}%
  \expandafter\x
    \csname \romannumeral\the\count@ pt\expandafter\endcsname
    \csname @\romannumeral\the\count@ pt\endcsname
  \csname #3\endcsname}%
\fi
\fi\endgroup
\begin{picture}(5424,1824)(439,-1963)
\put(2265,-319){\makebox(0,0)[lb]{\smash{\SetFigFont{10}{12.0}{rm}{\color[rgb]{0,0,0}$D'$}%
}}}
\put(2594,-1666){\makebox(0,0)[lb]{\smash{\SetFigFont{10}{12.0}{rm}{\color[rgb]{0,0,0}$M$}%
}}}
\put(910,-323){\makebox(0,0)[lb]{\smash{\SetFigFont{10}{12.0}{rm}{\color[rgb]{0,0,0}$\partial \Omega_\epsilon'$}%
}}}
\put(2431,-1876){\makebox(0,0)[lb]{\smash{\SetFigFont{10}{12.0}{rm}{\color[rgb]{0,0,0}{\bf Construction of $S'$}}%
}}}
\put(4857,-327){\makebox(0,0)[lb]{\smash{\SetFigFont{10}{12.0}{rm}{\color[rgb]{0,0,0}$D'$}%
}}}
\put(5168,-1681){\makebox(0,0)[lb]{\smash{\SetFigFont{10}{12.0}{rm}{\color[rgb]{0,0,0}$M$}%
}}}
\put(3232,-935){\makebox(0,0)[lb]{\smash{\SetFigFont{10}{12.0}{rm}{\color[rgb]{0,0,0}$S'$}%
}}}
\put(3502,-331){\makebox(0,0)[lb]{\smash{\SetFigFont{10}{12.0}{rm}{\color[rgb]{0,0,0}$\partial \Omega'$}%
}}}
\end{picture}

\end{center} 

Furthermore, a given fixed $\Omega_\epsilon'$ can be slightly isotoped
(translated) to a domain $\Omega'$ still strongly pseudoconvex and
having boundary transverse to $M$ so that $D'$ is precisely contained
in the isotoped 2-sphere $\partial \Omega' \cap M$. After a very
slight generic perturbation, we may insure that $S'$ has only elliptic
or hyperbolic complex tangencies (a part of $\partial \Omega'$ has
also to be perturbed). In sum:

\def\thelemma{5.21}\begin{lemma}
{\rm (\cite{ po2003})}
There exists a bounded domain $\Omega' \subset \C^2$ such that{\rm :}

\begin{itemize}

\smallskip\item[$\bullet$]
$\partial \Omega'$ is $\mathcal{ C}^\infty$, 
strongly pseudoconvex and diffeomorphic to 
a $3$-sphere{\rm ;}

\smallskip\item[$\bullet$]
$\partial \Omega'$ intersects $M$ transversally in a
two-sphere $S' := \partial \Omega' \cap M${\rm ;}

\smallskip\item[$\bullet$]
$S'$ has $k$ hyperbolic and $k+2$ elliptic points{\rm ;}

\smallskip\item[$\bullet$]
$\partial \Omega'$ contains the open $2$-disc $D' \supset K$.

\end{itemize}\smallskip
\end{lemma}

Then Theorem~5.17 applies in the strongly pseudoconvex boundary
$\partial \Omega'$, yielding a Levi-flat $3$-sphere $B' \subset
\Omega'$ with $\partial B' =S'$. However, the nonpseudoconvexity of
$M$ obstructs further insights in the position of $B'$ with respect to
$M$. In fact, $B'$ may change sides or even be partly contained in
$M$.

In the (simpler) case where $M = \partial \Omega$ is a strongly
pseudoconvex boundary, we introduce a foliation of a neighborhood of
$S'$ in $M$ by $\mathcal{ C}^\infty$ $2$-spheres $S_t'$ with $S_0'
=S'$. By filling them, we get a family of Levi-flat 3-balls $B_t'$
with $\partial B_t' =S_t'$. Denote $B_t' = \cup_s \, \Delta_{ t, s}'$
the foliation of $B_t'$ by holomorphic discs. For $t \neq 0$, each
$\Delta_{ t, s}'$ has boundary $\partial \Delta_{ t, s}' \subset S_t'
\subset M \backslash K$. Thus, by means of the continuity principle,
we may extend holomorphic functions in the neighborhood $\mathcal{ U}$
of $M \backslash K$ to a neighborhood of $B_t'$ in $\C^n$, for all
small $t\neq 0$. A final simple check shows that Theorem~2.30 {\bf
(rm5)} applies to remove $B_0'$, and we get holomorphic extension to the
union $\cup_t \, B_t'$, a set containing the strongly pseudoconvex
open local side of $\Omega$ at every point of $K$.

Without pseudoconvexity assumption on $M$, 
we can still consider a foliation
$S_t'$, but now the global geometry of $B_t'$ is no longer clear. If for
instance $M$ is Levi-flat near $K$ and the $S_t'$ are contained in the
Levi-flat part, then the $B_t'$ just form an increasing family whose
union is just a subdomain of $M$. Therefore it seems necessary to
deform $S'$ once again in order to gain transversality of $B'$ and $M$.
Since the global behavior of $B'$ is hard to control, a further
localization is advisable.

As in \cite{ me1997}, we consider the set $K_{\rm nr}$ of points $q
\in K$ such that $\mathcal{ O} \big( \mathcal{ V} (M \backslash K)
\big)$ does not extend holomorphically to a one-sided neighborhood of
$q$. So $\mathcal{ O} \big( \mathcal{ V} (M \backslash K) \big)$
extends holomorphically to a one-sided neighorhood $\mathcal{ V}^b
\big( K \backslash K_{\rm nr} \big)$. By deforming $M$ at points of
$K \backslash K_{ \rm nr}$, we come down to the same situation
with $K$ replaced with $K_{ \rm nr}$, except that no point
of $K_{ \rm nr}$ should be removable. Assuming $K_{ \rm nr} \neq
\emptyset$, to conclude by contradiction, it then
suffices to remove only one point of $K_{ \rm nr}$.

To begin with, assume that $K_{\rm nr}$ is contained in {\sl finitely}
many of the disc boundaries $\partial \Delta_{ 0, s}'$ which foliate
$S' = S_0'$. Then we claim that no $\partial \Delta_{ 0, s }'$ can be
contained in $K_{\rm nr}$. Otherwise, $\partial \Delta_{ 0, s }'
\subset K_{ \rm nr} \subset D' \subset D$ and the $2$-disc enclosed by
$\partial \Delta_{ 0, s }$ in $S_0'$ inside the totally real $2$-disc
$D'$ contain no complex tangencies, but the filling provided by
Theorem~7.17 excludes such a topological possibility.
So $K_{\rm nr}$ is properly
contained in a finite union of arcs, and hence removable by 
Theorem~4.9.

Therefore we may assume that $K_{\rm nr}$ has nonvoid intersection
with {\sl infinitely} many of $\partial \Delta_{ 0, s}'$. Since there
is only finitely many complex tangencies, there exists a $\partial
\Delta_{0, s_0}'$ with $\partial \Delta_{0, s_0}' \cap K_{ \rm nr}
\neq \emptyset$ not encountering them. The same argument as above
shows that $\partial \Delta_{0, s_0}' \not \subset K_{ \rm nr}$. Let
$p_0 \in \partial \Delta_{0, s_0}' \cap K_{ \rm nr}$.

If $\overline{ \Delta_{ 0, s_0}'}$ and $M$ meet {\sl transversally}
at $p_0'$, holomorphic extension to a one-sided neighborhood at
$p_0'$ proceeds as in the strongly pseudoconvex case, by
applying the continuity principle with discs $\Delta_{ t, s}'
\subset B_t'$ for $t\neq O$.

Assume now that $\overline{ \Delta_{ 0, s_0}'}$ is tangential to $M$
in $p_0'$ or equivalently, that $\partial \Delta_{0, s_0}'$ is
tangential to the characteristic leaf in $p_0'$. The idea is to change
the angle of the discs close to $\Delta_{0, s_0 }'$, and to apply the
above argument to the deformed disc passing through $p_0'$. Since
$\partial \Delta_{ 0, s_0}' \not \subset K_{ \rm nr}$, we may deform
slightly $S'$ near some point $q_0' \in \partial \Delta_{0, s_0}'
\backslash K_{ \rm nr}$ in the direction normal to $B'$. More
precisely, one deforms $S'$ slightly, so that Theorem 5.17 still
applies, and then picks up the disc of the deformed Levi-flat 3-ball
that passes through $p_0'$. In view of known results about normal
deformations of small discs (Proposition~2.21(V); \cite{ trp1990,
brt1994, tu1994a}), the turning of the angle for large discs (\cite{
fo1986, gl1994}) may also be established in such a way ({\it
see}~\cite{ po2003, po2004}).

There is one final point to be handled carefully. We have to be sure
that after turning the discs, the deformed disc boundary passing
through the point $p_0' \in K_{ \rm nr}$ is not entirely contained in
$K_{ \rm nr}$. 

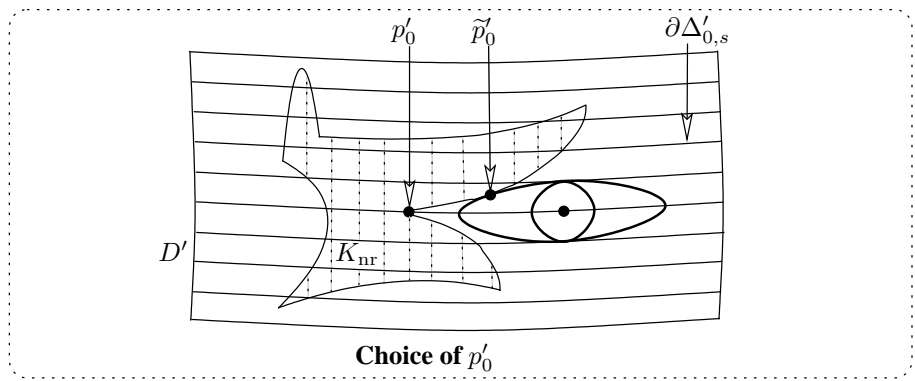
\begin{center}
\begin{picture}(0,0)%
\includegraphics{bazille.pstex}%
\end{picture}%
\setlength{\unitlength}{4144sp}%
\begingroup\makeatletter\ifx\SetFigFont\undefined
\def\x#1#2#3#4#5#6#7\relax{\def\x{#1#2#3#4#5#6}}%
\expandafter\x\fmtname xxxxxx\relax \def\y{splain}%
\ifx\x\y   
\gdef\SetFigFont#1#2#3{%
  \ifnum #1<17\tiny\else \ifnum #1<20\small\else
  \ifnum #1<24\normalsize\else \ifnum #1<29\large\else
  \ifnum #1<34\Large\else \ifnum #1<41\LARGE\else
     \huge\fi\fi\fi\fi\fi\fi
  \csname #3\endcsname}%
\else
\gdef\SetFigFont#1#2#3{\begingroup
  \count@#1\relax \ifnum 25<\count@\count@25\fi
  \def\x{\endgroup\@setsize\SetFigFont{#2pt}}%
  \expandafter\x
    \csname \romannumeral\the\count@ pt\expandafter\endcsname
    \csname @\romannumeral\the\count@ pt\endcsname
  \csname #3\endcsname}%
\fi
\fi\endgroup
\begin{picture}(5424,2229)(432,-1985)
\put(2393,-1272){\makebox(0,0)[lb]{\smash{\SetFigFont{10}{12.0}{rm}{\color[rgb]{0,0,0}$K_{\rm nr}$}%
}}}
\put(2514,-1890){\makebox(0,0)[lb]{\smash{\SetFigFont{10}{12.0}{rm}{\color[rgb]{0,0,0}{\bf Choice of $p_0'$}}%
}}}
\put(2732, 67){\makebox(0,0)[lb]{\smash{\SetFigFont{10}{12.0}{rm}{\color[rgb]{0,0,0}$p_0'$}%
}}}
\put(3210, 63){\makebox(0,0)[lb]{\smash{\SetFigFont{10}{12.0}{rm}{\color[rgb]{0,0,0}$\widetilde{p}_0'$}%
}}}
\put(1336,-1273){\makebox(0,0)[lb]{\smash{\SetFigFont{10}{12.0}{rm}{\color[rgb]{0,0,0}$D'$}%
}}}
\put(4366, 74){\makebox(0,0)[lb]{\smash{\SetFigFont{10}{12.0}{rm}{\color[rgb]{0,0,0}$\partial\Delta_{0,s}'$}%
}}}
\end{picture}

\end{center}

This can be assured by replacing $p_0'$ by another special
nearby point $\widetilde{ p}'_0\in K_{\rm nr}$
with a good transversality property as illustrated above.
\qed

\smallskip

Theorem~5.20 is not yet the complete generalization of
Theorem~5.15 to nonpseudoconvex hypersurfaces, since $D$ is assumed to
be totally real at every point. If $D$ has hyperbolic complex
tangencies, it is not clear whether a sphere $S'$ together with a
strongly pseudoconvex boundary $\partial \Omega' \supset S'$ as in the
above key lemma can be constructed. The recent Theorem~5.13 indicates
that this is possible if hyperbolic complex tangencies are
holomorphically flat, an assumption which would be rather {\it ad
hoc}\, for the removal of compact sets $K \subset D$.

In fact, assuming generally that $M$ is an arbitrary globally minimal
hypersurface, that a given surface $S \subset M$ has arbitrary
topology (not necessarily diffeomorphic to an open $2$-disc) and
possesses complex tangencies, the reduction to the filling
Theorem~5.17 seems to be impossible. Indeed, Forn{\ae}ss-Ma (\cite{
fm1995}) constructed an unknotted nonfillable $2$-sphere $S \subset
\C^2$ having only two elliptic complex tangencies. To the authors'
knowledge, the possibility of filling by Levi-flat $3$-spheres some
$2$-spheres lying in a {\it nonpseudoconvex}\, hypersurface is a
delicate open problem. In addition, for the higher codimensional
generalization of Theorem~1.2, the idea of global filling seems to be
irrelevant at present times, because no analog of the filling
Theorem~5.17 is known in dimension $n \geqslant 3$.

\subsection*{ 5.22.~Beyond this survey}
In the research article~\cite{ mp2006a} placed in direct continuation
to this survey, we consider surfaces $S$ having arbitrary topology and
we generalize Theorem~5.20 to arbitrary codimension, {\it
localizing}\, the removability arguments and using only {\it small}\,
analytic discs.

\vfill\end{document}